\RequirePackage{rotating}
\documentclass[a4paper]{amsart}

\makeatletter
\renewcommand\part{%
   \if@noskipsec \leavevmode \fi
   \par
   \addvspace{4ex}%
   \@afterindentfalse
   \secdef\@part\@spart}

\def\@part[#1]#2{%
    \ifnum \c@secnumdepth >\m@ne
      \refstepcounter{part}%
      \addcontentsline{toc}{part}{\thepart\hspace{1em}#1}%
    \else
      \addcontentsline{toc}{part}{#1}%
    \fi
    {\parindent \z@ \raggedright
     \interlinepenalty \@M
     \normalfont
\center{
     \ifnum \c@secnumdepth >\m@ne
       \large\bfseries \partname\nobreakspace\thepart: 
     \fi
     \large \bfseries #2}%
     \par}%
    \nobreak
    \vskip 3ex
    \@afterheading}
\def\@spart#1{%
    {\parindent \z@ \raggedright
     \interlinepenalty \@M
     \normalfont
     \huge \bfseries #1\par}%
     \nobreak
     \vskip 3ex
     \@afterheading}
\makeatother

\usepackage{amsmath,amstext,amssymb,mathrsfs,amscd,amsthm,indentfirst,rotating}
\usepackage{amsfonts}
\usepackage{mathtools}
\usepackage{stmaryrd}

\usepackage[sort=def,symbols,nomain,numberline]{glossaries}
\newcommand{\listofsymbolsname}{List of Symbols}
\makeglossaries

\usepackage{xr-hyper}
\usepackage[dvipsnames,svgnames,x11names,hyperref]{xcolor}
\usepackage[colorlinks,citecolor=Mahogany,linkcolor=Mahogany,urlcolor=Mahogany,filecolor=Mahogany]{hyperref}
\usepackage{bbm}

\usepackage{booktabs}
\usepackage{verbatim}
\usepackage{enumerate}
\usepackage{nicefrac}
\usepackage{lmodern}
\usepackage{microtype}
\usepackage{setspace}
\usepackage{url}

\usepackage{tikz,tikz-cd}
\usetikzlibrary{matrix,calc,positioning,arrows,decorations.pathreplacing,patterns,fit}
\usepackage{graphicx}

\newcommand\mvline[3][]{%
	\pgfmathtruncatemacro\hc{#3-1}
	\draw[#1]({$(#2-1-#3)!.5!(#2-1-\hc)$} |- #2.north) -- ({$(#2-1-#3)!.5!(#2-1-\hc)$} |- #2.south);
}
\newcommand\mhline[4][]{%
	\node[fit=(#2-#3-1),inner sep=0pt,outer sep=0pt](R){};
	\foreach \i in {1,...,#4}\node[fit=(R) (#2-#3-\i),inner sep=0pt,outer sep=0pt](R){};
	\draw[#1] (R.north -| #2.west) -- (R.north -| #2.east);
}

\numberwithin{equation}{section}

\newcommand{\Kleis}{\cat{Kleis}}

\newcommand{\cat}[1]{\mathsf{#1}}
\newcommand{\mr}[1]{{\rm #1}}
\newcommand{\bunit}{\mathbbm{1}}
\newcommand{\fS}{\mathfrak{S}}
\newcommand{\bterm}{\mathbbm{t}}
\newcommand{\binit}{\mathbbm{i}}
\newcommand{\sfG}{\mathsf{G}}
\newcommand{\sfH}{\mathsf{H}}
\newcommand{\sfC}{\mathsf{C}}
\newcommand{\sfCstar}{\sfC_*} 
\newcommand{\op}{\mathrm{op}}
\newcommand{\suspsplit}{T^{E_1}}
\renewcommand{\split}{S^{E_1}}
\newcommand{\splitinf}{S^{E_\infty}}

\newcommand{\Modk}{\cat{Mod}_{\bk}}

\newcommand{\indec}{Q}
\newcommand{\dec}{\mr{Dec}}
\newcommand{\ul}[1]{\underline{#1}}
\newcommand{\Hom}{\cH om}
\newcommand{\Sing}{\mathrm{Sing}}

\makeatletter
\newcommand{\boxedsymbol}[2]{%
	\begingroup
	\setlength{\fboxsep}{0pt}%
	\fbox{%
		$\m@th#1\mspace{-1.25mu}#2\mspace{-1.25mu}$%
	}%
	\endgroup
}
\makeatother

\DeclareRobustCommand\longtwoheadleftarrow
{\twoheadleftarrow\joinrel\relbar}

\makeatletter
\newcommand*{\obigwedge}{%
	\mathbin{%
		\mathpalette\@obigwedge{}%
	}%
}
\newcommand*{\@obigwedge}[2]{%
	\sbox0{$#1\bigoplus\m@th$}%
	\dimen2=.5\dimexpr\wd0-\ht0-\dp0\relax 
	\dimen@=\dimexpr\ht0+\dp0\relax
	\def\lw{.04}
	\def\radius{.5-\lw/2}%
	\kern\dimen2 
	\tikz[
	line width=\lw\dimen@,
	line join=round,
	x=\dimen@,
	y=\dimen@,
	baseline=\dimexpr-.5\dimen@+\dp0\relax,
	]
	\draw
	(0,0) circle[radius=\radius]
	(225:\radius) -- (0,.5-\lw) -- (-45:\radius)
	;%
	\kern\dimen2 
}
\makeatother

\newcommand{\CircNum}[1]{\ooalign{\hfil\raise .00ex\hbox{\scriptsize #1}\hfil\crcr\mathhexbox20D}}

\newcommand{\bk}{\mathbbm{k}}

\newcommand{\bF}{\mathbb{F}}

\newcommand{\bH}{\mathbb{H}}

\newcommand{\bL}{\mathbb{L}}

\newcommand{\bN}{\mathbb{N}}

\newcommand{\bQ}{\mathbb{Q}}
\newcommand{\bR}{\mathbb{R}}
\newcommand{\bS}{\mathbb{S}}

\newcommand{\bZ}{\mathbb{Z}}

\newcommand{\gA}{\bold{A}}
\newcommand{\gB}{\bold{B}}
\newcommand{\gC}{\bold{C}}
\newcommand{\gD}{\bold{D}}
\newcommand{\gE}{\bold{E}}
\newcommand{\gF}{\bold{F}}

\newcommand{\gL}{\bold{L}}
\newcommand{\gM}{\bold{M}}
\newcommand{\gN}{\bold{N}}

\newcommand{\gR}{\bold{R}}
\newcommand{\gS}{\bold{S}}
\newcommand{\gT}{\bold{T}}

\newcommand{\gV}{\bold{V}}

\newcommand{\gX}{\bold{X}}
\newcommand{\gY}{\bold{Y}}
\newcommand{\gZ}{\bold{Z}}

\newcommand{\cA}{\mathcal{A}}
\newcommand{\cB}{\mathcal{B}}
\newcommand{\cC}{\mathcal{C}}

\newcommand{\cE}{\mathcal{E}}

\newcommand{\cH}{\mathcal{H}}
\newcommand{\cI}{\mathcal{I}}

\newcommand{\cO}{\mathcal{O}}
\newcommand{\cP}{\mathcal{P}}

\newcommand{\cR}{\mathcal{R}}

\newcommand{\cV}{\mathcal{V}}

\newcommand{\cX}{\mathcal{X}}
\newcommand{\cY}{\mathcal{Y}}

\newcommand\lra{\longrightarrow}

\newcommand\Emb{\mathrm{Emb}}

\newcommand\const{\operatorname*{const}}
\newcommand\colim{\operatorname*{colim}}
\newcommand\hocolim{\operatorname*{hocolim}}
\newcommand\holim{\operatorname*{holim}}
\newcommand\limone{\operatorname*{lim^1}}

\newcommand{\hcoker}{/\!\!/}

\newcommand{\grr}{\mathrm{gr}}
\newcommand{\alg}{\mathrm{alg}}
\newcommand{\Dec}{\mathrm{Dec}}

\newcommand{\fgr}[1]{\| {#1} \|}
\newcommand{\gr}[1]{| {#1} |}

\newcommand{\N}{\mathbb{N}}
\newcommand{\Z}{\mathbb{Z}}

\renewcommand{\epsilon}{\varepsilon}

\newcommand{\id}{\mathrm{id}}

\newcommand{\Map}{\mathrm{Map}}

\newcommand{\Ob}{\mathrm{Ob}}

\newcommand{\nin}{\not\in}

\newcommand{\sfS}{\mathsf{S}}
\newcommand{\sfD}{\mathsf{D}}
\newcommand{\Alg}{\cat{Alg}}

\mathchardef\ordinarycolon\mathcode`\:
\mathcode`\:=\string"8000
\begingroup \catcode`\:=\active
  \gdef:{\mathrel{\mathop\ordinarycolon}}
\endgroup

\usepackage{amsthm}
\theoremstyle{plain}

\newtheorem{theorem}{Theorem}[section]

\newtheorem{proposition}[theorem]{Proposition}
\newtheorem{lemma}[theorem]{Lemma}

\newtheorem{corollary}[theorem]{Corollary}

\theoremstyle{definition}
\newtheorem{definition}[theorem]{Definition}

\newtheorem{notation}[theorem]{Notation}

\newtheorem{assum}[theorem]{Assumption}
\newtheorem{axiom}[theorem]{Axiom}

\theoremstyle{remark}

\newtheorem{remark}[theorem]{Remark}
\newtheorem{example}[theorem]{Example}
\newtheorem*{remark*}{Remark}

\title{Cellular $E_k$-algebras}

  \author{S{\o}ren Galatius}
  \address{Department of Mathematics\\
  	University of Copenhagen\\
  	Denmark}
\email{galatius@math.ku.dk}

 \author{Alexander Kupers}
  \address{Department of Computer and Mathematical Sciences\\
 	University of Toronto Scarborough \\
 	1265 Military Trail \\
 	Toronto, ON M1C 1A4 \\ Canada}
 \email{a.kupers@utoronto.ca}

 \author{Oscar Randal-Williams}
  \address{Centre for Mathematical Sciences\\
 	Wilberforce Road\\
 	Cambridge CB3 0WB\\
 	UK}
 \email{o.randal-williams@dpmms.cam.ac.uk}

\date{\today}
\subjclass[2020]{57R90, 18N40, 55R40, 55P35, 55P48, 57T30}
\keywords{$E_k$-algebras, cellular approximation, indecomposables, bar construction, homological stability}

\makeindex

\allowdisplaybreaks

\begin{document}
\begin{abstract}
We give a set of foundations for cellular $E_k$-algebras which are especially convenient for applications to homological stability. We provide conceptual and computational tools in this setting, such as filtrations, a homology theory for $E_k$-algebras with a Hurewicz theorem, CW approximations, and many spectral sequences, which shall be used for such applications in future papers.
\end{abstract}

\maketitle

\addtocontents{toc}{\protect\setcounter{tocdepth}{1}}
\tableofcontents

\section{Introduction}\label{sec:introduction}In this paper and its sequels we develop a \emph{multiplicative} approach to the study of automorphism groups $\mr{G}_n$ such as mapping class groups, automorphism groups of free groups, or general linear groups. The multiplicative structure on these groups is given by homomorphisms $\mr{G}_n \times \mr{G}_m \to \mr{G}_{n+m}$ which are appropriately associative and commutative. This can be made precise by saying that the disjoint union of classifying spaces $\gR \coloneqq \bigsqcup_{n \geq 1} B\mr{G}_n$ has the structure of a (non-unital) $E_k$-algebra, usually for $k=2$ or $\infty$.

This perspective on the groups $\mr{G}_n$ and their homology is fundamentally different than the traditional \emph{additive} approach. That approach focuses on homomorphisms $\mr{G}_n \to \mr{G}_{n+1}$ which induce \emph{stabilization maps} $B\mr{G}_n \to B\mr{G}_{n+1}$ on classifying spaces. This is akin to thinking of the space $\gR$ as a module over the monoid of natural numbers under addition, with $k \in \bN$ sending $B\mr{G}_n$ to $B\mr{G}_{n+k}$ by iterating the stabilization map. A typical result of the additive approach is \emph{homological stability}; the statement that the stabilization map $B\mr{G}_n \to B\mr{G}_{n+1}$ induces an isomorphism on the $d$th homology group when $d \ll n$.

Our multiplicative approach recovers many such homological stability results. However, it is also capable of producing qualitatively different results, such as what we call \emph{secondary homological stability} (or \emph{non-stability}) results. Our strategy for proving these results is to construct $\gR$ or related objects out of free $E_k$-algebras in a manner analogous to CW approximation. We use a homology theory for $E_k$-algebras to bound how many $E_k$-cells are needed. Such bounds together with explicit knowledge of the homology of free $E_k$-algebras are then used to deduce results about the homology of $\gR$. 

With an eye towards implementing this strategy, in this first paper we provide a robust set of foundations for a cellular theory of $E_k$-algebras. Later papers shall focus on the geometric and algebraic arguments relevant to particular examples.

\subsection{$E_k$-algebras} We start by explaining the notion of an $E_k$-algebra in a sufficiently nice category $\sfS$, e.g.~the category of simplicial sets, spectra, or simplicial $\bk$-modules. To talk about $E_k$-algebras in $\sfS$, it must be copowered over simplicial sets and have a monoidal structure. This monoidal structure needs to be braided if we wish to discuss $E_2$-algebras and symmetric if we wish to discuss $E_k$-algebras with $k > 2$. 

\begin{figure}[h]
	\centering
	\begin{tikzpicture}
	
	\draw[fill = Mahogany!10!white] (3,2) rectangle (4,3);
	\draw[fill = Mahogany!10!white] (1,2) rectangle (2,3);
	
	\draw (0,0) rectangle (5,5);
	
	\node at (1.5,2.5) {$e_1$};
	\node at (3.5,2.5) {$e_2$};
	
	\node[right] at (5.1,2.5) {$\in \cC_2(2)$};
	\end{tikzpicture}
	\caption{An element of $\cC_2(2)$.}
	\label{fig:introduction-c2}\end{figure}

Next we need to pick an $E_k$-operad in simplicial sets: we use the (singular simplicial set of the) little $k$-cubes operad $\cC_k$. The space of $i$-ary operations in this operad, $\cC_k(i)$, is the space of rectilinear embeddings of $i$ cubes $I^k$ into $I^k$ with disjoint interior, and composition of operations is given by composition of embeddings. 

Operads encode algebraic structures: a (unital) \emph{$E_k$-algebra} $\gR$ in $\sfS$ is an object $R$ of $\sfS$ equipped with maps
\[\cC_k(n) \times R^{\otimes n} \lra R\]
satisfying unit, associativity, and equivariance axioms. Here $\times$ denotes the copowering over simplicial sets and $\otimes$ denotes the monoidal product on $\sfS$.
For example, Figure \ref{fig:introduction-c2} depicts an element $\mu$ of $\cC_2(2)$ which encodes a particular operation which ``multiplies'' two elements of $E_2$-algebra. This operation $R \otimes R \to R$ is commutative in the sense that it is homotopic to the map with input switched because $\mu \in \cC_2(2)$ may be connected by a path to $\mu$ with labels $e_1$ and $e_2$ switched. There is no canonical choice of such a path and the existence of non-trivial families of multiplications gives rise to operations on the homology of $E_k$-algebras.

We write $\Alg_{E_k}(\sfS)$ for the category of $E_k$-algebras in $\sfS$, and $U^{E_k} \colon \Alg_{E_k}(\sfS) \to \sfS$ for the functor which forgets the $E_k$-algebra structure. To do homotopy theory with $E_k$-algebras, we define a morphism in $\Alg_{E_k}(\sfS)$ to be a weak equivalence if it becomes one after applying $U^{E_k}$, i.e.\ after forgetting the $E_k$-algebra structure. 

Let us give two examples of $E_k$-algebras:

\begin{figure}
	\centering
	\begin{tikzpicture}
	\draw[pattern=north west lines, pattern color=black!10!white] (0,0) rectangle (5,5);
	\node at (0.5,0.5) {$x_0$};
	
	\draw[fill = white] (1,.5) rectangle (1.5,3.5);
	\draw[fill = white] (2,2) rectangle (3,3);
	\draw[fill = white] (3.5,3.5) rectangle (4.5,4.5);
	
	\node at (1.25,2) {$f_2$};
	\node at (2.5,2.5) {$f_1$};
	\node at (4,4) {$f_3$};
	
	\node[right] at (5.1,2.5) {$\in \Omega^2 X$};
	\end{tikzpicture}
	\caption{The result of combining three elements $f_1,f_2,f_3 \in \Omega^2 X$ with an element of $\cC_2(3)$.}
	\label{fig:loop-space-example}
\end{figure}

\begin{example}\label{exam:loop-spaces} $E_k$-algebras were introduced to study $k$-fold loop spaces \cite{BoardmanVogt,GILS}, the prototypical examples of $E_k$-algebras. If $X$ is a pointed topological space and $\Omega^k X$ denotes the space of continuous maps $I^k \to X$ sending $\partial I^k$ to the basepoint $x_0$, then we may combine an element $e$ of $\cC_k(i)$ with $i$ elements $f_j$ of $\Omega^k X$ to form a single elements of $\Omega^k X$; insert $f_j$ into the image of the $j$th cube $e_j$, and extend by the constant map with value $x_0$. Section \ref{sec:group-completion-map} contains more details, and Figure \ref{fig:loop-space-example} gives an example in the case $k=2$ and $j=3$.
\end{example}

\begin{example}The example $\bigsqcup_{n \geq 1} B\mr{G}_n$ discussed above arises as follows: if $\sfG$ is a braided or symmetric monoidal groupoid whose objects are given by $\bN$ and such that the monoidal product of objects is addition, then the disjoint union $\bigsqcup_{n \geq 1} B\mr{G}_n$ of the automorphism groups $\mr{G}_n$ of the objects of $\sfG$ has a canonical structure of an $E_2$- or $E_\infty$-algebra. Section \ref{sec:EkAlgFromGpd} contains the details.
\end{example}

\subsection{Cellular $E_k$-algebras} Cell attachments for $E_k$-algebras are modeled after cell attachments for topological spaces. In the latter case, we start with the data of a topological space $X$ and an attaching map $\partial D^d \to X$. Taking a pushout
\[\begin{tikzcd}\partial D^d \rar \dar & X \dar \\
D^d \rar & X',\end{tikzcd}\]
we get a space $X'$ which we say is obtained by attaching a $d$-dimensional cell to $X$.

If we assume the category $\sfS$ has all colimits then so does $\Alg_{E_k}(\sfS)$, and we can form cell attachments in an analogous manner. That is, we start with the data of an $E_k$-algebra $\gR$, a pair of simplicial sets $(D^d,\partial D^d)$ whose geometric realization is homeomorphic to $(D^d,\partial D^d)$, and an attaching map $\partial D^d \to U^{E_k} \gR$ in $\sfS$. (Here the copowering is used to make sense of a map from a simplicial set to an object of $\sfS$.) Then we have a diagram
\begin{equation}\label{eqn:pushout-uek} \begin{tikzcd}\partial D^d \rar \dar & U^{E_k}(\gR)  \\
D^d & \qquad \end{tikzcd}\end{equation}
in $\sfS$, but its pushout does not in general admit the structure of an $E_k$-algebra.

To remedy this, we use that the forgetful functor $U^{E_k}$ participates in the adjunction
\begin{equation*}
\begin{tikzcd}
\sfS \arrow[shift left=.5ex]{r}{F^{E_k}} & \Alg_{E_k}(\sfS) \arrow[shift left=.5ex]{l}{U^{E_k}},
\end{tikzcd}
\end{equation*}
with left adjoint $F^{E_k}$ given by the free $E_k$-algebra functor. That is, every map $X \to U^{E_k}(\gR)$ gives rise to a unique map $F^{E_k}(X) \to \gR$ of $E_k$-algebras. So instead of taking the pushout of \eqref{eqn:pushout-uek}, we take the pushout of the adjoint diagram
\[\begin{tikzcd}F^{E_k}(\partial D^d) \rar \dar & \gR  \\
F^{E_k}(D^d) &  \end{tikzcd}\]
in the category $\Alg_{E_k}(\sfS)$, which exists because $\sfS$ has all colimits so in particular pushouts. By construction the result will be an $E_k$-algebra, and we say is obtained by attaching a $d$-dimensional cell to $\gR$. This is justified because it satisfies the same universal property in $\Alg_{E_k}(\sfS)$ that an ordinary cell attachment does in $\cat{Top}$.

An object $\gC \in \Alg_{E_k}(\sfS)$ is cellular if it can be obtained from the initial object by a (perhaps transfinite) sequence of such cell attachments. A slightly stronger version of this is a CW object, in which the cells are in particular attached in order of dimension. Just as in the category of topological spaces, an $\gR \in \Alg_{E_k}(\sfS)$ admits a CW approximation under mild conditions on both $\sfS$ and $\gR$, i.e.~a weak equivalence $\gC \overset{\sim}\to \gR$ in $\Alg_{E_k}(\sfS)$ where $\gC$ is a CW object.

Relative cellular algebras have always played an important role in the study of the homotopy theory of algebras over an operad, as they are the cell complexes in the projective model structure, e.g.\ \cite[Section 4]{BergerMoerdijk2}. They were applied to study homological stability in the context of factorization homology in \cite{KupersMiller}.

\subsection{$E_k$-homology} We want to obtain small CW approximations, i.e.\ ones with as few cells of each dimension as possible. When finding a CW approximation of a topological space $X$, a lower bound on the number of cells needed is given in terms of generators of its singular homology groups: if the abelian group $H_d(X)$ is generated by $a$ elements and the torsion subgroup of $H_{d-1}(X)$ is generated by $b$ elements, then no CW approximation has fewer than $a+b$ cells of dimension $d$. Furthermore, if $X$ is simply connected then this bound is realized. 

The analogous question for $E_k$-algebras has a similar answer, given in terms of a type of homology theory for algebras going back to Hochschild, Quillen, and Andr\'e \cite{HochschildCohomology,QuillenHomology,AndreHomology}: there are homology groups $H_*^{E_k}(\gR)$ which always give a lower bound on the number of cells in a cellular approximation, and under certain assumptions on both $\sfS$ and $\gR$ this bound may be realized. These $E_k$-homology groups are defined as the homology groups of an object $Q^{E_k}_\bL(\gR)$ of \emph{derived $E_k$-indecomposables of $\gR$}, constructed to satisfy $Q_\bL^{E_k}(F^{E_k}(X)) \simeq X_+$ and preserve homotopy colimits. 

The derived $E_k$-indecomposables may be thought of as measuring generators, relations, syzygies, etc.\ for $\gR$ \emph{as an $E_k$-algebra}, and are obtained by deriving the construction which takes the quotient of $\gR$ by the sub-object obtained by applying all possible operations of $E_k$ to $\gR$ of arity at least $2$. This is a version of $\mr{THH}$ and $\mr{TAQ}$ for the $E_k$-operad (those constructions corresponding to $k=1$ and $k=\infty$ respectively), and is related to factorization homology (also known as higher Hochschild homology or topological chiral homology) and to $k$-fold deloopings. For recent results on these, see e.g.\ \cite{BGR,BM,FrancisThesis,Fresse,FresseZiegenhagen,LivernetRicher,Mandell} and their references.

The technical tool underlying this result is a Hurewicz theorem for $E_k$-homology, see Section \ref{sec:hurewicz-cw}. This says that under certain assumptions on both $\sfS$ and $\gR$, the first non-vanishing relative $E_k$-homology group of a map of $E_k$-algebras coincides with the first non-vanishing relative homology group of that map. This is then used to deduce the existence of minimal CW approximations.

\subsection{Computing $E_k$-homology}

In general it is hard to compute $H_*^{E_k}(\gR)$ and we will often settle for establishing vanishing results in a range. These often take the form of vanishing lines with respect to the naturally present grading of $\gR = \bigsqcup_{n \geq 1} B\mr{G}_n$ by $n$. That is, we consider $\gR$ not as an $E_k$-algebra in $\cat{Top}$, but as an $E_k$-algebra in the category $\smash{\cat{Top}^{\bN}}$ of functors from the non-negative integers to $\cat{Top}$. More generally we study $E_k$-algebras in the category $\sfS^\cat{G}$ of functors $\cat{G} \to \sfS$, where $\sfS$ is a sufficiently nice model category. If $\cat{G}$ has monoidal structure then one can make sense of $E_1$-algebras in $\sfS^\cat{G}$, and if it is braided or symmetric monoidal then one can define $E_2$-algebras or $E_k$-algebras for any $k$. In this setting our cells take the form $\partial D^d \times \sfG(g,-) \to D^d \times \sfG(g,-)$, where $\sfG(g,-) \colon \sfG \to \sfS$ denotes the functor represented by an object $g \in \sfG$ considered as an object of $\sfS$ using the copowering. For $\sfG = (\bN, +, 0)$, this just endows the cells with an additional grading. Consequently each cell has both a \emph{geometric dimension} $d$ and a \emph{rank} $[g] \in \pi_0(\sfG)$, and we can ask for CW approximations with few cells in each bidegree $([g],d)$. The $E_k$-homology will be bigraded, with $([g],d)$-cells measured by $H^{E_k}_{g,d}(\gR)$.

To establish vanishing results for $E_k$-homology, we prove that it can be computed using a $k$-fold bar construction, see Section \ref{sec:BarConstr}. This is an observation going back to \cite{GetzlerJones}, and instances of this result appear in \cite{BM,Fresse,FrancisPaper}. It is of particular use in the examples we study in subsequent papers. In these applications, it is often  easier to compute $E_1$-homology, but more convenient to extract information out of the $E_{2}$- or $E_\infty$-homology. Using the description in terms of iterated bar constructions we prove that vanishing lines can be transferred upwards from $E_k$-homology to $E_{k+1}$-homology (they can also be transferred downwards using cellular methods), see Section \ref{sec:Transferring}.

\subsection{Towards applications}\label{sec:applications}

The technology developed in this paper will be applied in the sequels to prove new results about mapping class groups of oriented surfaces \cite{e2cellsII} and general linear groups \cite{e2cellsIII,e2cellsIV}.

The following is an outline of our basic strategy. The groups $\mr{G}_n$ of interest arise as automorphism groups in a braided or symmetric monoidal category $\cat{G}$. For example, the general linear group $\mr{G}_n = \mr{GL}_n(\bF)$ over a field $\bF$ is the automorphism group of  $\bF^n$ in the symmetric monoidal category of finite-dimensional vector spaces over $\bF$. The object $\gR^+ \coloneqq \bigsqcup_{n \geq 0} B\mr{G}_n$ is the derived pushforward $\bL r_*(\underline{*})$ of the terminal functor $\underline{*} \colon \cat{G} \to \cat{sSet}$ sending each object to $\ast$ and each morphism to the identity, along a monoidal ``rank'' functor $r \colon \cat{G} \to \bN_0$, which canonically is a (unital) $E_2$- or $E_\infty$-algebra. 

Using the aforementioned bar constructions, we prove that its derived $E_1$-inde\-com\-posables can be computed in terms of the $E_1$-splitting complexes described in Section \ref{sec:splittingcomplex}. In applications one must prove these are highly-connected using geometric or algebraic techniques particular to the category $\cat{G}$. By transferring vanishing lines up, we obtain vanishing lines for the $E_{2}$- or $E_\infty$-homology of these objects (it is helpful to establish these vanishing lines for the largest possible $k$, as this simplifies later computations). Combined with low-rank low-degree homology computations, we can use the Hurewicz theorem to build a small CW $E_k$-algebra $\gA$ with map $\gA \to \gR$ which is a good approximation in low degrees. 

One can not compute the absolute homology of $\gR$ using $\gA$, but one can compute relative homology of various ordinary and higher stabilization maps with it. These computations are done using skeletal and cell attachment spectral sequences, which we develop in Section \ref{sec:homology+ss}, and their interaction with the homology operations particular to $E_k$-algebras, which we explain in Section \ref{sec:Cohen}. Alternatively, one can use a comparison result proven in Section \ref{sec:comparing-algebra-and-module cells}: given a map $f \colon \gR \to \gS$ of $E_k$-algebras, this compares how to obtain $\gS$ from $\gR$ by attaching $E_k$-algebras cells, with the way to obtain it by attaching $\gR$-module cells. 

To make the application of this general theory as straightforward as possible, in Section \ref{sec:algebras-from-groupoids} we describe a general framework producing $E_k$-algebras from (perhaps braided or symmetric) monoidal groupoids. When their corresponding $E_1$-splitting complexes satisfy a standard connectivity result (a property related to Koszulity of the $E_1$-algebra $\underline{*}$, cf.~Section \ref{sec:koszul}), we give generic homological stability results with both constant and local coefficients in Sections \ref{sec:homological-stability-applications} and \ref{sec:local-coefficients}. These are ``multiplicative" analogues to the ``additive" generic homological stability theorems of \cite{RWW}. 

This is only a basic outline of the strategy used in the sequels to prove homological stability, secondary homological stability, and non-stability results, and the specifics of implementing it depend on the groups studied. For example, to obtain secondary stability for mapping class groups in \cite{e2cellsII} requires a study of the unstable homology of mapping class groups, and for general linear groups as in \cite{e2cellsIII,e2cellsIV} it requires new results about coinvariants of Steinberg modules. In this paper, however, we shall give examples that do not require any additional work: improved stability results for general linear groups over certain Dedekind domains (including $\bZ$) and finite fields, see Sections \ref{sec:outlook}, \ref{sec:FiniteFields} and \ref{sec:FiniteFieldsTwisted}.

\subsection{Guide for the reader} Supposing that the reader is interested in applications of cellular $E_k$-algebras to homological stability, their goal will be to understand Part 4 of this paper. Here we describe what should be picked up from the earlier parts of the paper, supposing that the reader is already categorically and homotopically sophisticated.

Much of the development in Parts \ref{part:category} and \ref{part:homotopy} of this paper is not intrinsically new, and such a reader will find little here to surprise them. They will be able to skim over these parts quite rapidly, though should pay some attention to 
\begin{enumerate}[(i)]
	\item the way we deal with $\sfG$-graded objects (described in Section \ref{sec:functor-cats}), 
	\item the way we deal with filtrations (Section \ref{sec:filtered-objects-1}) and hence cellular and CW objects (Section \ref{sec:cell-attachments}), 
	\item our notation for homology and the various spectral sequences we develop (Section \ref{sec:homology+ss}). 
\end{enumerate} Such a reader should also at least familiarise themselves with the statement of the CW approximation theorem (Theorem \ref{sec:additive-case}), which will be used often, and may want to read the rest of that section to understand the proof.

On a first reading of Part \ref{part:ek} such a reader should cover the definitions of Section \ref{sec:EkAlgAndMod}, but skip the technical Section \ref{sec:constr-adapt}. From the long Section \ref{sec:BarConstr} they should read Section \ref{sec:iterated-bar-def} and understand the statement of Theorem \ref{thm:BarHomologyIndec}, then understand the statement of Theorem \ref{thm:CalcFree} which should be familiar by analogy with the classical situation of $E_k$-algebras in $\cat{Top}$. It is then a good idea to work through Section \ref{sec:Transferring} in detail, to develop a taste of how the tools discussed so far can be deployed. One can then move on to the reminder in Sections \ref{sec:homology-operations} and \ref{sec:RelationsAmongOps} of the natural operations on the homology of $E_k$-algebras, and recall the description of the homology of free $E_k$-algebras given in Theorem \ref{thm.wkfreealg}. Section \ref{sec:homology-ekalg-coproduct} can be omitted, but  Sections \ref{sec:spectralsequences.skel} and \ref{sec:SSStr} should be understood, as they will be used in almost all calculations. The reader can now return to Section \ref{sec:comparing-algebra-and-module cells}, but omit the difficult proof of Theorem \ref{thm:AlgAsMod} in favour of convincing themselves that it is true (at the level of homology) when $\sfG$ is discrete using the results of Section \ref{sec:Cohen}.

\subsection*{Acknowledgments}

AK and SG were supported by the European Research Council (ERC) under the European Union’s Horizon 2020 research and innovation programme (grant agreement No 682922).  SG was also supported by NSF grant DMS-1405001. AK was also supported by the Danish National Research Foundation through the Centre for Symmetry and Deformation (DNRF92) and by NSF grant DMS-1803766. ORW was  supported by EPSRC grant EP/M027783/1, by the ERC under the European Union’s Horizon 2020 research and innovation programme (grant agreement No. 756444), and by a Philip Leverhulme Prize from the Leverhulme Trust.

SG would like to thank J.\ Lurie for helpful discussions about cell structures and $E_k$-homology when this project was in its infancy. We would like to thank I.\ Sierra, C.\ Bernard, and M.\ Fischer for corrections to an earlier version. We would also like to thank the two anonymous referees for their detailed comments.

\addtocontents{toc}{\protect\setcounter{tocdepth}{2}}

\part{Category theory of algebras over a monad} \label{part:category}

In the first part of this paper we will discuss the category theory of algebras over a sifted monad on a category $\sfS$. Most of these monads arise from operads. Our main goals for this part are the definitions of indecomposables, cell attachments, and CW-algebras, in Sections \ref{sec:simplicial-operad-augmentation}, \ref{sec:cell-attachments-monads}, and \ref{sec:cw-algebras}. In doing so, we will set up the machinery which will be used throughout the second and third parts of this paper.

\section{Contexts for category theory}\label{sec:cat-contexts}  To make the theory of sifted monads go through smoothly, the category $\sfS$ in which we work needs to be endowed with certain structures satisfying certain conditions, as discussed in this section. We shall also explain the contexts we shall work for the remainder of this paper.

\subsection{Axioms for categories} \label{sec:axioms-of-cats} We start by describing our axioms for $\sfS$. We shall assume the reader is familiar with elementary category theory, and refer to \cite{MacLane} and \cite{KellyEnriched} for more background material.

\newglossaryentry{sfs}{%
	name={\ensuremath{\sfS}},
	description={Convenient category},
	type=symbols
}

\subsubsection{Simplicial enrichment}\label{sec:simplicial-enrichment}  Our first assumption is that the category $\gls{sfs}$ is enriched in simplicial sets, i.e.\ it has a class of objects, for each pair $X,Y$ of objects has a simplicial set $\mr{Map}_\sfS(X,Y)$ of morphisms from $X$ to $Y$, a composition law for such morphisms, and for $X=Y$ an identity $0$-simplex in $\mr{Map}_\sfS(X,Y)$. This is the definition of a $\cV$-enriched category specialized to $\cV = \cat{sSet}$.\index{simplicial enrichment}

\begin{axiom}\label{axiom:simplicial-enriched} $\sfS$ is simplicially enriched.\end{axiom}

Restricting to the $0$-simplices of the simplicial sets of morphisms we obtain an ordinary category, which by abuse of notation we also denote $\sfS$. We thus regard the simplicial enrichment as an additional structure on this ordinary category.

We want $\sfS$ to be complete and cocomplete in the enriched sense, which means that not only should $\sfS$ be complete and cocomplete in the ordinary sense, but it should also have all $\cat{sSet}$-indexed colimits and limits. This means that for each simplicial set $K$, the functors 
\begin{equation}\label{eqn:indexed-colimits-limits} Y \mapsto \mr{Map}_{\cat{sSet}}(K,\mr{Map}_\sfS(X,Y)) \text{ and }Y \mapsto \mr{Map}_{\cat{sSet}}(K,\mr{Map}_\sfS(Y,X))\end{equation}
are representable. This implies that colimits and limits indexed by small simplicial categories exist \cite[Theorem 3.73]{KellyEnriched}.
		
\begin{axiom}\label{axiom:enriched-complete-cocomplete} The category $\sfS$ is complete and cocomplete in the enriched sense.\end{axiom}

\newglossaryentry{copowering}{%
	name={\ensuremath{-\times-}},
	description={Copowering},
	type=symbols
}

\newglossaryentry{powering}{%
	name={\ensuremath{(-)^-}},
	description={Powering},
	type=symbols
}
		
The representing objects for the functors (\ref{eqn:indexed-colimits-limits}) respectively give us a copowering\index{copowering} $\gls{copowering} \colon \cat{sSet} \times \sfS \to \sfS$ and a powering\index{powering} $\gls{powering} \colon \cat{sSet} \times \sfS \to \sfS$, which satisfy $K \times (L \times -) \cong (K \times L) \times -$ and $((-)^K)^L \cong (-)^{K \times L}$. It also follows that $K \times -$ is left adjoint to $(-)^K$ in the enriched sense; there is an isomorphism of simplicial sets
		\[\mr{Map}_\sfS(K \times X,Y) \cong \mr{Map}_\sfS(X,Y^K)\]
		 natural in $X$ and $Y$.

\newglossaryentry{tensor}{%
	name={\ensuremath{- \otimes -}},
	description={Monoidal tensor product},
	type=symbols
}

\newglossaryentry{bunit}{%
	name={\ensuremath{\bunit}},
	description={Monoidal unit},
	type=symbols
}

\subsubsection{Monoidal structure}
A \emph{simplicially enriched monoidal structure}\index{monoidal structure} on a category $\sfS$ consists of 
\begin{itemize}\item a \emph{tensor product} simplicial functor $\gls{tensor} \colon \sfS \times \sfS \to \sfS$,
	\item a \emph{unit} object $\gls{bunit} \in \sfS$,
\end{itemize}
together with \emph{associativity natural isomorphisms} and \emph{right and left unit natural isomorphisms} subject to associativity pentagon and unit triangle axioms. By passing to $0$-simplices one obtains a monoidal structure on the underlying ordinary category. 

\newglossaryentry{braiding}{%
	name={\ensuremath{\beta_{X,Y}}},
	description={Braiding},
	type=symbols
}
A \emph{braided}\index{monoidal structure!braided} monoidal structure has additional \emph{braiding natural isomorphisms}\index{braiding}, $\gls{braiding} \colon X \otimes Y \to Y \otimes X$, subject to additional associativity hexagon and unit triangle axioms. It is a \emph{symmetric}\index{monoidal structure!symmetric} monoidal structure if $\beta_{Y,X} \circ \beta_{X,Y} = \mr{id}_{X \otimes Y}$.

\begin{notation}\label{not:k-monoidal} Let $k \in \{1,2,3,\ldots,\infty\}$. A \emph{$k$-monoidal structure} is a monoidal structure if $k=1$, a braided monoidal structure if $k=2$, and a symmetric monoidal structure if $k>2$.\index{$k$-monoidal}
\end{notation}

\begin{remark}This should not be confused with more refined notion of a $k$-fold monoidal category, as in e.g.~\cite{BFSV}.\end{remark}

\newglossaryentry{internalhom}{%
	name={\ensuremath{\Hom_{\sfS}(-,-)}},
	description={Internal hom},
	type=symbols
}
A simplicially enriched monoidal structure is said to be \emph{closed}\index{monoidal structure!closed} if $- \otimes Y \colon \sfS \to \sfS$ has an enriched right adjoint $\Hom_{\sfS}(Y,-)$ for all $Y$; that is, there are isomorphisms of simplicial sets
\[\mr{Map}_\sfS(X \otimes Y,Z) \cong \mr{Map}_\sfS(X,\Hom_\sfS(Y,Z)).\]
natural in $X$ and $Z$, or in other words $\gls{internalhom}$ is the internal hom. The unit isomorphisms imply that $\Hom_{\sfS}(\bunit,-)$ is naturally isomorphic to the identity functor $\mr{id} \colon \sfS \to \sfS$. 

Using the braiding we see that when $k \geq 2$, the functor $X \otimes - \colon \sfS \to \sfS$ also has a right adjoint, which is naturally isomorphic to $\Hom_{\sfS}(X,-)$. This right adjoint is canonical if $k > 2$, but it has a $\bZ$-torsor's worth of isomorphisms to $\Hom_\sfS(X,-)$ if $k=2$. If $k=1$ it does not necessarily follow that $X \otimes -$ has an enriched right adjoint, and when we want to additionally impose this
we say the monoidal structure is \emph{closed on both sides}.

The internal hom makes $\sfS$ enriched over itself: the internal identity $\bunit \to \Hom_{\sfS}(X,X)$ is adjoint to the identity $X \to X$ in $\sfS$, and the enriched composition $\Hom_{\sfS}(X,Y) \otimes \Hom_{\sfS}(Y,Z) \to \Hom_{\sfS}(X,Z)$ is given by using twice the evaluation maps $X \otimes \Hom_{\sfS}(X,Y) \to Y$ adjoint to the identity $\Hom_{\sfS}(X,Y) \to \Hom_{\sfS}(X,Y)$.

\begin{axiom}\label{axiom.cat-monoidal} $\sfS$ is equipped with a simplicially enriched closed $k$-monoidal structure. If $k=1$, it should be closed on both sides.\end{axiom}

\begin{notation}\label{not:otimes-unit-s} If we want to indicate that $\otimes$ and $\bunit$ are part of a simplicially enriched monoidal structure \emph{on the category $\sfS$}, we will denote them $\otimes_{\sfS}$ and $\bunit_{\sfS}$.\end{notation}
		
These axioms imply that the simplicial copowering behaves ``centrally'' with respect to the monoidal structure, even if the latter is not braided monoidal. In particular, there are isomorphisms
\[(K \times X) \otimes Y \cong K \times (X \otimes Y) \cong X \otimes (K \times Y)\]
naturally in $K$, $X$ and $Y$. By specializing $X$ to $\bunit$, we see that there is a (tautologically simplicial) functor 
\newglossaryentry{simplicial}{%
	name={\ensuremath{s(-)}},
	description={Copowering of monoidal unit},
	type=symbols
}
\begin{align*}\gls{simplicial} \colon \cat{sSet} &\lra \sfS \\
	K &\longmapsto K \times \bunit,\end{align*}
so that $s(K) \otimes X$ is naturally isomorphic to $K \times X$. This functor $s$ preserves colimits because it has a right adjoint, and is strong monoidal because 
\[s(K \times L) \cong (K \times L) \times \bunit \cong K \times  (L \times \bunit) \cong (K \times \bunit) \otimes (L \times \bunit) \cong s(K) \otimes s(L).\]

\subsubsection{Monoidal functors} If $\sfS$ and $\sfS'$ are simplicially enriched categories equipped with a simplicially enriched monoidal structures, and $F \colon \sfS \to \sfS'$ is a simplicial functor, then a \emph{lax monoidality}\index{monoidality} (or just \emph{monoidality}, for brevity) on $F$ is the
data of an enriched natural transformation 
\[m \colon \otimes_{\sfS'} \circ (F \times
F) \Longrightarrow F \circ \otimes_\sfS \colon \sfS \times \sfS \lra \sfS',\]
and a morphism $e \colon \bunit_{\sfS'} \to F(\bunit_{\sfS})$ in $\sfS'$.  The natural transformation $m$ is subject to an associativity axiom, and the morphism $e$ is subject to a unitality axiom. This data is a \emph{strong monoidality} if $e$ is an isomorphism and $m$ is a natural isomorphism. An \emph{oplax monoidality} on $F$ is a lax monoidality on $F^\op \colon \sfS^\op \to (\sfS')^\op$.  

If $\sfS$ and $\sfS'$ are additionally equipped with braidings or symmetries, then (lax, strong, or oplax) monoidalities on a functor are \emph{braided} or \emph{symmetric} if they satisfy the additional property of being compatible with the braidings in $\sfS$ and $\sfS'$. Following Notation \ref{not:k-monoidal}, for $k \in \{1,2,3,\ldots,\infty\}$, a \emph{$k$-monoidal functor} means a monoidal functor if $k=1$, a braided monoidal functor if $k=2$, and a symmetric monoidal functor if $k>2$.

We emphasize that monoidal structures and monoidalities on functors are really extra structures, but we shall follow common usage of the ``monoidal'' adjective: a \emph{simplicially enriched $k$-monoidal category} is a category together with a simplicially enriched $k$-monoidal structure, and a lax/strong/oplax $k$-monoidal functor is a functor between the underlying simplically enriched categories together with a lax/strong/oplax $k$-monoidality. A $k$-monoidal natural transformation between $k$-monoidal functors is a natural transformation between the underlying functors, subject to two conditions (one about multiplication and one about units) but involves no additional data.

\subsubsection{Pointed categories and pointed simplicial sets}\label{sec:point-categ-point} A simplicially enriched category $\sfS$ which is also pointed, i.e.\ has an object $\ast$ that is simultaneously initial and terminal, is automatically enriched in pointed simplicial sets. The adjunctions in the unpointed setting in Section \ref{sec:simplicial-enrichment} imply those in the pointed setting. As a result, we obtain a refinement of copowering and powering to $\cat{sSet}_\ast$, which will be denoted $\wedge \colon \cat{sSet}_* \times \sfS \to \sfS$ and $(-)_\ast^{-} \colon \cat{sSet}_* \times \sfS \to \sfS$. We may recover the original simplicial enrichment by applying the forgetful functor $U^+ \colon \cat{sSet}_+ \to \cat{sSet}$ to the pointed simplicial enrichment, and we may recover the copowering and powering by precomposing with the left adjoint $(-)_+ \coloneqq F^+ \colon \cat{sSet} \to \cat{sSet}_\ast$ to $U^+$.
\newglossaryentry{halfwedge}{%
	name={\ensuremath{\rtimes}},
	description={Half wedge product},
	type=symbols
}
\newglossaryentry{wedge}{%
	name={\ensuremath{\wedge}},
	description={Wedge product},
	type=symbols
}
\newglossaryentry{vee}{%
	name={\ensuremath{\vee}},
	description={Coproduct in pointed category},
	type=symbols
}
\newglossaryentry{sqcup}{%
	name={\ensuremath{\sqcup}},
	description={Coproduct},
	type=symbols
}
If $\sfS$ is cocomplete, then associated to this we have $\gls{halfwedge} \colon \cat{sSet}_* \times \sfS \to \sfS_*$ given by
\[A \rtimes X \coloneqq \mathrm{colim}(* \leftarrow * \times X \rightarrow A \times X),\]
and $\gls{wedge} \colon \mathsf{sSet}_* \times \sfS_* \to \sfS_*$ given by
\[A \wedge X \coloneqq \mathrm{colim}(* \leftarrow * \times X \vee A \times * \rightarrow A \times X).\]
These are related by the familiar isomorphism $A \wedge X_+ \cong A \rtimes X$. Here, as later, we write $\gls{vee}$ rather than $\gls{sqcup}$ for the coproduct in a pointed category.

\subsubsection{Sifted colimits and geometric realization} Since the $k$-monoidal structure is closed (on both sides if $k=1$), $\otimes$ preserves colimits in each variable. While this does not imply that $\otimes \colon \sfS \times \sfS \to \sfS$ preserves all colimits, it still preserves sifted colimits. Recall that a \emph{sifted colimit}\index{sifted!colimit} is a colimit over a sifted category, and a diagram $\cat{I}$ is \emph{sifted} if the diagonal functor $\cat{I} \to \cat{I}\times \cat{I} $ is final. For example, the reflexive coequalizer diagram 
\[\begin{tikzcd} 
{[1]} \arrow[shift right=1ex]{r} \arrow[shift left=1ex]{r} & {[0]} \ar{l}
\end{tikzcd}\]
is sifted, as is $\Delta^\mr{op}$ and every filtered category. In fact, a functor between cocomplete categories preserves sifted colimits if and only if it preserves filtered colimits and reflexive coequalizers \cite{AdamekRosickyVitale}. The following is elementary, see e.g.\ \cite[Proposition 5.7]{HarperOperads}. 

\begin{lemma}\label{lem.cat-sifted} If $\sfS$ satisfies Axioms \ref{axiom:simplicial-enriched}, \ref{axiom:enriched-complete-cocomplete}, and \ref{axiom.cat-monoidal}, then the functors $\otimes \colon \sfS \times \sfS \to \sfS$ and $\times \colon \cat{sSet} \times \sfS \to \sfS$ commute with sifted colimits.\end{lemma}

We will most often use this to say that sifted colimits commute with $\otimes$-powers: $\colim_{i \in \cat{I}} (X_i^{\otimes n}) \cong (\colim_{i \in \cat{I}} X_i)^{\otimes n}$ as long as the diagram $\cat{I}$ is sifted.

\newglossaryentry{simpob}{%
	name={\ensuremath{\cat{sS}}},
	description={Simplicial objects in $\sfS$},
	type=symbols
}
\newglossaryentry{singsfs}{%
	name={\ensuremath{\mr{Sing}}},
	description={Singular simplicial object},
	type=symbols
}
\newglossaryentry{geomrel}{%
	name={\ensuremath{|-|}},
	description={Geometric realization},
	type=symbols
}
\newglossaryentry{thickgeomrel}{%
	name={\ensuremath{\fgr{-}}},
	description={Thick geometric realization},
	type=symbols
}
Axiom \ref{axiom:enriched-complete-cocomplete} has further consequences for geometric realizations of simplicial objects, which will play a large role in the later parts. We let $\gls{simpob}$ denote the category of \emph{simplicial objects} in $\sfS$, i.e.\ of functors $\Delta^\mr{op} \to \sfS$ and natural transformations between these. We may then define an internal singular simplicial object functor
\begin{align*}\mr{Sing} \colon \sfS &\lra \cat{sS} \\
	X &\longmapsto \mr{Sing}_\sfS(X) \coloneqq \big([n] \mapsto X^{\Delta^n}\big).\end{align*}
The functor $\gls{singsfs}$ has a left adjoint $\gls{geomrel} \colon \cat{sS} \to \sfS$, an internal version of geometric realization,\index{geometric realization} given on $X_\bullet \in \cat{sS}$ by the coend $\smallint^{n \in \Delta^\mr{op}} \Delta^n \times X_n$, which is isomorphic to the reflexive coequalizer of the two maps
\begin{equation}\label{eqn:geomrelcoeq} \begin{tikzcd}\bigsqcup_{[n] \to [m]} \Delta^n \times X_m \arrow[shift left=.5ex]{r} \arrow[shift left=-.5ex]{r} & \bigsqcup_{n} \Delta^n \times X_n.\end{tikzcd}\end{equation}

There is also a natural transformation $|X_\bullet \otimes Y_\bullet| \to |X_\bullet| \otimes |Y_\bullet|$, where on the left $\otimes$ denotes the levelwise tensor product. To construct it, we provide a natural map $\bigsqcup_{n \geq 0} \Delta^n \times (X_n \otimes Y_n) \to \gr{X_\bullet} \otimes \gr{Y_\bullet}$ coequalizing the two maps of (\ref{eqn:geomrelcoeq}), by taking the canonical map
\begin{equation}\label{eqn:geom-rel-nat} \begin{tikzcd}
\left(\bigsqcup_{n,m \geq 0} \Delta^n \times \Delta^m \times (X_n \otimes Y_m) \right) \cong \left(\bigsqcup_{n \geq 0} \Delta^n \times X_n\right) \otimes \left(\bigsqcup_{m \geq 0} \Delta^m \times Y_m\right) \arrow{d} \\
{\gr{X_\bullet} \otimes \gr{Y_\bullet}}
\end{tikzcd}\end{equation}
and mapping into the terms for $m=n$ using the diagonal map $\Delta^n \to \Delta^n \times \Delta^n$. The following is a consequence of \cite[Proposition A.3]{BergerMoerdijkBV}, and the fact that we have isomorphisms $\gr{\Delta(-,[n]) \times \Delta(-,[m])} \to \gr{\Delta(-,[n])} \times \gr{\Delta(-,[m])}$ satisfying natural associativity, unit and symmetry conditions.

\begin{lemma}\label{lem.monoidal-geomrel} If $\sfS$ satisfies Axioms \ref{axiom:simplicial-enriched}, \ref{axiom:enriched-complete-cocomplete}, and \ref{axiom.cat-monoidal}, then $\otimes$ commutes with geometric realization, i.e.\ \eqref{eqn:geom-rel-nat} is a natural isomorphism $|X_\bullet \otimes Y_\bullet| \overset{\cong}\to |X_\bullet| \otimes |Y_\bullet|$.\end{lemma}

Letting $\Delta_\mr{inj} \subset \Delta$ denote the subcategory with only injective morphisms, we define the category $\cat{ssS}$ of \emph{semi-simplicial objects} in $\sfS$ to be that of functors $\Delta^\mr{op}_{\mr{inj}} \to \cat{S}$ and natural transformations between these. Analogous to geometric realization of simplicial objects we have \emph{thick geometric realization}\index{geometric realization!thick} $\gls{thickgeomrel}$ of semi-simplicial objects, given by replacing $\Delta$ with $\Delta_\mr{inj}$ in the coend and coequalizer descriptions. We may left Kan extend along the inclusion $\sigma \colon  \Delta^\mr{op}_\mr{inj} \hookrightarrow \Delta^\mr{op}$ to construct a simplicial object out of a semi-simplicial object. This is given by freely adding degeneracies\index{adding degeneracies}; for $Z_\bullet \in \cat{ssS}$ we have that 
\begin{equation}\label{eqn:add-degeneracies} \sigma_* Z_n \cong \bigsqcup_{[n] \twoheadrightarrow [m]} Z_m,\end{equation}
from which it follows that there is a natural isomorphism $|\sigma_*(-)| \cong ||-||$. It is usually not true that thick geometric realization commutes with $\otimes$.

\subsection{Examples} The following is a non-exhaustive list of examples of monoidal categories which satisfy the axioms of Section \ref{sec:axioms-of-cats}:
\newglossaryentry{sset}{%
	name={\ensuremath{\cat{sSet}}},
	description={Category of simplicial sets},
	type=symbols
}
\newglossaryentry{ssetast}{%
	name={\ensuremath{\cat{sSet}_*}},
	description={Category of pointed simplicial sets},
	type=symbols
}
\newglossaryentry{top}{%
	name={\ensuremath{\cat{Top}}},
	description={Category of CGWH spaces},
	type=symbols
}
\newglossaryentry{topast}{%
	name={\ensuremath{\cat{Top}_*}},
	description={Category of pointed CGWH spaces},
	type=symbols
}
\newglossaryentry{smod}{%
	name={\ensuremath{\cat{sMod}_\bk}},
	description={Category of simplicial $\bk$-spaces},
	type=symbols
}
\newglossaryentry{spsigma}{%
	name={\ensuremath{\cat{Sp}^\Sigma}},
	description={Category of symmetric spectra},
	type=symbols
}
\newglossaryentry{rmod}{%
	name={\ensuremath{R\text{-}\cat{Mod}}},
	description={Category of $R$-module spectra},
	type=symbols
}
\begin{itemize}
	\item $\gls{sset}$, simplicial sets, with cartesian product.
	\item $\gls{ssetast}$, pointed simplicial sets, with smash product.
	\item $\gls{top}$, compactly generated weakly Hausdorff (CGWH) topological spaces, with cartesian product. 
	\item $\gls{topast}$, pointed CGWH topological spaces, with smash product.
	\item $\gls{smod}$, simplicial $\bk$-modules, with levelwise tensor product over $\bk$.
	\item $\gls{spsigma}$, symmetric spectra in the sense of \cite{HoveyShipleySmith}, with smash product.
	\item $\gls{rmod}$, $R$-module symmetric spectra over a commutative ring 
	spectrum $R$, with smash product over $R$.
\end{itemize}

\begin{remark}
We purposefully did not include (non-negatively graded) chain complexes over a commutative ring $\bk$, as there is no strong monoidal functor $s \colon \cat{sSet} \to \cat{Ch}_\bk$; neither the Eilenberg--Zilber nor the Alexander--Whitney map is an isomorphism. This is only mildly inconvenient, as by the Dold--Kan theorem \cite[Chapter III.2]{GoerssJardine} simplicial $\bk$-modules is an equivalent category which does satisfy our axioms.
\end{remark}

With the exception of the case $R\text{-}\cat{Mod}$, which will be covered by Section \ref{sec:module-cats}, we shall now verify the axioms for the aforementioned categories.

\subsubsection{Simplicial sets and pointed simplicial sets} The category $\cat{sSet}$ of \emph{simplicial sets} is the category $\cat{Fun}(\Delta^\mr{op},\cat{Set})$ of functors $\Delta^\mr{op} \to \cat{Set}$ and natural transformations between these. It is simplicially enriched by the simplicial mapping set with $p$-simplices $\mr{Map}_\cat{sSet}(X,Y)_p \coloneqq \mr{Hom}_\cat{sSet}(X \times \Delta^p,Y)$, establishing Axiom \ref{axiom:simplicial-enriched}. 

Limits and colimits are computed pointwise and hence $\cat{sSet}$ is complete and cocomplete in the ordinary sense. To establish Axiom \ref{axiom:enriched-complete-cocomplete}, it remains to note $Y \mapsto \mr{Map}_\cat{sSet}(K,\mr{Map}_\cat{sSet}(X,Y))$ and $Y \mapsto \mr{Map}_\cat{sSet}(K,\mr{Map}_\cat{sSet}(Y,X))$ are representable by $K \times X$ and $\mr{Map}_\cat{sSet}(K,X)$ respectively. 

The functor $\mr{Map}_\cat{sSet}(X,-)$ is the right adjoint to $- \times X$, with cartesian product $\times$ the monoidal product of a cartesian enriched closed symmetric monoidal structure, establishing Axiom \ref{axiom.cat-monoidal}. Pointed simplicial sets are similar. 

\subsubsection{CGWH topological spaces and pointed CGWH topological spaces} The category $\cat{Top}$ of \emph{compactly generated weakly Hausdorff topological spaces} is a standard convenient category of topological spaces, see \cite[Chapter 5]{MayConcise}. It is enriched in simplicial sets by taking the $p$-simplices of $\mr{Map}_\cat{Top}(X,Y)$ to be $\mr{Hom}_\cat{Top}(X \times \Delta^p,Y)$, establishing Axiom \ref{axiom:simplicial-enriched}. 

It is complete and cocomplete in the ordinary sense, and the functors $Y \mapsto \mr{Map}_\cat{sSet}(K,\mr{Map}_\cat{Top}(X,Y))$ and $Y \mapsto \mr{Map}_\cat{sSet}(K,\mr{Map}_\cat{Top}(Y,X))$ are representable by $|K| \times X$ and $\mr{Map}_\cat{Top}(|K|,X)$ respectively, where $\mr{Map}_\cat{Top}(-,-)$ is the mapping space in the (retopologized) compact-open topology. This finishes the verification of Axiom \ref{axiom:enriched-complete-cocomplete}. 

Finally, it is a cartesian enriched closed symmetric monoidal category, with $\otimes$ given by (retopologized) cartesian product, establishing Axiom \ref{axiom.cat-monoidal}. In particular, $s \colon \cat{sSet} \to \cat{Top}$ is given by geometric realization $|-| \colon \cat{sSet} \to \cat{Top}$, which has as right adjoint the singular simplicial set $\mr{Sing} \colon \cat{Top} \to \cat{sSet}$.

\subsubsection{Simplicial $\bk$-modules} Let $\bk$ be a commutative ring. The category $\cat{sMod}_\bk$ of \emph{simplicial $\bk$-modules} is the category of functors $\cat{Fun}({\Delta^\mr{op}}, \cat{Mod}_\bk)$. The simplicial enrichment is given by taking the $p$-simplices of $\mr{Map}_{\cat{sMod}_\bk}(X,Y)$ to the maps $X \otimes \bk[\Delta^p] \to Y$ of simplicial $\bk$-modules, establishing Axiom \ref{axiom:simplicial-enriched}.

As limits and colimits are computed pointwise, this is complete and cocomplete since $\cat{Mod}_\bk$ is. The functors $Y \mapsto \mr{Map}_\cat{sSet}(K,\mr{Map}_\cat{sMod_\bk}(X,Y))$ and $Y \mapsto \mr{Map}_\cat{sSet}(K,\mr{Map}_\cat{sMod_\bk}(Y,X))$ are representable, respectively by the levelwise tensor product $\bk[K] \otimes X$ and $\mr{Map}_\cat{sSet}(K,X)$, which inherits from $X$ the structure of a simplicial $\bk$-module. This finishes the verification of Axiom \ref{axiom:enriched-complete-cocomplete}. 

It has an enriched closed symmetric monoidal structure by using tensor product $\otimes_\bk$ of $\bk$-modules levelwise, establishing Axiom \ref{axiom:simplicial-enriched}. In particular, the map $s \colon \cat{sSet} \to \cat{sMod}_\bk$ is given by taking the free $\bk$-module levelwise, which is strong monoidal and has a right adjoint giving by forgetting the $\bk$-module structure.

\subsubsection{Symmetric spectra} \label{sec:symmetric-spectra} The category $\cat{Sp}^\Sigma$ of \emph{symmetric spectra} was introduced by Hovey, Shipley and Smith \cite{HoveyShipleySmith}, see also \cite{HoveySymmetricModel,Schwedebook}. A symmetric spectrum $E$ is a sequence $\{E_n\}_{n \geq 0}$ of pointed simplicial sets with $\fS_n$-actions and maps $E_n \wedge S^1 \to E_{n+1}$ so that the iterated suspension maps $E^n \wedge S^k \to E_{n+k}$ are $\fS_n \times \fS_k$-equivariant. A morphism of symmetric spectra $E \to E'$ is a sequence of maps $E_n \to E_n'$ compatible with the structure.

To define the simplicial enrichment, we start with the tensoring of $\cat{Sp}^\Sigma$ over $\cat{sSet}$ by $(K \times E)_n = E_n \wedge K_+$. Using this we may define $\mr{Map}_\cat{Sp^\Sigma}(E,F)$ by setting its $p$-simplices to be $\mr{Hom}_{\cat{Sp^\Sigma}}(\Delta^p \times E,F)$. This gives Axiom \ref{axiom:simplicial-enriched}.

As limits and colimits may be computed objectwise, $\cat{Sp}^\Sigma$ is complete and cocomplete. The functor $F \mapsto \mr{Map}_\cat{sSet}(K,\mr{Map}_\cat{Sp^\Sigma}(E,F))$ is represented by $K \times E$ with $n$th level given by $K_+ \wedge E_n$. Similarly the functor $F \mapsto \mr{Map}_\cat{sSet}(K,\mr{Map}_\cat{Sp^\Sigma}(F,E))$ is represented by $E^K$ with $n$th level given by $\mr{Map}_\ast(K_+,E_n)$. This verifies Axiom \ref{axiom:enriched-complete-cocomplete}.

The motivation for introducing symmetric group actions is to endow $\cat{Sp}^\Sigma$ with an enriched symmetric monoidal structure. This is called the smash product, denoted $\wedge$, and is given by setting $(E \wedge F)_n$ to be the coequalizer of the diagram 
\[\begin{tikzcd}
\bigvee_{p+1+q = n} (\fS_n)_+ \wedge_{\fS_p \times \fS_1 \times \fS_q} E_p \wedge S^1 \wedge F_q \arrow[shift left=.5ex]{r} \arrow[shift right=.5ex]{r}  &[-5pt] \bigvee_{p+q = n} (\fS_n)_+ \wedge_{\fS_p \times \fS_q} E_q \wedge F_q.\end{tikzcd}\]
This formula implies that the smash product commutes with colimits in each variable, and jointly commutes with sifted colimits, because the smash product of pointed simplicial sets does. This is closed with right adjoint to $- \wedge E$ given by the function spectrum $\mr{Fun}(E,-)$ of Definition 2.2.9 of \cite{HoveyShipleySmith}. The monoidal unit $\bunit$ is the sphere spectrum $\bS$ with $n$th entry given by $S^n \coloneqq (S^1)^{\wedge n}$. This completes the verification of Axiom \ref{axiom.cat-monoidal}.

The sphere spectrum is an example of a suspension spectrum: for any pointed simplicial set $K$ there is a symmetric spectrum $\Sigma^\infty K$ with $n$th entry given by $S^n \wedge K$. The strong monoidal functor $s \colon \cat{sSet} \to \cat{Sp}^\Sigma$ is the composition of $+ \colon \cat{sSet} \to \cat{sSet}_\ast$ and $\Sigma^\infty \colon \cat{sSet}_\ast \to \cat{Sp}^\Sigma$, the latter being strong monoidal by Proposition 1.3.1 of \cite{HoveyShipleySmith}. It has right adjoint given by $\mr{ev}_0$, see Proposition 2.2.6 ($F_0$ is their notation for $\Sigma^\infty$).

\subsection{Module categories}\label{sec:module-cats} 
\newglossaryentry{grmod}{%
	name={\ensuremath{\gR\text{-}\cat{Mod}}},
	description={Category of $\gR$-module objects},
	type=symbols
}
If $\sfS$ is symmetric monoidal, we can define commutative algebra objects in $\sfS$. Given a commutative algebra $\gR$ in $\sfS$, with underlying object $R \in \sfS$, we can define the category $\gls{grmod}$ of (left) $\gR$-modules. The following shows that this satisfies the axioms of Section \ref{sec:axioms-of-cats}. 

\begin{proposition}\label{prop:ModuleCatAxioms}
If $\sfS$ satisfies the axioms of Section \ref{sec:axioms-of-cats} and $\gR$ is a commutative algebra object of $\sfS$, then $\gR\text{-}\cat{Mod}$ also satisfies the axioms of Section \ref{sec:axioms-of-cats}.\end{proposition}

This satisfies Axiom \ref{axiom:simplicial-enriched}, as there is a simplicial enrichment given by taking $\mr{Map}_{\gR\text{-}\cat{Mod}}(\gM,\gN)$ to be the coreflexive equalizer in $\cat{sSet}$ of the two maps
\[\begin{tikzcd}\mr{Map}_{\gR\text{-}\cat{Mod}}(\gM,\gN) \rar & \mr{Map}_\sfS(M,N) \arrow[shift left=1ex]{r} \arrow[shift left=-1ex]{r} & \mr{Map}_\sfS(\gR \otimes M,N), \lar\end{tikzcd}\]
where $M$, $N$ denote the underlying objects of $\gM$, $\gN$ in $\sfS$.

\newglossaryentry{ugr}{%
	name={\ensuremath{U^\gR}},
	description={Forgetful functor on $\gR$-modules},
	type=symbols
}
\newglossaryentry{fgr}{%
	name={\ensuremath{F^\gR}},
	description={Free $\gR$-module functor},
	type=symbols
}
\newglossaryentry{otimesgr}{%
	name={\ensuremath{-\otimes_\gR -}},
	description={Tensor product over $\gR$-modules},
	type=symbols
}
\newglossaryentry{homgr}{%
	name={\ensuremath{\Hom_\gR(-,-)}},
	description={Internal hom of $\gR$-modules},
	type=symbols
}
There is a forgetful functor $\gls{ugr} \colon \gR\text{-}\cat{Mod} \to \sfS$ with left adjoint given by a free $\gR$-module functor $\gls{fgr} \colon \sfS \to \gR\text{-}\cat{Mod}$, explicitly given by $X \mapsto \gR \otimes X$. It is complete and cocomplete, as the forgetful functor to $\sfS$ creates both limits and colimits. The functor $\gN \mapsto \mr{Map}_\cat{sSet}(K,\mr{Map}_\sfS(\gM,\gN))$ is representable by $K \times M$, which inherits a $\gR$-module structure the $\gR$-module structure $a_\gM \colon \gR \otimes M \to M$, and the functor $\gN \mapsto \mr{Map}_\cat{sSet}(K,\mr{Map}_\sfS(\gN,\gM))$ is similarly representable by $M^K$. This verifies Axiom \ref{axiom:enriched-complete-cocomplete}. 

To verify Axiom \ref{axiom.cat-monoidal}, if $\gM$ and $\gN$ are $\gR$-modules, with underlying objects $M, N \in \sfS$, we define the tensor product $\gM \otimes_\gR \gN$ to have underlying object given by the reflexive coequalizer
\[\begin{tikzcd}M \otimes R \otimes N \arrow[shift left=1 ex]{r} \arrow[shift right=1ex]{r} & M \otimes N \lar \rar & U^\gR(\gM \otimes_\gR \gN).\end{tikzcd}\]
of the maps given by the $\gR$-module structure maps of $\gM$ and $\gN$, with reflection given by the unit of $\gR$. This may be endowed with the structure of an $\gR$-module, using the left $\gR$-module structure on $\gM$ (this uses the fact that $\gR$ is a commutative algebra object). This yields a functor $\gls{otimesgr} \colon \gR\text{-}\cat{Mod} \times \gR\text{-}\cat{Mod} \to \gR\text{-}\cat{Mod}$.

The resulting symmetric monoidal structure is closed. To describe the right adjoint, let us first explain how to endow the internal hom $\Hom_\sfS(X,N)$ with an $\gR$-module structure when $N$ is endowed with an $\gR$-module structure $a_\gN \colon R \otimes N \to N$. Namely, the structure map is adjoint to the map
\[R \otimes X \otimes \Hom_\sfS(X,N) \xrightarrow{\,R \otimes \text{ev}\,} R \otimes N \overset{a_\gN}\lra N.\]
We write $\Hom_\sfS(X,\gN)$ for the resulting $\gR$-module. Now we can describe the right adjoint $Hom_\gR(\gM,-)$ to $\gM \otimes_\gR -$ as the equalizer
\[\begin{tikzcd}\Hom_\gR(\gM,\gN) \rar & \Hom_\sfS(M,\gN) \arrow[shift left=.5ex]{r} \arrow[shift right=.5ex]{r} & \Hom_\sfS(R \otimes M,\gN)\end{tikzcd}\]
where the top map is induced by $R \otimes M \to M$ and the bottom map is the adjoint of $R \otimes M \otimes \cH om_\sfS(M,N) \to R \otimes N \to N$. A priori this is an object of $\sfS$, but the equalizer is one of $\gR$-modules and the forgetful functor to $\sfS$ creates limits so this describes $\Hom_R(M,N)$ as an $\gR$-module. This yields a functor $\gls{homgr} \colon \gR\text{-}\cat{Mod}^\op \times \gR\text{-}\cat{Mod} \to \gR\text{-}\cat{Mod}$.

In this case, the strong monoidal functor $s \colon \cat{sSet} \to \gR\text{-}\cat{Mod}$ is given by the composition of the strong monoidal functor $\cat{sSet} \to \sfS$ with $F^\gR$.

\subsection{Diagram categories}
\label{sec:functor-cats}
\newglossaryentry{sfc}{%
	name={\ensuremath{\sfC}},
	description={Diagram category, of the form $\sfS^\sfG$},
	type=symbols
}
\newglossaryentry{sfg}{%
	name={\ensuremath{\sfG}},
	description={Indexing category},
	type=symbols
}
If $\gls{sfg}$ is a small category then we may form the category $\gls{sfc} = \cat{Fun}(\sfG, \sfS) = \sfS^\sfG$ of $\sfG$-shaped diagrams in $\sfS$.\index{diagram category} When $\sfG$ is discrete, i.e.\ the only morphisms are identities, taking $\sfS^\sfG$ amounts to adding an additional grading to $\sfS$.\index{discrete category} We will denote the evaluation of an object $X$ of $\sfS^\sfG$ at an object $g \in \sfG$ by $X(g)$. The following explains when such a diagram category satisfies the axioms of Section \ref{sec:axioms-of-cats}.

\begin{proposition}\label{prop.functor-cats-axioms-of-cats}
If $\sfS$ satisfies the axioms of Section \ref{sec:axioms-of-cats}, so is in particular $k$-monoidal, and $\sfG$ is a $k$-monoidal category, then $\sfC = \sfS^\sfG$ also satisfies the axioms of Section \ref{sec:axioms-of-cats}.\end{proposition}

To prove Proposition \ref{prop.functor-cats-axioms-of-cats}, we begin by noting that it has a simplicial enrichment with $\mr{Map}_{\sfC}(X,Y)$ to be the sub-simplicial set of $\prod_{g \in \sfG} \mr{Map}_\sfS(X(g),Y(g))$ of natural transformations. This verifies Axiom \ref{axiom:simplicial-enriched}.

To verify Axiom \ref{axiom:enriched-complete-cocomplete}, we first observe that because colimits and limits are computed pointwise, the category $\sfC$ is complete and cocomplete. By taking the copowering and powering objectwise, we obtain a copowering and powering of $\sfC$ over $\cat{sSet}$ giving rise to the required representing objects.

\newglossaryentry{otimess}{%
	name={\ensuremath{-\otimes_\sfC-}},
description={Day convolution tensor product},
type=symbols
}
Next we claim that if $\sfS$ and $\sfG$ are both $k$-monoidal then $\sfC$ has a $k$-monoidal structure, which is closed if the $k$-monoidal structure on $\sfC$ is, verifying Axiom \ref{axiom.cat-monoidal}. The monoidal structure on $\sfC$ is given by \emph{Day convolution},\index{Day convolution} defined as follows. We shall introduce the convention that the tensor product of $\sfG$ is denoted $\oplus_\sfG$, even though we do \emph{not} assume it is symmetric. This is to prevent the symbol $\otimes$ from being overloaded. The monoidal structure on $\sfS$ induces an exterior product
\begin{align*}\bar{\otimes} \colon \sfC \times \sfC = \sfS^\sfG \times \sfS^\sfG &\lra \sfS^{\sfG \times \sfG} 
\\ (X,Y) &\longmapsto \left[(g,h) \mapsto X(g) \otimes_\sfS Y(h)\right],\end{align*}  
and we may then define $X \otimes_\sfC Y$ as the left Kan extension of $X \overline\otimes Y \colon \sfG \times \sfG \to \sfS$ along the monoidal structure $\oplus_\sfG \colon \sfG \times \sfG \to \sfG$, i.e.\
\begin{equation*}
\begin{tikzcd}
\sfG \times \sfG \rar{X \times Y} \dar[swap]{\oplus_\sfG} &
\sfS \times \sfS \dar{\otimes_\sfS} \arrow[Rightarrow]{ld} \\
\sfG \rar[swap]{X \otimes_\sfC Y} & \sfS.
\end{tikzcd}
\end{equation*} 
The result is a functor $X \otimes_\sfC Y \colon \sfG \to \sfS$, i.e. an object of $\sfC$. This construction gives a simplicially enriched functor $\gls{otimess} \colon \sfC \times \sfC \to \sfC$.

The following is often stated only when $\sfS$ in symmetric monoidal. The reason is that in the literature $\sfG$ is often taken to be enriched in $\sfS$ (e.g. \cite{Day}), and then checking associativity for coends requires the existence of a symmetric braiding. This is not necessary if $\sfG$ is an ordinary small category, and in fact, the centrality of $\cat{sSet}$ should allow one to take $\sfG$ to be a simplicial category.

\begin{theorem}The functor $\otimes_{\sfC} \colon \sfC \times \sfC \to \sfC$ is part of a simplicially enriched $k$-monoidal structure on $\sfC = \sfS^\sfG$ with unit $\bunit_\sfC$ given by $\sfG(\bunit_\sfG,-) \otimes \bunit_\sfS$. If the monoidal structure on $\sfS$ is closed (on both sides if $k=1$) then this monoidal structure is closed as well (on both sides if $k=1$).\end{theorem}

By definition of $\otimes_{\sfC}$ as a left Kan extension, specifying a morphism $\phi \colon X \otimes_\sfC Y \to Z$ is the same as specifying morphisms $\phi_{g,h} \colon X(g) \otimes_\sfS Y(h) \to Z(g \oplus_\sfG h)$ in $\sfS$ for all $g,h \in \sfG$, forming a natural transformation of functors $\sfG \times \sfG \to \sfS$. In particular, there is a universal map $X(g) \otimes_\sfS Y(h) \to (X \otimes_\sfC Y)(g \oplus_\sfG h)$. Given a braiding $\beta^\sfG_{g,h} \colon g \oplus_\sfG h \to h \oplus_\sfG g$ of $\oplus_\sfG$ we construct a braiding $\beta^\sfC_{X,Y} \colon X \otimes_\sfC Y \to Y \otimes_\sfC X$ of $\otimes_\sfC$ in terms of the maps (the choices of braidings used in the formula are just a convention, though using $(\beta^\sfG)^{-1}$ means that Yoneda is braided monoidal if $\sfG^\mr{op}$ is given the braiding $(\beta^\sfG)^{-1}$) 
\begin{align*}
  (\beta^\sfC_{X,Y})_{g,h} \colon X(g) \otimes_\sfS Y(h) &\xrightarrow{\beta^\sfS} Y(h) \otimes_{\sfS} X(g)
  \to (Y \otimes_\sfC X)(h \oplus_\sfG g) \\
  &\xrightarrow{(Y
    \otimes_{\sfC} X)((\beta^\sfG_{g,h})^{-1})} (Y \otimes_\sfC X)(g
  \oplus_\sfG h).
\end{align*}
If the braidings on $\sfG$ and $\sfS$ are symmetries, then so is the induced braiding on $\sfC$.

\index{change-of-diagram-category} 
\newglossaryentry{pup}{%
	name={\ensuremath{p^*}},
	description={Restriction along $p$},
	type=symbols
}
\newglossaryentry{pdown}{%
	name={\ensuremath{p_*}},
	description={Left Kan extension along $p$},
	type=symbols
}
How does this depend on $\sfG$? For any functor $p \colon \sfG \to \sfG'$, there is a simplicially enriched change-of-diagram-category functor $p_* \colon \sfS^{\sfG} \to \sfS^{\sfG'}$ obtained as the left adjoint to the functor $p^* \colon \sfS^{\sfG'} \to \sfS^\sfG$, given by enriched left Kan extension. The following is proven by a straightforward manipulation of coends.

\begin{lemma}\label{lem:oplax-lax-Day-convo}
	If $p \colon \sfG \to \sfG'$ is (strong or oplax) monoidal, then
	$p_* \colon \sfS^\sfG \to \sfS^{\sfG'}$ is (strong or lax) monoidal.
\end{lemma}

In particular, we can apply this to the inclusion of the monoidal unit $\bunit_\sfG \colon \ast \to \sfG$. This is a strong monoidal functor, so the result is a strong monoidal simplicially enriched functor $(\bunit_\sfG)_* \colon \sfS \to \sfS^\sfG$ which is a left adjoint and so preserves colimits. We then obtain a strong 
monoidal functor $s_\sfG \colon \cat{sSet} \to \sfS^\sfG$ as the composition $(\bunit_\sfG)_* \circ s \colon \cat{sSet} \to \sfS \to \sfS^\sfG$.

\begin{example}\label{exam:functor-cats-filtered}
	\newglossaryentry{zleq}{%
		name={\ensuremath{\bZ_\leq}},
	description={Integers as poset with usual order},
	type=symbols
}
	Although our main examples shall have $\sfG$ be a groupoid, it is convenient to allow also non-invertible morphisms.  In particular, in the discussion of filtrations in Section \ref{sec:filtered-objects-1} we will use $\gls{zleq}$, the set of integers considered as a poset with the usual order. An object $X$ of $\sfS^{\sfG \times \Z_{\leq}} \cong \sfC^{\Z_\leq}$ is then a filtered object of $\sfC = \sfC^{\sfG}$, though the maps $X(i) \to X(i+1)$ need not be injective in any way.
\end{example}
\section{Sifted monads and indecomposables}
\label{sec:sifted-monads}
As explained in the introduction, our eventual goal is a robust theory of (cellular) $E_k$-algebras and of (cellular) modules over an associative algebra, both in a category $\sfC = \sfS^\sfG$. This can be done in a uniform manner by considering both examples as instances of algebras over an operad in $\sfC$, and for many purposes all that is important is that both are examples of algebras over a monad on $\sfC$ which preserves sifted colimits (we shall explain in Section \ref{sec:operads} that the monad associated to an operad always preserves sifted colimits). The goal of this section is to explain that theory. Unless mentioned otherwise, we assume that $\sfS$ satisfies the axioms of Section \ref{sec:axioms-of-cats} and hence so does $\sfC$ by Proposition \ref{prop.functor-cats-axioms-of-cats}.

\subsection{Monads and adjunctions} \label{sec:monads}

The category of endofunctors $\cat{Fun}(\sfC,\sfC)$ is a monoidal category under composition, the monoidal unit being the identity functor $\mr{id} \colon \sfC \to \sfC$.

\newglossaryentry{T}{%
	name={\ensuremath{T}},
description={Monad},
type=symbols
}
\newglossaryentry{muT}{%
	name={\ensuremath{\mu^T}},
description={Monad multiplication},
type=symbols
}
\newglossaryentry{1T}{%
	name={\ensuremath{1^T}},
description={Monad unit},
type=symbols
}
\begin{definition}A \emph{monad}\index{monad} $\gls{T}$ is a unital monoid in $\cat{Fun}(\sfC,\sfC)$.\end{definition}

Concretely, this means that there are natural transformation $\gls{muT} \colon T^2 \Rightarrow T$ and $\gls{1T} \colon \mr{id} \Rightarrow T$, satisfying unit and associativity axioms saying that the following diagrams commute for all objects $X$ of $\sfC$:
\[\begin{tikzcd} TX \rar{1^T_{TX}} \ar[equals]{rd} &[5pt] T^2X \arrow{d}[description]{\mu^T_{TX}} &[5pt] TX \lar[swap]{T(1^T_X)} \ar[equals]{ld} \\[5pt]
& TX & \end{tikzcd} \qquad \begin{tikzcd} T^3X \rar{\mu^T_{TX}} \dar[swap]{T(\mu^T_X)} &[5pt] T^2X \dar{\mu^T_X} \\[5pt]
T^2X \rar{\mu^T_X} & TX.
\end{tikzcd}\]

\newglossaryentry{algT}{%
	name={\ensuremath{\Alg_T}},
description={Category of $T$-algebras},
type=symbols
}
\begin{definition}A \emph{$T$-algebra}\index{$T$-algebra} $\gX = (X, a^T_X)$ consists of an object $X \in \sfC$ and a structure map $\alpha^T_X \colon TX \to X$ such that the following diagrams commute:
\[\begin{tikzcd} X \rar{1^T_X} \ar[equals]{rd} & TX \dar{\alpha^T_X} \\
 & X\end{tikzcd} \qquad \begin{tikzcd} T^2 X \rar{\mu^T_X} \dar[swap]{T(\alpha^T_X)} & TX \dar{\alpha^T_X} \\
TX \rar{\alpha^T_X} & X.\end{tikzcd}\]

A \emph{morphism of $T$-algebras} $\gX \to \gY$ is a morphism $f \colon X \to Y$ in $\sfC$ such that $\alpha^T_Y \circ T(f) = f \circ \alpha^T_X$. The category of $T$-algebras is denoted $\gls{algT}(\sfC)$.\end{definition}

Monads are closely related to adjunctions. From any adjunction
\[\begin{tikzcd}
\sfC \arrow[shift left=.5ex]{r}{F} & \arrow[shift left=.5ex]{l}{G} \sfD
\end{tikzcd}\]
we can obtain a monad $GF \colon \sfC \to \sfC$, where the counit natural transformation $\epsilon \colon FG \to \mr{id}$ is used to define $\mu^{GF}$ as $G\epsilon_F \colon GFGF \Rightarrow GF$ and the unit natural transformation $\eta \colon \mr{id} \Rightarrow GF$ gives $1^{GF}$.

\newglossaryentry{uT}{%
	name={\ensuremath{U^T}},
	description={Forgetful functor on $T$-algebras},
	type=symbols
}
\newglossaryentry{fT}{%
	name={\ensuremath{F^T}},
	description={Free $T$-algebra functor},
	type=symbols
}
The monad $T$ is a special case of this. The forgetful functor $\gls{uT} \colon \Alg_T(\sfC) \to \sfC$ given by $U^T(\gX) = X$ has a left adjoint given by sending an object $X$ of $\sfC$ to the free $T$-algebra $\gls{fT}(X)$, having underlying object $TX$ and structure map $\mu^T_X \colon T^2 X \to TX$. We have $T \cong U^T F^T$ and the monadic adjunction
\[\begin{tikzcd}
\sfC \arrow[shift left=.5ex]{r}{F^T} & \arrow[shift left=.5ex]{l}{U^T} \Alg_T(\sfC).
\end{tikzcd}\]

\subsection{Sifted monads} \label{sec:sifted-complete-cocomplete} We shall restrict our attention to monads which are ``finitary'' in the following sense.

\begin{definition}
	A monad $T$ on $\sfC$ is \emph{sifted}\index{monad!sifted}\index{sifted!monad} if the underlying functor $T \colon \sfC \to \sfC$ preserves all sifted colimits.
\end{definition}

It is easier to work with a sifted monad for two reasons: (i) its category of algebras is well-behaved, and (ii) we can use a ``density argument'' to construct functors out of the category of $T$-algebras.

\subsubsection{Categorical properties}
We start by establishing basic properties of the categories of algebras over a sifted monad. Though the following lemmas are well-known, see e.g.\ Exercise II of Section VI.2 of \cite{MacLane} and Proposition II.7.4 of \cite{EKMM}, we believe it is helpful to give a proof.

\begin{lemma}\label{lem:ut-preserves-sifted} 
The category $\Alg_T(\sfC)$ has sifted colimits, which are preserved by the forgetful functor $U^T \colon \Alg_T(\sfC) \to \sfC$.\end{lemma}

\begin{proof}
Let $i \mapsto \gX_i \colon \cat{I} \to \Alg_T(\sfC)$ be a sifted diagram. Applying the forgetful functor to $\sfC$ gives the diagram $i \mapsto X_i$ in $\sfC$ with colimit $\colim_{i \in \cat{I}} X_i$, and as $T$ preserves sifted colimits we may endow this with a $T$-algebra structure via
	\[T(\colim\limits_{i \in \cat{I}} X_i) \cong \colim_{i \in \cat{I}} TX_i \lra \colim_{i \in \cat{I}} X_i.\]
This $T$-algebra satisfies the universal property for the colimit.
\end{proof}

The reflexive coequalizer diagram is sifted, and there is a reflexive coequalizer diagram
	\begin{equation*}
\begin{tikzcd}
	F^T(T(X)) \arrow[shift left=1ex]{r} \arrow[shift left=-1ex]{r} & F^T(X) \lar 
\end{tikzcd}
	\end{equation*}
in $\Alg_T(\sfC)$, where the maps are given by $F^T(\alpha^T_X)$ and $\mu_T$ and the reflection is given by $F^T(\iota_X)$. The counit $F^T(X) = F^T U^T(\gX) \to \gX$ coequalizes this diagram, and by the above lemma one can check whether this is the coequalizer after applying $U^T$: it is, because the resulting coequalizer diagram in $\sfC$ is split. The diagram 
	\begin{equation}\label{eq:CanonPres}
\begin{tikzcd}
	F^T(T(X)) \arrow[shift left=1ex]{r} \arrow[shift left=-1ex]{r} & F^T(X) \lar \rar & \gX
\end{tikzcd}
	\end{equation}
is called the \emph{canonical presentation}\index{canonical presentation} of $\gX$ in $\Alg_T(\sfC)$.

\begin{lemma}\label{lem:algt-complete-cocomplete} The category $\Alg_T(\sfC)$ is complete and cocomplete.\end{lemma}

\begin{proof}
Limits are calculated in $\sfC$, and inherit a $T$-algebra structure in a standard way.
	
	For colimits, first note that free diagrams,  i.e.\ those functors $\cat{J} \to \Alg_T(\sfC)$ factoring through $F^T \colon \sfC \to \Alg_T(\sfC)$, have colimits, as $F^T$ is a left adjoint and so preserves colimits. 
A general diagram $j \mapsto \gX_j$ in $\Alg_T(\sfC)$ may be reduced to reflexive coequalizers and free diagrams by means of the canonical presentation. 
 Explicitly, the colimit of $j \mapsto \gX_j \in \Alg_T(\sfC)$ is the coequalizer of the
	diagram
	\begin{equation*}
	\begin{tikzcd}
	F^T(\colim\limits_{j \in \cat{J}} TX_j) \arrow[shift left=.5ex]{r} \arrow[shift left=-.5ex]{r} &
	F^T(\colim\limits_{j \in \cat{J}} X_j).
	\end{tikzcd}\qedhere
	\end{equation*}
\end{proof}

We will occasionally add superscripts $T$ to colimits or limits when we want to emphasize that they are taken in $\Alg_T(\sfC)$. For example, we use the notation $\sqcup^T$ to denote the coproduct in the category $\Alg_T(\sfC)$, or $\vee^T$ if we want to stress that $\sfC$ is pointed.

\subsubsection{Functors out of $\Alg_T(\sfC)$}
\label{sec:funct-out-eilenb}

\newglossaryentry{kleis}{%
	name={\ensuremath{\Kleis_T}},
	description={Kleisli category of a monad $T$},
	type=symbols
}
The full subcategory of $\Alg_T(\sfC)$ on the image of $F^T$ is called the
category of \emph{Kleisli algebras}\index{Kleisli algebras} for $T$ and denoted $\gls{kleis}(\sfC)$. Equivalently, the Kleisli category has the same objects as $\sfC$, but the morphisms from $X$ to $Y$ are given by $\sfC(X,TY)$, with composition of morphisms defined using the monad structure on $T$. 

If $G \colon \Alg_T(\sfC) \to \sfD$ is a functor, then the composition $H \coloneqq G \circ F^T \colon \sfC \to \sfD$ inherits a natural transformation
\[\mu^H \colon H \circ T = (G \circ F^T) \circ (U^T \circ F^T) \Rightarrow G \circ
F^T = H\]
coming from the counit of the adjunction. The natural
transformation $\mu^H$ satisfies the axioms of
a \emph{right $T$-module functor}. This means it is a right module over $T$ in $\cat{Fun}(\sfC,\sfD)$, i.e.\ the following diagrams commutes for all objects $X$ of $\sfC$:
\[\begin{tikzcd} HT^2X \rar{\mu^H_{TX}} \dar[swap]{H \mu^T_X} & HTX \dar{\mu^H_{X}} \\
HTX \rar{\mu^H_{X}} & HX\end{tikzcd} \qquad \begin{tikzcd} HX \rar{H1^T_X} \arrow[equal] {rd}&HTX \dar{\mu^H_X} \\
& HX.\end{tikzcd} \]

Notice that the above construction only used the restriction of $G$ to the Kleisli category: if $G \colon \Kleis_T(\sfC) \to \sfD$ is a functor, then $H = G \circ F^T \colon \sfC \to \sfD$ acquires the structure of a right $T$-module functor.

\begin{lemma}\label{lem:Kleisli}
There is an equivalence of categories between the category of functors $\Kleis_T(\sfC) \to \sfD$ and the category of right $T$-module functors $\sfC \to \sfD$.
\end{lemma}

\begin{proof}We described the forwards direction above on objects of $\cat{Fun}(\Kleis_T(\sfC),\sfD)$. It clearly extends to morphisms, i.e. natural transformations.
	
For the converse direction, if $H \colon \sfC \to \sfD$ is given the structure of a right module functor for a monad $T$ on $\sfC$, then we define a functor $G \colon \Kleis_T(\sfC) \to \sfD$ on objects by $G(F^T(X)) = H(X)$. On a (possibly non-free) morphism $f \colon F^T(X) \to F^T(Y)$, with adjoint $f' \colon X \to TY$, we define a morphism $G(f) \in \sfD(HX,HY)$ as the composition
\begin{equation*}
H(X) \xrightarrow{H(f')} H(T(Y)) \overset{\mu_H}\lra H(Y).
\end{equation*}
This defines the functor $G$, and this construction clearly extends to natural transformations. These two constructions induce the stated equivalence of categories.
\end{proof}

Since the Kleisli category generates the category of $T$-algebras
under sifted colimits (or even just reflexive coequalizers), there is a similar description of functors out of $\Alg_T(\sfC)$.

\begin{proposition}\label{prop:extend-by-density} Suppose that $T$ is a sifted monad on $\sfC$ and $\sfD$ has all sifted colimits. Then there is an equivalence of categories between the category of functors $\Alg_T(\sfC) \to \sfD$ preserving sifted colimits, and the category of right $T$-module functors $\sfC \to \sfD$ preserving sifted colimits.\end{proposition}

Most relevant for us will be the reverse direction, so let us make it explicit. Given a right $T$-module functor $H \colon \sfC \to \sfD$ we define $G(\gX)$ as the reflexive coequalizer
\[\begin{tikzcd}
	H(T(X)) \arrow[shift left=1ex]{r} \arrow[shift left=-1ex]{r} & H(X) \rar \lar &
	G(\gX)
\end{tikzcd}\]
in $\Alg_T(\sfC)$, of the maps $\mu_H$ and $H(\alpha^T_X)$ with reflection $H(\iota_X)$. We call this technique of defining functors out of $\Alg_T(\sfC)$ \emph{extension by density under sifted colimits}.\index{extension by density}

\subsubsection{Simplicial monads}\label{sec:simplicial-monads} Suppose now that $\sfC$ is a category enriched in $\cat{sSet}$.

\begin{definition}A \emph{simplicial monad}  is a unital monoid in the category of simplicially enriched functors $T\colon \sfC \to \sfC$.\index{monad!simplicial}\end{definition}

We now show that if $T$ is sifted and simplicial, then $\Alg_T(\sfC)$ is simplicially enriched and copowered over $\cat{sSet}$. In particular, the copowering is constructed using extension by density under sifted colimits.

\begin{lemma}\label{lem:algt-simplicial} If $T$ is simplicial the category $\Alg_T(\sfC)$ is enriched over $\cat{sSet}$, and if additionally $T$ is sifted then $\Alg_T(\sfC)$ is copowered over $\cat{sSet}$. The copowering satisfies $K \times F^T(X) \cong F^T(K \times X)$ naturally in $K$ and $X$.
\end{lemma}

\begin{proof}
	For $\gX, \gY \in \Alg_T(\sfC)$ we define the simplicial set $\Map_{\Alg_T}(\gX,\gY)$ as the equalizer
	\begin{equation*}
	\begin{tikzcd}
	\Map_{\Alg_T}(\gX,\gY) \rar & \Map_{\sfC}(X, Y) \arrow[shift left=.5ex]{r} \arrow[shift left=-.5ex]{r} &
	\Map_{\sfC}(T(X),Y)
	\end{tikzcd}
	\end{equation*}
	in $\cat{sSet}$, where the top and bottom maps are given by 
	\[\alpha^T_Y \circ T \colon \Map_{\sfC}(X,Y) \lra \Map_{\sfC}(TX,TY)\lra \Map_{\sfC}(T(X),Y),\]
	\[\alpha^T_X \colon \Map_{\sfC}(X,Y) \lra \Map_{\sfC}(TX,Y).\]
	
	For a simplicial set $K$ and $\gX \in \Alg_T(\sfC)$ we define the copowering $K \times \gX$ by the reflexive coequalizer
	\begin{equation*}
	\begin{tikzcd}
	F^T(K \times T(X)) \arrow[shift left=1ex]{r} \arrow[shift left=-1ex]{r} & 
	F^T(K \times X) \rar \lar & K \times \gX
	\end{tikzcd}
	\end{equation*}
	in $\Alg_T(\sfC)$, of the map $F^T(K \times \alpha^T_X)$ and the map
	\[F^T(K \times T(X)) \xrightarrow{F^T(\nu)} F^T(T(K \times X)) \overset{\mu_T}\longrightarrow F^T(K \times X)\]
	where $\nu \colon K \times TX \to T(K \times X)$ is adjoint to 
	\[K \lra \Map_{\sfC}(X, K \times X) \lra \Map_{\sfC}(TX, T(K \times X)).\]
The reflection is given by $F^T(K \times \iota_X) \colon F^T(K \times X) \to F^T(K \times T(X))$.
	From this construction it follows that $K \times F^T(X) \cong F^T(K \times X)$.
	
	An elementary argument gives a natural isomorphism $\Map_{\Alg_T}(K \times \gX,\gY) \cong \Map_{\Alg_T}(\gX, \gY)^K$, showing that $\times \colon \mathsf{sSet} \times \Alg_T(\sfC) \to \Alg_T(\sfC)$ is indeed a copowering. 
\end{proof}

If $T$ is simplicial, then so is the Kleisli category $\Kleis_T(\sfC)$: the morphism spaces are given by $\Map_{\Kleis_T}(X,Y) = \Map_\sfC(X,TY)$ and composition uses the simplicial monad structure on $T$. If $\sfD$ is simplicially enriched, then the proof of Lemma \ref{lem:Kleisli} with morphism sets replaced by morphism spaces yields an equivalence between the category of simplicial functors $\Alg_T(\sfC) \to \sfD$ and right simplicial $T$-module functors $\sfC \to \sfD$. Using this, Proposition \ref{prop:extend-by-density} upgrades to an equivalence between simplicial functors preserving sifted colimits, and right simplicial $T$-module functors preserving sifted colimits.

\subsection{The basepoint monad}
\newglossaryentry{initial}{%
	name={\ensuremath{\binit}},
	description={Initial objects},
	type=symbols
}
\newglossaryentry{terminal}{%
	name={\ensuremath{\bterm}},
	description={Terminal objects},
	type=symbols
}
\newglossaryentry{sfcstar}{%
	name={\ensuremath{\sfCstar}},
	description={Category of pointed objects in $\sfC$},
	type=symbols
}
\newglossaryentry{uplus}{%
	name={\ensuremath{U^+}},
	description={Forgetful functor on pointed objects},
	type=symbols
}
\newglossaryentry{fplus}{%
	name={\ensuremath{F^+}},
	description={Free pointed object functor},
	type=symbols
}
\newglossaryentry{plus}{%
	name={\ensuremath{+}},
	description={Basepoint monad},
	type=symbols
}
\label{sec:basepoint-monad} \index{monad!basepoint} We have assumed that $\sfC$ is complete and cocomplete, so in particular it has an initial object $\gls{initial}$ and terminal object $\gls{terminal}$.  We have however not assumed that the initial and terminal objects coincide, i.e.\ that $\sfC$ is pointed. We can make it pointed by considering instead $\gls{sfcstar}$, the undercategory of $\bterm$. The functor $\gls{uplus}$ which forgets the reference map from the terminal object is part of an adjunction
\[\begin{tikzcd}
\sfC \arrow[shift left=.5ex]{r}{F^+} & \sfC_\ast \arrow[shift left=.5ex]{l}{U^+} 
\end{tikzcd}\]
where $\gls{fplus}$ is given by $X \mapsto (\bterm \to X \sqcup \bterm)$. The composition $U^+ F^+ \colon \sfC \to \sfC$ is thus given by taking coproduct with the terminal object, and this composition forms a monad on $\sfC$.  We denote this monad by $\gls{plus}$. The underlying functor of $(-)_+$ is defined in terms of colimits and finite products (the terminal object is the empty product), so preserves filtered colimits and reflexive coequalizers. The adjunction is easily seen to be monadic, i.e.\ there is an equivalence of categories $\sfC_{\ast} \cong \Alg_{(-)_+}(\sfC)$. If the category $\sfC$ is already pointed, i.e.\ the unique map $\binit \to \bterm$ is an isomorphism, then $\sfCstar \cong \sfC$ and $+$ is the identity functor.

In the remainder of this section, we will construct a tensor product on $\sfC_\ast$ from a tensor product on $\sfC$. Whenever we consider $\sfC_*$ as a monoidal category we shall endow it with this monoidal structure.

\begin{lemma}
	Suppose $(\sfC,\otimes,\bunit_\sfC)$ is a category with a $k$-monoidal structure. Then there is a $k$-monoidal structure $(\sfC_\ast,\owedge,\bunit_{\sfC_\ast})$ on $\sfC_\ast$ such that $(-)_+ \colon \sfC \to \sfC_\ast$ is strong monoidal and $U^+ \colon \sfC_\ast \to \sfC$ is lax monoidal. Moreover, if the monoidal structure on $\sfC$ is closed (on both sides), so is that on $\sfC_\ast$.
\end{lemma}

\begin{proof}We will use general results about \emph{monoidal monads}. If $(\sfC, \otimes, \bunit)$ is a monoidal category then a monad $T \colon \sfC \to \sfC$ is called \emph{monoidal} if it is equipped with a lax monoidality such that $T^2 \Rightarrow T$ and $\mr{id} \Rightarrow T$ are monoidal natural transformations. It is well-known \cite[p.\ 30]{DayII} that in this situation the category $\mathsf{Kleis}_T(\sfC)$ inherits a monoidal structure. Let us briefly recall the construction. Using Lemma \ref{lem:Kleisli}, we may define a functor
	\newglossaryentry{owedge}{%
		name={\ensuremath{-\owedge-}},
		description={Induced tensor product on pointed objects},
		type=symbols
	}
	\[- \otimes_T - \colon \mathsf{Kleis}_T(\sfC) \times \mathsf{Kleis}_T(\sfC) = \mathsf{Kleis}_{T \times T}(\sfC \times \sfC) \lra \mathsf{Kleis}_T(\sfC)\]
	by giving a right $T \times T$-module functor $F \colon \sfC \times \sfC \to \mathsf{Kleis}_T(\sfC)$. Taking $F(-, -) = F^T(- \otimes -)$ and equipping it with the right $T \times T$-module structure
	\[F^T(T(-) \otimes T(-)) \lra F^T(T(- \otimes -)) \lra F^T(- \otimes -),\]
	given by the lax monoidality of $T$ followed with the right $T$-module structure of $F^T$, therefore defines $- \otimes_T -$. A similar discussion produces associators and left and right unitors, verifies the triangle and pentagon axioms, and yields a natural transformation $U^T(-) \otimes U^T(-) \Rightarrow U^T(- \otimes_T -)$.
	
	If in addition $- \otimes -$ preserves sifted colimits, which it always does in cases we will consider by Lemma \ref{lem.cat-sifted}, then by Proposition \ref{prop:extend-by-density} this definition extends to a monoidal structure $\otimes_T$ on $\Alg_T(\sfC)$. By construction, the functor $F^T \colon \sfC \to \Alg_T(\sfC)$ has a strong monoidality and its right adjoint $U^T$ a lax monoidality. Finally, if the monoidality on $T$ is braided or symmetric, then $\otimes_T$ inherits a braiding or symmetry.
	
	\vspace{.5em}
	
	If $(\sfC, \otimes, \bunit)$ is a monoidal (or braided or symmetric monoidal) structure on $\sfC$, then the monad $+$ may be given a lax monoidality, using the morphisms
	\[(X \sqcup \bterm) \otimes (Y \sqcup \bterm) \cong (X \otimes Y) \sqcup ((X \otimes \bterm) \sqcup (\bterm \otimes Y) \sqcup (\bterm \otimes \bterm)) \lra (X \otimes Y) \sqcup \bterm\]
	which is the identity map on the first summand and the unique map to $\bterm$ on the second summand. By the above considerations, this defines a monoidal (or braided or symmetric monoidal) structure on $\sfC_{\ast} \cong \Alg_{(-)_+}(\sfC)$, which we shall denote as
	\[\gls{owedge} \colon \sfC_* \times \sfC_* \lra \sfC_*.\]
	If $\sfC$ is closed then $\sfC_\ast$ with $\owedge$ is also closed, with enriched right adjoint to $- \owedge X$ given by the equalizer
	\[\begin{tikzcd}
	\Hom_{\sfC_\ast}(X,Y) \rar & \Hom_{\sfC}(X_+, Y) \arrow[shift left=.5ex]{r} \arrow[shift left=-.5ex]{r} &
	\Hom_{\sfC}(X,Y). \qedhere
	\end{tikzcd}\]
\end{proof}

\begin{example}	When applied to $\cat{sSet}$ or $\cat{Top}$ with $\otimes$ the cartesian product, $\owedge$ is the usual smash product. When applied to a pointed category, $\owedge$ is $\otimes$. \end{example}

\subsection{Morphisms of monads and indecomposables} \label{sec:indecomposables} It is helpful to define the notion of a morphism between monads on different categories, rather than demand they are defined on the same category.

\begin{definition}Let $T$ be a monad on a category $\cat{C}$ and $T'$ be a monad on a category $\cat{C}'$. Then a \emph{morphism of monads}\index{monad!morphism of} $(F,\phi) \colon (\cat{C},T) \to (\cat{C}',T')$ is a pair of a functor $F \colon \cat{C}' \to \cat{C}$ and a natural transformation $\phi \colon TF \Rightarrow FT'$ so that the following diagrams commute for all objects $X$ of $\cat{C}$
	\[\begin{tikzcd} T^2FX \rar{T\phi_X} \dar{\mu^{T}_{FX}} & TFT'X \rar{\phi_{T'X}} & F(T')^2X \dar{F\mu^{T'}_X} \\
	TFX \arrow{rr}{\phi_X} & & FT'X,\end{tikzcd} \qquad \begin{tikzcd} FX \dar[swap]{1^{T}_{FX}} \rar{F1^{T'}_X} & FT'X\\
	TFX. \arrow{ru}[swap]{\phi_X} & \end{tikzcd}  \]
\end{definition}

This definition is made so that the following holds:

\begin{lemma}A morphism of monads $(F,\phi) \colon (\cat{C},T) \to (\cat{C}',T')$ induces a functor
	\[(F,\phi)^* \colon \Alg_{T'}(\sfC') \lra \Alg_{T}(\sfC)\]
satisfying $U^{T}(F,\phi)^* = F U^{T'}$.
\end{lemma}

\begin{proof}
If $\gX = (X,a^{T'}_X)$ is a $T'$-algebra, then $FX$ is a $T$-algebra, with $a^{T}_{FX}$ given by 
\[TFX \overset{\phi_X}\lra FT'X \xrightarrow{Fa^{T'}_X} FX.\]
To verify this, we observe that in the diagram
\[\begin{tikzcd} T^2FX \arrow{rr}{\mu^T_{FX}}  \dar{T\phi_X} & & TFX \dar{\phi_X} \\
TFT'X \rar{\phi_{T'X}}  \dar{TF\alpha^{T'}_X} & F(T')^2X \rar{F\mu^{T'}_X} \dar[description]{FT'a^{T'}_X} & FT'X \dar{Fa^{T'}_X} \\
TFX \rar{\phi_X} & FT'X \rar{Fa^{T'}_X} & FX \end{tikzcd}\]
the top rectangle commutes because $\phi$ is a morphism of monads, the bottom-left square commutes because $\phi$ is a natural transformation, and the bottom-right square commutes because $X$ is a $T'$-algebra. Similarly, the folllowing diagram commutes because $F$ is a functor and $\phi$ is a morphism of monads
\[\begin{tikzcd} FX \rar{1^T_{FX}} \dar[equal] \arrow{rd}[description]{F1^{T'}_X} & TFX \dar{\phi_X} \\
FX & FT'X \lar{F\alpha^{T'}_X}.\end{tikzcd} \]
This construction of a $T$-algebra structure on $FX$ is natural in $\gX$, and yields the required functor.
\end{proof}

If the functor $F$ is a right adjoint, then $(F,\phi)^*$ is a right adjoint too, as follows.

\newglossaryentry{Rphi}{%
	name={\ensuremath{(R,\phi)^*}},
	description={Right adjoint in change-of-category-and-monad adjunction},
	type=symbols
}
\newglossaryentry{Lphi}{%
	name={\ensuremath{(L,\phi)_*}},
	description={Left adjoint in change-of-category-and-monad adjunction},
	type=symbols
}

\begin{lemma}\label{lem:morphism-monads-adjunction} Suppose that $T'$ is a sifted monad and $\cat{C}'$ has all sifted colimits. Then if $(R,\phi) \colon (\cat{C},T) \to (\cat{C}',T')$ is a morphism of monads and $R$ has a left adjoint $L$, then there is a \emph{change-of-category-and-monad adjunction}
	\begin{equation*}
	\begin{tikzcd}\Alg_T(\sfC) \arrow[shift left=.5ex]{r}{\gls{Lphi}} & \Alg_{T'}(\sfC'). \arrow[shift left=.5ex]{l}{\gls{Rphi}}
	\end{tikzcd}
	\end{equation*}
satisfying $U^{T}(R,\phi)^* = R U^{T'}$ and $(L,\phi)_*F^T = F^{T'}L$. Moreover, if $T$ and $T'$ are simplicial and if $L \dashv R$ is a simplicial adjunction, then so is $(L,\phi)_* \dashv (R,\phi)^*$.
\end{lemma}

\begin{proof}The natural transformation
	\[LT \xrightarrow{LT\eta} LTRL \xrightarrow{L\phi L} LRT'L \xrightarrow{\epsilon T'L} T'L,\]
	induces a right $T$-module functor structure on $F^{T'}L$. As this functor also preserves sifted colimits, because it is the composition of two left adjoints, we may invoke Proposition \ref{prop:extend-by-density} to obtain a functor $(L,\phi)_*$ satisfying  $(L,\phi)_*F^T = F^{T'}L$. The statement about simplicial enrichments is justified by the considerations of Section \ref{sec:simplicial-monads}.
\end{proof}

\begin{example}[The change-of-monad adjunction]
\label{exam:change-monads}
\newglossaryentry{phiup}{%
	name={\ensuremath{\phi^*}},
	description={Change-of-monads along $\phi$},
	type=symbols
}
Let us specialize to $L = R = \mr{Id}_\cat{C}$. Suppose $(\mr{Id}_\cat{C},\phi) \colon (\cat{C},T) \to (\cat{C},T')$ is a morphism of monads, which we abbreviate to $\phi \colon T \to T'$, then we obtain a \emph{change-of-monad} adjunction
\begin{equation*}
\begin{tikzcd}\Alg_T(\sfC) \arrow[shift left=.5ex]{r}{\phi_*} & \Alg_{T'}(\sfC) \arrow[shift left=.5ex]{l}{\phi^*}
\end{tikzcd}
\end{equation*}
satisfying $U^T\phi^* = U^{T'}$ and $\phi_* F^{T} = F^{T'}$.\end{example}

\begin{example}[The monadic adjunction] Let $\mr{Id}$ denote the identity monad on $\sfC$, so that $\Alg_\mr{Id}(\sfC) = \sfC$.  Then the change-of-monad adjunction for the monad map $\mr{Id} \to T$ is precisely the monadic adjunction $F^T \dashv U^T$.\end{example}

\subsubsection{Augmented monads and indecomposables}\label{sec:AugMonad}
Recall that $+$ denotes the basepoint monad on $\sfC$, as described in Section \ref{sec:basepoint-monad}.

\newglossaryentry{augmentation}{%
	name={\ensuremath{\epsilon}},
	description={Augmentation},
	type=symbols
}

\begin{definition}\label{def:augmentation} An \emph{augmentation}\index{augmentation} on a monad $T$ on $\sfC$ is a morphism of monads $\gls{augmentation} \colon T \to +$.\index{monad!augmented} An \emph{augmented monad} is a monad with an augmentation.
\end{definition}

\begin{example}If $\sfC$ is pointed then $\sfC = \sfCstar$ and $+ = \mr{Id}$, so an augmentation is a morphism of monads $T \to \mr{Id}$.\end{example}

\newglossaryentry{QT}{%
	name={\ensuremath{Q^T}},
	description={$T$-algebra indecomposables functor},
	type=symbols
}
\newglossaryentry{ZT}{%
	name={\ensuremath{Z^T}},
	description={Trivial $T$-algebra functor},
	type=symbols
}

Invoking Example \ref{exam:change-monads}, we obtain:

\begin{definition}\label{defn:indecomposables}
	Let $T$ be a sifted monad on a category $\sfC$ which has sifted colimits and $\epsilon \colon T \to +$ be an augmentation. The \emph{indecomposables}\index{indecomposables} functor $\gls{QT} \colon \Alg_T(\sfC) \to \sfC_\ast$ is given by
	\begin{equation*}
	Q^T(\gX) \coloneqq \epsilon_*(\gX).
	\end{equation*}
\end{definition}

The composition of morphisms of monads $\mr{Id} \to T \to +$ on $\cat{C}$ gives rise to a pair of adjunctions
\begin{equation*}
\begin{tikzcd}
\sfC = \Alg_{\mr{Id}}(\sfC) \arrow[shift left=.5ex]{r}{F^T} & \Alg_T(\sfC) \arrow[shift left=.5ex]{l}{U^T} \arrow[shift left=.5ex]{r}{Q^T}
& \Alg_+(\sfC) = \sfCstar \arrow[shift left=.5ex]{l}{Z^T},
\end{tikzcd}
\end{equation*}
where we have written $\gls{ZT} = \epsilon^* \colon \sfCstar \to \Alg_T(\sfC)$ for the right
adjoint to $Q^T = \epsilon_*$. The right adjoint $Z^T$ is the \emph{trivial $T$-algebra} functor.

\begin{corollary}\label{cor:properties-of-indecomposables}
	The functors $F^T \colon \sfC \to \Alg_T(\sfC)$ and   $Q^T \colon \Alg_T(\sfC) \to \sfCstar$ preserve colimits.  The composition   $Q^T \circ F^T \colon \sfC \to \sfCstar$ is naturally isomorphic to the functor $F^+$ (which takes coproduct with the terminal object, regarded as the basepoint).
\end{corollary}

\begin{proof}
	The functors preserve colimits because they are left adjoints. The composition $Z^T \circ U^T$ is the functor which forgets the basepoint and hence its left adjoint is as asserted.
\end{proof}

Let us unravel the definition of $Q^T$. As $Q^T$ preserves colimits, applying it to the canonical presentation \eqref{eq:CanonPres} of a $T$-algebra $\gX$ gives a reflexive coequalizer
\[\begin{tikzcd}
+T(X) \arrow[shift left=1ex]{r} \arrow[shift left=-1ex]{r} & +X \rar \lar &
Q^T(\gX)
\end{tikzcd}\]
in $\sfCstar$. The maps are given by $\mu_+ \circ (+\epsilon)$ and $+\alpha^T_X$, and the reflection is given by $+\iota_X$. Roughly speaking, $Q^T(X)$ is obtained from $+X$ by collapsing to the basepoint everything obtained by applying a non-identity operation in $T$.

\subsubsection{Indecomposables and morphisms of monads}

Let $(R,\phi) \colon (\cat{C},T) \to (\cat{C}',T')$ be a morphism of sifted monads, and suppose that $R$ has a left adjoint $L$. Since $R$ is a right adjoint, it preserves $\bterm$ so the morphisms
\[R(X) \sqcup \bterm \cong R(X) \sqcup R(\bterm)  \lra R(X \sqcup \bterm)\]
induce a natural transformation $\upsilon \colon + \circ R \Rightarrow R \circ +$ of functors $\cat{C}' \to \cat{C}$, yielding a morphism of monads $(R,\upsilon) \colon (\cat{C},+) \to (\cat{C}',+')$. By Lemma \ref{lem:morphism-monads-adjunction} we obtain from this an adjunction
\begin{equation*}
\begin{tikzcd}\Alg_+(\sfC) = \sfC_\ast \arrow[shift left=.5ex]{r}{(L,\upsilon)_*} & \Alg_{+}(\sfC') = \sfC'_\ast. \arrow[shift left=.5ex]{l}{(R,\upsilon)^*}
\end{tikzcd}
\end{equation*}
Explicitly, we have
\[(R,\upsilon)^*(\bterm \to X') = \big(\bterm \cong R(\bterm) \to R(X')\big)\]
\[(L,\upsilon)_*(\bterm \to X) = \big(\bterm \to \colim(\bterm \leftarrow L(\bterm) \to L(X))\big).\]

\begin{lemma}\label{lem:ind-change-of-monads} Suppose we have a diagram
\[\begin{tikzcd} (\cat{C},T) \rar{(\mr{Id}_\cat{C},\epsilon)} \dar[swap]{(R,\phi)} & (\cat{C},+_{\sfC}) \dar{(R,\upsilon)} \\
(\cat{C}',T') \rar{(\mr{Id}_{\cat{C}'},\epsilon')} & (\cat{C}',+_{\sfC'})\end{tikzcd} \]
of morphisms of monads commuting up to natural isomorphism. Then there is a natural isomorphism
\[Q^{T'}((L,\phi)_*(-)) \cong (L,\upsilon)_*(Q^T(-)) : \Alg_T(\sfC) \lra \sfC_\ast.\]
\end{lemma}

\begin{proof}The diagram yields natural isomorphism of right adjoints $(R,\phi)^*(\mr{Id}_{\cat{C}'},\epsilon')^* \cong (\mr{Id}_{\cat{C}},\epsilon)^*(R,\upsilon)^*$, and hence a natural isomorphism of left adjoints.\end{proof}

\begin{example}\label{exam:morphism-over-id-indec}
Let us specialize to $L = R = \mr{Id}_\cat{C}$ and let $\phi \colon T \to T'$ be a morphism of monads. If $\epsilon' \colon T' \to +$  an augmentation of $T'$, then $\epsilon \coloneqq \epsilon' \circ \phi$ is an augmentation of $T$. Then there is a natural isomorphism
\[Q^{T'}(\phi_* (-)) \cong Q^T(-) \colon \Alg_T(\sfC) \lra \sfC_\ast.\]
\end{example}

\begin{example}[Indecomposables and change-of-diagram-category] \label{exam:change-groupoid-indec} We will often invoke Lemma \ref{lem:ind-change-of-monads} in the following setting. A functor $p \colon \sfG \to \sfG'$ induces by precomposition a functor $p^* \colon \sfC^{\sfG'} \to \sfC^\sfG$ which has a left adjoint given by left Kan extension $p_* \colon \sfC^{\sfG} \to \sfC^{\sfG'}$. Suppose we are given augmented sifted monads $T$ and $T'$ in $\sfC^{\sfG}$ and $\sfC^{\sfG'}$ and a commutative diagram of the form
\begin{equation*}
\begin{tikzcd}
(\cat{C},T) \rar{\epsilon_T} \dar[swap]{(p^*,\phi)} & (\cat{C},+_{\sfC}) \dar{(p^*,\upsilon)}\\
(\cat{C}',T') \rar{\epsilon_{T'}} & (\cat{C}',+_{\sfC'}).
\end{tikzcd}
\end{equation*}
This diagram of functors commutes up to natural isomorphism, so there is a natural isomorphism
\[Q^{T'}((p_*,\phi)_*(-)) \cong p_*(Q^T(-)) \colon \Alg_T(\sfC^{\sfG}) \lra \sfC_\ast^{\sfG'}.\]
When $\sfG'$ is the terminal category, so that $\sfC^{\sfG'} = \sfC$, this expresses $\indec^{T'}((p_*,\phi)_*(\gX)) \in \sfC_\ast$ as the colimit over $\sfG$ of the functor $\indec^T(\gX) \colon \sfG \to \sfC_*$. In contrast, there is usually no simple formula for $\indec^T((p^*,\phi)^*(\gX))$.\end{example}
\section{Monads associated to operads}\label{sec:operads} 

The most important class of monads on $\sfC$ to which we will apply the previous discussion are monads associated to operads in $\sfC$: these are always sifted, as we will show in Corollary \ref{cor:monad-from-operad-sifted}. We continue to suppose that $\sfC$ satisfies the axioms of Section \ref{sec:axioms-of-cats}, and that it is a $k$-monoidal category. In order to treat these varying levels of symmetry uniformly, we will consider operads as certain $k$-symmetric sequences, and we first define this notion. See \cite{FresseBook,HarperOperads} for further background.

\subsection{Symmetric, braided and ordered sequences} Operads will be objects of a category of \emph{$k$-symmetric sequences}\index{symmetric sequence}\index{$k$-symmetric sequence}, for $k \in \{1,2,\infty\}$. When describing modules over an associative algebra in terms of operads, there are no higher operations, and we may use $0$-symmetric sequences to encode these special operads.

\newglossaryentry{fsn}{%
	name={\ensuremath{\mathfrak{S}_n}},
	description={Symmetric group on $n$ elements},
	type=symbols
}
\newglossaryentry{braidn}{%
	name={\ensuremath{\beta_n}},
	description={Braid group on $n$ strands},
	type=symbols
}
\begin{definition}
Write $\ul{n} \coloneqq \{1,\ldots,n\}$ for $n \geq 0$, so $\ul{0}=\varnothing$. We define four groupoids:
	\begin{enumerate}[(i)]
		\item $\cat{FB}_\infty$ has objects $\ul{n}$ for $n \geq 0$, which have automorphism groups given by the symmetric groups $\gls{fsn}$. There are no other morphisms.
		\item $\cat{FB}_2$ has objects $\ul{n}$ for $n \geq 0$, which have automorphism groups given by the braid groups $\gls{braidn}$. There are no other morphisms.
		\item $\cat{FB}_1$ has objects $\ul{n}$ for $n \geq 0$, and no non-identity morphisms.
		\item $\cat{FB}_0$ has one object $\ul{1}$, and no non-identity morphisms.
	\end{enumerate}
\end{definition}

\newglossaryentry{FBk}{%
	name={\ensuremath{\cat{FB}_k}},
	description={Indexing categories for symmetric sequences},
	type=symbols
}
Let us first discuss the cases $k \in \{1,2,\infty\}$, as the case $k=0$ is slightly different. Recall that a $k$-monoidal category is a monoidal category if $k\geq 1$, a braided monoidal category if $k\geq 2$, and a symmetric monoidal category if $k >2 $. For $k \geq 1$, the categories $\gls{FBk}$ have a monoidal structure $\otimes$ given on objects by $\ul{n} \otimes \ul{m} \coloneqq \ul{n+m}$ and on morphisms by disjoint union of permutations or braids. Furthermore, $\cat{FB}_2$ has a braiding 
\[\beta_{m,n} \colon \ul{m} \otimes \ul{n} \lra \ul{n} \otimes \ul{m}\]
given by the braid which crosses the left $m$ strands in front of the right $n$ strands (see \cite[p.\ 36]{JoyalStreet}), and the corresponding permutation gives a braiding on $\cat{FB}_\infty$ which is actually a symmetry. Thus $\cat{FB}_k$ has the structure of a $k$-monoidal category, and is in fact the free strict $k$-monoidal category on one generator. In order to give a uniform treatment of these three cases, we will write $G_n = \mathrm{Aut}_{\cat{FB}_k}(\ul{n})$.

To any category $\sfC$ we may associate the category of functors $\sfC^{\cat{FB}_k}$. If $\sfC$ is closed $k$-monoidal then Day convolution endows $\sfC^{\cat{FB}_k}$ with a closed $k$-monoidal structure. Explicitly, $\cX \otimes \cY$ is given by
\[(\cX \otimes \cY)(\ul{r}) \coloneqq \colim_{(\ul{k_1},\ul{k_2},f \colon \ul{k_1}\otimes \ul{k_2} \overset{\sim}\to \ul{r})} \cX(\ul{k_1}) \otimes \cY(\ul{k_2})\]
which we may identify as
\[\bigsqcup_{\substack{k_1, k_2 \geq 0 \\k_1+k_2 = r}} G_r \times_{G_{k_1} \times G_{k_2}} \cX(\ul{k_1}) \otimes \cY(\ul{k_2}),\]
where the homomorphism $G_{k_1} \times G_{k_2} \to G_{k_1 + k_2} = G_r$ is given by the monoidal structure of $\cat{FB}_k$.

There is a further piece of structure available on the categories $\cat{FB}_k$; it is perhaps most familiar in the case $k=2$, where it is given by cabling braids.

\begin{definition}
We define three groupoids:
	\begin{enumerate}[(i)]
		\item $\cat{FB}^{(2)}_\infty$ has objects $(\ul{n}, \ul{k_1}, \ldots, \ul{k_n})$ for $n, k_i \geq 0$. A morphism from $(\ul{n}, \ul{k_1}, \ldots, \ul{k_n})$ to $(\ul{n'}, \ul{l_1}, \ldots, \ul{l_{n'}})$ exists only if $n=n'$, in which case it consist of a permutation $\sigma \in \mathfrak{S}_n$ such that $k_{\sigma(i)} = l_i$ for all $i$, as well as permutations $\tau_i \in \mathfrak{S}_{k_i}$. Composition is given by\[(\sigma'; \tau'_1, \ldots, \tau'_n) \circ (\sigma; \tau_1, \ldots, \tau_n) = (\sigma' \circ \sigma; \tau'_{\sigma(1)} \circ \tau_1, \ldots, \tau'_{\sigma(n)} \circ \tau_n).\]

		\item $\cat{FB}^{(2)}_2$ has objects $(\ul{n}, \ul{k_1}, \ldots, \ul{k_n})$ for $n, k_i \geq 0$. A morphism from $(\ul{n}, \ul{k_1}, \ldots, \ul{k_n})$ to $(\ul{n'}, \ul{l_1}, \ldots, \ul{l_{n'}})$ exists only if $n=n'$, in which case it consist of a braid $\sigma \in \beta_n$ such that $k_{\sigma(i)} = l_i$ for all $i$, as well as braids $\tau_i \in \beta_{k_i}$. Composition is given by the same formula as above.
		\item $\cat{FB}^{(2)}_1$ has objects $(\ul{n}, \ul{k_1}, \ldots, \ul{k_n})$ for $n, k_i \geq 0$, and no non-identity morphisms.
	\end{enumerate}
\end{definition}

In each case there is a functor $c \colon \cat{FB}^{(2)}_k \to \cat{FB}_k$ given on objects by sending $(\ul{n}, \ul{k_1}, \ldots, \ul{k_n})$ to $\ul{k_1} \otimes \cdots \otimes \ul{k_n} $, and (for $k \geq 2$) given on morphisms by \emph{cabling}: sending the morphism $(\sigma; \tau_1, \ldots, \tau_n) \colon (\ul{n}, \ul{k_1}, \ldots, \ul{k_n}) \to (\ul{n}, \ul{l_1}, \ldots, \ul{l_n})$ to
\begin{equation*}
\begin{tikzcd}
\ul{k_1} \otimes \cdots \otimes \ul{k_n} \rar{\tau_1 \otimes \ldots \otimes \tau_n}&[12 pt] \ul{k_1} \otimes \cdots \otimes \ul{k_n} \rar{\sigma} & \ul{k_{\sigma(1)}} \otimes \cdots \otimes \ul{k_{\sigma(n)}} = \ul{l_1} \otimes \cdots \otimes \ul{l_n}.
\end{tikzcd}
\end{equation*}

Objects $\cX, \cY \in \sfC^{\cat{FB}_k}$ can be combined to give a functor
\begin{align*}
\cY^\cX \colon \cat{FB}^{(2)}_k &\lra \sfC\\
(\ul{n}, \ul{k_1}, \ldots, \ul{k_n}) &\longmapsto \cX(\ul{n}) \otimes \cY(\ul{k_1}) \otimes \cdots \otimes \cY(\ul{k_n}),
\end{align*}
where the morphism $(\sigma; \tau_1, \ldots, \tau_n) \colon (\ul{n}, \ul{k_1}, \ldots, \ul{k_n}) \to (\ul{n}, \ul{l_1}, \ldots, \ul{l_n})$ is sent to
\begin{equation*}
\begin{tikzcd}
\cX(\ul{n}) \otimes \cY(\ul{k_1}) \otimes \cdots \otimes \cY(\ul{k_n}) \rar{\sigma \otimes \tau_1 \otimes \cdots \otimes \tau_n} \ar[draw=none]{d}[name=X, anchor=center]{} &[10pt] \cX(\ul{n}) \otimes \cY(\ul{k_1}) \otimes \cdots \otimes \cY(\ul{k_n})             
\ar[rounded corners,
                to path={ -- ([xshift=2ex]\tikztostart.east)
                	|- (X.center) \tikztonodes
                	-| ([xshift=-2ex]\tikztotarget.west)
                	-- (\tikztotarget)}]{dl}[very near end]{\cX(n) \otimes \sigma}\\
              \cX(\ul{n}) \otimes \cY(\ul{k_{\sigma(1)}}) \otimes \cdots \otimes \cY(\ul{k_{\sigma(n)}}) \ar[equal]{r} & \cX(\ul{n}) \otimes \cY(\ul{l_1}) \otimes \cdots \otimes \cY(\ul{l_n})
\end{tikzcd}
\end{equation*}

\newglossaryentry{circ}{%
	name={\ensuremath{-\circ-}},
	description={Composition product},
	type=symbols
}
\begin{definition}
We define the \emph{composition product}\index{composition product}
\[\gls{circ} : \sfC^{\cat{FB}_k} \times \sfC^{\cat{FB}_k} \lra \sfC^{\cat{FB}_k}\]
by setting $\cX \circ \cY$ to be 
%
the left Kan extension of $\cY^\cX \colon \cat{FB}^{(2)}_k \to \sfC$ along $c \colon  \cat{FB}^{(2)}_k \to \cat{FB}_k$.
\end{definition}
Concretely, the value of $\cX \circ \cY$ on $\ul{r}$ 
is given by
\[\bigsqcup_{n=0}^\infty \left(\cX(\ul{n}) \otimes_{G_n}  \bigsqcup_{\substack{k_1, \ldots, k_n \geq 0\\k_1+\cdots+k_n=r}}  G_r \times_{G_{k_1} \times \cdots \times G_{k_n}} \left(\bigotimes_{i=1}^n \cY(\ul{k}_i)\right)\right).\]
The functor $\iota \colon \cat{FB}_k \to \sfC$ which sends $\ul{1}$ to $\bunit_\sfC$ and all other objects to $\binit_\sfC$ satisfies $\iota \circ \cX \cong \cX \cong \cX \circ \iota$. As for any colimit, to map out of a composition product one may equivalently provide equivariant maps out of $\cX(\ul{n}) \otimes \bigotimes_{i=1}^n \cY(\ul{k}_i)$. 

If $\sfC$ is $(k+1)$-monoidal, then the composition product $\circ$ on $\sfC^{\cat{FB}_k}$ is associative up to isomorphism and its associator isomorphism satisfies the pentagon identity. This follows by considering the evident generalisations $\cat{FB}^{(3)}_k$ and $\cat{FB}^{(4)}_k$ and the several different cabling functors between them. See Section 2.2 of \cite{FresseBook} or Section 4.2 of \cite{HarperOperads} for a (partial) proof that the composition product on $\sfC^{\cat{FB}_k}$ gives a monoidal structure, under the assumption that $\sfC$ is $\infty$-monoidal. 

It remains to treat the much easier case $k=0$ (Example \ref{exam:module-over-algebra-operad} may be clarifying). In this case, the category $\sfC^{\cat{FB}_0}$ is canonically identified with $\sfC$, and we define the composition product to be the $1$-monoidal structure on $\sfC$; thus it is $1$-monoidal.

When $k=2$, in the applications we have in mind the operads arise through a strong monoidal functor $s \colon \cat{sSet} \to \sfC$, and as such their spaces of $k$-ary operations will lie in the ``center'' of $\sfC$ under the monoidal structure (known as ``transparent'' in the mathematical physics literature):

\begin{definition}An object $X$ of a braided monoidal category $(\sfC,\otimes,\bunit)$ is \emph{central}\index{central} if for all $Y$, the braidings $\beta_{X,Y} \colon X \otimes Y \to Y \otimes X$ satisfy $\beta_{Y,X}\circ\beta_{X,Y} = \mr{id}_{X \otimes Y}$.\end{definition}

If $\sfC$ is $2$-monoidal then the full subcategory of $\sfC^{\cat{FB}_2}$ on those $2$-symmetric sequences $\cX$ such that each $\cX(n)$ is central in $\sfC$, is a monoidal category under the composition product. 

Let us collect the previous discussion in a proposition:

\begin{proposition}If $k \geq 0$ and $\sfC$ is $k$-monoidal, then there is a composition product on $\sfC^{\cat{FB}_k}$. This is a monoidal structure if $\sfC$ is $(k+1)$-monoidal. Furthermore, if $k=2$ and $\sfC$ is $2$-monoidal, this is a monoidal structure when restricted to the full subcategory of $\sfC^{\cat{FB}_2}$ consisting of those $2$-symmetric sequences that are objectwise central.
\end{proposition}

To emphasize the difference between the Day convolution and composition products, from now on we will write $\cat{FB}_k(\sfC)$ for the functor category $\sfC^{\cat{FB}_k}$ equipped with the composition product. Let us record some properties of the composition product.

\begin{lemma}\label{lem:odot-props}
The composition product $\circ$ has the following properties:
	\begin{enumerate}[(i)]
		\item \label{enum:odot-props-sifted}$\circ \colon  \cat{FB}_k(\sfC) \times \cat{FB}_k(\sfC) \to \cat{FB}_k(\sfC)$ preserves sifted colimits,
		\item \label{enum:odot-props-geomrel} $\circ \colon \cat{FB}_k(\sfC) \times \cat{FB}_k(\sfC) \to \cat{FB}_k(\sfC)$ preserves geometric realization.
	\end{enumerate}
\end{lemma}

\begin{proof}
As left Kan extension is given by a colimit, it commutes with all colimits and with geometric realization. It is therefore enough to show that
\[(\cX, \cY) \mapsto \cY^\cX \colon \sfC^{\cat{FB}_k} \times \sfC^{\cat{FB}_k} \lra \sfC^{\cat{FB}^{(2)}_k}\]
commutes with sifted colimits and with geometric realization. This is evident from the formula for $\cY^\cX$, and the fact that $\otimes_\sfC$ commutes with
sifted colimits by Lemma \ref{lem.cat-sifted} and geometric realization by Lemma \ref{lem.monoidal-geomrel}.\end{proof}
	  
\begin{notation}To ease notation, we will drop the underline from the object of $\cat{FB}_k$, i.e. denote $\cX(\ul{n})$ by $\cX(n)$.\end{notation}

For each $n \geq 0$, change-of-diagram-category along the inclusion $\{n\} \hookrightarrow \cat{FB}_k$ gives rise to an adjunction
\begin{equation}\label{eqn:symseq0adj} \begin{tikzcd} 
\sfC \arrow[shift left=.5ex]{r}{n_*} & \cat{FB}_k(\sfC) \arrow[shift left=.5ex]{l}{n^*},
\end{tikzcd}\end{equation}
where explicitly we have $n^*(\cX) = \cX(n)$ and $n_*(X)$ is given by the $k$-symmetric sequence assigning $G_n \times X$ to $n$, where $G_n = \mr{Aut}_{\cat{FB}_k}(n)$ acts on the first factor by translation, and $\binit_\sfC$ to all other objects of $\cat{FB}_k$. In this notation, the monoidal unit of $\cat{FB}_k(\sfC)$ can be written as $1_*(\bunit_\sfC)$.

\begin{definition}\label{def:symm-seq-application} We define a bifunctor $\cat{FB}_k(\sfC) \times \sfC \to \sfC$ by
\[(\cY,X) \longmapsto \cY(X) \coloneqq 0^*(\cY \circ 0_*(X)).\]
\end{definition}

Concretely, this is given by the formula
\begin{equation}\label{eqn:symmetric-pairing}
\cY(X) = \bigsqcup_{n \geq 0} \cY(n) \otimes_{G_n}
X^{\otimes n}.
\end{equation}

\begin{lemma}\label{lem:sseq-sifted}
For any $k$-symmetric sequence $\cY$ the functor $X \mapsto \cY(X)$ preserves sifted colimits and geometric realization. 
\end{lemma}
\begin{proof}
This follows from Lemma \ref{lem:odot-props}, as both $0^*$ and $0_*$ commute with all colimits and geometric realization. To spell this out a little, $X \mapsto X^{\otimes n}$ commutes with sifted colimits by Lemma \ref{lem.cat-sifted} and geometric realization by Lemma \ref{lem.monoidal-geomrel}, and the remaining constructions in \eqref{eqn:symmetric-pairing} commute with all colimits and geometric realization. 
\end{proof}

\subsection{Operads} Let $\sfC$ be a $(k+1)$-monoidal category as above, so that there is a monoidal category $(\cat{FB}_k(\sfC), \circ, 1_*(\bunit_\sfC))$. Alternatively, if $k=2$ one may assume $\sfC$ is $2$-monoidal and restrict to those $2$-symmetric sequences that are objectwise central. 

\newglossaryentry{co}{%
	name={\ensuremath{\cO}},
	description={Operad},
	type=symbols
}
\newglossaryentry{1o}{%
	name={\ensuremath{1_\cO}},
	description={Operad unit},
	type=symbols
}
\newglossaryentry{muo}{%
	name={\ensuremath{\mu_\cO}},
	description={Operad multiplication},
	type=symbols
}
\begin{definition}An \emph{operad}\index{operad} $\gls{co}$ in $\sfC$ is a unital monoid in $\cat{FB}_k(\sfC)$, with unit $\gls{1o} \colon 1_*(\bunit_\sfC) \to \cO$ and multiplication $\gls{muo} \colon \cO \circ \cO \to \cO$. 
\end{definition}

(This notion of course depends on the choice of $k$-symmetric monoidal structure on $\sfC$, though the notation does not reflect this.) Unwinding the definitions, an operad $\cO$ in $\sfC$ consists of a sequence of objects $\cO(n)$ with $G_n$-actions, for $n \geq 0$, and morphisms
\[1_\cO(1) \colon \bunit \lra \cO(1),\]
\[\mu_{\cO}(n;k_1,\ldots,k_n) \colon \cO(n) \otimes \cO(k_1) \otimes \cdots \otimes \cO(k_n) \lra \cO(k_1+\cdots+k_n),\]
which satisfy unit, associativity and equivariance axioms. The objects $\cO(n)$ are called \emph{$n$-ary operations of $\cO$}. If $\sfC$ is symmetric monoidal, then for $k=\infty$ this is a \emph{symmetric operad}\index{operad!symmetric} as defined Chapter 1 of \cite{GILS}, for $k=2$ it is a \emph{braided operad}\index{operad!braided} as in \cite{FiedorowiczBraided}, and for $k=1$ it is a \emph{non-symmetric operad}\index{operad!non-symmetric}.

\newglossaryentry{endx}{%
	name={\ensuremath{\cE_X}},
	description={Endomorphism operad of $X$},
	type=symbols
}
\begin{example}\label{example:endomorphism-operad} 
The prototypical example of an operad is the \emph{endomorphism operad}\index{operad!endomorphism}. Let $\sfC$ be a $(k+1)$-monoidal category as above, then given an object $X \in \sfC$ we may form the $k$-symmetric sequence
	\[\gls{endx}(n) \coloneqq \Hom_\sfC(X^{\otimes n},X),\]
with composition $\cE_X \circ \cE_X \to \cE_X$ induced by the internal composition, and unit $\bunit_\sfC \to \Hom_\sfC(X,X)$ adjoint to the identity morphism of $X$. (One should take $X^{\otimes 0} = \bunit_\sfC$ in the definition of $\mathcal{E}_X(0)$.)
\end{example}

As $\cO$ is a unital monoid, the functor $\cX \mapsto \cO \circ \cX$ has the structure of a monad on $\cat{FB}_k(\sfC)$, which we also denote $\cO$. Using the adjunction (\ref{eqn:symseq0adj}), we may transfer this monad from $\cat{FB}_k(\sfC)$ to $\sfC$, as described in Section \ref{sec:monads}. The resulting functor is $X \mapsto \cO(X)$ on $\sfC$, given in Definition \ref{def:symm-seq-application} and so by the formula \eqref{eqn:symmetric-pairing}, and we continue to denote this monad by $\cO$. Unravelling definitions, an algebra $\gX$ for the monad $\cO$ is an object $X \in \sfC$ together with morphisms\index{$\cO$-algebra}
\[a_n \colon \cO(n) \otimes X^{\otimes n} \lra X\]
for all $n \geq 0$ satisfying unit, associativity and equivariance axioms. If $\sfC$ is symmetric monoidal, for $k = \infty$ this recovers the classical definition over an algebra over an operad, as in Section 2 of \cite{GILS}, and similarly for $k=2$ \cite{FiedorowiczBraided} and $k=1$ \cite{HarperOperads}. Lemma \ref{lem:sseq-sifted} implies the following.

\begin{corollary}\label{cor:monad-from-operad-sifted} 
A monad associated to an operad is sifted.
\end{corollary}

\begin{remark}\label{rem:alg-alternative-defn} There are various equivalent points of view on an $\cO$-algebra structure. If $\sfC$ is $(k+1)$-monoidal,  using the closedness of the monoidal structure, the data of an $\cO$-algebra structure on $X$ is the same as a morphism $\alpha \colon \cO \to \cE_X$ of operads in $\sfC$. An $\cO$-algebra structure on an object $X \in \sfC$ is the same as an $\cO$-algebra structure on the $k$-symmetric sequence $0_*(X)$, since $0_*(X)(n)$ is initial for $n > 0$.\end{remark}

\begin{example}\label{exam:module-over-algebra-operad} Let us discuss the prototypical example in the case $k=0$. If $\gR$ is a unital monoid in a monoidal category $\sfC$ (also known as an \emph{associative algebra object}\index{associative algebra})\index{monoid} with underlying object $R \in \sfC$ then a left $\gR$-module $\gM$ is an object $M \in \sfC$ with a map $R \otimes M \to M$ satisfying unit and associativity axioms. 

Left $\gR$-modules may be encoded by an operad, denoted $\cR$. This operad has underlying $0$-symmetric sequence given by $\cR(1) = R$. The operad structure involves a map $\cR \circ \cR \to \cR$, which is the same as a map $R \otimes R \to R$ and is given by multiplication, and a map $1_*(\bunit_\sfC) \to \cR$, which is the same as a map $\bunit_\sfC \to R$ and is given by the unit. Thus the associated monad is given by $X \mapsto R \otimes X$. We will denote $\Alg_\gR(\sfC)$ by $\gR\text{-}\cat{Mod}$. There is a similar definition of right $\gR$-modules, which form a category $\cat{Mod}\text{-}\gR$.
This example has some special properties that are worth pointing out. Firstly, the monad $R \otimes -$ preserves all colimits, not just sifted ones, because the monoidal structure is closed. Secondly, the forgetful functor $U^\gR \colon \gR\text{-}\cat{Mod} \to \sfC$ has an enriched right adjoint given by $X \mapsto \Hom_{\sfC}(\gR,X)$. To give $\Hom_{\sfC}(R,X)$ the structure of a left $\gR$-module, we produce a map
\[R \otimes \Hom_\sfC(R,X) \lra \Hom_{\sfC}(R,X)\]
as the adjoint of the map 
\[R \otimes R \otimes \Hom_\sfC(R, X) \lra X\]
given by multiplication followed by evaluation. To show that this is the right adjoint, we note that $\Hom_{\gR\text{-}\cat{Mod}}(\gM,\gN)$ is given by the equalizer of two maps $\Hom_{\sfC}(M,N) \to \Hom_\sfC(R \otimes M,N)$. In the case $\gN = \Hom_{\sfC}(\gR,X)$, this simplifies to the equalizer of two maps $\Hom_{\sfC}(R \otimes M, X) \to \Hom_{\sfC}(R \otimes R \otimes M,X)$, which one may compute as $\Hom_\sfC(M,X)$. As a consequence of this $U^\gR$ commutes with all colimits.\end{example}

\subsection{Operads and lax monoidal functors}\label{sec:operads-and-lax-monoidal-functors} Suppose that $\sfC$ and $\sfD$ are $k$-monoidal categories. Recall that a lax $k$-monoidality on a functor $F \colon \sfC \to \sfD$ is given by a natural transformation
$F(X) \otimes_\sfD F(Y) \to F(X \otimes_\sfC Y)$ and a morphism $\bunit_\sfD \to F(\bunit_\sfC)$, subject to certain conditions. The conditions guarantee the existence of well-defined maps of iterated tensor powers
\begin{equation}\label{eq:25}
F(X)^{\otimes_\sfD n} \lra F(X^{\otimes_\sfC n})
\end{equation}
for all $n \geq 0$, and that these are compatible with the various composition maps. 

When $\sfC$ and $\sfD$ are $(k+1)$-monoidal the categories $\cat{FB}_k(\sfC)$ and $\cat{FB}_k(\sfD)$ are monoidal. If $F$ is $(k+1)$-monoidal, then applying $F$ levelwise gives rise to a lax monoidal functor
\[F \colon \cat{FB}_k(\sfC) \lra \cat{FB}_k(\sfD).\]
Such a functor sends an operad $\cO$ in $\sfC$ to an operad $F(\cO)$ in $\sfD$, and sends left modules, right modules, and algebras over $\cO$ in $\sfC$ to such objects over $F(\cO)$ in $\sfD$, with underlying objects obtained by applying $F$. Furthermore, the lax-monoidality on $F$ induces a natural transformation $\phi_F \colon F(\cO) F \Rightarrow F \cO$ which is part of a morphism of monads 
\[(F,\phi_F) \colon (\cat{D},F(\cO)) \to (\cat{C},\cO).\] 

Suppose that $L \colon \cat{C} \to \cat{D}$ is strong $k$-monoidal, and has a right adjoint $R$. By the above arguments, given an operad $\cO$ on $\cat{C}$ we obtain an operad $L(\cO)$ on $\cat{D}$ and a morphism of monads $(L,\phi_L) \colon (\cat{D},L(\cO)) \to (\cat{C},\cO)$. 

\begin{lemma}$(L,\phi_L)^* \colon \cat{Alg}_\cO(\cat{C}) \to \cat{Alg}_{L(\cO)}(\cat{D})$ satisfies $(L,\phi_L)^* F^\cO \cong F^{L(\cO)} L$ and preserves colimits of $\cO$-algebras.\end{lemma}

\begin{proof}Because $L$ is strong monoidal and preserves all colimits, the natural map
\[\begin{tikzcd}L(\cO)(L(X)) = \bigsqcup_{n \geq 0} L(\cO(n)) \otimes_{\sfD,G_n} L(X)^{\otimes_\sfD n} \dar \\[-5pt] L(\cO(X)) = L\Big(\bigsqcup_{n \geq 0} \cO(n) \otimes_{\sfC,G_n} X^{\otimes_\sfC n}\Big),\end{tikzcd}\]
upgrades to an isomorphism between the $L(\cO)$-algebras $F^{L(\cO)}(X)$ and $(L,\phi_L)^* F^\cO(X)$. This proves the first claim.

For the second claim, we first observe that $(L,\phi_L)^*$ preserves colimits of free diagrams, because it is preserves free algebras and $L$ preserves colimits. It similarly sends the canonical presentation of $\gX$ to the canonical presentation of $L(\gX)$, because we may verify this on underlying objects. Since these two types of colimits generate all colimits in $\Alg_\cO(\sfC)$, $(L,\phi_L)^*$ preserves all colimits.\end{proof}

On the other hand, the morphisms
\[\begin{tikzcd} L\Big(\bigsqcup_{n \geq 0} \cO(n) \otimes_{\sfC,G_n} R(A)^{\otimes_\sfC n}\Big) \cong \bigsqcup_{n \geq 0} L(\cO(n)) \otimes_\sfC L(R(A))^{\otimes_\sfD n} \dar \\[-5pt] \bigsqcup_{n \geq 0} L(\cO(n)) \otimes_{\sfC,G_n} A^{\otimes_\sfC n}\end{tikzcd}\]
give by adjunction a natural transformation $\psi \colon \cO R \Rightarrow RL(\cO)$, which is part of a morphism of monads $(R,\psi) \colon (\cat{C},\cO) \to (\cat{D},L(\cO))$. So we can also form $(L,\psi)_* \colon \Alg_\cO(\sfC) \to \Alg_{L(\cO)}(\sfD)$.

\begin{lemma}\label{lem:l-left-right-iso} The functors $(L,\psi)_*,(L,\phi_L)^* \colon \Alg_\cO(\sfC) \to \Alg_{L(\cO)}(\sfD)$ are naturally isomorphic.\end{lemma}

\begin{proof}The functor $(L,\psi)_*$ is determined uniquely up to natural isomorphism by the fact that it preserves sifted colimits and sends the free $\cO$-algebra on $X$ to the free $L(\cO)$-algebra on $L(X)$. But $(L,\phi_L)^*$ has the same properties.\end{proof}

\subsection{Non-unitary operads and unitalization} \label{sec:non-unitary} 
The following discussion is only relevant when $k>0$. Later we consider the following class of operads:

\begin{definition}An operad $\cO$ in $\sfC$ is said to be \emph{non-unitary}\index{operad!non-unitary} if $\cO(0) \cong \binit$, the initial object.\end{definition}

\begin{remark}Note this is still a \emph{unital} operad, i.e.\ there is a map $\bunit \to \cO(1)$ inducing a map $\iota \to \cO$.\end{remark}

\newglossaryentry{oplus}{%
	name={\ensuremath{\cO^+}},
	description={Unitalization of operad $\cO$},
	type=symbols
}
Given an operad $\cO^+$, we may form a non-unitary operad $\cO$ by letting $\cO(0) \coloneqq \binit$ and $\cO(n) \coloneqq \cO^+(n)$ for $n > 0$, with the same composition in positive arity. 

The morphism of operads $\upsilon \colon \cO \to \cO^+$ induces a morphism of associated monads, and hence an adjunction
\begin{equation}\label{eqn:unitalization-adjunction} \begin{tikzcd} \Alg_{\cO^+}(\sfC) \arrow[shift left=.5ex]{r}{\upsilon^*} & \Alg_{\cO}(\sfC), \arrow[shift left=.5ex]{l}{\upsilon_*}
\end{tikzcd}\end{equation}
where the underlying object of $\upsilon_*(\gR)$ is naturally isomorphic to $\cO^+(0) \sqcup \gR$. We call this the \emph{unitalization}\index{unitalization} of $\gR$ and denote it $\gR^+$. 

\newglossaryentry{ecan}{%
	name={\ensuremath{\epsilon_\mr{can}}},
	description={Canonical augmentation},
	type=symbols
}
Observe that restriction the operad composition $\cO^+ \circ \cO^+ \to \cO^+$ to arity $0$ yields a canonical $\cO^+$-algebra structure on $\cO^+(0)$. If $\sfC$ is pointed then there is a canonical map of $\cO^+$-algebras $\gls{ecan} \colon \gR^+ \to \cO^+(0)$, which we call the \emph{canonical augmentation}.  On free $\cO^+$-algebras this may be defined by the map $\cO^+(X) \cong \cO^+(0) \sqcup \cO(X) \to \cO(0)^+$ induced by the map $\cO(X) \to \ast$, and on general $\cO^+$-algebras it is defined by density under sifted colimits using Proposition \ref{prop:extend-by-density}. This is an example of the following notion:

\begin{definition}\label{defn:AlgAug}
An \emph{augmentation of an $\cO^+$-algebra}\index{augmentation} $\gR$ is an $\cO^+$-algebra map $\epsilon \colon \gR \to \cO^+(0)$. An \emph{augmented $\cO^+$-algebra}\index{$\cO^+$-algebra!augmented} is a pair $(\gR,\epsilon)$ of an $\cO^+$-algebra $\gR$ and an augmentation $\epsilon \colon \gR \to \cO^+(0)$.
\end{definition}

Given an augmented $\cO^+$-algebra $(\gR,\epsilon)$, we can form its augmentation ideal.

\newglossaryentry{augmideal}{%
	name={\ensuremath{I(\gR)}},
	description={Augmentation ideal},
	type=symbols
}
\begin{definition}If $\sfC$ is pointed, the \emph{augmentation ideal}\index{augmentation!ideal} $\gls{augmideal}$ of an augmented $\cO^+$-algebra $I(\gR)$ is given by $\ast \times_{\cO^+(0)} \gR$, the pullback along the augmentation.\end{definition}

\begin{lemma}If $\sfC$ is pointed and $\epsilon \colon \gR \to \cO(0)^+$ is an augmented $\cO^+$-algebra, then $I(\gR)$ has a canonical structure of an $\cO$-algebra such that the map $I(\gR) \to \gR$ is a map of $\cO$-algebras.
\end{lemma}

\begin{proof}The category $\Alg_{\cO}(\sfC)$ is complete and the forgetful functor $U^{\cO} \colon \Alg_{\cO}(\sfC) \to \sfC$ preserves limits as it is a right adjoint. As the diagram
\[* \lra \cO^+(0) \overset{\epsilon}\longleftarrow \gR\]
of which $I(\gR)$ is the limit consists of $\cO$-algebras and $\cO$-algebra maps, $I(\gR)$ inherits an $\cO$-algebra structure such that the map $I(\gR) \to \gR$ is one of $\cO$-algebras. \end{proof}

\newglossaryentry{unitalization}{%
	name={\ensuremath{(-)^+}},
	description={Unitalization of algebras over operad},
	type=symbols
}
Since the unitalization functor $\gls{unitalization} \colon \Alg_\cO(\sfC) \to \Alg_{\cO^+}(\sfC)$ is a left adjoint to the forgetful functor $\Alg_{\cO^+}(\sfC) \to \Alg_\cO(\sfC)$, from the map $I(\gR) \to \gR$ of $\cO$-algebras we obtain a canonical map $I(\gR)^+ \to \gR$ of $\cO^+$-algebras. This will not be an isomorphism in general.

\begin{definition}\label{def:split-augmented} If $\sfC$ is pointed, we say that an augmentation $\epsilon \colon \gR \to \cO^+(0)$ is \emph{split}\index{augmentation!split} if the canonical map $I(\gR)^+ \to \gR$ is an isomorphism.\end{definition}

\begin{remark}The adjunction (\ref{eqn:unitalization-adjunction}) induces an equivalence of categories between $\cO$-algebras and the subcategory of split augmented $\cO^+$-algebras: the functor $\upsilon_*$ is fully faithful and its essentially image are the split augmented $\cO^+$-algebras.\end{remark}

\subsection{Augmentations of operads} \label{sec:simplicial-operad-augmentation} 

Our goal is to study $T$-indecomposables, when $T$ is the monad associated to a non-unitary operad $\cO$ in $\sfC$. (Recall that we also write $\cO$ for the associated monad.) In order to do so, as described in Section \ref{sec:AugMonad} this monad must be equipped with an augmentation $\epsilon \colon \cO \to +$.

In an operad a distinguished role is played by the unary operations $\cO(1)$, which the operad structure makes into a unital monoid. (On the other hand, any unitary monoid may be considered as an operad with only unary operations, as explained in Example \ref{exam:module-over-algebra-operad}.) We will first explain how a non-unitary operad has a canonical ``relative" augmentation, and then explain what is needed to promote this to an augmentation.

\subsubsection{The canonical relative augmentation and relative indecomposables} \label{sec:relative-augmentation}
\newglossaryentry{urel}{%
	name={\ensuremath{U^\cO_{\cO(1)}}},
	description={Relative forgetful functor on $\cO$-algebras},
	type=symbols
}
\newglossaryentry{freerel}{%
	name={\ensuremath{F^\cO_{\cO(1)}}},
	description={Relative free $\cO$-algebra functor},
	type=symbols
}
\newglossaryentry{qrel}{%
	name={\ensuremath{Q^\cO_{\cO(1)}}},
	description={Relative $\cO$-indecomposables functor},
	type=symbols
}
\newglossaryentry{zrel}{%
	name={\ensuremath{Z^\cO_{\cO(1)}}},
	description={Relative trivial $\cO$-algebra functor},
	type=symbols
}
There is a map of operads $\cO(1) \to \cO$, which induces a factorization of the monadic adjunction for $\cO$:
\[\begin{tikzcd}
 & \Alg_{\cO}(\sfC) \arrow[shift left=.5ex]{ld}{U^\cO_{\cO(1)}} \arrow[shift left=-.5ex,swap]{rd}{U^\cO} &[20pt] \\[20pt]
\Alg_{\cO(1)}(\sfC) \arrow[shift left=.5ex]{ru}{F^\cO_{\cO(1)}} \arrow[shift left=-.5ex,swap]{rr}{U^{\cO(1)}} & & \sfC, \arrow[shift left=-.5ex,swap]{ll}{F^{\cO(1)}} \arrow[shift left=-.5ex,swap]{lu}{F^\cO}
\end{tikzcd}\]
where $\gls{urel}$ is the relative forgetful functor and the relative free algebra functor $\gls{freerel}$ has underlying object given by
\[F^\cO_{\cO(1)}(X) \coloneqq \bigsqcup_{n \geq 0} \cO(n) \otimes_{G_n \wr \cO(1)} X^{\otimes n}.\]
Here $G_n \wr \cO(1)$ is the monoid, or equivalently associative algebra, with underlying object $G_n \times \cO(1)^{\otimes n}$ and composition induced the composition in $G_n$ and $\cO(1)$, and the action of $G_n$ on $\cO(1)^{\otimes n}$ by the braiding. The relative tensor product is that given in the proof of Proposition \ref{prop:ModuleCatAxioms}.

To define relative indecomposables, we need to replace the monad $+$ by one denoted $\cO(1)_+$, which we shall construct now. The monoid $\cO(1)$ yields a monad $\cO(1)_+$ on $\sfC$ whose value on an object $X$ is $(\cO(1) \otimes X)_+$ and whose multiplication is given by the map
\[(\cO(1) \otimes (\cO(1) \otimes X)_+)_+ \lra (\cO(1) \otimes \cO(1) \otimes X)_+ \lra (\cO(1) \otimes X)_+\]
which first uses that $\cO(1) \otimes -$ preserves coproducts and collapses the term $\cO(1) \otimes \bterm$ to $\bterm$, and then uses the monoid structure of $\cO(1)$. Equivalently, this is the monad obtained as the composition $+\cO(1)$ of the monad induced by $\cO(1)$ and the monad $+$, using the distributive law $\cO(1) + \Rightarrow + \cO(1)$ induced by the above collapse map (as in \cite{BeckDistributive}). We can interpret this in terms of the category $\sfC_\ast$: the equivalence $\Alg_\mr{Id}(\sfC_\ast) \cong \Alg_+(\sfC)$ generalises to an equivalence $\Alg_{\cO(1)_+}(\sfC) \cong \Alg_{\cO(1)}(\sfC_*)$ where the monad $\cO(1)$ on $\sfC_*$ has value on an object $X$ given by $\cO(1)_+ \owedge X$.

Let us now assume that the operad $\cO$ is non-unitary. In this case there is morphism of monads $\epsilon^\cO_{\cO(1)} \colon \cO \to \cO(1)_+$ given by sending the $k$-ary operations for $k \geq 2$ to the basepoint. We define the functor $\gls{qrel}$ of \emph{relative indecomposables}\index{indecomposables!relative} to be
\[(\epsilon^\cO_{\cO(1)})_* \colon \Alg_{\cO}(\sfC) \lra \Alg_{\cO(1)}(\sfC_\ast).\] That is, $Q^\cO_{\cO(1)}$ is the left adjoint in change-of-monad for the map of monads $\epsilon^\cO_{\cO(1)}$, which has a right adjoint $\gls{zrel} \colon \Alg_{\cO(1)}(\sfC_\ast) \to \Alg_{\cO}(\sfC)$ called the \emph{relative trivial algebra functor}. The trivial $\cO$-algebra on an $\cO(1)$-algebra has the same underlying object (where we forget the special status of the basepoint) and all 1-ary operations act by $\cO(1)$ and $k$-ary operations for $k \geq 2$ map to the basepoint.

\subsubsection{Relative decomposables}\label{sec:relative-decomposables} It is occasionally useful to express the relative indecomposables $Q^\cO_{\cO(1)}$ in terms of relative decomposables.

\newglossaryentry{dec}{%
	name={\ensuremath{\dec^\cO_{\cO(1)}}},
	description={Relative $\cO$-decomposables functor},
	type=symbols
}
\begin{definition}
	\label{def:decomposables} The \emph{relative decomposables}\index{decomposables!relative} of a free $\cO$-algebra $F^\cO(X)$ are 
	\[\gls{dec}(F^\cO(X)) \coloneqq \left(\bigsqcup_{n \geq 2} \cO(n) \otimes_{G_n} X^{\otimes n}\right)_+ \in \Alg_{\cO(1)}(\sfC_\ast).\]
	The right-hand side defines a right $\cO$-functor which commutes with sifted colimits, so there is a unique extension to $\dec^\cO_{\cO(1)} \colon \Alg_{\cO}(\sfC) \to \Alg_{\cO(1)}(\sfC_\ast)$ by Proposition \ref{prop:extend-by-density}.
\end{definition}

Intuitively, this takes quotient by all $n$-ary operations with $n \geq 2$, and remembers the action by the monoid $\cO(1)$.  Let us spell out these two functors more explicitly. The natural transformation given by the inclusion
\[\left(\bigsqcup_{n \geq 2} \cO(n) \otimes_{G_n} X^{\otimes n}\right)_+ \lra \left(\bigsqcup_{n \geq 1} \cO(n) \otimes_{G_n} X^{\otimes n}\right)_+\]
induces a natural transformation $\dec^\cO_{\cO(1)} \Rightarrow +(U^{\cO}_{\cO(1)})$. If $\gR \in \Alg_\cO(\sfC)$ we may form the following pushout in $\Alg_{\cO(1)}(\sfC_\ast)$
\[\begin{tikzcd}\dec^\cO_{\cO(1)}(\gR) \rar \dar & U^{\cO}_{\cO(1)}(\gR)_+ \dar \\
\ast \rar &  U^{\cO}_{\cO(1)}(\gR)_+/\dec^\cO_{\cO(1)}(\gR).\end{tikzcd}\]
On free $\cO$-algebras we have 
\[U^{\cO}_{\cO(1)}(F^\cO(X))_+/\dec^\cO_{\cO(1)}(F^\cO(X)) \cong (\cO(1) \times X)_+ \cong Q^\cO_{\cO(1)}(F^\cO(X)),\] 
so by Proposition \ref{prop:extend-by-density} we have:

\begin{lemma}\label{lem:indec-dec-quotient}
For $\gR \in \Alg_\cO(\sfC)$, there is a natural isomorphism 
	\[Q^\cO_{\cO(1)}(\gR) \cong U^{\cO}_{\cO(1)}(\gR)_+/\dec^\cO_{\cO(1)}(\gR) \in \Alg_{\cO(1)}(\sfC_\ast).\]
\end{lemma}

\subsubsection{Absolute augmentations and absolute indecomposables}

Suppose the monoid $\cO(1)$, which is an algebra for the unital associative operad, is equipped with an augmentation $\epsilon \colon \cO(1) \to \bunit$ in the sense of Definition \ref{defn:AlgAug}, i.e.\ a morphism of unital monoids. Then it defines an augmentation $\epsilon_+ \colon \cO(1)_+ \to \bunit_+ = +$ of the monad $\cO(1)_+$ and so composing it with the canonical relative augmentation $\epsilon^\cO_{\cO(1)} \colon \cO \to \cO(1)_+$ it defines an augmentation of the operad $\cO$.

Using this we may define the \emph{absolute $\cO$-indecomposables}\index{indecomposables!absolute}
\[Q^\cO \colon \Alg_\cO(\sfC) \lra \Alg_+(\sfC) = \sfC_*,\]
which by Lemma \ref{lem:ind-change-of-monads} satisfies $Q^\cO(\gX) \cong Q^{\cO(1)_+} (Q^\cO_{\cO(1)}(\gX))$. Using the Lemma \ref{lem:indec-dec-quotient}, we can express the absolute $\cO$-indecomposables as
\[Q^\cO(\gX) \cong Q^{\cO(1)_+} (U^{\cO}_{\cO(1)}(\gR)_+/\dec^\cO_{\cO(1)}(\gR)) \in \sfC_*.\]
As an $\cO(1)_+$-algebra is a pointed object with an $\cO(1)$-action, it is natural to think of the $\cO(1)_+$-indecomposables as taking orbits for the $\cO(1)$-action. (Although it is not reflected in the notation, this of course depends on the choice of augmentation $\epsilon \colon \cO(1) \to \bunit$.)

\subsubsection{Relative indecomposables and lax monoidal functors} Let $L \colon \sfC \to \sfD$ be a strong monoidal functor, with right adjoint $R$. Then $L(\cO)$ is non-unitary if $\cO$ is, because $L$ preserves initial objects. In particular, both $\cO$ and $L(\cO)$ have a canonical relative augmentations. Moreover, the counit $LRX \to X$ induces a map
\[L(\cO(1) \otimes RX) \cong L(\cO(1)) \otimes LRX \lra L(\cO(1)) \otimes X,\]
which upon adjunction and adding a disjoint basepoint yields a map $\cO(1)_+(R(X)) \to R(L(\cO(1))_+(X)$. These assemble to a morphism of monads $(R,\upsilon_{\cO(1)}) \colon (\cat{C},\cO(1)_+) \to (\cat{D},L(\cO(1))_+)$ which fits into a commutative diagram
\[\begin{tikzcd} (\cat{C},\cO) \dar[swap]{(R,\psi)} \rar{\epsilon^\cO_{\cO(1)}} &[10pt] (\cat{C},\cO(1)_+) \dar{(R,\upsilon_{\cO(1)})} \\
(\cat{D},L(\cO)) \rar{\epsilon^{L(\cO)}_{L(\cO)(1)}} & (\cat{D},L(\cO(1))_+) \end{tikzcd} \]
As a consequence of Lemmas \ref{lem:ind-change-of-monads} and \ref{lem:l-left-right-iso}, this implies that $L$ preserves relative indecomposables, as follows:

\begin{lemma}There is a natural isomorphism $(L,\phi_L)_+ Q^\cO_{\cO(1)} \cong Q^{L(\cO)}_{L(\cO(1))} (L,\phi_L)_*$ of functors $\Alg_\cO(\sfC) \to \Alg_{L(\cO(1))_+}(\sfD) = \Alg_{L(\cO(1))}(\sfD_\ast)$.\end{lemma}

\begin{remark}Since $L$ preserves free algebras and colimits, it also preserves relative decomposables and the formula of Lemma \ref{lem:indec-dec-quotient}.\end{remark}

Suppose now that the unital monoid $\cO(1)$ comes with an augmentation $\epsilon \colon \cO(1) \to \bunit$. Since $L$ is strong monoidal, we obtain from this an augmentation $L\epsilon \colon L(\cO(1)) \to L(\bunit) \cong \bunit$. As the above construction only used that $\cO(1)$ is a monoid, we can apply it equally well to $\bunit$ to get a morphism of monads $(R,\upsilon_\bunit) \colon (\sfC,+) \lra (\sfD,+)$, which fits into a commutative diagram
\[\begin{tikzcd} (\sfC,\cO(1)_+) \dar[swap]{(R,\upsilon_{\cO(1)})} \rar{\epsilon} & (\sfC,+) \dar{(R,\upsilon_\bunit)} \\
(\sfD,L(\cO(1))_+) \rar{L\epsilon} & (\sfD,+)\end{tikzcd}\]
This implies that $L$ preserves absolute indecomposables, in the following precise sense: 

\begin{lemma}\label{lem:l-pres-abs-indec}There is a natural isomorphism $L_+ \circ Q^\cO \cong Q^{L(\cO)} \circ (L,\phi_L)^*$ of functors $\Alg_\cO(\sfC) \to \sfD_\ast$.\end{lemma}

\subsection{Operads in simplicial sets} \label{sec:simplicial-operads} 

By the discussion in Section \ref{sec:operads-and-lax-monoidal-functors}, every operad in simplicial sets gives rise to an operad in $\sfC$ by applying $s \colon \cat{sSet} \to \sfC$ objectwise. Operads in $\sfC$ which arise in this way enjoy certain special properties, and some of our results will only hold for such operads. We shall generally write $\cC$ for an operad in simplicial sets (just as we write $\cO$ for an operad in $\sfC$), and continue to write $\cC$ for the operad $s(\cC)$ in $\sfC$.

Firstly, for $k=2$ the spaces of operations of such operads are always \emph{central}. This avoids any difficulties with the composition product of 2-symmetric sequences in a 2-monoidal category. 

Secondly, the monad associated to a simplicial operad is a simplicial monad, i.e.~it is a monoid in the category of simplicially enriched functors. 

Thirdly, they are compatible with change-of-diagram-category. If $p \colon \sfG \to \sfG'$ is strong $k$-monoidal, then by Lemma \ref{lem:oplax-lax-Day-convo} the functor $p_* \colon \sfS^\sfG \to \sfS^{\sfG'}$ is also strong $k$-monoidal. Furthermore, as $p(\bunit_\sfG) \cong \bunit_{\sfG'}$ we have $p_* \circ (\bunit_\sfG)_* \cong (\bunit_{\sfG'})_*$. Thus from the formula above we see that $p_*(\cC(X)) \cong \cC(p_*(X)),$ and hence there is an induced functor $p_* \colon \Alg_\cC(\sfS^\sfG) \to \Alg_\cC(\sfS^{\sfG'})$, which by construction satisfies $p_*(F^\cC(X)) \cong F^\cC(p_*(X))$. In fact, this also works when using an operad $\cO$ in $\sfS$ in place of $\cC$.

Fourthly, as the category $\cat{sSet}$ is cartesian, i.e.\ the monoidal product coincides with the categorical product, any monoid in simplicial sets has a unique augmentation: if $\cC$ is an operad in simplicial sets then $\epsilon \colon \cC(1) \to *$ is an augmentation, and so if it is non-unitary the operad $s(\cC)$ is augmented as described in the previous section. We call this the \emph{canonical augmentation}\index{augmentation!canonical}. There is therefore defined an absolute indecomposables functor $Q^\cC \colon \Alg_\cC(\sfC) \to \sfC_*$.

\begin{example}\label{ex:UnitIsCAlg}
For any operad $\cC$ in simplicial sets, the monoidal unit $\bunit_\sfC$ is canonically a $\cC$-algebra, via the map 
	\[\cC(\bunit_\sfC) \cong \left(\bigsqcup_{n \geq 0} \cC(n)/G_n \right) \times \bunit_\sfC \lra \ast \times \bunit_\sfC \cong \bunit_\sfC,\]
induced by the unique map of simplicial sets $\bigsqcup_{n \geq 0} \cC(n)/G_n \to *$. The isomorphism uses that the $G_n$-action on $\bunit_\sfC \cong \bunit_\sfC^{\otimes n}$ is trivial, which follows from the axioms of a braided monoidal category.
\end{example}

The last special property of operads in simplicial sets is that their algebras in $\sfC$ have a monoidal structure, coming from the diagonal maps in the category of simplicial sets which make every simplicial set into a cocommutative coalgebra.

\begin{proposition}\label{prop:algcc-symmetric-monoidal}
Let $k \in \{2, \infty\}$, $\cC \in \cat{FB}_k(\cat{sSet})$ be an operad, and $\sfC$ be $\infty$-monoidal. Then $\cat{Alg}_{\cC}(\cat{C})$ has a $(k-1)$-monoidal structure $\otimes_\cC$ such that the functor $U^\cC \colon \cat{Alg}_{\cC}(\cat{C}) \to \cat{C}$ is strong monoidal.\index{tensor product!of algebras over operad}
\end{proposition}

\begin{proof}
The key point is that, as $\cat{sSet}$ is cartesian, the $k$-symmetric sequence $\cC$ has the structure of an unital monoid in $\cat{FB}_k(\cat{Coalg}(\cat{sSet}))$ and hence, on applying $s$, in $\cat{FB}_k(\cat{Coalg}(\cat{C}))$, making the associated monad a \emph{Hopf monad} as described in \cite{MoerdijkHopf}. This is described in Example 3.1 of that paper if $k=\infty$, but goes through for $k \geq 2$. Furthermore, if $k =\infty$ it is actually a cocommutative Hopf monad. The result is then the combination of Propositions 1.4 and 3.2 of \cite{MoerdijkHopf}.
\end{proof}
\section{Filtered algebras}
\label{sec:filtered-objects-1}

For computational as well as for conceptual reasons we shall consider filtered objects in $\sfC$, and well as filtered $\cO$-algebras. For the latter we never mean filtered objects in the category of $\cO$-algebras in $\sfC$, but rather $\cO$-algebras in the category of filtered objects in $\sfC$. In this section we define these notions and present various constructions related to them. We continue to let $\sfC$ be a category satisfying the axioms of Section \ref{sec:axioms-of-cats}.

\subsection{Graded and filtered objects} The perspective we take on filtered objects is close to that of Gwilliam--Pavlov \cite{GwilliamPavlov}; this will continue when we later discuss the homotopy theory of filtered objects.
\newglossaryentry{zequal}{%
	name={\ensuremath{\bZ_=}},
	description={Integers as category with only identity morphisms.},
	type=symbols
}
\begin{definition}\mbox{}\label{def:filtered-graded-objects}
 \begin{itemize}
	\item Let $\gls{zequal}$ denote the category with objects $\Z$, and only identity morphisms.  A \emph{graded object}\index{graded object} in $\sfC$ is a functor $X \colon \Z_= \to \sfC$. The category of graded objects shall be denoted $\sfC^{\Z_=}$.
	\item Let $\Z_\leq$ denote the category associated to the partially ordered set $(\Z,\leq)$.  A \emph{filtered object}\index{filtered object}\index{filtration} in $\sfC$ is a functor $X \colon \Z_\leq \to \sfC$.  The category of filtered objects shall be denoted $\sfC^{\Z_\leq}$.
	\item A filtered object $X$ is called \emph{ascending}\index{filtration!ascending} if the morphism $X(i) \to X(i+1)$ is an isomorphism for $i < -1$. It is called \emph{descending}\index{filtration!descending} if the morphism $X(i) \to X(i+1)$ is an isomorphism for all $i \geq 0$.
	\end{itemize}
Both $\Z_=$ and $\Z_\leq$ are symmetric monoidal categories using addition of integers, and so Day convolution endows $\sfC^{\bZ_=}$ and $\sfC^{\bZ_\leq}$ with (symmetric) monoidal structures.
\end{definition}

For a graded object $X \in \sfC^{\bZ_=}$, we shall think of $\bigsqcup_{i \in \bZ} X(i) \in \sfC$ as the ``underlying'' object, and think of $X(i)$ as ``grading $i$.'' Similarly, for a filtered object $X \in \sfC^{\Z_\leq}$ we shall think of $\mr{colim}(X) \coloneqq \colim_{i \in \bZ_\leq} X(i) \in \sfC$ as the ``underlying'' unfiltered object of $X$, and think of $X(i)$ as ``filtration $i$.'' However, we emphasize that we do not require the maps $X(i) \to X(i+1)$ to be injective in any sense. For example, the terminal object $\bterm$ has many interesting filtrations in our sense.

\subsection{Functors between filtered and unfiltered objects}
\label{sec:from-filt-unfilt}
\newglossaryentry{grr}{%
	name={\ensuremath{\grr}},
description={Associated graded},
type=symbols
}
There are adjunctions
\[\mr{colim} \dashv \mr{const} \dashv \mr{lim},\]
which we shall describe momentarily. Furthermore, for each integer $a \in \bZ$, we shall describe adjunctions
\begin{equation*}
a_! \dashv a_* \dashv a^* \dashv a^!,
\end{equation*}
where $a^*X = X(a)$ and the leftmost adjunction need only exist when $\sfC$ is pointed. Finally, we shall describe  an ``associated graded" functor $\gls{grr}$ and an adjunction
\[\grr \dashv u.\]

\subsubsection{The functors $\mr{colim}$ and $\mr{const}$} \label{sec:filtered-colim} Because the functor $\mr{colim} \colon \sfC^{\bZ_{\leq}} \to \sfC$ may be described as forming left Kan extension $\pi_*$ along the functor $\pi \colon \Z_\leq \to \{\ast\}$, it has a right adjoint $\pi^*$ given by precomposing with that functor,
i.e.\ sending an object of $\sfC$ to the corresponding constant functor, which we shall denote $\mr{const} \colon \sfC \to \sfC^{\bZ_{\leq}}$. The functor $\mr{const}$ has a further right adjoint given by right Kan extension, sending $X \in \sfC^{\Z_\leq}$ to the limit $\mr{lim}(X) \coloneqq \lim_{i \in \bZ_{\leq}} X(i) \in \sfC$.

\subsubsection{The functors $a_!$, $a_*$, $a^*$ and $a^!$} \label{sec:filtered-integers}
The basic functor associated to an object $a \in \Z$ is given by evaluation of $X$ at $a$, giving a functor
\begin{align*}
  a^* \colon \sfC^{\Z_\leq} & \lra \sfC\\
  X & \longmapsto X(a).
\end{align*}
This may be identified with precomposition with the functor $a \colon \{\ast\} \hookrightarrow \Z_\leq$ given by $\ast \mapsto a$.

Since $\sfC$ is cocomplete and complete, $a^*$ has both left and right adjoints. The left adjoint is given by left Kan extension along $a \colon \{\ast\} \to \Z_\leq$:
  \begin{align*}a_* \colon \sfC &\lra \sfC^{\Z_\leq} \\
  Y &\longmapsto a_*Y\coloneqq \left(n \mapsto \begin{cases}
    \binit & \text{if $n<a$} \\
    Y & \text{if $n \geq a$}
  \end{cases}\right),\end{align*}
where the functoriality sends a morphism $n \leq m$ in $\bZ_\leq$ to the identity map of $Y$ if $a \leq n$. Similarly, the right adjoint is given by right Kan extension:
\begin{align*}
a^! \colon \sfC &\lra \sfC^{\Z_\leq} \\
  Y &\longmapsto a^!Y \coloneqq \left(n \mapsto
  \begin{cases}
    \bterm & \text{if $n>a$} \\
    Y & \text{if $n \leq a$}
  \end{cases}\right),\end{align*}
where the functoriality sends a morphism $n \leq m$ to the identity map of $Y$ if $a \geq m$.

If $\sfC$ is pointed then $a_* \colon \sfC \to \sfC^{\Z_\leq}$ admits a further left adjoint
\begin{align*}a_! \colon \sfC^{\Z_\leq} &\lra \sfC \\
X &\mapsto \mr{colim}\left({\begin{tikzcd}[ampersand replacement=\&] X(a-1) \rar \dar \& \mr{colim}(X) \\
\ast \& \end{tikzcd}}\right).\end{align*}
If $\sfC$ is not pointed we may still use this formula to define a functor $a_! \colon \sfC^{\Z_\leq} \to \sfC_*$, replacing $\ast$ by $\bterm$.

\begin{remark}If $\sfC$ is pointed then $a^!$ admits a further right adjoint, sending $X$ to the pullback of $\ast \to X(a+1) \leftarrow \lim_i X(i)$. It does not play a role in this paper, but its existence does imply that $a^!$ preserves colimits.\end{remark}

\subsubsection{The functors $\grr$ and $u$.}
\label{sec:associated-graded-2} The \emph{associated graded}\index{filtration!associated graded}\index{associated graded} functor is given by
\begin{align*}\grr \colon \sfC^{\Z_\leq} &\lra \sfC^{\Z_=}_* \\
X &\longmapsto \grr(X) \coloneqq \left(n \mapsto \mr{colim}\left({\begin{tikzcd}[ampersand replacement=\&] X(n-1) \rar \dar \& X(n) \\
	\bterm \& \end{tikzcd}}\right)\right),\end{align*}
where the pushout is regarded as a pointed object using the induced morphism from $\bterm$. We shall occasionally write this as $\grr(X)(n) = X(n)/X(n-1)$, but we emphasize again that $X(n-1) \to X(n)$ need not be injective in any sense. 

Recall that $\sfC_*$ is the category of pointed objects, and $U^+ \colon \sfC_\ast \to \sfC$ is the functor which forgets basepoint. This is part of an adjunction
\[\begin{tikzcd} \sfC^{\bZ_\leq} \arrow[shift left=.5ex]{r}{\grr} & \sfC_\ast^{\bZ_=} \arrow[shift left=.5ex]{l}{u}\end{tikzcd}\]
with right adjoint given by
\begin{align*}u \colon \sfC^{\Z_=}_* & \lra \sfC^{\Z_\leq} \\
X & \longmapsto \left(n \mapsto U^+ X(n)\right), \end{align*}
and extending the functoriality of the composition $U^+ \circ X \colon \Z_= \to \sfC$ by sending non-identity morphisms $a < b$ to the unique maps $X(a) \to X(b)$ factoring through the basepoint $\bterm \to X(b)$.

\subsection{Monoidality}
\label{sec:filtr-mult-struct}

Let us discuss the extent to which the various functors discussed in Section \ref{sec:from-filt-unfilt} preserve the monoidal structures on $\sfC^{\Z_=}$ and $\sfC^{\Z_\leq}$ given by Day convolution. Bemusingly, several of them admit only ``half'' the structure of a monoidal functor: they preserve the product in the lax sense, but not the monoidal unit.

\subsubsection{The functors $\mr{colim}$ and $\mr{const}$}
\label{sec:colimit-2}

Since $\pi \colon \Z_\leq \to \ast$ is strong monoidal, the left Kan extension $\colim \colon \sfC^{\Z_\leq} \to \sfC$ is also strong monoidal by Lemma \ref{lem:oplax-lax-Day-convo}. Hence its right adjoint $\mr{const}$ is lax monoidal.

\subsubsection{The functors $a_!$, $a^*$, $a_*$ and $a^!$} Because the monoidal structure on $\sfC^{\Z_\leq}$ is defined by Day
convolution, there are canonical morphisms
\begin{equation*}
  X(a) \otimes_\sfC Y(b) \lra (X \otimes_{\sfC^{\Z_\leq}} Y)(a + b),
\end{equation*}
which may be interpreted as a natural transformation $(a^* X) \otimes
(b^* Y) \to (a+b)^*(X \otimes Y)$. For $a = b \leq 0$ we have $a + b \leq a$, so we may compose to get morphisms
\begin{equation}\label{eq:24}
  (a^*X) \otimes_{\sfC^{\Z_\leq}} (a^*Y) \lra a^*(X \otimes_{\sfC^{\Z_\leq}} Y),
\end{equation}
compatible with associators. These are $k$-symmetric if $\sfC$ is $k$-monoidal.

\begin{remark}We warn reader that (\ref{eq:24}) does \emph{not} promote $a^*$ to a monoidal functor unless $a = 0$, because there might not be a morphism $\bunit_\sfC \to a^*(\bunit_{\sfC^{\Z_\leq}})$ with the required properties. For example, for $a < 0$ the functor $a^*$ will send a non-unital monoid in $\sfC^{\Z_\leq}$ to a non-unital monoid in $\sfC$, but need not send a unital monoid to a unital monoid.\end{remark}

Next we turn to the left adjoint $a_* \colon \sfC \to \sfC^{\Z_\leq}$ of $a^*$.  Here
we have a morphism
\begin{equation}\label{eq:0983}
  (a_*X) \otimes_{\sfC^{\bZ_\leq}} (b_* Y) \lra (a+b)_*(X \otimes_\sfC Y),
\end{equation}
coming from the canonical morphisms for $i \geq a$ and $j \geq b$
\[X(i) \otimes_\sfC Y(j) \lra (X \otimes_\sfC Y)(i + j) = ((a+b)_*(X
\otimes_\sfC Y))(i + j).\]

\begin{lemma}\label{lem:alower-tensor} The morphism \eqref{eq:0983} is an isomorphism. 
Furthermore, if $K \in \cat{sSet}$, then $K \times (a_* X) \cong a_*(K \times X)$.
\end{lemma}

\begin{proof}Recall that $(a_* X)(r)$ is $\binit$ if $r<a$ and $X$ otherwise. The tensor product $(a_* Y) \otimes_{\sfC^{\bZ_\leq}} (b_* Z)$ is given on $r \in \bZ$ by the colimit
	\[\underset{r_1+r_2 \leq r}{\colim} \left(a_*Y(r_1) \otimes_\sfC b_*Z(r_2)\right).\]
A term in this diagram is initial if $r_1<a$ or $r_2<b$, and is $Y \otimes_\sfC Z$ otherwise. The maps in the diagram are either the canonical map from the initial object or isomorphisms. Thus the value is initial when $r < a+b$, and otherwise is equivalent to the colimit over a constant diagram on $Y \otimes_\sfC Z$ having initial object $(a,b)$. 

The proof of the second claim is similar.
\end{proof}

For $a \geq 0$ the morphism $a \leq 2a$ gives a natural transformation $(2a)_* \Rightarrow a_*$ and hence we get a natural map
\[(a_*X) \otimes_{\sfC^{\leq}} (a_* Y) \cong (2a)_*(X \otimes_\sfC Y) \lra a_*(X \otimes_\sfC Y).\]
Due to issues with the unit, this is again not a monoidal functor unless $a = 0$, but it will preserve non-unital multiplicative structures. More precisely, the functor $a \colon \{\ast\} \to \Z_\leq$ satisfies the part of being oplax concerning the tensor product, but not the part concerning the unit.

Finally, we discuss the functor $a_! \colon \sfC^{\Z_\leq} \to \sfC_*$, when $\sfC_*$ is given the monoidal structure $\otimes_{\sfC_*} = \owedge$ constructed in Section \ref{sec:basepoint-monad}. In general there is no good map between $(a_! X) \otimes_{\sfC_{\ast}} (a_! Y)$ and $a_!(X \otimes_{\sfC^{\bZ_{\leq}}} Y)$. For descending filtered objects $X \in \sfC^{\Z_\leq}$ we have a natural isomorphism $a_!X \cong X(0)/X(a-1)$. For $a \leq 0$ we have natural maps
\begin{equation*}
  (X(0)/X(a-1)) \otimes_{\sfC_*} (Y(0)/Y(a-1)) \lra (X \otimes Y)(0)/(X \otimes
  Y)(a-1),
\end{equation*}
which may be promoted to the structure of a monoidal functor on 
\[X \mapsto X(0)/X(a-1) \colon \sfC^{\Z_\leq} \lra \sfC_*,\]
and thus $a_!$ becomes a monoidal functor when restricted to descending objects.

\subsubsection{The functors $\grr$ and $u$}\label{sec:MonoidalityOnGr} We
construct a strong monoidality on the functor $\grr \colon \sfC^{\Z_\leq} \to \sfC^{\Z_=}_*$,  i.e.\ natural isomorphisms
\begin{align*}
  \grr(X) \otimes_{\sfC^{\Z_=}_*} \grr(Y) &\lra \grr(X
  \otimes_{\sfC^{\Z_\leq}} Y) \\
  \bunit_{\sfC_{\ast}^{\bZ_=}} &\lra \grr(\bunit_{\sfC^{\bZ_\leq}}),
\end{align*}
satisfying the usual axioms of monoidal functors. In Section \ref{sec:basepoint-monad} we have described the monoidal structure $\otimes_{\sfC_*} = \owedge$ on $\sfC_*$ and given a strong monoidality on $F^+ \colon \sfC \to \sfC_*$, and this induces a strong monoidality on $F^+ \colon \sfC^{\Z_\leq} \to \sfC^{\Z_\leq}_*$ when both are equipped with monoidal structures by Day convolution. This reduces the question to the case where $\sfC$ is pointed, i.e.\ we want a strong monoidality on $\grr \colon \sfC^{\Z_\leq}_* \to \sfC^{\Z_=}_*$.

Let us first describe a lax monoidality on the right adjoint $u \colon \sfC^{\Z_=}_* \to \sfC^{\Z_\leq}_*$ to $\grr$.  By definition, $((uX) \otimes_{\sfC_\ast^{\bZ_\leq}} (uY))(n)$ is the colimit of $(uX(a)) \owedge (uY(b))$ over the poset consisting of $(a,b) \in \Z_\leq \times \Z_{\leq}$ with $a + b \leq n$.  That diagram sends any non-identity morphisms to the trivial morphisms, from which it is easily deduced that
\begin{equation*}
  \bigvee_{a + b = n} X(a) \owedge Y(b) \lra \left((uX) \otimes_{\sfC^{\Z_\leq}_*}
  (uY)\right)(n)
\end{equation*}
is an isomorphism.  This amounts to a natural isomorphism
\[(uX) \otimes_{\sfC_\ast^{\bZ_\leq}}  (uY) \to u(X \otimes_{\sfC_\ast^{\bZ_=}}  Y),\]
but nevertheless we only get a lax monoidality because the obvious map from the monoidal unit of
$\sfC^{\Z_\leq}_*$ to the object $S_0 \coloneqq \smash{u(\bunit_{\sfC_\ast^{\bZ_=}})}$ is rarely an isomorphism. Indeed, we have that $S_0$ and the unit of $\sfC_*^{\bZ_\leq}$ are given by
\[S_0(n) = \begin{cases} \bunit_{\sfC_*} & \text{if $n=0$,} \\
\ast & \text{otherwise,} \end{cases} \qquad \text{and} \qquad \bunit_{\sfC_\ast^{\bZ_\leq}}(n) = \begin{cases} \bunit_{\sfC_\ast} & \text{if $n \geq 0$,} \\
\ast & \text{otherwise,}\end{cases} \]
where we recall that $\ast$ denotes the initial and terminal object in a pointed category.

The object $S_0$ is canonically a unital monoid in the
monoidal category $\sfC^{\Z_\leq}_*$. Since $\bunit \to S_0$
is an epimorphism, being a module over $S_0$ is simply the \emph{property} of being in the essential image of the fully faithful functor $u$, i.e.\ that all maps $X(n-1) \to X(n)$ factor through the terminal object.  In fact, the (lax monoidal) functor $u$ gives an equivalence of categories from $\sfC^{\Z_=}_*$ to the full subcategory
of $\sfC^{\Z_\leq}_*$ consisting of $S_0$-modules.  Under this
equivalence, the functor $\grr$ becomes identified with
$X \mapsto S_0 \otimes_{\sfC_\ast^{\bZ_\leq}}  X_+$.  The monoid structure on $S_0$ gives a
\emph{lax} monoidality of $\grr$, and the fact that it is a strong
monoidality is just the fact that the multiplication map
$S_0 \otimes_{\sfC_\ast^{\bZ_\leq}}  S_0 \to S_0$ is an isomorphism.

\begin{remark}
This adjunction is formally quite similar to the adjunction between
abelian groups and $\bF_p$-modules: the left adjoint is $- \otimes_\Z
\bF_p$ and is strong monoidal; the right adjoint is the forgetful map
which ``preserves tensor product up to isomorphism'', but not the monoidal
unit.
\end{remark}

\subsubsection{Induced functors on algebras} \label{sec:filtrations-induced-functors-on-algebras} An operad $\cO$ in $\sfC$ produces operads in $\sfC^{\Z_\leq}$ and $\sfC^{\Z_=}_*$, using the strong monoidal functors $0_* \colon \sfC \to \sfC^{\Z_\leq}$ and $0_* \circ (-)_+ \colon \sfC \to \sfC^{\Z_=}_*$. We shall continue to call these operads $\cO$.

We get a diagram of functors
\begin{equation*}
  \begin{tikzcd}
    \Alg_\cO(\sfC^{\Z_\leq}) \arrow[shift left=.5ex]{r}{\grr} \arrow[shift left=.5ex]{d}{U^\cO} &
    \Alg_\cO(\sfC^{\Z_=}_*)  \arrow[shift left=.5ex]{d}{U^\cO} \arrow[shift left=.5ex]{l}{u} \\
    \sfC^{\Z_\leq} \arrow[shift left=.5ex]{r}{\grr}  \arrow[shift left=.5ex]{u}{F^\cO} &
    \sfC^{\Z_=}_* \arrow[shift left=.5ex]{l}{u} \arrow[shift left=.5ex]{u}{F^\cO} ,
  \end{tikzcd}  
\end{equation*}
commuting up to natural isomorphism, and similarly for $\colim \dashv \const$, and $0_* \dashv 0^*$. This is natural in the operad. We may apply this to an augmentation $\epsilon: \cO \to +$ (e.g.\ the canonical one when $\cO = s(\cC)$ is a non-unitary operad in simplicial sets). This induces a commutative diagram
\[\begin{tikzcd}\Alg_\cO(\sfC^{\Z_\leq}) \arrow[shift left=.5ex]{r}{\grr}   &
\Alg_\cO(\sfC^{\Z_=}_*)  \arrow[shift left=.5ex]{l}{u} \\
\sfC_{\ast}^{\bZ_\leq}  \arrow[shift left=.5ex]{r}{\grr} \arrow{u}{\epsilon^!} & \sfC_{\ast}^{\bZ_=} \arrow[shift left=.5ex]{l}{u} \arrow{u}[swap]{\epsilon^!}\end{tikzcd}\]
and since this diagram of right adjoints commutes, so does the corresponding diagram of left adjoints: this gives a natural isomorphism $Q^\cO\, \grr \cong \grr\, Q^\cO \colon \Alg_\cO(\sfC^{\bZ_\leq}) \to \sfC_{\ast}^{\bZ_=}$. We obtain similar natural isomorphisms $Q^\cO\, \colim \cong \colim\, Q^\cO$, and $Q^\cO\, 0_* \cong 0_*\, Q^\cO$ from the analogous argument with $\colim \dashv \const$, and $0_* \dashv 0^*$.

For a functor to induce a functor between categories of algebras over a non-unitary operad $\cO$, it suffices to produces a morphism of monads. We will do so for $a^*$ with $a \leq 0$, by observing that the morphisms \eqref{eq:24} assemble to a natural transformation
\[\psi \colon \cO a^* \Rightarrow a^* \cO.\]
(That $\cO$ is non-unitary is necessary because the left side has a term $\cO(0)$ but the right side only does so when $a \geq 0$.) The resulting morphism of monads $(a^*,\psi) \colon (\sfC,\cO) \to (\sfC^{\bZ_\leq},\cO) $ yields  
a diagram 
\begin{equation*}
\begin{tikzcd}
\Alg_\cO(\sfC^{\Z_\leq}) \rar{a^*}\dar[swap]{U^\cO} &
\Alg_\cO(\sfC) \dar{U^\cO} \\
\sfC^{\Z_\leq} \rar{a^*} & \sfC
\end{tikzcd}    
\end{equation*}
commuting up to natural isomorphism. Similar considerations apply to $a_* \colon \sfC \to \sfC^{\Z_\leq}$ for $a \geq 0$, but the resulting functor between categories of algebras does not appear to play any important role for $a > 0$.

\begin{remark}
If $M$ is a filtered monoid in sets, then the product of two elements of filtration $-1$ is in filtration $-2$ and so in particular in filtration $-1$, so filtration $-1$ becomes a non-unital monoid; in contrast the filtration +1 subset contains the unit of the monoid, but has no well defined multiplication.
\end{remark}

If $\gX \in \Alg_\cO(\sfC^{\Z_\leq})$ is descendingly filtered, in Section \ref{sec:filtered-integers} we defined $a_!(\gX) \in \Alg_\cO(\sfC_*)$ for each $a \leq 0$. This is intuitively given by ``taking quotient by filtration $a$''. These may be assembled into a pro-object
\begin{equation*}
\gX^\wedge = (a \mapsto a_!(\gX)) \in \cat{pro\text{-}Alg}_\cO(\sfC_*).
\end{equation*}

Both the pro-object $\gX^\wedge$ and the corresponding limit in
$\Alg_\cO(\sfC_*)$ may be regarded as a \emph{completion} of
$\gX(0) = \gX(\infty)$ with respect to the filtration.

\subsection{The canonical multiplicative filtration}
\label{sec:mult-struct-filtr} 
\newglossaryentry{canmult}{%
	name={\ensuremath{(-1)_\mr{alg}^*}},
description={Canonical multiplicative filtration},
type=symbols
}
A first example of an $\cO$-algebra in filtered objects is given by the canonical multiplicative filtration, which is defined when $\cO$ is a non-unitary operad. This is the analogue of the filtration of a non-unital ring (i.e.\ an ideal) by powers of itself. In the context of spectra this filtration has been studied by Harper--Hess \cite{HarperHess} (where it is called the (homotopy) completion tower) and Kuhn--Pereira \cite{KuhnPereira} (where it is called the augmentation ideal filtration).

\subsubsection{Extending $a_*$ to algebras}

For $a < 0$ the left adjoint $a_* \colon \sfC \to \sfC^{\Z_\leq}$ to
$a^*$ does not appear to preserve multiplicative structures in any
interesting way. Nevertheless, for a non-unitary operad $\cO$ the functor
\begin{equation*}
  \Alg_\cO(\sfC^{\Z_\leq})
  \xrightarrow{a^*}
  \Alg_\cO(\sfC)
\end{equation*}
\emph{does} admit a left adjoint
\begin{equation*}
  \Alg_\cO(\sfC) \xrightarrow{a_*^\mathrm{alg}}
  \Alg_\cO(\sfC^{\Z_\leq}),
\end{equation*}
which we now discuss. (Note that it will usually be the case that
$U^\cO a_*^\mathrm{alg} \gR$ and $a_* U^\cO \gR$ are
\emph{not} isomorphic.) If such a left adjoint exists, we must necessarily have a natural isomorphism $ a_*^\mathrm{alg} F^\cO (X) \cong F^\cO (a_* X)$, and conversely, by Proposition \ref{prop:extend-by-density}, this formula may be used to define $a_*^\mathrm{alg}$, by stipulating that it preserve sifted colimits and providing the functor $F^\cO a_* \colon \sfC \to \Alg_{\cO}(\sfC^{\Z_\leq})$
with the structure of a right $\cO$-module functor. We will now provide this right $\cO$-module functor structure. 

\begin{lemma}If $\cO$ is non-unitary and $a \leq 0$, then there is a natural isomorphism
	\[\cO(X) \cong (\cO(a_* X))(a)\]
\end{lemma}

\begin{proof}
We have $\cO(a_*(X)) = \bigsqcup_{n \geq 1}\cO(n) \otimes_{G_n} a_*(X)^{\otimes n}$.
By Lemma \ref{lem:alower-tensor}, we have that $(\cO(n) \otimes_{G_n} a_*(X)^{\otimes n}) \cong (na_*)(\cO(n) \otimes_{G_n} X^{\otimes n})$. Evaluating at $a$ thus gives $\cO(n) \otimes_{G_n} X^{\otimes n}$, since $a \geq na$ if $a \leq 0$ and $n \geq 1$.\end{proof}

There is a natural transformation $m_F \colon a_* \cO \Rightarrow \cO a_*$ adjoint to the isomorphism $\cO(X) \cong (\cO(a_* X))(a)$ described in the previous lemma. Combining this with the counit of the adjunction we obtain a natural transformation
\begin{equation}\label{eqn:muf}
 \mu_F \colon F^\cO a_* \cO \Rightarrow F^\cO \cO a_* = F^\cO U^\cO F^\cO a_* \Rightarrow F^\cO a_*.
\end{equation}

\begin{lemma}The natural transformation (\ref{eqn:muf}) endows $F^\cO a_*$ with the structure of a right $\cO$-module functor.\end{lemma}

\begin{proof}
Let us write $F \coloneqq F^\cO a_*$. We must check that $\mu_F(\mu_F \circ \cO)= \mu_F( F \circ \mu_\cO )$ as natural transformations $F\cO^2 \Rightarrow F$. For this it suffices to check that the following diagram commutes
\[\begin{tikzcd}a_* \cO^2 \arrow[Rightarrow]{rr}{a_*(\mu_\cO)} \dar[Rightarrow,swap]{\mu_F \circ \cO} & & a_* \cO \arrow[Rightarrow]{d}{\mu_F} \\
\cO a_* \cO \rar[Rightarrow]{\cO \circ \mu_F} & \cO^2 a_* \rar[Rightarrow]{\mu_\cO(a_*)} & \cO a_*.
\end{tikzcd} \]
By adjunction, this follows because the following two maps concide:
\begin{align*}\cO^2(X) &\lra \cO(X) \cong \cO(a_*(X))(a)),\\
\cO^2(X) &\cong \cO(a_*(\cO(X))(a)) \cong \cO^2(a_*(X))(a) \lra \cO(a_*(X))(a).\end{align*}
\end{proof}

Note that $0_* \colon \sfC \to \sfC^{\Z_\leq}$ is strong monoidal, and satisfies $0_*(F^\cO(X)) \cong F^\cO(0_*(X))$. It follows that $0_*^\mr{alg} = 0_*$.

\subsubsection{The canonical multiplicative filtration}
\label{sec:canon-mult-filtr}

For $\cO$ a non-unitary operad, in the previous section we have constructed a functor
$\gls{canmult} \colon \Alg_\cO(\sfC) \to
\Alg_\cO(\sfC^{\Z_\leq})$, left adjoint to ``evaluation at
$-1$.''  In this section we shall study this construction more
carefully.
\begin{lemma}\label{lem:-1alg-descending}
  For any $\gR \in \Alg_\cO(\sfC)$, the underlying object $U^\cO(-1)_*^\alg \gR \in \sfC^{\Z_\leq}$ is descendingly filtered.  For any $a \geq -1$ there are natural isomorphisms
\[(U^\cO_{\cO(1)}(-1)_*^\alg \gR)(a) \cong U^\cO_{\cO(1)} \gR.\]
\end{lemma}

\begin{proof} On free $\cO$-algebras we defined $(-1)_*^\alg(F^\cO X)$ to be $F^\cO((-1)_* X)$, so that
	\begin{align*}
	U^\cO_{\cO(1)}(-1)_*^\alg F^\cO X &= U^\cO_{\cO(1)}
	F^\cO (-1)_* X \\
	&= \bigsqcup_{n \geq 1} \cO(n)
	\otimes_{G_n} ((-1)_*X)^{\otimes n} \\
	&\cong \bigsqcup_{n \geq 1}
	(-n)_* (\cO(n)	\otimes_{G_n} X^{\otimes n}),
	\end{align*}
	where the last two objects are $\cO(1)$-modules by the action of $\cO(1)$ on $\cO(n)$.  This gives for any $a \in \Z_\leq$ a natural isomorphism
	\begin{equation*}
	(U^\cO_{\cO(1)}(-1)_*^\alg F^\cO X)(a) \cong \bigsqcup_{n \geq
		-a} \cO(n) \otimes_{G_n} X^{\otimes n},
	\end{equation*}
of $\cO(1)$-modules. This is isomorphic to $U^\cO_{\cO(1)}F^\cO(X)$ as long as $a \geq -1$.	As $(U^\cO_{\cO(1)}(-1)_*^\alg F^\cO(X))(a)$ and $U^\cO_{\cO(1)} F^\cO(X) \cong \cO(X)$ commute with sifted colimits, the conclusion follows by Proposition \ref{prop:extend-by-density}.
\end{proof}

Thus to any $\gR\in  \Alg_\cO(\sfC)$ there is associated a canonical descendingly filtered object $(-1)_*^\alg \gR$ given by
\begin{equation*}
  \gR \cong [(-1)_*^\alg \gR](0) \leftarrow [(-1)_*^\alg \gR](-1) \leftarrow
  [(-1)_*^\alg \gR](-2) \leftarrow \cdots.
\end{equation*}
For any $a \leq 0$, we have a new algebra $a_! (-1)_*^\alg \gR \in \Alg_\cO(\sfC_*)$, whose underlying object in $\sfC_*$ is the quotient $\gR/((-1)_*^\alg \gR)(a-1)$.  As above, these algebras assemble to a pro-object
\[((-1)_*^\alg \gR)^\wedge = (a \mapsto a_! (-1)_*^\alg \gR) \in
\cat{pro}\text{-}\Alg_\cO(\sfC_*),\]
canonically associated to $\gR$. 

\begin{remark}A rough analogy is that of a
local ring $R$, or better its maximal ideal $\mathfrak{m}$: it comes with a canonical
filtration forming a pro-object which we could denote $\mathfrak{m}^\wedge$; if the local ring is Artinian, then $\mathfrak{m}^\wedge$ is pro-constant and isomorphic to $\mathfrak{m}$.  If $R$ is a complete local ring, then $\mathfrak{m}$ is isomorphic to the inverse limit of the pro-object $\mathfrak{m}^\wedge$.\end{remark}

\subsubsection{Its associated graded}\label{sec:CanMultFiltAssocGr} To understand the canonical filtration construction $(-1)^\alg_* \colon \Alg_\cO(\sfC) \to \Alg_\cO(\sfC^{\Z_\leq})$, let us describe its associated graded in terms of the relative indecomposables functor $Q^\cO_{\cO(1)}$. We have
\begin{align*}
U^\cO_{\cO(1)}((-1)_*^\mr{alg}(F^\cO(X))) = U^\cO_{\cO(1)}(F^\cO((-1)_*(X)))
= \bigsqcup_{n \geq 1} (-n)_*(\cO(n) \otimes_{G_n} X^{\otimes n})
\end{align*}
and the $(-1)$st piece of the associated graded is defined by the pushout
\begin{equation*}
\begin{tikzcd} \bigsqcup_{n \geq 2} \cO(n) \otimes_{G_n} X^{\otimes n} \rar \dar & {\bigsqcup_{n \geq 1}} \cO(n) \otimes_{G_n} X^{\otimes n} \dar \\
\bterm \rar & {[\grr\, U^\cO_{\cO(1)}((-1)_*^\mr{alg}(F^\cO(X)))](-1)},
\end{tikzcd}
\end{equation*}
so that there is an identification $(\cO(1) \otimes X)_+ \cong [\grr\, U^\cO_{\cO(1)}((-1)_*^\mr{alg}(F^\cO(X)))](-1)$ as $G_1$ is always trivial. Because $(\cO(1) \otimes X)_+ \cong Q^\cO_{\cO(1)}(F^\cO(X))$, by adjunction we obtain a natural transformation
\[(-1)_*Q^\cO_{\cO(1)}F^\cO \Rightarrow \grr\, U^\cO_{\cO(1)}(-1)_*^\mr{alg}F^\cO \colon \sfC \lra \Alg_{\cO(1)}(\sfC_{\ast}^{\bZ_=}).\]
These functors both preserve sifted colimits, are right $\cO$-module functors, and the natural transformation is one of right $\cO$-module functors, so by Proposition \ref{prop:extend-by-density} this extends to a natural transformation
\[(-1)_*Q^\cO_{\cO(1)} \Rightarrow \grr \, U^\cO_{\cO(1)}(-1)_*^\mr{alg} \colon \Alg_{\cO}(\sfC) \lra \Alg_{\cO(1)}(\sfC_{\ast}^{\bZ_=}).\]
Commuting $U^\cO_{\cO(1)}$ and $\grr$, and taking a further adjoint gives a natural transformation 
  \begin{equation}\label{eqn:nat-trans-1alg}
    F^\cO_{\cO(1)} (-1)_*
    Q^\cO_{\cO(1)} \Rightarrow
    \grr (-1)^\alg_* \colon \Alg_{\cO}(\sfC) \to \Alg_\cO(\sfC^{\Z_=}_*).
  \end{equation}

\begin{proposition}\label{prop:associated-graded-of-canonical-filt} The natural transformation (\ref{eqn:nat-trans-1alg}) is a natural isomorphism.
\end{proposition}
\begin{proof}
Using Proposition \ref{prop:extend-by-density} it suffices to verify this on free $\cO$-algebras. Furthermore, as $U^\cO$ creates isomorphisms it is enough to check after applying this. As above we have
\[U^\cO\grr (-1)_*^\mr{alg}(F^\cO(X))(-k) = (\cO(k) \otimes_{G_k} X^{\otimes k})_+.\]
On the other hand we have
\[(-1)_* Q^\cO_{\cO(1)}(F^\cO(X)) = (-1)_*(\cO(1) \otimes X)_+ = \cO(1)_+ \wedge ((-1)_*(X)_+)\]
which is $F^{\cO(1)}((-1)_*(X)_+)$.  So, as $F^\cO_{\cO(1)} \circ F^{\cO(1)} \cong F^\cO$, we have
\[F^\cO_{\cO(1)}((-1)_* Q^\cO_{\cO(1)}(F^\cO(X))) \cong F^\cO((-1)_*(X)_+)\]
whose underlying object in grading $-k$ is $(\cO(k) \otimes_{G_k} X^{\otimes k})_+$. It is easy to see that (\ref{eqn:nat-trans-1alg}) induces the identity map under these identifications.
\end{proof} 
\section{Cell attachments}
\label{sec:cell-attachments}

In this section we explain how to attach a $T$-algebra cell to a $T$-algebra. When the monad comes from an operad $\cO$, we describe a canonical filtration on a cell attachment. After this, we define cellular and CW $\cO$-algebras, and the skeletal filtration of a CW algebra. As before, in this section we work in a category $\sfC = \sfS^\sfG$ where $\sfS$ satisfies the axioms of Section \ref{sec:axioms-of-cats}.

\subsection{Cell attachments for sifted monads} \label{sec:cell-attachments-monads}

\subsubsection{The definition of a cell attachment}\label{sec:cell-def} For $\gX_0 \in \Alg_T(\sfC)$ the data for a $T$-cell attachment to $\gX_0$ is given by:
\begin{itemize}
\item a cofibration of simplicial sets $\partial D^d \hookrightarrow D^d$, whose geometric realization is homeomorphic to the inclusion of the boundary of the $d$-disk,
	\item an object $g \in \sfG$, and
	\item a morphism $e \colon \partial D^d \to U^T(\gX_0)(g)$. (Here $\partial D^d$ is considered as an object of $\sfS$ via the functor $s \colon \cat{sSet} \to \sfS$.)
\end{itemize}

\newglossaryentry{dgd}{%
	name={\ensuremath{D^{g,d}}},
description={$d$-disk in grading $g$},
type=symbols
}
\newglossaryentry{partdgd}{%
	name={\ensuremath{\partial D^{g,d}}},
description={$(d-1)$-sphere in grading $g$},
type=symbols
}
There is an adjunction $g_* \dashv g^*$, and we will write
\begin{align*}
\gls{dgd} &\coloneqq g_*(D^d),\\
\gls{partdgd} &\coloneqq g_*(\partial D^d),
\end{align*}
and for later use, also define the pointed object $S^{g,d} \coloneqq D^{g,d}/\partial D^{g,d} \in \sfC_{\ast}$.

Using the adjunction $g_* \dashv g^*$ the morphism $e$ corresponds to a morphism  $\partial D^{g,d} \to U^T(\gX_0)$, and using the adjunction $F^T \dashv U^T$, this in turn corresponds to a morphism $F^T(\partial D^{g,d}) \to \gX_0$ in $\Alg_T(\sfC)$ which we shall also denote by $e$. We then define $\gX_1 \in \Alg_T(\sfC)$ to be the following pushout in $\Alg_T(\sfC)$
\begin{equation}\label{eq:6}
  \begin{aligned}
    \begin{tikzcd}
      F^T(\partial D^{g,d}) \dar \rar{e} & \gX_0 \dar \\
      F^T(D^{g,d}) \rar & \gX_1. \end{tikzcd}
  \end{aligned}
\end{equation}

\newglossaryentry{cellattachment}{%
	name={\ensuremath{\cup^T_e}},
description={$T$-cell attachment along $e$},
type=symbols
}
\begin{definition}
Given the pushout diagram (\ref{eq:6}), we say that \emph{$\gX_1$ is obtained from $\gX_0$ by attaching a $T$-cell of dimension $(g,d)$ along $e$}\index{cell attachment} and we shall often denote $\gX_1$ by $\gX_0 \gls{cellattachment} \gD^{g,n}$.
\end{definition}

The proof of existence of colimits in $\Alg_T(\sfC)$ in Lemma \ref{lem:algt-complete-cocomplete} was constructive, so we shall concretely describe the underlying object of the pushout~(\ref{eq:6}) in $\Alg_T(\sfC)$ for the benefit of the reader. Given a diagram 
\[D^{g,d} \hookleftarrow \partial D^{g,d} \overset{e}{\lra} X_0\]
as above, let us write $X_0 \cup_e D^{g,n}$ for the pushout in $\sfC$.  Then the $T$-algebra $\gX_1$ has underlying object of $\sfC$ given by the reflexive coequalizer
\begin{equation}\label{eq:9}
  \begin{tikzcd}
    T(T(X_0) \cup_e D^{g,d}) \arrow[shift left=1ex]{r} \arrow[shift left=-1ex]{r}& T(X_0 \cup_e D^{g,d}) \lar \rar &
    X_1,
  \end{tikzcd}
\end{equation}
where: the top arrow is obtained by applying $T$ to the induced map on pushouts $\mu^T_{X_0} \cup_{\partial D^{g,d}} D^{g,d} \colon T(X_0) \cup_e D^{g,d}\to X_0 \cup_e D^{g,d}$; the bottom arrow is obtained by applying $T$ to $i \colon T(X_0) \cup_e D^{g,d} \to T(X_0 \cup_e D^{g,d})$ and then composing with the component of the natural transformation $\mu^T \colon T^2 \Rightarrow T$ at $X_0 \cup_e D^{g,d}$; the reflection is obtained by applying $T$ to $i_{X_0} \cup_e D^{g,d}$, with $i_{X_0} \colon X_0 \to T(X_0)$ the unit of the monad.

\subsubsection{Cell attachments and change-of-monad}  If $\phi \colon T \to T'$ is a morphism of sifted monads and $\phi_* \colon \Alg_T(\sfC) \to \Alg_{T'}(\sfC)$ is the left adjoint in the corresponding change-of-monads adjunction, then $\phi_*$ preserves pushout diagrams like any left adjoint functor. Hence if we apply it
to the diagram (\ref{eq:6}) and use the natural isomorphism $\phi_* F^T \cong F^{T'}$ as discussed in Example \ref{exam:change-monads}, then we obtain 
\begin{equation*}
    \begin{tikzcd}
	F^{T'}(\partial D^{g,d}) \dar \rar{\phi_*(e)} & \phi_*(\gX_0) \dar \\
	F^{T'}(D^{g,d}) \rar & \phi_*(\gX_1),
	\end{tikzcd}
\end{equation*}
a pushout diagram in $\Alg_T(\sfC)$. That is, if $\gX_1$ is obtained from $\gX_0$ by attaching a cell along $e$ in $\Alg_T(\sfC)$, then $\phi_*(\gX_1)$ is obtained from $\phi_*(\gX_0)$ by attaching a cell along $\phi_*(e)$ in $\Alg_{T'}(\sfC)$.

\begin{lemma}If $\phi \colon T \to T'$ is a morphism of sifted monads, then $\phi_* \colon \Alg_T(\sfC) \to \Alg_{T'}(\sfC)$ preserves cell attachments.\end{lemma}

\subsubsection{Cell attachments and indecomposables}\label{sec:cell-attachments-and-indecomposables}

\index{cell attachment!effect on indecomposables}Since the $T$-indecomposables $Q^T$ is a special case of the change-of-monads construction, the effect on indecomposables of a cell attachment in $\Alg_T(\sfC)$
is a cell attachment in $\Alg_+(\sfC) = \sfCstar$. This is captured by the slogan that ``$Q^T$ transforms cell structures in $\Alg_T(\sfC)$ to cell structures in $\sfCstar$.''

Let us discuss this special case more explicitly. If we apply the left adjoint $Q^T$ to~(\ref{eq:6}), using the formula $Q^T F^T(-) \cong (-)_+$, we obtain a pushout diagram in $\sfCstar$
\begin{equation*}
\begin{tikzcd}
	\partial D^{g,d}_+ \rar \dar & Q^T(\gX_0)\dar \\
	D^{g,d}_+ \rar & Q^T(\gX_1).
\end{tikzcd}
\end{equation*}

For a general diagram category $\sfG$ such a cell attachment may be quite complicated, but if $\sfG$ is a groupoid then one has the following more concrete description. If $g \not\cong h \in \sfG$ then the map $Q^T(\gX_0)(h) \to Q^T(\gX_1)(h)$ is an isomorphism. If $g \cong h \in \sfG$ then the difference between $Q^T(\gX_0)(h)$ and $Q^T(\gX_1)(h)$ is described by  a pushout diagram in $\sfS_*^{\sfG(h,h)}$
\begin{equation}\label{eq:2}
\begin{tikzcd}
	\partial D^{d}_+ \wedge \sfG(h,h)_+ \rar \dar  & Q^T(\gX_0)(h)\dar \\
	D^{d}_+ \wedge \sfG(h,h)_+ \rar & Q^T(\gX_1)(h),
\end{tikzcd}
\end{equation}
that is, $Q^T(\gX_1)(h)$ is obtained from $Q^T(\gX_0)(h)$ by ``attaching a free $\sfG(h,h)$-equivariant $d$-cell.''

\subsubsection{Cell attachments for operads and change-of-diagram-category}
 
We now show that change-of-diagram-category preserves $\cO$-algebra cell attachments when $\cO$ is an operad. 

In Section \ref{sec:simplicial-operads} we saw that if $p \colon \sfG \to \sfG'$ is strong $k$-monoidal, there is a functor $p_* \colon \Alg_\cO(\sfS^\sfG) \to \Alg_\cO(\sfS^{\sfG'})$ which is a left adjoint and satisfies $p_*(F^\cO(X)) \cong F^\cO(p_*(X))$. Thus $p_*$ preserves pushouts, so applying it to the diagram (\ref{eq:6}) and using that $p_* g_* = p(g)_*$, we obtain a pushout diagram in $\Alg_\cO(\sfS^{\sfG'})$ 
\begin{equation*}
  \begin{aligned}
    \begin{tikzcd}
      F^\cO(\partial D^{p(g),d}) \dar \rar{p_*(e)} & p_*(\gX_0) \dar \\
      F^\cO(D^{p(g),d}) \rar & p_*(\gX_1). \end{tikzcd}
  \end{aligned}
\end{equation*}
That is, $p_*(\gX_1)$ is obtained from $p_*(\gX_0)$ by attaching a $\cO$-algebra  $(p(g), d)$-cell:

\begin{lemma}If $p \colon \sfG \to \sfG'$ is a strong $k$-monoidal functor and $\cO$ is an operad in $\sfS$, then $p_* \colon \Alg_\cO(\sfS^{\sfG}) \to \Alg_\cO(\sfS^{\sfG'})$ preserves cell attachments.\end{lemma}

\subsection{Ascending filtrations from cell attachments}
\label{sec:filtr-cell-attachm}

In this section we shall specialize to the case of $\cO$-algebras, where $\cO$ is an operad in $\sfC = \sfS^\sfG$. We shall study the filtration on a cell attachment.

\subsubsection{The filtration on a cell attachment} \label{sec:filtration-on-pushout-algebras} 

In Section \ref{sec:colimit-2} we saw that the functor $\colim \colon \sfC^{\Z_\leq} \to \sfC$ is (strong) monoidal and commutes with colimits, so from the formula
\[\cO(X) = \bigsqcup_{n \geq 0} 0_*(\cO(n)) \otimes_{G_n} X^{\otimes n}\]
for the monad $\cO$ on $\sfC^{\Z_\leq}$, and the fact that $\colim \circ\, 0_* = \mr{id}$, we conclude that $\colim \cO \cong \cO \colim$. Thus there is a functor $\colim \colon \Alg_\cO(\sfC^{\Z_\leq}) \to \Alg_\cO(\sfC)$ which commutes with $U^\cO$.

If $\gS \in \Alg_\cO(\sfC^{\Z_\leq})$ is a filtered $\cO$-algebra, we think of $\colim \gS \in \Alg_\cO(\sfC)$ as the ``underlying'' $\cO$-algebra and an isomorphism $\gR \overset{\sim}\to \colim \gS$ in $\Alg_\cO(\sfC)$ as specifying a multiplicative filtration on $\gR$.  We now describe such a multiplicative filtration on a cell attachment in $\Alg_\cO(\sfC)$. Let us return to the situation of Section \ref{sec:cell-def}: we have an $\gR_0 \in \Alg_\cO(\sfC)$, and a diagram 
\begin{equation}\label{eq:20}
\begin{tikzcd} F^\cO(\partial D^{g,d}) \dar \rar{e} & \gR_0 \\
F^\cO(D^{g,d}), & \end{tikzcd}
\end{equation}
whose pushout in $\Alg_\cO(\sfC)$ we denoted $\gR_1$. To obtain multiplicative filtration on $\gR_1$, i.e.\ lift it to $\Alg_\cO(\sfC^{\Z_\leq})$, it suffices to lift the defining pushout diagram~(\ref{eq:20}) to a diagram in $\Alg_\cO(\sfC^{\bZ_{\leq}})$. 

\newglossaryentry{veeo}{%
	name={\ensuremath{\vee^\cO}},
description={Coproduct in $\cO$-algebras in pointed category},
type=symbols
}
If we did this using the strong monoidal functor $0_*$, the result would be isomorphic to $0_*\gR_1$ since $0_*$ commutes with pushouts as a left adjoint. Instead, in \eqref{eq:20} we replace $\gR_0$ by $0_*\gR_0$, $D^{g,n}$ by $1_* D^{g,d}$, and $\partial D^{g,n}$ by $1_* \partial D^{g,d}$. The free algebra $F^\cO(1_* D^{g,d})$ has underlying object $F^\cO(D^{g,d})$, and its filtration is not concentrated in any particular degree. Consequently the pushout in $\Alg_\cO(\sfC^{\Z_\leq})$ of the diagram
\begin{equation}
\label{eq:21}
\begin{tikzcd}F^\cO(1_* \partial D^{g,n}) \dar \rar{e} & 0_*\gR_0 \\
F^\cO(1_* D^{g,n}), & \end{tikzcd}
\end{equation}
which we shall denote by $\mr{f}\gR_1$ and call the \emph{cell attachment filtration}\index{filtration!cell attachment}, also has a filtration which is not concentrated in any particular degree. This may be seen from the following description of its associated graded, where $\gls{veeo}$ denotes the coproduct in $\Alg_\cO(\sfC_\ast^{\bZ_=})$:

\begin{theorem}\label{thm:cell-attachment-associated-graded} In $\Alg_\cO(\sfC_{\ast}^{\bZ_=})$ there is an isomorphism  
	\[\grr(\mr{f}(\gR_1)) \cong 0_*(\gR_0)_+ \vee^{\cO} F^\cO(1_*(S^{g,d})).\]
\end{theorem}

\begin{proof}Since $\grr$ commutes with colimits and with $F^\cO$, we have a pushout diagram in $\Alg_{\cO}(\sfC^{\bZ_{=}}_*)$
	\[\begin{tikzcd}F^\cO(1_* \partial D^{g,d}_+) \rar \dar & 0_*(\gR_0)_+ \dar \\
	F^\cO(1_* D^{g,d}_+) \rar & \grr(\mr{f}(\gR_1)).
	\end{tikzcd}\]
	The top map is adjoint to the unique map $\partial D^{g,d} \to 1^*(0_*(\gR_0)_+) = *$ in $\sfC_*$, so it is the unique map to $* \to 0_*(\gR_0)_+$. Thus the pushout is isomorphic to the coproduct in $\Alg_{\cO}(\sfC_*^{\bZ_{=}})$ of $0_*(\gR_0)_+$ with the cofiber of the left map. As a left adjoint $F^\cO$  preserves cofibers, giving the asserted answer.
\end{proof}

\subsubsection{The stages of the filtration on powers of pushouts} \label{sec:filtration-powers} Next we describe the filtration steps of the cell attachment filtration. To do so, we first understand the filtration steps of a tensor power of a pushout.

Given a map $i \colon X_0 \to X_1$, we define a filtered object $f[i]$ by setting 
\[f[i](n) = \begin{cases} \binit & \text{if $n < 0$,} \\
X_0 & \text{if $n \leq 0$,} \\
X_1 & \text{if $n > 0$,} \end{cases}\]
with non-trivial structure maps induced by $i \colon X_0 \to X_1$. That is, $f[i]$ is the pushout
\[\begin{tikzcd} 1_* X_0 \rar \dar & 0_*X_0 \dar \\
1_* X_1 \rar & f[i].\end{tikzcd}\]

We shall describe the induced filtration on $f[i]^{\otimes n}$,\index{filtration!on powers of pushouts} with $\otimes$ given by Day convolution on $\sfC^{\bZ_{\leq}}$, in terms of pushout-products $\Box$ on $\sfC^{\bZ_{\leq}}$. 

In the case of a general cocomplete category $\sfD$ with monoidal structure, the \emph{pushout-product}\index{pushout-product} is a monoidal structure on the arrow category $\sfC^{[1]}$, where $[1]$ is the diagram category $0 \to 1$. It is given by Day convolution with respect to the symmetric monoidal functor $\min \colon [1] \times [1] \to [1]$, so inherits many of the properties of $(\sfD,\otimes,\bunit)$, e.g.\ it is symmetric monoidal if $\otimes$ is. Explicitly, for two morphisms $f \colon X_0 \to X_1$ and $g \colon Y_0 \to Y_1$, $f \Box g$ is the induced map in the following diagram:
\[\begin{tikzcd} X_0 \otimes Y_0 \rar{f \otimes \mr{id}} \dar{\mr{id} \otimes g} & X_1 \otimes Y_0 \dar \arrow[bend left = 25]{rdd}[description]{\mr{id} \otimes g} & \\
X_0 \otimes Y_1 \rar \arrow[bend left = -15]{drr}[description]{f \otimes \mr{id}} & (X_1 \otimes Y_0) \sqcup_{X_0 \otimes Y_0} (X_0 \otimes Y_1) \arrow{rd}[dashed,description]{f \Box g}  & \\
& & X_1 \otimes Y_1, \end{tikzcd}\]
where we shall usually shorten $(X_1 \otimes Y_0) \sqcup_{X_0 \otimes Y_0} (X_0 \otimes Y_1)$ to $X_1 \Box Y_1$ (suppressing from notation the fact that it depends on $f$ and $g$).

Let us return to our study of $f[i]^{\otimes n}$. In filtration degree $c$, this is given by
\[f[i]^{\otimes n}(c) \cong \underset{c_1+\ldots+c_n \leq c}{\mr{colim}} \left(f[i](c_1) \otimes \ldots \otimes f[i](c_n)\right).\]
Since the maps $f(i)[c] \to f(i)[c+1]$ are isomorphisms when $c \leq -1$ or $c \geq 1$, the maps $f[i]^{\otimes n}[c] \to f[i]^{\otimes n}[c+1]$ are isomorphism for $c \leq -1$ or $c \geq n$, and we may compute the non-trivial filtration steps by restricting from $\bZ_{\leq}$ to the subcategory $[1]$ in the colimit, as for $0 \leq c < n$, the inclusion of $(i_1,\ldots,i_n) \in \{0,1\}^n$ such $i_1+\ldots+i_n \leq c$ into the $(c_1,\ldots,c_n) \in (\bZ^{\leq})^n$ such that $c_1+\ldots+c_n \leq c$ is final functor. Thus it suffices to give for $I = (i_1,\ldots,i_n) \in \{0,1\}^n$ an inductive construction of
\[\underset{i_1+\ldots+i_n \leq c}{\mr{colim}} \left(X_{i_1} \otimes \ldots \otimes X_{i_n}\right).\]

To simplify notation, we shall write $X_I^{\otimes n} \coloneqq X_{i_1} \otimes \ldots \otimes X_{i_n}$. To obtain the inductive construction we use that $[1]^n$ is a Reedy category by setting $I = (i_1,\ldots,i_n)$ to have degree $|I| = i_1+\ldots+i_n$ (see Section \ref{sec:geom-rel-generalities} for a discussion of this theory in the context of geometric realization). Our filtration is the skeletal filtration of the colimit over this Reedy category:
\[\mr{sk}_c(X_\bullet^{\otimes n}) \cong \underset{\substack{I \in \{0,1\}^n \\ i_1+\ldots+i_n \leq c}}{\mr{colim}} X_I^{\otimes n}.\]

Then we can compute $\mr{sk}_{c}(X_\bullet^{\otimes n})$ from $\mr{sk}_{c-1}(X_\bullet^{\otimes n})$ in terms of $\bigsqcup_{|I|=c} X_I^{\otimes n}$ and the $c$th latching object $L_c(X_\bullet^{\otimes n})$:
\begin{equation}\label{eqn:general-tensor-pushout} \begin{tikzcd} L_c(X_\bullet^{\otimes n}) \rar \dar & \mr{sk}_{c-1}(X_\bullet^{\otimes n}) \dar \\
\underset{|I|=c}{\bigsqcup} X_I^{\otimes n} \rar & \mr{sk}_c(X_\bullet^{\otimes n}).\end{tikzcd}\end{equation}
This latching object may be computed as
\[L_c(X_\bullet^{\otimes n}) \coloneqq \underset{\substack{I \in \{0,1\}^n \\ i_1+\ldots+i_n = c}}{\bigsqcup} \left(\underset{\substack{J \in \{0,1\}^n,\, j_r \leq i_r \\ j_1+\ldots+j_r \leq c-1}}{\mr{colim}} X_J^{\otimes n}\right),\]
and if the monoidal structure on $\sfC$ has a braiding then $L_c(X_\bullet^{\otimes n})$ is a ${n \choose c}$-fold disjoint union of terms isomorphic to $X^{\otimes n-c}_0 \otimes X_1^{\Box c}$, using the canonical isomorphism
\[\underset{\substack{J \in \{0,1\}^c \\ j_1+\ldots+j_c \leq c-1}}{\colim} X_J^{\otimes c} \cong X_1^{\Box c}.\]

If we let $G_{n-c,c} \leq G_n$ denote the setwise stabilizer of the set $\{n-c+1,\ldots,n\}$ of the last $c$ elements (this is just $\fS_{n-c} \times \fS_c$ if $k=\infty$, but more complicated in the case $k=2$), then (\ref{eqn:general-tensor-pushout}) may then be written as
\begin{equation}\label{eqn:general-tensor-pushout-2} \begin{tikzcd} G_{n} \times_{G_{n-c,c}} \left( X_0^{\otimes n-c} \otimes X_1^{\Box c} \right) \rar \dar[swap]{G_{n} \times_{G_{n-c,c}}\left( X_0^{\otimes n-c} \otimes i^{\Box c} \right)} & f[i]^{\otimes n}(c-1) \dar \\
G_{n} \times_{G_{n-c,c}} \left( X_0^{\otimes n-c} \otimes X_1^{\otimes c} \right) \rar & f[i]^{\otimes n}(c).\end{tikzcd}\end{equation}

\begin{lemma}\label{lem:fi-sym-seq} For $c \geq 1$, and any symmetric sequence $\cX$ in $\sfC$, there is a pushout diagram
\[\begin{tikzcd}\bigsqcup_{n \geq c} \cX_n \otimes_{G_{n-c,c}} \left(X_0^{\otimes n-c} \otimes X_1^{\Box c}\right) \dar \rar & 0_*(\cX)(f[i])(c-1) \dar \\
\bigsqcup_{n \geq c} \cX_n \otimes_{G_{n-c,c}} \left(X_0^{\otimes n-c} \otimes X_1^{\otimes c}\right) \rar & 0_*(\cX)(f[i])(c). \end{tikzcd}\]\end{lemma}
\begin{proof}
We use the formula $(0_*\cX)(Y)(c) = \bigsqcup_{n \geq 0} \cX_n \otimes_{G_n} (Y^{\otimes n}(c))$. Applying $\cX_n \otimes_{G_n} -$ to the pushout square \eqref{eqn:general-tensor-pushout-2} in $\sfC^{G_n}$ we obtain a pushout square in $\sfC$ (as $- \otimes -$ preserves colimits in each variable). Next we take the coproduct over $n \geq c$, using that $f[i]^{\otimes n}(c-1) \to f[i]^{\otimes n}(c)$ is an isomorphism for $n < c$.
\end{proof}

\subsubsection{The stages of the filtration on a pushout in algebras}\label{sec:filt-stages-pushout}

We now describe the stages of the cell attachment filtration $\mr{f}(\gR_1)$ of Section \ref{sec:filtration-on-pushout-algebras} in the case that $k=2,\infty$. To do so, we give an alternative expression for the underlying object of the coproduct $\gR \sqcup^\cO F^\cO(Y_1)$ in $\cO$-algebras. 

Let $\tilde{G}_{n-c,c}$ be the kernel of the homomorphism $G_{n-c,c} \to G_c$ obtained by deleting the first $n-c$ strands if $k=2$ or the action on the first $n-c$ elements when $k=\infty$ (so that $\tilde{G}_{n-c,c} \cong \fS_{n-c}$ if $k=\infty$). A $G_c$-action on $\cO(n) \otimes_{\tilde{G}_{n-c,c}} X^{\otimes n-c}$ remains. We then define the following right $\cO$-module functor
\begin{align}
\label{eqn:eco}
\begin{aligned}
\mr{Env}_c(\cO) \colon \sfC &\lra \sfC^{G_c} \\
X & \longmapsto \bigsqcup_{n \geq c} \cO(n) \otimes_{\tilde{G}_{n-c,c}} X^{\otimes n-c} 
\end{aligned}
\end{align}
which visibly commutes with sifted colimits. Extending by density under sifted colimits as in Section \ref{sec:funct-out-eilenb} we obtain functors $\mr{Env}_c^\cO \colon \Alg_T(\sfC) \to \sfC^{G_c}$ satisfying $\mr{Env}_c^\cO(F^\cO(X)) = \mr{Env}_c(\cO)(X)$ for each $c \geq 0$, which can be assembled into a single functor
\begin{align*} \mr{Env}^\cO \colon \Alg_\cO(\sfC) &\lra \cat{FB}_k(\sfC) \\
\gR &\longmapsto \left(\mr{Env}^\cO(\gR) \colon c \mapsto \mr{Env}_c^\cO(\gR)\right).\end{align*}

\begin{remark}If $\sfC$ is $(k+1)$-monoidal this forms an operad, and the category of $\mr{Env}^\cO(\gR)$-algebras is equivalent to the category of $\cO$-algebras under $\gR$ \cite[Lemma 1.7]{BergerMoerdijkDerived}.\end{remark}

\newglossaryentry{sqcupo}{%
	name={\ensuremath{\sqcup^\cO}},
description={Coproduct in $\cO$-algebras},
type=symbols
}
\begin{lemma}If $k=2,\infty$, there is an isomorphism
	\[\gR \gls{sqcupo} F^\cO(X_1) \cong \mr{Env}^\cO(\gR)(X_1),\]
which is natural in $\gR$ and $X_1$.\end{lemma}

\begin{proof}
	Since both
	\[X \longmapsto F^\cO(X) \sqcup^\cO F^\cO(X_1) \qquad \text{ and } \qquad X \longmapsto \mr{Env}^\cO(F^\cO(X))(X_1)\]
	are right $\cO$-module functors which commute with sifted colimits, by Proposition \ref{prop:extend-by-density} it suffices to establish this for free algebras. In that case both sides are naturally isomorphic to
	\[\bigsqcup_{n_1,n_2 \geq 0} \cO(n_1+n_2) \otimes_{G_{n_1,n_2}} X^{\otimes n_1} \otimes X_1^{\otimes n_2}.\qedhere\]
\end{proof}

Let us now consider the defining pushout diagram for the cell attachment filtration
\[\begin{tikzcd} F^\cO(1_* X_0) \rar \dar[swap]{F^\cO(i)} & 0_* \gR_0 \dar \\
F^\cO(1_* X_1) \rar & \mr{f}(\gR_1).\end{tikzcd}\]
It may factored as a composition of two pushout diagrams
\[\begin{tikzcd} F^\cO(1_* X_0) \rar \dar[swap]{F^\cO(i)} & F^\cO(0_* X_0) \dar \rar &0_* \gR_0 \dar \\
F^\cO(1_* X_1) \rar & F^\cO(f[i]) \rar & \mr{f}(\gR_1),\end{tikzcd}\]
and we may restrict our attention on the right pushout square. Using the results of Section \ref{sec:filtration-powers}, the following proposition gives the filtration steps. This is the starting point of the homotopical analysis of algebras over operads, and versions of it appear in the unpublished work of Spitzweck, \cite[Chapter 18]{FresseBook}, \cite[Section 7.3]{HarperOperads}, and \cite[Section 5]{HarperHess}.

\begin{proposition}\label{prop:filt-stages-pushout} For all $c \geq 1$ there is a pushout diagram
	\begin{equation}\label{eqn:filt-stages-pushout} \begin{tikzcd}\mr{Env}_c^\cO(\gR_0) \otimes_{G_c} X_1^{\Box c} \rar \dar & \mr{f}(\gR_1)(c-1) \dar \\
	\mr{Env}_c^\cO(\gR_0) \otimes_{G_c} X_1^{\otimes c} \rar & \mr{f}(\gR_1)(c). \end{tikzcd}\end{equation}
\end{proposition}

\begin{proof}By Lemma \ref{lem:fi-sym-seq}, the following two functors $\sfC^{X_0 \downarrow -} \to \sfC$ are naturally isomorphic; the first is given by
	\[G_1(X_0 \to X) \coloneqq \cO(0_* X \sqcup_{1_* X_0} 1_* X_1)(c)\]
	and the second, $G_2$, defined by the pushout square
	\[\begin{tikzcd}\mr{Env}_c(\cO)(X) \otimes_{G_c} X_1^{\Box c} \rar \dar & \cO(0_* X \sqcup_{1_* X_0} 1_* X_1)(c-1) \dar \\
	\mr{Env}_c(\cO)(X) \otimes_{G_c} X_1^{\otimes c} \rar & G_2(X_0 \to X).\end{tikzcd}\]
	Both are right $\cO$-module functor preserving sifted colimits. The extensions of $G_1$ and $G_2$ to functors $\Alg_\cO(\sfC)^{F^\cO(X_0) \downarrow -} \to \sfC$ by density are $\mr{f}(\gR_1)(c)$ and the pushout of (\ref{eqn:filt-stages-pushout}) respectively, and we may conclude that these are naturally isomorphic as well.\end{proof}

\subsection{Cellular algebras and CW-algebras}
\label{sec:cw-algebras}

We shall define two notions of $\cO$-algebras built using cells: (i) cellular $\cO$-algebras, which are obtained by iterated cell attachments starting at the initial object $\binit$, and (ii) CW $\cO$-algebras, which are obtained by iterated cell attachments respecting a skeletal filtration. This not only means that the cells are attached in order of dimension, but also imposes restrictions on the possible attaching maps. 

\subsubsection{Cellular maps and CW-structures on maps}
For later use, we shall discuss the more general notions of a cellular map and a CW-structure on a map. Here one starts with an $\cO$-algebra $\gR$ instead of the initial object $\binit$.

\begin{definition}
A map $f \colon \gR \to \gS$ of $\cO$-algebras is said to \emph{cellular} if it is the transfinite composition of cell attachments. More precisely, we mean that there exists a diagram
	\[\begin{tikzcd}\gR = \gR_{-1} \arrow{d}[description]{f} \rar & \gR_0 \arrow{dl}[description]{f_0} \rar & \gR_1 \arrow{dll}[description]{f_1} \rar & \cdots \\
	\gS\end{tikzcd}\]
indexed by some ordinal $\kappa$, such that (i) $\mr{colim}_{i \in \kappa} f_i$ is an isomorphism, and (ii) for each successor ordinal $i \in \kappa$ there is a pushout diagram in $\Alg_{\cO}(T)$
		\[\begin{tikzcd}F^\cO\left(\bigsqcup_{\alpha \in I_i} \partial D^{g_\alpha,d_\alpha} \right) \rar \dar & \gR_{i-1} \dar \\
		F^\cO\left(\bigsqcup_{\alpha \in I_i} D^{g_\alpha,d_\alpha}\right) \rar & \gR_i, \end{tikzcd}\]
for some set of morphisms $\{h_\alpha \colon \partial D^{g_\alpha,d_\alpha} \to R_{i-1}\}_{\alpha \in I_{i}}$, while for each limit ordinal $i \in \kappa$, we have that $f_i \colon \gR_i \to \gS$ is the colimit of $f_{i'} \colon \gR_{i'} \to \gS$ for $i'<i$.
\end{definition}

\begin{definition}An $\cO$-algebra $\gR$ is said to be \emph{cellular}\index{$\cO$-algebra!cellular} if the map $\binit \to \gR$ is cellular.\end{definition}

Cellular $\cO$-algebras do not admit a useful filtration, even if one demands the cells are attached in increasing order of dimension. This is because by definition they only give a filtered object in $\cO$-algebras, i.e.\ an object of $\Alg_\cO(\sfC)^{\bZ_\leq}$, and its associated graded need \emph{not} be an $\cO$-algebra.

This defect is addressed by the notion of a CW $\cO$-algebra. This will be defined in terms of their skeletal filtration, which is a filtered $\cO$-algebra, i.e.\ an object of $\Alg_\cO(\sfC^{\bZ_\leq})$. It shall be obtained by attaching $d$-dimensional cells in filtration $d$, along attaching maps into filtration $d-1$. To make this precise, given a cofibration $\partial D^d \hookrightarrow D^d$ of simplicial sets, whose geometric realization is homeomorphic to the $d$-disc and its boundary, we form a filtered simplicial set as follows: we put the source $\partial D^d$ in filtration $d-1$
\[\partial D^d[d-1] \coloneqq (d-1)_*(\partial D^d),\]
and let the target be filtered by putting the subset $\partial D^d$ in filtration $(d-1)$ and the remainder in filtration $d$. That is, $D^d[d]$ is the pushout in $\cat{sSet}^{\bZ_\leq}$
\[\begin{tikzcd}d_*(\partial D^d) \rar \dar & (d-1)_*(\partial D^d) \dar \\
d_*(D^d) \rar & D^d[d].\end{tikzcd}\]

We will consider objects $X$ of the category $\sfC^{\Z_\leq} = (\sfS^\sfG)^{\Z_\leq} = \sfS^{\sfG \times \Z_\leq}$, and we write their values as $X(g,n)$ for $(g,n) \in \sfG \times \Z_\leq$. As usual, we will implicitly apply $\cat{sSet} \to \sfS \to \sfC$ to consider simplicial sets as objects of $\sfS$ or $\sfC$, and hence in particular we consider $D^d[d]$ as an object of $\sfC^{\Z_\leq}$.

\begin{definition}\label{def:filtered-cw-attachment-dim-d} A structure of \emph{filtered CW attachment of dimension $d$} on a morphism  $f \colon \gR \to \gS$ in  $\Alg_\cO(\sfC^{\bZ_{\leq}})$ consists of the following data:
	\begin{enumerate}[(i)]
		\item a set $I_d$,
		\item a collection of cofibrations of simplicial sets $\{\partial D^{d}_\alpha \hookrightarrow D^d_\alpha\}_{\alpha \in I}$ each of whose geometric realizations is homeomorphic to the inclusion of the boundary of the $d$-disk,
		\item a collection of objects $\{g_\alpha\}_{\alpha \in I_d}$ of $\sfG$,
		\item a collection of morphisms $e_\alpha \colon \partial D^d_\alpha \to \gR(g_\alpha,d-1)$ in $\sfS$, adjoint to morphisms $\partial D^{g_\alpha d}_\alpha[d-1] \to R$ in $\sfC^{\Z_\leq}$,
		\item a pushout diagram
		\[\begin{tikzcd}F^\cO\left(\bigsqcup_{\alpha \in I_d} \partial D^{g_\alpha,d}_\alpha[d-1] \right) \rar \dar & \gR \dar{f} \\
		F^\cO\left(\bigsqcup_{\alpha \in I_d} D^{g_\alpha,d}_\alpha[d]\right) \rar & \gS. \end{tikzcd}\]
	\end{enumerate}
\end{definition}

\newglossaryentry{skeletal}{%
	name={\ensuremath{\mr{sk}}},
description={Skeletal filtration on CW-object},
type=symbols
}
\begin{definition}\label{def:relative-cw-structure} 
A \emph{relative CW-structure}\index{CW-!structure} on a morphism $f \colon \gR \to \gS$ consists of the following data: 
\begin{enumerate}[(i)]
	\item a diagram indexed by the poset $\bN \cup \{-1\}$
	\[\begin{tikzcd}
	0_*(\gR) = \mr{sk}_{-1}(f) \rar & \mr{sk}_0(f) \rar & \mr{sk}_{1}(f) \rar &\cdots \end{tikzcd}\]
	in $\Alg_{\cO}(\sfC^{\bZ_{\leq}}),$ 
	\item  for $d \geq 0$, the data of a filtered CW attachment of dimension $d$ on the morphism $f_d \colon \mr{sk}_{d-1}(f) \to \mr{sk}_d(f)$, in particular a pushout diagram
	\begin{equation}\label{eqn:skeletal-step}
\begin{tikzcd}
F^\cO\left(\bigsqcup_{\alpha \in I_d} \partial D^{g_\alpha,d}_\alpha[d-1] \right) \rar \dar & \mr{sk}_{d-1}(f) \dar{f} \\
	F^\cO\left(\bigsqcup_{\alpha \in I_d} D^{g_\alpha,d}_\alpha[d]\right) \rar & \mr{sk}_d(f), \end{tikzcd}
\end{equation}
	\item using the notation $\gls{skeletal}(f) \coloneqq \mr{colim}_d\, \mr{sk}_d(f)$, a commutative diagram
	\[\begin{tikzcd} \gR \rar{f} \dar{\cong} & \gS \dar{\cong} \\
	\mr{colim}(\mr{sk}_{-1}(f)) \rar & \mr{colim}(\mr{sk}(f)). \end{tikzcd}\] 
\end{enumerate}
\end{definition}
	
\begin{definition}
  A \emph{CW-algebra structure}\index{CW-!algebra} on $\gR \in \Alg_\cO(\sfC)$ is a relative CW-structure on the initial map $\binit \to \gR$.
\end{definition}

Note that $\mr{sk}(f)(d) \cong \mr{sk}_{-1}(f)(d) \cong \binit$ for all $d<0$, so the filtration $\mr{sk}(f)$ is ascending in the sense of Definition \ref{def:filtered-graded-objects}. We will think of $\mr{sk}(f)$ as the skeletal filtration on $f$. Since the left adjoint $\colim$ commutes with colimits, if $f \colon \gR \to \gS$ admits a relative CW-structure then it is cellular.

When we eventually construct CW-structures inductively, it is helpful to have a notion of map between CW-algebras. A map $\gR \to \gS$ of CW-algebras is a \emph{CW-map}\index{CW-!map} if the CW-structure on $\gR$ may be obtained from $\gS$ by taking a subset of the cells, and $\gR \to \gS$ is induced by the inclusion of these cells upon applying $\colim$.

\subsubsection{The associated graded of the skeletal filtration} 

In this section we describe the associated graded of the skeletal filtration; heuristically passing to the associated graded ``filters away'' the attaching maps for the CW-algebra structure.

\begin{theorem}\label{thm:associated-graded-skeletal}
Using the notation of Definition \ref{def:relative-cw-structure}, there is an isomorphism
	\[\grr(\mr{sk}(f)) \cong 0_*(\gR_+) \vee^{\cO} F^\cO\left(\bigvee_{d \geq 0} \bigvee_{\alpha \in I_d} d_*(S^{g_\alpha, d}_\alpha)\right)\]
	 in $\Alg_\cO(\sfC^{\bZ_{=}}_\ast)$.
\end{theorem}

\begin{proof}
On taking quotients the map of pairs
\[(\gR, \binit) = ( \mr{sk}_{-1}(f)(0), \mr{sk}(f)_{-1}(-1)) \lra (\mr{sk}(f)(0), \mr{sk}(f)(-1))\]
gives a morphism $\gR_+ = \grr(\mr{sk}_{-1}(f))(0) \to  \grr(\mr{sk}(f))(0)$ in $\Alg_\cO(\sfC_\ast)$ and hence by adjunction a morphism $\phi \colon 0_*(\gR_+) \to \grr(\mr{sk}(f))$ in $\Alg_\cO(\sfC^{\bZ_{=}}_\ast)$. For each cell we have a characteristic map 
	\[i_\alpha \colon (D^{g_\alpha,d}_\alpha,\partial D^{g_\alpha,d}_\alpha) \lra (\mr{sk}(f)(d),\mr{sk}(f)(d-1))\]
	 which on quotients gives a pointed morphism $j_\alpha \colon S^{g_\alpha, d}_\alpha \to \grr(\mr{sk}(f))(d)$ and hence by adjunction a morphism $d_*(j_\alpha) \colon d_*(S^{g_\alpha, d}_\alpha) \to \grr(\mr{sk}(f))$. Freely extending this to a map of $\cO$-algebras we obtain a morphism
	\[\varphi \colon 0_*(\gR_+) \vee^{\cO} F^\cO\left(\bigvee_{d \geq 0} \bigvee_{\alpha \in I_d} d_*(S^{g_\alpha, d}_\alpha)\right) \xrightarrow{\phi \vee^{\cO} F^\cO\left(\bigvee_{d \geq 0} \bigvee_{\alpha \in I_d} d_*(j_\alpha)\right)\,} \grr (\mr{sk}(f))\]
	in $\Alg_\cO(\sfC_\ast^{\bZ_{=}})$, which we claim is an isomorphism.
	
	We shall prove by induction over $k$ that
	\begin{equation*}
\varphi_k \colon 0_*(\gR_+) \vee^{\cO} F^\cO\left(\bigvee_{d \leq k} \bigvee_{\alpha \in I_d} d_*(S^{g_\alpha, d}_\alpha)\right) \xrightarrow{\phi \vee^{\cO} F^\cO\left(\bigvee_{d \leq k} \bigvee_{\alpha \in I_d} d_*(j_\alpha)\right)\,} \grr (\mr{sk}_k(f))
\end{equation*}
	is an isomorphism given that $\varphi_{k-1}$ is: as $\grr$ and $F^\cO$ both commute with colimits, it then follows by taking the colimit as $k \to \infty$ that $\varphi$ is an isomorphism. The initial case $k=-1$ is obvious, because in that case both sides are $0_*(\gR_+)$. For the inductive step from $k-1$ to $k$, we apply $\grr$ to the pushout diagram (\ref{eqn:skeletal-step}): as $\grr$ commutes with colimits and with $F^\cO$ this gives a pushout diagram
	\[\begin{tikzcd} F^\cO\left(\bigvee_{\alpha \in I_k} (k-1)_*(\partial D^{g_\alpha,k-1}_\alpha)\right) \dar \rar & \grr(\mr{sk}_{k-1}(f)) \dar \\
		F^\cO\left(\bigvee_{\alpha \in I_k} (k-1)_*(\partial D^{g_\alpha,k-1}_\alpha) \vee \bigvee_{\alpha \in I_k} k_*(S^{g_\alpha,k}_\alpha)\right) \rar & \grr(\mr{sk}_{k}(f)). \end{tikzcd}\]
	Omitting the corner $\grr(\mr{sk}_{k}(f))$, this is a coproduct	in $\Alg_\cO(\sfC_\ast^{\bZ_{=}})$ of the two diagrams
	\begin{equation*}
	\begin{tikzcd} F^\cO\left(\bigvee_{\alpha \in I_k} (k-1)_*(\partial D^{g_\alpha,k-1}_\alpha)\right) \dar \rar &[-5pt] \grr(\mr{sk}_{k-1}(f)) \\
		F^\cO\left(\bigvee_{\alpha \in I_k} (k-1)_*(\partial D^{g_\alpha,k-1}_\alpha)\right) & \qquad \end{tikzcd}\end{equation*}
	and
		\[\begin{tikzcd} \ast \dar \rar &[-10pt] \ast \\
	F^\cO\left(\bigvee_{\alpha \in I_k} k_*(S^{g_\alpha,k}_\alpha)\right). & \qquad \end{tikzcd}\]
	Since the left vertical arrow of the first diagram is an isomorphism, the expression for $\grr(\mr{sk}_k(f))$ simplifies to
	\[\grr(\mr{sk}_k(f)) \cong \grr(\mr{sk}_{k-1}(f)) \vee^\cO F^\cO\left(\bigvee_{\alpha \in I_k} k_*(S^{g_\alpha,k}_\alpha)\right),\]
	under which the map $\varphi_k$ is identified with $\varphi_{k-1} \vee^{\cO} F^\cO\left(\bigvee_{\alpha \in I_k} k_*(S^{g_\alpha,k}_\alpha)\right)$, so is an isomorphism.
\end{proof}


\part{Homotopy theory of algebras over a monad} \label{part:homotopy}

In this second part we will add homotopy-theoretic considerations to the theory developed in Part \ref{part:category}. We shall suppose that $\sfS$ is given a model structure, and from this produce model structures on $\sfC = \sfS^\sfG$ and $\Alg_T(\sfC)$. Using these we can derive many of the constructions made so far, such as change-of-diagram-category, change-of-monads, indecomposables, decomposables, or cell attachments, which we do in Section \ref{sec:derived-functors}. The main construction of interest is the functor of derived $T$-indecomposables, $Q^T_\bL$, whose homology groups are defined in Section \ref{sec:homology} and called $T$-homology. We will often be specifically interested in the case that $T$ is the monad associated to an operad $\cO$.

The main technical tools we will develop are simplicial formulae for computing various derived functors, spectral sequences for computing with filtered objects, a Hurewicz theorem for $\cO$-homology, and a cellular approximation theorem for $\cO$-algebras. These appear in Sections \ref{sec:simpl-resol}, \ref{sec:homotopy-theory-filtered}, \ref{sec:an-estimate-abs}, and \ref{sec:additive-case} respectively.

Throughout this part, we will assume that the axioms of Section \ref{sec:axioms-of-cats} hold for $\sfS$ (and hence for $\sfC = \sfS^\sfG$) unless mentioned otherwise:
\begin{itemize}
	\item Axiom \ref{axiom:simplicial-enriched}: $\cat{S}$ is simplicially enriched.
	\item Axiom \ref{axiom:enriched-complete-cocomplete}: $\cat{S}$ is complete and cocomplete in an enriched sense.
	\item Axiom \ref{axiom.cat-monoidal}: $\cat{S}$ has a simplicially enriched closed $k$-monoidal structure, closed on both sides if $k=1$.
\end{itemize}

\section{Contexts for homotopy theory}

In Section \ref{sec:cat-contexts} we discussed axioms on a category $\sfS$ necessary for a good theory of algebras over a sifted monad. As this part concerns homotopy theory, we will require $\sfS$ to be a model category and discuss model categorical axioms necessary for a good homotopy theory of algebras over a sifted monad.

\subsection{Axioms for convenient contexts} \label{sec:axioms-of-model-cats} Model categories are a convenient setting for homotopy theory and several good references exist, e.g.\ \cite{QuillenHA,GoerssSchemmerhorn,Hirschhorn,Hovey}. A \emph{model category structure}\index{model category} on a complete cocomplete category $\sfS$ consists of three classes of morphisms; \emph{weak equivalences}\index{weak equivalence}, \emph{cofibrations}\index{cofibration}, and \emph{fibrations}\index{fibration}. These classes should be closed under retracts and 2-out-of-3. Morphisms that are both weak equivalences and cofibrations are called \emph{trivial cofibrations} and similarly morphisms that are both weak equivalences and fibrations are called \emph{trivial fibrations}. The trivial cofibrations should have the left lifting property with respect to fibrations and cofibrations should have the left lifting property with respect to trivial fibrations. There should further exist two functorial factorizations of a morphism $f \colon X \to Y$ into $X \to Z \to Y$ with either (i) $X \to Z$ a trivial cofibration and $Y \to Z$ a fibration, or (ii) $X \to Y$ a cofibration and $Y \to Z$ a trivial fibration. Applying this to $\binit \to X$ or $X \to \bterm$ we obtain functorial cofibrant or fibrant replacements. We call a weak equivalence $f \colon Y \to X$ with $Y$ cofibrant a \emph{cofibrant approximation of $X$} and say it is a \emph{cofibrant replacement} if $f$ is a trivial fibration. Similarly, we call a weak equivalence $g \colon X \to Z$ with $Z$ fibrant a \emph{fibrant approximation of $X$} and say it is a \emph{fibrant replacement} if $g$ is a trivial cofibration. A model category $\sfS$ has a homotopy category $\cat{Ho}(\sfS)$, which can be constructed as a localization but also has a concrete model as the category with objects the cofibrant-fibrant objects in $\sfS$ and morphisms the homotopy classes of morphisms in $\sfS$.

\newglossaryentry{bl}{%
	name={\ensuremath{\bL}},
	description={Left-derived functor},
type=symbols
}
\newglossaryentry{br}{%
	name={\ensuremath{\bR}},
	description={Right-derived functor},
type=symbols
}
In an adjunction between model categories, the left adjoint $F$ is a \emph{left Quillen functor} and the right adjoint $G$ a \emph{right Quillen functor}\index{Quillen!functor} if one of the following equivalent conditions hold: (i) $F$ preserves cofibrations and trivial cofibrations, (ii) $G$ preserves fibrations and trivial fibrations, (iii) $F$ preserves cofibrations and $G$ fibrations, (iv) $F$ preserves trivial cofibrations and $G$ trivial fibrations. Left, resp.\ right, Quillen functors may be derived by composing them with a functorial cofibrant, resp.\ fibrant, replacement functor, which we shall denote $c$, resp.\ $f$. We use the notation $\gls{bl} F \coloneqq L \circ c$ and $\gls{br} G \coloneqq G \circ f$ for these derived functors\index{derived functor}: they come with a natural transformation $\bL F \Rightarrow F$ (which is a weak equivalence on cofibrant objects) and $G \Rightarrow \bR G$ (which is a weak equivalence on fibrant objects). The functor $\bL F$ preserves weak equivalences between cofibrant objects and $\bR G$ preserves weak equivalences between fibrant objects, by Ken Brown's Lemma. These derived functors induce a pair of adjoint functors on the homotopy categories, which is independent of the choice of $c$ or $f$, and will also be denoted $\bL F \dashv \bR G$. It is possible to derive functors with weaker properties: a functor $F$ may be left derived if it takes trivial cofibrations between cofibrant objects to weak equivalences, and a functor $G$ may be right derived if it takes trivial fibrations between fibrant objects to weak equivalences, see Proposition 8.4.8 of \cite{Hirschhorn}.

For categories of algebras or diagram categories, we define so-called projective model structures. The construction of these uses that the model category structure on $\sfS$ is \emph{cofibrantly generated}\index{model category!cofibrantly generated}, see Section 2.1 of \cite{Hovey}. Roughly, this means that the cofibrations and trivial cofibrations are generated by \emph{sets} of morphisms under pushouts, transfinite composition and retracts. These sets are called \emph{generating cofibrations} and \emph{generating trivial cofibrations}. 

\begin{axiom}\label{axiom.CofGen}$\sfS$ is equipped with a cofibrantly generated model category structure.\end{axiom}
	
As the category $\sfS$ is required to satisfy the axioms in Section \ref{sec:cat-contexts}, it comes with a simplicial enrichment and a closed $k$-monoidal structure (using the latter, the former can be reconstructed from a functor $s \colon \cat{sSet} \to \sfC$). We will require the model category structure to be compatible with these two structures, in the following precise sense (see Section 4.2 of \cite{Hovey}).

A functor $F \colon \sfC \times \sfD \to \cat{E}$ is called a \emph{Quillen bifunctor}\index{Quillen!bifunctor} if it is an adjunction of two variables and has the following two properties: firstly, if $i \colon K \to L$ and $j \colon X \to Y$ are cofibrations then the pushout $F(K,Y) \sqcup_{F(K,X)} F(L,X) \to F(L,Y)$ is a cofibration. Secondly, this should be a weak equivalence if at least one of $i$, $j$ is a weak equivalence. If $F$ is a Quillen bifunctor, fixing cofibrant objects $C \in \sfC$ or $D \in \sfD$ we obtain a left Quillen functor $F(C,-) \colon \sfD \to \cat{E}$ or $F(-,D) \colon \sfC \to \cat{E}$.

A model category is \emph{monoidal}\index{model category!monoidal} if $\otimes \colon \sfS \times \sfS \to \sfS$ is a Quillen bifunctor and tensoring a cofibrant $X$ with a cofibrant replacement $c(\bunit) \to \bunit$ gives a weak equivalence. Note that the second condition is automatically satisfied when $\bunit$ is cofibrant, which we shall soon assume in Axiom \ref{axiom.model-monoidalsimplicial}. It has the following consequences: (i) if $X$, $Y$ are cofibrant, then so is $X \otimes Y$, (ii) if $X \to X'$ is a weak equivalence between cofibrant objects and $Y$ is cofibrant, then $X \otimes Y \to X' \otimes Y$ is a weak equivalence (and similarly in the other variable).

A left Quillen functor $F$ between monoidal model categories is strong/oplax/lax monoidal if it is a strong/oplax/lax monoidal functor and $F(c(\bunit)) \to F(\bunit)$ is a weak equivalence. If $\bunit$ is cofibrant the latter is automatic.

The model category $\sfS$ is \emph{simplicial}\index{model category!simplicial} if the copowering $\times \colon \cat{sSet} \times \sfS \to \sfS$ is a Quillen bifunctor. It has the following consequences: (i) if $X$ is cofibrant, then so is $K \times X$, (ii) if $K \to K'$ is a weak equivalence and $X$ cofibrant, then $K \times X \to K' \times X$ is a weak equivalence, (iii) if $X \to X'$ is a weak equivalence between cofibrant objects, then $K \times X \to K' \times X$ is a weak equivalence. Using $s(K) = K \times \bunit_\sfS$ and $K \times X = s(K) \otimes X$, it follows that if $\sfS$ satisfies the axioms in Section \ref{sec:axioms-of-cats} and is a monoidal model category with $\bunit$ cofibrant, then $\sfS$ is simplicial if and only if $s \colon \cat{sSets} \to \sfS$ is a left Quillen functor if and only if $\mr{Sing} \colon \sfS \to \cat{sSets}$ is a right Quillen functor. 

Though we believe that with appropriate modifications one can state all our results when the unit is not cofibrant (e.g.~\cite[Theorem 1]{MuroUnit}), it will be quite convenient to make this assumption.

\begin{axiom}\label{axiom.model-monoidalsimplicial} The model category structure on $\sfS$ is monoidal and simplicial. The unit $\bunit$ of the monoidal structure is cofibrant.\end{axiom}

We do not demand any compatibility between fibrant replacements and the monoidal structure. Sometimes there is, through the existence of a lax $k$-monoidal fibrant approximation (which need \emph{not} be a functorial fibrant replacement):

\begin{definition}\label{def.model-laxfibrant} A \emph{lax $k$-monoidal fibrant approximation} on $\sfC$ is a functor $R \colon \sfC \to \sfC$ such that $R(X)$ is fibrant for all $X$ which comes with a lax $k$-monoidality and a lax $k$-monoidal natural weak equivalence $X \to R(X)$.\end{definition}

Finally, it shall be useful to know that homotopy equivalences in the following sense are weak equivalences. This is Proposition 9.5.16 of \cite{Hirschhorn}.

\begin{definition}\label{def:homotopy-equiv} Given two maps $f_0,f_1 \colon X \to Y$, a \emph{homotopy from $f_0$ to $f_1$} is a map $H \colon \Delta^1 \times X \to Y$ such that the composite $H \circ i_0 \colon \{0\} \times X \to \Delta^1 \times X \to Y$ equals $f_0$ and the composite $H \circ i_1 \colon \{1\} \times X \to \Delta^1 \times X \to Y$ equals $f_1$.

A map $f \colon X \to Y$ is a \emph{homotopy equivalence} if there exists a map $g \colon Y \to X$ such that $f \circ g$ is homotopic to $\mr{id}_Y$ and $g \circ f$ is homotopic to $\mr{id}_X$.\end{definition}

\begin{remark}It may be helpful to note which axioms in this section we consider necessary, and which we consider merely helpful (and likely avoidable with more technical work). The ``helpful'' assumption is the part of Axiom \ref{axiom.model-monoidalsimplicial} stating that the monoidal unit $\bunit$ is cofibrant. The remaining axioms seem to be ``necessary.''\end{remark}
	
\subsection{Examples}\label{sec:examples-of-model-cats}
We explain how to endow the following examples of $\sfS$ with model structures satisfying the axioms of Section \ref{sec:axioms-of-model-cats}.

\subsubsection{Simplicial sets} 

The Quillen model structure has cofibrations the monomorphisms and fibrations the Kan fibrations, and these determine the remaining data. The weak equivalences are those maps that induce weak equivalences of topological spaces upon geometric realization. It is cofibrantly generated by generating cofibrations $\{\partial \Delta^n \hookrightarrow \Delta^n\}$ and generating trivial cofibrations $\{\Lambda^n_i \hookrightarrow \Delta^n\}$ \cite[Chapter 3.2]{Hovey}. Thus the Quillen model structure satisfies Axiom \ref{axiom.CofGen}.

For Axiom \ref{axiom.model-monoidalsimplicial}, the Quillen model structure is monoidal for the Cartesian product $\times$ \cite[Proposition 4.2.8]{Hovey}. Clearly $\mr{id} \colon \cat{sSet} \to \cat{sSet}$ is a left Quillen functor, so it is also simplicial. It has a lax symmetric monoidal fibrant replacement functor given by $\mr{Ex}^\infty$, which preserves finite products as it is a filtered colimit of functors with a left adjoint.

\subsubsection{CGWH topological spaces} The Serre model structure on CGWH topological spaces has fibrations the Serre fibrations and weak equivalences the weak homotopy equivalences, and these determine the remaining data. It is cofibrantly generated by the generating cofibrations $\{S^{n-1} \hookrightarrow D^n\}$ and generating trivial cofibrations $\{D^n \times \{0\} \hookrightarrow D^n \times [0,1]\}$ \cite[Chapter 2.4]{Hovey}.

For Axiom \ref{axiom.model-monoidalsimplicial}, it is monoidal for the Cartesian product $\times$ \cite[Proposition 4.2.11]{Hovey}. The functor $\gr{-} \colon \cat{sSet} \to \cat{Top}$ is a left Quillen functor, as it sends a monomorphism of simplicial sets to a relative CW-complex and creates weak equivalence by definition. As all spaces are fibrant, the identity is a lax (in fact strong) symmetric monoidal fibrant replacement.

\subsubsection{Simplicial $\bk$-modules} \label{sec.model-simplicial-rmods} Let $\cat{sC}$ denote the category of simplicial objects in  $\sfC$, and suppose there is a functor $F \colon \cat{sSet} \to \cat{sC}$ with right adjoint $U \colon \cat{sC} \to \cat{sSet}$. There exists a simplicial model structure on $\cat{sC}$ in which $f$ is a fibration or weak equivalence if $Uf$ is, if the following conditions are satisfied \cite[Theorem II.5.4 and Lemma II.6.1]{GoerssJardine}: (i) $U$ commutes with filtered colimits and (ii) every object of $\cat{sC}$ has fibrant underlying simplicial set. This will then be cofibrantly generated by generating cofibrations $\{F(\partial \Delta^n) \to F(\Delta^n)\}$ and generating trivial cofibrations $\{F(\Lambda^n_i) \to F(\Delta^n)\}$. 

This result in particular applies to the category of simplicial $\bk$-modules $\cat{sMod}_\bk$, with $F \colon \cat{sSet} \to \cat{sMod}_\bk$ given by taking the levelwise free $\bk$-module, as $U$ preserves filtered colimits and every simplicial abelian group is Kan. Because every object of $\cat{sMod}_\bk$ is fibrant, the identity is a lax symmetric monoidal fibrant replacement functor. For Axiom \ref{axiom.model-monoidalsimplicial}, \cite[Corollary 4.2.5]{Hovey} implies it suffices to check the properties of a Quillen bifunctor for $\otimes$ only on generating (trivial) cofibrations. This is clear, as all maps involved can be expressed as $F$ applied to colimits of generating (trivial) cofibrations.

\subsubsection{Symmetric spectra} There is a number of model structures on the category of symmetric spectra, all of which are cofibrantly generated, monoidal and simplicial \cite{Schwedebook}, verifying Axioms \ref{axiom.model-monoidalsimplicial} and \ref{axiom.CofGen}. We shall take the \emph{absolute projective stable model structure} (called the stable model structure in \cite{HoveyShipleySmith}).

In the absolute projective stable model structure, the cofibrations are those maps that have the left lifting property with respect to the level trivial fibrations, i.e.\ the maps $f \colon X \to Y$ such that $f_n \colon X_n \to Y_n$ is a Kan fibration of underlying simplicial sets. The weak equivalences are the stable equivalences; these are the map $f \colon X \to Y$ such that the induced map $f^* \colon [Y,A] \to [X,A]$ is a bijection for all injective $\Omega$-spectra $A$. It is worth pointing out that no lax fibrant approximation functor can exist on symmetric spectra, due to Lewis's argument \cite{LewisSpectra}.

\begin{remark}An alternative model structure on symmetric spectra is the \emph{positive projective stable model structure}. To obtain this, only use the generating (trivial) cofibrations for $n > 0$. This has advantages when dealing with operads that are not $\Sigma$-cofibrant, see e.g.\ \cite{HarperSymSpec,HarperHess,PavlovScholbach3}, but $\bS$ is no longer cofibrant.\end{remark}

\subsubsection{Pointed categories}\label{sec:model-pointed-categories} If $\sfS$ satisfies the axioms of Section \ref{sec:axioms-of-model-cats}, then $\sfS_*$ does too. Firstly, if $\sfS$ is a cofibrantly generated model category it follows from \cite{HirschhornUnder} that $\sfS_*$ has the structure of a cofibrantly generated model category, where weak equivalences, cofibrations, and fibrations are all created by $U^+ \colon \sfS_* \to \sfS$. The generating (trivial) cofibrations are $f \sqcup \bterm$ where $f \colon A \to B$ is a generating (trivial) cofibration of $\sfS$.

We equip $\sfS_*$ with the monoidal structure $\owedge$ described in Section \ref{sec:basepoint-monad}. The unit is $\bunit_\sfS \sqcup \bterm$, so is cofibrant in $\sfS_*$ as $\bunit_\sfC$ is cofibrant in $\sfS$. As $\owedge$ participates in an adjunction of two variables, to prove that $\owedge$ is a Quillen bifunctor it suffices to verify the condition on generating (trivial) cofibrations by Corollary 4.2.5 of \cite{Hovey}. As $F^+(X) \owedge F^+(Y) \cong F^+(X \otimes Y)$, the pushout-product of generating cofibrations $F^+(f)$ and $F^+(g)$ is obtained by applying $F^+$ to the $\otimes$ pushout-product of $f$ and $g$ in $\sfS$, so is a cofibration, and is a trivial cofibration if one of $f$ or $g$ is. This verifies Axiom \ref{axiom.model-monoidalsimplicial}. That $\times \colon \cat{sSet} \times \sfS_\ast \to \sfS_\ast$ is a Quillen bifunctor may similarly be verified on generating (trivial) cofibrations using $K \times F^+(X) \cong F^+(K \times X)$.

If $\sfS$ is pointed, the fact that $\times \colon \cat{sSet} \times \sfS_\ast \to \sfS_\ast$ is a Quillen bifunctor implies that $\wedge \colon \cat{sSet}_\ast \times \sfS_\ast \to \sfS_\ast$ is also a Quillen bifunctor. We again use Corollary 4.2.5 of \cite{Hovey} and that $F^+(K) \wedge F^+(X) \cong F^+(K \times X)$.

\subsection{Model category structures transferred along adjunctions} \label{sec.model-adjunctions} The construction of a model category structure on $\cat{sMod}_\bk$ and $\sfS_\ast$ are examples of the transfer of a model structure along an adjunction. In this section we explain this technique, and apply it to diagram categories and module categories.

\subsubsection{The projective model structure}\label{sec:projective-model-structure-general} 
Given a model category $\sfC$ and a functor $F \colon \sfC \to \sfD$ with right adjoint $U$, one may try to construct a model structure on $\sfD$ by declaring a morphism $f$ in $\sfD$ to be a fibration or weak equivalence if $Uf$ is in $\sfC$. If it exists, this is called the \emph{projective model structure}\index{model structure!projective} (also known as the \emph{right-induced model structure}), and with respect to this model structure $F$ becomes a left Quillen functor since $U$ preserves fibrations and trivial fibrations. If $\sfC$ is cofibrantly generated by sets $I$ and $J$ of generating cofibrations and trivial cofibrations, then this is equivalent to declaring $FI \coloneqq \{F(i) \mid i \in I\}$ to be the set of generating cofibrations and $FJ \coloneqq \{F(j) \mid j \in J\}$ to be the set of generating trivial cofibrations. There are Theorem 11.1.13 of \cite{FresseBook} or Theorem 11.3.2 of \cite{Hirschhorn} give general conditions under which this indeed defines a cofibrantly generated model structure (these results essentially go back to Quillen \cite{QuillenHA}).

\begin{theorem}\label{thm:general-proj-model-structure} The projective model structure on $\sfD$ transferred along the adjunction $F \dashv U$ exists if 
	\begin{enumerate}[(i)]
		\item \label{enum:proj-i} $FI$ and $FJ$ admit a small object argument.
		\item \label{enum:proj-ii} $U$ takes relative $FJ$-cell complexes to weak equivalences.
	\end{enumerate}
\end{theorem}

In condition \eqref{enum:proj-i}, that $FI$ admits a small object argument means that the domains of morphisms of the set $FI$ are small relative to the class of relative $FI$-cell complexes \cite[Definition 10.5.15]{Hirschhorn}. The following observation appears in Section 2.5 of \cite{BergerMoerdijk}:

\begin{lemma}\label{lem:small-object-condition} Theorem \ref{thm:general-proj-model-structure} \eqref{enum:proj-i} holds if
	\begin{enumerate}[(i$^\prime$)]
		\item \label{enum:proj-iprime} The domains of $I$ and $J$ are small relative to all morphisms, and $U$ preserves filtered colimits.
	\end{enumerate}
\end{lemma}

Theorem 3.8 of \cite{GoerssSchemmerhorn} gives a way to verify \eqref{enum:proj-ii} (see also Section 2.6 of \cite{BergerMoerdijk}):

\begin{proposition}\label{prop:proj-gs} Theorem \ref{thm:general-proj-model-structure} \eqref{enum:proj-ii} holds if
\begin{enumerate}[(i$^\prime$)] \setcounter{enumi}{1}
	\item \label{enum:proj-iiprime} $\sfD$ is simplicial (to form natural path objects) and there exists an \emph{fibrant approximation}, i.e.\ a functor $R \colon \sfD \to \sfD$ such that $R(X)$ is fibrant for all $X$ (i.e.\ fibrant upon applying $U$) and a natural transformation $\mr{id}_\sfD \to R$ that is a natural weak equivalence (i.e.\ a natural transformation that becomes a weak equivalence upon applying $U$).
\end{enumerate}
\end{proposition}

\subsubsection{Diagram categories}\label{sec:FunctorCategories}

These techniques can be applied to diagram categories. If $\sfG$ is a small category, then the inclusion $\mr{ob}(\sfG) \to \sfG$ with $\mr{ob}(\sfG)$ the set of objects in $\sfG$ defines a functor $U \colon \sfS^{\sfG} \to \prod_{\mr{ob}(\sfG)} \sfS$ by restriction, with left adjoint $F$ given by \[F((X_g)_{g \in \mr{ob}(\sfS)}) = \bigsqcup_{g \in \mr{ob}(\sfG)} \sfG(g,-) \times X_g.\]
In \cite[Theorem 11.6.1]{Hirschhorn} the assumptions of Theorem \ref{thm:general-proj-model-structure} are verified for $F \dashv U$ assuming that $\sfS$ is cofibrantly generated. Thus assuming Axiom \ref{axiom.CofGen}, the projective model structure on $\sfS^\sfG$ exists. By \cite[Theorem 11.7.3]{Hirschhorn} it will be a simplicial model structure if $\sfS$ is a simplicial cofibrantly generated model category.

If the category $\sfG$ is in addition monoidal then $\sfS^\sfG$ is equipped with a Day convolution product, and we will now verify that this makes it a monoidal model category \cite[Section 2.2]{IsaacsonThesis}.

\begin{lemma}\label{lem:model-structure-functor-cats}
If $\sfG$ is a small closed monoidal category and $\sfS$ satisfies the axioms of Section \ref{sec:axioms-of-model-cats}, then $\sfS^\sfG$ with the projective model structure also satisfies the axioms of Section \ref{sec:axioms-of-model-cats}. If $\sfS$ has a lax monoidal fibrant approximation functor, then so does $\sfS^\sfG$.
\end{lemma}

\begin{proof}
We have explained above that $\sfS^\sfG$ has the projective model structure, which is by definition cofibrantly generated, and that it is simplicial, so to verify Axiom \ref{axiom.model-monoidalsimplicial} it remains to show that it is monoidal with respect to Day convolution and that the unit is cofibrant.

The Day convolution product is closed, so participates in an adjunction of two variables and hence by \cite[Corollary 4.2.5]{Hovey} it suffices to check that Day convolution is a Quillen bifunctor only on generating (trivial) cofibrations. Since $F$ is a left adjoint, and $\otimes$ distributes over colimits by closedness, it suffices to check this for morphisms of the form $\sfG(g,-) \times i$ or $\sfG(h,-) \times j$, where $i \in I$ is a generating cofibration and $j \in J$ is a generating cofibration. For example, in the case $\sfG(g,-) \times i$, $\sfG(g',-) \times i'$ for $i,i' \in I$, the pushout-product
\[\sfG(g \oplus g',-) \times ((X \otimes Y') \sqcup_{X' \otimes Y'} (X' \otimes Y)) \lra \sfG(g \oplus g',-) \times (X' \otimes Y')\]
is a cofibration, as it is obtained by applying the left Quillen functor $(g \oplus g')_* \colon \sfS \to \sfS^\sfG$ to a cofibration. Furthermore, if one of $i$ or $i'$ is a trivial cofibration then it is obtained by applying $(g \oplus g')_* \colon \sfS \to \sfS^\sfG$ to a trivial cofibration, so is a trivial cofibration. The unit of $\sfS^\sfG$ is $\bunit_{\sfS^\sfG} = \sfG(\bunit_\sfG,-) \times \bunit_\sfS$, obtained by applying
\[\sfS \xrightarrow{(\bunit_\sfG)_*} \prod_{\mr{ob}(\sfG)} \sfS \overset{F}\lra \sfS^\sfG\]
to $\bunit_\sfS$. These are both left Quillen functors, so as $\bunit_\sfS$ is cofibrant so is $\bunit_{\sfS^\sfG}$.

Any fibrant approximation functor $f \colon \sfS \to \sfS$ induces a functor $f^\sfG \colon \sfS^\sfG \to \sfS^\sfG$ by applying $f$ objectwise. This will be lax monoidal with respect to Day convolution if $f$ is lax monoidal.\end{proof}

\begin{example}In Example \ref{exam:change-groupoid-indec} we explained that the functor $p^* \colon \sfS^{\sfG'} \to \sfS^\sfG$ induced by precomposition by functor $p \colon \sfG \to \sfG'$ has a left adjoint $p_* \colon \sfS^{\sfG} \to \sfS^{\sfG'}$ given by left Kan extension. This is a Quillen adjunction, since it is clear that $p^*$ preserves (trivial) fibrations. Thus we may derive either of these functors.\end{example}

\begin{remark}\label{rem:hocolim}
At this point we can also mention homotopy colimits\index{homotopy!colimit} (see \cite{DwyerHirschhornKanSmith} for more background). If $\cat{I}$ is a small category and $\sfC$ is a cofibrantly generated model category, then as described above the projective model structure on $\sfC^\cat{I}$ exists. The adjunction
\[\begin{tikzcd}
\sfC^\cat{I} \arrow[shift left=.5ex]{r}{\colim} & \arrow[shift left=.5ex]{l}{\const} \sfC
\end{tikzcd}\]
is a Quillen adjunction, as $\const$ clearly preserves fibrations and trivial fibrations by definition of the projective model structure, and the homotopy colimit functor $\hocolim$ is the derived functor $\bL \colim$.

If $\sfD$ is another such model category and $F \colon \sfC \to \sfD$ is a functor which takes weak equivalences between cofibrant objects to weak equivalences then it has a left derived functor $\bL F \colon \cat{Ho}(\sfC) \to \cat{Ho}(\sfD)$. Furthermore $F^\cat{I} \colon \sfC^\cat{I} \to \sfD^\cat{I}$ also preserves weak equivalences between cofibrant objects, because weak equivalences are objectwise and cofibrant objects of $\sfC^\cat{I}$ are in particular objectwise cofibrant, by e.g.\  \cite[Proposition 11.6.3]{Hirschhorn}, so $F^\cat{I}$ also has a left derived functor. Hence we may ask if
\[\begin{tikzcd}
\cat{Ho}(\sfC^\cat{I}) \arrow{r}{\bL\colim} \ar{d}[swap]{\bL F^\cat{I}} &[5pt] \cat{Ho}(\sfC) \ar{d}{\bL F}\\
\cat{Ho}(\sfD^\cat{I}) \arrow{r}{\bL\colim} & \cat{Ho}(\sfD)
\end{tikzcd}\]
commutes (up to natural isomorphism). If it does then we say that $F$ \emph{preserves homotopy colimits}. A sufficient condition is that $F$ is a left Quillen functor, but this is not a necessary condition.
\end{remark}

\subsubsection{Module categories}\label{sec:model-module-categories}

The second application of Section \ref{sec:projective-model-structure-general} is to the category $\gR\text{-}\cat{Mod}$ of modules over a commutative algebra object $\gR$ in $\sfS$. We shall give two conditions under which the projective model structure transferred along the free-forgetful adjunction exists. The second of these involves the \emph{monoid axiom} \cite[Definition 3.3]{SchwedeShipley}: if a morphism is a transfinite composition of pushouts of tensor products of trivial cofibrations with any object, it is a weak equivalence. 

\begin{theorem}\label{thm:model-structure-on-modules}
	Suppose that $\sfS$ satisfies the axioms of Section \ref{sec:axioms-of-model-cats} and $\gR$ is a commutative algebra in $\sfS$ such that either
	\begin{enumerate}[(i)]
		\item the underlying object of $\gR$ is cofibrant in $\sfS$, or
		\item $\sfS$ satisfies Schwede--Shipley's monoid axiom,
	\end{enumerate}
	then the category $\gR\text{-}\cat{Mod} \coloneqq \Alg_{\gR \otimes -}(\sfC)$ of left $\gR$-modules has a projective model structure. In case (i), $U^\gR \colon \gR\text{-}\cat{Mod} \to \sfC$ preserves (trivial) cofibrations.
\end{theorem}

\begin{proof}In case (i), we use \cite[Proposition 11.1.4]{FresseBook}, which says that sufficient conditions to apply Theorem \ref{thm:general-proj-model-structure} are that (a) $U^\gR$ preserves sifted colimits, (b) in each pushout diagram
	\begin{equation}\label{eqn:module-pushout} \begin{tikzcd} F^{\gR}(Y_0) \rar \dar[swap]{F^\cO(i)} & \gM_0 \dar{f} \\
	F^\gR(Y_1) \rar & \gM_1, \end{tikzcd}\end{equation}
	the morphism $U^\gR(f)$ is a (trivial) cofibration if $i$ is.
	
	Part (a) follows from Lemma \ref{lem:ut-preserves-sifted}. For part (b), we note that $\gR \otimes -$ preserves all colimits as a left adjoint, so that $U^\gR$ preserves all colimits as well by an argument analogous to Lemma \ref{lem:ut-preserves-sifted}; if $i \mapsto \gM_i$ is a diagram of $\gR$-modules with underlying objects $M_i$, we can endow $\colim_{i \in I} M_i$ with a $\gR$-module structure satisfying the universal property. Thus it suffices to consider the following pushout diagram in $\sfC$
	\[\begin{tikzcd} R \otimes Y_0 \rar \dar[swap]{R \otimes i} & M_0 \dar{f} \\
	R \otimes Y_1 \rar & M_1. \end{tikzcd}\]
	Under assumptions of case (i), $R \otimes -$ preserves (trivial) cofibrations by the axioms of a monoidal model category. As trivial cofibrations these are preserved by pushouts, (b) is satisfied. Since every (trivial) cofibration $g$ in $\gR\text{-}\cat{Mod}$ is a retract of transfinite compositions of morphisms of the form of the right hand vertical map (\ref{eqn:module-pushout}), whose underlying morphisms is also a (trivial) cofibration $\sfC$, $U^\gR(g)$ is also a (trivial) cofibration.
	
	For case (ii), we refer to \cite[Theorem 4.1]{SchwedeShipley} (note Remark 2.4 weakens the smallness assumptions). 
\end{proof}

\begin{remark} In case (i), we may replace the assumption that $R$ is cofibrant with the assumption that $\bunit \to R$ is a cofibration. One follows the proof given above, but to show that $R \otimes i$ is a (trivial) cofibration one applies the pushout-product axiom to the cofibration $\bunit \to R$ and the (trivial) cofibration $Y_0 \to Y_1$ to see that the map $g$ in
	\[\begin{tikzcd} \bunit \otimes Y_0 \cong Y_0 \rar \dar & R \otimes Y_0 \dar{f} \rar[equal] & R \otimes Y_0 \dar{R \otimes i} \\
		\bunit \otimes Y_1 \cong Y_1 \rar & Y_1 \sqcup_{Y_0} (R \otimes Y_1) \rar{g} & R \otimes Y_1 \end{tikzcd}\]
	is a (trivial) cofibration. As (trivial) cofibrations are closed under pushouts, $f$ is a (trivial) cofibration and hence so is $R \otimes i$.
	
	In case (ii), $U^\gR \colon \gR\text{-}\cat{Mod} \to \sfC$ might not preserve (trivial) cofibrations; if $\gR$ is not cofibrant but the monoidal unit $\bunit$ is, then $\binit \to \gR$ is a cofibration in $\gR\text{-}\cat{Mod}$ whose underlying morphism in $\sfC$ is not a cofibration.\end{remark}

For our applications, we need to verify some additional axioms on the model structure on $\gR\text{-}\cat{Mod}$, listed in Section \ref{sec:axioms-of-model-cats}:

\begin{lemma}\label{lem:ModuleCategories}
Suppose that $\sfS$ satisfies the axioms of Section \ref{sec:axioms-of-model-cats} and $\gR$ is a commutative algebra in $\sfS$ such that either
\begin{enumerate}[(i)]
\item the underlying object of $\gR$ is cofibrant in $\sfS$, or
\item $\sfS$ satisfies Schwede--Shipley's monoid axiom,
\end{enumerate}
then $\gR\text{-}\cat{Mod}$ also satisfies the axioms of Section \ref{sec:axioms-of-model-cats}.
\end{lemma}

\begin{proof}First, recall that in Proposition \ref{prop:ModuleCatAxioms} we showed that $\gR\text{-}\cat{Mod}$ satisfies the axioms of Section \ref{sec:axioms-of-cats}. We start with case (i), and the existence of the projective model structure was established in Theorem \ref{thm:model-structure-on-modules}. However, in Section \ref{sec:axioms-of-model-cats} we have imposed additional axioms beyond the mere existence of a model category structure, which we shall now verify for the category $\gR\text{-}\cat{Mod}$. 

Axiom \ref{axiom.CofGen} holds as $\gR\text{-}\cat{Mod}$ is by definition cofibrantly generated. Axiom \ref{axiom.model-monoidalsimplicial} first requires $\gR\text{-}\cat{Mod}$ to be a simplicial model category, which means that
\begin{align*}
- \times - \colon \cat{sSet} \times \gR\text{-}\cat{Mod} &\lra \gR\text{-}\cat{Mod}\\
 (K, \gM) &\longmapsto (\gR \otimes s(K)) \otimes_\gR \gM
\end{align*}
should be a Quillen bifunctor, which can be identified with $(K, \gM) \mapsto K \times \gM \coloneqq s(K) \otimes \gM$. Since $- \times -$ participates in an adjunction of two variables, by \cite[Corollary 4.2.5]{Hovey} it suffices to check that it is a Quillen bifunctor on generating (trivial) cofibrations. So, given generating cofibrations $f \colon K \to L \in \cat{sSet}$ and $g \colon A \to B \in \sfS$ the pushout-product $f \square (\gR \otimes g)$ may be identified with
\[\gR \otimes (f \square g) \colon \gR \otimes (L \times A \sqcup_{K \times A} K \times B) \lra \gR \otimes (L \times B)\]
which is $\gR \otimes -$ applied to a cofibration in $\sfS$ (as $\times \colon \cat{sSet} \times \sfS \to \sfS$ is a Quillen bifunctor), so is a cofibration; if $f$ or $g$ is a trivial cofibration then so $f \square g$ and because $\gR \otimes -$ is a left Quillen functor, so is $\gR \otimes (f \square g)$.

Axiom \ref{axiom.model-monoidalsimplicial} next requires $\gR\text{-}\cat{Mod}$ to be a monoidal model category, and as we have said we will show in Section \ref{sec:AssocModules} that there is a tensor product $\otimes_\gR$, which is a Quillen bifunctor by Lemma \ref{lem:otimes-gr-quillen-bifunctor}. As $\bunit_\sfS$ is cofibrant in $\sfS$ by assumption and $\gR \otimes -$ is a left Quillen functor, $\bunit_{\gR\text{-}\cat{Mod}} = \gR \otimes \bunit_\sfS$ is cofibrant in $\gR\text{-}\cat{Mod}$.

In case (ii), \cite[Theorem 4.1(2)]{SchwedeShipley} shows that the projective model structure on $\gR\text{-}\cat{Mod}$ exists, and that it is monoidal. Verifying the remaining axioms may be done as above.
\end{proof}

\begin{example}\label{exam:hk-cofibrant} 
The category $\sfS = \cat{Sp}^\Sigma$ with the absolute projective stable model structure satisfies Schwede--Shipley's monoid axiom \cite[Section 5]{SchwedeShipley}. The Eilenberg--Maclane object $H\bk$ has $(H\bk)_k$ given by the underlying simplicial set of $\bk[S^k] \in \cat{sMod}_\bk$, and is a commutative ring spectrum. The category of $H\bk$-module spectra is Quillen equivalent to chain complexes over $\bk$ by \cite{Shipley} (see also \cite{ShipleyCorrectionNocorrection}).
\end{example}

\section{Homotopy theory of monads and indecomposables} \label{sec:homotopy-theory-monads}

In this section we discuss the interaction of homotopy theory with the theory of monads of Section \ref{sec:sifted-monads}, and thus we assume that the axioms of Section \ref{sec:axioms-of-cats} and \ref{sec:axioms-of-model-cats} are satisfied. We will explain how to obtain a model structure on the category of algebras over a monad, and hence how to derive the functors defined in that section, most importantly indecomposables of an augmented monad. We also give simplicial formulae to make these derived functors computable.

\subsection{Homotopy theory of algebras over a monad}
\label{sec:homotopy-theory-algebras}

Fix a sifted and simplicial monad $T$. If it exists, the {projective model structure}\index{model structure!projective} on $\Alg_T(\sfC)$ is obtained by transferring the model structure on $\sfC$ along the adjunction 
\begin{equation}\label{eqn:ftut-adj}\begin{tikzcd} \sfC \arrow[shift left=.5ex]{r}{F^T} & \Alg_T(\sfC). \arrow[shift left=.5ex]{l}{U^T}\end{tikzcd}\end{equation}
Section \ref{sec.model-adjunctions} explains conditions which guarantee the existence of such a transferred model structure, and we make the existence of the projective model structure an assumption on $\sfC$ and $T$.

\begin{axiom}\label{axiom:monad-proj} The projective model structure on $\Alg_T(\sfC)$ exists and the forgetful functor $U^T \colon \Alg_T(\sfC) \to \sfC$ preserves (trivial) cofibrations between cofibrant objects.\end{axiom}

\begin{lemma}\label{lem:monad-cofibrant} Assuming Axiom \ref{axiom:monad-proj}, the monad $T$ preserves (trivial) cofibrations between cofibrant objects and the projective model structure on $\Alg_T(\sfC)$ is simplicial.\end{lemma}

\begin{proof}By definition of the projective model structure, if it exists, the adjunction (\ref{eqn:ftut-adj}) is a Quillen adjunction. Thus $F^T(f)$ is a (trivial) cofibration between cofibrant objects if $f$ is, and applying $U^T$ and using the second part of Axiom \ref{axiom:monad-proj} we see that $Tf$ is a (trivial) cofibration between cofibrant objects.

	As $\times \colon \cat{sSet} \times \Alg_T(\sfC) \to \Alg_T(\sfC)$ participates in an adjunction of two variables, that it is a Quillen bifunctor may be verified on generating (trivial) cofibrations by \cite[Corollary 4.2.5]{Hovey}. That is, given generating cofibrations $f \colon K \to L \in \cat{sSet}$ and $g \colon A \to B$ in $\sfC$, we must check that $f \square F^T(g)$, determined from the pushout diagram
	\[\begin{tikzcd} K \times F^T(A) \rar \dar & L \times F^T(A) \dar \arrow[bend left = 25]{rdd} & \\
	K \times F^T(B) \rar \arrow[bend left = -15]{drr} & L \times F^T(A)  \sqcup_{K \times F^T(A)} K \times F^T(B) \arrow{rd}{f \square F^T(g)} & \\
	& & L \times F^T(B), \end{tikzcd} \]
	is a cofibration, which is a trivial cofibration if $f$ or $g$ is a trivial cofibration. By Lemma \ref{lem:algt-simplicial} and the fact that $F^T$ is a left adjoint, we may identify $f \square F^T(g)$ with $F^T(f \square g)$. Since $\times \colon \cat{sSet} \times \sfC \to \sfC$ is a Quillen bifunctor, $f \square g$ is a cofibration and hence so is $F^T(f \square g)$. Similarly, if $f$ or $g$ is a trivial cofibration, then so is $f \square g$ and hence $F^T(f \square g)$.
\end{proof}

We note that it is usually \emph{not} the case that the monad $T$ preserves cofibrations between objects that are not cofibrant.
	
To compute derived functors we will use simplicial resolutions of $T$-algebras, and hence we will impose the condition that $T$ preserves geometric realizations. This assumption is not of a model-categorical nature, but only becomes relevant when doing homotopy theory. It should be thought of as analogous to preserving sifted colimits (indeed, preserving enriched or $\infty$-categorical sifted colimits means preserving filtered colimits and geometric realizations).

The identity morphism of $\Delta^n \times X$ yields an $n$-simplex of $\Map_\sfC(X,\Delta^n \times X)$. Applying a simplicial functor $F$ we obtain an $n$-simplex of $\Map_\sfC(F(X),F(\Delta^n \times X))$ which in turn yields a morphism $\Delta^n \times F(X) \to F(\Delta^n \times X)$. This is natural in $n$ and $X$, so given a simplicial object $X_\bullet$ we get a map
\[|F(X_\bullet)| =  \int^{n \in \Delta^\op} \Delta^n \times F(X_n) \lra \int^{n \in \Delta^\op} F(\Delta^n \times X_n).\]
Upon applying $F$ to the natural maps $\Delta^n \times X_n \to |X_\bullet|$, we obtain a map from the right side to $F(|X_\bullet|)$. Taking the functor to be $T$, we can ask for the compositon of these two maps to be an isomorphism:

\begin{axiom}\label{axiom:monad-geomrel} $T$  and preserves geometric realization, in the sense that the natural map $|TX_\bullet| \to T|X_\bullet|$ is an isomorphism.\end{axiom}

Axioms \ref{axiom:monad-proj} and \ref{axiom:monad-geomrel} are quite restrictive, but hold under reasonable assumptions when the monad is obtained from an operad $\cO$, as we will show in Section \ref{sec:model-str-alg-over-operads}.

\subsection{Deriving constructions on $T$-algebras} \label{sec:derived-functors} Now that we have a model structure on $\Alg_T(\sfC)$, we can derive all the important constructions on $T$-algebras from the first part of this paper.

\subsubsection{Derived cell attachments}\label{sec:derived-cell-attachment} A model structure on $T$-algebras in $\sfC$ allows us to define derived cell attachments in $\Alg_T(\sfC)$. In Section \ref{sec:cell-def} we defined a cell attachment in $\Alg_T(\sfC)$ for $\sfC = \sfS^\sfG$, and it depended on the data of an $\gX_0 \in \Alg_T(\sfS^\sfG)$, a cofibration of simplicial sets $\partial D^{d} \hookrightarrow D^d$, an object $g$ of $\sfG$, and a morphism $e \colon \partial D^d \to X_0(g)$. The cell attachment was then defined to be a pushout
\[\begin{tikzcd}
	F^T(\partial D^{g,d}) \dar \rar{e} & \gX_0 \dar \\
	F^T(D^{g,d}) \rar & \gX_1, \end{tikzcd}\]
in $\Alg_T(\sfS^\sfG)$ as in Diagram (\ref{eq:6}). This is not necessarily homotopy invariant, and to remedy this, we should replace the pushout by a homotopy pushout. This can be done replacing the diagram by one where all three objects are cofibrant and one of the two maps is a cofibration. Since 
\[\cat{sSet} \overset{s}\lra \sfS \overset{g_*}\lra \sfS^\sfG \overset{F^T}\lra \Alg_T(\sfS^\sfG)\]
is a composition of left Quillen functors, the left-hand map is a cofibration between cofibrant objects, this pushout diagram is a homotopy pushout if $\gX_0$ is cofibrant. Thus we may derive the cell attachment by taking a cofibrant approximation $c\gX_0 \overset{\sim}\to \gX_0$. The attaching map $e$ lifts because the map $cX_0(g) \to X_0(g)$ on underlying object is a trivial fibration because $c\gX_0 \to \gX_0$ is, and the source of $e$ is cofibrant because $s$ is a left Quillen functor. 

Since $F^T$ is a left Quillen functor, if $X \hookrightarrow Y$ is a cofibration in $\sfC$ then $F^T(X) \to F^T(Y)$ is a cofibration in $\Alg_T(\sfC)$, and similarly for trivial cofibrations. Hence a (transfinite) composition of pushouts along such  maps is a cofibration, and similarly for trivial cofibrations, so that in particular a cellular map is a cofibration.

\subsubsection{Derived change-of-monads}
\label{sec:deriv-change-monads} We can similarly derive change-of-monad functors, lifting the adjunction of Lemma \ref{lem:morphism-monads-adjunction} to a Quillen adjunction.

\begin{lemma}
Suppose that $(R,\phi) \colon (\sfC,T) \to (\sfC',T')$ is a map of monads such that $R$ has a left adjoint $L$ forming a Quillen adjunction, and $T$, $T'$ are sifted monads satisfying Axiom \ref{axiom:monad-proj}. Then the adjunction 
	\[\begin{tikzcd}
	\Alg_T(\sfC) \arrow[shift left=.5ex]{r}{(L,\phi)_*} & \Alg_{T'}(\sfC') \arrow[shift left=.5ex]{l}{(R,\phi)^*}\end{tikzcd}\]
	is a Quillen adjunction.
\end{lemma}

\begin{proof}
	It suffices to check that $(R,\phi)^*$ preserves fibrations and trivial fibrations, but this is clear since both are verified by forgetting the algebra structure and $R$ is a right Quillen functor.
\end{proof}

Using functorial cofibrant and fibrant replacement functors we obtain derived functors $\bL(L,\phi)_* \colon \Alg_T(\sfC) \to \Alg_{T'}(\sfC')$ and $\bR(R,\phi)^* \colon \Alg_{T'}(\sfC)' \to \Alg_{T}(\sfC)$.

\subsubsection{Derived indecomposables}
\label{sec:deriv-indec-1} Since indecomposables are an example of a change-of-monad functor, the previous section provides a derived functor of indecomposables. Due to its importance we spell out the details.

\newglossaryentry{qlt}{%
	name={\ensuremath{Q^T_\bL}},
	description={Derived $T$-algebra indecomposables},
type=symbols
}
\begin{definition}\label{def:derived-indecomposables} Let $T$ be a sifted monad on $\sfC$ and $\epsilon \colon T \to +$ be an augmentation. The \emph{derived indecomposables}\index{indecomposables!derived} of a  $T$-algebra $\gX$ with respect to the augmentation $\epsilon$ are
  \[
    \gls{qlt}(\gX) \coloneqq \bL\epsilon_*(\gX) \in \sfCstar.
  \]
That is, $Q^T_\mathbb{L}(\gX) \simeq Q^T(c\gX)$ for a cofibrant approximation $c\gX \overset{\sim}{\to} \gX$.
\end{definition}

The derived indecomposables functor $Q^T_\mathbb{L}$ inherits many properties from $Q^T$. In particular, Lemma \ref{lem:ind-change-of-monads} says that if $(R,\phi) \colon (\sfC,T) \to (\sfC',T')$ is a morphism of augmented monads with $R$ the right adjoint in a Quillen adjunction $L \dashv R$ then we have
\[Q^{T'}_\bL(\bL(L,\phi)_* \gX) \simeq \bL(L,\upsilon)_*(Q^T_\bL(\gX)),\]
where $(L,\upsilon)_* \colon \sfC_* \to \sfC_*$ is the functor induced by $L$ on pointed objects.

\begin{example}The analogue of Example \ref{exam:morphism-over-id-indec} is that if $\phi \colon T \to T'$ is a map of monads augmented over $+$, then there is a natural weak equivalence
\[Q^{T'}_\mathbb{L}(\bL\phi_*(\gX)) \simeq Q^T_\mathbb{L}(\gX).\]
\end{example}

\begin{example}Given the data in Example \ref{exam:change-groupoid-indec}, for $\gX \in  \Alg_T(\sfC^{\sfG})$ change-of-indexing-category induces a natural weak equivalence
\[Q^{T'}_\bL(\bL\phi_*(\gX)) \simeq \bL p_*(Q^T_\bL(\gX)).\]
\end{example}

\subsection{Simplicial formulae}
 \label{sec:simpl-resol}

When we extend a right $T$-module functor $H \colon \sfC \to \sfD$ preserving sifted colimits to a functor $G \colon \Alg_T(\sfC) \to \sfD$ using the techniques in Section \ref{sec:funct-out-eilenb}, its value on $\gX \in \Alg_T(\sfC)$ is given by a reflexive coequalizer. While we could compute $\bL G(\gX)$ by applying $G$ to a cofibrant approximation of $\gX$, the result usually is rather inexplicit. The relevant reflexive coequalizer diagram is the truncations of a simplicial objets and in this section we discuss why under mild conditions $\bL G(\gX)$ can be described as the geometric realization of this simplicial object.

\subsubsection{Geometric realization}\label{sec:geom-rel-generalities} \index{geometric realization!homotopy theory of} When using simplicial objects to produce resolutions, we intend to eventually take their geometric realization. For this to be homotopically well-behaved, our resolutions need to be sufficiently cofibrant. One could ask for cofibrancy in the projective model structure on simplicial objects, but this is too strong a condition. Instead it suffices to ask for Reedy cofibrancy.

For a simplicial object $X_\bullet \in \cat{sC}$, the \emph{$n$th latching object} $L_n(X_\bullet)$ is 
\begin{equation} \label{eqn:nth-latching} L_n(X_\bullet) \coloneqq \colim_{[q] \twoheadleftarrow [n],[q] \neq [n]} X_q.\end{equation}
There is a natural map $L_n(X_\bullet) \to X_n$, by adding the identity $[n] \to [n]$ to the diagram indexing the colimit, which is terminal. This is called the \emph{latching map}. 

\begin{definition}
\label{def:reedy-cofibrant} A simplicial object in $\sfC$ is \emph{Reedy cofibrant}\index{simplicial object!Reedy cofibrant} if all latching maps are cofibrations in $\sfC$.
\end{definition}

This is an example of a general construction, see \cite{RiehlVerityReedy}. The category $\Delta^\mr{op}$ has a \emph{Reedy structure}. A Reedy structure on a category $\sfD$ consists of subcategories $\sfD^+,\sfD^- \subset \sfD$ and a degree functor from $\sfD$ to some ordinal such that every non-identity morphism in $\sfD^+$ strictly increases degree, every non-identity morphism in $\sfD^-$ strictly decreases degree, and every morphism $f$ factors uniquely as $f = gh$ with $g \in \sfD^+$ and $h \in \sfD^-$. For example, on $\Delta$ the degree function is given by the cardinality, $\Delta^-$ consists of surjections and $\Delta^+$ of injections. This data also induces a Reedy structure on $\Delta^\mr{op}$. For a functor $F \colon \sfD \to \sfC$ the $n$th latching object $L_n F$ is given on $D \in \sfD$ by the colimit over the subcategory of $\sfD^+ \downarrow D$ of objects of strictly lower degree than $D$. For $\Delta^\mr{op}$, the $n$th latching object is given by the formula in (\ref{eqn:nth-latching}).

If $\sfD$ has a Reedy structure and $\sfC$ is cofibrantly generated, then $\sfC^\sfD$ has a \emph{Reedy model structure}\index{model category!Reedy}, where weak equivalences are levelwise and the cofibrations are the maps $f \colon X_\bullet \to Y_\bullet$ such that the maps $X_n \sqcup_{L_n(X_\bullet)} L_n(Y_\bullet) \to Y_n$ are cofibrations in $\sfC$. Definition \ref{def:reedy-cofibrant} describes cofibrant objects in the Reedy model structure on simplicial objects. See \cite[Chapter VII]{GoerssJardine} or \cite[Chapter 15]{Hirschhorn} for more information on Reedy model structures. The following collates Theorems 18.6.6 (1) and 18.6.7 (1) of \cite{Hirschhorn}.

\begin{lemma}\label{lem:geom-rel-cofibrations} Geometric realization $\gr{-} \colon \cat{sC}\to \cat{C}$ preserves Reedy cofibrations and Reedy weak equivalences between Reedy cofibrant simplicial objects. In particular, it sends Reedy cofibrant simplicial objects to cofibrant objects.\end{lemma}

We shall give a criterion for a simplicial object to be Reedy cofibrant.

\begin{definition}A simplicial object $X_\bullet$ \emph{has split degeneracies}\index{simplicial object!split degeneracies} if there are objects $N_p(X_\bullet)$ and morphisms $N_p(X_\bullet) \to X_p$ for all $p \geq 0$ such that the induced map
	\[\bigsqcup_{[p] \twoheadrightarrow [q]} N_q(X_\bullet) \lra X_p\]
is an isomorphism.\end{definition}

\begin{lemma}\label{lem:reedy-cofibrant-split-degeneracies} A simplicial object $X_\bullet$ with split degeneracies is Reedy cofibrant if and only if each $N_p(X_\bullet)$ is cofibrant.\end{lemma}

\begin{proof}Note that 
	\[L_n(X_\bullet) = \bigsqcup_{[n] \twoheadrightarrow [q],[n] \neq [q]} N_q(X_\bullet)\]
	so that the map $L_n(X_\bullet) \to X_n$ is given by taking the coproduct of $L_n(X_\bullet)$ with $N_n(X_\bullet)$. This is a cofibration if and only if $N_n(X_\bullet)$ is cofibrant.
\end{proof}

The subcategory $\Delta_\mr{inj} \subset \Delta$ has a trivial Reedy structure: $\Delta^+ = \Delta_\mr{inj}$, and $\Delta^-_\mr{inj}$ consists of identities only. In this Reedy model structure both weak equivalences and cofibrations are levelwise. 

The inclusion $\sigma \colon \Delta_\mr{inj}^\mr{op}\to \Delta^\mr{op}$ defines a restriction functor $\sigma^* \colon \cat{sC} \to \cat{ssC}$ by forgetting degeneracy maps, which has a left adjoint $\sigma_*$  given by freely adding degeneracies. The adjunction
\[\begin{tikzcd}
	\cat{ssC} \arrow[shift left=.5ex]{r}{\sigma_*} & \arrow[shift left=.5ex]{l}{\sigma^*} \cat{sC}
\end{tikzcd}\] 
is a Quillen adjunction if we endow both sides with the Reedy model structures. It is evident that $\sigma_*$ preserves weak equivalences. By the formula \eqref{eqn:add-degeneracies} for $\sigma_*$, the latching maps are given by the map
\[X_n \sqcup \bigsqcup_{[n] \twoheadrightarrow [q],[n] \neq [q]} Y_q \lra \bigsqcup_{[n] \twoheadrightarrow [q]} Y_n\]
induced by $X_n \to Y_n$. These are hence cofibrations if all maps $X_p \to Y_p$ are. 

This has two useful consequences worth recording separately. Recalling that $||X_\bullet|| \cong |\sigma_* X_\bullet|$, Lemma \ref{lem:geom-rel-cofibrations} yields: 

\begin{lemma}\label{lem:thick-geom-rel-cofibrations} Thick geometric realization $\fgr{-} \colon \cat{ssC} \to \cat{C}$ preserves Reedy cofibrations and Reedy weak equivalences between Reedy cofibrant semi-simplicial objects. In particular, it sends Reedy cofibrant semi-simplicial objects to cofibrant objects.
\end{lemma}

Note that a Reedy cofibrant semi-simplicial object is simply one that is levelwise cofibrant, so we conclude:

\begin{lemma}\label{lem:sigma-semi-reedy} If $X_\bullet$ is a levelwise cofibrant semi-simplicial object, then $\sigma_* X_\bullet$ is Reedy cofibrant.\end{lemma}

By writing a simplicial object as a colimit of representables and switching two coends, one proves that the thick geometric realization $\fgr{\sigma^* X_\bullet}$, of the restriction of a simplicial object $X_\bullet$ to a semi-simplicial object, may also be computed as the coend $\int^{n \in \Delta^\mr{op}} \fgr{\sigma^* \Delta^n} \times X_n$. The canonical map $\fgr{\sigma^* \Delta^n} \to \Delta^n$ then induces a natural transformation $\fgr{\sigma^*-} \Rightarrow \gr{-}$ of functors $\cat{sC} \to \sfC$.

\begin{lemma}\label{lem:thick-to-thin} If $X_\bullet$ is a Reedy-cofibrant simplicial object, then $\fgr{\sigma^* X_\bullet} \to \gr{X_\bullet}$ is a weak equivalence.\end{lemma}

\begin{proof}By \cite[Corollary 18.4.14]{Hirschhorn}, it suffices to prove that $||\Delta^{\bullet}|| \to \Delta^\bullet$ is a weak equivalence between Reedy cofibrant cosimplicial objects. It is clearly a weak equivalence, and its target is Reedy cofibrant \cite[Corollary 15.9.11]{Hirschhorn}, which is proven using \cite[Corollary 15.9.10]{Hirschhorn} saying that it suffices to verify that the maximal augmentation is empty. This is also true for $||\Delta^{\bullet}||$.\end{proof}

\begin{remark}We remark that when working in $\cat{Top}$ one may weaken the required cofibrancy conditions of Lemma \ref{lem:geom-rel-cofibrations} and \ref{lem:sigma-semi-reedy} by playing the Str\o m and Serre model structures off each other. In particular, in $\cat{Top}$ the functor $\fgr{-}$ sends levelwise weak equivalences to weak equivalences, and for $\gr{-}$ this is true when the simplicial spaces involved are proper.
\end{remark}

\subsubsection{Extra degeneracies} We next give a useful condition for a map out of the geometric realization of a simplicial object to be a homotopy equivalence. 

An \emph{augmentation}\index{simplicial object!augmentation} of simplicial object $X_\bullet$ is a morphism $\epsilon \colon X_0 \to X_{-1}$ coequalizing $d_0,d_1 \colon X_1 \to X_0$. As the geometric realization $\gr{X_\bullet}$ maps to the coequalizer of 
\[\begin{tikzcd}X_1 \arrow[shift left=.5ex]{r}{d_0} \arrow[shift left=-.5ex,swap]{r}{d_1} & X_0,\end{tikzcd}\] and hence to $X_{-1}$, there is a canonical map $\epsilon \colon \gr{X_\bullet} \to X_{-1}$. The data of an augmentation is the same as an extension from $\Delta^\mr{op}$ to the category of possibly empty finite ordered sets.

An augmented simplicial object is said to have an \emph{extra degeneracy}\index{simplicial object!extra degeneracy} if there are maps $s_{-1} \colon X_i \to X_{i+1}$ for $i \geq -1$ satisfying the additional simplicial identities $\epsilon s_{-1} = \mr{id}$, $d_0 s_{-1} = \mr{id}$, $d_{i+1}s_{-1} = s_{-1}d_i$, $s_{j+1}s_{-1} = s_{-1}s_j$. There is a similar definition for semi-simplicial objects, omitting the final additional simplicial identity. Equivalently, a semi-simplicial object $X_\bullet$ admits an extra degeneracy when the simplicial object $\sigma_* X_\bullet$ does. 

\begin{lemma}\label{lem:extra-degeneracy}
	If $X_\bullet \in \cat{sC}$ has an augmentation to $X_{-1}$ with an extra degeneracy, then $\epsilon \colon \gr{X_\bullet} \to X_{-1}$ is a weak equivalence in $\sfC$.  The same is true for semi-simplicial objects and the thick geometric realization.
\end{lemma}

\begin{proof}It suffices to prove this for simplicial objects. By \cite[\S 6]{MeyerBar}, an extra degeneracy is equivalent to giving a pair of simplicial maps
	\[\mr{const}(X_{-1})_\bullet \xrightarrow{s_{-1}} X_\bullet \overset{d_0}\lra \mr{const}_{\bullet}(X_{-1})\]
so that $d_0 \circ s_{-1}$ is the identity and $s_{-1} \circ d_0$ admits a simplicial homotopy to the identity. Upon geometric realization, these respectively yield the identity and a homotopy equivalence in the sense of Definition \ref{def:homotopy-equiv}, hence weak equivalences.
\end{proof}

\subsubsection{Geometric realization of algebras} 
\newglossaryentry{grt}{%
	name={\ensuremath{\gr{-}_T}},
	description={Geometric realization in category of $T$-algebras},
type=symbols
}
The category $\Alg_T(\sfC)$ is simplicial and has colimits, so it has its own geometric realization functor $\gr{-}_T \colon \cat{s}\Alg_T(\sfC) \to \Alg_T(\sfC)$. We learned the following from Johnson--Noel \cite{JohnsonNoel}. Recall we assumed that Axiom \ref{axiom:monad-geomrel} holds.

\begin{lemma}\label{lem:monad-functors-geomrel} The following are natural isomorphisms:
	\begin{enumerate}[(i)]\item For $X_\bullet \in \cat{sC}$, the natural map $F^T \gr{X_\bullet} \to \gr{F^TX_\bullet}_T$ in $\Alg_T(\sfC)$.
		\item For $\gX_\bullet \in \cat{s}\Alg_T(\sfC)$, the natural map $Q^T \gr{\gX_\bullet}_T \to \gr{Q^T \gX_\bullet}$ in $ \sfC_\ast$.
		\item For $\gX_\bullet \in \cat{s}\Alg_T(\sfC)$, the natural map $\gr{U^T\gX_\bullet} \to U^T\gr{\gX_\bullet}_T$.
	\end{enumerate}
In particular, for $\gX_\bullet \in \cat{s}\Alg_T(\sfC)$ the structure map $TX_\bullet \to X_\bullet$ makes $\gr{U^T\gX_\bullet}$ into a $T$-algebra, and it follows that this is $\gr{\gX_\bullet}_T$.
\end{lemma}

\begin{proof}Parts (i) and (ii) follow because $F^T \colon \sfC \to \Alg_T(\sfC)$ and $Q^T \colon \Alg_T(\sfC) \to \sfC$ are left adjoints and preserve the copowering by simplicial sets by construction of the copowering of $\Alg_T(\sfC)$ in Lemma \ref{lem:algt-simplicial}.

For part (iii), the following commutative diagram has top row expressing $U^T|\gX_\bullet|_T$ as a reflexive coequalizer
\[\begin{tikzcd} U^T|F^T T U^T\gX_\bullet|_T \arrow[shift left=1ex]{r} \arrow[shift left=-1ex]{r} & U^T|F^T U^T \gX_{\bullet}|_T \rar \lar & U^T|\gX_\bullet|_T \\
U^TF^T|TU^T \gX_{\bullet}|\arrow[shift left=.5ex]{r} \arrow[shift left=-.5ex]{r} \uar{\cong} \ar[equals]{d} & U^TF^T|U^T\gX_\bullet| \uar{\cong} \ar[equals]{d}  & \\
T|TU^T \gX_\bullet| \arrow[shift left=.5ex]{r} \arrow[shift left=-.5ex]{r}  & T|U^T\gX_\bullet| & \\
{|T^2 U^T \gX_\bullet|} \uar{\cong} \arrow[shift left=1ex]{r} \arrow[shift left=-1ex]{r} & {|TU^T \gX_\bullet|} \rar \lar \uar{\cong} & {|U^T \gX_\bullet|}
\end{tikzcd}\]
while the bottom row expresses $|U^T\gX_\bullet|$ as a reflexive coequalizer. The top maps are isomorphisms by (i) and bottom maps by Axiom \ref{axiom:monad-geomrel}.
\end{proof}

Because we only used that $F^T$ and $Q^T$ are left adjoints and preserve the copowering by simplicial sets, parts (i) and (ii) of this lemma hold more generally, e.g.\ for the thick geometric realization of semi-simplicial objects and without Axiom \ref{axiom:monad-geomrel}. However, part (iii) might not hold for thick geometric realization, as it might not commute with $T$ because $T$ often does \emph{not} commute with the functor $\sigma_*$ which freely adds degeneracies.

\subsubsection{Free resolutions}\label{sec:FreeRes} One use of simplicial objects is to resolve $T$-algebras by free $T$-algebras.

\begin{definition}\label{def:Free-Res}
An augmented simplicial $T$-algebra $\epsilon \colon \gX_\bullet \to \gX$ is a \emph{free simplicial resolution}\index{resolution!free} of $\gX$ if
\begin{enumerate}[(i)]
\item the map $\epsilon \colon \gr{\gX_\bullet}_T \to \gX$ is a weak equivalence,
\item each $\gX_p$ is a free $T$-algebra,
\item $\gX_\bullet \in \cat{s}\Alg_T(\sfC)$ is Reedy cofibrant.
\end{enumerate}
\end{definition}

By Lemma \ref{lem:geom-rel-cofibrations}, if $\gX_\bullet$ is a free simplicial resolution then $\gr{\gX_\bullet}_T \in \Alg_T(\sfC)$ is cofibrant. Because the weak equivalence $\epsilon \colon \gr{\gX_\bullet}_{T} \to \gX$ need not be a fibration, this is not necessarily a cofibrant replacement but only a cofibrant approximation. However, as explained in Section \ref{sec:axioms-of-model-cats} we can still use it to compute left derived functors: if $F \colon \Alg_T(\sfC) \to \sfD$ is a functor which commutes with geometric realization and preserves trivial cofibrations between cofibrant objects, there are natural weak equivalences
\[\bL F(X) \simeq F(|\gX_{\bullet}|_T) \cong |F(\gX_\bullet)|.\]

\subsubsection{The monadic bar resolution}\label{sec:monadic-bar-resolution} 
\newglossaryentry{Bbullet}{%
	name={\ensuremath{B_\bullet}},
	description={Bar construction},
type=symbols
}
Recall that we write $X = U^T(\gX)$ for the underlying object in $\sfC$ of a $T$-algebra $\gX$. There is an explicit source of free simplicial resolutions coming from the \emph{monadic bar resolution}\index{resolution!monadic bar} $\gls{Bbullet}(F^T,T,\gX) \to \gX$ given by
\[[p] \longmapsto B_p(F^T,T,\gX) = F^T T^p X,\]
with face maps given by the $T$-algebra structure map of $\gX$ and the right $T$-module functor structure of $F^T$, and with degeneracy maps given by the unit of the monad $T$.

This is visibly a free $T$-algebra in each degree, establishing (ii) of Definition \ref{def:Free-Res}. The unit transformation $1^T \colon \mr{id} \to T$ of the monad gives a map $X \to T(X)$ which is \emph{not} a $T$-algebra map, but nonetheless equips the underlying augmented simplicial object $U^T B_\bullet(F^T,T,\gX) \to X$ in $\sfC$ with an extra degeneracy, so
\[U^T\gr{B_\bullet(F^T,T,\gX)} \lra X\]
is a weak equivalence by Lemma \ref{lem:extra-degeneracy}. This establishes (i) by noting that Lemma \ref{lem:monad-functors-geomrel} (iii) identifies this map with $U^T\gr{B_\bullet(F^T,T,\gX)}_T \to X$.

In general it need not be the case that $B_\bullet(F^T,T,\gX)$ is Reedy cofibrant, but we do have the following lemma. We can form the simplicial $T$-algebra $\sigma_* \sigma^* B_\bullet(F^T,T,\gX)$ by freely adding degeneracies to the semi-simplicial $T$-algebra $\sigma^* B_\bullet(F^T,T,\gX)$. Here it is important that we are freely adding degeneracies in the category of $T$-algebras, as $\sigma_*$ does not commute with $F^T$ (though $\sigma^*$ does). We call the simplicial object $\sigma_* \sigma^* B_\bullet(F^T,T,\gX)$ the \emph{thick monadic bar construction}\index{resolution!thick monadic bar}.

\begin{lemma}\label{lem:adding-degeneracies-reedy} If $X$ is cofibrant in $\sfC$, then $\sigma_* \sigma^* B_\bullet(F^T,T,\gX)$ is Reedy cofibrant and hence a free simplicial resolution of $\gX$.
\end{lemma}

\begin{proof}The following computation
	\begin{align*}\sigma_* \sigma^* B_p(F^T,T,\gX) &\cong \bigsqcup^{\Alg_T(\sfC)}_{\alpha \colon [p] \twoheadrightarrow [q]} B_q(F^T,T,X) \\
	&\cong \bigsqcup^{\Alg_T(\sfC)}_{\alpha \colon [p] \twoheadrightarrow [q]} F^T T^q(X) \\
	&\cong F^T\left(\bigsqcup^{\sfC}_{\alpha \colon [p] \twoheadrightarrow [q]} T^q(X)\right)\end{align*}
shows that $\sigma_* \sigma^* B_\bullet(F^T,T,\gX)$ is levelwise free and has split degeneracies. Lemma \ref{lem:sigma-semi-reedy} says it is Reedy cofibrant if all $F^T(T^p(X))$ are cofibrant in $\Alg_T(\sfC)$. These are cofibrant if $T^p(X)$ is cofibrant in $\sfC$, which follows from the assumption that $T$ preserves cofibrant objects.
\end{proof}

\subsubsection{A simplicial formula for derived cell attachments} \label{sec:derived-cell-simplicial} 
The derived cell attachment as discussed in Section \ref{sec:derived-cell-attachment} may be described simplicially, by extending the reflexive coequalizer diagram~(\ref{eq:9}) to a simplicial diagram. \index{cell attachment!derived}

To derive the attachment of a cell to $\gX \in \Alg_T(\sfC^\sfG)$, we replace it by a cofibrant $T$-algebra $c\gX$, lift the attaching map to $c\gX$ (which we have explained is possible), and form the pushout. If $F^T$ commutes with geometric realization then the geometric realization of a free simplicial resolution $\epsilon \colon \gr{\gX_{\bullet}}_T \overset{\sim}\to \gX$ provides a cofibrant approximation, and we need that the attaching map factors through $\epsilon$. 

If $T$ preserves cofibrations and $X$ is cofibrant in $\sfC$, then we may take $\gX_\bullet = \sigma_*\sigma^*B_\bullet(F^T, T, \gX)$, in which case the unit $i_X \colon X \to T(X) =U^TF^T(X) = U^T \gX$ shows that the attaching map does factor through $\epsilon$. The geometric realization of $\gX_\bullet$ is the thick geometric realization of $B_\bullet(F^T, T, \gX)$.  As thick geometric realization commutes with pushouts, an explicit simplicial formula for derived cell attachment is given by
\[\fgr{B_\bullet(F^T,T,\gX)}_T \cup^T_{F^T(\partial D^{g,d})} F^T(D^{g,d}) \cong \fgr{B_\bullet(F^T,T,\gX) \cup^T_{F^T(\partial D^{g,d})} F^T(D^{g,d})}_T.\]
Furthermore, the free algebra functor $F^T$ is a left adjoint and commutes with colimits, so we may write this as the thick geometric realization of the simplicial object
\[[p] \longmapsto F^T(T^p(\gX) \cup_{\partial D^{g,d}} D^{g,d}).\]
This is the formula that appears in \cite{KupersMiller}.

\subsubsection{A simplicial formula for derived indecomposables} \label{sec:simplicial-formula-indecomposables} As $Q^T \colon \Alg_T(\sfC) \to \sfC_\ast$ is a left Quillen functor, the discussion above implies that if $F^T$ commutes with geometric realization we have
\[Q^T_\bL(\gX) \simeq Q^T(\gr{\gX_\bullet}_T) \cong \gr{Q^T(\gX_\bullet)},\]
using Lemma \ref{lem:monad-functors-geomrel} (ii) for the second isomorphism. If $T$ preserves cofibrations and $X$ is cofibrant in $\sfC$ then by Lemma \ref{lem:adding-degeneracies-reedy} we have that $\sigma_* \sigma^* B_\bullet(F^T,T,\gX) \to \gX$ is a free simplicial resolution, so using the monadic bar resolution and the fact that $Q^T F^T \cong +$ we get the following explicit simplicial formula:\index{indecomposables!derived}
	\[Q^T_\bL(\gX) \simeq Q^T(\fgr{B_\bullet(F^T,T,\gX)}_T) \cong \fgr{ B_\bullet(+,T,\gX)}.\]

More generally, if $\phi \colon T \to T'$ is a morphism of monads which commute with geometric realization, and $\gX_{\bullet} \to \gX$ is a free simplicial resolution of $\epsilon \colon \gX \in \Alg_T(\sfC)$, then we have in $\Alg_{T'}(\sfC)$ that
\[\bL \phi_*(\gX) \simeq \phi_*(\gr{\gX_\bullet}_T) \cong \gr{\phi_*(\gX_\bullet)}_{T'}.\]
As before, when $X$ is cofibrant this may be made more explicit by choosing the monadic bar resolution.
\section{Homotopy theory of operads, algebras, and modules} \label{sec:homotopy-theory-operads} 

In this section we verify the axioms of Section \ref{sec:homotopy-theory-algebras} in the case of monads $T \colon \sfC \to \sfC$ arising from an operad $\cO$ in $\sfC$ as discussed in Section \ref{sec:operads}; this will mean that we can do homotopy theory in the category $\Alg_\cO(\sfC)$. As before, we assume Axioms \ref{axiom.CofGen} and \ref{axiom.model-monoidalsimplicial} throughout.

\subsection{Symmetric sequences} \label{sec:homotopy-sym-seq}  In Section \ref{sec:axioms-of-model-cats} we discussed the projective model category structure on diagram categories. The category of $k$-symmetric sequences in $\sfC$, $\cat{FB}_k(\sfC)$, is a diagram category and so has a projective model structure as in Section \ref{sec:FunctorCategories}, in which a $k$-symmetric sequence $\cX$ is cofibrant if and only if each $\cX(n)$ is a cofibrant object in the projective model structure on $\sfC^{G_n}$. It is a monoidal model category for the Day convolution product, but not necessarily for the the composition product. However, the following often suffices. 

\begin{lemma}\label{lem.symmetric-sequence-preserving} 
The bifunctor 
		\begin{align*}\cat{FB}_k(\sfC) \times \sfC &\lra \sfC \\
	(\cY,X) &\longmapsto \cY(X),\end{align*}
described in Definition \ref{def:symm-seq-application}, has the property that the pushout-product
\[f \square g \colon \cY(X') \sqcup_{\cY(X)} \cY'(X) \lra \cY'(X')\]
of $f \colon \cY \to \cY'$ and $g \colon X \to X'$ is a cofibration if $f$ and $g$ are cofibrations and $X$ is cofibrant, and a trivial cofibration if additionally $f$ or $g$ (or both) are trivial cofibrations. In particular, it enjoys the following properties:
\begin{enumerate}[(i)]
\item If $\cY \in \cat{FB}_k(\sfC)$ is cofibrant then $X \mapsto \cY(X) \colon \sfC \to \sfC$ preserves cofibrations, trivial cofibrations with cofibrant domain, and weak equivalences between cofibrant objects.
\item If $X \in \sfC$ is cofibrant then $\cY \mapsto \cY(X) \colon \cat{FB}_k(\sfC) \to \sfC$ preserves cofibrations, trivial cofibrations,  and weak equivalences between cofibrant objects.
\end{enumerate}
\end{lemma}

\begin{proof}
Using the formula $\cY(X) = 0^*(\cY \circ 0_*(X))$, this follows from the pushout-product property established in Theorem 6.2 (a) of \cite{HarperOperads} for the composition product, in the cases $k=1$ and $k=\infty$. Harper's argument also applies in the case $k=2$ (the fact that the automorphism groups $G_n = \beta_n$ in $\cat{FB}_2$ are infinite has no effect on the argument, as Harper never uses any of the finiteness hypotheses in Section 6 of \cite{HarperOperads}). See \cite[Lemma 11.5.1]{FresseBook} for a related discussion.
\end{proof}

\subsection{Algebras over operads} \label{sec:model-str-alg-over-operads} Next we discuss the homotopy theory of algebras over operads, verifying the axioms of Section \ref{sec:homotopy-theory-algebras} using the results in Sections \ref{sec.model-adjunctions} and \ref{sec:homotopy-sym-seq}. We start by verifying Axiom \ref{axiom:monad-geomrel}, which does not concern the model structure.

\begin{lemma}\label{lem:operad-geomrel} The monad associated to an operad commutes with geometric realization.\end{lemma}
\begin{proof}
This follows from part (\ref{enum:odot-props-geomrel}) of Lemma \ref{lem:odot-props}.
\end{proof}

Axiom \ref{axiom:monad-proj} requires the existence of the projective model structure transferred along the adjunction $F^\cO \dashv U^\cO$, and that $U^\cO$ preserves (trivial) cofibrations between cofibrant objects. We first show that this second part follows from the existence of projective model structure under a mild condition on $\cO$.

\begin{definition}\label{def:sigma-cofibrant}
	An operad $\cO$ is called \emph{$\Sigma$-cofibrant}\index{$\Sigma$-cofibrant} if its underlying symmetric sequence is cofibrant in $\cat{FB}_k(\sfC)$.
\end{definition}

\begin{remark}It is better to say that $\cO$ is $\Sigma$-cofibrant if the unit map $\bunit \to \cO$ is a cofibration. This implies our definition, as we assumed that $\bunit$ is cofibrant.\end{remark}

\begin{lemma}\label{lem:proj-lax-monoidal-axiom} 
If $\Alg_\cO(\sfC)$ has a projective model structure and $\cO$ is $\Sigma$-cofibrant, then $U^\cO$ preserves (trivial) cofibrations between cofibrant objects. That is, Axiom \ref{axiom:monad-proj} holds for $\Alg_\cO(\sfC)$.\end{lemma}

\begin{proof}We may restrict our attention to the cofibrations, as by definition of the projective model structure $U^\cO$ preserves weak equivalences. The projective model structure on $\Alg_\cO(\sfC)$ is a cofibrantly generated with set of generating cofibrations given by $F^\cO(I) \coloneqq \{F^\cO(i) \mid i \in I\}$ with $I$ the set of generating cofibrations of $\sfC$. Hence cofibrations in $\Alg_\cO(\sfC)$ are retracts of relative $F^\cO(I)$-cell complexes, which are (potentially transfinite) compositions of pushouts
	\[\begin{tikzcd}F^\cO(Y_0) \rar \dar[swap]{F^\cO(i)} & \gR_0 \dar{f} \\
	F^\cO(Y_1) \rar & \gR_1\end{tikzcd}\]
	in $\Alg_\cO(\sfC)$ with $i$ a generating cofibration of $\sfC$ \cite[Proposition 2.1.18]{Hovey}. 
		
	If $\gR_0$ is cofibrant, then the small object argument implies it is a retract of a $F^\cO(I)$-cell complex. Since cofibrations are closed under retracts and transfinite composition, in order to show that $U^\cO$ preserves cofibrations between cofibrant objects it suffices to show that given a pushout diagram 
	\[\begin{tikzcd}F^\cO(Y_0) \rar \dar[swap]{F^\cO(i)} & \gR_0 \dar{f} \\
	F^\cO(Y_1) \rar & \gR_1\end{tikzcd}\]
	where $\gR_0$ is a $F^\cO(I)$-cell complex, $U^\cO(f)$ is cofibration if $i$ is. This is proven in \cite[Proposition 5.1]{BergerMoerdijk} or \cite[Lemma 20.1.A]{FresseBook}.\end{proof}

It remains to discuss the existence of the projective model structure on $\Alg_\cO(\sfC)$. Most of the example model categories we have described in Section \ref{sec:examples-of-model-cats} have a lax $k$-monoidal fibrant approximation, so we treat that easier case first. For symmetric spectra with the absolute stable projective model structure we will need to give a different argument.

\subsubsection{Categories with lax $k$-monoidal fibrant approximation} Under standard smallness hypotheses, a sufficient condition for the existence of the projective model structure is the existence of a lax monoidal fibrant approximation as in Definition \ref{def.model-laxfibrant}. This is essentially \cite[Theorem 2.1]{BergerMoerdijk2}.

\begin{lemma}\label{lem:operad-fibrant-rep} If there exists a lax $k$-monoidal fibrant approximation on $\sfC$, then there exists a fibrant approximation as in Proposition \ref{prop:proj-gs} \eqref{enum:proj-iiprime} on $\Alg_\cO(\sfC)$.
\end{lemma}

\begin{proof}
Given a lax $k$-monoidal fibrant approximation $R \colon \sfC \to \sfC$ as in Definition \ref{def.model-laxfibrant} and an $\cO$-algebra $\gX$ with underlying object $X$, $R(X)$ may be endowed with the structure of an $\cO$-algebra using the following structure maps:
\[\cO(n) \otimes R(X)^{\otimes n} \lra R(\cO(n)) \otimes R(X)^{\otimes n} \lra R(\cO(n) \otimes X^{\otimes n}) \lra R(X),\]
where the first map uses the natural transformation $\mr{id}_\sfC \Rightarrow R$, the second map uses the lax monoidality of $R$, and the third map uses the $\cO$-algebra structure map of $\gX$. The natural transformation $\mr{id}_\sfC \Rightarrow R$ then induces a weak equivalence $\gX \to R(\gX)$ of $\cO$-algebras, and thus $R$ lifts to a fibrant approximation $R \colon \Alg_{\cO}(\sfC) \to \Alg_{\cO}(\sfC)$.\end{proof}

\begin{corollary}\label{cor:proj-lax-monoidal}
Suppose $\sfC$ is cofibrantly generated model category such that 
\begin{enumerate}[(i)]
	\item \label{enum:proj-lax-monoidal-i} $F^\cO(I)$ and $F^\cO(J)$ admit a small object argument.
	\item \label{enum:proj-lax-monoidal-ii} There exists a lax $k$-monoidal fibrant approximation on $\sfC$.
\end{enumerate}
Then the projective model structure on $\Alg_\cO(\sfC)$ exists.
\end{corollary}

\begin{proof}
We apply Theorem \ref{thm:general-proj-model-structure}. Condition \eqref{enum:proj-i} of that theorem holds by assumption \eqref{enum:proj-lax-monoidal-i}. To verify condition \eqref{enum:proj-ii} of that theorem, we check condition \eqref{enum:proj-iiprime} of Proposition \ref{prop:proj-gs}: Lemma \ref{lem:algt-simplicial} says $\Alg_\cO(\sfC)$ is simplicial, and by assumption \eqref{enum:proj-lax-monoidal-ii}, Lemma \ref{lem:operad-fibrant-rep} implies the existence of the appropriate fibrant approximation functor.\end{proof}

Note, however, that this model structure need not satisfy the second part of Axiom \ref{axiom:monad-proj} unless $\cO$ is $\Sigma$-cofibrant, cf.~Lemma \ref{lem:proj-lax-monoidal-axiom}.

\begin{example} If the domains of $I$ and $J$ are small with respect to all morphisms in $\sfC$, then condition \eqref{enum:proj-lax-monoidal-i} of Corollary \ref{cor:proj-lax-monoidal} follows from Lemmas \ref{lem:ut-preserves-sifted} and \ref{lem:small-object-condition}. This assumption holds for $\sfC = \cat{sSet}$ by \cite[Lemma 3.1.1]{Hovey} and $\sfC = \cat{sMod}_\bk$ by a similar argument, as well as $\sfC_\ast$ and $\sfC^\sfG$.
\end{example}

\begin{example}For $\sfC = \cat{Top}$ condition \eqref{enum:proj-lax-monoidal-i} of Corollary \ref{cor:proj-lax-monoidal} is satisfied because $I$ and $J$ consist of closed inclusions, $F^\cO$ preserves closed inclusions (as these are preserved by disjoint unions, products, and quotients by group actions), and CGWH topological spaces are small relative to closed inclusions (by \cite[Lemma 2.4.1]{Hovey} and the remarks following \cite[Proposition 2.4.22]{Hovey}). A similar argument applies to $\sfC_\ast$ and $\sfC^\sfG$.
\end{example}

\subsubsection{Symmetric spectra with the absolute stable model structure}The previous section applies to all $\sfS$ in Section \ref{sec:examples-of-model-cats} with the exception of the category $\cat{Sp}^\Sigma$ of symmetric spectra. 

\begin{remark}In general, the projective model structure of $\cO$-algebras in symmetric spectra transferred from the absolute stable model structure does \emph{not} exist, by an similar argument as in \cite{LewisSpectra}. For a contradiction take $\cO = \mr{Comm}$, the commutative operad, and suppose that the projective model structure on commutative algebras in symmetric spectra exists. A fibrant replacement of $\bS$ in this model structure would also be fibrant in $\cat{Sp}^\Sigma$. Taking its $0$th space we would obtain a strictly commutative model of $QS^0$, which would imply $QS^0$ is a product of Eilenberg--Maclane spaces (which of course it is not).\end{remark}

However, one might expect that the projective model structure \emph{does} exist if $\cO$ is a $\Sigma$-cofibrant operad. Such a result was sketched by Schwede, but Hornbostel has pointed out in the remark following Conjecture 3.7 of \cite{Hornbostel} that a crucial step (our Lemma \ref{lem:sigma-cofibrant-free}) is hard to verify without additional conditions on $\cO$, and mentions a counterexample of Fresse for operads in chain complexes over a field of positive characteristic.

We shall show that Schwede's argument works given a $\Sigma$-cofibrant operad $\cC$ in simplicial sets satisfying one of two conditions on its $0$-ary operations. Our goal is to apply Theorem \ref{thm:general-proj-model-structure}, where as usual \eqref{enum:proj-ii} is hardest to verify. Our main tool for this is the existence of the \emph{absolute injective stable model structure} on symmetric spectra: the cofibrations are levelwise cofibrations and the weak equivalences the stable equivalences (in this model structure every object is cofibrant). The existence of such a model structure for $\cat{Sp}^\Sigma$ is established in Section 5.3 of \cite{HoveyShipleySmith} (see also \cite[p.\ 367]{Schwedebook}).  When referring to (trivial) cofibrations, etc., in this model structure, we add the adjective ``injective''. Similarly, when referring to the (trivial) cofibrations, etc., in the absolute stable projective model structure, we add ``projective''.

\begin{proposition}\label{prop:existence-proj-model-structure-spectra} Let $\cC$ be a $\Sigma$-cofibrant operad in simplicial sets such that either 
	\begin{enumerate}[(i)]
		\item $\cC(0) \neq \varnothing$, or 
		\item $\cC$ is non-unitary and admits a unitalization $\cC^+$.
	\end{enumerate}
	Then the projective model structure on $\Alg_\cC(\cat{Sp}^\Sigma)$ transferred from the absolute stable projective model structure on $\cat{Sp}^\Sigma$ exists.
\end{proposition}

\begin{proof}We apply Theorem \ref{thm:general-proj-model-structure}. For (i), we remark that symmetric spectra are combinatorial, so that this is tautological because all objects are small.
	
	To verify (ii) we shall use the absolute stable injective model structure. It suffices to prove that if $Y_0 \to Y_1$ is a projective trivial cofibration of symmetric spectra and $\gR_0 \in \Alg_\cC(\cat{Sp}^\Sigma)$, then in a pushout diagram 
	\begin{equation}\label{eqn:pushout-2} \begin{tikzcd} F^\cC(Y_0) \rar \dar[swap]{F^\cC(i)} & \gR_0 \dar{f} \\
	F^\cC(Y_1) \rar & \gR_1, \end{tikzcd}\end{equation}
	the map $U^\cC(f)$ is a stable equivalence. We shall prove the stronger statement that it is an injective trivial cofibration. The proof involves a number of technical lemmas about symmetric spectra, which we postpone for the moment.
	
	As explained in Section \ref{sec:filt-stages-pushout}, the underlying map of the right hand side of (\ref{eqn:pushout-2}) admits a filtration $\mr{f}(\gR_1)$ in $\cat{Sp}^\Sigma$ with filtration steps given by pushout diagrams
	\begin{equation}\label{eqn:filt-stages-pushout-2} \begin{tikzcd}\mr{Env}_c^\cC(\gR_0) \otimes_{G_c}  Y_1^{\Box c} \rar \dar & \mr{f}(\gR_1)(c-1) \dar \\
	\mr{Env}_c^\cC(\gR_0) \otimes_{G_c} Y_1^{\otimes c} \rar & \mr{f}(\gR_1)(c), \end{tikzcd}\end{equation}
	where $\mr{Env}^\cC_c(\gR_0)$ is defined on free $\cC$-algebras in (\ref{eqn:eco}) and in general is extended by density using Proposition \ref{prop:extend-by-density}.	  
	
	Since injective trivial cofibrations are closed under transfinite composition, it suffices to show that for each $c > 0$ the right vertical map of (\ref{eqn:filt-stages-pushout-2}) is an injective trivial cofibration. As pushouts preserve trivial cofibrations in any model structure, it suffices to prove that the left vertical map of (\ref{eqn:filt-stages-pushout-2}) is an injective trivial cofibration. The iterated pushout-product $i^{\Box c} \colon Y_1^{\Box c} \to Y_1^{\otimes c}$ is a projective trivial cofibration because this is a monoidal model structure. By Lemmas \ref{lem:sym-spectra-levelwise-free} and \ref{lem:sigma-cofibrant-free}, the functor $\mr{Env}^\cC_c(\gR_0) \wedge_{\fS_c} -$ sends projective trivial cofibrations to injective trivial cofibrations.
\end{proof}

Let us now prove the remaining technical lemmas. The following is a variation of \cite[Proposition 6.12]{Schwedebook}.

\begin{lemma}\label{lem:sym-spectra-levelwise-free} If $E$ is a symmetric spectrum with a $G$-action which is levelwise free away from the basepoint, and $f \colon X \to Y$ is a projective trivial cofibration between symmetric spectra with $G$-actions which happens to be $G$-equivariant, then $E \wedge_{G} X \to E \wedge_{G} Y$ is an injective trivial cofibration.\end{lemma}

\begin{proof}The map $E \wedge X \to E \wedge Y$ is an injective trivial cofibration by Theorem 5.3.7 (5) of \cite{HoveyShipleySmith}, because $E$ is injective cofibrant (every object is), $f$ is an $S$-cofibration as in Definition 5.3.6 of \cite{HoveyShipleySmith} because it is a projective cofibration, and $f$ is a stable equivalence. This implies $E \wedge_G X \to E \wedge_G Y$ is an injective cofibration because taking the quotient by $G$ levelwise preserves monomorphisms.
	
	It is however not clear that taking the quotient by $G$ levelwise preserves stable equivalences. Hence, to verify it is also a stable equivalence, we consider the diagram
	\[\begin{tikzcd} EG_+ \wedge_G (E \wedge X) \rar{\simeq} \dar & \ast_+ \wedge_G(E \wedge X) \cong E \wedge_G X \dar \\
	EG_+ \wedge_G (E \wedge Y) \rar{\simeq} & \ast_+ \wedge_G (E \wedge Y) \cong E \wedge_G Y,\end{tikzcd}\]
	where the horizontal maps are levelwise weak equivalences (because the $G$-action on $E$ is levelwise free away from the basepoint) and hence stable equivalences. Hence it suffices to prove that the left vertical map is a stable equivalence.
	
	To do so, we filter $EG$ by skeleta $EG^{(k)}$ to obtain a filtration with filtration quotients given by $\Delta^k/\partial \Delta^k \wedge G^{\wedge k} \wedge G_+$ (with $G$ based at the identity and $G$ acting on the right term $G_+$). There is then a map of cofiber sequences of symmetric spectra
	\[\begin{tikzcd} EG^{(k-1)}_+ \wedge_G (E \wedge X) \dar \rar &[-10pt] EG^{(k)}_+ \wedge_G (E \wedge X) \rar \dar &[-10pt] \Delta^k/\partial \Delta^k \wedge G^{\wedge k} \wedge E \wedge X \dar \\
	EG^{(k-1)}_+ \wedge_G (E \wedge Y)  \rar & EG^{(k)}_+ \wedge_G (E \wedge Y) \rar  & \Delta^k/\partial \Delta^k \wedge G^{\wedge k} \wedge E \wedge Y.  \end{tikzcd}\]
	
	In this diagram, the right vertical map is an injective trivial cofibration as a consequence of Theorem 5.3.7 (5) of \cite{HoveyShipleySmith} because $\Delta^k/\partial \Delta^k \wedge G^{\wedge k} \wedge E$ is injective cofibrant (every object is), and $f$ is an $S$-cofibration because it is a projective cofibration. Because the two left horizontal maps are injective cofibrations between injective cofibrant objects, both rows are in fact homotopy cofiber sequences. Thus it follows by induction over $k$ that each map $EG^{(k)}_+ \wedge_G (E \wedge X) \to EG^{(k)}_+ \wedge_G (E \wedge Y)$ is a stable equivalence, the initial case $k = -1$ being trivial. It is also an injective cofibration, since it is levelwise injective by a similar argument as for $E \wedge_G X \to E \wedge_G Y$. Injective trivial cofibrations are closed under transfinite composition, so we conclude $EG_+ \wedge_G (E \wedge X) \to EG_+ \wedge_G (E \wedge Y)$ is an injective trivial cofibration.
\end{proof}

In the proof of Proposition \ref{prop:existence-proj-model-structure-spectra}, we want to apply Lemma \ref{lem:sym-spectra-levelwise-free} to the left vertical map of \eqref{eqn:filt-stages-pushout-2}. Hence we need to prove the following:

\begin{lemma}\label{lem:sigma-cofibrant-free} Let $\cC$ be $\Sigma$-cofibrant operad in simplicial sets satisfying either
	\begin{enumerate}[(i)]
		\item $\cC(0) \neq \varnothing$, or 
		\item $\cC$ is non-unitary and admits a unitalization $\cC^+$.
	\end{enumerate}
	If $\gR$ is an $\cC$-algebra in $\cat{Sp}^\Sigma$, then the $G_c$-action on $\mr{Env}^\cC_c(\gR)$ is levelwise free.\end{lemma}

\begin{proof}The symmetric spectrum with $G_c$-action $\mr{Env}^\cC_c(\gR)$ is defined by a reflexive coequalizer
	\[\begin{tikzcd}
	\mr{Env}_c(\cC)(\cC(\gR)) \arrow[shift left=1ex]{r} \arrow[shift left=-1ex]{r} & 	\mr{Env}_c(\cC)(\gR) \rar \lar &
	\mr{Env}^\cC_c(\gR).
	\end{tikzcd}\] 
	Since colimits of symmetric spectra are computed levelwise and the freeness of the $G_c$-action is a property of each of the individual levels, we may restrict our attention to the $k$th level $\mr{Env}^\cC_c(\gR)_k \in \cat{sSet}_\ast$, which is given by the coequalizer in $\cat{sSet}_\ast$ of the diagram
	\begin{equation}\label{eqn:level-pushout-env} \begin{tikzcd}\bigvee_{n \geq 0} \cC(n+c)_+ \wedge_{G_{n,c}} (\cC(\gR)^{\wedge n})_{k}  \arrow[shift left=.5ex]{r} \arrow[shift left=-.5ex]{r} & \bigvee_{n \geq 0} \cC(n+c)_+ \wedge_{G_{n,c}} (\gR^{\wedge n})_{k}.\end{tikzcd}\end{equation}
	
	The $G_c$-action on each of the right-hand terms is free (away from the basepoint) because on $p$-simplices there is $G_c$-equivariant function
	\[f_{p,k} \colon (\cC(n+c)_+ \wedge_{G_{n,c}} (\gR^{\wedge n})_k)_p \lra (\cC(c)_+)_p,\]
	given by sending the basepoint to the basepoint and an equivalence class $(o;r_1,\ldots,r_n)$ to the element of $\cC(c)_p$ obtained by inserting a fixed element of either (i) $\cC(0)$ or (ii) $\cC^+(0)$, into the first $n$ slots of $o \in \cC(n+c)$ to get a $p$-simplex of $\cC(c)$. By construction, this only maps the basepoint to the basepoint. Since $\cC(c)_p$ is a free $G_c$-set, by $\Sigma$-cofibrancy of $\cC$, and the map is $G_c$-equivariant, we conclude that the $G_c$-action on the domain of $f_{p,k}$ must be free away from the basepoint as well.
	
	We may compute the underlying pointed sets of $p$-simplices of the coequalizer of (\ref{eqn:level-pushout-env}) as the quotient of its right hand side by the equivalence relation $\sim$ generated by \[(o;o_1(r^1_1,\ldots,r^1_{k_1}),\ldots,o_n(r^n_1,\ldots,r^n_{k_n})) \sim (o \circ (o_1,\ldots,o_n);r^1_1,\ldots,r^n_{k_n}).\]
	The map $f_{p,k}$ is not compatible with $\sim$ because it may happen that one of the $o_i(r^i_1,\ldots,r^i_{k_i})$ is the basepoint without any of the $r^i_j$ being the basepoint. However, note that the subset of elements of $(\cC(n+c)_+ \wedge_{G_{n,c}} (\gR^{\wedge n})_k)_p$ that are identified with the basepoint under the equivalence relation is closed under the $G_c$-action: if $x$ is equivalent to the basepoint via $x \sim x_1 \sim \ldots \sim x_r \sim \ast$, then $gx$ is equivalent to the basepoint via $gx \sim gx_1 \sim \ldots \sim gx_r \sim g\ast = \ast$. Hence we may restrict $f_{p,k}$ to the complement of this subset to get a $G_c$-equivariant map of pointed sets $(\mr{Env}^\cC_c(\gR)_k)_p \to (\cC(c)_+)_p$ which sends only the basepoint to the basepoint. As before, because the codomain has free $G_c$-action away from the basepoint the same is true for the domain.
\end{proof}

Finally, we explain how to generalize this to the diagram categories $\sfC = (\cat{Sp}^\Sigma)^\sfG$. We want to prove the existence of the model structure on $\Alg_\cC(\sfC)$ transferred from the projective model structure on $\sfC$ induced by the absolute stable projective model structure on $\cat{Sp}^\Sigma$, the ``projective projective model structure.'' Our tool shall be the injective model structure on $\sfC$ induced by the absolute stable injective model structure on $\cat{Sp}^\Sigma$, the ``injective injective model structure.'' The injective model structure on $\sfC$ exists by Proposition A.2.8.2 of \cite{LurieHTT}, as the absolute stable injective model structure on $\cat{Sp}^\Sigma$ is combinatorial. As before, the advantage of this injective injective model structure is that the cofibrations are those maps $E \to F$ that for each object $g \in \sfG$ the map $E(g) \to F(g)$ of symmetric spectra is levelwise a cofibration of pointed simplicial sets, and thus, as before, every object is cofibrant in this model structure.

The previous argument then needs to be modified in Lemma \ref{lem:sym-spectra-levelwise-free}: in particular, we need to explain why $F \otimes_{\sfG} -$ sends (trivial) projective projective cofibrations to (trivial) injective injective cofibrations. It suffices to verify this only for generating projective projective trivial cofibrations $\sfG(g,-)_+ \wedge f$, where $f$ is a projective trivial cofibration. In that case, at a fixed object $h \in \sfG$ we have that $(F \otimes (\sfG(g,-)_+ \wedge f))(h) = F(g \oplus h) \wedge f$. Now we may apply \cite[Theorem 5.3.7]{HoveyShipleySmith} again.

\begin{proposition}Let $\cC$ be a $\Sigma$-cofibrant operad in simplicial sets, satisfying either
	\begin{enumerate}[(i)]
		\item $\cC(0) \neq \varnothing$, or 
		\item $\cC$ is non-unitary and admits a unitalization $\cC^+$. 
	\end{enumerate}
	Let $\sfC = (\cat{Sp}^\Sigma)^\sfG$, then the projective model structure on $\Alg_\cC(\sfC)$ transferred from the projective model structure on $\sfC$ transferred from absolute stable projective model structure on $\cat{Sp}^\Sigma$, exists.
\end{proposition}

\begin{remark}We believe the analogous statement holds when one replaces $\cat{Sp}^\Sigma$ with the category $R\text{-}\cat{Mod}$ of module spectra over a commutative ring spectrum $R$, given the existence of an injective model structure on $R\text{-}\cat{Mod}$ transferred from the absolute stable projective model structure on $\cat{Sp}^\Sigma$.\end{remark}

\subsubsection{Further results on the existence of the projective model structure on algebras over operads} \label{sec:further-model-results} There is a large literature on the existence of the projective model structure of $\cO$-algebras. We shall survey the main results, all of which are obtained by studying the filtration of Section \ref{sec:filt-stages-pushout} to verify the conditions in Theorem \ref{thm:general-proj-model-structure}. 

In general, for $\Sigma$-cofibrant operad $\cO$, the lifting and factorization axioms only hold when the domain is cofibrant and thus one only has a projective \emph{semi-model structure} on $\Alg_\cO(\sfC)$, see \cite[Theorem 12.3.A]{FresseBook}.

However, there are properties which guarantee the existence of the projective model structure on algebras over \emph{any} operad. In the symmetric monoidal case, such a condition is given in \cite{PavlovScholbach} (generalizing \cite{bataninberger}): all maps of the form of the left vertical morphism in (\ref{eqn:filt-stages-pushout-2}) have to be trivial $h$-cofibrations in the sense of \cite{bataninberger} (also called ``flat maps'', e.g. \cite{HillHopkinsRavenel}), which is guaranteed by requiring that the model structure on $\sfC$ is \emph{symmetric $h$-monoidal}.
Theorem 5.10 of \cite{PavlovScholbach} implies that if $\sfC$ is a symmetric monoidal model category that is symmetric $h$-monoidal, then projective model structure on $\Alg_\cO(\sfC)$ exists for any operad $\cO$ under mild additional model-categorical conditions on $\sfC$. Many model structures of interest are symmetric $h$-monoidal and satisfy the additional conditions, including $\cat{Top}$, $\cat{sSet}$, and $\cat{sMod}_\bk$. This gives an alternative approach to the existence of projective model structures in these categories. Similar results hold on weaker conditions in the monoidal setting, see  \cite{MuroI,MuroII,MuroCorrection}.

\subsection{Simplicial formulae revisited} We now discuss how some of the results in Section \ref{sec:simpl-resol} can be improved if the monad $T$ comes from an operad $\cO$.

\subsubsection{Derived indecomposables and decomposables}
\label{sec:deriv-decomp}

If $\cO$ is an augmented operad, Section \ref{sec:deriv-indec-1} describes how to derive its indecomposables and under mild conditions Section \ref{sec:simplicial-formula-indecomposables} gives a simplicial formula for $Q^\cO_\bL$.

As explained in Section \ref{sec:relative-augmentation}, there is a canonical relative augmentation $\epsilon^\cO_{\cO(1)}$ and if $\cO(1)$ is augmented, we get an augmentation of $\cO$ which factors as
\[\cO \lra \cO(1)_+ \lra (-)_+.\]
Hence the indecomposables functor $Q^\cO$ factors as the composition over a partial indecomposables functor $Q^\cO_{\cO(1)} \colon \Alg_\cO(\sfC) \to \Alg_{\cO(1)}(\sfC_*)$ to the category of pointed objects with $\cO(1)$-action. Under the assumption that projective model structures on $\Alg_\cO(\sfC)$ and $\Alg_{\cO(1)}(\sfC_*)$ exist, $Q^{\cO}_{\cO(1)}$ is a left Quillen functor.

Let us now assume that $\cO$ is $\Sigma$-cofibrant and $\gR$ is a $\cO$-algebra whose underlying object $R = U^\cO(\gR)$ is cofibrant in $\sfC$. By Lemma \ref{lem:operad-geomrel} the monad $\cO$ commutes with geometric realization and thus we can apply the discussion in Section \ref{sec:simplicial-formula-indecomposables} to derive $Q^\cO_{\cO(1)}$ using the explicit simplicial formula\index{indecomposables!derived}
\begin{equation*}
\mathbb{L}Q^\cO_{\cO(1)}(\gR) \simeq \fgr{B_\bullet(Q^\cO_{\cO(1)}F^\cO,\cO,\gR)} \cong \fgr{B_\bullet(\cO(1)_+,\cO,\gR)}.
\end{equation*}

In Section \ref{sec:relative-decomposables} we described $Q^\cO_{\cO(1)}$ as the cofiber of a natural transformation 
\[\Dec^\cO_{\cO(1)} \Rightarrow U^{\cO}_{\cO(1)} \circ (-)_+ \colon \Alg_{\cO}(\sfC) \lra \Alg_{\cO(1)}(\sfC_\ast).\]
The relative decomposables functor $\Dec^\cO_{\cO(1)}$ as defined in Definition \ref{def:decomposables} is not necessarily a left Quillen functor because it is not necessarily a left adjoint, being the composition of the left adjoint $(-1)_*^\alg$ and the right adjoint $(-2)^* U^\cO$. However, we can still left derive it as long as it preserves trivial cofibrations between cofibrant objects. For this it suffices that it sends generating (trivial) cofibrations $F^\cO(f)$ to (trivial) cofibrations. To prove this, recall that on free algebras $\Dec^\cO$ is given by the formula
\begin{equation*}
\Dec^\cO_{\cO(1)}(F^\cO X) \cong \left(\bigsqcup_{n \geq 2} \cO(n)
\times_{G_n} X^{\otimes n}\right)_+ \eqqcolon \cO^{\geq 2}(X)_+,
\end{equation*}
As $\cO$ is $\Sigma$-cofibrant, $\cO^{\geq 2}$ is also cofibrant in $\cat{FB}_k(\sfC)$. It then follows from Lemma \ref{lem.symmetric-sequence-preserving} (i) that the functor $X \mapsto \cO^{\geq 2}(X)_+$ preserves trivial cofibrations between cofibrant objects. We conclude that we may left derive the functor $\Dec^\cO_{\cO(1)}$.

\index{decomposables!derived}We can resolve $\gR$ by the thick monadic bar resolution, and compute the derived functor $\mathbb{L}\Dec^\cO_{\cO(1)} \colon \Alg_\cO(\sfC) \to \Alg_{\cO(1)}(\sfC_+)$ by the explicit simplicial formula
\begin{equation*}
\mathbb{L}\Dec^\cO_{\cO(1)}(\gR) \coloneqq \fgr{B_\bullet(\Dec^\cO F^\cO,\cO,\gR)} \cong \fgr{B_\bullet(\cO^{\geq 2}_+,\cO,\gR)}.\end{equation*}
As the underlying object $R$ of $\gR$ is cofibrant, and $X \mapsto \cO(X)$ preserves cofibrant objects by Lemma \ref{lem.symmetric-sequence-preserving} (i), it follows that $\cO^p(R)$ is cofibrant in $\sfC$. Thus
\[\cO^{\geq 2}(\cO^p(R)) \lra \cO(\cO^p(R))\]
is a cofibration by Lemma \ref{lem.symmetric-sequence-preserving} (ii), and so $B_\bullet(\cO^{\geq 2}_+,\cO,\gR) \to B_\bullet(\cO_+,\cO,\gR)$ is a levelwise cofibration, so by Lemma \ref{lem:thick-geom-rel-cofibrations} we obtain a cofibration sequence
\begin{equation}\label{eq:DecCofSeq}
\mathbb{L}\Dec^\cO_{\cO(1)}(\gR) \lra U^\cO_{\cO(1)}\fgr{B_\bullet(F^\cO,\cO,\gR)}_+ \lra \mathbb{L}Q^\cO_{\cO(1)}\gR
\end{equation}
in $\Alg_{\cO(1)}(\sfC_\ast)$.

\subsubsection{Reedy cofibrancy and operads in simplicial sets}

For later use, we establish the following lemma about Reedy cofibrancy of two-sided bar constructions.

\begin{lemma}\label{lem:reedy-cofibrancy-simplicial} Let $\cC$ be a $\Sigma$-cofibrant operad in simplicial sets, $F \colon \sfC \to \sfD$ be a right $\cC$-module functor satisfying one of the following properties:
	\begin{enumerate}[(i)]
		\item $F$ preserves colimits and (trivial) cofibrations between cofibrant objects.
		\item $F = \cX(-)$ where $\cX \in \cat{FB}_k(\cat{sSet})$ is a $\Sigma$-cofibrant right $\cO$-module in $k$-symmetric sequences.
	\end{enumerate} 
	Then $B_\bullet(F,\cC,\gR)$ is Reedy cofibrant if $\gR$ is cofibrant in $\sfC$. More generally, $B_\bullet(F,\cC,-)$ sends (trivial) cofibrations between objects in $\Alg_\cC(\sfC)$ cofibrant in $\sfC$ to (trivial) cofibrations between Reedy cofibrant objects in $\cat{sD}_\cat{Reedy}$.\end{lemma}

\begin{proof}It suffices to prove the second statement. We recall that a map $X_\bullet \to Y_\bullet$ of simplicial objects is (trivial) Reedy cofibration if $L_n(Y_\bullet) \sqcup_{L_n(X_\bullet)} X_n \to Y_n$ is a (trivial) cofibration. We start with a discussion of the latching objects $L_n(B_\bullet(F,\cC,\gR))$ and the latching morphism $L_n(B_\bullet(F,\cC,\gR)) \to B_n(F,\cC,\gR)$. For $i \in \{0,1\}$, let $\cC^i$ denote the identity functor if $i=0$ and the functor $\cC$ if $i=1$. If we let $[1]$ denote the category $0 \to 1$, then the $n$th latching object $L_n$ is given by the colimit over the punctured cubical diagram 	\[[1]^n \setminus \{1,\ldots,1\} \ni I \longmapsto F(\cC^{I}(U^\cC(\gR))) \]
	where $\cC^I \coloneqq \cC^{i_1} \circ \cdots \circ \cC^{i_n}$. The $n$th latching map from $L_n(B_\bullet(F,\cC,\gR))$ to $B_n(F,\cC,\gR)$ is induced by adding the missing corner, which is final. 
	
	In case (i) it suffices to prove that 
	\[\mr{colim}_{I} \cC^I(\gS) \cup_{\mr{colim}_{I} \cC^I(\gR)} \cC^n(\gR) \lra \cC^n(\gS)\]
	is a (trivial) cofibration if the map $\gR \to \gS$ is a (trivial) cofibration and $\gR$ is cofibrant in $\sfC$. This uses the fact that $F$ commutes with colimits and preserves (trivial) cofibrations between cofibrant objects.
	
	The map $\mr{colim}_{I} \cC^I \to \cC^n$ of $k$-symmetric sequences is an inclusion with codomain a $k$-symmetric sequence of simplicial sets with levelwise free action, because $\cC$ is $\Sigma$-cofibrant, and hence a cofibration. Using Lemma \ref{lem.symmetric-sequence-preserving}, we conclude that in
	\[\begin{tikzcd}(\mr{colim}_{I} \cC^I)(\gR) \rar \dar & \cC^n(\gR) \dar \arrow[bend left = 25]{rdd} &  \\
	(\mr{colim}_{I} \cC^I)(\gS) \arrow[bend left = -15]{rrd} \rar & (\mr{colim}_{I} \cC^I)(\gS) \sqcup_{(\mr{colim}_{I} \cC^I)(\gR)} \cC^n(\gR) \arrow[dotted]{rd} \\
	 & & \cC^n(\gS)
	\end{tikzcd}\]
	the dotted map is a (trivial) cofibration if $\gR \to \gS$ is a (trivial) cofibration and $\gR$ is cofibrant in $\sfC$.
	
	In case (ii), one repeats the above argument to prove that
		\[\mr{colim}_{I} \cX(\cC^I) \lra \cX(\cC^n)\]
	is a cofibration of $k$-symmetric sequences, and applies Lemma \ref{lem.symmetric-sequence-preserving} once more.
	\end{proof}

\subsection{Modules over associative algebras}\label{sec:AssocModules}

Let us now take $\gR$ to be a unital associative algebra (i.e.\ a monoid) in $\sfC$, and consider the associated monad $\gR \otimes -$, whose algebras are by definition left $\gR$-modules. This is the monad associated to an operad with only $1$-ary operations, given by $\gR$. This operad is $\Sigma$-cofibrant if the underlying object of $\gR$ is cofibrant in $\sfC$. We may then apply the general theory of Section \ref{sec:model-str-alg-over-operads}, which was in fact already explained in Section \ref{sec:model-module-categories}: Theorem \ref{thm:model-structure-on-modules} says that the projective model structure on $\gR\text{-}\cat{Mod}$ exists if either (i) the underlying object of $\gR$ is cofibrant in $\sfS$, or (ii) $\sfS$ satisfies Schwede--Shipley's monoid axiom. Furthermore, in case (i), $U^\gR \colon \gR\text{-}\cat{Mod} \to \sfC$ preserves (trivial) cofibrations.

Of course, all of the above may be repeated with the monad $- \otimes \gR$, describing the category $\cat{Mod}\text{-}\gR \coloneqq \Alg_{- \otimes \gR}(\sfC)$ of right $\gR$-modules and its model structure. This can be interpreted in terms of the operad with only $1$-ary operations, given by $\gR^\mr{op}$.

\subsubsection{Tensor product of modules}\label{sec:TensorOfModules}

There is a functor $- \otimes_\gR - \colon \cat{Mod}\text{-}\gR  \times \gR\text{-}\cat{Mod} \to \sfC$ defined by the reflexive coequalizer\index{tensor product!of modules}
	\begin{equation*}
	\begin{tikzcd}
	\gN \otimes \gR \otimes \gM \arrow[shift left=1ex]{r} \arrow[shift left=-1ex]{r} &	\gN \otimes \gM \lar\arrow{r} & \gN \otimes_\gR \gM.
	\end{tikzcd}
	\end{equation*}
where the maps are given by the left and right $\gR$-actions and the reflection is given by the unit of $\gR$. There are natural isomorphisms
\[(X \otimes \gR) \otimes_\gR \gM \cong X \otimes \gM, \quad\quad \gN \otimes_{\gR} (\gR \otimes X) \cong \gN \otimes X.\]

\begin{lemma}\label{lem:otimes-gr-quillen-bifunctor} If $\gR$ is cofibrant in $\sfC$, then $- \otimes_\gR - \colon \cat{Mod}\text{-}\gR  \times \gR\text{-}\cat{Mod} \to \sfC$ is a Quillen bifunctor.
\end{lemma}
\begin{proof}
Similarly to Section \ref{sec:module-cats}, one proves that $-\otimes_\gR -$ participates in an adjunction of two variables. We can then use \cite[Corollary 4.2.5]{Hovey}, which says that it suffices to check the properties of a Quillen bifunctor for $\otimes_\gR$ only on generating (trivial) cofibrations. The generating (trivial) cofibrations in $\cat{Mod}\text{-}\gR$ are $f \otimes \gR$ for $f$ generating (trivial) cofibrations in $\sfC$, and similarly in $\gR\text{-}\cat{Mod}$. If $f \colon A \to B$ and $f' \colon A' \to B'$ are cofibrations in $\sfC$ then the first natural isomorphism above shows that $(X \otimes \gR) \otimes_\gR (\gR \otimes Y) \cong X \otimes R \otimes Y$, and hence the pushout-product $(f \otimes \gR) \Box_\gR (\gR \otimes f')$ is identified with the pushout-product $f \Box (\gR \otimes f')$ in $(\sfC, \otimes, \bunit)$. As $- \otimes -$ is a Quillen bifunctor and $\gR$ is cofibrant in $\sfC$, this is a cofibration if both $f$ and $f'$ are and is a trivial cofibration if in addition one of $f$ and $f'$ is a weak equivalence.
\end{proof}

We write $- \otimes^\bL_\gR -$ for the derived functor. We wish to explain how this may be computed by a two-sided bar construction under favourable circumstances. Recall that if $\gN$ is a right $\gR$-module and $\gM$ is a left $\gR$-module, then the \emph{two-sided bar construction}\index{bar construction!two-sided} $B_\bullet(\gN, \gR, \gM)$ is the simplicial object with
\[B_p(\gN, \gR, \gM) = \gN \otimes \gR^{\otimes p} \otimes \gM,\]
face maps given by the multiplication on $\gR$ and the module structures, and degeneracies given by the unit of $\gR$. We shall generally consider this as a semi-simplicial object, and write $B(\gN, \gR, \gM) \coloneqq \fgr{B_\bullet(\gN, \gR, \gM)}$ for its thick geometric realization.

\begin{lemma}\label{lem:DerivedTensorBar}
If $\gN$, $\gR$, and $\gM$ are cofibrant in $\sfC$, then there is an equivalence
\[\gN \otimes^\bL_\gR \gM \simeq B(\gN, \gR, \gM).\]
\end{lemma}
\begin{proof}
There is an augmentation $\epsilon \colon B_\bullet(\gN, \gR, \gM) \to \gN \otimes_\gR \gM$, as the coequalizer defining $\gN \otimes_\gR \gM$ is the 1-skeleton of the semi-simplicial object $\sigma^*B_\bullet(\gN, \gR, \gM)$.

Suppose first that $\gM$ is a cofibrant left $\gR$-module. If $\gN = \gX \otimes \gR$ is a free right $\gR$-module then the augmented simplicial object $\epsilon \colon B_\bullet(\gN, \gR, \gM) \to \gN \otimes_\gR \gM$ is identified with $\gX \otimes -$ applied to the augmented simplicial object $B_\bullet(\gR, \gR, \gM) \to \gM$. After applying $U^{\gR}(-)$ this has an extra degeneracy, so after applying $\gX \otimes -$ it does too, so by Lemma \ref{lem:extra-degeneracy} the map $\epsilon \colon B(\gN, \gR, \gM) \to \gN \otimes_\gR \gM$ is a weak equivalence.

Now let $\gN_\bullet = \sigma^*B_\bullet(\gN, \gR, \gR) \to \gN$, an augmented semi-simplicial object with an extra degeneracy, so a weak equivalence on thick geometric realization as above. We then have a bi-semi-simplicial object
\[([p], [q]) \longmapsto B_p(\gN_q, \gR, \gM)\]
augmented in the $q$ direction to $B_p(\gN, \gR, \gM)$, and augmented in the $p$ direction to $\gN_q \otimes_\gR \gM$. The maps $\fgr{B_\bullet(\gN_q, \gR, \gM)} \to \gN_q \otimes_\gR \gM$ are weak equivalences for each $q$ by the above, as each $\gN_q$ is a free right $\gR$-module. The augmented semi-simplicial object $B_p(\gN_\bullet, \gR, \gM) \to B_p(\gN, \gR, \gM)$ has an extra degeneracy so is a weak equivalence on thick geometric realization as above. This gives weak equivalences
\[\fgr{B_\bullet(\gN, \gR, \gM)} \overset{\sim}\longleftarrow \fgr{\fgr{B_\bullet(\gN_\bullet, \gR, \gM)}} \overset{\sim}\lra \fgr{\gN_\bullet \otimes_\gR \gM} \cong \fgr{\gN_\bullet} \otimes_\gR \gM.\]
Now, as $\gN$ is cofibrant in $\sfC$, $\gN_\bullet$ is levelwise a cofibrant right $\gR$-module, and so by Lemma \ref{lem:thick-geom-rel-cofibrations}, $\fgr{\gN_\bullet}$ is a cofibrant right $\gR$-module. As we have supposed that $\gM$ is a cofibrant left $\gR$-module, the rightmost term is a model for $\gN \otimes^\bL_\gR \gM$, which proves the lemma under the assumption that $\gM$ is a cofibrant left $\gR$-module.

If $c\gM \overset{\sim}\to \gM$ is a cofibrant approximation as a left $\gR$-module, then
\[B_\bullet(\gN, \gR, c\gM) \lra B_\bullet(\gN, \gR, \gM)\]
is a levelwise weak equivalence as $\gN$ and $\gR$ are cofibrant in $\sfC$, so if $\gM$ is also cofibrant in $\sfC$ then both objects are levelwise cofibrant and so this map is a weak equivalence on thick geometric realization by Lemma \ref{lem:thick-geom-rel-cofibrations}. This proves the lemma in general.
\end{proof}

\subsubsection{Derived indecomposables}\label{sec:derived-indec-modules}  \index{indecomposables!derived} Let us now suppose that there is given an augmentation $\epsilon \colon \gR \to \bunit_\sfC$. This gives a map of monads $\epsilon \colon \gR \otimes - \to \bunit_\sfC \otimes - = \mr{Id}$, and hence there is defined a relative indecomposables functor
\[Q^{\gR}_\mr{Id} \colon \Alg_{\gR \otimes -}(\sfC) \lra \sfC.\]

If $\sfC$ is pointed then $\mr{Id}=+$, and so $\epsilon$ defines an augmentation in the sense of Definition \ref{def:augmentation}.\index{indecomposables!derived} In this case the relative indecomposables $Q^{\gR}_\mr{Id}$ are simply the indecomposables $Q^{\gR}$. \emph{By abuse of notation, we shall write $Q^{\gR}$ for $Q^{\gR}_\mr{Id}$, even if $\sfC$ is not pointed.}

With the model structures discussed above we may form the derived indecomposables functor $Q_\bL^\gR$. This derived functor can often be computed by a bar construction, following Section \ref{sec:simpl-resol}.

\begin{corollary}\label{cor:DerivModIndecBar}
For a left $\gR$-module $\gM$, if $\gR$ and $\gM$ are cofibrant in $\sfC$ then there is an equivalence
\[Q_\bL^\gR(\gM) \simeq B(\bunit, \gR, \gM).\]
\end{corollary}

\begin{proof}
The monadic resolution of a left $\gR$-module $\gM$ is $B_\bullet(\gR, \gR, \gM) \to \gM$, and by Lemma \ref{lem:adding-degeneracies-reedy} if $\gM$ and $\gR$ are cofibrant in $\sfC$ then $\sigma_*\sigma^*B_\bullet(\gR, \gR, \gM)$ is Reedy cofibrant. This is levelwise a free left $\gR$-module, and the underlying simplicial object has an extra degeneracy so $\gr{\sigma_*\sigma^*B_\bullet(\gR, \gR, \gM)} \to \gM$ is a weak equivalence by Lemma \ref{lem:extra-degeneracy}. By Lemma \ref{lem:monad-functors-geomrel} (iii), geometric realization in $\Alg_{\gR \otimes -}(\sfC)$ has the same underlying object as geometric realization in $\sfC$, so $\gr{\sigma_*\sigma^*B_\bullet(\gR, \gR, \gM)}_{\gR} \to \gM$ is also a weak equivalence, and hence this simplicial object is a free simplicial resolution. Then as in Section \ref{sec:simplicial-formula-indecomposables} we get $Q_\bL^\gR(\gM) \simeq \gr{\sigma_*\sigma^*B_\bullet(\bunit, \gR, \gM)} = B(\bunit, \gR, \gM)$.
\end{proof}

Of course, the entire discussion above goes through for right modules, giving an equivalence $Q_\bL^\gR(\gN) \simeq B(\gN, \gR, \bunit)$ under the same hypotheses.

\section{Homology and spectral sequences}\label{sec:homology+ss}

In this section we shall discuss the filtrations of Section \ref{sec:filtered-objects-1} from a homotopical point of view. Of particular interest is the construction of their attendant spectral sequences, and this requires a discussion of homology in our contexts. We suppose throughout that $\sfS$, and hence $\sfC= \sfS^\sfG$, satisfies the assumptions of of Section \ref{sec:axioms-of-model-cats}. If there is a monad $T$ discussed, or the monad associated to an operad $\cO$, then we suppose that monad satisfies the axioms of Section \ref{sec:homotopy-theory-algebras}.

\subsection{Homology}\label{sec:homology} 

We shall discuss homology with coefficients in a $\bk$-module, for a commutative ring $\bk$.

\subsubsection{Singular chain functors}
\newglossaryentry{cata}{%
	name={\ensuremath{\cat{A}}},
	description={Target of singular chain functor},
type=symbols
}
\newglossaryentry{cast}{%
	name={\ensuremath{C_*}},
	description={Singular chain functor},
type=symbols
}
We shall define the homology of either chain complexes of $\bk$-modules or $H\bk$-modules, and so we write
\[\gls{cata} = \begin{cases}
\cat{Ch}_\bk, & \\
H\bk\text{-}\cat{Mod},
\end{cases}\]
with their symmetric monoidal projective model category structures. Here $H\bk$ is an Eilenberg--Maclane spectrum that is a commutative ring spectrum in $\cat{Sp}^\Sigma$ as in Example \ref{exam:hk-cofibrant}, and such a model structure exists by the discussion in Section \ref{sec:model-module-categories}. We write $H_i(X)$ for the homology of a chain complex $X$, and $H_i(X) \coloneqq \pi_i(X)$ for the ``homology" of an $H\bk$-module $X$.

In order to discuss homology of objects of $\sfS$, we shall ask for a \emph{singular chain functor}\index{singular chain functor}
\[\gls{cast} \colon \sfS \lra \cat{A},\]
by which we mean a functor such that
\begin{enumerate}[(i)]
\item there is a lax symmetric monoidal structure
\[C_*(X) \otimes C_*(Y) \lra C_*(X \otimes Y)\]
which is a weak equivalence when $X$ and $Y$ are cofibrant,

\item the composition $C_* \circ s \colon  \cat{sSet} \to\sfS \to \cat{A}$ is naturally weakly equivalent to either $C_*(-;\bk)$ or $H\bk \wedge \Sigma^\infty (-)_+$,

\item $C_*$ preserves cofibrant objects and weak equivalences between cofibrant objects, and preserves homotopy colimits (as described in Remark \ref{rem:hocolim}).
\end{enumerate}

Given such a functor, we define the associated \emph{reduced singular chain functor} $\widetilde{C}_* \colon \sfS_* \to \cat{A}$ as $\widetilde{C}_*(X) = C_*(X)/C_*(\bterm)$, with homology $\widetilde{H}_*(X)$. As $C_*(\bterm) \to C_*(X)$ has a retraction $C_*(X \to \bterm)$, it is a cofibration, and so
\[C_*(\bterm) \lra C_*(X) \lra C_*(X)/C_*(\bterm) = \widetilde{C}_*(X)\]
is a split cofibration sequence, and hence there is a canonical decomposition
\[H_*(C_*(X)) \cong {H}_*(\widetilde{C}_*(X)) \oplus H_*(C_*(\bterm)).\]
Note that $\binit = s(\varnothing)$ so by (ii) we have $C_*(\binit) \simeq 0$. Thus if the category $\sfS$ is pointed then reduced and unreduced singular chains agree.

\begin{lemma}
If ${C}_* \colon \sfS \to \cat{A}$ is a singular chain functor and $\bterm \in \sfS$ is cofibrant then $\widetilde{C}_* \colon \sfS_* \to \cat{A}$ is a singular chain functor.
\end{lemma}
\begin{proof}
For $X, Y \in \sfS_*$ the map
\[C_*(X) \otimes C_*(Y) \lra C_*(X \otimes Y) \lra C_*(X \owedge Y) \lra \widetilde{C}_*(X \owedge Y)\]
is trivial when restricted to $C_*(X) \otimes C_*(\bterm)$ and $C_*(\bterm) \otimes C_*(Y)$, so descends to a map $\widetilde{C}_*(X) \otimes \widetilde{C}_*(Y) \lra \widetilde{C}_*(X \owedge Y)$ which is a lax symmetric monoidality. One may verify that it is a weak equivalence when $X$ and $Y$ are cofibrant, using the canonical decomposition above, verifying (i).

The composition $\widetilde{C}_* \circ s \colon \cat{sSet} \to \sfS_* \to \cat{A}$ sends $X$ to $\widetilde{C}_*(s(X) \sqcup \bterm) \cong C_*(s(X))$, verifying (ii). 

If we suppose that $\bterm \in \sfS$ is cofibrant, then $U^+ \colon \sfS_* \to \sfS$ preserves cofibrant objects. Thus if $\bterm \to X \overset{f}\to Y$ is a weak equivalence between cofibrant objects of $\sfS_*$ then $f$ is a weak equivalence between cofibrant objects in $\sfS$, and so $C_*(f)$ is a weak equivalence between cofibrant objects. By the functorial split cofibration sequence above, $\widetilde{C}_*(f)$ is a weak equivalence too. If $F \colon \cat{I} \to \sfS_*$ is a diagram, its homotopy colimit may be formed as the homotopy cofibre of
\[\hocolim_\cat{i \in I} \bterm \lra \hocolim_\cat{i \in I} U^+F(i).\]
This identifies $\widetilde{C}_*(\hocolim_\cat{i \in I} F(i))$ with the homotopy cofibre of
\[\hocolim_\cat{i \in I} C_*(\bterm) \lra \hocolim_\cat{i \in I} C_*(F(i)),\]
which is $\hocolim_\cat{i \in I} \widetilde{C}_*(F(i))$. This verifies (iii).
\end{proof}

It is also easy to see that given a singular chain functor $\widetilde{C}_* \colon \sfS_* \to \cat{A}$, precomposing with $F^+ \colon \sfS \to \sfS_*$ defines a singular chain functor on $\sfS$.

\subsubsection{Examples}

For $\sfS = (\cat{sMod}_\bk, \otimes, \bk)$, we take $\cat{A} = \cat{Ch}_\bk$ and let $C_*(X) \coloneqq N(X)$ be normalized chains. The Eilenberg--Zilber map gives the required natural weak equivalence $N(X) \otimes N(Y) \to N(X \otimes Y)$. The normalized chains functor $N \colon \cat{sMod}_\bk \to \cat{Ch}_\bk$ is a left Quillen functor \cite[\S 4.1]{SchwedeShipley} so satisfies (iii), and the composition $C_* \circ s \colon \cat{sSet} \to \cat{Ch}_\bk$ is the functor of normalized simplicial chains which is naturally weakly equivalent to the functor $C_*(-;\bk)$ of simplicial chains.

For $\sfS = (\cat{sSet}_*, \wedge, S^0)$, we take $\cat{A} = \cat{Ch}_\bk$ and use the strong symmetric monoidal ``reduced free $\bk$-module" functor $\tilde{\bk}[-] \coloneqq \bk[-]/\bk[*] \colon \cat{sSet}_* \to \cat{sMod}_\bk$, composed with the construction for simplicial $\bk$-modules above. This satisfies (iii) as $\tilde{\bk}[-]$ is a left Quillen functor. For $(\cat{sSet}, \times, *)$ we precompose this with the strong symmetric monoidal left Quillen functor $F^+ \colon \cat{sSet} \to \cat{sSet}_*$.

For $\sfS = (\cat{Top}, \times, *)$ or $(\cat{Top}_*, \wedge, S^0)$, we take $\cat{A} = \cat{Ch}_\bk$ and use the strong symmetric monoidal singular simplices functor $\Sing$ composed with the construction above for (pointed) simplicial sets. The composition $C_* \circ s$ sends a simplicial set $K$ to $C_*(\Sing\vert K \vert;\bk)$ which has a canonical weak equivalence from $C_*( K ;\bk)$, so this satisfies (ii). For (iii), first note that the functor $\Sing$ always produces cofibrant objects, and preserves weak equivalences between all (not just cofibrant) objects. Thus it has a left derived functor $\bL\Sing$ (even though $\Sing$ is a \emph{right} Quillen functor) given by $(\bL \Sing)(X) = \Sing(cX)$ for $c$ a cofibrant replacement functor (which we may take to be $cX = \vert \Sing(X)\vert$). As $\Sing$ participates in a Quillen equivalence, it is enough to show that its adjoint, $\gr{-} \colon \cat{sSet} \to \cat{Top}$ preserves homotopy colimits, but this is clear as it is a left Quillen functor.

For $\sfS = (\cat{Sp}^\Sigma, \wedge, S^0)$ we take $\cat{A} = H\bk\text{-}\cat{Mod}$ and let $C_*(X) = H\bk \wedge X$. This is strong symmetric monoidal and precomposed with $s \colon \cat{sSet} \to \cat{Sp}^\Sigma$ it is equal to $H\bk \wedge \Sigma^\infty (-)_+$. It is a left Quillen functor by definition of the projective model structure on $H\bk\text{-}\cat{Mod}$, so satisfies (iii).

\subsubsection{Homology of objects of $\mathsf{S}$}\label{sec:homology-objects-sfs}

Given such a singular chain functor $C_* \colon \sfS \to \cat{A}$, we define the \emph{homology}\index{homology} $H_*(X;\bk)$ of $X \in \sfS$ as the homology of $\bL C_*(X)$. More generally, if $A$ is a left $\bk$-module then we define 
\[C_*(X;A) \coloneqq \begin{cases}
\bL C_*(X) \otimes_\bk^\bL A & \text{ if $\cat{A} = \cat{Ch}_\bk$,}\\
\bL C_*(X) \otimes_{H\bk}^\bL HA & \text{ if $\cat{A} = H\bk\text{-}\cat{Mod}$,}
\end{cases}\]
with homology $H_*(X;A)$, where $HA$ is the Eilenberg--Maclane spectrum for $A$ defined by replacing $\bk$ with $A$ in Example \ref{exam:hk-cofibrant}.

If $f \colon X \to Y$ is a morphism in $\sfS$, then we define the relative chains $C_*(f;A)$ as the mapping cone in $\cat{A}$ of $f_* \colon C_*(X;A) \to C_*(Y;A)$, whose homology we denote by $H_*(f;A)$. There is then a long exact sequence
\[\cdots \lra H_i(X;A) \overset{f_*}\lra H_i(Y;A) \lra H_i(f;A) \overset{\partial}\lra H_{i-1}(X;A) \lra \cdots.\]
If the map $f$ is understood, we write $H_*(Y, X;A)$ for $H_*(f;A)$.

If $f \colon X \to Y$ is a cofibration in $\sfS$ with cofibre $Y/X = Y \cup_X \bterm$, then the square
\begin{equation*}
\begin{tikzcd}
	C_*(X;A) \dar \rar{f_*}& C_*(Y;A) \dar\\
	C_*(\bterm;A) \rar  & C_*(Y/X;A)
\end{tikzcd}
\end{equation*}
is a homotopy pushout, giving an equivalence $C_*(Y,X;A) \simeq \widetilde{C}_*(Y/X;A)$, where the homology of the latter is denoted $\widetilde{H}_*(Y/X;A)$. Using this notation, 
there is a long exact sequence
\[\cdots \lra H_i(X;A) \overset{f_*}\lra H_i(Y;A) \lra \widetilde{H}_i(Y/X;A) \overset{\partial}\lra H_{i-1}(X;A) \lra \cdots.\]

\subsubsection{Homology of objects of $\mathsf{C}$}\label{sec:homology-objects-sfc}

For objects $X \in \sfC = \sfS^\sfG$ we consider homology to give $\sfG$-graded $\bk$-modules, as follows.
\newglossaryentry{hgd}{%
	name={\ensuremath{H_{g,d}}},
	description={Homology in degree $d$ at $g \in \sfG$},
type=symbols
}
\begin{definition}
	Let $\bk$ be a commutative ring, $A$ be a $\bk$-module, and $X \in \sfC = \sfS^\sfG$. The \emph{homology groups of $X$ with coefficients in $A$}\index{homology!with coefficients} are defined to be the $\bk$-module
		\begin{equation*}
		H_{g,d}(X;A) \coloneqq H_d(X(g);A).
		\end{equation*}
	We consider the collection of these groups as giving a functor 
\begin{align*}
\sfG &\lra \cat{Mod}_\bk^\bZ\\
g &\longmapsto H_{g,*}(X;A) 
\end{align*}
to $\bZ$-graded $\bk$-modules, obtained as the homology of $g \mapsto C_*(X(g);A) \colon \sfG \to \cat{A}$.
\end{definition}

As usual, it is convenient to also have available relative homology. For a morphism $f \colon X \to Y$ in $\sfC = \sfS^\sfG$ the \emph{relative	homology groups}\index{homology!relative} $H_{g,d}(f;A)$ are defined as the homology groups of the mapping cones of $f_* \colon C_*(X(g);A) \to C_*(Y(g);A)$. 
As usual, we shall write these groups as $H_{g,d}(Y,X;A)  \coloneqq H_{g,d}(f;A)$ when $f$ is understood. As in the previous section, if $f \colon X \to Y$ is a cofibration in $\sfC$ with cofibre $Y/X$ then there is an identification $\widetilde{H}_{g,d}(Y/X;A) \cong H_{g,d}(f;A)$, and an associated long exact sequence.

\begin{remark}
There is a more general notion of coefficients for objects of $\sfS^\sfG$, which we will not have need for but which some readers may find clarifying. If $\cA \colon \sfG^\mathrm{op} \to \cat{Mod}_\bk$ is a functor then we can define $C_*(X;\cA)$ to be the homotopy coend of the functor
		\begin{align*}
		\sfG \times \sfG^\op & \lra \cat{A}\\
		(g,g') &\longmapsto C_*(X(g);\cA(g')).
		\end{align*}
One can then define homology of $X$ with coefficients in $\cA$ as the $\bk$-module
		\begin{equation*}
		H_d(X;\cA) \coloneqq H_d(C_*(X;\cA)).
		\end{equation*}
 
Let $A$ be a $\bk$-module and $\cA(g) = A \otimes_\Z \Z[\sfG(-,g)]$, the representable functor
$\sfG(-,g) \colon \sfG^\op \to \cat{Set}$ composed with the free $\bZ$-module
functor, composed with $A \otimes_\Z- \colon \cat{Ab} \to \cat{Mod}_\bk$. Then $C_*(X;\cA)$ is equivalent to $C_*(X(g);A)$, i.e.\ ordinary chains of the
object $X(g) \in \sfS$ with coefficients in $A$.  Hence the two notions of homology are related
by a natural isomorphism $H_{g,d}(X;A) \cong H_d(X; \cA(g))$.
\end{remark}

\subsubsection{K{\"u}nneth theorems}\label{sec:Kunneth}\index{K{\"u}nneth theorem}

For cofibrant objects $A, B \in \sfC = \sfS^\sfG$ and any commutative ring $\bk$ of coefficients, the lax monoidal structure on $C_* \colon \sfS \to \cat{A}$ gives weak equivalences
\[C_*(A) \otimes C_*(B ) \lra C_*(A \otimes B)\]
in $\cat{A}^\sfG$. Thus for the purposes of establishing K{\"u}nneth-type theorems, we may as well work entirely in $\cat{A}^\sfG$. For $U, V \in \cat{A}^\sfG$ we have, by definition of the Day convolution product,
\[(U \otimes V)(x) = \colim_{(a, b, f) \in \sfH_x} U(a) \otimes V(b)\]
where the category $\sfH_x$ has objects triples $(a, b, f)$ with $a, b \in \sfG$ and $f \colon a \oplus b \to x$ a morphism in $\sfG$, and $\sfH_x((a', b', f'),(a, b, f))$ is given by morphisms $g \colon a' \to a$ and $h \colon b' \to b$ in $\sfG$ such that $f' = f \circ (g \oplus h)$.

\begin{lemma}\label{lem:KanExtCofibrant}
If $U, V \in \cat{A}^\sfG$ are cofibrant then the functor
\[(a, b, f) \mapsto U(a) \otimes V(b) \colon \sfH_x \lra \cat{A}\]
is cofibrant. Thus there is a strongly convergent spectral sequence
\[E^2_{x, s,t} = \bL_s\colim_{(a, b, f) \in \sfH_x} H_t(U(a) \otimes V(b)) \Longrightarrow H_{x, s+t}(U \otimes V)\]
with differentials $d^r \colon E^r_{x,s,t} \to E^r_{x,s-r,t+r-1}$.
\end{lemma}

\begin{proof}
It suffices to prove that the functor
	\[-\ul{\otimes}- \colon \cat{A}^\sfG \times \cat{A}^\sfG \lra \cat{A}^{\sfH_x},\]
given by exterior product to $\cat{A}^{\sfG \times \sfG}$ followed by restriction along the functor $(a,b,f) \mapsto (a,b) \colon \sfH_x \to \sfG \times \sfG$, is a Quillen bifunctor. 

By Corollary 4.2.5 of \cite{Hovey} it suffices to verify that $-\ul{\otimes}-$ satisfies the property of a Quillen bifunctor only on generating (trivial) cofibrations. If $f_i \colon A_i \to B_i$ are the generating (trivial) cofibrations in $\cat{A}$, then the morphisms $\sfG(z, -) \times f_i$ are the generating (trivial) cofibrations of $\cat{A}^\sfG$. The pushout-product $(\sfG(y, -) \times f_1) \square (\sfG(z, -) \times f_2)$ evaluated at $(a,b,f)$ is identified with the map
\[\sfG(y,a) \times \sfG(z,b) \times \left(A_1 \otimes B_2 \bigsqcup_{A_1 \otimes A_2} B_1 \otimes A_2\right) \lra \sfG(y,a) \times \sfG(z,b) \times (B_1 \otimes B_2)\]
which is the identity on the first two factors and the pushout-product in $\cat{A}$ on the second. The pushout-product is a cofibration in $\cat{A}$ and is a trivial cofibration if $f_1$ or $f_2$ is, as $-\otimes-$ on $\cat{A}$ is a Quillen bifunctor. On the other hand the functor
\begin{align*} \sfH_x &\lra \cat{Set} \\
(a,b,f) &\longmapsto \sfG(y,a) \times \sfG(z,b)\end{align*}
is naturally isomorphic to the coproduct of representable functors 
\[(a,b,f) \longmapsto \bigsqcup_{g \colon y \oplus z \to x} \sfH_x((y,z,g),(a,b,f)),\]
so $(\sfG(y, -) \times f_1) \square (\sfG(z, -) \times f_2)$ is a cofibration in $\cat{A}^{\sfH_x}$, trivial if either $f_i$ is, as required.
	
Since the functor $(a, b, f) \mapsto U(a) \otimes V(b) \colon \sfH_x \to \cat{A}$ is cofibrant, its colimit is also a homotopy colimit, and the Bousfield--Kan spectral sequence for a homotopy colimit is given in \cite[Section XII.5.7]{BousfieldKan}. As it is in particular the spectral sequence of a simplicial object, it is a half-plane spectral sequence with exiting differentials and $A^\infty=0$ in the sense of Boardman, so converges strongly by \cite[Theorem 6.1]{Boardman}.
\end{proof}	

This result will typically be used in conjunction with a K{\"u}nneth theorem---or K{\"u}nneth spectral sequence---in $\cat{A}$, which may be described as follows.

\begin{lemma}\label{lem:KunnethSS}
If $U, V \in \cat{A}^\sfG$ are cofibrant then for each $a, b \in \sfG$ there is a strongly convergent K{\"u}nneth spectral sequence
\[E^2_{p,q} = \bigoplus_{q'+q''=q}\mathrm{Tor}_p^{\bk}(H_{q'}(U(a)), H_{q''}(V(b))) \Longrightarrow H_{p+q}(U(a) \otimes V(b))\]
with differentials $d^r \colon E^r_{p,q} \to E^r_{p-r,q+r-1}$. The edge homomorphism gives an external product map
\[H_*(U(a)) \otimes H_*(V(b))\lra H_*(U(a) \otimes V(b)),\]
which is an isomorphism if $H_*(U(a))$ is a flat $\bk$-module.
\end{lemma}

\begin{proof}
As $U$ is cofibrant in $\cat{A}^\sfG$ each $U(a)$ is cofibrant in $\cat{A}$. If $\cat{A} = \cat{Ch}_\bk$ then this means it is a chain complex of projective $\bk$-modules, so by e.g.\ \cite[Ex.\ 5.7.5]{Weibel} there is a K{\"u}nneth spectral sequence
\[E^2_{p,q} = \bigoplus_{q'+q''=q}\mathrm{Tor}_p^{\bk}(H_{q'}(U(a)), H_{q''}(V(b))) \Longrightarrow H_{p+q}(U(a) \otimes V(b)),\]
where the abutment is as claimed because $U(a)$ is a chain complex of flat modules, so hyper-Tor is simply given by tensor product in this case.  It is a half-plane spectral sequence with exiting differentials and $A^\infty=0$ in the sense of Boardman, so converges strongly by \cite[Theorem 6.1]{Boardman}.

If $\cat{A} = H\bk\text{-}\cat{Mod}$ there is completely analogous spectral sequence, yielding the same conclusion. This spectral sequence is developed in \cite[IV Theorem 4.1]{EKMM} (for a different model of spectra, but a similar analysis gives it in $\cat{Sp}^\Sigma$).

Under the further assumption that $H_*(U(a))$ is a flat $\bk$-module, we have $E^2_{p,q}=0$ for $p>0$ so the spectral sequence degenerates to give the claimed isomorphism.
\end{proof}

The previous two results will be used in Lemma \ref{lem:connectivity-and-tensor-products} to estimate the homological connectivity of $U \otimes V$ in terms of those of $U$ and $V$, and similarly for maps. Here we give a finer result than that, a K{\"u}nneth isomorphism in $\cat{A}^\sfG$, valid when $\sfG$ is a groupoid satisfying a mild hypothesis. We use the notation $G_x \coloneqq \mr{Aut}_{\sfG}(x) = \sfG(x,x)$ for the automorphisms of an object $x$.

\begin{lemma}\label{lem:KunnethFormula}
If $U, V \in \cat{A}^\sfG$ are cofibrant and either
\begin{enumerate}[(i)]
\item $\sfG$ is a groupoid such that the map $- \otimes - \colon G_x \times G_y \to G_{x \oplus y}$ is injective for all $x, y \in \sfG$, or

\item $U(a)=\binit$ for $a \not\cong \bunit_\sfG$,
\end{enumerate}
then
\[H_{x,*}(U \otimes V) \cong \colim_{(a, b, f) \in \sfH_x} H_*(U(a) \otimes V(b)).\]

The edge homomorphism of Lemma \ref{lem:KunnethSS} then defines an external product map
\[H_*(U) \otimes H_*(V) \lra H_*(U \otimes V),\]
which is an isomorphism if $H_*(U(a))$ is a flat $\bk$-module for all $a \in \sfG$.
\end{lemma}
\begin{proof}
Under assumption (i) on $\sfG$, the category $\sfH_x$ is filtered (because it is equivalent to a discrete category) so, as taking homology in $\cat{A}$ commutes with filtered colimits, the claimed formula holds. 
Under assumption (ii) on $U$ we have $(U \otimes V)(x) \cong U(\bunit_\sfG) \otimes V(x)$, so the same formula holds. The second part follows from Lemma \ref{lem:KunnethSS}.
\end{proof}

If the $H_*(U(a))$ are not flat, but $\bk$ is a PID, then the above can also be used to develop a K{\"u}nneth short exact sequence for $H_{*,*}(U \otimes V)$ involving $\mathrm{Tor}^\bk_1$.

\subsubsection{$T$-homology}\label{sec:t-homology}

Let $\sfS$ satisfy the axioms of Section \ref{sec:axioms-of-model-cats} and let $T$ be an augmented monad on $\sfC = \sfS^\sfG$ satisfying the Axioms of Section \ref{sec:homotopy-theory-algebras}. In Section \ref{sec:deriv-indec-1} we defined the derived $T$-indecomposables $\indec_\mathbb{L}^T \gX \in \sfC_*$ for a $\gX \in \Alg_T(\sfC)$. The $T$-homology is given by the homology of the $T$-indecomposables:

\newglossaryentry{thomology}{%
	name={\ensuremath{H^T_{g,d}}},
	description={$T$-homology groups},
type=symbols
}
\begin{definition}\label{defn:relative-T-homology}
	Let $f \colon \gX \to \gY$ be a morphism in $\Alg_T(\sfC)$. We define the \emph{$T$-homology groups}\index{homology!$T$-} to be
	\begin{equation*}
	\gls{thomology}(\gY,\gX;A) \coloneqq H_{g,d}(Q^T_\mathbb{L}\gY,Q^T_\mathbb{L} \gX;A).
	\end{equation*}

	The absolute $T$-homology is then defined to be $T$-homology relative to the
	initial $T$-algebra $F^T(\binit)$, whose underlying object is $\binit$. As $Q^T_\bL(F^T(\binit)) \simeq \binit_+ = \bterm$, this means that
\[H^T_{g,d}(\gX;A) \coloneqq \widetilde{H}_{g,d}(Q^T_\mathbb{L} \gX;A).\]
\end{definition}

	Of course, the absolute and relative $T$-homology groups fit into a long exact sequence
\begin{equation}\label{eqn:les-t-homology} 
  \cdots \to H_{g,d}^T(\gX;A) \to H_{g,d}^T(\gY;A) \to H^T_{g,d}(\gY,\gX;A) \to H^T_{g,d-1}(\gX;A) \to \cdots.
\end{equation}
If $\bk$ is a field then we can extract numerical invariants of $f$, the \emph{relative $T$-Betti numbers}\index{Betti number} $b_{g,d}^{T}(\gY,\gZ) \coloneqq \dim_\bk H_{g,d}^T(\gY,\gX;\bk) \in \bN \cup\{\infty\}$. Similarly defined are the absolute \emph{$T$-Betti numbers}.
	
If the monad $T$ is obtained from an operad $\cO$, we shall replace $H^T$ by the notation $H^\cO$ and $b^T$ by the notation $b^\cO$.

\subsection{Homotopy theory of filtered objects}\label{sec:homotopy-theory-filtered} Recall from Section~\ref{sec:filtered-objects-1} that we have adjunctions
\begin{equation}\label{eqn:filtered-adj}
\begin{tikzcd}
	{\sfC^{\Z_{=}}_\ast} \arrow[shift left=-.5ex,swap]{r}{u} & {\sfC^{\Z_\leq}}
	\arrow[shift left=-.5ex,swap]{l}{\mr{gr}} 
	\arrow[shift left=.5ex]{r}{\mr{colim}} & \sfC. \arrow[shift left=.5ex]{l}{\mr{const}}
\end{tikzcd}
\end{equation}

Following the discussion in Section \ref{sec.model-adjunctions}, as $\sfC$ satisfies the assumption of Section \ref{sec:axioms-of-model-cats} it follows from Lemma \ref{lem:model-structure-functor-cats} that the model structure on $\sfC$ induces projective model structures on $\sfC^{\Z_{=}}_\ast$ and $\sfC^{\Z_\leq}$, making them into cofibrantly generated $k$-monoidal (with respect to Day convolution) simplicial model categories when $\sfC$ is $k$-monoidal.

\begin{proposition}\label{prop:ModStrFilt}
	With these model structures, the adjunctions (\ref{eqn:filtered-adj}) are Quillen adjunctions.
\end{proposition}

\begin{proof}
  The right adjoint $\mr{const} \colon \sfC \to \sfC^{\bZ_\leq}$ sends $X$ to the
  constant functor $n \mapsto X$. Because weak equivalences and
  fibrations are objectwise in $\sfC^{\bN_\leq}$ this preserves
  fibrations and trivial fibrations, so $\mr{colim} \dashv \mr{const}$ is a
  Quillen adjunction.

  The right adjoint $u \colon \sfC^{\bZ_{=}}_* \to \sfC^{\bZ_\leq}$ sends
  $X \colon \bZ_{=} \to \sfC_*$ to the functor $u(X) \colon \bZ_\leq \to \sfC$
  which sends $n$ to $X(n)$ all non-identity morphisms to the
  constant map to the basepoint. In particular on objects it is given
  by ($\bZ$ copies of) the forgetful functor $U^+ \colon \sfC_* \to \sfC$. But
  $U^+$ is a right Quillen functor, as its left adjoint
  $+ \colon \sfC \to \sfC_*$ preserves weak equivalences and
  cofibrations. Thus $U^+$ preserves weak equivalences and fibrations,
  and so $u$ does too,  so $\mr{gr} \dashv u$ is a Quillen adjunction.
\end{proof}

It is possible to (partially) characterize cofibrant graded and filtered objects.

\begin{lemma}\label{lem:CofOb}
  An object $X \in \sfC^{\Z_=}$ is cofibrant if and only if each
  $X(n) \in \sfC$ is cofibrant. If
  $X \in \sfC^{\Z_\leq}$ is cofibrant, then each $X(n)$ is cofibrant and each structure map
  $X(n) \to X(n+1)$ is a cofibration. If $X \in \sfC^{\Z_\leq}$ is ascending
  and each structure map $X(n) \to X(n+1)$ is a cofibration, then $\mr{const}(X(-1)) \to X$ is a cofibration.
\end{lemma}
\begin{proof}
The case of $\sfC^{\Z_=}$ is immediate from the definition of the projective model structure. If $X \in \sfC^{\Z_\leq}$ is cofibrant then each $X(n)$ is cofibrant by Proposition 11.6.3 of \cite{Hirschhorn}. To see that $X(n) \to X(n+1)$ is a cofibration choose a trivial fibration $f \colon A \to B$ in $\sfC$ and consider the lifting problem
\begin{equation*}
\begin{tikzcd}
	X(n) \dar \rar{g}& A \dar{f}\\
	X(n+1) \rar{h} \arrow[dashed]{ru} & B
\end{tikzcd}
\end{equation*}
in $\sfC$. This gives rise to a lifting problem
\begin{equation*}
\begin{tikzcd}
	{\binit} \dar \rar & (\cdots \to A \to A \to A  \to \bterm \to \bterm \to \cdots) \dar{(\ldots, \mathrm{id}_A, f, \mr{id}_{\bterm},\ldots)}\\
	X \rar  \arrow[dashed]{ru}  & (\cdots \to A \to A \overset{f}\to B \to \bterm \to \bterm \to \cdots)
\end{tikzcd}
\end{equation*}
in $\sfC^{\Z_{\leq}}$. The right hand map is a trivial fibration, as these are defined objectwise and $f$ is a trivial fibration, so this lifting problem can be solved as $X$ is cofibrant. A solution to this determines in particular a solution to the original lifting problem: hence $X(n) \to X(n+1)$ has the left lifting property with respect to all trivial fibrations, so is a cofibration in $\sfC$. 

For the converse, suppose $X \in \sfC^{\Z_\leq}$ is ascending and that there is given a trivial fibration $f \colon A \to B \in \sfC^{\bZ_\leq}$ and a lifting problem
 \begin{equation*}
	\begin{tikzcd}
 \mr{const}(X(-1)) \dar \rar & A \dar{f}\\
 X \rar \arrow[dashed]{ru} & B.
 \end{tikzcd}
 \end{equation*}
 We will build a lift inductively. Since the filtration is ascending, it suffices to begin with
 \begin{equation*}
 \begin{tikzcd}
 X(-1) \dar \rar& A(0) \dar{f(0)}\\
 X(0) \rar \arrow[dashed]{ru} & B(0).
 \end{tikzcd}
 \end{equation*}
This has a solution $L_0$, as $f(0)$ is a trivial fibration and the left-hand map is a cofibration by assumption. Supposing compatible maps $L_i \colon X(i) \to A(i)$ have been constructed for $i < n$, consider
 \begin{equation*}
 \begin{tikzcd}
 X(n-1) \dar \rar{L_{n-1}}& A(n-1) \rar& A(n) \dar{f(n)}\\
 X(n) \arrow{rr} \arrow[dashed]{rru} && B(n)
 \end{tikzcd}
 \end{equation*}
 which has a solution $L_n$ as $f(n)$ is a trivial fibration and the left-hand map is a cofibration by assumption. These inductively determined $L_n$ give a lift $L \colon X \to A$ in the original lifting problem.
\end{proof}

\index{spectral sequence!filtered object}

\begin{theorem}\label{thm:SSAsc}
If $X \in \sfC^{\Z_{\leq}}$ is cofibrant then there is a spectral sequence
  \begin{equation*}
    E^1_{g,p,q} = \widetilde{H}_{g,p+q,p}(\grr(X);A) \Longrightarrow H_{g,p+q}(\mr{colim}(X);A)
  \end{equation*}
with differentials $d^r \colon E^r_{g,p,q} \to E^r_{g, p-r, q+r-1}$, which is conditionally convergent if
  \begin{equation*}
    \lim\limits_{p \to \infty} H_{*,*}(X(p);A)=0 = \limone\limits_{p \to \infty} H_{*,*}(X(p);A).
  \end{equation*}
\end{theorem}

This does not relate the different $g \in \sfG$ at all, so we may equivalently think of this as one spectral sequence for each $g$. The indexing of each these is that of the homological Serre spectral sequence.
	
\begin{proof}
The long exact sequences of the homology of the pairs $(X(q),X(q-1))$ assemble into an exact couple as usual, with
\[E^1_{g, p,q} = H_{g,p+q}(X(p),X(p-1);A).\]
As $X$ is cofibrant, by Lemma \ref{lem:CofOb} all maps $X(n) \to X(n+1)$ are cofibrations and so the homology of the pair $(X(p),X(p-1))$ is the same as the reduced homology of $\grr(X)(p)(g)$. This gives a spectral sequence with the claimed $E^1$-page, and with $A^1_{g,p,q} = H_{g,p+q}(X(p);A)$. The colimit $\colim(X)$ is a homotopy colimit, as $X$ is cofibrant. As taking derived singular simplicial chains commutes with homotopy colimits, and taking homology in $\cat{A}$ commutes with filtered colimits, the spectral sequence abuts to $H_{g,p+q}(\colim(X);A)$. Conditional convergence is by definition of that term, cf.\ \cite[Definition 5.10]{Boardman}.
\end{proof}

We give two easy applications of this general existence result for spectral sequences. We will later give more delicate applications.

\subsubsection{The geometric realization spectral sequence}\label{sec:geom-rel-ss}\index{spectral sequence!geometric realization}

The geometric realization $\gr{X_\bullet}$ of a simplicial object $X_\bullet$ has a canonical ascending filtration by skeleta: for $k \in \bZ$ the $k$-skeleton is the coend
\[{\gr{X_\bullet}^{(k)}} = \int^{n \in \Delta^\mr{op}_{\leq k}} \Delta^n \times X_n\]
over the full subcategory $\Delta^\mr{op}_{\leq k}$ of $\Delta^\mr{op}$ on those ordered finite sets of cardinality $\leq k$. There is a pushout diagram
\begin{equation}\label{eqn:geom-rel-skeletal-filt-pushout}
\begin{tikzcd} {\Delta^k \times L_k(X_\bullet)  \sqcup_{\partial \Delta^k \times L_k(X_\bullet) }  \partial \Delta^k \times X_k} \dar \rar & {\gr{X_\bullet}^{(k-1)}} \dar \\
\Delta^k \times X_k \rar	&  {\gr{X_\bullet}^{(k)}},
\end{tikzcd}
\end{equation}
As $\times \colon \cat{sSet} \times \sfC \to \sfC$ is a Quillen bifunctor, the left vertical map is a cofibration as long as $L_k(X_\bullet) \to X_k$ is, i.e.\ as long as $X_\bullet$ is Reedy cofibrant. In this case, as $\const(\gr{X_\bullet}^{(-1)}) = \binit$ is cofibrant, it follows from Lemma \ref{lem:CofOb} that this filtered object is cofibrant.

\begin{theorem}\label{thm:geom-rel-ss} If $X_\bullet$ is a Reedy cofibrant simplicial object, there is a spectral sequence
	\[E^1_{g, p, q} = H_{g,q}(X_p;A) \Longrightarrow H_{g,p+q}(\gr{X_\bullet};A),\]
which converges strongly, has $d^1$-differential given by $\sum_i (-1)^{i} (d_i)_*$, and has differentials $d^r \colon E^r_{g,p,q} \to E^r_{g,p-r,q+r-1}$.\end{theorem}

\begin{proof}
Applying Theorem \ref{thm:SSAsc} to the above cofibrant filtered object, the colimit of the skeletal filtration is $\gr{X_\bullet}$ and its associated graded has $p$th term given by $\Delta^p/\partial \Delta^p \wedge X_p/L_p(X_\bullet)$. Thus, using the suspension isomorphism (which is a consequence of Lemma \ref{lem:KunnethFormula} under assumption (ii)), there is a spectral sequence
	\[F^1_{g,p, q} = \widetilde{H}_{g,q}(X_p / L_p(X_\bullet);A) \Longrightarrow H_{g,p+q}(\gr{X_\bullet};A).\]
	As explained in the proof of \cite[Proposition 5.1]{SegalClassifyingSpaces}, $F^2_{g, p,q}$ agrees with the homology of $E^1_{g, p, q} = H_{g,q}(X_p;A)$ with respect to the differential $\sum_i (-1)^{i} (d_i)_*$.

It remains to verify the convergence condition, for which we note that it is a half-plane spectral sequence with exiting differentials and $A^\infty=0$ in the sense of Boardman, so converges strongly by \cite[Theorem 6.1]{Boardman}.
\end{proof}

A similar but easier argument applies to Reedy cofibrant (i.e. levelwise cofibrant) semi-simplicial objects, with $\partial \Delta^n \times X_n \to \Delta^n \times X_n$ replacing the left vertical arrow in (\ref{eqn:geom-rel-skeletal-filt-pushout}).

\begin{theorem}\label{thm:geom-rel-ss-thick} 
If $X_\bullet$ is a Reedy cofibrant semi-simplicial object, there is a spectral sequence
	\[E^1_{g, p, q} = H_{g,q}(X_p;A) \Longrightarrow H_{g, p+q}(\fgr{X_\bullet};A),\]
	which converges strongly, has $d^1$-differential given by $\sum_i (-1)^{i} (d_i)_*$, and has differentials $d^r \colon E^r_{g,p,q} \to E^r_{g,p-r,q+r-1}$.\end{theorem}

\subsubsection{The bar spectral sequence} \label{sec:bar-ss} An example of this type of spectral sequence is the \emph{bar spectral sequence}\index{spectral sequence!bar} associated to the two-sided bar construction described in Section \ref{sec:AssocModules}. Recall that for a unital associative algebra (i.e.\ a monoid) $\gR$ in $\sfC$ and right and left $\gR$-modules $\gN$ and $\gM$, the two-sided bar construction $B(\gN, \gR, \gM)$ is the thick geometric realization of the (semi-)simplicial object with $p$-simplices $B_p(\gN, \gR, \gM) = \gN \otimes \gR^{\otimes p} \otimes \gM$ and face maps given by the multiplication on $\gR$ and its action on $\gN$ and $\gM$. In this situation Theorem \ref{thm:geom-rel-ss-thick} gives a strongly convergent spectral sequence 
	\[E^1_{g, p, q} = H_{g, q}(\gN \otimes \gR^{\otimes p} \otimes \gM ;A) \Longrightarrow H_{g, p+q}(B(\gN, \gR, \gM);A)\]
with $d^1$-differential given by the alternating sum of the face maps.

Often the $E^1$-page of this spectral sequence may be simplified using a version of the K\"unneth formula. In particular, if a K\"unneth theorem as in Lemma \ref{lem:KunnethFormula} applies, then $H_*(\gR;\bk)$ is an augmented associative algebra object in the diagram category $\cat{GrMod}_\bk^\sfG$, where $\cat{GrMod}_\bk$ denotes the category of graded $\bk$-modules with tensor product involving a Koszul sign as normal, and $\cat{GrMod}_\bk^\sfG$ is given the Day convolution monoidal structure. 

Via the K\"unneth isomorphism $E^1_{*, p,*}  \cong H_{*, *}(\gN;\bk) \otimes H_{*, *}(\gR;\bk)^{\otimes p} \otimes H_{*, *}(\gM ;A)$, we may identify the $p$th column of the $E^1$-page with
\[(H_{*, *}(\gN;\bk) \otimes H_{*, *}(\gR;\bk)^{\otimes p+1}) \otimes_{H_{*, *}(\gR;\bk)} H_{*, *}(\gN ;A).\]
Under this identification the $d^1$-differential may be identified with that of the complex associated to the bar resolution of $H_{*, *}(\gN;\bk)$ by free $H_{*,*}(\gR;\bk)$-modules, so that we may identify $E^2_{*, p,\ast}$ with $\mathrm{Tor}_p^{H_{*, *}(\gR;\bk)}(H_{*, *}(\gN;\bk),H_{*,*}(\gM;A))$ with $\mr{Tor}$-groups taken in $\cat{GrMod}_\bk^\sfG$.

\subsubsection{The homotopy orbit spectral sequence} \label{sec:homotopy-orbit-ss} 
An special case of the above spectral sequence is the \emph{homotopy orbit spectral sequence}\index{spectral sequence!homotopy orbit}. Let $M$ be a unital monoid in $\sfC$, defining a monad $M \otimes -$ on $\sfC$ whose algebras are left $M$-spaces, and suppose that there is an augmentation $\epsilon \colon M \to \bunit$ of monoids. If $M$ arises from a monoid in simplicial sets then there is a canonical such augmentation.

If $X$ is an algebra for this monad it is acted upon by $M$ on the left, and the orbits $X/M$ is given by the reflexive coequalizer
\[\begin{tikzcd} M \otimes X \rar[shift left=1ex] \rar[shift left=-1ex] & \lar X \rar & X/M\end{tikzcd}\]
of the action map $M \otimes X \to X$ and the map $M \otimes X \to \bunit \otimes X \cong X$ induced by $\epsilon$. In other words, these are the indecomposables of the monad $M \otimes -$ with respect to the augmentation $M \otimes - \Rightarrow \mr{Id}$ induced by the augmentation $\epsilon$. 

\newglossaryentry{homotopyorbits}{%
	name={\ensuremath{\sslash}},
	description={Homotopy orbits},
type=symbols
}
We then define \emph{homotopy orbits}\index{homotopy!orbits} $X \gls{homotopyorbits} M$ to be its left derived functor. If $M$ is cofibrant in $\sfC$, then the monad $M \otimes -$ satisfies the axioms of Section \ref{sec:homotopy-theory-algebras}. By Section \ref{sec:simplicial-formula-indecomposables}, if $X$ is cofibrant in $\sfC$, this derived functor may be computed as the thick geometric realization of the two-sided bar construction
\[X \sslash M \coloneqq \fgr{B_\bullet(\bunit,M,X)}.\]
As in the previous section there is a strongly convergent spectral sequence 
\[E^1_{g, p, q} = H_{g, q}(M^{\otimes p} \otimes X ;A) \Longrightarrow H_{g, p+q}(X \sslash M;A)\]
with $d^1$-differential given by an alternating sum of the maps induced by the augmentation, the multiplication of $M$, and the action of $M$ on $X$. If there is a K\"unneth theorem available then we may identify $E^2_{*, p,\ast}$ with $\mathrm{Tor}_p^{H_{*, *}(M;\bk)}(\bk[\bunit],H_{*,*}(X;A))$.

\subsubsection{The change-of-diagram-category spectral sequence}\label{sec:change-of-diagram-ss}

Just before stating Lemma \ref{lem:oplax-lax-Day-convo}, we explained how to associate to a (strong $k$-monoidal) functor $f \colon \sfG \to \sfG'$ a (strong $k$-monoidal) functor $f_* \colon \sfS^\sfG \to \sfS^{\sfG'}$, which is left adjoint to restriction along $f$, and in Section \ref{sec:FunctorCategories} we explained that the model structures on these categories are such that $f_*$ is a left Quillen functor.\index{spectral sequence!change-of-diagram-category}

\begin{theorem}\label{thm:derived-change-of-diagram-ss}
Let $f \colon \sfG \to \sfG'$ be a functor between groupoids which induces surjective maps $G_g \to G'_{f(g)}$ on automorphism groups, with kernels $K_g$. If $X \in \sfS^\sfG$ is cofibrant then there is a spectral sequence
	\[E^2_{g',p,q} = \bigoplus_{\substack{[g] \in \pi_0\sfG\\ f(g) \cong g'}}H_{p}(K_g;H_{g,q}(X;A)) \Longrightarrow H_{g', p+q}(f_* X;A),\]
	which converges strongly, and has differentials $d^r \colon E^r_{g',p,q} \to E^r_{g',p-r,q+r-1}$.
\end{theorem}
\begin{proof}
By definition we have 
\[(f_*(X))(g') = \colim\limits_{g \in f/g'} X(g).\]
Under the stated conditions there is an equivalence of groupoids
\[f/g' \simeq \coprod_{\substack{[g] \in \pi_0\sfG\\ f(g) \cong g'}} K_g,\]
and so
\[(f_*(X))(g') \cong \coprod_{\substack{[g] \in \pi_0\sfG\\ f(g) \cong g'}} X(g)/K_g.\]

As $X$ is cofibrant in $\sfS^\sfG$, $X(g)$ is cofibrant in $\sfS^{K_g}$ and so $EK_g \times X(g) \to X(g)$ is a weak equivalence between cofibrant objects, hence $EK_g \times_{K_g} X(g) \to X(g)/K_g$ is a weak equivalence too. Considering $EK_g = \fgr{E_\bullet K_g}$ as the realization of a semi-simplicial object in $\cat{sSet}$, the above is $\fgr{E_\bullet K_g \times_{K_g} X(g)}$. Applying Theorem \ref{thm:geom-rel-ss-thick}, we obtain a spectral sequence which coincides with the direct sum of the homotopy orbit spectral sequences of $K_g$ acting on $X(g)$. As explained in Sections \ref{sec:bar-ss} and \ref{sec:homotopy-orbit-ss}, using Lemma \ref{lem:KunnethFormula} (ii) we may identify its $E^1$-page as 
\[E^1_{g',p,q} = \bigoplus_{\substack{[g] \in \pi_0\sfG\\ f(g) \cong g'}} \bk[K_g]^{\otimes p+1} \otimes_{\bk[K_g]} H_q(X(g);A),\]
where the $d^1$-differential on the $g$th summand may be identified with that of the bar complex computing $\mathrm{Tor}_{p}^{\bk[K_g]}(\bk, H_q(X(g);A))$, which is one definition of the group homology $H_{p}(K_g ; H_q(X(g);A))$.
\end{proof}

As described in Section \ref{sec:simplicial-operads}, for an operad $\cO$ in $\sfS$ the functor $f_*$ induces a functor $f_* \colon \Alg_\cO(\sfS^\sfG) \to \Alg_\cO(\sfS^{\sfG'})$. If the operad $\cO$ is augmented then the absolute $\cO$-indecomposables functor $Q^\cO$ is defined, and if both categories of algebras admit the projective model structure then the derived $\cO$-indecomposables and $\cO$-homology are defined too. In this case the following lemma allows one to compute the $\cO$-homology of $f_*(\gX)$ in terms of that of $\gX$.

\begin{corollary}
If $\gX \in \Alg_\cO(\sfS^\sfG)$ has underlying object cofibrant in $\sfS^\sfG$ then there is a spectral sequence
	\[E^2_{g',p,q} = \bigoplus_{\substack{[g] \in \pi_0\sfG\\ f(g) \cong g'}}H_p(K_g;H_{g,q}^\cO(\gX;A)) \Longrightarrow H_{g',p+q}^\cO(f_* \gX;A),\]
	which converges strongly, and has differentials $d^r \colon E^r_{g',p,q} \to E^r_{g',p-r,q+r-1}$.
\end{corollary}
\begin{proof}
Let $c\gX \overset{\sim}\to \gX$ be a cofibrant approximation in $\Alg_\cO(\sfS^\sfG)$, so $Q^\cO_\bL(\gX) \simeq Q^\cO(c\gX)$. The functor $f_* \colon \Alg_\cO(\sfS^\sfG) \to \Alg_\cO(\sfS^{\sfG'})$ is a left Quillen functor, because it sends the generating (trivial) cofibrations of $\Alg_\cO(\sfS^\sfG)$ to (trivial) cofibrations, as there is a natural isomorphism $f_* F^\cO \cong F^\cO f_*$ and $f_* \colon \sfS^\sfG \to  \sfS^{\sfG'}$ is a left Quillen functor. Hence the object $f_*(c\gX) \in \Alg_\cO(\sfS^{\sfG'})$ is cofibrant, and as the underlying object of $\gX$ is cofibrant the map $f_*(c\gX) \to f_*(\gX)$ is a weak equivalence. Thus we have $Q^\cO_\bL(f_*(\gX)) \simeq Q^\cO(f_*(c\gX)) \cong f_*(Q^\cO(c\gX))$.  Furthermore, as $Q^\cO \colon \Alg_\cO(\sfS^{\sfG'}) \to \sfS^{\sfG'}_*$ is a left Quillen functor the object $Q^\cO(c\gX)$ is cofibrant in $\sfS^{\sfG'}_*$. Applying Theorem \ref{thm:derived-change-of-diagram-ss} to $Q^\cO(c\gX)$ and using the above identifications gives the required spectral sequence.
\end{proof}

\subsection{Multiplicative filtrations of $\cO$-algebras}\label{sec:MultFilt}

Let $\cO$ be a $\Sigma$-cofibrant operad in $\sfC$. By the discussion in Section \ref{sec:model-str-alg-over-operads} and Proposition \ref{prop:ModStrFilt}, we obtain model structures on $\Alg_\cO(\sfC_*^{\bZ_{=}})$ and $\Alg_\cO(\sfC^{\bZ_\leq})$. As discussed in Section \ref{sec:filtrations-induced-functors-on-algebras}, the Quillen adjunctions
\begin{equation*}
\begin{tikzcd}
{\sfC^{\Z_{=}}_\ast} \arrow[shift left=-.5ex,swap]{r}{u} & {\sfC^{\Z_\leq}}
\arrow[shift left=-.5ex,swap]{l}{\grr} 
\arrow[shift left=.5ex]{r}{\mr{colim}} & \sfC, \arrow[shift left=.5ex]{l}{\mr{const}}
\end{tikzcd}
\end{equation*}
induce Quillen adjunctions
\begin{equation*}
\begin{tikzcd}
{\Alg_\cO(\sfC^{\Z_{=}}_\ast)} \arrow[shift left=-.5ex,swap]{r}{u} & {\Alg_\cO(\sfC^{\Z_\leq})}
\arrow[shift left=-.5ex,swap]{l}{\grr} 
\arrow[shift left=.5ex]{r}{\mr{colim}} & \Alg_\cO(\sfC) \arrow[shift left=.5ex]{l}{\mr{const}}
\end{tikzcd}
\end{equation*}
between the associated categories of $\cO$-algebras.

If $\cO$ is equipped with an augmentation $\epsilon \colon \cO \to +$ then there is an associated $\cO$-indecomposables functor $Q^\cO$, which commutes with the left adjoints $\grr$ and $\colim$ by Section \ref{sec:filtrations-induced-functors-on-algebras}.

\subsubsection{The derived indecomposables spectral sequence} 
Assuming that $\cO$ is an augmented operad so that $Q^\cO$ is defined, since it takes filtered algebras to filtered pointed objects, under suitable cofibrancy assumptions a filtration on an $\cO$-algebra induces a spectral sequence on $\cO$-homology.\index{spectral sequence!derived indecomposables}

\begin{theorem}\label{thm:DerIndecSS}
If $\gX \in \Alg_\cO(\sfC^{\bZ_{\leq}})$ has underlying object cofibrant in $\sfC^{\bZ_{\leq}}$ then
there is a spectral sequence
\[E^1_{g,p,q} \cong H_{g,p+q, p}^{\cO}(\grr(\gX);A) \Longrightarrow H_{g,p+q}^{\cO}(\colim \gX;A)\]
with differentials $d^r \colon E^r_{g,p,q} \to E^r_{g,p-r,q+r-1}$, which is conditionally convergent if
  \begin{equation*}
    \lim\limits_{p \to \infty} H_{*,*,p}^\cO(\gX;A)=0 = \limone\limits_{p \to \infty} H_{*,*,p}^\cO(\gX;A).
  \end{equation*}
\end{theorem}

\begin{proof}Let $c\gX \overset{\sim}\to\gX$ be a cofibrant approximation in $\Alg_\cO(\sfC^{\bZ_{\leq}})$, so that $Q_\bL^{\cO}(\gX) \coloneqq Q^\cO(c\gX) \in \sfC^{\bZ_{\leq}}_*$ is cofibrant as $Q^\cO$ is a left Quillen functor. Thus we may apply the spectral sequence of Theorem \ref{thm:SSAsc} to it, which takes the form
\[E^1_{g,p,q} = \widetilde{H}_{g,p+q,p}(\grr(Q^{\cO}(c\gX));A) \Longrightarrow \widetilde{H}_{g,p+q}(\colim  (Q^{\cO}(c\gX));A),\]
and is conditionally convergent under the given assumptions.

Let us identify the abutment. By definition of $Q_\bL^{\cO}(\gX)$, we have $\colim (Q_\bL^{\cO}(\gX)) \simeq \colim (Q^\cO(c\gX))$. As discussed in Section \ref{sec:filtrations-induced-functors-on-algebras}, $\colim $ commutes with $Q^\cO$ so we have that $\colim (Q^\cO(c\gX)) \cong Q^\cO(\colim (c\gX))$. Since $\colim $ is a left Quillen functor by Section \ref{sec:homotopy-theory-filtered}, and the map $c\gX \overset{\sim}\to \gX$ is a weak equivalence between objects which are cofibrant in $\sfC^{\bZ_{\leq}}$, the map $\colim (c\gX) \to \colim (\gX)$ is a weak equivalence. Furthermore, the source is a cofibrant $\cO$-algebra, again as $\colim $ is a left Quillen functor, so this may be used to compute $Q^\cO_\bL(\colim (\gX))$ as $Q^\cO(\colim (c\gX))$. This may be summarized as
\[(Q_\bL^{\cO}(\gX)) \simeq \colim (Q^\cO(c\gX)) \cong Q^\cO(\colim (c\gX)) \simeq Q^\cO_\bL(\colim (\gX)).\]

Similarly, to identify the $E^1$-page we note that, as discussed in Section \ref{sec:filtrations-induced-functors-on-algebras}, the functor $\grr$ commutes with $Q^\cO$, giving $\grr(Q^{\cO}(c\gX)) = \grr(Q^\cO(c\gX)) \cong Q^\cO(\grr(c\gX))$, and as $\grr$ is a left Quillen functor by Section \ref{sec:homotopy-theory-filtered} and $\gX$ is cofibrant in $\sfC^{\bZ_{\leq}}$, we also have $\grr(c\gX) \overset{\sim}\to \grr(\gX)$. This is a weak equivalence from a cofibrant $\cO$-algebra, so $Q^\cO_\bL(\grr(\gX))$ may be computed as $Q^\cO(\grr(c\gX))$.
\end{proof}

\subsubsection{The cell attachment spectral sequence} 

In Section \ref{sec:filtr-cell-attachm} we discussed the canonical ascending filtration associated to a cell attachment. In $\sfC = \sfS^\sfG$, given the data of an $\gX_0 \in \Alg_\cO(\sfC)$, a cofibration of simplicial sets $\partial D^d \hookrightarrow D^d$, an element $g \in \sfG$, and a map $e \colon \partial D^{d} \to \gX_0(g)$, this was defined as the pushout $\mr{f}\gX_1$ in $\Alg_\cO(\sfC^{\bZ_\leq})$ of the diagram
\[F^\cO(1_* D^{g,d}) \longleftarrow F^\cO(1_* \partial D^{g,d}) \lra 0_*\gX_0\]
which one should think of as putting $\gX_0$ is filtration degree $0$ and the cell in filtration degree $1$. The underlying object $\gX_1 = \colim \mr{f}\gX_1$ was also denoted $\gX_0 \cup^{\cO}_{e} \gD^{g,d}$. In Theorem \ref{thm:cell-attachment-associated-graded} we identified its associated graded as
\begin{equation}\label{eqn:cell-attachment-associated-graded}
\grr(\mr{f}\gX_1) \cong 0_*(\gX_0) \vee^{\cO} F^\cO(1_*(S^{g,d})).
\end{equation}
\index{spectral sequence!cell attachment}

\begin{lemma}\label{lem:cell-attachment-filtration-cofibrant} 
If $\gX_0 \in \Alg_\cO(\sfC)$ is cofibrant, then the filtered object $\mr{f}\gX_1$ is cofibrant in $\sfC^{\bZ_\leq}$.\end{lemma}

\begin{proof}
The filtered object $\mr{f}{\gX_1}$ is given by the pushout
	\[\begin{tikzcd} F^\cO(1_* \partial D^{g,d}) \rar \dar & 0_* \gX_0 \dar \\
	F^\cO(1_* D^{g,d}) \rar & \mr{f}\gX_1. \end{tikzcd}\]
Here $0_* \gX_0$ is cofibrant in $\Alg_\cO(\sfC^{\bZ_{\leq}})$ using Lemma \ref{lem:CofOb} and the assumption that $\gX_0$ is cofibrant in $\gX_0$. Since cofibrations are closed under pushouts, the map $0_* \gX_0 \to \mr{f}\gX_1$ is a cofibration and hence $\mr{f} \gX_1$ is cofibrant as well. Finally, we use that $U^\cO \colon \Alg_\cO(\sfC^{\bZ_{\leq}}) \to \sfC^{\bZ_\leq}$ preserves cofibrant objects by Lemma \ref{lem:proj-lax-monoidal-axiom}.
\end{proof}

\begin{corollary}\label{cor:cell-attachment-spectral-sequence}
If $\gX_0 \in \Alg_\cO(\sfC)$ is cofibrant, then there is a strongly convergent spectral sequence
\[E^1_{g,p,q} \cong \widetilde{H}_{g,p+q,p}\left(0_*(\gX_0) \vee^{\cO} F^\cO(1_*(S^{g,d});A\right) \Longrightarrow H_{g,p+q}(\gX_0 \cup^{\cO}_{e} \gD^{g,d};A)\]
with differentials $d^r \colon E^r_{g,p,q} \to E^r_{g,p-r,q+r-1}$.
\end{corollary}

\begin{proof}
  We apply the spectral sequence of Theorem \ref{thm:SSAsc} to the filtered object $U^\cO(\mr{f}{\gX_1}) \in \sfC^{\bZ_{\leq}}$, which is cofibrant by Lemma \ref{lem:cell-attachment-filtration-cofibrant}. As $\mr{f}\gX_1(d) = \binit$ for $d<0$ the assumptions for conditional convergence are satisfied, but even further this becomes a half-plane spectral sequence with exiting differentials and $A^\infty=0$ in the sense of Boardman, so converges strongly by \cite[Theorem 6.1]{Boardman}. The $E^1$-page is identified by (\ref{eqn:cell-attachment-associated-graded}).
\end{proof}

\subsubsection{The skeletal spectral sequence} 

In Section \ref{sec:cw-algebras}, we described the skeletal filtration of a relative CW algebra $f \colon \gR \to \gS$. By definition, $f \colon \gR \to \gS$ the underlying object of an ascendingly filtered object $\mr{sk}(f)$ and in Theorem \ref{thm:associated-graded-skeletal} we showed that its associated graded is given by
\begin{equation}\label{eqn:skeletal-associated}
\grr(\mr{sk}(f)) \simeq 0_*(\gR_+) \vee^\cO F^\cO\left(\bigvee_{d \geq 0} \bigvee_{\alpha \in I_d} d_*(S^{g_\alpha, d}_\alpha)\right).
\end{equation}
\index{spectral sequence!skeletal}

\begin{lemma}\label{lem:skeletal-filtration-cofibrant} 
If $\gR \in \Alg_\cO(\sfC)$ is cofibrant, then the filtered object $\mr{sk}(f)$ is cofibrant in $\sfC^{\bZ_{\leq}}$.
\end{lemma}
	
\begin{proof}
The map $0_*(\gR_+) \to \mr{sk}(f)$ is a cofibration in $\Alg_\cO(\sfC^{\bZ_{\leq}})$, as it is the transfinite composition of the maps $\mr{sk}_{d-1}(f) \to \mr{sk}_{d}(f)$ which are obtained as pushouts along cofibrations
\[F^\cO\left(\bigsqcup_{\alpha \in I_d} \partial D^{g_\alpha,d}_\alpha[d-1] \right) \lra F^\cO\left(\bigsqcup_{\alpha \in I_d} D^{g_\alpha,d}_\alpha[d]\right).\]
If $\gR \in \Alg_\cO(\sfC)$ is cofibrant then so is $0_*(\gR_+)$, as $0_*$ and $(-)_+$ are left Quillen functors, and hence $0_*(\gR_+) \to \mr{sk}(f)$ is a cofibration between cofibrant objects in $\Alg_\cO(\sfC^{\bZ_{\leq}})$. It follows from Lemma \ref{lem:proj-lax-monoidal-axiom} that its underlying map is a cofibration between cofibrant objects in $\sfC^{\bZ_{\leq}}$; in particular $\mr{sk}(f)$ is cofibrant in $\sfC^{\bZ_{\leq}}$.
\end{proof}

\begin{corollary}\label{cor:skeletalss}
If $\gR \in \Alg_\cO(\sfC)$ is cofibrant and $f \colon \gR \to \gS$ is a relative CW algebra, then there is a strongly convergent spectral sequence
\[E^1_{g,p,q} \cong \widetilde{H}_{g,p+q,p}\left(0_*(\gR_+) \vee^\cO F^\cO\left(\bigvee_{d \geq 0} \bigvee_{\alpha \in I_d} d_*(S^{g_\alpha, d}_\alpha)\right);A\right) \Longrightarrow H_{g,p+q}(\gS, \gR;A)\]
with differentials $d^r \colon E^r_{g,p,q} \to E^r_{g,p-r,q+r-1}$.
\end{corollary}
\begin{proof}
We apply Theorem \ref{thm:SSAsc} to the object $U^\cO (\mr{sk}(f)) \in \sfC^{\bZ_\leq}_*$, which is cofibrant by Lemma \ref{lem:skeletal-filtration-cofibrant}. Strong convergence is as in the proof of Corollary \ref{cor:cell-attachment-spectral-sequence}. The $E^1$-page is identified by (\ref{eqn:skeletal-associated}).
\end{proof}

\subsubsection{The canonical multiplicative filtration spectral sequence}

In Section \ref{sec:canon-mult-filtr} we have associated to a non-unitary operad $\cO$ and an $\cO$-algebra $\gR$ a canonical multiplicative filtration $(-1)_*^\alg(\gR)$. By Lemma \ref{lem:-1alg-descending} this is a descending filtration. This canonical filtration and its associated spectral sequence has been studied by Harper--Hess \cite{HarperHess} and Kuhn--Pereira \cite{KuhnPereira}. 

The spectral sequence for the canonical multiplicative filtration has quite subtle convergence properties, as the filtration is descending. In fact, to state them we shall have to borrow some terms which will be introduced in Section \ref{sec:hurewicz-cw}, and to prove them we will have to use some results from that section too. There is no risk of circularity, as this result is not necessary for anything which follows. We start by defining an abstract connectivity $c$ (cf.\ Definition \ref{defn:HomConnective}) by
\begin{equation}\label{eqn:connectivity-c} c(g) \coloneqq \begin{cases}
1 & \text{ if $g \in \sfG^\times$,}\\
0 & \text{ otherwise}.
\end{cases}\end{equation}
We shall assume in the next theorem that $\gR$ is homologically $c$-connective. We remark that this is implied by $\gR$ being homologically 0-connective (cf.\ Definition \ref{def:Connective}) and $\gR$ being reduced (cf.\ Definition \ref{defn:Reduced}).\index{spectral sequence!canonical multiplicative filtration}

\begin{theorem}\label{thm:CanMultFiltSS}
If $\gR \in \Alg_\cO(\sfC)$ is cofibrant then there is a spectral sequence
\[E^1_{g,p,q} = \widetilde{H}_{g, p+q, p}(F^\cO_{\cO(1)}(-1)_*Q^\cO_{\cO(1)}(\gR); A) \Longrightarrow H_{g, p+q}(\gR; A)\]
with differentials $d^r \colon E^r_{g,p,q} \to E^r_{g,p-r,q+r-1}$. If $\gR$ is homologically $c$-connective, each $\cO(n)$ is homologically 0-connective (cf.\ Definition \ref{defn:HomConnective}), and $\sfG$ is Artinian (cf.\ Definition \ref{defn:Artinian}), then this converges strongly.
\end{theorem}

\begin{proof}
Recall that $(-1)^\alg_*$ is left adjoint to the evaluation map
\[(-1)^* \colon \Alg_\cO(\sfC^{\bZ_\leq}) \lra \Alg_\cO(\sfC),\]
and this preserves fibrations and trivial fibrations (as these are defined pointwise on underlying objects) so is a right Quillen functor, and hence $(-1)_*^\alg$ is a left Quillen functor. Thus $(-1)_*^\alg(\gR) \in \Alg_\cO(\sfC^{\bZ_\leq})$ is cofibrant, so by Theorem \ref{thm:SSAsc} it has an associated spectral sequence. We have $\colim (-1)_*^\alg(\gR) \cong \gR$ which identifies the abutment. By Proposition \ref{prop:associated-graded-of-canonical-filt}, we may identify the associated graded as $\grr (-1)_*^\alg(\gR) \cong F^\cO_{\cO(1)}(-1)_*Q^\cO_{\cO(1)}(\gR)$, which identifies the $E^1$-page.

This is a half-plane spectral sequence with entering differentials in the sense of Boardman, so to show that it converges strongly we will show that it converges conditionally and then verify the hypothesis of \cite[Theorem 7.3]{Boardman}, i.e.\ that the derived $E^\infty$-page vanishes. By the $\limone$ exact sequence, conditional convergence is the same as asking for
\[\holim_{a \in \bZ_\leq} \bL C_*((-1)_*^\alg \gR(a);A) \simeq *.\]
It is enough to show this with $\bk$-coefficients. We will prove that:

\vspace{1ex}

\noindent\textbf{Claim:} {\it $\bL C_*((-1)^\mr{alg}_* \gR(-a);\bk)$ is $c^{\ast a}$-connective.}

\begin{proof}[Proof of claim]
We first prove this in the case that $\gR = F^\cO(X)$ is a free algebra with $X$ homologically $c$-connective. Then in the proof of Lemma \ref{lem:-1alg-descending}, which uses the assumption that $\cO$ is non-unitary, we saw that
\begin{equation*}
	U^\cO (-1)_*^\alg F^\cO(X)(-a) = \bigsqcup_{n \geq a} \cO(n) \otimes_{G_n} X^{\otimes n}.
\end{equation*}
As the derived functor $\bL C_*$ of the singular chain functor preserves arbitrary coproducts up to weak equivalence, the natural map
\[\bigoplus_{n \geq a} \bL C_*(\cO(n) \otimes_{G_n} X^{\otimes n}; \bk) \lra \bL C_*(U^\cO (-1)_*^\alg F^\cO(X)(-a); \bk)\]
is a weak equivalence. It thus suffices to prove that $\cO(n) \otimes_{G_n} X^{\otimes n}$ is $c^{* a}$-connective whenever $n \geq a$. 

As $X$ is homologically $c$-connective, by Lemma \ref{lem:connectivity-and-tensor-products} (i) the object $X^{\otimes n}$ is homologically $c^{* n}$-connective.  By the homotopy orbits spectral sequence of Section \ref{sec:homotopy-orbit-ss}, the fact that each $\cO(n)$ is homologically 0-connective, and Lemma \ref{lem:connectivity-and-tensor-products} (i), it follows that $\cO(n) \otimes_{G_n} X^{\otimes n}$ is $c^{*n}$-connective as well. As $c*c \geq c$, if $n \geq a$ then $c^{*n} \geq c^{* a}$. This has proved the claim for free $\cO$-algebras on homologically $c$-connected objects.

Let us now suppose that $\gR$ is a general cofibrant $\cO$-algebra. In particular $\gR$ is cofibrant in $\sfC$, so by the discussion in Section \ref{sec:monadic-bar-resolution}  there is a free simplicial resolution $\epsilon \colon \gR_\bullet = \sigma_*\sigma^*B_\bullet(F^\cO, \cO, \gR) \to \gR$ given by the thick monadic bar construction. As $(-1)_*^\alg$ is a left Quillen functor and preserves geometric realization (as it commutes with the copowering by simplicial sets), we have a weak equivalence $\gr{(-1)_*^\alg \gR_\bullet}_\cO \overset{\sim}\to (-1)_*^\alg(\gR)$ and so, on underlying objects, a weak equivalence 
\[\fgr{U^\cO (-1)_*^\alg B_\bullet(F^\cO, \cO, \gR)} \cong \vert U^\cO(-1)_*^\alg \gR_\bullet \vert \overset{\sim}\lra U^\cO(-1)_*^\alg(\gR)\]
using Lemma \ref{lem:monad-functors-geomrel} (iii). As $\gR$ is homologically $c$-connective, $\cO$ is non-unitary, and each $\cO(n)$ is homologically 0-connective, each $U^\cO (-1)_*^\alg B_p(F^\cO, \cO, \gR)$ is obtained by applying $U^\cO(-1)_*^\alg$ to a free $\cO$-algebra on a homologically $c$-connective object. Hence for fixed $g \in \sfG$ and $a \in \bZ$, $(U^\cO (-1)_*^\alg B_p(F^\cO, \cO, \gR))(-a)$ is homologically $c^{\ast a}$-connective.

We will prove by induction over $p$ that $\fgr{U^\cO (-1)_*^\alg B_\bullet(F^\cO, \cO, \gR)}^{(p)}(-a)$ is homologically $c^{\ast a}$-connective. For $p=0$, we have that
\[\fgr{U^\cO (-1)_*^\alg B_\bullet(F^\cO, \cO, \gR)}^{(0)}(-a) = (U^\cO (-1)_*^\alg F^\cO(U^\cO \gR))(-a)\]
is homologically $c^{\ast a}$-connective by the case proved above. For the induction step, we use the homotopy cofibre sequence
\[\gr{U^\cO(-1)_*^\alg \gR_\bullet}^{(p-1)} \lra \gr{U^\cO(-1)_*^\alg \gR_\bullet}^{(p)} \lra S^p \wedge (U^\cO (-1)_*^\alg B_p(F^\cO, \cO, \gR))_+\]
in $\sfC^{\bZ_\leq}$, which on applying $\bL C_*$ give homotopy cofibre sequence in $(\cat{A}^\sfG)^{\bZ_\leq}$, which are thus homotopy fibre sequences as $\cat{A}$ is stable. Since $(-a)^*$ is a right Quillen functor, it preserves homotopy fibre sequences. Applying $(-a)^*$, by the inductive hypothesis 
\[\bL C_*(\gr{U^\cO(-1)_*^\alg \gR_\bullet}^{(p-1)};\bk)(-a)\] is $c^{\ast a}$-connective, and by the case proved above
\[\bL C_*(S^p \wedge (U^\cO (-1)_*^\alg B_p(F^\cO, \cO, \gR))_+;\bk)(-a)\]
 is $(p+ c^{\ast a})$-connective. It follows that 
 \[\bL C_*(\gr{U^\cO(-1)_*^\alg \gR_\bullet}^{(p)};\bk)(-a)\] is $c^{\ast a}$-connective as well. 
This completes the proof of the induction step. The claim follows as $\bL C_*$ preserves sequential homotopy colimits.
\end{proof}

Let us make the connectivities more explicit. Since $\sfG$ is Artinian, it has a rank functor $r \colon \sfG \to \bN_{\leq}$. Observe that $c^{* n}(g) \geq n - r(g)$, as a decomposition $g \cong a_1 \oplus \cdots \oplus a_n$ must have at least $(n -r(g))$ $a_i$'s $\oplus$-invertible, by definition of a rank functor. We conclude that for each fixed $g \in \sfG$ and $d \in \bZ$, the inverse system $H_{g,d}(\bL C_*((-1)_*^\alg\gR(a);\bk))$ is eventually constantly equal to $0$. This implies that the homotopy limit is contractible, proving conditional convergence.

To see that the derived $E^\infty$-page vanishes, we establish a vanishing line. As $\gR$ is $c$-connective, $Q^\cO_{\cO(1)}(\gR)$ is too (this may be seen, for example, by using $Q^\cO_{\cO(1)}(\gR) \simeq B(\cO(1)_+, \cO, \gR)$ and working simplicially as above). Now the spectral sequence has
\begin{align*}
E^1_{g,p,q} &= \widetilde{H}_{g, p+q, p}(F^\cO_{\cO(1)}(-1)_*Q^\cO_{\cO(1)}(\gR); A) \\
&= \widetilde{H}_{g,p+q}( \cO(-p) \otimes_{G_{-p} \wr \cO(1)} Q^\cO_{\cO(1)}(\gR)^{\otimes (-p)}; A).
\end{align*}
As $Q^\cO_{\cO(1)}(\gR)$ is a cofibrant $\cO(1)$-module, and $G_{-p}$ acts freely on $\cO(-p)$, the quotient of $ \cO(-p) \otimes Q^\cO_{\cO(1)}(\gR)^{\otimes (-p)}$ by ${G_{-q} \wr \cO(1)}$ is in fact a homotopy quotient. Thus, using the homotopy orbits spectral sequence and the fact that each $\cO(n)$ is homologically 0-connective in the same way as above, we see that $E^1_{g,p,q}$ vanishes for $p+q \leq c^{* (-p)}(g)$. In particular using the estimate above it vanishes for $p+q \leq -p-r(g)$. Thus, holding $(g,p,q)$ fixed, the target of the differential
\[d^r \colon E^r_{g,p,q} \lra E^r_{g, p-r, q+r-1}\]
vanishes if $(p-r) + (q+r-1) \leq -(p-r) -r(g)$, i.e.\ if $p+q-1 \leq r-p-r(g)$, which is satisfied for all $r \gg 0$. Therefore $E^r_{g,p,q}$ is independent of $r$ for $r \gg 0$, so the derived $E^\infty$-page vanishes.
\end{proof}

\begin{remark}
Suppose that $f \colon \gR \to \gS$ is a morphism of cofibrant $\cO$-algebras, and that $\cO$, $\gR$, and $\gS$ satisfy the assumptions of Theorem \ref{thm:CanMultFiltSS} for strong convergence. If $f$ induces a homology equivalence
\[Q^\cO_{\cO(1)}(f) \colon Q^\cO_{\cO(1)}(\gR) \lra Q^\cO_{\cO(1)}(\gS),\]
then it induces an isomorphism between $E^1$-pages of the spectral sequence of that theorem (as $F^\cO_{\cO(1)}$ and $(-1)_*$ preserve homology equivalences), so by \cite[Theorem 5.3]{Boardman} 
it follows that $f \colon \gR \to \gS$ is a homology equivalence too. This is related to Theorem 1.12 (c) of \cite{HarperHess}, which studies the case of $\cO$-algebras in modules over a commutative symmetric ring spectrum with the positive projective stable model structure and calls the canonical multiplicative filtration spectral sequence the ``homotopy completion spectral sequence.'' We will study such questions in more detail, and by different methods, in Section \ref{sec:hurewicz-cw} (see especially Section \ref{sec:whitehead}).
\end{remark}

\subsection{Filtrations of associative algebras and their modules}

If $\gR$ is a unital associative algebra in $\sfC$, then in Example \ref{exam:module-over-algebra-operad} we described an operad in $\sfC$ whose algebras are (left) $\gR$-modules. All of the discussion so far can therefore be applied to manipulate filtered $\gR$-modules. 

However, we can also consider a filtered unital associative algebra $\gR$, i.e.\ a unital associative algebra in the category $\sfC^{\bZ_\leq}$. In this case, if $\gN$ and $\gM$ are right and left $\gR$-modules respectively, then as described in Section \ref{sec:AssocModules} we can form the two sided bar construction $B_\bullet(\gM, \gR, \gN) \in \mathsf{s}\sfC^{\bZ_\leq}$ and hence its thick geometric realization
\[B(\gM, \gR, \gN) \coloneqq \fgr{ B_\bullet(\gM, \gR, \gN)} \in \sfC^{\bZ_\leq}.\]
\index{spectral sequence!filtered module}

\begin{lemma}
If the underlying objects of $\gR$, $\gM$, and $\gN$ are cofibrant in $\sfC^{\bZ_\leq}$ then there is a spectral sequence
\[E^1_{n,p,q} \cong \widetilde{H}_{n,p+q,p}(B(\grr(\gM), \grr(\gR), \grr(\gN))) \Longrightarrow H_{n,p+q}(\colim B(\gM, \gR, \gN))\]
with differentials $d^r \colon E^r_{n,p,q} \to E^r_{n,p-r,q+r-1}$. If the filtrations on  $\gR$, $\gM$, and $\gN$ are ascending then it converges strongly.
\end{lemma}

\begin{proof}
The semi-simplicial object $B_\bullet(\gM, \gR, \gN) \in \cat{ss}\sfC^{\bZ_\leq}$ is levelwise cofibrant (as its $p$-simplices $\gM \otimes \gR^{\otimes p} \otimes \gN$ are a tensor product of cofibrant objects) so it is Reedy cofibrant, and hence $B(\gM, \gR, \gN) \in \sfC^{\bZ_\leq}$ is cofibrant by Lemma \ref{lem:thick-geom-rel-cofibrations}. We may thus apply Theorem \ref{thm:SSAsc} to it, to obtain a spectral sequence. To identify the $E^1$-page we use the isomorphism
\[\grr(B(\gM, \gR, \gN)) \cong B(\grr(\gM), \grr(\gR), \grr(\gN)) \in \sfC^{\bZ_{=}}_*,\]
which holds as $\grr$ is symmetric monoidal, and preserves (thick) geometric realizations as it is objectwise given by a pushout and commutes with the simplicial copowering.

It converges strongly under the given assumption by \cite[Theorem 6.1]{Boardman}, as then the filtration on $B(\gM, \gR, \gN)$ is also ascending, so it becomes a half-plane spectral sequence with exiting differentials and $A^\infty=0$.
\end{proof}

A similar spectral sequence has been studied by Angelini-Knoll--Salch \cite{AngeliniKnollSalch}: they consider (descendingly) filtered $E_\infty$-algebras (in fact, strictly commutative monoids), the induced filtration of the cyclic bar construction (i.e.\ topological Hochschild homology), and its associated spectral sequence.

\section{Hurewicz theorems and CW approximation}
\label{sec:hurewicz-cw}

In Sections \ref{sec:an-estimate} and \ref{sec:an-estimate-abs} we will establish Hurewicz theorems for $\cO$-homology. This culminates in Corollary \ref{cor:absolute-Hurewicz} which under certain conditions identifies the first non-trivial $\cO$-homology group in terms of the corresponding ordinary homology group. 

Using these, in Section \ref{sec:whitehead} we will establish conditions under which $\cO$-homology can be used to detect (homology or) weak homotopy equivalences. In Section \ref{sec:additive-case} we will develop a theory of (minimal) CW approximations, which makes use of the ordinary Hurewicz theorem comparing homology and homotopy groups. 

\subsection{Connectivity functors}
\label{sec:conn-funct} We begin with a discussion of the appropriate method to keep track of (homological) connectivity of objects in diagram categories, as well as establish the behavior of connectivity under tensor products.

For simplicial sets or topological spaces, it is common to say that a map $X \to Y$ is \emph{homologically $c$-connective} if $H_d(Y,X) = 0$ for all $d < c$. 
We would like to point out the difference between ``connective'' and ``connected'': the term ``$(c-1)$-connected'' usually means $c$-connective. In this section we shall discuss how to encode vanishing conditions on the homology of each of the values $X(g)$ of an object $X \in \sfC = \sfS^\sfG$.

\begin{definition}\label{defn:HomConnective}
Let $[-\infty,\infty]_\geq$ be the category with objects given by the set $[-\infty,\infty]$ of extended real numbers, and a unique morphism $x \to y$ if and only if $x \geq y$. We endow it with a symmetric monoidal structure given by addition, with the convention
\[(\infty)+(-\infty) = (\infty).\]
The category $[-\infty,\infty]_\geq$ has all colimits, and this monoidal structure preserves colimits in each variable.

An \emph{abstract connectivity}\index{connectivity!abstract} for $\sfG$ is a functor $c \colon \sfG \to [-\infty,\infty]_\geq$.
\end{definition}

\begin{definition}\label{def:Connective} For an abstract connectivity $c$ and a commutative ring $\bk$, a morphism $f \colon X \to Y$ in $\sfC = \sfS^\sfG$ is \emph{homologically $c$-connective}\index{connective!homologically} if $H_{g,d}(Y,X;\bk) = 0$ for $d < c(g)$.  An object $X \in \sfC$ is \emph{homologically $c$-connective} if $H_{g,d}(X;\bk) = 0$ for $d < c(g)$.
\end{definition}

It follows from the properties of a singular chains functor that $H_{*,*}(\binit;\bk)=0$, so we can equivalently define $X$ to be homologically $c$-connective if the morphism $\binit \to X$ is. We warn the reader that when $\sfS$ is not pointed, this condition may not be what they expect: the condition is about the map from the initial object, not the map to the terminal object.  For example, a
functor $X \colon \sfG \to \cat{sSet}$ is homologically $c$-connective precisely
when $X(g) = \varnothing$ whenever $c(g) > 0$.  This seemingly unusual
definition is in fact desirable, since it makes connectivity be
``additive under tensor product'', cf.\ Lemma~\ref{lem:connectivity-and-tensor-products} (i) below.

\newglossaryentry{conv}{%
	name={\ensuremath{-\ast-}},
	description={Convolution of abstract connectivities},
type=symbols
}
\newglossaryentry{1conn}{%
	name={\ensuremath{\bunit_\mr{conn}}},
	description={Monoidal unit of abstract connectivities},
type=symbols
}
The $k$-monoidal structure on $\sfG$ induces a $k$-monoidal structure on abstract connectivities with tensor product \gls{conv}, as follows. Let $c$ and $c'$ be abstract connectivities for $\sfG$, then their \emph{convolution} is the abstract connectivity $c \ast c'$, defined by Day convolution of functors $c,c' \colon \sfG \to [-\infty,\infty]_\geq$, using that the target category has all colimits. Explicitly, we have
\begin{equation}\label{eqn:abstract-connectivity-unit}
(c\ast c')(g) = \inf \{c(a) + c'(a') \mid \sfG(a \oplus a',g) \neq
\varnothing \} \in [-\infty,\infty].
\end{equation}
Recognize $\sfG(a \oplus a',g) \neq \varnothing$ if and only if $a \oplus a' \cong g$, as $\sfG$ is a groupoid. This endows the category of functors $\sfG \to [-\infty,\infty]_{\geq}$ with a $k$-monoidal structure, with unit given by 
\begin{equation}\label{eqn:1-conn}\gls{1conn}(g) = \begin{cases} 0  & \text{if $\sfG(\bunit_\sfG,g) \neq \varnothing$,}\\ \infty & \text{otherwise.}\end{cases}\end{equation}

Recall that the cofiber of a cofibration $f \colon X \to Y$ is the pushout $Y/X = Y \cup_X \bterm$, which comes with a canonical map from $\bterm$ and is hence considered as an object of $\sfC_*$. The following is a direct consequence of the definitions:

\begin{lemma}\label{lem:ConnectivityViaCofibre}
Let $f \colon X \to Y$ be a cofibration and let $c$ be an abstract connectivity.  Then $f$ is a homologically $c$-connective morphism in $\sfC$ if and only if $\bterm \to Y/X$ is a homologically $c$-connective morphism.
\end{lemma}

\begin{proof}
  There is an identification $H_*(Y,X;\bk) \cong H_*(Y/X, \bterm ; \bk)$ as in Section \ref{sec:homology-objects-sfc}.
\end{proof}

The next lemma is preparation for Corollary \ref{cor:connectivity-of-tensor-power} below, giving an estimate on the connectivity of iterated tensor products of maps.

\begin{lemma}\label{lem:connectivity-and-tensor-products}\mbox{}
\begin{enumerate}[(i)]
\item Let $X, X' \in \sfC$ be cofibrant, and assume that $X$ is homologically
  $c$-connective and $X'$ is homologically $c'$-connective.  Then
  $X \otimes X'$ is homologically $(c \ast c')$-connective.

\item Let $f \colon X'' \to X'$ be a homologically $c_f$-connective morphism between
  cofibrant objects of $\sfC$ and let $X \in \sfC$ be a homologically $c$-connective
  cofibrant object. Then $X \otimes f \colon X \otimes X'' \to X \otimes X'$
  is homologically $(c * c_f)$-connective, as is $f \otimes X$.
\end{enumerate}
\end{lemma}

\begin{proof}
The first part follows from the second, applied to the morphism $X'' = \binit \to X'$. To prove the second part we may as well suppose that $f$ is a cofibration, so $X \otimes f$ is too and, as $X \otimes -$ preserves pushouts, the cofibre of $X \otimes f$ is isomorphic to the cofibre of $X \otimes \bterm \to X \otimes (X'/X'')$. By Lemma \ref{lem:ConnectivityViaCofibre} the map $\bterm \to X'/X''$ is homologically $c_f$-connective.

The strongly convergent K\"unneth spectral sequence of Lemma \ref{lem:KunnethSS} has an obvious relative elaboration,
\[\begin{tikzcd}\bigoplus_{q'+q''=q}\mathrm{Tor}_p^{\bk}(H_{q'}(X(a)), H_{q''}((X'/X'')(a'), \bterm)) \dar[Rightarrow] \\ H_{p+q}(C_*(X(a)) \otimes C_*((X'/X'')(a'), \bterm)),\end{tikzcd}\]
and we have $H_{q'}(X(a))=0$ for $q' < c(a)$ and $H_{q''}((X'/X'')(a'), \bterm)=0$ for $q'' < c_f(a')$, so the target vanishes in degrees $ < c(a) + c_f(a')$. Therefore it vanishes in degrees $ < (c \ast c_f)(g)$, whenever there exists a morphism $f \colon a \oplus a' \to g$ in $\sfG$. It then follows from the spectral sequence
\[\begin{tikzcd}E^2_{g, s,t} = \bL_s\colim_{(a, a', f) \in \sfH_g} H_t(C_*(X(a)) \otimes C_*((X'/X'')(a'), \bterm)) \dar[Rightarrow] \\ H_{g, s+t}(X \otimes (X'/X''), X \otimes \bterm),\end{tikzcd}\]
of Lemma \ref{lem:KanExtCofibrant} that the map $X \otimes \bterm \to X \otimes (X'/X'')$ is homologically $(c*c_f)$-connective.
\end{proof}

\begin{corollary}\label{cor:connectivity-under-tensor2}
  Let $f \colon X \to X'$ and $g \colon Y \to Y'$ be morphisms between cofibrant
  objects.  Assume $X$ is homologically $c_X$-connective, $f$ is
  homologically $c_f$-connective, and so on.  Then
  $f\otimes g \colon X \otimes Y \to X' \otimes Y'$ is homologically
  \begin{equation*}
    \max (\min(c_X \ast c_g,c_{f } \ast c_{Y'}), \min(c_f \ast
    c_Y,c_{X'}\ast c_g))\text{  -connective.}
  \end{equation*}
\end{corollary}
\begin{proof}
  The four convolutions in the formula are the connectivities of the four arrows in the diagram
  \begin{equation*}
    \begin{tikzcd}
      X \otimes Y \rar{X \otimes g} \dar[swap]{f \otimes Y} & X
      \otimes Y' \dar{f \otimes Y'}\\
      X' \otimes Y \rar{X' \otimes g} & X' \otimes Y'.
    \end{tikzcd}
  \end{equation*}
  The corollary follows because the connectivity of a composition is at least the minimum of the two connectivities, by the five lemma, and we may pick the way around the diagram which results in the maximal connectivity.
\end{proof}

Applying the previous corollary $(n-1)$ times, we obtain:

\begin{corollary}\label{cor:connectivity-of-tensor-power}
  If $f \colon X \to Y$ is a homologically $c_f$-connective morphism between
  cofibrant objects which are homologically $c_X$ and $c_Y$-connective, then
  $f^{\otimes n} \colon X^{\otimes n} \to Y^{\otimes n}$ is homologically
  \begin{equation*}
    \min \{c_X^{\ast a} \ast c_f \ast c_Y^{\ast b}\mid a + b = n-1\}\text{-connective.}
  \end{equation*}
\end{corollary}

\subsection{Hurewicz theorems for relative indecomposables}\label{sec:an-estimate}

In this section, we prove the first Hurewicz theorem comparing ordinary homology to $\cO$-homology. To do so, we will additionally assume that $\cO$ is a non-unitary $\Sigma$-cofibrant operad and that all $\cO(n)$ are homologically $0$-connective: we will abbreviate the latter condition by saying that $\cO$ is homologically $0$-connective. We start by applying the results of Section \ref{sec:conn-funct} to the monad associated to this operad.

\begin{lemma}\label{lem:HomologicalDecVsConnectivities}
Let $\cO$ be a non-unitary homologically 0-connective $\Sigma$-cofibrant operad. Let $c, c_f \colon \sfG \to  [-\infty,\infty]_\geq$ be abstract connectivities such that $c \ast c \geq  c$, $c \ast c_f \geq c_f$, and $c_f \ast c \geq c_f$.  Let $X$ and $Y$ be homologically $c$-connective cofibrant objects of $\sfC$ and let $f \colon X \to Y$ be a homologically $c_f$-connective map.  Then 
\begin{enumerate}[(i)]
\item $\cO(X)$ is homologically $c$-connective, 
\item $\cO(f)$ is homologically $c_f$-connective, and
\item $\Dec^\cO_{\cO(1)}(F^\cO f)$ is homologically $\min\{c*c_f, c_f * c\}$-connective.
\end{enumerate}
\end{lemma}
\begin{proof}
Note that (i) follows from (ii) applied to $f \colon \binit \to X$. For (ii), we have
  $\cO(X) = \bigsqcup_{n \geq 1} \cO(n) \times_{G_n} X^{\otimes n}$.
  The map $f^{\otimes n} \colon X^{\otimes n}\to Y^{\otimes n}$ is homologically
  $\min \{c^{\ast a} \ast c_f \ast c^{\ast b}\mid a + b = n-1\}$-connective by Corollary \ref{cor:connectivity-of-tensor-power}, so as $c \ast c_f \geq c_f$ and $c_f \ast c \geq c_f$ it is homologically $c_f$-connective for all
  $n \geq 1$. That $\cO(n) \times_{G_n} f^{\otimes n}$ is homologically $c_f$-connective then follows from the fact that $\cO(n)$ is homologically 0-connective, Lemma \ref{lem:connectivity-and-tensor-products} (ii), and the map of homotopy orbit spectral sequences
\[E^2_{g, p,q} = \mathrm{Tor}^{\bk[G_n]}_p(\bk, H_{g,q}(\cO(n) \times f^{\otimes n} ; \bk)) \Longrightarrow H_{g, p+q}(\cO(n) \times_{G_n} f^{\otimes n}; \bk)\]
of Section \ref{sec:homotopy-orbit-ss}, using that the $G_n$-action on $\cO(n)$ is free.

For (iii) a similar argument applies to $\Dec^\cO_{\cO(1)}(F^\cO f) = \bigsqcup_{n \geq 2} \cO(n)\times_{G_n} f^{\otimes n}$, except now each summand is homologically $\min\{c*c_f, c_f * c\}$-connective.
\end{proof}

\begin{lemma}\label{lem:DecConnectivity}
Let $\cO$ be a non-unitary homologically 0-connective $\Sigma$-cofibrant operad. Let $c, c_f \colon \sfG \to [-\infty,\infty]_\geq$ be abstract connectivities such that $c \ast c \geq  c$, $c \ast c_f \geq c_f$, and $c_f \ast c \geq c_f$. Let $f \colon \gR \to \gS$ be a homologically $c_f$-connective morphism of homologically $c$-connective $\cO$-algebras. Then
\[\bL\Dec^\cO_{\cO(1)}(f)\colon  \bL\Dec^\cO_{\cO(1)}(\gR) \lra \bL\Dec^\cO_{\cO(1)}(\gS)\]
is homologically $\min\{c*c_f, c_f * c\}$-connective.
\end{lemma}

\begin{proof}
As the desired conclusion is about derived decomposables, without loss of generality we may suppose that $\gR$ and $\gS$ are cofibrant $\cO$-algebras, and hence (by Axiom \ref{axiom:monad-proj}) that their underlying objects are cofibrant in $\sfC$. Thus as described in Section \ref{sec:deriv-decomp} we may compute their derived decomposables using the monadic bar resolution. 

  The $p$-simplices of $B_\bullet(\Dec^\cO_{\cO(1)} F^\cO,\cO,\gR)$ are
  $\Dec^\cO_{\cO(1)} F^\cO \cO^p(\gR)$.  By iteratedly applying Lemma \ref{lem:HomologicalDecVsConnectivities} we see that
  $\cO^p f$ is homologically $c_f$-connective for all $p \geq 0$ and hence that
  $\Dec^\cO_{\cO(1)} F^\cO \cO^p f$ is homologically $\min\{c*c_f, c_f * c\}$-connective for all $p$.  The
  geometric realization in $\sfC$ is calculated objectwise, and so
  preserves homological connectivity, so the geometric realization is also
  homologically $\min\{c*c_f, c_f * c\}$-connective.
\end{proof}

\begin{proposition}\label{prop:UsefulHurewicz} 
Let $\cO$ be a non-unitary homologically 0-connective $\Sigma$-cofibrant operad. Let $\gR, \gS \in \Alg_\cO(\sfC)$ be cofibrant and homologically $c$-connective for some abstract connectivity $c$ such that $c * c \geq c$. Let $f \colon \gR \to \gS$ be a map such that $U^\cO f$ is homologically $c_f$-connective. Then the square
  \begin{equation*}
      \begin{tikzcd}
        (U^\cO \gR)_+ \rar\dar & Q^\cO_{\cO(1)} \gR \dar\\
        (U^\cO \gS)_+ \rar & Q^\cO_{\cO(1)} \gS 
      \end{tikzcd}
  \end{equation*}
  is homologically $(1 + \min\{c*c_f, c_f * c\})$-cocartesian, i.e.\ the induced
  map
  \begin{equation*}
    H_{g,d}(\gS, \gR) \lra H_{g,d}(Q^\cO_{\cO(1)} 
    \gS,Q^\cO_{\cO(1)} \gR)
  \end{equation*}
  is an epimorphism for $d < (1 + \min\{c*c_f, c_f * c\})(g)$ and an isomorphism
  for $d < (\min\{c*c_f, c_f * c\})(g)$.
\end{proposition}

\begin{proof}
As $\gR$ and $\gS$ are cofibrant $\cO$-algebras, we may identify the square in question with the right-hand square of
  \begin{equation*}
      \begin{tikzcd}
        \bL\Dec^\cO_{\cO(1)}(\gR) \dar \rar & U^\cO_{\cO(1)} \fgr{B_\bullet(F^\cO , \cO, \gR)}_+ \rar\dar & \bL Q^\cO_{\cO(1)} \gR \dar\\
        \bL\Dec^\cO_{\cO(1)}(\gS) \rar& U^\cO_{\cO(1)} \fgr{B_\bullet(F^\cO , \cO, \gS)}_+ \rar & \bL Q^\cO_{\cO(1)} \gS,
      \end{tikzcd}
  \end{equation*}
where the rows are cofibration sequences described in Section \ref{sec:deriv-decomp}. Establishing the result in question means showing that
\[\bL\Dec^\cO_{\cO(1)}(f) \colon \bL\Dec^\cO_{\cO(1)}(\gR) \lra \bL\Dec^\cO_{\cO(1)}(\gS)\]
is homologically $\min\{c*c_f, c_f * c\}$-connective, which is the content of Lemma \ref{lem:DecConnectivity}.
\end{proof}

The next result is our Hurewicz theorem for relative $\cO$-indecomposables, the main result of this section. For  a map $f \colon \gR \to \gS$ of $\cO$-algebras we wish to make a statement about ``the lowest" degree in which $H_{*,*}(\gS, \gR)$ is not known to vanish. In order to make sense of this, we make the following definition.

\begin{definition}\label{defn:Artinian}
A $k$-monoidal groupoid $(\sfG,\oplus,\bunit)$ is \emph{Artinian}\index{Artinian} if there exists a lax $k$-monoidal functor $r \colon \sfG \to \bN_\leq$, such that $r(g) > 0$ when $g\in\sfG$ is not $\oplus$-invertible. We let $\sfG^\times$ denote the full subcategory of $g$ that are $\oplus$-invertible.
\end{definition}

When $\sfG$ is Artinian, $r(g)$ gives an upper bound on the number of summands in a decomposition $g \cong g_1 \oplus \dots \oplus g_n$ with no $g_i$ invertible under $\oplus$.  In fact, if such an $r$ exists, then we may define $\omega \colon \sfG \to \N_\leq$ by
\[\omega(g) \coloneqq \sup \{n \mid g \cong g_1 \oplus \dots \oplus g_n \text{ with no $g_i$ $\oplus$-invertible}\},\]
using $r$ to show that the supremum is attained, so is in fact a maximum. Then $\omega$ is in fact lax monoidal and satisfies $\omega(g) \leq r(g)$ for any rank functor $r$. Let us call $\omega$ the \emph{canonical rank functor}\index{canonical rank functor}.

\begin{definition}\label{defn:Reduced}
An object $X \in \sfC$ is \emph{reduced}\index{reduced} if $H_{g,0}(X)=0$ for all $g \in \sfG$ that are $\oplus$-invertible.
\end{definition}

To state the main result, we introduce the relation $\lneqq$ on $\bZ^2$: $(\omega',d') \lneqq (\omega,d)$ if $(\omega',d')$  belongs to the set
\begin{equation*}
\Z_{\leq \omega} \times \Z_{\leq d -1} \cup \Z_{\leq
	\omega -1} \times \Z_{\leq d},
\end{equation*}
as in Figure \ref{fig:omega-d-diagram}. 

\begin{figure}[h]
	\centering
	\begin{tikzpicture}

	\matrix (m) [matrix of nodes,nodes={minimum width=1cm,minimum height=1cm}]
	{ 	3 & \quad & \quad & \quad & \quad & \quad \\
		2 & $\bullet$ & $\bullet$ & \quad & \quad & \quad \\
		1 & $\bullet$ & $\bullet$ & $\bullet$ & \quad & \quad \\
		0 & $\bullet$ & $\bullet$ & $\bullet$ & \quad & \quad  \\
		$\nicefrac{d}{\omega}$ & 0 & 1 & 2 & 3 & 4 \\
	};
	
	\foreach \r in {3, 4, 5, 6} {\mvline[dotted]{m}{\r}}
	\mvline{m}{2}
	\foreach \r in {1, 2, 3, 4} {\mhline[dotted]{m}{\r}{5}}
	\mhline{m}{5}{2}
	
	\fill[pattern=north west lines,pattern color=black!50!white] ($(m-2-2.north west)+(-.25,.3)$) rectangle ($(m-4-3.south east)+(.25,-.25)$);
	\fill[pattern=north west lines,pattern color=black!50!white] ($(m-3-2.north west)+(-.25,.3)$) rectangle ($(m-4-4.south east)+(.25,-.25)$);

	\end{tikzpicture}
	\caption{The set $\Z_{\leq \omega} \times \Z_{\leq d -1} \cup \Z_{\leq
			\omega -1} \times \Z_{\leq d}$ for $d=2$ and $\omega=2$.}
	\label{fig:omega-d-diagram}
\end{figure}

\begin{corollary}\label{cor:homology-Hurewicz}
Let $\cO$ be a non-unitary homologically 0-connective $\Sigma$-cofibrant operad. Let $\sfG$ be an Artinian groupoid, and $\gR, \gS \in \Alg_\cO(\sfC)$ be cofibrant, reduced, and homologically 0-connective. Let $f \colon \gR \to \gS$ be a morphism such that $H_{g',d'}(\gS,\gR) = 0$ whenever $(\omega(g'),d') \lneqq (\omega(g),d)$, for some $d \in \bZ$ and $g \in \sfG$. Then the induced map
  \begin{equation*}
    H_{g,i}(\gS, \gR) \lra H_{g,i}(Q^\cO_{\cO(1)} \gS, Q^\cO_{\cO(1)} \gR)
  \end{equation*}
  is an isomorphism for $i=d$, and a surjection for $i=d+1$.
\end{corollary}

\begin{proof}
Recall that $\sfG^\times$ denotes those objects that are $\oplus$-invertible. Define an abstract connectivity $c$ by
\[c(g) \coloneqq  \begin{cases}
1 & \text{ if $g \in \sfG^\times$,}\\
0 & \text{ otherwise}.
\end{cases}\]
Then $c*c \geq c$, and as $\gS$ and $\gR$ are homologically 0-connective and reduced, they are homologically $c$-connective.

Let $c_f$ denote the homological connectivity of the map $f$. If there is a morphism $a \oplus b \overset{\sim}\to g \in \sfG$ then $\omega(a) + \omega(b) \leq \omega(g)$. If $b \nin \sfG^\times$ then $\omega(b) \geq 1$ and so $\omega(a) < \omega(g)$, and hence $H_{a, d'}(\gS, \gR)=0$ for $d' < d+1$, and so $c_f(a) \geq d+1$; therefore $c_f(a) + c(b) \geq d+1$. On the other hand if $b \in \sfG^\times$ then $a \oplus b = g$ so $\omega(a) \leq \omega(g)$ and $g \oplus b^{-1} = a$ so $\omega(g) \leq \omega(a)$, which imply that $\omega(a) = \omega(g)$, and hence $H_{a, d'}(\gS, \gR)=0$ for $d' < d$, so $c_f(a) \geq d$, but also $c(b) \geq 1$ as $g \in \sfG^\times$, so $c_f(a) + c(b) \geq d+1$ in this case too. By the formula for $(c_f * c)(g)$ as an infimum it follows that
\[(c_f * c)(g) \geq d+1.\]
Similarly for $c * c_f$. Applying Proposition \ref{prop:UsefulHurewicz} gives the required result.
\end{proof}

For the little $n$-cubes operad in $S$-modules this has been proved by Basterra--Mandell \cite[Theorem 3.7]{BMBP}. A similar result has been obtained by Harper--Hess \cite[Theorem 1.8]{HarperHess} for operads in symmetric spectra or chain complexes (what they call $\mathsf{TQ}$ or $\mathsf{Q}$ is what we call $\bL Q^\cO_{\cO(1)}$, see Definitions 3.15 and 8.4 of their paper), and by Basterra \cite[Lemma 8.2]{Basterra} for the commutative operad in $S$-algebras.

\subsection{Hurewicz theorems for absolute indecomposables}\label{sec:an-estimate-abs} Corollary \ref{cor:homology-Hurewicz} concerned the relative indecomposables $Q^\cO_{\cO(1)}$. If $\cO$ is a non-unitary operad equipped with an augmentation $\epsilon \colon \cO(1) \to \bunit_\sfC$, then we may also form the absolute indecomposables $Q^\cO$ as described in Section \ref{sec:simplicial-operad-augmentation}. As $Q^\cO = Q^{\cO(1)} Q^\cO_{\cO(1)}$, the derived functors are related by the homotopy orbit construction of Section \ref{sec:homotopy-orbit-ss}, as
\begin{equation}\label{eq:AbsIndec}
\bL Q^\cO(\gR) \simeq (\bL Q^\cO_{\cO(1)}(\gR)) \sslash \cO(1).
\end{equation}

The homology groups $H_{*,0}(\cO(1);\bk) \in \cat{Mod}_\bk^\sfG$ form an (augmented) algebra in this category, and if $\gR$ is a $\cO(1)$-algebra then in particular each $H_{*,d}(\gR;\bk)$ has the structure of a $H_{*,0}(\cO(1);\bk)$-module. As the composition $(\cO(1) \otimes \gR)_+ \to \gR_+ \to Q^{\cO(1)}(\gR)$ canonically factors through the augmentation, if $c\gR \overset{\sim}\to\gR$ is a cofibrant approximation then the map on homology induced by 
\[\gR_+ \overset{\sim}\longleftarrow c\gR_+ \lra Q^{\cO(1)}(c\gR) \simeq \gR \sslash \cO(1)\]
descends to a map
\[\bk[\bunit] \otimes_{H_{*,0}(\cO(1);\bk)}H_{*,d}(\gR ; \bk) \lra H_{*,d}(\gR \sslash \cO(1);\bk) = H_{*,d}^{\cO(1)}(\gR;\bk),\]
where $\bk[\bunit]$ is the functor $(\bunit_\sfG)_\ast(\bk)$ given by $g \mapsto \bk[\sfG(\bunit_\sfG,-)]$.

\begin{lemma}\label{lem:CoinvariantsOnFirstNonvanishing}
Suppose that $\cO(1)$ is homologically 0-connective. If $f \colon \gR \to \gS$ is a morphism of $\cO(1)$-algebras which satisfies $H_{*,d'}(\gS, \gR ; \bk) = 0$ whenever $d' < d$, for some $d \in \bZ$, then 
  \begin{equation*}
     \bk[\bunit] \otimes_{H_{*,0}(\cO(1);\bk)}H_{*,i}( \gS,  \gR ; \bk) \lra H_{*,i}^{\cO(1)}(\gS, \gR ; \bk)
  \end{equation*}
is an isomorphism for $i \leq d$. If in addition $\epsilon \colon H_{*,0}(\cO(1);\bk) \to \bk[\bunit]$ is an isomorphism, then it is also a surjection for $i=d+1$.
\end{lemma}
\begin{proof}
Without loss of generality we may suppose that $\gS$ and $\gR$ are cofibrant $\cO(1)$-algebras and $f$ is a cofibration, so in particular $\gS$ and $\gR$ are cofibrant in $\sfC$. We shall apply the obvious relative analogue of the homotopy orbit spectral sequence of Section \ref{sec:homotopy-orbit-ss}, which takes the form
\[E^1_{g, p,q} = H_{g, q} (\cO(1)^{\otimes p} \otimes \gS, \cO(1)^{\otimes p} \otimes \gR ; \bk) \Longrightarrow H_{g,p+q}^{\cO(1)}(\gS,\gR ; \bk).\]
We have $E^1_{g, 0, d} = H_{g,d}(\gS,\gR;\bk)$, and we wish to identify $E^1_{g, 1, d}$. The spectral sequence of Lemma \ref{lem:KanExtCofibrant} shows that the natural map
\[\colim_{(a, b, f) \in \sfH_g} H_d(C_*(\cO(1))(a) \otimes C_*(\gS,\gR)(b)) \lra H_{g, d}(\cO(1) \otimes \gS, \cO(1) \otimes \gR ; \bk)\]
is an isomorphism. As $\cO(1)$ is homologically 0-connective the K{\"u}nneth spectral sequence of Lemma \ref{lem:KunnethSS} shows that the natural map
\[H_0(\cO(1)(a);\bk) \otimes_\bk H_{b,d}(\gS,\gR;\bk) \lra H_d(C_*(\cO(1))(a) \otimes C_*(\gS,\gR)(b))\]
is an isomorphism. Combining these two isomorphisms shows that the natural map
\[H_{*,0}(\cO(1);\bk) \otimes H_{*,d}(\gS,\gR;\bk) \lra H_{*, d}(\cO(1) \otimes \gS, \cO(1) \otimes \gR ; \bk)\]
is an isomorphism, which identifies $E^1_{*, 1, d} \cong H_{*,0}(\cO(1);\bk) \otimes H_{*,d}(\gS, \gR;\bk)$.

Under this isomorphism the differential $d^1 \colon E^1_{*, 1, d} \to E^1_{*, 0, d}$ is identified with the difference of the $\cO(1)$-action map on $H_{*,d}(\gS,\gR;\bk)$ and the augmentation. Thus
\[E^2_{*, 0, d} \cong \bk[\bunit] \otimes_{H_{*,0}(\cO(1);\bk)}H_{*,d}( \gS, \gR ; \bk).\]
This is the only term in total degree $d$, and there are none in total degree less than $d$, which gives the claimed isomorphisms in degrees $i \leq d$.

\vspace{.5em}

To obtain a surjection in degree $i=d+1$, we need to analyse the $E^1$-page in one degree further. Firstly, the argument used above to identify $E^1_{*,1,d}$ generalizes to show that the natural map
\[H_{*,0}(\cO(1);\bk)^{\otimes p} \otimes H_{*,d}(\gS,\gR;\bk) \lra H_{*,d}(\cO(1)^{\otimes p} \otimes \gS,\cO(1)^{\otimes p} \otimes \gR;\bk) = E^1_{*,p,d}\]
is an isomorphism, and to identify the $d^1$-differential with that of the bar complex. Thus we obtain 
\[E^2_{*,p,d} \cong \mr{Tor}^{H_{*,0}(\cO(1);\bk)}_p(\bk[\bunit], H_{*,d}(\gS,\gR;\bk)).\]
Using the assumption that $\epsilon \colon H_{*,0}(\cO(1);\bk) \to \bk[\bunit]$ is an isomorphism, it follows that this vanishes for $p \geq 1$. 

We still have $E^1_{*,0,d+1} = H_{*,d+1}(\gS,\gR;\bk)$, but the entry $E^1_{*,1,d+1} = H_{*,d+1}(\cO(1) \otimes \gS,\cO(1) \otimes \gR;\bk)$ is more complicated. The edge homomorphism in the spectral sequence of Lemma \ref{lem:KanExtCofibrant} gives a natural map
\[\colim_{(a, b, f) \in \sfH_g} H_{d+1}(C_*(\cO(1))(a) \otimes C_*(\gS,\gR)(b)) \lra H_{g, d+1}(\cO(1) \otimes \gS, \cO(1) \otimes \gR ; \bk).\]
Similarly, the K{\"u}nneth spectral sequence of Lemma \ref{lem:KunnethSS} gives a natural map
\[\bigoplus_{j=0,1} H_j(\cO(1)(a);\bk) \otimes_\bk H_{b,d+j}(\gS,\gR;\bk) \lra H_{d+1}(C_*(\cO(1))(a) \otimes C_*(\gS,\gR)(b)).\]
As before, the composition 
\[\bigoplus_{j=0,1} H_j(\cO(1)(a);\bk) \otimes_\bk H_{b,d+j}(\gS,\gR;\bk) \longrightarrow E^1_{*,1,d+1} \overset{d_1}{\lra} E^1_{*,0,d+1}\]
coincides with the difference of the $\cO(1)$-action map and the augmentation, and so $E^2_{*,0,d+1}$ is a quotient of $\bk[\bunit] \otimes_{H_{*,0}(\cO(1);\bk)}H_{*,d+1}( \gS, \gR ; \bk)$. With the vanishing of $E^2_{*, p, d}$ for $p>0$ established above, in degree $d+1$ we have $H_{*,d+1}(\gS,\gR;\bk) \cong E^2_{*, 0, d+1}$, which gives the claimed surjectivity.
\end{proof}

If $\gR$ is a cofibrant $\cO$-algebra, then we can form an absolute Hurewicz map by composing $\bk[\bunit] \otimes_{H_{*,0}(\cO(1);\bk)} -$ applied to the relative Hurewicz map with
\[\bk[\bunit] \otimes_{H_{*,0}(\cO(1);\bk)}H_{*,d}(Q^\cO_{\cO(1)}\gR ; \bk) \lra H_{*,d}(Q^\cO_{\cO(1)}\gR \sslash \cO(1);\bk) \cong H_{*,d}^{\cO}(\gR;\bk).\]
This definition extends easily to relative homology, and to non-cofibrant $\cO$-algebras by taking cofibrant replacements.

\begin{corollary}\label{cor:absolute-Hurewicz}
Let $\cO$ be an augmented non-unitary homologically 0-connective $\Sigma$-cofibrant operad. Let $\gR, \gS \in \Alg_\cO(\sfC)$ be homologically 0-connective. Let $f \colon \gR \to \gS$ be a morphism such that either
\begin{enumerate}[(i)]
\item $\sfG$ is an Artinian groupoid, $\gR$ and $\gS$ are reduced, and $H_{g',d'}(\gS,\gR) = 0$ whenever $(\omega(g'),d') \lneqq (\omega(g),d)$, for some $d \in \bZ$ and $g \in \sfG$, or

\item $\cO$ is the operad associated to an associative ring, and $H_{*,d'}(\gS,\gR) = 0$ whenever $d' < d$, for some $d \in \bZ$.
\end{enumerate}
Then the induced map
  \begin{equation*}
    (\bk[\bunit] \otimes_{H_{*,0}(\cO(1);\bk)} H_{*,i}(\gS, \gR))(g) \lra H_{g,i}^\cO(\gS,  \gR)
  \end{equation*}
  is an isomorphism for $i=d$. If in addition $\epsilon \colon H_{*,0}(\cO(1);\bk) \to \bk[\bunit]$ is an isomorphism, then it is also a surjection for $i=d+1$.
\end{corollary}
\begin{proof}
Without loss of generality we may suppose that $\gS$ and $\gR$ are cofibrant $\cO$-algebras.
In case (i), we apply the right-exact functor $\bk[\bunit] \otimes_{H_{*,0}(\cO(1);\bk)} -$ to the conclusion of Corollary \ref{cor:homology-Hurewicz}, then use Lemma \ref{lem:CoinvariantsOnFirstNonvanishing} on the map $\bL Q^\cO_{\cO(1)}(f)$ and \eqref{eq:AbsIndec} to identify the target with $H_{*,*}^\cO(\gS,  \gR)$. In case (ii), we note that $Q^\cO_{\cO(1)}$ is the identity, so this is just Lemma \ref{lem:CoinvariantsOnFirstNonvanishing}.
\end{proof}

\subsection{Whitehead theorems}\label{sec:whitehead} Having established Hurewicz theorems comparing ordinary homology to $\cO$-homology, we next wish to describe conditions under which the homology of the derived relative or absolute indecomposables detects homology equivalences between $\cO$-algebras. With the methods developed in the last two sections this question could be studied quite generally, but we restrict ourselves to those situations which will be important for our applications.

\subsubsection{Whitehead theorem for relative indecomposables}

The following is immediate from Corollary \ref{cor:homology-Hurewicz}.

\begin{proposition}\label{prop:relwhitehead}
Let $\cO$ be a non-unitary homologically 0-connective $\Sigma$-cofibrant operad. Let $\sfG$ be an Artinian groupoid, and $\gR, \gS \in \Alg_\cO(\sfC)$ be reduced and homologically 0-connective. If $f \colon \gR \to \gS$ is a morphism such that $\bL Q^\cO_{\cO(1)}(f)$ is a homology equivalence, then $f$ is also a homology equivalence.
\end{proposition}

\subsubsection{Whitehead theorem for absolute indecomposables}

For absolute indecomposables, the discussion of Section \ref{sec:an-estimate} means that we need to understand when the functor $\bk[\bunit] \otimes_{H_{*,0}(\cO(1);\bk)} - \colon H_{*,0}(\cO(1);\bk)\text{-}\cat{Mod} \to \cat{Mod}_\bk^\sfG$ detects trivial objects. This seems like a difficult question to answer in general, and we content ourselves with the following condition, which covers all the examples we have in mind.

\begin{lemma}\label{lem:TensorDetectNull}
Suppose that $\sfG$ is Artinian and either
\begin{enumerate}[(i)]
\item $\epsilon \colon H_{*,0}(\cO(1);\bk) \to \bk[\bunit]$ is an isomorphism, or
\item the ideal $\mr{Ker}(\epsilon \colon H_{\bunit_\sfG, 0}(\cO(1);\bk) \to \bk)$ is nilpotent and all $\oplus$-invertible objects are isomorphic to $\bunit_\sfG$.
\end{enumerate}
If $M$ is a $H_{*,0}(\cO(1);\bk)$-module such that $(\bk[\bunit] \otimes_{H_{*,0}(\cO(1);\bk)} M)(g)=0$ for all $g \in \sfG$ satisfying $\omega(g) \leq r$, then $M(g)=0$ for all such $g$ too.
\end{lemma}
\begin{proof}
Define a $H_{*,0}(\cO(1);\bk)$-module $I$ by the exact sequence
\[0 \lra I \lra H_{*,0}(\cO(1);\bk) \lra \bk[\bunit] \lra 0.\]
It follows from the assumption that
\[(I \otimes_{H_{*,0}(\cO(1);\bk)} M)(g) \lra M(g)\]
is surjective for all $g \in \sfG$ such that $\omega(g) \leq r$. 

Under hypothesis (i) we have $I=0$ so it follows that $M(g)=0$, as required. Under hypothesis (ii), suppose for a contradiction that $M(g) \neq 0$ for some $g$ with $\omega(g) \leq r$, and that $\omega(g)$ is minimal with this property. If $a \oplus b \cong g$ with $a$ not a unit then $\omega(b) < \omega(g)$ so $M(b) =0$. On the other hand, if $a$ is a unit then by hypothesis it is isomorphic to $\bunit_\sfG$. Thus we have
\[(I \otimes_{H_{*,0}(\cO(1);\bk)} M)(g) \cong I(\bunit_\sfG) \otimes_{H_{\bunit_\sfG,0}(\cO(1);\bk)} M(g)\]
so $I(\bunit_\sfG) \otimes_{H_{\bunit_\sfG,0}(\cO(1);\bk)} M(g) \to M(g)$ is surjective. As $I(\bunit_\sfG)$ is a nilpotent ideal in $H_{\bunit_\sfG,0}(\cO(1);\bk)$ it follows that $M(g)=0$, a contradiction.
\end{proof}

\begin{proposition}\label{prop:AbsWhitehead}
Let $\cO$ be an augmented non-unitary homologically 0-connective $\Sigma$-cofibrant operad, such that $\sfG$ and $\cO$ satisfy the hypotheses of Lemma \ref{lem:TensorDetectNull}. Let $f \colon \gR \to \gS$ be a morphism between homologically 0-connective $\cO$-algebras such that $H^\cO_{*,d'}(\gS, \gR)=0$ for all $d' \leq d$.
If either
\begin{enumerate}[(i)]
\item  $\gR$ and $\gS$ are reduced, or

\item $\cO$ is the operad associated to an associative ring,
\end{enumerate}
then $H_{*,d'}(\gS, \gR)=0$ for all $d' \leq d$ too.
\end{proposition}
\begin{proof}
Suppose for a contradiction that $H_{g',d'}(\gS, \gR;\bk) \neq 0$ for some $g'$ and $d' \leq d$; we may further suppose that $(\omega(g'), d')$ is minimal with respect to the partial order $\lneqq$, so that $H_{g'',d''}(\gS, \gR;\bk)=0$ for all $(\omega(g''),d'') \lneqq (\omega(g'),d')$. By Corollary \ref{cor:absolute-Hurewicz} we then have that
  \begin{equation*}
    (\bk[\bunit] \otimes_{H_{*,0}(\cO(1);\bk)} H_{*,d'}(\gS, \gR;\bk))(g') \lra H_{g',d'}^\cO(\gS,  \gR;\bk) = 0
  \end{equation*}
  is an isomorphism. The same holds for all $g''$ with $\omega(g'') \leq \omega(g')$, so by Lemma \ref{lem:TensorDetectNull} we have that $H_{g',d'}(\gS, \gR;\bk)=0$, a contradiction.
\end{proof}

In particular, taking $d=\infty$ it follows that $\cO$-homology detects homology equivalences, as follows.

\begin{corollary}\label{cor:AbsWhiteheadInfty}
Under the same hypotheses, if $\bL Q^\cO(f)$ is a homology equivalence then $f$ is also a homology equivalence.
\end{corollary}

\subsection{CW approximation in the semistable case}\label{sec:additive-case}

In this section we use the Hurewicz results of Section \ref{sec:an-estimate} to prove the existence of (minimal) CW approximations. Until now we have had no need to consider homotopy groups\index{homotopy!group} of objects of $\sfS$, but it will now be essential to do so. For this we shall suppose that $\sfS = \sfS_*$ is pointed. Then, using the copowering $- \wedge - \colon \cat{sSet}_* \times \sfS \to \sfS$ we obtain for any pointed simplicial set $X$ an object $s_+(X) \coloneqq X \wedge \bunit$, giving a pointed version of the usual map $s \colon \cat{sSet} \to \sfS$ satisfying $s_+(X) \cong s(X)/s(*)$. Write $i_d \colon S^{d-1} \to D^d$ for the inclusion. To define relative homotopy groups we work in the category $\sfS^{[1]}$ of arrows in $\sfS$, which by abuse of notation we consider as pairs, and for a morphism $f: X \to Y$ in $\sfS$ and a $d \in \bN$ we set
\[\pi_d(Y, X) = \pi_d(f) \coloneqq \cat{Ho}(\sfS^{[1]})(s_+(i_d), f).\]
We define absolute homotopy groups as $\pi_d(X) \coloneqq \pi_d(X, *)$, which may be identified with $\cat{Ho}(\sfS)(s_+(S^{d}), X)$. As usual the homotopy cogroup structure on $S^d \in \cat{sSet}_*$ makes $\pi_d(X)$ into a group for $d \geq 1$, which is abelian for $d \geq 2$, and similarly for relative homotopy groups. There is a long exact sequence for relative homotopy, developed in the usual way.

The functor $\Sigma(-) = s_+(S^1) \otimes - = S^1 \wedge - \colon \sfS \to \sfS$ has a right adjoint $\Omega(-) = \Hom_{\sfS}(s_+(S^1),-) \colon \sfS \to \sfS$ and these form a Quillen adjunction: following Heller \cite{HellerSemistab} we call $\sfS$ \emph{semistable}\index{semistable} if the derived unit of this adjunction is a natural \emph{isomorphism} between the identity functor and
\[\cat{Ho}(\sfS) \overset{\bL \Sigma}\lra \cat{Ho}(\sfS) \overset{\bR \Omega}\lra \cat{Ho}(\sfS).\]
In particular we have weak equivalences $Y \simeq (\bR\Omega^n)(\bL \Sigma^n(Y))$ for every $n \in \bN$ and so weak equivalences of derived mapping spaces
\[\mr{Map}_\sfS(X, Y) \simeq \Omega^n \mr{Map}_\sfS(X, \bL \Sigma^n(Y)).\]
This can probably be used to enrich $\sfS$ in infinite loop spaces or even spectra, but we shall settle for observing that it yields an enrichment of $\cat{Ho}(\sfS)$ in abelian groups. Furthermore semistability implies that $\bL \Sigma  \colon \cat{Ho}(\sfS) \to \cat{Ho}(\sfS)$ is full and faithful, so in particular that the maps
\[\bL \Sigma(-) \colon \pi_d(X) \lra \pi_{d+1}(\bL \Sigma(X))\]
are bijections; the same follows for relative homotopy groups.

The first point endows each $\pi_d(f)$ with the structure of an abelian group (which agrees with the old structure when it is defined, by the Eckmann--Hilton argument), and the second allows us to extend the definition of relative homotopy groups above to
\[\pi_{-d}(f) \coloneqq \pi_0(\bL\Sigma^d(f)) \quad \text{ for } d \in \bN.\]
We then say that a morphism $f \colon X \to Y$ in $\sfS$ is $c$-connective if $\pi_d(f)=0$ for all $d < c$, and an object $X$ is $c$-connective if $\pi_d(X)=0$ for all $d < c$.  As with homology groups, for a map $f \colon X \to Y$ in $\sfC=\sfS^\sfG$ we can then define
\[\pi_{g,d}(Y, X) \coloneqq \pi_d(Y(g), X(g))\]
for any $d \in \bZ$, and so define $c$-connectivity for any abstract connectivity $c$ on $\sfG$.

As we wish to use homology rather than homotopy to detect cells, it is vital that we work in a context where homology may be used to detect homotopical connectivity of maps. By the second axiom of a singular chain functor on $\sfS$, the homology groups ${H}_i(D^d \wedge \bunit, S^{d-1} \wedge \bunit;\bk)$ are naturally isomorphic to the ordinary homology groups ${H}_i(D^d, S^{d-1};\bk)$, and in particular there is a canonical generator $u_d \in {H}_d(D^d \wedge \bunit, S^{d-1} \wedge \bunit;\bk)$. Given a morphism $X \to Y$ in $\sfS$, functoriality defines a relative Hurewicz map\index{Hurewicz map}
\begin{equation}\label{eq:RElHurewiczMap}
h \colon \pi_d(Y, X) \lra {H}_d(Y, X;\bk)
\end{equation}
for each $d \in \bN$. As homology has a suspension isomorphism (this follows from Lemma \ref{lem:KunnethFormula} (ii)) this definition extends to all $d \in \bZ$.

\begin{axiom}\label{axiom:Hurewicz}
The category $\sfS$ is pointed and semistable, and the weak equivalences are precisely those maps which induce a bijection on $\pi_d(-)$ for all $d \in \bZ$. Furthermore for any $d \in \bZ$ and any map $f \colon X \to Y$ in $\sfS$ such that $\pi_i(Y, X)=0$ for all $i < d$, the relative Hurewicz map \eqref{eq:RElHurewiczMap} is a bijection.
\end{axiom}

This property holds in $\cat{sMod}_\bk$ (with $\bk$ coefficients), in $\cat{Sp}^\Sigma$ (with $\bZ$ coefficients), and in the category of $R$-modules in $\cat{Sp}^\Sigma$ for a fixed commutative ring spectrum $R$ (with $\pi_0(R)$ coefficients). In the latter two cases we shall write $\bk$ for $\bZ$ and $\pi_0(R)$ respectively.

\subsubsection{CW approximation}
We use much of the terminology from Section \ref{sec:cw-algebras}, but recall some definitions for the convenience of the reader. Definition \ref{def:relative-cw-structure} defined a relative CW-structure on a map $f \colon \gR \to \gS$ to be an object $\mr{sk}(f)$ in $\Alg_\cO(\sfC^{\bZ_{\leq}})$ which is the colimit in $\Alg_\cO(\sfC^{\bZ_{\leq}})$ of a diagram
\[0_*(\gR) = \mr{sk}_{-1}(f) \lra \mr{sk}_0(f) \lra \mr{sk}_1(f) \lra \cdots\]
where $f_d \colon \mr{sk}_{d-1}(f) \to \mr{sk}_d(f)$ comes with the structure of a filtered CW attachment of dimension $d$, and a factorisation
\[f \colon \gR \lra \colim \mr{sk}(f) \overset{\cong}\lra \gS.\]
Here the second map is, crucially, an isomorphism. A relative CW approximation is the homotopical analogue of this definition.

\begin{definition}
A \emph{relative CW approximation}\index{CW-!approximation} of a map $f \colon \gR \to \gS$ of $\cO$-algebras is an object $\mr{sk}(f)$ in $\Alg_\cO(\sfC^{\bZ_{\leq}})$ as above and a factorization 
\[f \colon \gR \lra \colim \mr{sk}(f) \overset{\sim}\lra \gS\]
as a relative CW algebra followed by a weak equivalence.
\end{definition}

The main theorem of this section we give conditions under which a map admits a relative CW approximation, and furthermore show that then the dimensions of the cells involved can be constrained by the derived absolute $\cO$-indecomposables. 

\begin{theorem}\label{thm:MinCellStr-additive}
Let $\sfS$ be a pointed category satisfying the axioms of Section \ref{sec:axioms-of-model-cats}, and Axiom \ref{axiom:Hurewicz}. Let $\cO$ be an augmented non-unitary homologically 0-connective $\Sigma$-cofibrant operad in $\sfC = \sfS^\sfG$, such that $\sfG$ and $\cO$ satisfy the hypotheses of Lemma \ref{lem:TensorDetectNull}. Let $f \colon \gR \to \gS$ be a morphism between homologically 0-connective $\cO$-algebras, such that either
\begin{enumerate}[(i)]
\item\label{it:MinCellStr-additive:1}  $\gR$ and $\gS$ are reduced, or
\item\label{it:MinCellStr-additive:2} $\cO$ is the operad associated to an associative ring.
\end{enumerate}
Let $c \colon \sfG \to [-\infty,\infty]_\geq$ be an abstract connectivity such that $\smash{H_{g,d}^\cO}(\gS, \gR;\bk)=0$ for $d < c(g)$. Then there exists a relative CW approximation $f \colon \gR \to \colim \mr{sk}(f) \xrightarrow{\sim} \gS$ where $\mr{sk}(f)$ has no $(g,d)$-cells with $d < c(g)$.
\end{theorem}

The construction in the proof below will in fact give a \emph{minimal cell structure}, i.e.\ one having the smallest possible number of cells in a given bidegree. 
To make this precise, suppose that $\bk$ is a field and recall that we have defined the $\cO$-Betti numbers $b_{g,d}^{\cO}(\gS,\gR) \coloneqq \dim_\bk H_{g,d}^\cO(\gS,\gR;\bk) \in \bN \cup\{\infty\}$. Then the relative CW $\cO$-algebra $\gR \to \colim \mr{sk}(f)\overset{\sim}\to \gS$ produced by this theorem will have precisely $b_{g,d}^{\cO}(\gS,\gR)$ $(g,d)$-cells.

\begin{proof}
In case (\ref{it:MinCellStr-additive:1}), by Corollary \ref{cor:absolute-Hurewicz} we may assume that $c(g)>0$ for $g \in \sfG^\times$, as $\gS$ and $\gR$ are reduced. We shall prove by induction over $\epsilon$, the \emph{dimension}, starting at $\epsilon=-1$, the following statement:
	\begin{quote}
		There exists a factorization of $0_*(f) \colon 0_*(\gR) \to 0_*(\gS)$:
		\[0_*(\gR) = \mr{sk}_{-1}(f) \overset{h_0}\lra \mr{sk}_0(f) \overset{h_1}\lra \cdots \overset{h_\epsilon}\lra  \mr{sk}_\epsilon(f) \overset{f_\epsilon}\lra 0_*(\gS)\]
		with the following properties for all $0 \leq e \leq \epsilon$:	
		\begin{enumerate}[\indent (a)]
             \item \label{enum:ind-range} $H^\cO_{*,d}(\gS,\colim \mr{sk}_e(f)) = 0$ for all $d$ satisfying $d \leq e$,
			\item \label{enum:ind-cells} $h_e \colon \mr{sk}_{e-1}(f) \to \mr{sk}_e(f)$ comes with the structure of a filtered CW attachment of dimension $e$, and only has cells attached to those $g$ with $c(g) \leq e$.
		\end{enumerate}
	\end{quote}
Supposing we have done so, then 
\[H_{g,d}^\cO(\gS, \colim \mr{sk}(f);\bk) = \colim_\epsilon H_{g,d}^\cO(\gS, \colim \mr{sk}_\epsilon(f);\bk)\]
vanishes for all $g$ and $d$ by (\ref{enum:ind-range}), so the induced map $f_\infty \colon \colim \mr{sk}(f) \to \gS$ induces an isomorphism on $\cO$-homology. In case (\ref{it:MinCellStr-additive:1}) it follows from (\ref{enum:ind-cells}) that $\colim \mr{sk}(f)$ is obtained from $\gR$ by attaching cells with no 0-cells attached to $g \in \sfG^\times$: as $\gR$ is reduced so is $\colim \mr{sk}(f)$. Thus in case (\ref{it:MinCellStr-additive:1}) or (\ref{it:MinCellStr-additive:2}), Corollary \ref{cor:AbsWhiteheadInfty} applies, and shows that the map $f_\infty$ is a homology isomorphism, so also a $\pi_*$-isomorphism and hence a weak equivalence by Axiom \ref{axiom:Hurewicz}.

It remains then to prove the above statements by induction; the case $\epsilon=-1$ is immediate. For the induction step let us write $\gZ_e = \colim \mr{sk}_{e}(f)$ to ease notation. Assuming the statement holds for $(\epsilon-1)$, we have $H^\cO_{*,d}(\gS, \gZ_{\epsilon-1}) = 0$ for all $d$ satisfying $d \leq \epsilon-1$, and hence by Proposition \ref{prop:AbsWhitehead} we have $H_{*,d}(\gS,\gZ_{\epsilon-1}) = 0$ for all such $d$.

\vspace{1ex}

\noindent\textbf{Claim:} The Hurewicz map
\[\pi_{*,\epsilon}(\gS, \gZ_{\epsilon-1}) \lra H^\cO_{*, \epsilon}(\gS, \gZ_{\epsilon-1};\bk)\]
is surjective.
\begin{proof}[Proof of claim]
It follows from Corollary \ref{cor:homology-Hurewicz} that the Hurewicz map
\begin{equation}\label{eq:Hurewicz55}
H_{*,d}(\gS, \gZ_{\epsilon-1};\bk) \lra H_{*,d}(\bL Q^\cO_{\cO(1)}\gS, \bL Q^\cO_{\cO(1)}\gZ_{\epsilon-1};\bk)
\end{equation}
is surjective for $d \leq \epsilon$, and an isomorphism for $d \leq \epsilon-1$. Thus the target also vanishes for $d \leq \epsilon-1$, and so by Lemma \ref{lem:CoinvariantsOnFirstNonvanishing} applied to this pair the map
\[\bk[\bunit] \otimes_{H_{*,0}(\cO(1);\bk)} H_{*,\epsilon}(\bL Q^\cO_{\cO(1)}\gS, \bL Q^\cO_{\cO(1)}\gZ_{\epsilon-1};\bk) \lra H^\cO_{*, \epsilon}(\gS, \gZ_{\epsilon-1};\bk)\]
is an isomorphism. Applying $\bk[\bunit] \otimes_{H_{*,0}(\cO(1);\bk)} -$ to the map \eqref{eq:Hurewicz55}, it follows that
\[\bk[\bunit] \otimes_{H_{*,0}(\cO(1);\bk)} H_{*,\epsilon}(\gS, \gZ_{\epsilon-1};\bk) \lra H^\cO_{*, \epsilon}(\gS, \gZ_{\epsilon-1};\bk)\]
is a surjection. Combining this with the surjections
\begin{align*}
\pi_{*,\epsilon}(\gS, \gZ_{\epsilon-1}) & \lra H_{*,\epsilon}(\gS, \gZ_{\epsilon-1};\bk)\\
& \lra \bk[\bunit] \otimes_{H_{*,0}(\cO(1);\bk)} H_{*,\epsilon}(\gS, \gZ_{\epsilon-1};\bk),
\end{align*}
where the first map is surjective (in fact, an isomorphism) by Axiom \ref{axiom:Hurewicz}, proves the claim.
\end{proof}

Thus for each $g \in \sfG$, we may choose a set of maps
\[\{E_\alpha \colon (D^\epsilon, \partial D^\epsilon) \to (\gS(g), \gZ_{\epsilon-1}(g))\}_{\alpha \in I_{g, \epsilon,}}\]
whose images under the above Hurewicz map generate $H^\cO_{g, \epsilon}(\gS, \gZ_{\epsilon-1};\bk)$ as a $\bk$-module. 

In the exact sequence
\[\cdots \lra H^\cO_{g, \epsilon}(\gS, \gR;\bk) \lra H^\cO_{g, \epsilon}(\gS, \gZ_{\epsilon-1};\bk)\lra H^\cO_{g, \epsilon-1}(\gZ_{\epsilon-1}, \gR;\bk) \lra \cdots\]
of the triple $\gR \subset \gZ_{\epsilon-1} \subset \gS$, the left-hand term vanishes for $\epsilon < c(g)$ by assumption, and the pair $(\gZ_{\epsilon-1}, \gR)$ only has relative $(g, \epsilon-1)$-cells if $c(g) \leq \epsilon-1$ by part (\ref{enum:ind-cells}) of the inductive assumption, so also vanishes for $\epsilon < c(g)$. Thus the middle term vanishes for $\epsilon < c(g)$, so we can take $I_{g, \epsilon} = \varnothing$ for $\epsilon < c(g)$.

In order to use the maps $e_\alpha \coloneqq E_\alpha\vert_{\partial D^\epsilon}$ to attach $(g, \epsilon)$-CW-cells to the filtered object $\mr{sk}_{\epsilon-1}(f)$, we must lift them along
\[\mr{sk}_{\epsilon-1}(f)(g, \epsilon-1) \lra \colim \mr{sk}_{\epsilon-1}(f)(g) = \gZ_{\epsilon-1}\]
up to homotopy.

\vspace{1ex}

\noindent\textbf{Claim:} This map is $(\epsilon-1)$-connected.
\begin{proof}[Proof of claim]
By Axiom \ref{axiom:Hurewicz} it is enough to check that it is homologically $(\epsilon-1)$-connected, and we may do this by analysing the associated graded of the filtered object $\mr{sk}_{\epsilon-1}(f)$. By Theorem \ref{thm:associated-graded-skeletal} there is an isomorphism
	\[\grr(\mr{sk}_{\epsilon-1}(f)) \cong 0_*(\gR_+) \vee^\cO F^\cO\left(\bigvee_{d \leq \epsilon-1} \bigvee_{\alpha \in I_{d}} d_*(S^{g_\alpha, d}_\alpha)\right),\]
	and we need to show that in each grading $n \geq \epsilon$ it is homologically $\epsilon$-connective, i.e.\ homologically $(\epsilon-1)$-connected.
	
	To do so, using that $\gR \in \Alg_\cO(\sfC)$ is cofibrant, so $0_*(\gR_+)$ is too, we can consider the above coproduct of $\cO$-algebras as a derived coproduct, and use the analogue of the simplicial formula in Section \ref{sec:derived-cell-simplicial} to describe it. This gives
	\[\grr(\mr{sk}_{\epsilon-1}(f)) \simeq \left\| [p] \mapsto F^\cO\left(\cO^p(0_*(\gR_+)) \vee \bigvee_{d \leq \epsilon-1} \bigvee_{\alpha \in I_{d}} d_*(S^{g_\alpha, d}_\alpha)\right) \right\|.\]
	Let us consider the abstract connectivity \emph{in the $\bZ_{=}$-grading direction} given by $\ell(n)=n$, which satisfies $\ell *\ell \geq \ell$. That is, an object $X \in \cat{C}^{\bZ_{=}}$ is homologically $\ell$-connective if $H_{g,d,n}(X)=0$ whenever $d < \ell(n)=n$. First observe that if $X \in \sfC^{\bZ_{=}}$ is a homologically $\ell$-connective object then $\cO(X) = \bigvee_{k \geq 1} \cO(k)_+ \wedge_{G_k} X^{\otimes k}$ is too, using that $X^{\otimes k}$ is homologically $(\ell^{* k} \geq \ell)$-connective and then applying the homotopy orbit spectral sequence. (Here we have used that $\cO$ is non-unitary, homologically 0-connective, and $\Sigma$-cofibrant.)
	
	The object $0_*(\gR_+)$ is supported in grading 0, and $\gR$ is homologically 0-connective by assumption, so $0_*(\gR_+)$ is homologically $\ell$-connective and hence by the observation above $\cO^p(0_*(\gR_+))$ is homologically $\ell$-connective too. The object $d_*(S^{g_\alpha, d}_\alpha)$ is supported in grading $d$ and $S^{g_\alpha, d}_\alpha$ is $d$-connective, so $d_*(S^{g_\alpha, d}_\alpha)$ is homologically $\ell$-connective. Thus $\cO^p(0_*(\gR_+)) \vee \bigvee_{d \leq \epsilon-1} \bigvee_{\alpha \in I_{d}} d_*(S^{g_\alpha, d}_\alpha)$ is homologically $\ell$-connective, and by the observation above the free $\cO$-algebra on it is too. Thus the semi-simplicial object is levelwise homologically $\ell$-connective, so it follows from the geometric realization spectral sequence (Theorem \ref{thm:geom-rel-ss-thick}) that $\grr(\mr{sk}_{\epsilon-1}(f))$ is homologically $\ell$-connective too, which is precisely what we required.
\end{proof}

Thus we may lift the maps $e_\alpha$ (up to homotopy) to $\mr{sk}_{\epsilon-1}(f)(g, \epsilon-1)$, use them to attach filtered $(g, \epsilon)$-cells to form $\mr{sk}_{\epsilon}(f)$, and use the corresponding $E_\alpha$ to extend the map $f_{\epsilon-1}$ to a map $f_\epsilon \colon \mr{sk}_{\epsilon}(f) \to \const(\gS)$. This satisfies property (\ref{enum:ind-cells}) by construction. The long exact sequence of the triple $\gZ_{\epsilon-1} \to \gZ_{\epsilon} \to \gS$ takes the form
\begin{equation*}
\begin{tikzcd}
     & 0 \rar & H_{g, \epsilon+1}^{\cO}(\gS, \gZ_{\epsilon-1};\bk) \rar \ar[draw=none]{d}[name=X, anchor=center]{} & H_{g, \epsilon+1}^{\cO}(\gS, \gZ_{\epsilon};\bk)
                \ar[rounded corners,
                to path={ -- ([xshift=2ex]\tikztostart.east)
                	|- (X.center) \tikztonodes
                	-| ([xshift=-2ex]\tikztotarget.west)
                	-- (\tikztotarget)}]{dll} \\
             & {\bigoplus\limits_{\alpha \in I_{g, \epsilon}}\mathbbm{k}\{E_\alpha\}} \rar[two heads] & H_{g, \epsilon}^{\cO}(\gS, \gZ_{\epsilon-1};\bk) \rar \ar[draw=none]{d}[name=Y, anchor=center]{} & H_{g, \epsilon}^{\cO}(\gS, \gZ_{\epsilon};\bk)
                \ar[rounded corners,
                to path={ -- ([xshift=2ex]\tikztostart.east)
                	|- (Y.center) \tikztonodes
                	-| ([xshift=-2ex]\tikztotarget.west)
                	-- (\tikztotarget)}]{dll} \\
             & 0 \rar & H_{g, \epsilon-1}^{\cO}(\gS, \gZ_{\epsilon-1};\bk)=0 \rar & H_{g, \epsilon-1}^{\cO}(\gS, \gZ_{\epsilon};\bk) \rar& 0
\end{tikzcd}
\end{equation*}
as the pair $(\gZ_{\epsilon}, \gZ_{\epsilon-1})$ only has relative $\cO$-algebra cells of dimension $\epsilon$. It follows that $H_{g, d}^{\cO}(\gS, \gZ_{\epsilon};\bk)=0$ for all $d \leq \epsilon$, which verifies property (\ref{enum:ind-range}). 
\end{proof}



\part{$E_k$-algebras} \label{part:ek} In this third part we will apply the results of Parts \ref{part:category} and \ref{part:homotopy} to $E_k$-algebras. Next we develop a number of tools particular to this setting, most importantly a result relating $Q^{E_k}_\bL$ to $k$-fold bar constructions.

In this part we will assume that the axioms of Sections \ref{sec:axioms-of-cats} and \ref{sec:axioms-of-model-cats} hold unless mentioned otherwise:
\begin{itemize}
	\item Axiom \ref{axiom:simplicial-enriched}: $\cat{S}$ is simplicially enriched.
	\item Axiom \ref{axiom:enriched-complete-cocomplete}: $\cat{S}$ is complete and cocomplete in an enriched sense.
	\item Axiom \ref{axiom.cat-monoidal}: $\cat{S}$ has a simplicially enriched closed $k$-monoidal structure, closed on both sides if $k=1$.
	\item Axiom \ref{axiom.CofGen}: $\cat{S}$ has a cofibrantly generated model structure.
	\item Axiom \ref{axiom.model-monoidalsimplicial}: this model structure is monoidal, simplicial, and $\bunit$ is cofibrant.
\end{itemize}

The simplicial monads associated to the little $k$-cubes operads will also satisfy the axioms of Section \ref{sec:homotopy-theory-algebras}:
\begin{itemize}
	\item Axiom \ref{axiom:monad-proj}: the projective model structure on $\Alg_T(\sfC)$ exists and $U^T$ preserves (trivial) cofibrations between cofibrant objects.
	\item Axiom \ref{axiom:monad-geomrel}: $T$ preserves geometric realization.
\end{itemize}

\section{$E_k$-algebras and $E_1$-modules}\label{sec:EkAlgAndMod}

In Section \ref{sec:ekoperads} we specialize the theory of Parts \ref{part:category} and \ref{part:homotopy} to the monad associated to the little $k$-cubes operad. There are two instances of this operad: unitary and non-unitary. For most of our work we shall consider algebras for the non-unitary little $k$-cubes operad, but we shall occasionally have cause to consider algebras over the unitary operad obtained by formally adjoining a unit. In Section \ref{sec:modules} we discuss some technical points of the theory of $E_1$-modules over $E_k$-algebras.

\subsection{$E_k$- and $E_k^+$-operads} \label{sec:ekoperads} In this section we define the operads which are the main subject of this paper.

\newglossaryentry{ck}{%
	name={\ensuremath{\cC_k}},
	description={Non-unital little $k$-cubes operad},
type=symbols
}
\newglossaryentry{ckplus}{%
	name={\ensuremath{\cC_k^+}},
	description={Little $k$-cubes operad},
type=symbols
}
\newglossaryentry{ek}{%
	name={\ensuremath{E_k}},
	description={Monad associated to $\cC_k$},
type=symbols
}
\newglossaryentry{ekplus}{%
name={\ensuremath{E_k^+}},
description={Monad associated to $\cC_k^+$},
type=symbols
}
\newglossaryentry{embrect}{%
	name={\ensuremath{\mr{Emb}^\mr{rect}}},
	description={Space of rectilinear embeddings},
type=symbols
}
\subsubsection{The $E_k$- and $E_k^+$-operads in the symmetric monoidal case} \label{sec:ekoperadssym} We start by defining the little $k$-cubes operad $\gls{ckplus}$, our choice of a unitary $E_k$-operad. We will write down an operad in $\cat{Top}$, but implicitly take the singular simplicial set to obtain an operad in $\cat{sSet}$ denoted the same way. The following definition uses the notion of an rectilinear embedding. An embedding $I^k \hookrightarrow I^k$ is \emph{rectilinear} if it is of the form 
\[(t_1,\ldots,t_k) \longmapsto \left((b_1-a_1)t_1+a_1,\ldots,(b_k-a_k)t_k+a_k\right)\]
for $a_i<b_i$ (note that necessarily we must have $0 \leq a_i < b_i \leq 1$ in this case, but later we shall encounter cubes of other sizes).  The set of rectilinear embeddings is topologized as a subspace of $\bR^{2k}$.

\begin{definition}\label{def:ekplus} The underlying symmetric sequence of the \emph{unitary little $k$-cubes operad}\index{operad!little cubes} $\cC_k^+$ is given by
	\begin{align*} \cC_k^+(n) \coloneqq \gls{embrect}(\sqcup_n I^k, I^k) \end{align*}
	where $\mr{Emb}^\mr{rect}(\sqcup_n I^k, I^k)$ denotes the space of $n$-tuples of rectilinear embeddings $e_1, e_2, \ldots, e_n \colon I^k \to I^k$ with disjoint interiors (see Figure \ref{fig:little-cube-example} for an example). This is a monoid for the composition product using the map induced by composition of rectilinear embeddings (see Figure \ref{fig:little-cube-comp-example} for an example)
	\[\cC_k(n) \times \cC_k(k_1) \times \cdots \times \cC_k(k_n) \lra \cC_k(k_1+\cdots+k_n),\]
	and unit $\ast \to \cC_k(1)$ given by the identity map $I^k \hookrightarrow I^k$.
\end{definition}

\begin{figure}[h]
	\centering
\begin{tikzpicture}

\draw[fill = Mahogany!10!white] (3,2) rectangle (4,4);
\draw[fill = Mahogany!10!white] (0.5,1) rectangle (2,3);
\draw[fill = Mahogany!10!white] (2.5,0.2) rectangle (4.5,1.8);

\draw (0,0) rectangle (5,5);

\node at (1.25,2) {$e_1$};
\node at (3.5,3) {$e_2$};
\node at (3.5,1) {$e_3$};

\node[right] at (5.1,2.5) {$\in \cC_2(3)$};
\end{tikzpicture}
	\caption{An element of $\cC_2(3)$.}
	\label{fig:little-cube-example}
\end{figure}

\begin{figure}
	\centering
	\begin{tikzpicture}
	
	\draw (-3.5,0) rectangle (-.5,3);
	\draw (0,0) rectangle (3,3);
	\draw (3.5,0) rectangle (6.5,3);
	
	\draw[fill = Mahogany!10!white] (-3.25,0.5) rectangle (-1.75,2.5);
	\draw[fill = Mahogany!10!white] (-1.75,0.25) rectangle (-.75,2.25);
	\node at (-2.5,1.5) {$e_1$};
	\node at (-1.25,1.25) {$e_2$};
	
	\draw[fill = Mahogany!10!white] (0.5,1) rectangle (2,3);
	\draw[fill = Mahogany!10!white] (3.75,0.5) rectangle (5,1.5);
	\draw[fill = Mahogany!10!white] (4,2) rectangle (6,3);
	\draw[fill = Mahogany!10!white] (5.25,0.5) rectangle (6.25,1.5);

	\node at (1.5,3.6) {$\cC_2(2) \times \cC_2(1) \times \cC_2(3)$};
	\node at (1.5,3.2) {\rotatebox{90}{$\in$}};
	
	\draw[decorate,decoration={brace,amplitude=5pt}] (-3.75,-0.25) -- (-3.75,3.25);
	\draw[decorate,decoration={brace,amplitude=5pt,mirror}] (6.75,-0.25) -- (6.75,3.25);
	
	\draw [|->] (1.5,-0.5) --  ++(0,-.5);
	
	\draw (0,-4.5) rectangle (3,-1.5);
	\fill[pattern=north west lines, pattern color=Mahogany!30!white] (0.25,-4) rectangle (1.75,-2);
	\fill[pattern=north west lines, pattern color=Mahogany!30!white] (1.75,-4.25) rectangle (2.75,-2.25);
	\draw[fill = Mahogany!10!white] (0.5*0.5+0.25,1*0.667+0.5-4.5) rectangle (2*0.5+0.25,3*0.666+0.5-4.5);
	\draw[fill = Mahogany!10!white] (0.25*0.333+1.75,0.5*0.666+.25-4.5) rectangle (1.5*0.333+1.75,1.5*0.666+.25-4.5);
	\draw[fill = Mahogany!10!white] (0.5*0.333+1.75,2*0.666+.25-4.5) rectangle (2.5*0.333+1.75,3*0.666+.25-4.5);
	\draw[fill = Mahogany!10!white] (1.75*0.333+1.75,0.5*0.666+.25-4.5) rectangle (2.75*0.333+1.75,1.5*0.666+.25-4.5);
	
	\node[right] at (3.15,-3) {$\in \cC_2(4)$};
	
	\end{tikzpicture}
	\caption{An example of composition in $\cC_2$. We have left out the labels on the inner cubes for readability.}
	\label{fig:little-cube-comp-example}
\end{figure}

As the notation suggests, this operad is isomorphic to the unitalization of a non-unitary operad $\gls{ck}$, as described in Section \ref{sec:non-unitary}. The operad $\cC_k$ will be our choice of an $E_k$-operad.

\begin{definition}
The underlying symmetric sequence of the \emph{non-unitary little $k$-cubes operad} $\cC_k$ is given by
\begin{align*} 
\cC_k(n) \coloneqq \begin{cases} \varnothing &  \text{if $n=0$,} \\
	\mr{Emb}^\mr{rect}(\sqcup_n I^k, I^k) & \text{otherwise,} \end{cases}
\end{align*}
with operad structure restricted from that of $\cC_k^+$.
\end{definition}

\begin{remark}\label{rem:comparing-ek-operads} We have specified that $\cC_k^+$ and $\cC_k$ are our choices of unitary and non-unitary $E_k$-operads; these operads are $\Sigma$-cofibrant. Another operad $\cE_k$ is said to be an \emph{$E_k$-operad}\index{operad!$E_k$-} if it is $\Sigma$-cofibrant and there is a zigzag of weak equivalences of operads between $\cE_k$ and $\cC_k$ (this zigzag may always be taken to consist of $\Sigma$-cofibrant operads, using the existence of a model structure on operads in $\cat{sSet}$ \cite[Theorem 3.2]{BergerMoerdijk}). The notion of an \emph{$E_k^+$ operad}\index{operad!$E_k^+$-} is defined analogously.
	
Some other choices $\cE_k$ of $E_k$- and $E_k^+$-operads are just minor variations of the definitions, e.g.~the operad of little $k$-disks, while some seem quite different, e.g.~the McClure--Smith $E_2$-operad \cite{McClure-Smith}. As long as $\cE_k$ is $\Sigma$-cofibrant, an algebra over $\cE_k$ is naturally weakly equivalent to an algebra over the operad $\cC_k$ described above. In fact, the categories of $\cC_k$- and $\cE_k$-algebras with the projective model structures are Quillen equivalent \cite[Theorem 4.1]{BergerMoerdijk2}.\end{remark}

\newglossaryentry{gek}{%
	name={\ensuremath{\gE_k}},
	description={Free $E_k$-algebra functor},
type=symbols
}
\newglossaryentry{gekplus}{%
name={\ensuremath{\gE_k^+}},
description={Free $E_k^+$-algebra functor},
type=symbols
}
The operads $\cC_k^+$ and $\cC_k$ may be used to define monads $\gls{ekplus}$ and $\gls{ek}$ on any symmetric monoidal category $\sfC$ satisfying the axioms of Section \ref{sec:axioms-of-cats}:
\begin{align*}
E_k^+ &\colon X \longmapsto \bigsqcup_{n \geq 0} \cC_k^+(n) \times_{\fS_n} X^{\otimes n},\\
E_k &\colon X \longmapsto \bigsqcup_{n \geq 1} \cC_k(n) \times_{\fS_n} X^{\otimes n}.
\end{align*}
In order to align with the notational convention we have employed that structured objects such as algebras over a monad are displayed in bold font, and in order to distinguish free algebras from the values of the monad, we adopt the notation
\[\gls{gekplus}(X) \coloneqq F^{E_k^+}(X) \quad\text{ and }\quad \gls{gek}(X) \coloneqq F^{E_k}(X).\]

\newglossaryentry{qek}{%
	name={\ensuremath{Q^{E_k}_\bL}},
	description={Derived $E_k$-indecomposables},
type=symbols
}
\newglossaryentry{decek}{%
name={\ensuremath{\Dec_\bL^{E_k}}},
description={Derived $E_k$-decomposables},
type=symbols
}
Both of these monads are sifted by Corollary \ref{cor:monad-from-operad-sifted}. The latter, $E_k$, has a canonical augmentation $\epsilon \colon E_k \to +$ as in Section \ref{sec:simplicial-operads}. Hence we may define the (absolute) derived indecomposables of $\gR \in \Alg_{E_k}(\sfC)$ by the derived functor
\[\gls{qek}(\gR) \coloneqq \bL\epsilon_*(\gR) \in \sfC_\ast,\]
as described in Section \ref{sec:deriv-indec-1}. We may also define the derived relative indecomposables $\bL Q^{E_k}_{E_k(1)}(\gR)$, but as $E_k(1) \simeq *$ the natural map $\bL Q^{E_k}_{E_k(1)}(\gR) \to Q^{E_k}_\bL(\gR)$ is a weak equivalence. For the same reason we abbreviate $\gls{decek}(\gR) \coloneqq \bL \Dec_{E_k(1)}^{E_k}(\gR)$.

\subsubsection{The $E_\infty$-operad in the symmetric monoidal case} \label{sec:e-infty-operad} 
\newglossaryentry{cinfty}{%
	name={\ensuremath{\cC_\infty}},
	description={Non-unital $E_\infty$-operad},
	type=symbols
}
\newglossaryentry{cinftyplus}{%
	name={\ensuremath{\cC_\infty^+}},
	description={$E_\infty$-operad},
	type=symbols
}
\newglossaryentry{einfty}{%
	name={\ensuremath{E_\infty}},
	description={Monad associated to $\cC_\infty$},
	type=symbols
}
\newglossaryentry{einftyplus}{%
	name={\ensuremath{E_\infty^+}},
	description={Monad associated to $\cC_\infty^+$},
	type=symbols
}
The $E_\infty$-operad is obtained by letting $k$ go to $\infty$ in the definition of $\cC_k$. Sending a rectilinear embedding $e \colon \sqcup_i I^k \hookrightarrow I^k$ to $e \times \mr{id}_I \colon \sqcup_i I^{k+1} \hookrightarrow I^{k+1}$ defines a map $\cC_k \to \cC_{k+1}$ of operads in $\cat{sSets}$.\index{operad!$E_\infty^+$-}\index{operad!$E_\infty$-}

\begin{definition}We define the unitary $E_\infty$-operad and non-unitary $E_\infty$-operad as $\gls{cinftyplus} \coloneqq \underset{k \to \infty}{\colim}\, \cC^+_k$ and $\gls{cinfty} \coloneqq \underset{k \to \infty}{\colim}\, \cC_k$ respectively.
\end{definition}

As above, we get sifted monads $\gls{einftyplus}$ and $\gls{einfty}$ on any category $\sfC$ satisfying the axioms of Section \ref{sec:axioms-of-cats}, and indecomposables $Q^{E_\infty}$ with its derived functor. Remark \ref{rem:comparing-ek-operads} simplifies because we can compare $\cC_{\infty}$ and $\cE_\infty$ using the $E_\infty$-operad $\cC_{\infty} \times \cE_\infty$.

\subsubsection{Modification for $k$-monoidal categories with $k=1,2$}\label{sec:EkAlgForSmallk}

\newglossaryentry{c1fb1}{%
	name={\ensuremath{\cC_1^\cat{FB_1}}},
	description={Non-symmetric non-unital little $1$-cubes operad},
	type=symbols
}
\newglossaryentry{e1fb1}{%
	name={\ensuremath{E_1^{\cat{FB_1}}}},
	description={Monad associated to $\cC_1^\cat{FB_1}$},
	type=symbols
}
\newglossaryentry{c2fb2}{%
	name={\ensuremath{\cC_2^\cat{FB_2}}},
	description={Braided non-unital little $2$-cubes operad},
	type=symbols
}
\newglossaryentry{e2fb1}{%
	name={\ensuremath{E_2^{\cat{FB_2}}}},
	description={Monad associated to $\cC_2^\cat{FB2}$},
	type=symbols
}

\begin{definition}
The underlying $1$-symmetric sequence of the \emph{non-unitary non-symmetric little $1$-cubes operad} $\gls{c1fb1}$ is given by
\begin{align*} \cC^\cat{FB_1}_1(n) \coloneqq \begin{cases} \varnothing &  \text{if $n=0$,} \\
\mr{Emb}^\mr{rect,\cat{FB_1}}(\sqcup_n I, I) & \text{otherwise,} \end{cases}\end{align*}
where $\mr{Emb}^\mr{rect,\cat{FB_1}}(\sqcup_n I, I) \subset \mr{Emb}^\mr{rect}(\sqcup_n I, I)$ is the path component consisting of those $(e_1, \ldots, e_n)$ such that $e_1(0) < e_2(0) < \cdots < e_n(0)$. This is a monoid for the composition product with composition induced by composition of rectilinear embeddings.

For any monoidal category $\sfC$ satisfying our axioms (in particular, being enriched and copowered over simplicial sets) we then get a monad $E_1^{\cat{FB_1}}$ whose underlying functor is given by
\begin{equation*}
  X \longmapsto \bigsqcup_{n \geq 1} \cC^\cat{FB_1}_1(n) \times X^{\otimes n}.
\end{equation*}
\end{definition}
If the monoidal structure on $\sfC$ comes from a symmetric monoidal structure, this monad agrees with the previously defined monad $E_1$, up to natural isomorphism of monads.  Hence we shall often drop the superscript $\cat{FB}_1$ from the notation.

For each $n \geq 1$ there is an injective map of spaces
\begin{equation}\label{eq:30}
  \cC^\cat{FB_1}_1(n) \lra \cC_2(n)
\end{equation}
given by taking cartesian product with the identity map of $I$.  As $n \in \cat{FB}_1$ varies, these maps form a morphism of non-symmetric operads.  Since the domain is contractible for $n \geq 1$, this provides a ``basepoint subspace'' for $\cC_2(n)$.
\begin{definition}\label{def:e2-2-monoidal} The underlying $2$-symmetric sequence of the \emph{non-unitary braided little $2$-cubes operad} $\gls{c2fb2}$ is given by
\begin{align*} \cC^{\cat{FB_2}}_2(n)\coloneqq \begin{cases} \varnothing &  \text{if $n = 0$,} \\
\mr{Emb}^\mr{rect,\cat{FB_2}}(\sqcup_n I^2, I^2) & \text{otherwise,}\end{cases}\end{align*}
where $\mr{Emb}^\mr{rect,\cat{FB_2}}(\sqcup_n I^2, I^2)$ consists of pairs $(e,\eta)$ of an element $e \in \mr{Emb}^\mr{rect}(\sqcup_n I^2,I^2)$ and a homotopy class $\eta$ of path from $e$ to an element in the image of~(\ref{eq:30}).  This is a monoid for the composition product with composition induced by composition of rectilinear embeddings and concatenation of homotopy classes of path.

For any braided monoidal category $\sfC$ satisfying our axioms we then get a monad $\gls{e2fb1}$ whose underlying functor is given by
\begin{equation*}
  X \longmapsto \bigsqcup_{n \geq 1} \cC^\cat{FB_2}_2(n) \times_{\beta_n} X^{\otimes_\sfC n}.
\end{equation*}
\end{definition}

For each $n$, the action of $\beta_n$ on $\cC^\cat{FB_2}_2(n)$ is free, and the quotient space is homeomorphic to $\cC_2(n)/\fS_n$.  The quotient maps $\cC^\cat{FB_2}_2(n) \to \cC_2(n)$ are universal covering maps, and as $n \in \cat{FB_2}$ varies they define a morphism of braided operads.  If the braiding on $\sfC$ is a symmetry, so that the previously discussed monad $E_2$ is defined, it follows that the monads $E_2$ and $E_2^{\cat{FB}_2}$ agree up to natural isomorphism.  Hence we shall often drop the superscripts $\cat{FB_2}$ from the notation.

Discussion of unitalizations, monads, indecomposables and comparison of different $E_k$-algebras in Section \ref{sec:ekoperadssym} may easily be adapted to both of these cases, and we shall not bother the reader with this.

\subsection{Modules over $E_1$-algebras}\label{sec:modules}
\newglossaryentry{roverline}{%
	name={\ensuremath{\overline{R}}},
	description={Unital associative replacement of $E_1$-algebra $\gR$},
	type=symbols
}
\newglossaryentry{adapter}{%
	name={\ensuremath{A(\gR)}},
	description={Adapter for an $E_2$-algebra $\gR$},
	type=symbols
}

There are several homotopically equivalent descriptions of modules over an $E_1$-algebra, see Remark \ref{rem:OtherModuleModels}, but for our applications the following is most convenient. For each 
$E_1$-algebra $\gR$ we shall construct a unital associative algebra (i.e.\ a monoid) $\gls{roverline}$, reminiscent of the \emph{Moore loops construction}. We will then consider a (left or right) $\gR$-module to be a module over the unital associative algebra $\overline{\gR}$, and appeal to Section \ref{sec:AssocModules} for the homotopy theory of such. We will also show that $\gR$ canonically has the structure of a (left or right) $\overline{\gR}$-module, and that there is an equivalence $\overline{\gR} \simeq \bunit \sqcup \gR$ of (left or right) $\overline{\gR}$-modules (at least when $\gR$ is cofibrant in $\sfC$).

If $\sfC$ is pointed we will construct an augmentation $\epsilon \colon \overline{\gR} \to \bunit$, so as in Section \ref{sec:AssocModules} there are indecomposables functors
\[Q^{\overline{\gR}} \colon \overline{\gR}\text{-}\cat{Mod} \lra \sfC_* \quad\text{ and }\quad Q^{\overline{\gR}} \colon \cat{Mod}\text{-}\overline{\gR} \lra \sfC_*\]
(we omit leftness or rightness of the module from the notation), and hence derived functors $Q^{\overline{\gR}}_\bL$.

Finally, if the $E_1$-algebra $\gR$ is obtained by neglect of structure from an $E_2$-algebra then we will show that there is an \emph{adapter}: an object $\gls{adapter}$ having two homotopic left $\overline{\gR}$-module structures and a right $\overline{\gR}$-module structure, all three of which commute strictly, and such that $A(\gR) \simeq \overline{\gR}$ as an $\overline{\gR}$-$\overline{\gR}$-bimodule.  (See below for the precise definition of what that means.)  If $\gM$ is a left $\overline{\gR}$-module then the tensor product
\[A(\gR) \otimes_{\overline{\gR}}\gM\]
is weakly equivalent to $\gM$ as a left $\overline{\gR}$-module, but has an additional special left $\overline{\gR}$-module structure which strictly commutes with the first. 

Eventually, the special module structure can be used to cone off multiplication by a map $f \colon \partial D^{g,d+1} \to \gR$ in the $\gR$-module structure on $\gM$ and have the result $\gM/f$ again be a left $\gR$-module.
This construction will play a central role in our applications.

\subsubsection{The unital associative replacement of an $E_1$-algebra}

For an $E_1$-algebra $\gR$, we explain how to obtain a unital associative algebra $\overline{\gR}$. This construction was made in \cite[Section 2]{MandellSmash} and \cite[Section 3]{BM} for specific categories of spectra, and there is no essential difficulty in adapting it to our setting (though we also change the construction slightly).

Firstly, we define the object in $\sfC$ underlying the unital associative algebra $\overline{\gR}$ as
the coproduct
\begin{equation}\label{eq:TildeUnderlying}
\overline{\gR} \coloneqq (\{0\} \times \bunit) \sqcup ((0,\infty) \times \gR),
\end{equation}
where $\{0\}$ and $(0,\infty)$ denote (singular simplicial sets of) these subset of $\bR$ in the Euclidean topology.  To define an associative multiplication, it suffices to give $\{0\} \times \bunit$ the structure of a unital associative algebra, $(0,\infty) \times \gR$ the structure of a $\{0\} \times \bunit$-bimodule, and $(0,\infty) \times \gR$ the structure of a non-unital associative algebra. 

We have
\[(\{0\} \times \bunit) \otimes (\{0\} \times \bunit) \cong (\{0\} \times \{0\}) \times (\bunit \otimes \bunit)\]
and addition on the first factor and the unit structure on the second make this into an associative algebra. Furthermore, we have
\begin{align*}
(\{0\} \times \bunit) \otimes ((0,\infty) \times \gR) &\cong (\{0\} \times (0,\infty)) \times (\bunit \otimes \gR)\\
 ((0,\infty) \times \gR) \otimes (\{0\} \times \bunit) &\cong ((0,\infty) \times \{0\}) \times (\gR \otimes \bunit)
\end{align*}
and addition on the first factor and the unit structure on the second define a $\{0\} \times \bunit$-bimodule structure. Finally, we have
\[((0,\infty) \times \gR) \otimes ((0,\infty) \times \gR) \cong ((0,\infty) \times (0,\infty)) \times (\gR \otimes \gR)\]
whereupon we use the map
\begin{align*}
\Gamma \colon (0,\infty) \times (0,\infty) &\lra (0,\infty) \times \cC_1^{\cat{FB}_1}(2)\\
(s, t) &\mapsto \left(s+t, \left(x \mapsto \frac{s \cdot x}{s+t}, x \mapsto \frac{s+t \cdot x}{s+t} \right)\right)
\end{align*}
and the action map $\cC_1^{\cat{FB}_1}(2) \times (\gR \otimes \gR) \to \gR$ to get to $(0,\infty) \times \gR$. The associativity of addition makes $(0,\infty) \times \gR$ into a non-unital associative algebra, and in total we have produced a unital associative algebra $\overline{\gR}$. Letting $\mr{Ass}^+$ denote the operad for unital associative algebras, this defines a functor
\[\overline{(-)} \colon \Alg_{E_1}(\sfC) \lra \Alg_{\mr{Ass}^+}(\sfC),\]
which has the following properties:

\begin{lemma}\label{lem:tilde-properties} We have that
	\begin{enumerate}[(i)]
		\item $\overline{\gR}$ is cofibrant in $\sfC$ if and only if $\gR$ is cofibrant in $\sfC$,
		\item $\overline{(-)}$ preserves weak equivalences between objects which are cofibrant in $\sfC$,
		\item $\overline{(-)}$ commutes with the functor $\grr \colon \Alg_{E_1}(\sfC^{\bZ_{\leq}}) \to \Alg_{E_1}(\sfC_*^{\bZ_{=}})$ and the functor $\colim \colon \Alg_{E_1}(\sfC^{\bZ_{\leq}}) \to \Alg_{E_1}(\sfC)$.
	\end{enumerate}
\end{lemma}

\begin{proof}For (i), recall that we are working under the assumption that $\bunit$ is cofibrant, cf.~Axiom \ref{axiom.model-monoidalsimplicial}, so this follows from (\ref{eq:TildeUnderlying}). Property (ii) similarly follows from that formula (note that coproducts in general only preserve weak equivalences between cofibrant objects, hence the cofibrancy assumption). Finally, (iii) follows from the fact that both $\grr$ and $\colim$ are simplicial functors. 
\end{proof}

We first show that when applied to an $E_1$-algebra that is already associative, the result is weakly equivalent (as a unital associative algebra) to the unitalization.

\begin{lemma}\label{lem:OverlineOnAssoc}
If $\gR$ is a non-unital associative algebra, considered as an $E_1$-algebra, then the morphism
\[\overline{\gR} \cong (\{0\} \times \bunit) \sqcup ((0,\infty) \times \gR) \lra \bunit \sqcup \gR = \gR^+\]
induced by the projections $\{0\} \to *$ and $(0,\infty) \to *$ is a natural weak equivalence of unital associative algebras. 
\end{lemma}

\begin{proof}
That this morphism is a homotopy equivalence is clear, and so it is a weak equivalence. That it is a morphism of unital associative algebras follows by unravelling the definition of the multiplication on $\overline{\gR}$.
\end{proof}

We will now show that for a $E_1$-algebra $\gR$, $\overline{\gR}$ is naturally weakly equivalent to the unitalization $\gR^+$ of $\gR$ as a $E_1^+$-algebra. To do so, we introduce the notation $\mr{Ass}$ for the non-unital associative operad,\index{operad!$\mr{Ass}$-} and $\mr{Ass}^+$\index{operad!$\mr{Ass}^+$-} for the unital associative operad.

\begin{proposition}\label{prop:OverlineEqPlus}
There is a zig-zag of natural transformations
\[\overline{(-)} \Longleftarrow \cdots \Longrightarrow (-)^+ \colon \Alg_{E_1}(\sfC) \lra \Alg_{E_1^+}(\sfC),\]
which are weak equivalences on those objects which are cofibrant in $\sfC$.
\end{proposition}
\begin{proof}
Let us first show that there is such a zig-zag of functors between the identity on $\Alg_{E_1}(\sfC)$ and a functor giving a non-unital associative algebra. To do this, note that there is map of operads $\cC_1 \to \mr{Ass}$ to the non-unital associative operad, by which every associative unital algebra becomes an $E_1$-algebra. There is therefore a zig-zag
\begin{equation}\label{eq:ZZ}
\gR \longleftarrow B_\bullet(F^{E_1}, E_1, \gR) \lra B_\bullet(F^\mr{Ass}, E_1, \gR),
\end{equation}
where the left map is an augmentation and the right one is a semi-simplicial map, and the thick geometric realization $B(F^{\mr{Ass}}, E_1, \gR)$ is a non-unital associative algebra.

Using Lemma \ref{lem:OverlineOnAssoc} we may therefore form the zig-zag
\begin{equation}\label{eq:LongZigZag}
\begin{tikzcd}
\gR^+ & \lar B(F^{E_1}, E_1, \gR)^+ \rar & B(F^{\mr{Ass}}, E_1, \gR)^+\\
{\overline{\gR}} & \lar  \overline{B(F^{E_1}, E_1, \gR)} \rar & {\overline{B(F^{\mr{Ass}}, E_1, \gR)}.} \uar
\end{tikzcd}
\end{equation}
Now suppose that $\gR$ is cofibrant in $\sfC$. The left map in \eqref{eq:ZZ} is a weak equivalence on geometric realization by Lemma \ref{lem:extra-degeneracy}, as (after neglecting the $E_1$ structure) the augmented semi-simplicial object has an extra degeneracy. Using Lemma \ref{lem.symmetric-sequence-preserving} and the assumption that $\gR$ is cofibrant in $\sfC$, the two semi-simplicial objects in \eqref{eq:ZZ} are Reedy cofibrant. Since the maps $\cC_1(n) \to \mr{Ass}(n)$ are weak equivalences of free $\fS_n$-spaces, the map of simplicial objects is a levelwise weak equivalence, so by Lemma \ref{lem:thick-geom-rel-cofibrations} its geometric realization is a weak equivalence between objects which are cofibrant in $\sfC$. Because $\bunit$ was assumed cofibrant, $(-)^+ = \bunit \sqcup -$ preserves weak equivalences between objects which are cofibrant in $\sfC$, so the maps in the top row of \eqref{eq:LongZigZag} are weak equivalences. Similarly, by Lemma \ref{lem:tilde-properties} (ii) the maps in the bottom row of \eqref{eq:LongZigZag} are weak equivalences. 
The vertical map in \eqref{eq:LongZigZag} is a weak equivalence by Lemma \ref{lem:OverlineOnAssoc}.
\end{proof}

\begin{remark}\label{rem:OtherModuleModels} There are two other approaches to producing a unital associative monoid out of a non-unital $E_1$-algebra, weakly equivalent under suitable cofibrancy conditions. 
	
Firstly, as in the proof of Proposition \ref{prop:OverlineEqPlus}, we may apply the bar construction to the morphism of operads $\cC_1 \to \mr{Ass}^+$ and define a functor
\[\gR \longmapsto B(F^{\mr{Ass}^+},E_1,\gR),\]
which is isomorphic to $B(F^{\mr{Ass}},E_1,\gR)^+$. By the proof of Proposition \ref{prop:OverlineEqPlus} this is weakly equivalent to $\overline{\gR}$ as long as $\gR$ is cofibrant in $\sfC$.

Secondly, we may use infinitesimal modules over an $\cO$-algebra $\gR$. An infinitesimal module is an object $M$ with maps $\cO(n) \times_{G_{n-1}} (R^{\otimes n-1} \otimes M) \to M$ satisfying suitable associativity and unit axioms. Equivalently, there exists a unital associative monoid $\mr{Inf}^\cO(\gR)$ such that an infinitesimal $\cO$-algebra is exactly a left $\mr{Inf}^\cO(\gR)$-module. In the case of the $E_1$-operad in a 1-monoidal category and $\gR = \gE_1(X)$ a free $E_1$-algebra, it is given by
\[\mr{Inf}^{E_1}(\gE_1(X)) \coloneqq \bigsqcup_{n \geq 1} \cC_1^{\cat{FB}_1}(n) \times X^{\otimes n-1},\]
which is weakly equivalent to $\gE_1^+(X) \simeq \overline{\gE_1(X)}$, because $\cC_1^{\cat{FB}_1}(n)$ is contractible and so in particular is weakly equivalent to $\cC_1^{\cat{FB}_1}(n-1)$. By taking a free simplicial resolution of $\gR$ one proves as $\gR$ is cofibrant in $\sfC$.
\end{remark}

\subsubsection{Adapters,  and bi- and tri-modules over $\overline{\gR}$}\label{sec:adapter}

\setcounter{secnumdepth}{5}

\newcommand{\ogR}{\overline{\gR}}

When the ambient category $\sfC$ is at least braided monoidal, the
tensor product $\ogR \otimes \ogR \in \sfC$ inherits the structure of
an associative algebra, with multiplication
\begin{equation*}
  (\ogR \otimes \ogR) \otimes (\ogR \otimes \ogR) \xrightarrow{\ogR
    \otimes \beta_{\ogR, \ogR} \otimes \ogR}   (\ogR \otimes \ogR) \otimes (\ogR
  \otimes \ogR) \xrightarrow{\mu \otimes \mu} \ogR \otimes \ogR.
\end{equation*}
There are maps of associative algebras
\begin{align}
  \ogR & \cong \ogR \otimes \mathbbm{1} \lra
         \ogR \otimes \ogR\label{eq:3}\\
  \ogR & \cong \mathbbm{1} \otimes \ogR \lra
          \ogR \otimes \ogR\label{eq:17},
\end{align}
obtained from the unit isomorphism in $\sfC$ and the unit of $\ogR$, giving rise to two (\emph{a priori} distinct) ways of viewing an $(\ogR \otimes \ogR)$-module as an $\ogR$-module.

If $\gR$ is an $E_k$-algebra and $k \geq 2$, the associative algebra
$\ogR$ depends only on the $E_1$-algebra obtained from $\gR$ by neglect
of structure. The forgotten structure may of course be used to
perform further constructions.  We shall discuss a method
for turning a left $\ogR$-module $\gM$ into a left $(\ogR\otimes\ogR)$-module,
provided $k \geq 2$, based on the following notion.
\begin{definition}
  Let $\gR$ be an $E_1$ algebra which is cofibrant in $\sfC$ and let
  $\ogR$ be the corresponding associative algebra.  An
  \emph{adapter}\index{adapter} for $\gR$ is a cofibrant
  $(\ogR \otimes \ogR)$-$\ogR$-bimodule $A(\gR)$, together with a zig-zag
  of weak equivalences
  \begin{equation}\label{eq:26}
    A(\gR) \overset{\simeq}\longleftarrow \cdots \overset{\simeq}\longrightarrow \ogR
  \end{equation}
  of $\ogR$-$\ogR$-bimodules when $A(\gR)$ is viewed as an
  $\ogR$-$\ogR$-bimodule via~(\ref{eq:3}), satisfying moreover that
  the diagram
  \begin{equation*}
    \begin{tikzcd}
      {(\ogR \otimes \mathbbm{1}) \otimes A(\gR)} \arrow{r}{} \arrow{d}{\cong}[swap]{\beta_{\ogR, \bunit}\otimes A(\gR)} &{A(\gR)}\\
      {(\mathbbm{1}\otimes \ogR) \otimes A(\gR)} \arrow{ur}{}
    \end{tikzcd}
  \end{equation*}
  in the category of right $\ogR$-modules, arising from the two left
  $\ogR$-module structures~(\ref{eq:3}) and~(\ref{eq:17}), is homotopy commutative (i.e.\ becomes commutative in the homotopy category of right
  $\ogR$-modules).
\end{definition}

\begin{remark}
Since $\gR$ is cofibrant in $\sfC$ so is $\ogR$ (by Lemma \ref{lem:tilde-properties} (i)) and hence $\smash{\ogR^{\otimes 3}}$, so the cofibrancy of $A(\gR)$ as a $(\ogR\otimes \ogR)$-$\ogR$-bimodule implies that it is cofibrant in $\sfC$ (the forgetful functor preserves cofibrant objects just like in Theorem \ref{thm:model-structure-on-modules}).  By a similar argument it is cofibrant as a right $\ogR$-module.

Conversely, if one has an $(\ogR\otimes \ogR)$-$\ogR$-bimodule $A'(\gR)$ satisfying all conditions above except for cofibrancy, then any cofibrant approximation $A(\gR) \to A'(\gR)$ as an $(\ogR\otimes \ogR)$-$\ogR$-bimodule will be an adapter.
\end{remark}

The defining properties of adapters depend only on the associative
algebra $\ogR$ and hence only on the $E_1$-algebra structure on $\gR$
(as well as the braided monoidal structure on the ambient category
$\sfC$).  But they most naturally arise when $\gR$ is obtained from an
$E_2$-algebra by neglect of structure.  In
Section~\ref{sec:constr-adapt} below, we shall construct adapters
$\gR \mapsto A(\gR)$, functorially in the $E_2$-algebra $\gR$,
together with a zig-zag of natural transformations~(\ref{eq:26}) which
are weak equivalences when $\gR$ is cofibrant in $\sfC$.
Functoriality in $\gR$ means that a map $c \colon \gR \to \gR'$ of
$E_2$-algebras induces a map $A(c) \colon A(\gR) \to A(\gR')$ of
$(\ogR \otimes \ogR)$-$\ogR$-bimodules.  There is a left adjoint map
of $(\ogR' \otimes \ogR)$-$\ogR$-bimodules
\begin{equation}\label{eq:27}
  (\ogR' \otimes \ogR) \otimes_{\ogR \otimes \ogR} A(\gR) \lra A(\gR'),
\end{equation}
and it follows from~(\ref{eq:26}) that this map is a weak equivalence when both $\gR$ and $\gR'$ are cofibrant in $\sfC$.

\subsubsection{Applications of adapters}\label{sec:ApplyingAdapters}

If $\gR$ is an $E_1$-algebra, $\gM$ is a left $\ogR$-module, and $f \colon
\partial D^{g,d+1} \to \ogR$ is a map, then we may form the ``left multiplication
by $f$'' map
\begin{equation*}
  f \cdot - \colon (\partial D^{g,d+1}) \otimes \gM \xrightarrow{f \otimes \gM} \ogR \otimes \gM
  \overset{\mu}\lra \gM.
\end{equation*}
As a preliminary definition, to be updated below, we define
an object $\gM/f$ by taking the
homotopy pushout diagram
\begin{equation}\label{eq:29}
  \begin{tikzcd}
    {(\partial D^{g,d+1}) \otimes \gM} \arrow{d}{} \arrow{r}{f \cdot -} & {\gM} \arrow{d}{}\\
    {D^{g,d+1} \otimes \gM} \arrow{r}{} & {\gM/f},
  \end{tikzcd}
\end{equation}
where $\partial D^{g,d+1}$ and $D^{g,d+1}$ are as defined in
Section~\ref{sec:cell-def}.  The objects
$\partial D^{g,d+1} \otimes \gM$, $\gM$, and $D^{g,d+1} \otimes \gM$
are all left $\ogR$-modules, but unfortunately $f \cdot -$ is not a
module map, so the homotopy pushout may only be formed in $\sfC$.

\begin{definition}\label{defn:R-mod-f}
  Let $A(\gR)$ be an adapter for $\gR$.  Given a map $f\colon \partial D^{g,d+1} \to \ogR$, define a map $\phi(f)\colon \partial D^{g,d+1} \otimes A(\gR) \to A(\gR)$ of $\ogR$-$\ogR$-bimodules as the composition
  \begin{equation*}
    (\partial D^{g,d+1}) \otimes A(\gR) \xrightarrow{f \otimes A(\gR)} \ogR \otimes A(\gR) \overset{\eqref{eq:17}}\lra (\ogR \otimes \ogR) \otimes A(\gR)
    \overset{\mu}\lra A(\gR),
  \end{equation*}
  and define $\ogR/f$ as the pushout
  \begin{equation}\label{eq:33}
    \begin{tikzcd}
      {(\partial D^{g,d+1}) \otimes A(\gR)} \arrow{d}{} \arrow{r}{\phi(f)} & {A(\gR)} \arrow{d}{}\\
      {D^{g,d+1} \otimes A(\gR)} \arrow{r}{} & {\ogR/f}
    \end{tikzcd}
  \end{equation}
  in the category of $\ogR$-$\ogR$-bimodules.
\end{definition}

Since $A(\gR)$ is assumed cofibrant, diagram~(\ref{eq:33}) is also homotopy pushout.  After forgetting down to a diagram in right $\ogR$-modules the map $\phi(f)$ is homotopic to left multiplication by $f$ in the bimodule structure, so the underlying map of right $\ogR$-modules $\ogR \simeq A(\gR) \to \ogR/f$
is equivalent to a cell attachment along $f$ in this category.

\begin{definition}
  Let $R$ and $f$ be as in definition~\ref{defn:R-mod-f}, and let $\ogR/f$ be the $\ogR$-$\ogR$-module defined there. For a left $\ogR$-module $\gM$ define an $\ogR$-module $\gM/f$ as
  \begin{equation*}
    \gM/f = (\ogR/f) \otimes_{\ogR} \gM.
  \end{equation*}
\end{definition}

The defining properties of adapters, and cofibrancy of $\ogR/f$ as an $\ogR$-module,
imply a homotopy pushout diagram
\begin{equation*}
  \begin{tikzcd}
    {\partial D^{g,d+1} \otimes \gM} \arrow{d}{} \arrow{r}{f \cdot -} &
    {\gM} \arrow{d}{}\\
    {D^{g,d+1} \otimes \gM} \arrow{r}{} & {U^{\ogR} (\gM/f)}, 
  \end{tikzcd}
\end{equation*}
in $\sfC$.  In case $\sfC$ is pointed, we have an object $S^{g,d} \in \sfC$ as in Section~\ref{sec:cell-def}, and a map $c \colon \partial D^{g,d+1} \to S^{g,d}$ obtained by identifying a point in $\partial D^{g,d+1}$ with the basepoint. If $f$ factors as 
\[f \colon \partial D^{g,d+1} \overset{c} \lra S^{g,d} \overset{f'}\lra \ogR\]
then there is a cofiber sequence
\begin{equation*}
  S^{g,d} \wedge \gM \xrightarrow{f' \wedge -} \gM \lra U^{\ogR} (\gM/f),
\end{equation*}
expressing informally that $U^{\ogR}(\gM/f)$ is the ``cofiber of left multiplication by $f'$''.

\begin{remark}
Upon working in a pointed setting, we can make sense of homotopy groups $\pi_{g,d}(\ogR) = \pi_{d}(\ogR(g))$ as in Section \ref{sec:additive-case}.
  When $\Ob(\sfG) = \bN$, $\bigoplus_{g,d} \pi_{g,d}(\ogR)$ is a bigraded ring and $\bigoplus_{g,d} \pi_{g,d}(\gM)$ is a bigraded module over it. If further $\sfS$ is semistable so that $S^{g,d} \wedge -$ only has the effect of shifting bigrading, then the homotopy groups of $\gM/f$ fit into a long exact sequence with the homotopy groups of $\gM$ and multiplication by $[f] \in \pi_{g,d}(\ogR)$ in the bigraded module structure.
\end{remark}

Since the resulting object $\gM/f$ is again an $\ogR$-module, this
operation may be iterated: given pointed maps $f_i' \colon S^{g_i,d_i} \to \ogR$ in $\sfC$, we may form a left
$\ogR$-module
\begin{equation*}
  \gM/(f_1,\dots, f_n) = (\dots((\gM/f_1)/f_2)/ \dots)/f_n.
\end{equation*}

Let us point out that this process of course depends on the choice of adapter, although we have omitted that from the notation.  We shall not discuss the extent to which it is independent up to homotopy, except point out that at least the homotopy type of the object $U^{\ogR} \gM/f$ of $\sfC$ is independent of this choice.

\subsubsection{Base change for the $\gM/f$ construction}

Let us also spell out a naturality property of this construction with
respect to maps $h\colon \gR \to \gR'$ between cofibrant $E_1$-algebras.  Given
$f\colon S^{g,d} \to \gR$ we may apply this construction to $f$ or to
$h \circ f\colon S^{g,d} \to \gR'$.  The weak equivalence~(\ref{eq:27})
then implies a weak equivalence 
\begin{equation}\label{eqn:adapter-quotient}
  \ogR' \otimes_{\ogR} (\gM/f) \simeq (\ogR' \otimes_{\ogR}
  \gM)/(h \circ f)
\end{equation}
in the category of left $\ogR'$-modules, and hence an isomorphism in the
homotopy category.

\subsubsection{Construction of adapters}
\label{sec:constr-adapt}

As promised, we shall now explain a construction $\gR \mapsto A(\gR)$ of adapters, functorial in $E_2$-algebras $\gR$. There are surely many ways to construct such an object and establish its basic properties; we take a hands-on approach.  It does not seem clear that the construction below results in a cofibrant $(\ogR \otimes \ogR)$-$\ogR$-bimodule, so $A(\gR)$ should be some functorial cofibrant approximation to the object resulting from the geometric construction.

\begin{figure}[h]
	\centering
\begin{tikzpicture}
\begin{scope}[scale=0.5]
\draw[fill = Mahogany!10!white] (1,0) -- (1,2.5) -- (2.5,2.5) -- (2.5,5) -- (2, 5) -- (2, 7.5)  -- (7.5,7.5) -- (7.5,2.5) -- (8,2.5) -- (8,0) -- (2,0);

\draw (0,0) rectangle (10,10);

\draw [dashed] (2,7.5) -- (2,10);

\node at (1,-0.5) {$-1-a$};
\node at (2,10.5) {$-1-c$};
\node at (8,-0.5) {$1+b$};

\node at (5,-0.5) {$0$};

\node at (-0.5,0) {$-1$};
\node at (-0.5,2.5) {$0$};
\node at (-0.5,5) {$1$};
\node at (-0.5,7.5) {$2$};

\node at (5,2.5) {$V_{a,b}^c$};
\end{scope}

	\begin{scope}[scale=0.5,xshift=13cm]
\draw  (1,0) -- (1,2.5) -- (2.5,2.5) -- (2.5,5) -- (2, 5) -- (2, 7.5)  -- (7.5,7.5) -- (7.5,2.5) -- (8,2.5) -- (8,0) -- (2,0);

\draw (0,0) rectangle (10,10);

\node at (5,-0.5) {$0$};

\node at (-0.5,0) {$-1$};
\node at (-0.5,2.5) {$0$};
\node at (-0.5,5) {$1$};
\node at (-0.5,7.5) {$2$};

\draw[fill = Mahogany!10!white] (3,2) rectangle (4,4);
\draw[fill = Mahogany!10!white] (5,1) rectangle (7,3);
\draw[fill = Mahogany!10!white] (2.5,0.2) rectangle (4.5,1.8);

\node at (6,2) {$e_1$};
\node at (3.5,3) {$e_2$};
\node at (3.5,1) {$e_3$};
\end{scope}
\end{tikzpicture}
	\caption{An example of the region $V_{a,b}^c$ and an element of $V(3)$.}
	\label{fig:Uabc}
\end{figure}

\begin{definition}\label{def:vabc}
For a triple $a, b, c \in [0,\infty)$, define a subset $V_{a,b}^c \subset \bR^2$ by
\begin{align*}
V_{a,b}^c &\coloneqq ([-1-a, 1+b] \times [-1,0])\\
 & \,\,\,\,\,\,\,\,\,\, \cup ([-1, 1] \times [-1,2])\\
 & \,\,\,\,\,\,\,\,\,\,  \cup ([-1-c,1] \times [1,2]),
\end{align*}
as depicted in Figure \ref{fig:Uabc}. We obtain a symmetric sequence of topological spaces by
\[V(n) \coloneqq \{(a,b,c; e) \in [0,\infty)^3 \times \Emb^\mr{rect}(\sqcup_n I^2, \bR^2) \, \vert \, \mr{im}(e) \subset V_{a,b}^c\},\]
where the symmetric group $\fS_n$ acts as usual on $\Emb^\mr{rect}(\sqcup_n I^2, \bR^2)$.

The associated functor is $X \mapsto V(X) = \bigsqcup_{n \geq 0} V(n) \times_{\fS_n} X^{\otimes n}$, and it inherits a right action of the monad $E_2$.
\end{definition}

\begin{remark}
  If the braiding on $\sfC$ is not a symmetry then the above definition should be rephrased, since the action of $\beta_n$ on $\gR^{\otimes n}$ does not factor through $\beta_n \to \fS_n$.  The following modification works for braided monoidal $\sfC$ and gives a result isomorphic to the above in case the braiding is a symmetry.  Replace the spaces $V(n)$ by their universal covering spaces $\tilde V(n)$, taken with respect to the contractible ``basepoint subspace'' defined by the condition on rectilinear embedding $e$ that the composition
  \begin{equation*}
    \{1, \dots, n\} \hookrightarrow \{1, \dots, n\} \times I^2 \overset{e}\lra \bR^2 \overset{\pi_1}\lra \bR
  \end{equation*}
  is an order-preserving injection.  This universal cover comes with an action of the braid group $\beta_n$ and the covering map $\tilde V(n) \to V(n)$ is equivariant for the usual homomorphism $\beta_n \to \fS_n$, and the resulting functor
  \begin{equation*}
    \tilde V\colon X \longmapsto {V}(X)=\bigsqcup_{n \geq 0} \tilde V(n) \times_{\beta_n} X^{\otimes n}
  \end{equation*}
  inherits a right action of the monad associated to the braided operad $\cC_2^{\cat{FB}_2}$.

If the braiding on $\sfC$ is a symmetry, this agrees (up to canonical isomorphism) with the functor $X \mapsto V(X)$ from Definition~\ref{def:vabc}.
\end{remark}

\begin{definition}
The functor $X \mapsto V(X)$ has the structure of a right $E_2$-functor and by Corollary \ref{cor:monad-from-operad-sifted} it preserves sifted colimits. Thus we may define $A \colon \Alg_{E_2}(\sfC) \to \sfC$ on free algebras by $A(\gE_2(X)) \coloneqq V(X)$ and extend this to general $E_2$-algebras by density under sifted colimits using Proposition \ref{prop:extend-by-density}. Explicitly, for an $E_2$-algebra $\gR$, $A(\gR)$ is the coequalizer of the two maps
\[\begin{tikzcd}V(E_2(\gR)) \arrow[shift left=.5ex]{r} \arrow[shift left=-.5ex]{r} & V(\gR) \rar & A(\gR)\end{tikzcd}\]
given by the right $E_2$-functor structure of $V$ and the $E_2$-algebra structure of $\gR$.

There is a sub-symmetric sequence $V^0 \subset V$ (and in the braided case, $\tilde V^0 \subset \tilde V$) consisting of tuples of the form $(a, b, 0 ; e)$ where $\mr{im}(e) \subset [-1-a,1+b] \times [-1,0]$. The associated functor $X \mapsto V^0(X)$ is again a right $E_2$-functor and preserves sifted colimits, and one defines $A^0(\gR)$ as above; by construction it has a natural map to $A(\gR)$. 

The empty collection of little cubes determines a map $\iota \colon \bunit \to A^0(\gR) \to A(\gR)$.
\end{definition}

\begin{figure}
	\centering
	\begin{tikzpicture}
		\begin{scope}[scale=0.5]
			\draw[fill = Mahogany!10!white] (1,0) -- (1,2.5) -- (2.5,2.5) -- (2.5,5) -- (2.5, 5) -- (2.5, 7.5)  -- (7.5,7.5) -- (7.5,2.5) -- (8,2.5) -- (8,0) -- (2,0);
			
			\draw (0,0) rectangle (10,10);
		
			\draw [dashed] (0,2.5) -- (10,2.5);
			
			\node at (1,-0.5) {$-1-a$};
			\node at (8,-0.5) {$1+b$};
			
			\node at (5,-0.5) {$0$};
			
			\node at (-0.5,0) {$-1$};
			\node at (-0.5,2.5) {$0$};
			\node at (-0.5,5) {$1$};
			\node at (-0.5,7.5) {$2$};
			
			\node at (5,5) {$V_{a,b}^0$};
		\end{scope}
		
		\begin{scope}[scale=0.5,xshift=13cm]
			\draw  (1,0) -- (1,2.5) -- (2.5,2.5) -- (2.5,5) -- (2.5, 5) -- (2.5, 7.5)  -- (7.5,7.5) -- (7.5,2.5) -- (8,2.5) -- (8,0) -- (2,0);
			
			\draw (0,0) rectangle (10,10);
			
			\draw [dashed] (0,2.5) -- (10,2.5);
						
			\node at (5,-0.5) {$0$};
			
			\node at (-0.5,0) {$-1$};
			\node at (-0.5,2.5) {$0$};
			\node at (-0.5,5) {$1$};
			\node at (-0.5,7.5) {$2$};
			
			\draw[fill = Mahogany!10!white] (5,.8) rectangle (7,1.8);
			\draw[fill = Mahogany!10!white] (2.5,0.2) rectangle (4.5,1.8);
			
			\node at (6,1.3) {$e_1$};
			\node at (3.5,1) {$e_2$};
		\end{scope}
	\end{tikzpicture}
	\caption{An example of the region $V_{a,b}^0$ and an element of $V^0(2)$. The requirement for lying in $V^0$ is that all cubes lies below the dashed line in $V_{a,b}^0$.}
	\label{fig:v0}
\end{figure}

We now wish to explain how $A(\gR)$ has an $(\ogR \otimes \ogR)$-$\ogR$-bimodule structure.  Equivalently, this is encoded by three morphisms $\mu_\mathrm{ll}\colon \ogR \otimes A(\gR) \to A(\gR)$, $\mu_\mathrm{ul}\colon \ogR \otimes A(\gR) \to A(\gR)$, and $\mu_\mathrm{r}\colon A(\gR) \otimes \ogR \to A(\gR)$ in $\sfC$ which ``commute'' in the sense that the six orders in which multiplication can be performed result in equal morphisms
\begin{equation*}
  \ogR \otimes \ogR \otimes A(\gR) \otimes \ogR \lra A(\gR)
\end{equation*}
in $\sfC$.  Two of these involve the braiding in $\sfC$, namely the compositions
\begin{equation*}
  \begin{tikzcd}[row sep=large,column sep=large]
    \ogR \otimes \ogR \otimes A(\gR) \otimes \ogR \arrow{rr}{\beta_{\ogR, \ogR}\otimes A(\gR) \otimes \ogR } & &  \ogR \otimes \ogR \otimes A(\gR) \otimes \ogR \arrow{d}{\ogR \otimes \mu_\mathrm{ll} \otimes \ogR} \\
    A(\gR) & & \ogR \otimes A(\gR) \otimes \ogR. \arrow[bend left=10]{ll}{\mu_\mathrm{ul} \circ (\ogR \otimes \mu_\mathrm{r})} \arrow[bend right=10]{ll}[swap]{\mu_\mathrm{r} \circ (\mu_\mathrm{ul} \otimes \ogR)} 
\end{tikzcd}
\end{equation*}
The subscripts ``$\mathrm{ll}$'' and ``$\mathrm{ul}$'' stand for ``lower left'' and
``upper left''; some readers may find it helpful to imagine the left and right tensor factors in $(\ogR \otimes \ogR)$ typeset instead as ``lower'' and ``upper'' tensor factors, respectively.

For the right $\overline{\gR}$-module structure, we start with the map
\[V(\gR) \otimes \overline{\gR} \cong \bigsqcup_{n \geq 0} (\tilde V(n) \times_{\beta_n} \gR^{\otimes n}) \otimes ((\{0\} \times \bunit) \sqcup ((0,\infty) \times \gR)) \lra V(\gR)\]
given heuristically by the formula
\[((a,b,c; e), (r_1, \ldots, r_n), (t, r)) \longmapsto ((a, b+t, c ; e'), (r_1, \ldots, r_n, r))\]
where $e' \colon \sqcup_{n+1} I^2 \to V_{a, b+t}^c$ is the rectilinear embedding given by $e$ on the first $n$ copies of $I^2$, and by
\[(x,y) \longmapsto (1+b + x \cdot t, y-1)\]
on the final copy. It is easy to see that this map descends to a map $A(\gR) \otimes \overline{\gR} \to A(\gR)$, and defines a right $\overline{\gR}$-module structure.

The ``lower left'' module structure $\mu_\mathrm{ll}\colon \ogR \otimes A(\gR) \to A(\gR)$ is defined similarly. We start with the map
\[\overline{\gR} \otimes V(\gR) \cong \bigsqcup_{n \geq 0} ((\{0\} \times \bunit) \sqcup ((0,\infty) \times \gR)) \otimes (\tilde V(n) \times_{\beta_n} \gR^{\otimes n})  \lra V(\gR)\]
given heuristically by the formula
\[((t, r), (a,b,c; e), (r_1, \ldots, r_n)) \mapsto ((a+t, b, c ; e'), (r_1, \ldots, r_n, r))\]
where $e' \colon \sqcup_{n+1} I^2 \to V_{a+t, b}^c$ is the rectilinear embedding given by $e$ on the first $n$ copies of $I^2$, and by
\[(x,y) \longmapsto (-1-a-t + x \cdot t, y-1)\]
on the final copy. As above, this descends to a map $\overline{\gR} \otimes A(\gR)  \to A(\gR)$, and defines a left $\overline{\gR}$-module structure. These two module structures clearly restrict to module structures on $A^0(\gR)$.

The ``upper left'' module structure $\mu_\mathrm{ul}\colon \ogR \otimes A(\gR) \to A(\gR)$ is also similar, but let us spell out the details. We start with the map
\[\overline{\gR} \otimes V(\gR) \cong \bigsqcup_{n \geq 0} ((\{0\}\times \bunit) \sqcup ((0,\infty) \times \gR)) \otimes (\tilde V(n) \times_{\beta_n} \gR^{\otimes n})  \lra V(\gR)\]
given heuristically by the formula
\[((t, r), (a,b,c; e), (r_1, \ldots, r_n)) \longmapsto ((a, b, c+t ; e'), (r_1, \ldots, r_n, r))\]
where $e' \colon \sqcup_{n+1} I^2 \to V_{a, b}^{c+t}$ is the rectilinear embedding given by $e$ on the first $n$ copies of $I^2$, and by
\[(x,y) \longmapsto (-1-c-t + t\cdot x, 1+y)\]
on the final copy. As above, this descends to a map $\overline{\gR} \otimes A(\gR)  \to A(\gR)$, and defines a left $\overline{\gR}$-module structure. 

It is clear that these three module structures commute with each other in the sense described; Figure
\ref{fig:example-special-left} may clarify these module structures and this fact.

\begin{figure}[h]
	\centering
	\begin{tikzpicture}
		\begin{scope}[scale=0.5]
		\draw (1,0) -- (1,2.5) -- (2.5,2.5) -- (2.5,5) -- (2, 5) -- (2, 7.5)  -- (7.5,7.5) -- (7.5,2.5) -- (8,2.5) -- (8,0) -- (2,0);
		
		\draw (0,0) rectangle (10,10);
		
		\draw [dashed] (2,7.5) -- (2,10);

		\node at (1,-0.5) {$-1-a$};
		\node at (2,10.5) {$-1-c$};
		\node at (8,-0.5) {$1+b$};

		\node at (5,-0.5) {$0$};

		\node at (-0.5,0) {$-1$};
		\node at (-0.5,2.5) {$0$};
		\node at (-0.5,5) {$1$};
		\node at (-0.5,7.5) {$2$};
		
		\draw[fill = Mahogany!10!white] (3,2) rectangle (4,4);
		\draw[fill = Mahogany!10!white] (5,1) rectangle (7,3);
		\draw[fill = Mahogany!10!white] (2.5,0.2) rectangle (4.5,1.8);
		
		\node at (6,2) {$r_1$};
		\node at (3.5,3) {$r_2$};
		\node at (3.5,1) {$r_3$};
	\end{scope}

	\draw [|->] (5.5,2.5) -- (6.5,2.5);

	\begin{scope}[scale=0.5,xshift=14cm]
	\draw (1,0) -- (1,2.5) -- (2.5,2.5) -- (2.5,5) -- (2, 5) -- (2, 7.5)  -- (7.5,7.5) -- (7.5,2.5) -- (8,2.5) -- (8,0) -- (2,0);
	
	\draw (0,0) rectangle (10,10);

		\draw [dashed] (0.5,7.5) -- (0.5,10);

		\node at (1,-0.5) {$-1-a$};
		\node at (0.5,10.5) {$-1-c-t_\mr{ul}$};
		\node at (9,-0.5) {$1+b+t_r$};

		\node at (5,-0.5) {$0$};

		\node at (-0.5,0) {$-1$};
		\node at (-0.5,2.5) {$0$};
		\node at (-0.5,5) {$1$};
		\node at (-0.5,7.5) {$2$};

	\draw[fill = Mahogany!10!white] (3,2) rectangle (4,4);
	\draw[fill = Mahogany!10!white] (5,1) rectangle (7,3);
	\draw[fill = Mahogany!10!white] (2.5,0.2) rectangle (4.5,1.8);

	\draw[fill = Mahogany!10!white] (0.5,5) rectangle (2,7.5);
	\draw[fill = Mahogany!10!white] (8,0) rectangle (9,2.5);
	
	\node at (6,2) {$r_1$};
	\node at (3.5,3) {$r_2$};
	\node at (3.5,1) {$r_3$};
	\node at (1.25,6.25) {$r_\mr{ul}$};
	\node at (8.5,1.25) {$r_r$};
\end{scope}
	\end{tikzpicture}
\caption{Heuristically, the result of using the ``upper left" $\overline{\gR}$-module structure with $(t_\mr{ul},r_\mr{ul})$ and the right $\overline{\gR}$-module structure with $(t_r,r_r)$. 
}
\label{fig:example-special-left}
\end{figure}

\begin{lemma}\label{lem:AEqR}
There is a zig-zag of weak equivalences of $\overline{\gR}$-bimodules between $A(\gR)$ and $\overline{\gR}$. 
\end{lemma}
\begin{proof}
Firstly, the inclusion $A^0(\gR) \to A(\gR)$ is a morphism of $\overline{\gR}$-$\ogR$-bimodules where the left $\gR$-module structure on $A(\gR)$ is the ``lower left'' one~(\ref{eq:3}), and we claim that it is a weak equivalence. To see this, we will show that the inclusion $V^0 \subset V$ of right $E_2$-functors has a homotopy inverse as such. This then induces a homotopy inverse of $A^0(\gR) \to A(\gR)$, and a homotopy equivalence is a weak equivalence. 

The homotopy inverse is in fact given by a weak deformation retraction of $V$ into $V^0$, as follows. There is a path of self-embeddings of $V_{a,b}^c$ given by first scaling horizontally until it fits into $[-1,1] \times [-1,2]$, and then scaling vertically until it fits into $[-1,1] \times [-1,0]$; more formally, it is given by the 1-parameter family:
\begin{align*}
\rho_{a,b}^c(t) \colon V_{a,b}^c &\lra V_{a,b}^c\\
(x,y) &\longmapsto \begin{cases}
\big((1+ 2t (\frac{1}{1+\max(a,b,c)}-1))\cdot x, y\big) & t \in [0,\tfrac{1}{2}],\\
\big(\frac{1}{1+\max(a,b,c)} \cdot x, (2-2t+(2t-1)\tfrac{1}{3})\cdot(y+1)-1\big) & t \in [\tfrac{1}{2},1].
\end{cases}
\end{align*}
This expression makes clear it depends continuously on $a,b,c$. We use this family of embeddings to define the homotopy
\begin{align*}
h \colon [0,1] \times V &\lra V\\
(t, (a,b,c;e)) &\longmapsto \big(a,b,\max((1+ 2t (\tfrac{1}{1+\max(a,b,c)}-1))\cdot c,0), \rho_{a,b}^c(t)\circ e\big).
\end{align*}
As each $\rho_{a,b}^c(t)$ is given by a vertical and horizontal scaling and a translation, this is a homotopy through morphisms of right $E_2$-functors. It starts at the identity, and $h(1,-)$ has image in $V^0$. Furthermore, each $h(t,-)$ sends $V^0$ into $V^0$, so $h_{\vert [0,1] \times V^0}$ is a homotopy from the identity map of $V^0$ to $h(1,-)_{\vert V^0}$.  As required, $h(1,-)\colon V \to V^0$ is therefore homotopy inverse to the inclusion, as right $E_2$-functors.

It remains to compare the $\overline{\gR}$-bimodules $A^0(\gR)$ and $\overline{\gR}$. By definition we have 
\[V^0(\gR) = \bigsqcup_{n \geq 0} \tilde V^0(n) \times_{\beta_n} \gR^{\otimes n},\]
and a coequalizer diagram
\[\begin{tikzcd}V^0(E_2(\gR)) \arrow[shift left=.5ex]{r} \arrow[shift left=-.5ex]{r} & V^0(\gR) \arrow{r} & A^0(\gR).\end{tikzcd}\]

There is a morphism
\[V^0(\gR) \lra [0,\infty)^2 \times \bunit \sqcup ([0,\infty)^2 \setminus \{(0,0)\}) \times \gR\]
given heuristically by
\[((a,b,0; e), r_1, \ldots, r_n) \mapsto \begin{cases}
((a, b), \bunit) & \text{if $n=0$,}\\
((a, b), \mu(e; r_1, \ldots, r_n)) & \text{if $n>0$,}
\end{cases}\]
where $\mu(e; r_1, \ldots, r_n)$ is given by considering $e \colon \sqcup_n I^2 \to [-1-a, 1+b] \times [-1,0]$ as an embedding of $n$ little cubes in a large cube and (identifying $[-1-a, 1+b] \times [-1,0] \cong [0,1]^2$ linearly) applying the operadic multiplication. This descends to a morphism
\[A^0(\gR) \lra [0,\infty)^2 \times \bunit \sqcup ([0,\infty)^2 \setminus \{(0,0)\}) \times \gR,\]
which is an isomorphism. Under this isomorphism, and \eqref{eq:TildeUnderlying}, there is a morphism
\[A^0(\gR) \lra \overline{\gR} = (\{0\} \times \bunit) \sqcup ((0,\infty) \times \gR)\]
given by the map $[0,\infty)^2 \times \bunit \to \{0\} \times \bunit$ induced by the projection $[0,\infty)^2 \to \{0\}$ on the first term and $(a,b;r) \mapsto (a+b+2; r)$ on the second term. This is easily seen to be a morphism of $\overline{\gR}$-bimodules, and is a homotopy equivalence and hence weak equivalence.
\end{proof}

\begin{proof}[Proof that $A(\gR)$ is an adapter]
  We have constructed the $(\ogR \otimes \ogR)$-$\ogR$-bimodule structure and proved that it is equivalent to $\ogR$ as an $(\ogR \otimes \mathbbm{1})$-$\ogR$-bimodule.  It remains to see that the two maps $\mu_\mathrm{ll}, \mu_\mathrm{ul}\colon \ogR \otimes A(\gR) \to A(\gR)$ are homotopic as maps of right $\ogR$-modules.
As $A(\gR) \simeq \overline{\gR}$ as right $\overline{\gR}$-modules, we have a bijection $[\overline{\gR} \otimes A(\gR), A(\gR)]_{\cat{mod}\text{-}\overline{\gR}} \cong [\overline{\gR}, A(\gR)]_\sfC$. Under this bijection, the two left $\overline{\gR}$-module structure maps are given heuristically by sending $(t, r)$ to $(t, 0, 0 ; e_1)$ and $(0,0,t,e_2)$ where the $e_i \colon I^2 \to \bR^2$ are two rectilinear embeddings. These maps are evidently homotopic, by sliding one such little cube to the other inside $V_{t,0}^t$.  As already mentioned, this construction does not obviously satisfy the cofibrancy requirements, so we redefine $A(\gR)$ by cofibrantly replacing the result of the geometric construction above.
\end{proof}
\section{Indecomposables and the bar construction}\label{sec:BarConstr}

As before, we work in the category $\sfC = \sfS^\sfG$ with $\sfS$ satisfying the axioms of Sections \ref{sec:axioms-of-cats} and \ref{sec:axioms-of-model-cats}. In Section \ref{sec:derived-functors} we defined the derived indecomposables $Q^T_\bL(-)$ for a monad $T$, and described two ways of computing it. Firstly, given a $T$-algebra cell structure on $\gX \in \Alg_T(\sfC)$ there is a associated ordinary cell structure on $Q^T_\bL(\gX) \in \sfC_*$. Secondly, we may choose a free simplicial resolution $\epsilon \colon \gX_\bullet \to \gX$, so that $Q^T_\bL(\gX) \simeq \gr{Q^T(\gX_\bullet)} \in \sfC_*$. These methods work for quite general monads $T$. In this section we shall describe a third way of computing the derived indecomposables particular to the case that $T=E_k$ is the non-unital little $k$-cubes monad: the $k$-fold iterated bar construction. 

We also discuss a number of related results; the effect of bar constructions on maps, $E_k$-algebra structures on iterated indecomposables, and group completion.

\subsection{The iterated bar construction}\label{sec:iterated-bar-def} 

The $k$-fold iterated bar construction is a flexible version of the ordinary bar construction, applied in $k$ directions at the same time. For simplicity of exposition we will first describe this construction under the assumption that $\sfC$ is symmetric monoidal; at the end of this section we will explain the mild changes to be made if this category is only monoidal or braided monoidal.

The iterated bar construction will be described in terms of grids in the $k$-dimensional cube.
\newglossaryentry{grids}{%
	name={\ensuremath{\cP_k(p_1,\ldots,p_k)}},
	description={Space of grids used in iterated bar construction},
	type=symbols
}
\begin{definition}\label{def:pk} 
Let us write $\cP_k(p_1,\ldots,p_k) \subset \smash{\prod_{j=1}^k} \bR^{p_j+1}$ for the subspace of $k$-tuples $\{t_i^j\}_{1 \leq j \leq k}$ of sequences $0<t_0^{j}<\ldots<t_{p_j}^j<1$. The assignment
\[[p_1,\ldots,p_k] \mapsto \cP_k(p_1,\ldots,p_k)\]
forms a $k$-fold semi-simplicial space, if we define the $i$th face map in the $j$th direction $d_i^j$  to be given by forgetting $t_i^j$.
\end{definition}

\begin{figure}[h]
	\centering
	\begin{tikzpicture}
	
	\foreach \x in {.5,2}
	\foreach \y in {1,2.5}
	{
		\draw[black] (\x,0) -- (\x,3);
		\draw[black] (0,\y) -- (3,\y);
	}
	
	\draw[dotted] (0,0) rectangle (3,3);
	
	\node[below] at (0,0) {$0$};
	\node[below] at (.5,0) {$t_0^1$};
	\node[below] at (2,0) {$t_1^1$};
	\node[below] at (3,0) {$1$};
	
	\node[left] at (0,3) {$1$};
	\node[left] at (0,2.5) {$t_1^2$};
	\node[left] at (0,1) {$t_0^2$};
	\node[left] at (0,0) {$0$};

	\draw [|->] (3.25,1.5) --  node[anchor=north] {$d^1_0$} ++(1,0);
	
	\foreach \x in {7}
	\foreach \y in {1,2.5}
	{
		\draw[black] (\x,0) -- (\x,3);
		\draw[black] (5,\y) -- (8,\y);
	}
	
	\draw[dotted] (5,0) rectangle (8,3);
	
	\node[below] at (5,0) {$0$};
	\node[below] at (7,0) {$t_0^1$};
	\node[below] at (8,0) {$1$};
	
	\node[left] at (5,3) {$1$};
	\node[left] at (5,2.5) {$t_1^2$};
	\node[left] at (5,1) {$t_0^2$};
	\node[left] at (5,0) {$0$};
	\end{tikzpicture}
	\caption{The face map of $d^1_0 \colon \cP(1,1) \to \cP(0,1)$.}
	\label{fig:cp-face-maps}
\end{figure}

It will be helpful to think of an element of $\cP_k(p_1,\ldots,p_k)$ as a collection of hyperplanes $\bR^{j-1} \times \{t^j_i\} \times \bR^{k-j}$ cutting $I^k$ into $\prod_{j=1}^k (p_j+2)$ $k$-cubes, as depicted in Figure \ref{fig:cp-face-maps}. Of these, the $p_1 \cdots p_k$ cubes given by 
\[\prod_{j=1}^k [t^j_{q_j-1},t^j_{q_j}] \quad \text{for} \quad (q_j)_{j=1}^k \in \prod_{j=1}^k \{1,\ldots,p_j\}\]
will play a major role (these are the shaded cubes in Figure \ref{fig:bar-construction}). The face map $d_i^j$ for $0<i<p_j$ merges the cubes $(q_1,\ldots,q_{j-1},i,q_{j+1},\ldots,q_k)$ and $(q_1,\ldots,q_{j-1},i+1,q_{j+1},\ldots,q_k)$. This merging may be interpreted as elements $\delta_i^j(q_1,\ldots,\hat{q}_j,\ldots,q_k) \in \cC_k(2)$, as in Figure \ref{fig:cp-delta-ij}.

\begin{figure}[h]
	\centering
	\begin{tikzpicture}
	
	\fill[Mahogany!10!white] (0,1) rectangle (2,2.5);
	
	\foreach \x in {.5,2}
	\foreach \y in {1,2.5}
	{
		\draw[black] (\x,0) -- (\x,3);
		\draw[black] (0,\y) -- (3,\y);
	}
	
	\draw[dotted] (0,0) rectangle (3,3);
	
	\node[below] at (0,0) {$0$};
	\node[below] at (.5,0) {$t_0^1$};
	\node[below] at (2,0) {$t_1^1$};
	\node[below] at (3,0) {$1$};
	
	\node[left] at (0,3) {$1$};
	\node[left] at (0,2.5) {$t_1^2$};
	\node[left] at (0,1) {$t_0^2$};
	\node[left] at (0,0) {$0$};
	
	\filldraw[fill=Mahogany!10!white, draw=black] (9,0) rectangle (12,3);
	
	\draw[black] (5,2) rectangle (8,5);
	\draw[black] (5,-2) rectangle (8,1);
	
	\draw[right hook->] (6.5,3.5) -- node[anchor=north, pos=0.25] {$e_1$} ++(3,-2);
	\draw[right hook->] (6.5,-.5) -- node[anchor=south, pos=0.25] {$e_2$} ++(4.5,2);
	
	\node[below] at (9,0) {$0$};
	\node[left] at (9,0) {$0$};
	\node[below] at (9.75,0) {$\frac{t_0^1}{t_0^1+t_1^1}$};
	\node[below] at (12,0) {$1$};
	\node[left] at (9,3) {$1$};
	
	\draw[black] (9.75,0) -- (9.75,3);
	\end{tikzpicture}
	\caption{The element $\delta_0^1(1)$ for the face map as in Figure \ref{fig:cp-face-maps}.}
	\label{fig:cp-delta-ij}
\end{figure}

More precisely the two rectilinear embeddings $e_1,e_2 \colon I^k \hookrightarrow I^k$ forming $\delta_i^j  \in \cC_k(2)$ are given by 
\begin{equation}\label{eqn:bar-embedding-deltaij}\begin{aligned} e_1(x_1,\ldots,x_k) &= \left(x_1, \ldots, x_{j-1},\frac{t^j_i-t^j_{i-1}}{t^j_{i+1}-t^j_{i-1}}x_j, x_{j+1}, \ldots,x_k\right),\\
e_2(x_1,\ldots,x_k) &= \left(x_1,\ldots,x_{j-1},\frac{(t^j_{i+1}-t^j_{i})x_j + (t^j_{i}-t^j_{i-1})}{t^j_{i+1}-t^j_{i-1}}, x_{j+1},\ldots,x_k\right).\end{aligned}\end{equation}

For later use, we remark that there is also an element $\delta \in \cC_k(p_1 \cdots p_k)$ consisting of the embeddings $e_{q_1,\ldots,q_k} \colon I^k \hookrightarrow I^k$ for $(q_j)_{j=1}^k \in \prod_{j=1}^k \{1,\ldots,p_j\}$ given by
\begin{equation}\label{eqn:bar-embedding-delta} e_{q_l,\ldots,q_k} \colon (x_1,\ldots,x_k) \mapsto \left((t^1_{q_1+1}-t^1_{q_1})x_1+t^1_{q_1},\ldots,(t^k_{q_k+1}-t^k_{q_k})x_k+t^k_{q_k}\right),\end{equation}
whose image consists of the $p_1 \cdots p_k$ inner cubes, where these embeddings are ordered lexicographically by $(q_1,\ldots,q_k)$.

We wish to define the $k$-fold bar construction for an $E^+_k$-algebra $\gR$ with augmentation $\epsilon \colon \gR \to \bunit$, but it is no more difficult to make a definition for an arbitrary morphism of $E_k^+$-algebras.

\begin{definition} \label{def:kfold-bar-augmented}
	Let $f \colon \gR \to \gS$ be a morphism of $E_k^+$-algebras. Then $B^{E_k}_{\bullet,\ldots,\bullet}(f)$ is the $k$-fold semi-simplicial object with $B^{E_k}_{p_1,\ldots,p_k}(f) \coloneqq \cP_k(p_1,\ldots,p_k) \times G_{p_1, \ldots, p_k}(f)$, where
\[G_{p_1, \ldots, p_k}(f) \coloneqq \bigotimes_{q_1=0}^{p_1+1} \cdots \bigotimes_{q_k=0}^{p_k+1} B_{p_1, \ldots, p_k}^{q_1,\ldots,q_k}\]
and $B_{p_1, \ldots, p_k}^{q_1, \ldots, q_k}$ is $\gR$ if $1 \leq q_j \leq p_j$ for all $j$, and $\gS$ otherwise (see Figure \ref{fig:bar-construction} for an example). 

The $i$th face map $d_i^j$ in the $j$th direction
	\[d_i^j \colon B^{E_k}_{p_1,\ldots,p_k}(f) \lra B^{E_k}_{p_1,\ldots, p_{j-1}, p_j-1, p_{j+1},p_k}(f)\]
	is given by the face map of Definition \ref{def:pk} on the first factor and then, by adjunction, by the map of simplicial sets
\begin{align*}
\cP_k(p_1,\ldots,p_k) &\lra \cC_k(2) \overset{\alpha}\lra \mr{Map}_\sfC(G_{p_1, \ldots, p_k}(f), G_{p_1,\ldots, p_{j-1}, p_j-1, p_{j+1},p_k}(f))\\
\{t_i^j\} &\longmapsto \delta^j_i 
\end{align*}
where $\alpha$ is given as follows: 
\begin{enumerate}[(i)]
\item For $i=0$, the maps
\begin{multline*}
\cC_k(2) \lra \cE_{\gS}(2) = \mr{Map}_\sfC(\gS \otimes \gS, \gS) \overset{(\gS \otimes \epsilon)^*}\lra \mr{Map}_\sfC(\gS \otimes \gR, \gS)\\
 = \mr{Map}_\sfC(B_{p_1, \ldots, p_k}^{q_1,\ldots,q_{j-1},0,q_{j+1},\ldots,q_k} \otimes B_{p_1, \ldots, p_k}^{q_1,\ldots,q_{j-1},1,q_{j+1},\ldots,q_k}, B_{p_1,\ldots, p_{j-1}, p_j-1, p_{j+1},p_k}^{q_1,\ldots,q_{j-1},0,q_{j+1},\ldots,q_k})
\end{multline*}
and the evident identity maps on the remaining factors.

\item For $0 < i < p_j$, the maps
\begin{multline*}
\cC_k(2) \lra \cE_{\gR}(2) = \mr{Map}_\sfC(\gR \otimes \gR, \gR) \\
 = \mr{Map}_\sfC(B_{p_1, \ldots, p_k}^{q_1,\ldots,q_{j-1},i,q_{j+1},\ldots,q_k} \otimes B_{p_1, \ldots, p_k}^{q_1,\ldots,q_{j-1},i+1,q_{j+1},\ldots,q_k}, B_{p_1,\ldots, p_{j-1}, p_j-1, p_{j+1},p_k}^{q_1,\ldots,q_{j-1},i,q_{j+1},\ldots,q_k})
\end{multline*}
and the evident identity maps on the remaining factors.

\item For $i=p_j$, the maps
\begin{multline*}
\cC_k(2) \lra \cE_{\gS}(2) = \mr{Map}_\sfC(\gS \otimes \gS, \gS) \overset{(\epsilon \otimes \gS)^*}\lra \mr{Map}_\sfC(\gR \otimes \gS, \gS)\\
 = \mr{Map}_\sfC(B_{p_1, \ldots, p_k}^{q_1,\ldots,q_{j-1},p_j,q_{j+1},\ldots,q_k} \otimes B_{p_1, \ldots, p_k}^{q_1,\ldots,q_{j-1},p_j+1,q_{j+1},\ldots,q_k}, B_{p_1,\ldots, p_{j-1}, p_j-1, p_{j+1},p_k}^{q_1,\ldots,q_{j-1},p_j,q_{j+1},\ldots,q_k})
\end{multline*}
and the evident identity maps on the remaining factors.
\end{enumerate}

	We write $B^{E_k}(f) \coloneqq \fgr{B_{\bullet, \ldots, \bullet}^{E_k}(f)} \in \sfC$ and call this the \emph{$k$-fold iterated bar construction}\index{bar construction!iterated}. It is natural in commutative squares of morphisms of $E_k^+$-algebras:
	\[\begin{tikzcd} \gR \rar{f} \dar & \gS \dar \\
		\gR' \rar{f'} & \gS' \end{tikzcd} \quad\text{yields}  \quad B^{E_k}(f) \to B^{E_k}(f').\]
	
\end{definition}

\begin{figure}[h]
	\centering
\begin{tikzpicture}

\fill[Mahogany!2!white] (0,0) rectangle (6,6);
\fill[Mahogany!10!white] (0.5,1) rectangle (5,5);

\foreach \x in {.5,2,4,5}
\foreach \y in {1,2.5,3.5,5}
{
	\draw[black] (\x,0) -- (\x,6);
	\draw[black] (0,\y) -- (6,\y);
}

\draw[dotted] (0,0) -- (0,6);
\draw[dotted] (6,0) -- (6,6);
\draw[dotted] (0,0) -- (6,0);
\draw[dotted] (0,6) -- (6,6);

\node at (0.25,.5) {$\gS$};
\node at (1.25,.5) {$\gS$};
\node at (3,.5) {$\gS$};
\node at (4.5,.5) {$\gS$};
\node at (5.5,.5) {$\gS$};

\node at (0.25,1.75) {$\gS$};
\node at (1.25,1.75) {$\gR$};
\node at (3,1.75) {$\gR$};
\node at (4.5,1.75) {$\gR$};
\node at (5.5,1.75) {$\gS$};

\node at (0.25,3) {$\gS$};
\node at (1.25,3) {$\gR$};
\node at (3,3) {$\gR$};
\node at (4.5,3) {$\gR$};
\node at (5.5,3) {$\gS$};

\node at (0.25,4.25) {$\gS$};
\node at (1.25,4.25) {$\gR$};
\node at (3,4.25) {$\gR$};
\node at (4.5,4.25) {$\gR$};
\node at (5.5,4.25) {$\gS$};

\node at (0.25,5.5) {$\gS$};
\node at (1.25,5.5) {$\gS$};
\node at (3,5.5) {$\gS$};
\node at (4.5,5.5) {$\gS$};
\node at (5.5,5.5) {$\gS$};

\node[below] at (0,0) {$0$};
\node[below] at (.5,0) {$t_0^1$};
\node[below] at (2,0) {$t_1^1$};
\node[below] at (4,0) {$t_2^1$};
\node[below] at (5,0) {$t_3^1$};
\node[below] at (6,0) {$1$};

\node[left] at (0,6) {$1$};
\node[left] at (0,1) {$t_0^2$};
\node[left] at (0,2.5) {$t_1^2$};
\node[left] at (0,3.5) {$t_2^2$};
\node[left] at (0,5) {$t_3^2$};
\node[left] at (0,0) {$0$};
\end{tikzpicture}
\caption{An illustration of $B^{E_2}_{3,3}(f)$.}
\label{fig:bar-construction}
\end{figure}
\newglossaryentry{bek}{%
	name={\ensuremath{B^{E_k}(\gR,\epsilon)}},
	description={Iterated bar construction of augmented $E_k$-algebra $\gR$},
	type=symbols
}
In particular, if $\epsilon \colon \gR \to \bunit$ is an augmented $E_k^+$-algebra then we write
\[B_{\bullet, \ldots, \bullet}^{E_k}(\gR, \epsilon) \coloneqq B_{\bullet, \ldots, \bullet}^{E_k}(\epsilon)\]
for the associated $k$-fold bar construction, and similarly for its geometric realization. The unit map $1 \colon \bunit \to \gR$ and augmentation $\epsilon \colon \gR \to \bunit$ are maps of augmented $E_k^+$-algebras, where $\bunit$ has augmentation given by the identity. Thus we get maps
\[B^{E_k}(\bunit,\epsilon_\bunit) \lra \gls{bek} \lra B^{E_k}(\bunit,\epsilon_\bunit),\]
whose composition is the identity.

\begin{lemma}\label{lem:bek-unit}  $B^{E_k}(\bunit,\epsilon_\bunit) \simeq \bunit$.
\end{lemma}

\begin{proof}In terms of the copowering over simplicial sets $\bunit \cong \ast$, so the thick geometric realization of $B^{E_k}(\bunit,\epsilon_\bunit)$ in $\sfC$ is isomorphic to $\bunit$ copowered with the thick geometric realization of the $k$-fold semi-simplicial simplicial set whose $(p_1,\ldots,p_k)$-simplices are given by $\cP_k(p_1,\ldots,p_k)$. This is contractible in each multisimplicial degree, hence so is its thick geometric realization. This induces the desired weak equivalence.
\end{proof}

\newglossaryentry{bekred}{%
	name={\ensuremath{\tilde{B}^{E_k}(\gR,\epsilon)}},
	description={Reduced iterated bar construction},
	type=symbols
}
\newglossaryentry{bekredcam}{%
	name={\ensuremath{\tilde{B}^{E_k}(\gR)}},
	description={Reduced iterated bar construction of canonical augmentation},
	type=symbols
}
\begin{definition}\label{def:reduced-kfold-bar}
We define \emph{reduced $k$-fold bar construction}\index{bar construction!reduced} $\gls{bekred} \in \sfC_\ast$ to be the cofiber of the map $B^{E_k}(\bunit,\epsilon_\bunit) \to B^{E_k}(\gR,\epsilon)$ induced by the unit.
\end{definition}

In Definition \ref{def:split-augmented}, we called $(\gR,\epsilon)$ split augmented if the induced map $(I(\gR))^+ \to \gR$ is an isomorphism, where the augmentation $\epsilon \colon \gR \to \bunit$ is used to define the augmentation ideal $I(\gR)$.

\begin{lemma}\label{lemma:bek-cofibrant} If $(\gR,\epsilon)$ is split augmented, then we have that 
	\[B^{E_k}(\gR,\epsilon)_+ \cong B^{E_k}(\bunit,\epsilon_\bunit)_+ \vee \tilde{B}^{E_k}(\gR,\epsilon).\]
If additionally $\gR$ is cofibrant in $\sfC$, then $B^{E_k}(\bunit,\epsilon_\bunit) \to B^{E_k}(\gR,\epsilon)$ is a cofibration.\end{lemma}

\begin{proof}If $\gR$ is split augmented, the map $B^{E_k}(\bunit,\epsilon_\bunit) \to B^{E_k}(\gR,\epsilon)$ is induced by a levelwise inclusion of a term into a coproduct, and the isomorphism follows from a levelwise isomorphism of pointed $k$-fold semi-simplicial objects.
	
If $U^{E_k}(\gR)$ is cofibrant, then terms of $B^{E_k}(\gR,\epsilon)$ which are not in $B^{E_k}(\bunit,\epsilon_\bunit)$, are cofibrant and thus the inclusion is a cofibration. By Lemma \ref{lem:thick-geom-rel-cofibrations} the map $B^{E_k}(\bunit,\epsilon_\bunit) \to B^{E_k}(\gR,\epsilon)$ is then a cofibration.\end{proof}

As explained in Section \ref{sec:non-unitary}, if $\cat{C}$ is pointed then the unitalization $\gR^+$ of an $E_k$-algebra $\gR$ may be endowed with the canonical augmentation $\epsilon_\mr{can}$, and this is always split augmented. In that case we will simplify notation and write
\begin{equation*}
\gls{bekredcam} \coloneqq \tilde{B}^{E_k}(\gR^+,\epsilon_\mr{can}) \in \sfC.
\end{equation*}
More generally if $\cat{C}$ is not pointed then we can instead consider $\gR_+$ as an $E_k$-algebra in the pointed category $\cat{C}_*$, and set
\begin{equation}\label{def:tildebek} 
\tilde{B}^{E_k}(\gR) \coloneqq \tilde{B}^{E_k}(\gR^+_+,\epsilon_\mr{can}) \in \sfC_\ast.
\end{equation}
The previous lemmas imply that if $\gR$ is cofibrant in $\sfC$ then $B^{E_k}(\gR^+_+,\epsilon_\mr{can}) \simeq \bunit_+ \vee \tilde{B}^{E_k}(\gR)$,
so we do not lose any homotopy-theoretical information when computing $\tilde{B}^{E_k}$ instead of $B^{E_k}$. In arguments later in this section, we shall use the following description of $\tilde{B}^{E_k}(\gR)$, which does not make reference to augmented $E_k^+$-algebras.
	
\begin{lemma}\label{lem:reduced-kfold-bar}
Let $\gR$ be a non-unital $E_k$-algebra which is cofibrant in $\sfC$. Then we may compute $\tilde{B}^{E_k}(\gR)$ as the geometric realization of the following $k$-fold semi-simplicial object $\tilde{B}^{E_k}_{\bullet,\ldots,\bullet}(\gR)$ in $\sfC_\ast$. It has the $(p_1,\ldots,p_k)$-simplices $\tilde{B}^{E_k}_{p_1,\ldots,p_k}(\gR)$ given by $\ast$ if any $p_i=0$, and otherwise by the quotient of
	\[\cP_k(p_1,\ldots,p_k) \times \bigotimes_{q_1 =1}^{p_1} \cdots \bigotimes_{q_k=1}^{p_k} (\bunit \sqcup \gR)\]
by the subobject \[\cP_k(p_1,\ldots,p_k) \times \bigotimes_{q_1 =1}^{p_1} \cdots \bigotimes_{q_k=1}^{p_k} \bunit.\]
The $i$th face map $d_i^j$ in the $j$th direction
\[d_i^j \colon \tilde{B}^{E_k}_{p_1,\ldots,p_k}(\gR) \lra \tilde{B}^{E_k}_{p_1,\ldots, p_{j-1}, p_j-1, p_{j+1}, \ldots, p_k}(\gR),\]
is given as in Definition \ref{def:kfold-bar-augmented}, with the variation that $d_0^j$ is given by applying the augmentation $\epsilon \colon \bunit \sqcup \gR \to \bunit$ to those factors with $q_j=1$, and $d^j_{p_j}$ is given by applying the augmentation $\epsilon \colon \bunit \sqcup \gR \to \bunit$ to those factors with $q_j=p_j$.
\end{lemma}

The following theorem is the main result of this section. It says one can compute the derived $E_k$-indecomposables in terms of the reduced $k$-fold bar construction. Instances of this result are due to Getzler--Jones \cite{GetzlerJones}, Basterra--Mandell \cite{BM}, Fresse \cite{Fresse}, and Francis \cite{FrancisPaper}.

\begin{theorem}\label{thm:BarHomologyIndec}
There is a zig-zag (\ref{eqn:zigzag-bk-qek}) of natural transformations 
	\[\tilde{B}^{E_k}(-) \Leftarrow \cdots \Rightarrow S^k \wedge Q_\bL^{E_k}(-)\]
of functors $\Alg_{E_k}(\sfC) \to \sfC_\ast$, which are weak equivalences when evaluated on objects which are cofibrant in $\sfC$.
\end{theorem}

We will prove this theorem in Section \ref{sec:ProofBEqQ} after some preparation.

\subsubsection{Modification for $k$-monoidal categories with $k=1,2$}\label{sec:bar-modification}

For concreteness we explain how Definition \ref{def:kfold-bar-augmented} must be modified; there are completely analogous modifications to Lemma \ref{lem:reduced-kfold-bar}.

If the category $\sfC$ is only 1-monoidal then it only makes sense to consider $E_1$-algebras, which must be done as described in Section \ref{sec:EkAlgForSmallk} using the monad $E_1^{\cat{FB_1}}$. In this case there is no difficulty in following the construction in Definition \ref{def:kfold-bar-augmented}. For a given grid $(t^1_0, \ldots, t^1_{p_1}) \in \mathcal{P}_1(p_1)$ the embeddings $e_1, e_2 \colon I \hookrightarrow I$ forming $\delta^1_i \in \mathcal{C}_1(2)$ satisfy $e_1(0)=0 < \tfrac{t_i^1 - t^1_{i-1}}{t^1_{i+1} - t^1_{i-1}}= e_2(0)$ and so $\delta^1_i$ lies in $\cC^\cat{FB_1}_1(2)$. Then the face map $d^1_i \colon B_{p_1}^{E_1}(f) \to B_{p_1-1}^{E_1}(f)$ is given in the same way, using $\{t_i^1\} \mapsto \delta^1_i \colon \cP_1(p_1) \to \cC^\cat{FB_1}_1(2)$, the map $\cC^\cat{FB_1}_1(2) \to \mr{Map}_\sfC(B_{p_1}^i \otimes B_{p_1}^{i+1}, B_{p_1-1}^i)$ given by the $E_1^{\cat{FB_1}}$-algebra structure on $\gR$ and $\gS$, and the map of simplicial sets
\begin{multline*}
\prod_{q_1=0}^{i-1} \mr{Map}_\sfC(B_{p_1}^{q_1}, B_{p_1-1}^{q_1}) \times \mr{Map}_\sfC(B_{p_1}^i \otimes B_{p_1}^{i+1}, B_{p_1-1}^i) \times \prod_{q_1=i+1}^{p_1} \mr{Map}_\sfC(B_{p_1}^{q_1+1}, B_{p_1-1}^{q_1})\\
 \lra \mr{Map}_\sfC\left(\bigotimes_{q_1=0}^{p_1+1} B_{p_1}^{q_1}, \bigotimes_{q_1=0}^{p_1} B_{p_1-1}^{q_1}\right)
\end{multline*}
given by multiplication. (In Definition \ref{def:kfold-bar-augmented} we had employed this construction, but implicitly followed it by a permutation of the factors in the target in order to identify the target with $G_{p_1,\ldots, p_{j-1}, p_j-1, p_{j+1},p_k}(f)$. In this case no permutation is necessary.)

The remaining case is when the category $\sfC$ is 2-monoidal and we consider $E_2$-algebras, which must be done as described in Section \ref{sec:EkAlgForSmallk} using the monad $E_2^{\cat{FB_2}}$. In this case a little care must be taken, because in defining $G_{p_1, p_2}(f)$ we have needed to impose a linear ordering of the terms $B_{p_1, p_2}^{q_1, q_2}$ when they are more naturally arranged in a 2-dimensional grid. The linear ordering we have chosen is a convention, but this convention does dictate a choice of lift
\[\cP_2(p_1, p_2) \lra \cC_2^{\cat{FB_2}}(2)\]
of the map $\{t^j_i\} \mapsto \delta^j_i \colon \cP_2(p_1, p_2) \to \cC_2(2)$. Namely: the point $\delta^1_i \in \cC_2(2)$ already lies in the subspace $\cC_1^{\cat{FB_1}}(2) \subset \cC_2(2)$ so equipping it with the constant path defines an element of $\cC_2^{\cat{FB_2}}(2)$; the point $\delta^2_i \in \cC_2(2)$ should be equipped with the path where the top cube moves rightwards and then down, to give a point in $\cC_1^{\cat{FB_1}}(2) \subset \cC_2(2)$. Using this choice, the face maps $d^2_i$ are as in Definition \ref{def:kfold-bar-augmented} using the map $\cC^\cat{FB_2}_1(2) \to \mr{Map}_\sfC(B_{p_1, p_2}^{q_1,i} \otimes B_{p_1, p_2}^{q_1,i+1}, B_{p_1-1, p_2}^{q_1,i})$ given by the $E_2^{\cat{FB_2}}$-algebra structure on $\gR$ and $\gS$. There is no need to explicitly use the braiding, as the terms which are being multiplied together are adjacent in the linear ordering of the terms in $G_{p_1, p_2}(f)$. On the other hand, for the face map $d^1_i$ the same construction naturally defines a map from
\[\cP_2(p_1, p_2) \times \left(\bigotimes_{q_1=0}^{i-1}\bigotimes_{q_2=0}^{p_2+1} B_{p_1, p_2}^{q_1, q_2}\right) \otimes \bigotimes_{q_2=0}^{p_2+1} (B_{p_1, p_2}^{i, q_2} \otimes B_{p_1, p_2}^{i+1, q_2}) \otimes \left(\bigotimes_{q_1=i+2}^{p_1+1}\bigotimes_{q_2=0}^{p_2+1} B_{p_1, p_2}^{q_1, q_2}\right)\]
to $\cP_2(p_1-1, p_2) \times \bigotimes_{q_1=0}^{p_1}\bigotimes_{q_2=0}^{p_2+1} B_{p_1-1, p_2}^{q_1, q_2} = \cP_2(p_1-1, p_2) \times G_{p_1, p_2}(f)$.
To identify the source with $\cP_2(p_1, p_2) \times G_{p_1, p_2}(f)$ it is necessary to choose a way to braid all the terms $B_{p_1, p_2}^{i+1, q_2}$ to the right of all the terms $B_{p_1, p_2}^{i, q_2}$. This choice is again dictated by our convention: in order for this $d^1_i$ to commute with the $d^2_{i'}$, we must precompose the map constructed so far with the isomorphism
\[\left(\bigotimes_{q_2=0}^{p_2+1} B_{p_1, p_2}^{i, q_2}\right) \otimes \left(\bigotimes_{q_2=0}^{p_2+1} B_{p_1, p_2}^{i+1, q_2}\right) \lra \bigotimes_{q_2=0}^{p_2+1} (B_{p_1, p_2}^{i, q_2} \otimes B_{p_1, p_2}^{i+1, q_2})\]
which braids the $B_{p_1, p_2}^{i+1, q_2}$'s \emph{in front of} the $B_{p_1, p_2}^{i, q_2}$'s. (This is with the convention that $\beta_{X,Y} \colon X \otimes Y \to Y \otimes X$ is represented by the braid where the strand labeled $X$ crosses in front of the strand labeled $Y$.)

\subsection{The bar construction for modules}Let us briefly draw a parallel. In Section \ref{sec:modules} we described how to associate a unital associative algebra $\overline{\gR}$ to an $E_1$-algebra $\gR$. A (left or right) module over $\gR$ is then a (left or right) $\overline{\gR}$-module in the usual sense. Furthermore, if $\sfC$ is pointed then we showed that there is a canonical augmentation $\epsilon \colon \overline{\gR} \to \bunit$, and so by Section \ref{sec:AssocModules} there is a notion of derived $\overline{\gR}$-module indecomposables $Q^{\overline{\gR}}_\bL$. In this case Corollary \ref{cor:DerivModIndecBar} is the analogue of Theorem \ref{thm:BarHomologyIndec}.

\subsection{The bar construction on free algebras}\label{sec:bar-on-free-algebras}

In order to prove Theorem \ref{thm:BarHomologyIndec} we will need to compute $\tilde{B}^{E_k}(\gE_k(X))$. In fact, for later use and because it is no more difficult, we shall explain how to compute $\tilde{B}^{E_k}(\gE_{n+k}(X))$. The following for $n=0$ has also been proved by Lurie \cite[Proposition 5.2.3.15]{LurieHA}. The reader familiar with $E_n$-algebras in spaces should compare Theorem \ref{thm:CalcFree} to the results in Section \ref{sec:GpCompletion} on group completion, where a \emph{different} augmentation is used. This highlights the role played by the augmentation.

\begin{theorem}\label{thm:CalcFree} There is a zig-zag (\ref{eqn:zig-zag-bek-fnk}) of natural transformations
\[\tilde{B}^{E_k}(\gE_{n+k}(-)) \Leftarrow \cdots \Rightarrow E_n(S^k \wedge (-)_+)\]
of functors from $\sfC$ to $\sfC_\ast$, which are weak equivalences on cofibrant objects.
\end{theorem}

The intermediate constructions involve spaces of little cubes in $I^n \times \bR^k$. It is easier to define such objects, maps and homotopies in $\cat{Top}$ rather than in $\cat{sSet}$. If so, we implicitly form the singular simplicial set before using copowering, and do not distinguish in notation between topological spaces and their singular simplicial sets.

For simplicity of exposition we will give the proof when $\sfC$ is symmetric monoidal, but the variations required to treat the 1- or 2-monoidal cases are routine (the following spaces of cubes should be defined in terms of $\mr{Emb}^{\mr{rect}, \cat{FB_k}}$ instead of $\mr{Emb}^{\mr{rect}}$).

\newglossaryentry{fnk}{%
	name={\ensuremath{F_{n,k}}},
	description={Symmetric sequence of $(n+k)$-cubes in $I^n \times \bR^k$},
	type=symbols
}
\newglossaryentry{partialfnk}{%
	name={\ensuremath{\partial F_{n,k}}},
	description={Sub-symmetric sequence of $F_{n,k}$},
	type=symbols
}
\newglossaryentry{fnkquotient}{%
	name={\ensuremath{F_{n,k}/\partial F_{n,k}}},
	description={Quotient of $F_{n,k}$ by $\partial F_{n,k}$},
	type=symbols
}
\begin{definition}\label{def:fnk-and-friends} We define the following three symmetric sequences:
	\begin{itemize}
		\item Let $\gls{fnk}$ denote the symmetric sequence in $\cat{Top}$ given by
		\begin{align*}F_{n,k}(i) \coloneqq \begin{cases} \varnothing & \text{if $i=0$,} \\
		\mr{Emb}^\mr{rect}(\sqcup_i I^{n+k},I^n \times \bR^k) & \text{otherwise,} \end{cases}\end{align*}
		see Figure \ref{fig:example-fnk} for an example.
		\item Let $\gls{partialfnk}(i)$ be the subspace of $F_{n,k}(i)$ where at least one cube lies entirely outside the interior of $I^n \times I^k$. This is preserved by the action of $\fS_i$ and hence we obtain a symmetric sequence $\partial F_{n,k}$ in $\cat{Top}$.
		\item Let $\gls{fnkquotient}$ denote the quotient symmetric sequence in $\cat{Top}_*$, whose value at $i$ is the pointed space $F_{n,k}(i)/\partial F_{n,k}(i)$.
	\end{itemize}
\end{definition}

\begin{figure}[h]
	\centering
	\begin{tikzpicture}
	
	\draw[dashed] (0,-1) -- (0,6);
	\draw[dotted] (0,0) -- (5,0);
	\draw[dotted] (0,5) -- (5,5);
	\draw (0,0) -- (0,5);
	\draw[dashed] (5,-1) -- (5,6);
	\draw (5,0) -- (5,5);
	
	\fill[Mahogany!10!white] (3,2) rectangle (4,5.5);
	\fill[Mahogany!10!white] (0.2,1) rectangle (0.8,3);
	\fill[Mahogany!10!white] (2.5,0.2) rectangle (4.5,1.8);
	
	\draw (3,2) rectangle (4,5.5);
	\draw (0.2,1) rectangle (0.8,3);
	\draw (2.5,0.2) rectangle (4.5,1.8);
	
	\node at (0.5,2) {$e_1$};
	\node at (3.5,3.75) {$e_2$};
	\node at (3.5,1) {$e_3$};
	
	\node[right] at (5.1,2.5) {$\in  F_{1,1}(3)$};

	\end{tikzpicture}
	\caption{An element of $F_{1,1}(3)$ which does not lie in $\partial F_{1,1}(3)$.}
	\label{fig:example-fnk}
\end{figure}

There is an inclusion, resp.\ quotient map, of symmetric sequences
\[\partial F_{n,k} \lra F_{n,k} \lra F_{n,k}/\partial F_{n,k}.\]

The symmetric sequence $F_{n,k}$ is a left module over $\cC_n$. The structure maps
\begin{equation}
\label{eqn:fnk-left-cn-module-structure}\mathcal{C}_n(j) \times F_{n,k}(i_1) \times \cdots \times F_{n,k}(i_j) \lra F_{n,k}(i_1+\cdots+i_j),\end{equation}
for $j \geq 1$ are given by sending $e \in \mathcal{C}_n(j) = \mr{Emb}^\mr{rect}(\sqcup_j I^n,I^n)$ to $e \times \mr{id}_{\bR^k}$ and composing rectilinear embeddings. The left $\cC_n$-module structure preserves $\partial F_{n,k}$ (here it is important that $\cC_n$ is non-unitary so no cubes are forgotten), so that $\partial F_{n,k}$ is also a left $\cC_n$-module. This implies that the quotient $F_{n,k}/\partial F_{n,k}$ inherits a left $\cC_n$-module structure, that is, the maps (\ref{eqn:fnk-left-cn-module-structure}) descend to maps
\[\mathcal{C}_n(j)_+ \wedge \frac{F_{n,k}(i_1)}{\partial F_{n,k}(i_1)} \wedge \cdots \wedge \frac{F_{n,k}(i_j)}{\partial F_{n,k}(i_j)} \lra \frac{F_{n,k}(i_1+\cdots+i_j)}{\partial F_{n,k}(i_1+\cdots+i_j)}.\]

Similarly, the symmetric sequence $F_{n,k}$ is a right module over $\cC_{n+k}$. In this case the structure maps 
\begin{equation}\label{eqn:fnk-right-cnk-module-structure} F_{n,k}(j) \times \mathcal{C}_{n+k}(i_1) \times \cdots \times \mathcal{C}_{n+k}(i_j) \lra F_{n,k}(i_1 + \cdots + i_j)\end{equation}
for $i_l \geq 1$ for each $l$ are given by composition of rectilinear embeddings of $(n+k)$-cubes. The subsequence $\partial F_{n,k}$ is preserved by this (here it is again important that $\cC_n$ is non-unitary), so both $\partial F_{n,k}$ and $F_{n,k}/\partial F_{n,k}$ inherit right $\cC_{n+k}$-module structures. That is, the maps (\ref{eqn:fnk-right-cnk-module-structure}) descend to maps
\[\frac{F_{n,k}(j)}{\partial F_{n,k}(j)} \wedge \mathcal{C}_{n+k}(i_1)_+ \wedge \cdots \wedge \mathcal{C}_{n+k}(i_j)_+ \lra \frac{F_{n,k}(i_1 + \cdots + i_j)}{\partial F_{n,k}(i_1 + \cdots + i_j)}.\]

\newglossaryentry{barfnk}{%
	name={\ensuremath{\bar{F}_{n,k}}},
	description={Functor associated to symmetric sequence $F_{n,k}/\partial F_{n,k}$},
	type=symbols
}
\begin{definition}\label{def:functor-fnk} 
Define the functor $\bar{F}_{n,k}\colon \sfC \to \sfC_*$ by
\[\bar{F}_{n,k} \colon X \longmapsto \bigsqcup_{i \geq 1} \frac{F_{n,k}(i)}{\partial F_{n,k}(i)} \rtimes_{\fS_i} X^{\otimes i},\]
where $\rtimes$ denotes smash product with the associated pointed object, as in Section~\ref{sec:point-categ-point}, and the subscript ${\fS_i}$ here denotes the quotient of that smash product by the evident symmetric group action.
\end{definition}

As explained above, this is both a left $E_n$-module functor and right $E_{n+k}$-module functor. Hence it naturally lifts to a functor $\bar{\gF}_{n,k} \colon \sfC \to \cat{Alg}_{E_n}(\sfC_*)$. In the following lemma, we construct a ``scanning map.'' In its definition, one has a choice whether to translate by $v$ or $-v$.  We prefer the latter.

\begin{lemma}\label{lem.fnksigmacomparison} 
There is a natural transformation 
	\[\varphi_{(-)} \colon \gE_n(S^k \rtimes -) \Longrightarrow \bar{\gF}_{n,k}(-)\]
of functors $\sfC \to \cat{Alg}_{E_n}(\sfC_*)$, which is a weak equivalence on cofibrant objects.
\end{lemma}

\begin{proof}
Consider the map $\phi \colon \bR^k \to F_{n,k}(1)$ sending $v \in \bR^k$ to the translation of the unit cube $I^{n+k} \subset I^n \times \bR^k$ by $-v$ in $\bR^k$:
\begin{align*}\phi(v) \colon I^{n+k} &\lra I^n \times \bR^k \\
	x &\longmapsto x-(0,v)\end{align*}
      where we write $(0,v) \in \bR^n \times \bR^k$. Let us temporarily write $\partial \bR^k \subset \bR^k$ for the subspace consisting of those points for which at least one of the $\bR$-coordinates lies outside $(-1,1)$, and identify $S^k$ with the quotient $\bR^k/\partial \bR^k$.  Then $\phi(v) \in \partial F_{n,k}(1)$ for $v \in \partial \bR^k$, so there is an induced map
\[\varphi \colon S^k = \frac{\bR^k}{\partial \bR^k} \lra \frac{F_{n,k}(1)}{\partial F_{n,k}(1)}\]
and hence $\varphi \rtimes \mathrm{id}_X \colon S^k \rtimes X \to \frac{F_{n,k}(1)}{\partial F_{n,k}(1)} \rtimes X\subset \bar{\gF}_{n,k}(X)$.  Since the target is an $E_n$-algebra, $\phi \rtimes \mathrm{id}_X$ extends uniquely to an $E_n$-algebra map
\begin{equation}\label{eq:fnksigmacomparison}
\varphi_X \colon \gE_n(S^k \rtimes X) \lra \bar{\gF}_{n,k}(X)
\end{equation}
which is clearly natural in $X$.

To show that this is a weak equivalence, we may forget the $E_n$-algebra structure, and then unravelling the above definition gives the following. Consider the maps 
\begin{align*}\phi_i \colon \cC_n(i) \times (\bR^k)^i &\lra F_{n,k}(i) \\
	((e_1,\ldots,e_i),(v_1,\ldots,v_i)) &\longmapsto (e_1 \times \mr{id}_{I^k} - (0,v_1),\ldots,e_i \times \mr{id}_{I^k}-(0,v_i)).\end{align*}
That is, we cross each of the embeddings of the $n$-cubes with the identity map on $I^k$, and then translate the $i$th resulting embedding of an $(n+k)$-cube by $-v_i \in \bR^k$. Let us write $\partial(\cC_n(i) \times (\bR^k)^i) \subset \cC_n(i) \times (\bR^k)^i$ for the subspace where at least one of the $\bR$-coordinates lies outside $(-1,1)$.  This subspace is sent into $\partial F_{n,k}(i)$ by $\phi_i$, which determines a map
\[\varphi_i \colon \cC_n(i)_+ \wedge (S^k)^{\wedge i} \cong \frac{\cC_n(i) \times (\bR^k)^i}{\partial(\cC_n(i) \times (\bR^k)^i)} \lra \frac{F_{n,k}(i)}{\partial F_{n,k}(i)}.\]
As $i$ varies, these form a map of symmetric sequences in $\cat{Top}_*$, and the map \eqref{eq:fnksigmacomparison} is induced by this map of symmetric sequences.

We will show that each $\varphi_i$ is a weak homotopy equivalence, and as $\fS_i$ acts freely (away from the basepoint) on the domain and codomain of $\varphi_i$ both symmetric sequences are cofibrant. It then follows from Lemma \ref{lem.symmetric-sequence-preserving} (ii) that \eqref{eq:fnksigmacomparison} is a weak equivalence whenever $X$ is cofibrant in $\sfC$.

To see that $\varphi_i$ is a weak homotopy equivalence, first observe that this map is a homeomorphism onto its image, which is the subspace $S$ of those points which may be represented by configurations of cubes in $I^n \times \bR^k$ which have edge length 1 in each of the $\bR^k$-directions and which remain disjoint when projected to $I^n$. (The map $e$ consists of $i$ embeddings of a cube, and here we identify each of these with its image for the sake of simplicity.) Let $S'$ be the larger subspace where we omit the ``edge length 1'' condition: it consists of those points which may be represented by configurations of cubes in $I^n \times \bR^k$ which remain disjoint when projected to $I^n$. We will show that the inclusions $S \hookrightarrow S' \hookrightarrow \frac{F_{n,k}(i)}{\partial F_{n,k}(i)}$ are weak homotopy equivalences.

For the inclusion $S \to S'$ we obtain a homotopy inverse by a simple scaling of the edge lengths of representative cubes in the $\bR^k$-direction, as follows. First consider the homotopy $\sigma_t \colon {F_{n,k}(i)} \to {F_{n,k}(i)}$, $t \in [0, 1]$, sending a tuple of disjoint rectilinear cubes $e_1, e_2, \ldots, e_i \colon I^{n+k} \to I^n \times \bR^k$ to the tuple $\sigma_t(e_1), \sigma_t(e_2), \ldots, \sigma_t(e_i)$, given by
\[\sigma_t(e_j) \coloneqq \big(\mathrm{id}_{I^n} \times (1+ t \cdot \max(0,\tfrac{1}{\epsilon}-1))\mathrm{id}_{\bR^k}\big) \circ e_j\]
where $\epsilon$ is the minimum of the edge lengths of the $e_j$ in the $\bR^k$-direction. This homotopy preserves the subspace $\partial {F_{n,k}(i)}$, so induces a homotopy of the same name on $\frac{F_{n,k}(i)}{\partial F_{n,k}(i)}$. 

As $\sigma_t$ does not change the projections of the cubes to $I^n$, this homotopy preserves the subspace $S'$: it gives a deformation retraction to the subspace $S'' \subset S'$ of those cubes which have edge length $\geq 1$ in each of the $\bR^k$-directions. Now we define a deformation retraction from $S''$ to its subspace $S$ by shrinking each cube linearly in each of the $\bR^k$-directions, fixing their centers, until they have edge length precisely 1 in each of these directions. This is well-defined, because shrinking cubes about their centers preserves disjointness, and if some cube $e_j$ has image outside of $I^n \times I^k$ then it still does after shrinking it about its center.

For the inclusion $S' \hookrightarrow \frac{F_{n,k}(i)}{\partial F_{n,k}(i)}$,
first consider the 1-parameter family of self-maps $\rho_t \colon {F_{n,k}(i)} \to {F_{n,k}(i)}$, $t \in [0, \infty)$, sending a tuple of disjoint rectilinear cubes $e_1, e_2, \ldots, e_i \colon  I^{n+k} \to I^n \times \bR^k$ to the tuple $\rho_t(e_1), \rho_t(e_2), \ldots, \rho_t(e_i)$ given by translating in the $\bR^k$-direction via
\[\rho_t(e_j)(x_1, \ldots, x_{n+k}) \coloneqq e_j(x_1, \ldots, x_{n+k}) + t \cdot \mathrm{proj}_{\bR^k}(e_j(\tfrac{1}{2}, \ldots, \tfrac{1}{2})).\]
This is again a rectilinear embedding, of the same edge lengths as $e_j$, and the $\rho_t(e_j)$ are disjoint from each other: for each pair of cubes the absolute value of the difference of the $\ell$th coordinates of their centers is non-decreasing, but their size remains the same. This 1-parameter family of maps has the following crucial properties: 
\begin{enumerate}[(i)]
\item if $e_j$ has image outside of $I^n \times I^k$ then so does $\rho_t(e_j)$ for all $t \geq 0$, and 
\item if $\mathrm{proj}_{\bR^k}(e_j(\tfrac{1}{2}, \ldots, \tfrac{1}{2})) \neq 0$ then $\rho_t(e_j)$ has image outside of $I^n \times I^k$ for all $t \gg 0$.
\end{enumerate} 
By property (i) $\rho_t$ descends to a 1-parameter family of self-maps of $\frac{F_{n,k}(i)}{\partial F_{n,k}(i)}$. This 1-parameter family preserves the subspace $S'$, as it does not change the projections of cubes to $I^n$. If an equivalence class $[e_1, e_2, \ldots, e_i] \in \frac{F_{n,k}(i)}{\partial F_{n,k}(i)}$ is not in $S'$ then some pair of cubes $\{e_j, e_\ell\}$ do not have disjoint projections to $I^n$, and so, as these cubes are disjoint in $I^n \times \bR^k$, $\mathrm{proj}_{\bR^k}(e_j(\tfrac{1}{2}, \ldots, \tfrac{1}{2}))$ and $\mathrm{proj}_{\bR^k}(e_\ell(\tfrac{1}{2}, \ldots, \tfrac{1}{2}))$ cannot both be 0, and hence by property (ii) either $\rho_t(e_j)$ or $\rho_t(e_\ell)$ lies outside of $I^n \times I^k$ for all $t \gg 0$. But then $\rho_t([e_1, e_2, \ldots, e_i])$ is the basepoint for all $t \gg 0$, so in particular lies in $S'$. As the necessary $t$'s can be chosen continuously, for any compact subset $K$ of $\frac{F_{n,k}(i)}{\partial F_{n,k}(i)}$ there is a $t$ such that $\rho_t(K) \subset S'$, and hence $S' \hookrightarrow \frac{F_{n,k}(i)}{\partial F_{n,k}(i)}$ is a weak homotopy equivalence.
\end{proof}

We now wish to relate $\bar{F}_{n,k}(X)$ to $\tilde{B}^{E_k}(\gE_{n+k}(X))$. In order to do this we will construct a $k$-fold semi-simplicial resolution of the functor $\bar{F}_{n,k}$. This $k$-fold semi-simplicial resolution contains the additional data of grids of hyperplanes between the cubes.

\begin{definition}Let $F_{n,k}(i)_{p_1, \ldots, p_k} \subset \cP(p_1,\ldots,p_k) \times F_{n,k}(i)$ be the space of pairs $(\{t^j_i\},e)$ of an element $\{t^j_i\}$ of $\cP(p_1,\ldots,p_k)$ consisting of
	\[0 < t^j_0 < t^j_1 < \cdots < t^j_{p_j} < 1\]
for $j=1,2,\ldots, k$, and an $e \in F_{n,k}(i)$, such that each hyperplane $I^n \times \bR^{j-1} \times \{t^j_{i}\} \times \bR^{k-j}$ is disjoint from the interior of all cubes of $e$ (see Figure \ref{fig:fnk-res-example} for an example).

From the $k$-fold semi-simplicial structure of $\cP(\bullet,\ldots,\bullet)$, this inherits the structure of a $k$-fold semi-simplicial $k$-symmetric sequence $F_{n,k}(i)_{\bullet,\ldots,\bullet}$ in $\cat{Top}$ augmented over $F_{n,k}(i)$. In particular, the $i$th face map in the $j$th simplicial direction $d^j_i$ forgets $t^j_i$ and the augmentation $\epsilon$ sends $(\{t^j_i\},e)$ to $e$.\end{definition}

\begin{figure}[h]
	\centering
	\begin{tikzpicture}
	\draw[dashed] (0,-1) -- (0,6);
	\draw[dotted] (0,0) -- (5,0);
	\draw[dotted] (0,5) -- (5,5);
	\draw (0,0) -- (0,5);
	\draw[dashed] (5,-1) -- (5,6);
	\draw (5,0) -- (5,5);
	
	\draw[fill=Mahogany!10!white] (3,4) rectangle (4,5.5);
	\draw[fill=Mahogany!10!white] (0.2,3) rectangle (0.8,3.8);
	\draw[fill=Mahogany!10!white] (2.5,0.2) rectangle (3.5,1.8);
	\draw[fill=Mahogany!10!white] (.5,-1) rectangle (1.5,-0.2);
	\draw[fill=Mahogany!10!white] (3,2.8) rectangle (3.8,3.8);
	
	\node at (0.5,3.4) {$e_1$};
	\node at (3.5,4.75) {$e_2$};
	\node at (3,1) {$e_3$};
	\node at (1,-0.6) {$e_4$};
	\node at (3.4,3.3) {$e_5$};
	
	\node[left] at (0,0) {0};	
	\draw (0,1.8) -- (5,1.8);
	\node[left] at (0,1.8) {$t_0^1$};
	\draw (0,2.6) -- (5,2.6);
	\node[left] at (0,2.6) {$t_1^1$};
	\draw (0,3.9) -- (5,3.9);
	\node[left] at (0,3.9) {$t_2^1$};
	\node[left] at (0,5) {1};	
	
	\node[right] at (5.1,2.5) {$\in  F_{1,1}(4)_{2}$};
	\end{tikzpicture}
	\caption{An example of an element of $F_{1,1}(4)_{2}$ which also lies in $\partial F_{1,1}(3)_{2}$, as $e_4$ lies outside $I \times I$.}
	\label{fig:fnk-res-example}
\end{figure}

As in Definition \ref{def:fnk-and-friends}, we let the $k$-fold semi-simplicial $k$-symmetric sequence $\partial F_{n,k}(i)_{\bullet, \ldots, \bullet}$ in $\cat{Top}$ augmented over $\partial F_{n,k}(i)$ be the sub-object consisting of $(\{t^j_i\},e)$ such that $e \in \partial F_{n,k}(i)$. We shall use the following convenient lemma, which uses the quotient map $q \colon \bigsqcup_{n \geq 0} \Delta^n \times X_n \to ||X_\bullet||$.

\begin{lemma}\label{lem:fat-gr-homotopy-to-compact-image} Let $X_\bullet$ be a semi-simplicial space. Then each map $f \colon S^i \to ||X_\bullet||$ is homotopic to a map $\tilde{f}$ such that there exist compact subsets $K_j \subset X_j$ for $0 \leq j \leq i$, so that $\tilde{f}$ has image in $q\left(\bigcup_{0 \leq j \leq i} \Delta^j \times K_j\right) \subset ||X_\bullet||$.\end{lemma}

\begin{proof}The counit natural weak equivalences $\epsilon_n \colon \gr{\mr{Sing}(X_n)} \to X_n$ give a levelwise weak equivalence of semi-simplicial spaces 
	\[\epsilon_\bullet \colon \gr{\mr{Sing}(X_\bullet)} \lra X_\bullet,\]
which induces a weak equivalence upon thick geometric realization. Thus up to homotopy we may lift $f$ to a map $f' \colon S^i \to \fgr{\gr{\mr{Sing}(X_\bullet)}}$. The latter is homeomorphic to the thin geometric realization of the diagonal of the bisimplicial set 
\[Y_{\bullet,\bullet} \colon [p,q] \longmapsto \bigsqcup_{[p] \twoheadrightarrow [r]} \mr{Sing}_q(X_r).\]

Hence, by the simplicial approximation theorem (e.g. Corollary 4.8 of 
\cite{JardineApprox}), there exists a simplicial triangulation $L_\bullet$ of $S^i$ so that $f'$  is homotopic to $\gr{\tilde{f}_\bullet}$ for a simplicial map $\tilde{f}_\bullet \colon L_\bullet \to \mr{diag}(Y_{\bullet,\bullet})$. Under $\tilde{f}_\bullet$ each non-degenerate $q$-simplex $\sigma$ maps to some continuous map $\Delta^q \to X^r$ for $r \leq q \leq i$. As the simplicial set $L_\bullet$ has finitely many non-degenerate simplices, necessarily of dimension $j \leq i$, for each $j \leq i$ we obtain a finite collection of continuous maps $\Delta^q \to X_j$ so that we may take $K_j$ to be the image in $X_j$ of these maps. This is a finite union of compact subsets and hence compact.
\end{proof}

\begin{lemma}\label{lem:CubesResolution}
	The augmentations induce maps
	\[\fgr{ F_{n,k}(i)_{\bullet, \ldots, \bullet} } \lra F_{n,k}(i) \quad \text{ and } \quad \fgr{ \partial F_{n,k}(i)_{\bullet, \ldots, \bullet} } \lra \partial F_{n,k}(i)\]
	which are weak equivalences of $k$-symmetric sequences. Thus, taking the pointed geometric realization, the map
	\[\fgr{ F_{n,k}(i)_{\bullet, \ldots, \bullet}/\partial F_{n,k}(i)_{\bullet, \ldots, \bullet} }_* \lra F_{n,k}(i)/\partial F_{n,k}(i)\]
	is a weak equivalence of $k$-symmetric sequences in $\cat{Top}_*$.
\end{lemma}

\begin{proof}We start with the easier proof that $\fgr{ F_{n,k}(i)_{\bullet, \ldots, \bullet} } \to F_{n,k}(i)$ is a weak equivalence. We define the space $C_{n,k}(i)$ of an ordered configuration of $i$ disjoint points in $I^n \times \bR^k$ as
	\[C_{n,k}(i) \coloneqq \mr{Emb}(\sqcup_i \ast,I^n \times \bR^k).\]
	Let $C_{n,k}(i)_{p_1,\ldots,p_k} \subset \cP(p_1,\ldots,p_k) \times C_{n,k}(i)$ be the space of pairs $(\{t^j_{i'}\},x)$ of a grid in $\cP(p_1,\ldots,p_k)$ and a configuration $x \in C_{n,k}(i)$ such that each hyperplane $I^n \times \bR^{j-1} \times \{t^j_{i'}\} \times \bR^{k-j}$ is disjoint from all points in the configuration $x$. We may assemble these into a $k$-fold semi-simplicial $k$-symmetric sequence $C_{n,k}(i)_{\bullet,\ldots,\bullet}$ in $\cat{Top}$ augmented over $C_{n,k}(i)$.
	
	Evaluating at the centers of cubes gives the maps of augmented $k$-fold semi-simplicial spaces
	\[\begin{tikzcd}F_{n,k}(i)_{\bullet,\ldots,\bullet} \dar \rar{\simeq} & C_{n,k}(i)_{\bullet,\ldots,\bullet} \dar \\
	F_{n,k}(i) \rar{\simeq} & C_{n,k}(i)\end{tikzcd}\] 
	which is easily seen to be a levelwise weak equivalence. Hence to prove that $\fgr{F_{n,k}(i)_{\bullet, \ldots, \bullet} } \to F_{n,k}(i)$ is a weak equivalence, it suffices to prove that the map $\pi \colon \fgr{ C_{n,k}(i)_{\bullet, \ldots, \bullet} } \to C_{n,k}(i)$ induced by the augmentation of the augmented $k$-fold semi-simplicial object is a weak equivalence.
	
We will prove this using the notion of a \emph{Serre microfibration}, see \cite[p.\ 190]{WeissClassify},\index{Serre microfibration} and in particular the result that a Serre microfibration with weakly contractible fibers is in fact a Serre fibration \cite[Lemma 2.2]{WeissClassify}, and hence a weak equivalence. Since a hyperplane disjoint from a finite configuration of points $x$ stays disjoint under a small perturbation of $x$, the map $\pi$ is a Serre microfibration. 

The fiber of $\pi$ over $x$ is given by thick geometric realization of the $k$-fold semi-simplicial space $C_{n,k}(x)_{\bullet,\ldots,\bullet}$ with $(p_1,\ldots,p_k)$-simplices $C_{n,k}(x)_{p_1,\ldots,p_k}$ given by the subspace of $\cP(p_1,\ldots,p_k)$ consisting of grids $\{t_{i'}^j\}$ such that each hyperplane $I^n \times \bR^{j-1} \times \{t^j_{i'}\} \times \bR^{k-j}$ is disjoint from all points in the configuration $x$. Equivalently, grids which are disjoint from the projection $\mr{proj}_k(x)$ of $x$ to $\bR^k$. The conditions on hyperplanes in each coordinate are independent, so we recognize this as a $k$-fold product of semi-simplicial spaces $\smash{\prod_{j=1}^k X_\bullet(c_1^j, \ldots, c_i^j)}$ where $c_1^j \leq \cdots \leq c_i^j$ are elements in $\bR$, and $X_\bullet(c_1^j, \ldots, c_i^j)$ is the nerve of the topological poset of real numbers in $(0,1)$ distinct from the $c_{i'}^j$.

It is therefore enough to show that $\fgr{X_\bullet(c_1, \ldots, c_i)}$ is weakly contractible for any real numbers $c_1 \leq \cdots \leq c_i$, so let $f \colon S^m \to \fgr{X_\bullet(c_1, \ldots, c_i)}$ be a continuous map, which we shall show is homotopic to a constant map. By Lemma \ref{lem:fat-gr-homotopy-to-compact-image}, we may assume that $f$ has image in $q\left(\bigsqcup_{0 \leq j \leq i'} \Delta^j \times K_j\right)$ for $K_j \subset X_j(c_1, \ldots, c_i)$ compact. Let $c$ be the smallest strictly positive $c_i$, or $1$ if there is none. By compactness there exists an $0 < \epsilon < c$ such that $\epsilon < t_0$ for all $\{t_0 < t_1 < \cdots < t_j\} \in K_j$ and all $j \leq i'$. Then the inclusion 
\[q\left(\bigsqcup_{0 \leq j \leq i'} \Delta^j \times K_j\right) \subset \fgr{X_\bullet(c_1, \ldots, c_i)}\]
is nullhomotopic, as it extends over the cone to the vertex $\epsilon$.
\vspace{1ex}

	For $\fgr{ \partial F_{n,k}(i)_{\bullet, \ldots, \bullet} } \to \partial F_{n,k}(i)$, consider instead the subspace $\partial^\ast F_{n,k}(i) \subset  F_{n,k}(i)$ such that at least one cube has center outside of $I^n \times I^k$. The inclusion $\partial F_{n,k}(i) \hookrightarrow \partial^\ast F_{n,k}(i)$ is a weak equivalence. Similarly, we may define $\partial^\ast F_{n,k}(i)_{p_1,\ldots,p_k}$ as the subspace of $F_{n,k}(i)_{p_1,\ldots,p_k}$ where at least one cube has center outside of $I^n \times I^k$, and the inclusion $\partial F_{n,k}(i)_{\bullet,\ldots,\bullet} \hookrightarrow \partial^\ast F_{n,k}(i)_{\bullet,\ldots,\bullet}$ is a levelwise weak equivalence. Hence it suffices to prove that \[\fgr{ \partial^\ast F_{n,k}(i)_{\bullet, \ldots, \bullet} } \lra \partial^\ast F_{n,k}(i)\]
	is a weak equivalence. This follows by specializing the previous proof to these subspaces.
\end{proof}

We define functors $(\bar{F}_{n,k})_{p_1, \ldots, p_k}$ analogously to $\bar{F}_{n,k}$, as
\[(\bar{F}_{n,k})_{p_1, \ldots, p_k}(X) \colon X \longmapsto \bigsqcup_{i \geq 1} \frac{F_{n,k}(i)_{p_1, \ldots, p_k}}{\partial F_{n,k}(i)_{p_1, \ldots, p_k}} \rtimes_{\fS_i} X^{\otimes i}.\]
This is a right $E_{n+k}$-functor for the same reason that $\bar{F}_{n,k}$ is. However, it is \emph{not} a left $E_n$-functor, as attempting to use elements of $\cC_n$ to combine different collections of cubes with grids might result in the grids of one collection intersecting the cubes of the other collection. 
These assemble into a $k$-fold augmented semi-simplicial object $(\bar{F}_{n,k})_{\bullet, \ldots, \bullet} \to \bar{F}_{n,k}$.

\begin{lemma}\label{lem:fnk-resolution}The map $\fgr{ (\bar{F}_{n,k})_{\bullet, \ldots, \bullet}(X)} \to \bar{F}_{n,k}(X)$ is a weak equivalence for $X \in \sfC$ cofibrant.\end{lemma}
\begin{proof}
There is an isomorphism
\[\fgr{ (\bar{F}_{n,k})_{\bullet, \ldots, \bullet}(X)} \cong \bigsqcup_{i \geq 1} \fgr{{F_{n,k}(i)_{\bullet, \ldots, \bullet}}/{\partial F_{n,k}(i)_{\bullet, \ldots, \bullet}}}_* \rtimes_{\fS_i} X^{\otimes i},\]
so, as $X$ is cofibrant, by Lemma \ref{lem.symmetric-sequence-preserving} (ii) it is enough to show that the augmentation
\[\fgr{{F_{n,k}(i)_{\bullet, \ldots, \bullet}}/{\partial F_{n,k}(i)_{\bullet, \ldots, \bullet}}}_* \lra {F_{n,k}(i)}/{\partial F_{n,k}(i)}\]
is a weak equivalence of cofibrant symmetric sequences. It is a weak equivalence by Lemma \ref{lem:CubesResolution}, and the $\fS_i$-action is free away from the basepoint by observation.
\end{proof}

In $F_{n,k}(i)_{\bullet,\ldots,\bullet}$ we have grids of hyperplanes as in the bar construction which avoid the interior of the $(n+k)$-cubes. If we take the quotient by $\partial F_{n,k}(i)_{\bullet,\ldots,\bullet}$, any collection of $(n+k)$-cubes with some cube lying outside $I^n \times I^k$ is identified with the basepoint. However, in the reduced $k$-fold bar construction as in Definition \ref{def:reduced-kfold-bar} and Lemma \ref{lem:reduced-kfold-bar} a collection should already be collapsed to the basepoint when some $(n+k)$-cube is in the outer parts of the grid. To remedy this discrepancy, we make the following definition: 

\newglossaryentry{partialfnkcirc}{%
	name={\ensuremath{\partial^\circ F_{n,k}}},
	description={Variation of $\partial F_{n,k}$},
	type=symbols
}
\begin{definition}Let $\gls{partialfnkcirc}(i)_{p_1, \ldots, p_k}$ be the subspace of $F_{n,k}(i)_{p_1, \ldots, p_k}$ of pairs $(\{t^j_i\},e)$ such that some cube of $e$ lies outside the interior of $I^n \times [t_0^1, t_{p_1}^1] \times \cdots \times [t_0^k, t_{p_k}^k]$. This is a collection of path components and defines a sub-object of $p$-fold semi-simplicial $k$-symmetric sequences
\[\partial^\circ F_{n,k}(i)_{\bullet, \ldots, \bullet} \subset F_{n,k}(i)_{\bullet, \ldots, \bullet}.\]
\end{definition}
The inclusions $\partial F_{n,k}(i)_{p_1, \ldots, p_k} \hookrightarrow \partial^\circ F_{n,k}(i)_{p_1, \ldots, p_k}$ induce a map of augmented $k$-fold semi-simplicial spaces 
	\[\partial F_{n,k}(i)_{\bullet, \ldots, \bullet} \lra \partial^\circ F_{n,k}(i)_{\bullet, \ldots, \bullet},\] 
which is easily seen to be a levelwise homotopy equivalence by scaling coordinates.

Using this variant, we define a $k$-fold simplicial functor
\[(\bar{F}^\circ_{n,k})_{p_1, \ldots, p_k}(X) \colon X \longmapsto \bigsqcup_{i \geq 1} \frac{F_{n,k}(i)_{p_1, \ldots, p_k}}{\partial^\circ F_{n,k}(i)_{p_1, \ldots, p_k}} \rtimes_{\fS_i} X^{\otimes i}\]
which comes with a natural transformation $(\bar{F}_{n,k})_{\bullet, \ldots, \bullet} \Rightarrow (\bar{F}^\circ_{n,k})_{\bullet, \ldots, \bullet}$ because $\partial^\circ F_{n,k}(i)_{p_1, \ldots, p_k}$ contains $\partial F_{n,k}(i)_{p_1, \ldots, p_k}$.

\begin{lemma}\label{lem:fnk-fnkcirc} The natural transformation $\fgr{(\bar{F}_{n,k})_{\bullet, \ldots, \bullet}(-)} \Rightarrow \fgr{(\bar{F}^\circ_{n,k})_{\bullet, \ldots, \bullet}(-)}$ is a weak equivalence on cofibrant objects.\end{lemma}

\begin{proof}
Let $X \in \sfC$ be cofibrant. By Lemma \ref{lem:thick-geom-rel-cofibrations} it is enough to show that each $(\bar{F}_{n,k})_{p_1, \ldots, p_k}(X) \to (\bar{F}^\circ_{n,k})_{p_1, \ldots, p_k}(X)$ is a weak equivalence between cofibrant objects of $\sfC_*$. As $X$ is cofibrant in $\sfC$, by Lemma \ref{lem.symmetric-sequence-preserving} (ii) it is enough to show that
\[\frac{F_{n,k}(i)_{p_1, \ldots, p_k}}{\partial F_{n,k}(i)_{p_1, \ldots, p_k}} \lra \frac{F_{n,k}(i)_{p_1, \ldots, p_k}}{\partial^\circ F_{n,k}(i)_{p_1, \ldots, p_k}}\]
is a weak equivalence between cofibrant symmetric sequences. The $\fS_i$-action on both spaces is free away from the basepoint, so they are cofibrant symmetric sequences; the map is a weak equivalence as $\partial F_{n,k}(i)_{p_1, \ldots, p_k} \to \partial^\circ F_{n,k}(i)_{p_1, \ldots, p_k}$ is.
\end{proof}

The following lemma connects the functor $\fgr{(\bar{F}^\circ_{n,k})_{\bullet, \ldots, \bullet}(-)}$ to the reduced $k$-fold bar construction $\tilde{B}^{E_k}$ of Definition \ref{def:reduced-kfold-bar} and Lemma \ref{lem:reduced-kfold-bar}.

\begin{lemma}\label{lem:bek-fnk-iso}
There is a natural isomorphism $\tilde{B}^{E_k}_{\bullet, \ldots, \bullet}(\gE_{n+k}(X)) \cong (\bar{F}^\circ_{n,k})_{\bullet, \ldots, \bullet}(X)$ of $k$-fold semi-simplicial objects.
\end{lemma}

\begin{proof}
There are two cases to consider. Firstly, if some $p_j$ is $0$ then we have isomorphisms
	\[\tilde{B}^{E_k}_{p_1, \ldots,p_k}(\gE_{n+k}(X)) \cong \ast \cong (\bar{F}^\circ_{n,k})_{p_1, \ldots, p_k}(X),\] 
	in the latter case because $\partial^\circ F_{n,k}(i)_{p_1, \ldots, p_k} = F_{n,k}(i)_{p_1, \ldots, p_k}$ for all $i \geq 0$.
	
The second case is when $p_j>0$ for all $j$. In this case we have that the pointed object $\tilde{B}^{E_k}_{p_1, \ldots,p_k}(\gE_{n+k}(X))$ is given by the quotient of
\[\cP(p_1,\ldots,p_k) \times \bigotimes_{q_1=1}^{p_1} \cdots \bigotimes_{q_k=1}^{p_k} \left(\bigsqcup_{i \geq 0} \cC_{n+k}(i) \times_{\fS_i} X^{\otimes i}\right)\]
by the sub-object corresponding to the terms $i=0$, given by $\cP(p_1,\ldots,p_k) \times \bigotimes_{q_1=1}^{p_1} \cdots \bigotimes_{q_k=1}^{p_k} \bunit$.

We may describe this quotient as $\cB(X)$ for a symmetric sequence $\cB \in \cat{FB}_\infty(\cat{Top}_\ast)$ applied to $X$. A non-basepoint element of $\cB(i)$ is an element in the space of grids $\cP(p_1,\ldots,p_k)$, together with for each $(q_j)_{j=1}^k \in \prod_{j=1}^k \{1,\ldots,p_j\}$ a collection $e[(q_j)_{j=1}^k]$ of $(n+k)$-cubes with interiors disjoint from the grid, and a bijection of $\{1,\ldots,i\}$ with the set of all of these $(n+k)$-cubes. 

This is isomorphic to the symmetric sequence $\cB' \in \cat{FB}_\infty(\cat{Top}_\ast)$ with a non-basepoint element of $\cB'(i)$ given by an element in the space of grids $\cP(p_1,\ldots,p_k)$ together with an element $e$ of $\cC_{n+k}(i)$ with image in $[t_0^1,t_{p_1}^1] \times \ldots [t^k_0,t^k_{p_k}]$ whose interior is disjoint from the grid. The isomorphism $\cB \to \cB'$ is given by composing with $\delta \in \cC_k(p_1\cdots p_k)$ of (\ref{eqn:bar-embedding-delta}).

The space $\cB'(i)$ may be described as adding a disjoint basepoint to those path components of $F_{n,k}(i)$ that are not in $\partial^\circ F_{n,k}(i)$. Since all the path components of $F_{n,k}(i)$ that are in $\partial^\circ F_{n,k}(i)$ are collapsed to the basepoint, we have described the desired isomorphism 
\[\tilde{B}^{E_k}_{p_1, \ldots,p_k}(\gE_{n+k}(X)) \cong (\bar{F}^\circ_{n,k})_{p_1, \ldots,p_k}(X).\]

We need to check this isomorphism is simplicial. Upon applying $d_i^j$ for $i=0$ or $i=p_j$, both sides get mapped to $\ast$. Otherwise, the face maps act identically on the grids $\{t_i^j\}$. On the rectilinear embeddings, the action on $e$ is the identity, and $\delta$ coequalizes the identity and the application of $d_i^j$ by associativity of composition of rectilinear embeddings.
\end{proof}

We may now complete the proof of Theorem \ref{thm:CalcFree}.

\begin{proof}[Proof of Theorem \ref{thm:CalcFree}] The result follows from  the following zig-zag of natural transformations, each of which has been shown to be a weak equivalence when $X$ is cofibrant:
\begin{equation}
\label{eqn:zig-zag-bek-fnk}
\begin{tikzcd}\tilde{B}^{E_k}(\gE_{n+k}(X)) \arrow["\cong","\text{Lemma \ref{lem:bek-fnk-iso}}"']{d} &[10pt] & \\
\fgr{ (\bar{F}^\circ_{n,k})_{\bullet, \ldots, \bullet}(X)} & \arrow["\simeq"',"\parbox{1 cm}{\centering \tiny Lemma \ref{lem:fnk-fnkcirc}}"]{l} \fgr{ (\bar{F}_{n,k})_{\bullet, \ldots, \bullet}(X)} \arrow["\simeq","\parbox{1 cm}{\centering \tiny Lemma \ref{lem:fnk-resolution}}"']{r} & \bar{F}_{n,k}(X)  \\
 & & E_n(S^k \ltimes X). \arrow["\simeq","\text{Lemma \ref{lem.fnksigmacomparison}}"']{u}\end{tikzcd}
\end{equation}
\end{proof}

\subsection{Proof of Theorem \ref{thm:BarHomologyIndec}}\label{sec:ProofBEqQ} We will now prove Theorem \ref{thm:BarHomologyIndec}, which describes how to compute derived indecomposables using iterated bar constructions. 

Let $S^k_{\bullet,\ldots,\bullet}$ be the $k$-fold semi-simplicial object in $\cat{sSet}_*$ given by taking the quotient simplicial set $\Delta^1/\partial \Delta^1$, taking its $k$-fold smash product considered as a $k$-fold pointed simplicial set and remembering only the $k$-fold semi-simplicial structure.

\begin{definition}\label{def:qek-bullets} We define a $k$-fold semi-simplicial object
\[Q^{E_k}_{\bullet, \ldots, \bullet}(\gR) \coloneqq S^k_{\bullet,\ldots,\bullet} \wedge Q^{E_k}(\gR).\]
\end{definition}

Let $\bar{S}^k$ denote thick geometric realization $\fgr{S^k_{\bullet,\ldots,\bullet}}$, and recall that the quotient map $\bar{S}^k \to S^k$ from the thick to the thin geometric realization is a homotopy equivalence. By construction we have that $\fgr{Q^{E_k}_{\bullet, \ldots, \bullet}(\gR)}_+ \cong \fgr{S^k_{\bullet,\ldots,\bullet}} \wedge Q^{E_k}(\gR)$, and thus the natural transformation
\[\fgr{S^k_{\bullet,\ldots,\bullet}} \wedge Q^{E_k}(\gR) \lra S^k \wedge Q^{E_k}(\gR)\] 
is a weak equivalence.

We claim that there is a map of $k$-fold semi-simplicial objects 
\[\tilde{B}^{E_k}_{\bullet,\ldots,\bullet}(\gR) \lra Q^{E_k}_{\bullet,\ldots,\bullet}(\gR)\]
induced by projection of $\cP(p_1,\ldots,p_k)$ to a point, the canonical map $\gR \to Q^{E_k}(\gR)$ if only one entry is $\gR$ and the canonical map to the terminal object otherwise. To see this is indeed semi-simplicial, recall that by Lemma \ref{lem:indec-dec-quotient} the object $Q^{\cC_k}_{\cC_k(1)}(\gR)$ can be obtained as the quotient of $U^{\cC_k}_{\cC_k(1)}(\gR)_+$ by $\Dec^{\cC_k}_{\cC_k(1)}(\gR)$, so the map $\cC_k(n) \times_{\fS_n} \gR^{\otimes n} \to Q^{E_k}(\gR)$ for $n \geq 2$ factors over the terminal object. 

Upon geometric realization of this $k$-fold semi-simplicial map we obtain a pair of natural transformations of functors $\Alg_{E_k}(\sfC) \to \sfC_\ast$
\[\begin{tikzcd}\tilde{B}^{E_k} \arrow[Rightarrow,swap]{dr}{\upsilon} \arrow[Rightarrow]{r}{\tilde{\upsilon}} & \bar{S}^k \wedge Q^{E_k} \arrow[Rightarrow]{d} \\
& S^k \wedge Q^{E_k}.\end{tikzcd}\]

\begin{lemma}\label{lem:EtaOnFrees}
If $\gR \cong \gE_k(X)$ with $X \in \sfC$ cofibrant, then $\upsilon_\gR$ is a weak equivalence.
\end{lemma}

\begin{proof}It suffices to prove that $\tilde{\upsilon} \colon \tilde{B}^{E_k} \Rightarrow \bar{S}^k \wedge Q^{E_k}$ is a natural weak equivalence.
	
Firstly, the canonical morphism $X_+ \to E_k(X)_+ \to Q^{E_k}(\gE_k(X))$ is an isomorphism in $\sfC_*$ by Corollary \ref{cor:properties-of-indecomposables}. This gives us an isomorphism $\bar{S}^k \wedge X_+ \cong \fgr{S^k_{\bullet,\ldots,\bullet} \wedge Q^{E_k}(\gE_k(X))}$, which forms the right column of (\ref{eqn:eta-on-frees}). Secondly, we obtain the middle column of (\ref{eqn:eta-on-frees}) from the natural transformations between $\tilde{B}^{E_k}(\gE_k(X))$ and $\bar{F}_{0,k}(X)$ which appear in the diagram (\ref{eqn:zig-zag-bek-fnk}) in the case $n=0$. These are weak equivalences under the assumption that $X$ is cofibrant.

\begin{equation}\label{eqn:eta-on-frees}
\begin{tikzcd}
	\bar{S}^k \wedge X_+ \ar[equals]{rr} &[-5pt] & \bar{S}^k \wedge X_+ \\
	\fgr{S^k_{\bullet, \ldots, \bullet} \wedge X_+} \arrow[equals]{u} &\fgr{\tilde{B}^{E_k}_{\bullet, \ldots, \bullet}(\gE_k(X))} \rar{\tilde{\upsilon}_{\gE_k(X)}} & \fgr{S^k_{\bullet,\ldots,\bullet} \wedge Q^{E_k}(\gE_k(X))} \arrow[equals]{u} \\
	\fgr{ (\bar{G}^\circ_{0,k})_{\bullet, \ldots, \bullet}(X)} \arrow["\fgr{j_\bullet}"',"\simeq"]{u} \rar & \fgr{(\bar{F}^\circ_{0,k})_{\bullet, \ldots, \bullet}(X)} \arrow{u}{\cong} & \\
	\fgr{ (\bar{G}_{0,k})_{\bullet, \ldots, \bullet}(X)} \rar \arrow[swap]{d}{\simeq} \uar{\simeq}& \fgr{ (\bar{F}_{0,k})_{\bullet, \ldots, \bullet}(X)} \arrow[swap]{d}{\simeq} \uar{\simeq}\\
	\bar{G}_{0,k}(X) \rar{\simeq} & \bar{F}_{0,k}(X).
\end{tikzcd}
\end{equation}

The left column of (\ref{eqn:eta-on-frees}) remains to be defined. Let $\bar{G}_{0,k}(X) \subset \bar{F}_{0,k}(X)$ be the sub-object which consists of $\leq 1$ cubes labeled by $X$, giving rise to a functor
\[\bar{G}_{0,k}(X) \colon X \longmapsto \frac{F_{0,k}(1)}{\partial F_{0,k}(1)} \rtimes X = S^k \rtimes X.\]
We may form the $k$-fold semi-simplicial objects $(\bar{G}_{0,k})_{\bullet, \ldots, \bullet}(X)$ and $(\bar{G}^\circ_{0,k})_{\bullet, \ldots, \bullet}(X)$ in analogy with those for $\bar{F}_{0,k}(X)$. As in the proof of Lemma \ref{lem:CubesResolution}, the maps between them become weak equivalences upon geometric realization. Furthermore, the proof of Lemma \ref{lem.fnksigmacomparison} shows that the inclusion $\bar{G}_{0,k}(X) \to \bar{F}_{0,k}(X)$ is a weak equivalence.

Let us now define the $k$-fold semi-simplicial map $j_\bullet$. We have that
\[(\bar{G}^\circ_{0,k})_{1, \ldots, 1}(X) \cong \mathcal{P}_k(1, \ldots, 1) \times \cC_k(1) \times X \subset \mathcal{P}_k(1, \ldots, 1) \times E_k(X)\]
and the map $j_\bullet$ is given by taking connected components of $(\bar{G}_{0,k})_{\bullet, \ldots, \bullet}$ and identifying $\pi_0((\bar{G}^\circ_{0,k})_{p_1, \ldots, p_k})$ with $S^k_{p_1,\ldots,p_k}$ by recording which of the subsets cut out by the hyperplanes of grid contains the unique cube. This is a levelwise weak equivalence.

We claim that the entire diagram (\ref{eqn:eta-on-frees}) commutes. The two bottom squares commute since they are induced by a commutative diagram of symmetric sequences. The two maps $\fgr{ (\bar{G}^\circ_{0,k})_{\bullet, \ldots, \bullet}(X)} \to \bar{S}^k \wedge X_+$ coincide, as both ways around record the label $X$ and which of the subsets cut out by the grid hyperplanes contains the unique cube. From this we conclude that $\upsilon_{\gE_k(X)}$ is a weak equivalence.
\end{proof}

We can now finish the proof of Theorem \ref{thm:BarHomologyIndec}.

\begin{proof}[Proof of Theorem \ref{thm:BarHomologyIndec}] To make use of Lemma \ref{lem:EtaOnFrees}, choose a simplicial resolution $\gR_\bullet \to \gR$ as in Section \ref{sec:FreeRes}. Such a resolution can be picked naturally in $\gR$ by taking the thick monadic bar resolution. Then we have a commutative diagram
\begin{equation}\label{eqn:zigzag-bk-qek}
\begin{aligned}
\begin{tikzcd} \tilde{B}^{E_k}(\gR) & \quad & \quad \\
\tilde{B}^{E_k}(\fgr{\gR_\bullet}_{E_k}) \rar{\upsilon_{\fgr{\gR_\bullet}}} \arrow["\epsilon"',"\simeq"]{u} \arrow[equal]{d} & S^k \wedge Q^{E_k}(\fgr{\gR_\bullet}_{E_k}) \arrow[equal]{r} \arrow[equal]{d} & S^k \wedge Q^{E_k}_\bL(\gR)\\
\fgr{\tilde{B}^{E_k}(\gR_\bullet)} \arrow["\fgr{\upsilon_{\gR_\bullet}}","\simeq"']{r} & \fgr{S^k \wedge Q^{E_k}(\gR_\bullet)},
\end{tikzcd} \end{aligned}
\end{equation}
where the right vertical equality comes from the fact that $Q^{E_k}$ commutes with geometric realization, and the left vertical equality follows from Lemma \ref{lem:monad-functors-geomrel} and commuting two geometric realizations. The map induced by $\epsilon$ is a weak equivalence by Lemma \ref{lem:thick-geom-rel-cofibrations}. The lower map is a weak equivalence by Lemma \ref{lem:EtaOnFrees}, as each $\gR_p$ is a free $E_k$-algebra. This proves Theorem \ref{thm:BarHomologyIndec}.\end{proof}

\subsection{The bar construction on maps between free algebras}\label{sec:FreeAlgMaps}
Given a map $f \colon X \to Y$ in $\cat{C}$, we obtain a map $\gE_{n+k}(f) \colon \gE_{n+k}(X) \to \gE_{n+k}(Y)$ and Theorem \ref{thm:CalcFree} identifies $\tilde{B}^{E_k}(\gE_{n+k}(f))$ with $E_n(S^k \wedge f_+) \colon E_n(S^k \wedge X_+) \to E_n(S^k \wedge Y_+)$ up to natural weak equivalence. However, a general map $F \colon \gE_{n+k}(X) \to \gE_{n+k}(Y)$ of $E_{n+k}$-algebras need not be of the form $\gE_{n+k}(f)$. We wish to describe the map $\tilde{B}^{E_k}(F)$ in terms of free $E_n$-algebras. To do so, note that $F$ is determined by a map $f \colon X \to E_{n+k}(Y)$, and may be factored as

\begin{equation*}\begin{tikzcd}
\gE_{n+k}(X) \arrow[swap]{r}{\gE_{n+k}(f)} \arrow[bend left=25]{rr}{F} &[20pt] \gE_{n+k}(E_{n+k}(Y)) \arrow[swap]{r}{\mu_Y}  &  \gE_{n+k}(Y)
\end{tikzcd}\end{equation*}
where the first map is of the form we have treated already, and the second is given by the monadic structure map. Thus it is enough to describe $\tilde{B}^{E_k}(\mu_Y)$.

\newglossaryentry{eta}{%
	name={\ensuremath{\eta_{(-)}}},
	description={Natural transformation defined using scanning},
	type=symbols
}
The description we shall give will be in terms of a natural transformation
\[\gls{eta} \colon S^k \wedge E_{n+k}(-) \Longrightarrow E_n(S^k \wedge -)\]
of functors $\cat{C}_\ast \to \cat{C}_\ast$.

\begin{remark}\label{rem:may-natural-transformation}
This is related to a construction due to May \cite[Proposition 5.4]{GILS}. May's construction is different from ours because it is built for the ``group completion'' augmentation, which does not exist in general in our setup, and even if it does may not be equal to the canonical augmentation. It is discussed in Section \ref{sec:GpCompletion} for $\sfC = \cat{Top}$.
\end{remark}

We construct $\eta_{(-)}$ in the case $n=0$, generalizing to $n > 0$ afterwards. The natural transformation $\eta_{(-)} \colon S^k \wedge E_k(-) \Rightarrow S^k \wedge -$ is induced by a map $\eta$ of symmetric sequences in $\cat{Top}_*$ given by $\fS_i$-equivariant maps
\[\eta_i \colon S^k \wedge \cC_k(i)_+ \lra (S^k)^{\wedge i}.\] 
The map $\eta_i$ is the constant map to the basepoint if $i>1$. To define $\eta_1$, suppose we are given $v = (v_1, \ldots, v_k) \in (0,1)^k \subset ((0,1)^k)^+ = S^k$ and a single cube $e \in \cC_{k}(1)$. 
Then we define
\begin{align*} \eta_1(v,e) \coloneqq \begin{cases} \ast & \text{if $v \notin \mr{im}(e)$,} \\
		e^{-1}(v) \in (0,1)^k \subset ((0,1)^k)^+ \cong S^k & \text{otherwise}.\end{cases}\end{align*}

For the generalization to $n>0$, we shall similarly record the intersection of cubes with an $n$-dimensional linear subspace, and map to the basepoint if any of these intersections is empty. The natural transformation $\eta_{(-)} \colon S^k \wedge E_{n+k} \Rightarrow E_n(S^k \wedge -)$ is induced by a map $\eta$ of symmetric sequences in $\cat{Top}_*$ given by $\fS_i$-equivariant maps
\[\eta_i \colon S^{k} \wedge \cC_{n+k}(i)_+ \lra \cC_n(i)_+ \wedge (S^k)^{\wedge i}.\]
To define $\eta_i$ on $v \in (0,1)^k \subset ((0,1)^k)^+ \cong S^k$ and a collection of cubes $e = (e_1, \ldots, e_i) \in \cC_{k+n}(i)$, we do the following. We write each $e_j$ as a product $e_j^n \times e_j^k$ of a rectilinear embedding $e_j^n: I^n \hookrightarrow \bR^n$ and a rectilinear embedding $e_j^k : I^k \hookrightarrow \bR^k$, that is, an $n$-cube and a $k$-cube. Firstly, if $v$ is not in $\mr{im}(e_j^k)$ for all $1 \leq j \leq i$ then $\eta_i(v,e) = \ast$.  Otherwise, we let $\eta_i(v,e)$ be the point of $\cC_n(i)_+ \wedge (S^k)^{\wedge i}$ given by
\begin{equation*}
\eta_i(v,e) \coloneqq \left((e_1^n,\ldots,e_i^n),((e_1^k)^{-1}(v),\ldots,(e_i^k)^{-1}(v))\right),
\end{equation*}
noting that the $n$-cubes $e_j^n$ for $1 \leq j \leq i$ are disjoint, because the $(n+k)$-cubes $e_j$ are disjoint and their projections to $(0,1)^k$ all contain $v$. For $n=0$ this definition agrees with the one given above. See Figure \ref{fig:phii} for an example.

\begin{figure}[h]
	\centering
	\begin{tikzpicture}
	
	\draw[fill=Mahogany!10!white] (3,2) rectangle (4,4);
	\draw[fill=Mahogany!10!white] (0.5,1) rectangle (2,3);

	\draw (0,0) rectangle (5,5);
	
	\node at (1.25,2) {$e_1$};
	\node at (3.5,3) {$e_2$};
	
	\draw[thick] (0,2.5) node[anchor=north west] {$v$} -- ++ (5.0,0);
	
	\draw[|->] (5.25,2.5) -- ++(0.5,0);
	
	\draw (6.5,2.75) -- (11.5,2.75);
	\draw[Mahogany!40!white,line width=1.6pt] (7,2.75) -- node[black,anchor=south] {$e_1^n$} ++ (1.5,0);
	\draw[Mahogany!40!white,line width=1.6pt] (9.5,2.75) -- node[black,anchor=south] {$e_2^n$} ++ (1,0);
	\node at (7.65,2.3) {\tiny $(e_1^k)^{-1}(v) \in (0,1)$};
	\node at (10.2,2.3) {\tiny $(e_2^k)^{-1}(v) \in (0,1)$};
	
	\draw[decorate,decoration={brace,amplitude=5pt}] (6.25,2) -- (6.25,3);
	\draw[decorate,decoration={brace,amplitude=5pt,mirror}] (11.75,2) -- (11.75,3);
	
	\end{tikzpicture}
	\caption{The map $\eta_i$ for $k = n = 1$ and $i=2$ assigns to $(v,(e_1,e_2))$ the pair of $1$-cubes $(e_1^n,e_2^n)$ obtained by projection to the $x$-axis, and elements $(e_1^k)^{-1}(v)$ and $(e_2^k)^{-1}(v)$ in $(0,1) \subset S^1$ obtained by recording the intersections of either cube with $v$. If $e_2$ is moved slightly upwards (or $v$ is decreased) so that the line at height $v$ no longer intersects the image of $e_2$, then we map to the basepoint $\ast$.}
	\label{fig:phii}
\end{figure}

These maps of pointed spaces form a map of symmetric sequences in $\cat{Top}_*$
\[\eta \colon S^k \wedge (\cC_{n+k})_+ \lra (\cC_n)_+ \wedge (S^k)^{\wedge -}\]
which induces a natural transformation given by
\begin{equation}\label{eqn:eta-commute}\begin{tikzcd} (S^k \wedge (\cC_{n+k})_+)(Y) \rar{\eta_Y} \dar{\cong} &[38pt] ((\cC_n)_+ \wedge (S^k)^{\wedge -})(Y) \dar{\cong} \\
\bigsqcup_{i \geq 1} S^k \wedge \cC_n(i)_+ \wedge_{\fS_i} Y^{\otimes i} \rar[swap]{ \bigsqcup_i \eta_i \wedge_{\fS_i} Y^{\otimes i}} & \bigsqcup_{i \geq 1} \cC_n(i)_+ \wedge_{\fS_i} (S^k \wedge Y)^{\otimes i}.\end{tikzcd}\end{equation}
On the bottom line we have implicitly used the symmetric monoidal structure to $\fS_i$-equivariantly write $(S^k \wedge Y)^{\otimes i} \cong (S^k)^{\wedge i} \wedge Y^{\otimes i}$.

\begin{theorem}
The map $\tilde{B}^{E_k}(\mu_Y)$ is weakly equivalent to the map
\[E_n(\eta_Y) \colon E_n(S^k \wedge E_{n+k}(Y)) \lra E_n(S^k \wedge Y)\]
induced by $\eta_Y \colon S^k \wedge E_{n+k}(Y) \to E_n (S^k \wedge Y)$.
\end{theorem}
\begin{proof}
Consider the diagram
\begin{equation}\label{eqn:tildebek-muy}
\begin{tikzcd}
\tilde{B}^{E_k}(\gE_{n+k}(E_{n+k}( Y))) \arrow{r}{\tilde{B}^{E_k}(\mu_Y)}&[10pt] \tilde{B}^{E_k}(\gE_{n+k}( Y))\\
\fgr{(\bar{F}_{n,k})_{\bullet, \ldots, \bullet}(E_{n+k}(Y))} \dar[swap]{\simeq} \uar{\simeq} \arrow{r} & \fgr{(\bar{F}_{n,k})_{\bullet, \ldots, \bullet}(Y)} \dar{\simeq} \uar[swap]{\simeq}\\
\bar{F}_{n,k}(E_{n+k}(Y)) \arrow{r}& \bar{F}_{n,k}(Y)\\
E_n(S^k \wedge E_{n+k}(Y)) \uar{\simeq}[swap]{\varphi_{E_{n+k}(Y)}} \rar & E_n(S^k \wedge Y), \uar{\varphi_Y}[swap]{\simeq}
\end{tikzcd}
\end{equation}
where the top horizontal map is given by functoriality of $\tilde{B}^{E_k}$ with respect to maps of $E_k$-algebras, the middle two horizontal maps are given by the right $\cC_{n+k}$-module structures on the functors $\bar{F}_{n,k}$ and $(\bar{F}_{n,k})_{p_1, \ldots, p_k}$, and the lower map is essentially the $E_n$-map in the statement of the theorem, induced by $\eta_Y \colon S^k \wedge E_{n+k}(Y) \to E_n (S^k \wedge Y)$.  The vertical maps are weak equivalences by the proof of Theorem \ref{thm:CalcFree}. In particular, the natural transformation $\varphi_{(-)}$, which was defined in Lemma \ref{lem.fnksigmacomparison}, is a weak equivalence by that lemma.

The reason that the bottom map of \eqref{eqn:tildebek-muy} is not equal to $E_n(\eta_Y)$ is that the bottom square involves two different identifications of $S^k$: in the construction of $\eta_Y$ it is given by $((0,1)^k)^+$ and in the construction of $\varphi_Y$ we identified this with $\bR^k/\partial \bR^k$ with $\partial \bR^k$ the complement of $(-1,1)^k$. As it occurs more often, we opt to write the former as $S^k$ in this proof. Using the homomorphism $\rho \colon \bR^k/\partial \bR^k \to S^k$ given in each coordinate by $x \mapsto \tfrac{x+1}{2}$, the bottom square is given by
\[\begin{tikzcd} \bar{F}_{n,k}(E_{n+k}(Y)) \rar &[20pt] \bar{F}_{n,k}(Y)\\
	E_n(\bR^k/\partial \bR^k \rtimes E_{n+k}(Y)) \uar{\varphi_{E_{n+k}(Y)}}[swap]{\simeq} \dar{\cong}[swap]{E_n(\rho \rtimes \id)} & 	\uar{\simeq}[swap]{\varphi_Y} E_n(\bR^k/\partial \bR^k  \rtimes Y) \\
	E_n(S^k \rtimes E_{n+k}(Y)) \rar{E_n(\eta_Y)} &  E_n(S^k \rtimes Y) \uar{\cong}[swap]{E_n(\rho^{-1} \rtimes \id)}. \end{tikzcd}\]

The top two squares of (\ref{eqn:tildebek-muy}) commute, as the horizontal maps are given by right $\cC_{n+k}$-module structures and the vertical maps are induced by natural transformation of right $E_{n+k}$-module functors. The bottom square does not commute, but we claim that it does commute up to homotopy. The proof of this occupies the remainder of this section.

\vspace{1em}

All maps involved in bottom square are $E_n$-maps, so it suffices to show that it commutes up to homotopy after restriction to $\bR^k/\partial \bR^k \rtimes E_{n+k}(Y)$. We may thus restrict our attention to the diagram of pointed spaces
\begin{equation}\label{eqn:delooping-map-bottom} \begin{tikzcd}F_{n,k}(1)/\partial F_{n,k}(1) \rtimes E_{n+k}(Y) \rar[hook] &[-10pt] \bar{F}_{n,k}(E_{n+k}(Y)) \rar &[-15pt]  \bar{F}_{n,k}(Y)\\
\bR^k/\partial \bR^k \rtimes E_{n+k}(Y) \uar{\varphi \rtimes \id} \dar{\cong}[swap]{\rho \rtimes \id}  &  & E_n(\bR^k/\partial \bR^k  \rtimes Y) \uar{\varphi_Y} \\
S^k \rtimes E_{n+k}(Y) \arrow{rr}{\eta_Y} & & E_n(S^k  \rtimes Y) \uar{\cong}[swap]{E_n(\rho^{-1} \rtimes \id)} 
\end{tikzcd}\end{equation}
The left-top composition of (\ref{eqn:delooping-map-bottom}) is induced by the map of symmetric sequences
\begin{align*}\bR^k \times \cC_{n+k}(i) &\lra F_{n,k}(i) \\
	(v,(e_1,\ldots,e_i)) &\longmapsto (e_1-(0,v),\ldots,e_i-(0,v))\end{align*}
where $(0,v) \in \bR^n \times \bR^k$. This sends $\partial \bR^k \times \cC_{n+k}(i)$ into $\partial F_{n+k}(i)$ so induces a map of symmetric sequences
\[\mr{glob} \colon \frac{\bR^k}{\partial \bR^k} \rtimes \cC_{n+k}(i) \lra \frac{F_{n,k}(i)}{\partial F_{n,k}(i)},\] 
whose notation suggests its informal description as a ``global translation''.

For the sake of writing explicit homotopies later, we let $F^\ast_{n,k}(i)$ denote the subspace of $\Emb^\mr{rect}(I^{n+k},I^n \times I^k)^i$ of $i$ cubes whose interiors are disjoint when none of the cubes lies entirely outside the interior of $I^n \times I^k$. We let $\partial F^\ast_{n,k}(i) \subset F^\ast_{n,k}(i)$ denote the subspace where at least one of the cubes lies entirely outside the interior of $I^n \times I^k$. Then the induced map
\[\frac{F_{n,k}(i)}{\partial F_{n,k}(i)} \lra \frac{F^\ast_{n,k}(i)}{\partial F^\ast_{n,k}(i)}\]
is a homeomorphism. 

Let $\vec{1} \in \bR^k$ denote the vector with all entries equal to $1$. If the bottom composition does not map the basepoint, then it translates the $j$th cube by the negative of
\[\rho^{-1}((\overline{e}_j^k)^{-1}(\rho(v))) = 2(\overline{e}_j^k)^{-1}(\tfrac{\vec{1}+v}{2})-\vec{1} = (\overline{e}_j^k)^{-1}(v)+(\overline{e}_j^k)^{-1}(\vec{1})-\vec{1},\]
where $\overline{e}_j^k$ is the unique extension of $e_j^k$ to an affine-linear map $\bR^k \to \bR^k$.  Thus the bottom-right composition of (\ref{eqn:delooping-map-bottom}) is induced by the map of symmetric sequences
\begin{align*} \bR^k \times \cC_{n+k}(i) &\lra F^\ast_{n,k}(i) \\
	(v,(e_1,\ldots,e_i)) &\longmapsto \big(e_1^n \times \id_{I^k}-(0,(\overline{e}_1^k)^{-1}(v)+(\overline{e}_1^k)^{-1}(\vec{1})-\vec{1}),\ldots\big).\end{align*}
This sends $\partial \bR^k \times \cC_{n+k}(i)$ into $\partial F^\ast_{n+k}(i)$ so induces a map of symmetric sequences
\[\mr{loc} \colon \frac{\bR^k}{\partial \bR^k} \rtimes \cC_{n+k}(i) \lra \frac{F_{n,k}(i)}{\partial F_{n,k}(i)},\] 
whose notation suggests its informal description as a ``local translation'', i.e.~each cube is translated by an amount depending on the cube in question.

It suffices to prove that the maps $\mr{glob}$ and $\mr{loc}$ are $\fS_i$-equivariantly homotopic when restricted to $\bR^k/\partial \bR^k \rtimes \cC^{=}_{n+k}(i)$, where $\cC^{=}_{n+k}(i) \subset \cC_{n+k}(i)$ denotes the subspace where all cubes have sides of the same length, because the inclusion $\cC^{=}_{n+k}(i) \hookrightarrow \cC_{n+k}(i)$ is a $\fS_i$-equivariant homotopy equivalence. We will denote the common length of all sides of all cubes in $e = (e_1,\ldots,e_i)$ by $\ell(e) \in (0,1]$. The homotopy will be a concatenation of two homotopies: starting at $\mr{glob}$, the first homotopy makes the cubes have the same edge-lengths as in $\mr{loc}$, and the second homotopy linearly makes them be translated by the same amount (see Figure \ref{fig:composition}).

\begin{figure}
	\centering
	\begin{tikzpicture}
		\begin{scope}[scale=.7]
			\draw[fill=Mahogany!10!white] (0.5,1) rectangle (2,3);
			\draw[fill=Mahogany!10!white] (3,2) rectangle (4,4);

			\draw (0,0) rectangle (5,5);
			
			\node at (1.25,2) {$e_1$};
			\node at (3.5,3) {$e_2$};
			\node at (2.5,-1) {(A)};
		\end{scope}
		
		\begin{scope}[scale=.7,xshift=6cm]
			\draw[fill=Mahogany!10!white] (0.5,-2.5) rectangle (2,2.5);
			\draw[fill=Mahogany!10!white] (3,2.5) rectangle (4,7.5);

			\draw (0,0) rectangle (5,5);
			
			\node at (1.25,0.5) {$e_1$};
			\node at (3.5,5.5) {$e_2$};
			\node at (2.5,-1) {(B)};
			
		\end{scope}
	\end{tikzpicture}
	\caption{Take $k=n=1$ and $i=2$, and consider the cubes in Figure \ref{fig:phii}. We have pictured the elements of $F^*_{n,k}(i)$ under both compositions in \eqref{eqn:delooping-map-bottom} for $v=0$, so that $\rho(v) = \tfrac{1}{2}$ as depicted in Figure \ref{fig:phii}: (A) is the top composition and (B) the bottom composition. More generally, if we let $v \in \bR^k$ run from $-1$ to $1$, on the one hand the top composition translates both $e_1$ and $e_2$ by $-v$. On the other hand the bottom composition replaces $e_1$ by $e_1^{1} \times \mr{id}_I$ and $e_2$ by $e_2^{1} \times \mr{id}_I$ and translates them by $-2(\bar{e}_1^{2})^{-1}(\tfrac{v+1}{2})+1$ and $-2(\bar{e}_1^{1})^{-1}(\tfrac{v+1}{2})+1$. In either case, both cubes start in $\partial F^*_{n,k}(i)$ above the pictured square $I^n \times I^k$, move downwards into the square at a rate linear in $v$, and end in $\partial F^*_{n,k}(i)$ below the square. The difference is only in the amount of scaling and translation in the $\bR^k$-direction.}
	\label{fig:composition}
\end{figure}

For $e=(e_1,\ldots,e_i) \in \cC^{=}_{n+k}(i)$ define a 1-parameter family of embeddings 
\[\lambda_{e,t} \coloneqq \id_{I^n} \times \frac{1}{t \ell(e)+(1-t)} \id_{\bR^k} \colon I^n \times \bR^k \hookrightarrow I^n \times \bR^k, \quad\quad t \in [0,1].\]
We start at the map $\mr{glob} \colon \bR^k/\partial \bR^k \rtimes \cC^{=}_{n+k}(i) \to F_{n,k}(i)/\partial F_{n,k}(i)$, and consider the homotopy of such maps induced by
\begin{align*} [0,1] \times \bR^k \times \cC^{=}_{n+k}(i) &\lra F_{n,k}(i) \\
	(t,v,(e_1,\ldots,e_i)) &\longmapsto (\lambda_{e,t} \circ e_1-(0,v),\ldots,\lambda_{e,t} \circ e_i-(0,v)).\end{align*}
This completes the construction of the first of the two homotopies.

The endpoint of this homotopy is induced by the map given at $(v,(e_1,\ldots,e_i))$ by
\[\big(e_1^n \times \id_{I^k}+(0,-v+\tfrac{1}{\ell(e)} e_1^k(\tfrac{1}{2},\ldots,\tfrac{1}{2})),\ldots,e_i^n \times \id_{I^k}+(0,-v+\tfrac{1}{\ell(e)} e_i^k(\tfrac{1}{2},\ldots,\tfrac{1}{2}))\big),\]
while the map $\mr{loc}$ is induced by
\[\big(e_1^n \times \id_{I^k} + (0,-(\overline{e}_1^k)^{-1}(v)-(\overline{e}_1^k)^{-1}(\vec{1})+\vec{1}),\ldots,e_i^n \times \id_{I^k} +(0,-(\overline{e}_i^k)^{-1}(v)-(\overline{e}_i^k)^{-1}(\vec{1})+\vec{1}))\big).\]
As they differ only by translation, describing a homotopy between them may be done by interpolating between the translation vectors in $\bR^k$ given by
\[\tau_j(0) \coloneqq -v+\tfrac{1}{\ell(e)} e_j^k(\tfrac{1}{2},\ldots,\tfrac{1}{2}) \qquad \text{and} \qquad
\tau_j(1) \coloneqq -(\overline{e}_j^k)^{-1}(v)-(\overline{e}_j^k)^{-1}(\vec{1})+\vec{1}.\]
We do so linearly:
\[\tau_j(s) \coloneqq(1-s)\tau_j(0)+s\tau_j(1).\]
This completes the construction of the second of the two homotopies, but it remains to check that it lies in $F^\ast_{n,k}(i)$, i.e.~the interiors of two cubes can only intersect if at least one is entirely outside the interior of $I^n \times I^k$.

To do so, let us examine two cubes: $e_j$ with image given by $\prod_{m=1}^{n+k} [a_m,b_m] \subset \bR^{n+k}$ and $e_{j'}$ with image given by $\prod_{m=1}^{n+k} [a'_m,b'_m] \subset \bR^{n+k}$. As these subsets have disjoint interiors, we have $b_m \leq a'_m$ for some $1 \leq m \leq n+k$. If $1 \leq m \leq n$ then the cubes remain disjoint under the linear interpolation of translation vectors, because the projections of the cubes to their first $n$ coordinates are unchanged. If $m=n+r$ with $1 \leq r \leq k$, then for $\mr{proj}_{r} \colon \bR^{k} \to \bR$ the projection onto the $r$th coordinate we have
\[\mr{proj}_{r}(\tau_j(0)) = -v_{r}+\tfrac{1}{2}+\tfrac{a_m}{\ell(e)}, \qquad \mr{proj}_{r}(\tau_j(1)) = -\tfrac{v_{r}}{\ell(e)}+\tfrac{a_m}{\ell(e)}-\tfrac{1}{\ell(e)}+\tfrac{a_m}{\ell(e)}+1.\]
As the same computation for $e_{j'}$ in place of $e_j$ replaces $a_m$ by $a'_m$, we conclude that \begin{align*}\mr{proj}_{r}(\tau_{j'}(0))-\mr{proj}_{r}(\tau_j(0)) &= \tfrac{a'_m-a_m}{\ell(e)} \\
	\mr{proj}_{r}(\tau_{j'}(1))-\mr{proj}_{r}(\tau_j(1)) &= 2\tfrac{a'_m-a_m}{\ell(e)}.\end{align*}
	As $a'_m \geq a_m$, the difference $\mr{proj}_{r}(\tau_{j'}(s))-\mr{proj}_{r}(\tau_j(s))$ is increasing as $s \in [0,1]$. Since the projections of the images of the cubes to the $m$th coordinate have disjoint interiors for $s=0$ and their distance in the $m$th coordinate is increasing while their sizes remain constant, they must have disjoint interiors for all $s \in [0,1]$.
\end{proof}

\subsection{Iterated indecomposables} \label{sec:iterated-indecomp}Using the discussion and results of Section \ref{sec:bar-on-free-algebras} we can describe a model for the derived $E_k$-indecomposables on an $E_{n+k}$-algebra which still has the structure of an $E_n$-algebra. While we find this to be a useful conceptual tool for thinking about derived $E_k$-indecomposables, it is not necessary for any of our applications. The reader can skip this section on a first reading.

\begin{theorem}\label{thm:iterated-decomposables}
	There is a functor $\gM_{n,k} \colon \Alg_{E_{n+k}}(\sfC) \to \Alg_{E_n}(\sfC_*)$ such that there are zig-zags of natural transformations of functors $\Alg_{E_{n+k}}(\sfC) \to\sfC_\ast$
	\[U^{E_n}\gM_{n,k}(-) \Longleftarrow \cdots \Longrightarrow S^k \wedge Q^{E_k}_\bL(-)\]
	\[S^n \wedge Q^{E_n}_\bL(\gM_{n,k}(-)) \Longleftarrow \cdots \Longrightarrow S^{n+k} \wedge Q^{E_{n+k}}_\bL(-)\]
	which are weak equivalences on $E_{n+k}$-algebras that are cofibrant in $\sfC$.
\end{theorem}

The functor $\gM_{n,k}$ is constructed using the functor $\bar{F}_{n,k}$ from Definition \ref{def:functor-fnk}, which consists of little $(n+k)$-cubes in $I^n \times \bR^k$ that can disappear at infinity. This is a right $\cC_{n+k}$-module and a left $\cC_n$-module. We define $\gM_{n,k}$ by
\[\gM_{n,k} \colon \gR \longmapsto \fgr{B_\bullet(\bar{\gF}_{n,k}, E_{n+k}, \gR)},\]
which has the structure of an $\cC_n$-algebra using the left $\cC_n$-module structure on $\bar{F}_{n,k}$.

\begin{proof}[Proof of Theorem \ref{thm:iterated-decomposables}] For the first part, we construct the following diagram, where $\simeq$ denotes that an arrow is a weak equivalence when $U^{E_{n+k}} \gR$ is cofibrant:
	\[\begin{tikzcd}\fgr{B_\bullet(\bar{F}_{n,k}, E_{n+k}, \gR)} \arrow[leftrightarrow,"\simeq","\parbox{1cm}{\centering \tiny (\ref{eqn:zig-zag-bek-fnk})}"']{d} &[10 pt] U^{E_n}\gM_{n,k}(\gR) \arrow[equals]{l} \\[10pt]
	   \fgr{B_\bullet(\tilde{B}^{E_k} \gE_{n+k}, E_{n+k}, \gR)} \dar{\simeq} & \\[10pt]
	 \tilde{B}^{E_k}(\gR) \arrow[leftrightarrow,"\simeq"',"\text{Theorem \ref{thm:BarHomologyIndec}}"]{r} & S^k \wedge Q^{E_k}_\bL(\gR).
	\end{tikzcd}\]
The top left arrow denotes (part of) the zig-zag (\ref{eqn:zig-zag-bek-fnk}) applied levelwise, which consists of weak equivalences as the functor $E_{n+k} \colon \sfC \to \sfC$ preserves cofibrant objects by Lemma \ref{lem.symmetric-sequence-preserving} (i). The bottom left arrow is obtained by applying $\tilde{B}^{E_k}(-)$ to the augmented simplicial object $\sigma_* \sigma^* B_\bullet(\gE_{n+k}, E_{n+k}, \gR) \to \gR$ and geometrically realizing. As this is a free resolution (by Lemma \ref{lem:adding-degeneracies-reedy}, because $\gR$ is cofibrant in $\sfC$) the map $\fgr{B_\bullet(\gE_{n+k}, E_{n+k}, \gR)} \to \gR$ is a weak equivalence, and so, by commuting geometric realizations, the bottom left arrow is also a weak equivalence. The bottom arrow is the zig-zag given by Theorem \ref{thm:BarHomologyIndec}.\\

For the second part, we start by noting that the first part of this theorem implies
\[S^n \wedge Q^{E_n}_\bL(\gM_{n,k}(\gR)) \simeq M_{0,n}(\gM_{n,k}(\gR)) = \fgr{B_\bullet(\bar{F}_{0,n},E_n,\gM_{n,k}(\gR))}.\]
The right hand side may be written as the thick geometric realization of the following 2-fold semi-simplicial object
\[([p],[q]) \longmapsto \bar{F}_{0,n} E_n^p \bar{F}_{n,k} E_{n+k}^q \gR.\]
Let us realize in the $p$-direction first, and therefore consider the semi-simplicial functor $B_\bullet(\bar{F}_{0,n},E_n,\bar{F}_{n,k}(-))\colon \sfC \to \cat{sC}$. This has an augmentation to $\bar{F}_{0, n+k}(-)$, defined as follows. First, there is a map of symmetric sequences $F_{0,n} \circ F_{n,k} \to F_{0,n+k}$ given by taking the product of the rectilinear embeddings $I^n \hookrightarrow \bR^n$ of the first term with $\bR^k$ and composing these with the rectilinear embeddings $I^{n+k} \hookrightarrow I^n \times \bR^k$ of the second term. This map sends $\partial F_{0,n} \circ F_{n,k}$, as well as the images of the terms 
\[F_{0,n}(r) \times \Big(F_{n,k}(s_1) \times \cdots F_{n,k}(s_{i-1}) \times \partial F_{n,k}(s_i) \times F_{n,k}(s_{i+1}) \times \cdots \times F_{n,k}(s_r)\Big)\]
for $1 \leq i \leq r$, into $\partial F_{0,n+k}$. Hence it induces a map of symmetric sequences \[\frac{{F}_{0,n}}{\partial {F}_{0,n}} \circ \frac{{F}_{n,k}}{\partial{F}_{n,k}} \lra \frac{{F}_{0, n+k}}{\partial{F}_{0, n+k}}\] and hence a natural transformation $\pi \colon \bar{F}_{0,n} \bar{F}_{n,k} \Rightarrow \bar{F}_{0,n+k}$. Since it is defined by composition, it equalizes the two maps $\bar{F}_{0,n} E_n \bar{F}_{n,k} \Rightarrow \bar{F}_{0,n} \bar{F}_{n,k}$ and thus is indeed an augmentation. It induces a map of simplicial objects
\[\fgr{[p] \mapsto B_p(\bar{F}_{0,n}, E_n, \bar{F}_{n,k} E^q_{n+k}\gR)} \lra B_q(\bar{F}_{0,n+k}, E_{n+k}, \gR)\]
and so, after geometrically realizing in the $q$-direction too, a natural map
\[\epsilon_\gR \colon \fgr{B_\bullet(\bar{F}_{0,n},E_n,\gM_{n,k}(\gR))} \lra \fgr{B_\bullet(\bar{F}_{0,n+k},E_{n+k},\gR)}.\]

By the first part of the theorem, the right-hand term has a zig-zag of natural transformations to $S^{n+k} \wedge Q^{E_{n+k}}_\bL(\gR)$ which are weak equivalences when $\gR$ is cofibrant in $\sfC$. Note that $E_{n+k}^q \gR$ is cofibrant in $\sfC$ when $\gR$ is, by Lemma \ref{lem.symmetric-sequence-preserving} (i), so by Lemma \ref{lem:thick-geom-rel-cofibrations} to show $\epsilon_\gR$ is a weak equivalence it suffices to prove that the augmentation 
\begin{equation}\label{eq:iterated-decomposables}
\pi \colon \fgr{B_\bullet(\bar{F}_{0,n}, E_n, \bar{F}_{n,k}(X))} \lra \bar{F}_{0,n+k}(X)
\end{equation}
is a weak equivalence between cofibrant objects of $\sfC$ as long as $X \in \sfC$ is cofibrant. The symmetric sequences of simplicial sets defining the functors $\bar{F}_{0,n}$, $E_n$, $\bar{F}_{n,k}$, and $\bar{F}_{0,n+k}$ are all cofibrant, so by Lemma \ref{lem.symmetric-sequence-preserving} (i) these functors preserve cofibrant objects. Thus $B_\bullet(\bar{F}_{0,n}, E_n, \bar{F}_{n,k}(X))$ is a Reedy cofibrant semi-simplicial object and so the geometric realization $\fgr{B_\bullet(\bar{F}_{0,n}, E_n, \bar{F}_{n,k}(X))}$ is cofibrant by Lemma \ref{lem:thick-geom-rel-cofibrations}. Thus when $X$ is cofibrant the map \eqref{eq:iterated-decomposables} is indeed between cofibrant objects, and it remains to show that it is a weak equivalence.

As we have done before, we focus on symmetric sequences and show that the augmentation $\fgr{B_\bullet(\bar{F}_{0,n}, E_n, \bar{F}_{n,k})} \to \bar{F}_{0,n+k}$ is a weak equivalence of symmetric sequences. Our proof of this is similar to the proof of Lemma \ref{lem:CubesResolution}. The spaces of $p$-simplices $B_p(F_{0,n},E_n,F_{n,k})(i)$ consists of sequences of embeddings of length $(n+2)$. It is convenient think of the $i$ elements of $F_{n,k}$ as cubes in $\bR^{n+k}$, as the image under composition of all these $(p+2)$ embeddings. We refer to them as the ``innermost cubes.'' Let us define $\partial B_p(F_{0,n},E_n,F_{n,k}) \subset B_p(F_{0,n},E_n,F_{n,k})$ as the subspace where at least one of the innermost cube lies outside the interior of $I^n \times I^k$, that is,
\[\partial B_p(F_{0,n},E_n,F_{n,k})(i) \coloneqq (\partial F_{0,n} \circ E_n^p \circ F_{n,k})(i) \cup (F_{0,n} \circ E_n^p \circ \partial F_{n,k})(i),\]
as a cube lies outside $I^n \times I^k$ if and only if it lies outside $I^n \times \bR^k$ or $\bR^n \times I^k$. As in Lemma \ref{lem:CubesResolution}, it suffices to prove that the two maps
\[\fgr{B_\bullet(F_{0,n},E_n,F_{n,k})} \lra F_{0,n+k}\]
\[\fgr{\partial B_\bullet(F_{0,n},E_n,F_{n,k})}\lra \partial F_{0,n+k}\]
are weak equivalences.  

To do so, as in the proof of Lemma \ref{lem:CubesResolution} we replace cubes by ordered configurations. Consider the symmetric sequence $C_{n,k}$ given by $C_{n,k}(i) \coloneqq \mr{Emb}(\sqcup_i \ast,I^n \times \bR^k)$. There is a map $F_{n,k} \to C_{n,k}$ of symmetric sequences which records the center of the $(i+j)$-cubes, and this is a weak equivalence of left $E_n$-functors. It and its analogue $F_{0,n+k} \to C_{0,n+k}$ yields a commutative diagram
\[\begin{tikzcd}\fgr{B_\bullet(F_{0,n},E_n,F_{n,k})} \rar \dar[swap]{\simeq} & F_{0,n+k} \dar{\simeq} \\
\fgr{B_\bullet(F_{0,n},E_n,C_{n,k})} \rar{\pi} & C_{0,n+k}, \end{tikzcd}\]
and hence it suffices to prove that the bottom map is a weak equivalence. 

The map $\pi$ is a Serre microfibration, so by \cite[Lemma 2.2]{WeissClassify} it suffices to prove that the fibers are weakly contractible. To see this, we observe that the map $\pi \colon B_p(F_{0,n},E_n,C_{n,k})(i) \to C_{0,n+k}(i)$ records the image in $\bR^{n+k}$ of the configuration under the sequence of embeddings of length $(p+1)$, and the embeddings as cubes $I^{n+k}$ or $I^n \times \bR^n$ around these points in $\bR^{n+k}$. Let $\mr{proj}_n \colon \bR^n \times \bR^k \to \bR^n$ denote the projection. Over a configuration $x \in C_{0,n+k}$, the fiber of the augmentation $\pi$ consists of the geometric realization of the semi-simplicial space $X_\bullet(x)$ with $p$-simplices given by the subspace of $F_{0,n} E_n^p C_{n,k}$ consisting of those elements such that the image of the configuration in $\bR^{n+k}$ is $x$ and each point of $\mr{proj}_n(x)$ is contained in the product of $\bR^k$ with the image of some $n$-cube of the innermost layer.

Let $f \colon S^m \to \fgr{X_\bullet(x)}$ be a continuous map. By Lemma \ref{lem:fat-gr-homotopy-to-compact-image}, it may be homotoped so that it factors through the image of compact subspaces $K_j \subset B_j(F_{0,n},E_n,\{x\})$ for $j \leq m$, and thus there exists an $\epsilon_0>0$ such that for all $y \in S^m$ the element $f(y)$ is represent by an configuration of cubes such that the sides of the innermost cubes have distance $>\epsilon_0$ from $\mr{proj}_n(x)$. By picking a collection of $n$-cubes with sides of length $\epsilon_0/2$ around the points in $\mr{proj}_n(x)$, we can cone off $f$ and hence conclude that $\fgr{X_\bullet(x)}$ is weakly contractible.

The argument for $\fgr{\partial B_\bullet(F_{0,n},E_n,F_{n,k})}\to \partial F_{0,n+k}$ is similar. As in the proof of Lemma \ref{lem:CubesResolution} we use a weakly equivalent versions $\partial^\ast F_{0,n+k}$ of $\partial F_{0,n+k}$ and $\partial^\ast B_p(F_{0,n},E_n,F_{n,k})$ of $\partial B_p(F_{0,n},E_n,F_{n,k})$ where the center of at least one innermost cube is outside $I^n \times I^k$. Then the above proof goes through with appropriate modifications.
\end{proof}

\subsection{The $E_\infty$-case and infinite bar spectra}\label{sec:InfiniteBarSpectra}
In Section \ref{sec:e-infty-operad}, we defined the operad $\cC_\infty$ as the colimit of the $\cC_k$. Hence it is not surprising that we may compute the derived $E_\infty$-indecomposables as a homotopy colimit of the derived $E_k$-indecomposables.

\begin{theorem}\label{thm:qeinfy-hocolim} 
There is a zig-zag of natural weak equivalences
	\[\underset{k \to \infty}{\mr{hocolim}}\, Q^{E_k}_{\bL}(\gR) \Leftarrow \cdots \Rightarrow Q^{E_\infty}_{\bL}(\gR) \colon \Alg_{E_\infty}(\sfC) \lra \sfC_*.\]
\end{theorem}

\begin{proof}
Let us write $c\gR \overset{\sim}\to \gR$ for a functorial cofibrant replacement in $\Alg_{E_\infty}(\sfC)$. This is in particular cofibrant in $\sfC$ by Axiom \ref{axiom:monad-proj}, so by Lemma \ref{lem:adding-degeneracies-reedy} we have a free resolution
\[\fgr{B_\bullet(\gE_k, E_k, c\gR)}_{E_k} \lra c\gR\]
in $\Alg_{E_k}(\sfC)$, and as in Section \ref{sec:simplicial-formula-indecomposables} we can compute $Q^{E_k}_\bL(\gR)$ using $\fgr{B_\bullet(+, E_k, c\gR)}$. It therefore remains to compare $\mr{hocolim}_{k \to \infty}\, \fgr{B_\bullet(+, E_k, c\gR)}$ with $\fgr{B_\bullet(+, E_\infty, c\gR)}$.

The $p$th level $B_p(+,E_\infty,c\gR)$ of the semi-simplicial object is obtained by a $p$-fold application of $\cC_\infty$ to $c\gR$, followed by adding a basepoint. Since $\cC_\infty$ is the sequential colimit of the $\cC_k$, and as a consequence of Lemma \ref{lem:odot-props} we have an isomorphism $B_p(+,E_\infty,c\gR) \cong \mr{colim}_{k \to \infty} B_p(+,E_k,c\gR)$. Since $c\gR$ is cofibrant in $\sfC$, by Lemma \ref{lem.symmetric-sequence-preserving} (i) each $B_p(+,E_k,\gR)$ is cofibrant. Similarly, $\cC_k \to \cC_{k+1}$ is a cofibration of symmetric sequences and hence the maps $B_p(+,E_k,c\gR) \to B_p(+,E_{k+1},c\gR)$ are cofibrations. Finally, thick geometric realization preserves cofibrant objects and cofibrations by Lemma \ref{lem:thick-geom-rel-cofibrations}. We thus conclude that
\[\fgr{B_\bullet(+,E_\infty,c\gR)} \cong \underset{k \to \infty}{\colim} \fgr{B_\bullet(+,E_k,c\gR)} \overset{\sim}\longleftarrow \underset{k \to \infty}{\hocolim} \fgr{B_\bullet(+,E_k,c\gR)} \]
where the second map is an equivalence as a sequential colimit is equivalent to the homotopy colimit if all objects are cofibrant and all morphisms are cofibrations.
\end{proof}

Since the $k$-fold suspension of the derived $E_k$-indecomposables was computed by a $k$-fold bar construction, we are similarly able to describe the derived $E_\infty$-indecomposables in terms of an iterated bar construction.

To make this precise, we will define the \emph{infinite bar symmetric spectrum}. This will be an object of the category $\cat{Sp}^\Sigma(\sfC)$ of symmetric spectra in $\sfC$, as defined in \cite{HoveySymmetricModel} (see also \cite{PavlovScholbach3}). (This coincides with $\cat{Sp}^\Sigma$ as in Section \ref{sec:symmetric-spectra} in the case $\sfC=\cat{sSet}$.) An object $E$ of $\cat{Sp}^\Sigma(\sfC)$ consists of a sequence $\{E_n\}_{n \geq 0}$ of objects of $\sfC_\ast$ and a $\fS_n$-action on $E_n$, along with structure maps $E_n \wedge S^1 \to E_{n+1}$ such that the iterated structure maps $E_n \wedge S^k \to E_{n+k}$ are $\fS_n \times \fS_k$-equivariant. For example, given an object $X \in \sfC_\ast$ we may form the suspension spectrum $\Sigma^\infty X \in \cat{Sp}^\Sigma(\sfC)$ by taking $(\Sigma^\infty X)_n \coloneqq X \wedge S^n$ with standard suspension maps. 

The circle $S^1$ is isomorphic to the thick geometric realization of the pointed semi-simplicial set $S^1_\bullet$ with a single $0$-simplex and $1$-simplex. We therefore obtain a $(k+1)$-fold semi-simplicial map
\[\tilde{B}^{E_k}(\gR)_{\bullet,\ldots,\bullet} \wedge S^1_\bullet \lra \tilde{B}^{E_{k+1}}(\gR)_{\bullet,\ldots,\bullet},\]
which yield structure maps 
\[b_k \colon \tilde{B}^{E_k}(\gR) \wedge S^1 \lra \tilde{B}^{E_{k+1}}(\gR)\]
upon geometric realization. The object $\tilde{B}^{E_k}(\gR) = \fgr{\tilde{B}^{E_k}(\gR)_{\bullet,\ldots,\bullet}}$ has an $\fS_k$-action by permuting the $k$ semi-simplicial directions, and this makes the iterated structure maps appropriately equivariant.

\newglossaryentry{binfty}{%
	name={\ensuremath{\tilde{B}^\infty(\gR)}},
	description={Infinite bar spectrum},
	type=symbols
}
\begin{definition}The \emph{infinite bar construction symmetric spectrum}\index{bar construction!infinite} in $\cat{Sp}^\Sigma(\sfC)$ is given by $\gls{binfty} \coloneqq \{\tilde{B}^{E_k}(\gR),b_k\}_{k \geq 0}$, where we set $\tilde{B}^{E_0}(\gR) \coloneqq \ast$.
\end{definition}

We can construct a closely related symmetric spectrum out of the derived $E_k$-indecomposables after picking an explicit model: we may take a cofibrant approximation $c\gR \to \gR$ in the category of $E_\infty$-algebras, and as in the proof of Theorem \ref{thm:qeinfy-hocolim} take $Q^{E_k}_\bL(\gR) = \fgr{B_\bullet(+,E_k,c\gR)}$. 

The inclusion of $E_k$ into $E_{k+1}$ induces a map $Q^{E_k}_\bL(\gR) \to Q^{E_{k+1}}_\bL(\gR)$. We can then define a symmetric spectrum $\tilde{Q}^\infty(\gR) \coloneqq \{ Q^{E_k}_{\bL}(\gR) \wedge S^k,\beta_k\}_{k \geq 0}$ with structure maps $\beta_k \colon Q^{E_k}_{\bL}(\gR) \wedge S^{k} \wedge S^1 \to Q^{E_{k+1}}_{\bL}(\gR) \wedge S^{k+1}$ given by smashing the map $Q^{E_k}_{\bL}(\gR) \to Q^{E_{k+1}}_{\bL}(\gR)$ with the identification $S^k \wedge S^1 \cong S^{k+1}$. In the proof of Theorem \ref{thm:BarHomologyIndec} we exhibited the zig-zag (\ref{eqn:zigzag-bk-qek}) of maps $S^k \wedge Q^{E_k}_{\bL}(\gR) \leftarrow \cdots \to \tilde{B}^{E_k}(\gR)$ which are weak equivalences in $\gR$ is cofibrant in $\sfC$, and when precomposed with the symmetry $Q^{E_k}_\bL(\gR) \wedge S^k \cong S^k \wedge Q^{E_k}_\bL(\gR)$ these assemble to a zig-zag of maps of symmetric spectra:

\begin{lemma}\label{lem:inf-q-sp} There is a natural zig-zag of morphisms
	\[\tilde{Q}^{E_\infty}(\gR) \longleftarrow \cdots \longrightarrow \tilde{B}^\infty(\gR)\]
in $\cat{Sp}^\Sigma(\sfC)$, which are levelwise weak equivalences if $\gR$ is cofibrant in $\sfC$.
\end{lemma}

Under mild conditions on $\sfC$ (namely, that it is left proper cellular, e.g.\ $\sfS = \cat{sSet}$, $\cat{sMod}_\bk$ or $\cat{Sp}^\Sigma$, a property which is preserved by transferring to the projective model structure on $\sfS^\sfG$), Hovey has shown \cite[Theorem 8.2]{HoveySymmetricModel} that there is a projective model structure on $\cat{Sp}^\Sigma(\sfC)$, with a localization called the stable model structure. Levelwise weak equivalences are weak equivalences in either of these model structures. Hovey also proves that if the functor $S^1 \wedge - \colon \sfC_* \to \sfC_*$ is already a Quillen equivalence, then $\sfC_*$ is Quillen equivalent to $\cat{Sp}^\Sigma(\sfC)$ with the stable model structure \cite[Theorem 9.1]{HoveySymmetricModel}. Thus in the case $\sfC = \cat{Sp}^\Sigma$ we obtain a Quillen equivalent category. 

In the case $\sfC =\cat{sSet}$, \cite[Theorem 3.1.11]{HoveyShipleySmith} says that stable homotopy equivalences are stable equivalences (see also \cite{SchwedeHomotopy}). Using this fact we may deduce the following.

\begin{corollary}\label{cor:InfBarSpectra}
Suppose that $\sfC = \cat{sSet}^\sfG$ with $\sfG$ discrete, then there is a natural zig-zag of morphisms
	\[\Sigma^\infty Q^{E_\infty}_{\bL}(\gR) \longleftarrow \cdots \longrightarrow \tilde{B}^\infty(\gR)\]
in $\cat{Sp}^\Sigma(\sfC)$, which are weak equivalences if $\gR$ is cofibrant in $\sfC$.\end{corollary}

\begin{proof}By Lemma \ref{lem:inf-q-sp} it suffices to show that $\Sigma^\infty Q^{E_\infty}_{\bL}(\gR)$ and $\tilde{Q}^{E_\infty}(\gR)$ are stable equivalent for each $g \in \sfG$. The inclusions $E_k \hookrightarrow E_\infty$ induce maps $Q^{E_k}_{\bL}(\gR) \wedge S^{k} \to Q^{E_\infty}_{\bL}(\gR) \wedge S^{k}$ which assemble into a map of symmetric spectra $f \colon \tilde{Q}^\infty(\gR) \to \Sigma^\infty Q^{E_\infty}_\bL(\gR)$. 
	
Since $\sfG$ is discrete, a map $X \to Y$ in $\cat{Sp}^\Sigma(\sfC)$ is a stable equivalence if and only if each $X(g) \to Y(g)$ is a stable equivalence in $\cat{Sp}^\Sigma$. Thus to see that $f$ is a stable equivalence, it suffices to prove that $\tilde{Q}^\infty(\gR)(g) \to \Sigma^\infty Q^{E_\infty}_\bL(\gR)(g)$ is a stable homotopy equivalence. This follows if the map of pointed simplicial sets $Q^{E_k}_{\bL}(\gR)(g) \to Q^{E_\infty}_{\bL}(\gR)(g)$ is $(k-1)$-connected. The thick geometric realization of a levelwise $(k-1)$-connected map between semi-simplicial simplicial sets is $(k-1)$-connected, and thus it suffices to prove that for each $p \geq 0$ the map $B_p(+,E_k,c\gR)(g) \to B_p(+,E_\infty,c\gR)(g)$ is $(k-1)$-connected. To prove this, we use that both functors $\cC_k(-)$ and $\cC_\infty(-)$ preserve connectivity of maps because $\cC_k$ and $\cC_\infty$ are $\Sigma$-cofibrant, and that $\cC_k(X) \to \cC_\infty(X)$ is $(k-1)$-connected for all $X \in \cat{sSet}_\ast$ because the map of $\Sigma$-cofibrant operads $\cC_k \to \cC_\infty$ is $(k-1)$-connected.
\end{proof}

\begin{remark}We expect this corollary to be true more generally; in $\sfS^\sfG$ for $\sfS = \cat{sSet}$, $\cat{sMod}_\bk$ or $\cat{Sp}^\Sigma$, and $\sfG$ any diagram category.\end{remark}

\subsection{Group completion}\label{sec:GpCompletion}

We shall discuss group completion in the setting of $E_k$-algebras in the category $\cat{Top}$. The discussion goes through unchanged for $\cat{sSet}$, but we make no claims about the other examples of categories we have discussed (one requirement seems to be that the category is \emph{semi-cartesian}: the morphism $\bunit \to \bterm$ is an isomorphism).

\subsubsection{The group completion augmentation} 
\newglossaryentry{egc}{%
	name={\ensuremath{\epsilon_\mr{gc}}},
	description={Group completion augmentation},
	type=symbols
}
In the monoidal category $(\cat{Top}, \times, *)$ the terminal object $\bterm = *$ is also the monoidal unit $\bunit = *$, and so any $E_k^+$-algebra $\gR$ in $\cat{Top}$ has a canonical map $\gls{egc} \colon \gR \to \bunit$ of $E_k^+$-algebras, which we call the \emph{group completion augmentation}\index{augmentation!group-completion}.\index{group completion!augmentation}

Suppose we have a filtered non-unital $E_k$-algebra $\gR \in \Alg_{E_k}(\cat{Top}^{\bN_{\leq}})$, with underlying $E_k$-algebra $\mr{colim}(\gR) \in \Alg_{E_k}(\cat{Top})$. If we unitalize $\gR$ then we obtain an $E_k^+$-algebra $\gR^+ \in \Alg_{E_k^+}(\cat{Top}^{\bN_\leq})$ satisfying $\mr{colim}(\gR^+) \cong \mr{colim}(\gR)^+$. 

Now $\mr{colim}(\gR^+)$ has the group completion augmentation $\epsilon_\mr{gc} \colon \mr{colim}(\gR^+) \to *$ described above, which by adjunction gives a map $\epsilon_\mr{gc}' \colon \gR^+ \to 0_*(*)$ in $\smash{\Alg_{E_k^+}(\cat{Top}^{\bN_\leq})}$. (Note that $0_*(*)$ is the terminal object of $\cat{Top}^{\bN_\leq}$, so is canonically a $E_k^+$-algebra.) The associated graded of $0_*(*)$ is the object $0_*(S^0) \in \cat{Top}_*^{\bN_{=}}$, which is the monoidal unit \emph{when $\cat{Top}_*$ is given a monoidal structure via smash product}, and $\cat{Top}_*^{\bN_{=}}$ is given the induced monoidal structure by Day convolution. Therefore the colimit of the map $\epsilon_\mr{gc}'$ is the map $\epsilon_\mr{gc}$, and its associated graded is the canonical augmentation $\epsilon_\mr{can} \colon \grr(\gR)^+ \to \grr(0_*(*)) = 0_*(S^0)$.

\begin{proposition}\label{prop:group-completion-ss}\index{spectral sequence!group completion}
If $\gR$ is cofibrant in $\cat{Top}^{\bN_\leq}$ then there is a {group completion spectral sequence}
\[E^1_{p,q} = H_{p+q-k,p}^{E_k}(\grr(\gR)) \Longrightarrow \tilde{H}_{p+q}({B}^{E_k}(\mr{colim}(\gR)^+,\epsilon_\mr{gc}))\]
which converges strongly, with differentials $d^r \colon E^r_{p,q} \to E^r_{p-r,q+r-1}$.
\end{proposition}

\begin{proof}
Consider the filtered object ${B}^{E_k}(\epsilon_\mr{gc}')$ given by the bar construction of Definition \ref{def:kfold-bar-augmented} applied to the morphism $\epsilon_\mr{gc}' \colon \gR^+ \to 0_*(*)$ of filtered $E_k^+$-algebras. As $\grr$ is strongly monoidal (by Section \ref{sec:MonoidalityOnGr}) and ${B}^{E_k}$ commutes with $\grr$ (as it is a geometric realization of iterated tensor products), we have an isomorphism
\[\grr({B}^{E_k}(\epsilon_\mr{gc}')) \cong {B}^{E_k}(\grr(\gR)^+,\epsilon_\mr{can})\]
in $\cat{Top}_*^{\bN_{=}}$. For the same reason ${B}^{E_k}$ commutes with $\mr{colim}$, so there is an isomorphism
\[\mr{colim}({B}^{E_k}({\epsilon}_\mr{gc}')) \cong {B}^{E_k}(\mr{colim}(\gR)^+,\epsilon_\mr{gc})\]
in $\cat{Top}$. If $\gR$ is cofibrant in $\cat{Top}^{\bN_{\leq}}$ then ${B}^{E_k}(\gR^+,{\epsilon}'_\mr{gc})$ is cofibrant in $\cat{Top}^{\bN_\leq}$. Thus we may apply Theorem \ref{thm:SSAsc} to obtain a spectral sequence 
\[E^1_{p,q} = \widetilde{H}_{p+q,p}({B}^{E_k}(\grr(\gR)^+,\epsilon_\mr{can})) \Longrightarrow {H}_{p+q}({B}^{E_k}(\mr{colim}(\gR)^+,\epsilon_\mr{gc})),\]
which in this case converges strongly by \cite[Theorem 6.1]{Boardman}. 

We have ${B}^{E_k}(\grr(\gR)^+,\epsilon_\mr{can}) \simeq S^0 \vee \tilde{B}^{E_k}(\grr(\gR)^+,\epsilon_\mr{can})$ by Lemma \ref{lemma:bek-cofibrant}. As $\grr(\gR)$ is cofibrant in $\cat{Top}_*^{\bN_{=}}$ (as $\grr$ is a left Quillen functor), by Theorem \ref{thm:BarHomologyIndec} we have $\tilde{B}^{E_k}(\grr(\gR)^+,\epsilon_\mr{can}) \simeq S^k \wedge Q^{E_k}_\bL(\gR)$. Thus the induced spectral sequence on homology relative to $*$ has the indicated $E^1$-page.
\end{proof}

\begin{remark}
In the case $k=\infty$, as in Section \ref{sec:InfiniteBarSpectra} the collection of spaces $\{{B}^{E_k}(\mr{colim}(\gR)^+,\epsilon_\mr{gc})\}_{k \geq 0}$ assemble into a symmetric spectrum $B^\infty(\mr{colim}(\gR)^+,\epsilon_\mr{gc})$, and as above (but using Corollary \ref{cor:InfBarSpectra}) we obtain a spectral sequence
\[E^1_{p,q} = H_{p+q,p}^{E_\infty}(\grr(\gR)) \Longrightarrow {H}_{p+q}^\mr{spec}({B}^{\infty}(\mr{colim}(\gR)^+,\epsilon_\mr{gc}))\]
converging to the spectrum homology of this spectrum.\index{spectral sequence!group completion}
\end{remark}

In the remainder of this section, our goal will be to show that ${B}^{E_k}(-,\epsilon_\mr{gc})$ coincides up to natural weak equivalence with classical delooping constructions in topology. This is a folklore result, but we do not know a reference.

\subsubsection{A variation of the $k$-fold iterated bar construction} 
To study group completions, it will be convenient to use a variation of Definition \ref{def:kfold-bar-augmented}. 

\begin{definition}\label{def:cik} Let us write $\cI_k(p_1,\ldots,p_k)$ for the space of collections of $k$-tuples $\{[a_i^j,b_i^j]\}$ of intervals for $1 \leq j \leq k$ and $0 \leq i \leq p_j$ with endpoints $0<a_i^j<1/2<b_i^j<1$ such that $[a_i^j,b_i^j] \supsetneq [a_{i+1}^j,b_{i+1}^j]$. This is a $k$-fold semi-simplicial space where the $i$th face map in the $j$th direction $d_i^j$ forgets $[a_i^j,b_i^j]$.\end{definition}

\begin{figure}
	\centering
	\begin{tikzpicture}
	
	\fill[Mahogany!2!white] (0,0) rectangle (6,6);
	\fill[Mahogany!8!white] (0.5,1) rectangle (5,5);
	\fill[Mahogany!16!white] (2,1) rectangle (4,5);

	\foreach \x in {.5,2,4,5}
	\foreach \y in {1,5}
	{
		\draw[black] (\x,0) -- (\x,6);
		\draw[black] (0,\y) -- (6,\y);
	}
	
	\draw[dotted] (0,0) -- (0,6);
	\draw[dotted] (6,0) -- (6,6);
	\draw[dotted] (0,0) -- (6,0);
	\draw[dotted] (0,6) -- (6,6);
	
	\node at (0.25,.5) {$\gS$};
	\node at (1.25,.5) {$\gS$};
	\node at (3,.5) {$\gS$};
	\node at (4.5,.5) {$\gS$};
	\node at (5.5,.5) {$\gS$};
	
	\node at (0.25,3) {$\gS$};
	\node at (1.25,3) {$\gR$};
	\node at (3,3) {$\gR$};
	\node at (4.5,3) {$\gR$};
	\node at (5.5,3) {$\gS$};

	\node at (0.25,5.5) {$\gS$};
	\node at (1.25,5.5) {$\gS$};
	\node at (3,5.5) {$\gS$};
	\node at (4.5,5.5) {$\gS$};
	\node at (5.5,5.5) {$\gS$};
	
	\node[below] at (0,0) {$0$};
	\node[below] at (.5,0) {$a_0^1$};
	\node[below] at (2,0) {$a_1^1$};
	\node[below] at (4,0) {$b_1^1$};
	\node[below] at (5,0) {$b_0^1$};
	\node[below] at (6,0) {$1$};
	
	\node[left] at (0,6) {$1$};
	\node[left] at (0,1) {$a_0^2$};
	\node[left] at (0,5) {$b_0^2$};
	\node[left] at (0,0) {$0$};
	\end{tikzpicture}
	\caption{An illustration of $\cB^{E_2}_{1,0}(f)$.}
	\label{fig:bar-construction-variation}
\end{figure}

\begin{definition}\label{def:bek-augmented-variation} Let $f \colon \gR \to \gS$ be a morphism of $E_k^+$-algebras in $\cat{Top}$. Then $\cB^{E_k}_{\bullet,\ldots,\bullet}(f)$ is the $k$-fold semi-simplicial object with $\cB^{E_k}_{p_1,\ldots,p_k}(f) \coloneqq \cI_k(p_1,\ldots,p_k) \times \mathcal{G}_{p_1, \ldots, p_k}(f)$, where
\[\mathcal{G}_{p_1, \ldots, p_k}(f) \coloneqq \prod_{q_1=0}^{2p_1+2} \cdots \prod_{q_k=0}^{2p_k+2} \mathcal{B}_{p_1, \ldots, p_k}^{q_1,\ldots,q_k}\]
and $\mathcal{B}_{p_1, \ldots, p_k}^{q_1,\ldots,q_k}$ is $\gR$ if $1 \leq q_j \leq 2p_j+1$ for all $j$, and is $\gS$ otherwise.

The $i$th face map $d_i^j$ in the $j$th direction
\[d_i^j \colon \mathcal{B}^{E_k}_{p_1,\ldots,p_k}(f) \lra \mathcal{B}^{E_k}_{p_1,\ldots, p_{j-1}, p_j-1, p_{j+1},p_k}(f)\]
is given by the face map of Definition \ref{def:cik} on the first factor and then, on the term $\mathcal{G}_{p_1, \ldots, p_k}(\gR)$, there are two cases. The precise description is analogous to Definition \ref{def:kfold-bar-augmented} (or Section \ref{sec:bar-modification}), but here we settle for a heuristic description:\index{bar construction!variation of iterated}
\begin{enumerate}[(i)]
	\item For $0 < i \leq p_j$ the map $d_i^j$ is induced by the $E_k^+$-algebra structure on $\gR$ by applying elements analogous to (\ref{eqn:bar-embedding-deltaij}) to the terms $\mathcal{B}^{q_1,\ldots,q_{j-1},i,q_{j+1},\ldots,q_k}_{p_1, \ldots, p_k}$ and $\mathcal{B}^{q_1,\ldots,q_{j-1},i+1,q_{j+1},\ldots,q_k}_{p_1, \ldots, p_k}$, and to the terms $\mathcal{B}^{q_1,\ldots,q_{j-1},2p_i-i,q_{j+1},\ldots,q_k}_{p_1, \ldots, p_k}$ and $\mathcal{B}^{q_1,\ldots,q_{j-1},2p_i-i-1,q_{j+1},\ldots,q_k}_{p_1, \ldots, p_k}$.
	\item The map $d^j_0$ is given by first applying the $f$ to each of the entries $\mathcal{B}_{p_1, \ldots, p_k}^{q_1,\ldots,q_{j-1},1,q_j,\ldots,q_k}$ and $\mathcal{B}_{p_1, \ldots, p_k}^{q_1,\ldots,q_{j-1},2p_j-1,q_j,\ldots,q_k}$ and then applying elements analogous to (\ref{eqn:bar-embedding-deltaij}).
\end{enumerate}
	
	We write $\cB^{E_k}(f) \coloneqq \fgr{\cB_{\bullet, \ldots, \bullet}^{E_k}(f)} \in \sfC$. This construction is natural in commutative diagrams of maps of $E_k^+$-algebras. 
\end{definition}

\newglossaryentry{bekvar}{%
	name={\ensuremath{\cB^{E_k}(\gR,\epsilon)}},
	description={Variation on iterated bar construction},
	type=symbols
}
\newglossaryentry{bekvarred}{%
	name={\ensuremath{\tilde{\cB}^{E_k}(\gR,\epsilon)}},
	description={Variation on reduced iterated bar construction},
	type=symbols
}
We shall be interested in the case that $f$ is an augmentation $\epsilon \colon \gR \to \bunit$. Then, as before, we write $\gls{bekvar}$ for $\cB^{E_k}(\epsilon)$ and denote the cofiber of $\cB^{E_k}(\bunit,\epsilon_\bunit) \to \cB^{E_k}(\gR,\epsilon)$ by $\gls{bekvarred}$. 

To compare $B^{E_k}(\gR,\epsilon)$ and $\cB^{E_k}(\gR,\epsilon)$, we shall use that $\cB^{E_k}(f)$ is obtained from $B^{E_k}(f)$ by $k$-fold \emph{edgewise subdivision} of \cite[Appendix I]{SegalConfiguration} up to weak equivalence. Recall that the $[n]$ denotes the ordered finite set $\{0<1<\ldots<n\}$. A single edgewise subdivision is obtained by precomposing a functor $\Delta^\mr{op} \to \cat{Top}$ or $\Delta^\mr{op}_\mr{inj} \to \cat{Top}$ with (the opposite of) the functor
\begin{align*}\mr{esd} \colon \Delta &\lra \Delta \\
[n] &\longmapsto [n]^\mr{op} \ast [n],\end{align*}
with $\ast$ denoting the join of finite ordered sets, so that $[n]^\mr{op} \ast [n]$ is an ordered set with $2n+2$ elements, i.e. isomorphic to $[2n+1]$. 

There is a natural transformation $\mr{id} \Rightarrow \mr{esd}$ given by the inclusion of the second $[n]$, and only in the case of $\Delta^\mr{op}$, also a natural transformation $\mr{esd} \Rightarrow \mr{id}$ by collapsing the first copy of $[n]$ onto the first element of the second copy of $[n]$. Since natural transformations induce simplicial homotopies, this implies that the thick geometric realization of an edgewise subdivision of a simplicial object $X_\bullet$ is homotopy equivalent to the thick geometric realization of $X_\bullet$, hence is weakly equivalent. We use these observations to prove the next lemma. 

\begin{lemma}$B^{E_k}(\gR,\epsilon)$ and $\cB^{E_k}(\gR,\epsilon)$ are naturally weakly equivalent.\end{lemma}

\begin{proof}The $k$-fold semi-simplicial space $B^{E_k}_{\bullet,\ldots,\bullet}(\gR,\epsilon)$ is levelwise weakly equivalent to a variation $B^{E_k,0}_{\bullet,\ldots,\bullet}(\gR,\epsilon)$. Let $\cP_k^0(p_1,\ldots,p_k)$ be as in Definition \ref{def:pk} but replace the condition that $0<t_0^j<\cdots<t^j_{p_j}<1$ by $0 \leq t^j_0 \leq \ldots \leq t^j_{p_j} \leq 1$. Using the notation of Definition \ref{def:kfold-bar-augmented}, \[B^{E_k,0}_{p_1,\ldots,p_k}(\gR,\epsilon) \subset \cP_k^0(p_1,\ldots,p_k) \times G_{p_1,\ldots,p_k}(\epsilon)\]
is the subspace where the element of the term $B^{q_1,\ldots,q_k}_{p_1,\ldots,p_k}$ of $G_{p_1,\ldots,p_k}(f)$ is required to be the unit of $\gR$ when $\smash{t^j_{q_j-1}} = \smash{t^j_{q_j}}$ for some $j$ (by convention $\smash{t^j_{-1}} = 0$). Informally, grids are allowed to contain cubes of volume $0$ but only when labeled by the unit. We may can the face maps of $B^{E_k}_{\bullet,\ldots,\bullet}(\gR,\epsilon)$ to $B^{E_k,0}_{\bullet,\ldots,\bullet}(\gR,\epsilon)$ as follows: $d^j_i$ is as before when $\smash{t^j_{i-1}} \neq \smash{t^j_i}$ and otherwise projects away the corresponding terms (necessarily given by units of $\gR$).
	
This admits a system of degeneracy maps by defining $s_i^j$ to duplicate $t^j_i$'s and adding the unit of $\gR$. We thus have weak equivalences
\[B^{E_k}(\gR,\epsilon) = \fgr{B^{E_k}_{\bullet,\ldots,\bullet}(\gR,\epsilon)} \simeq \fgr{B^{E_k,0}_{\bullet,\ldots,\bullet}(\gR,\epsilon)} \simeq \fgr{\mr{esd}\,B^{E_k,0}_{\bullet,\ldots,\bullet}(\gR,\epsilon)}.\]
The semi-simplicial space $\cB^{E_k}_{\bullet,\ldots,\bullet}(\gR,\epsilon)$ is levelwise weakly equivalent to the edgewise subdivision $\mr{esd}\,B^{E_k,0}_{\bullet,\ldots,\bullet}(\gR,\epsilon)$, differing only in the intervals allowed. As in $\cat{Top}$ the construction $- \times \gR$ preserves weak equivalences even if $\gR$ is not cofibrant, this induces a weak equivalence upon thick geometric realization.\end{proof}

\subsubsection{The group completion map}\label{sec:group-completion-map} In Section \ref{sec:FreeAlgMaps} we described a natural transformation $\eta \colon S^k \wedge E_k(-) \Rightarrow S^k \wedge -$ of functors $\sfC_\ast \to \sfC_\ast$ which is related to a natural transformation $\eta_M \colon \Sigma^k E_k \Rightarrow \Sigma^k$ of functors $\cat{Top} \to \cat{Top}_*$ due to May \cite{GILS} which we will now describe. Here we write $\Sigma^k(-) = S^k \wedge (-)_+$. The essential feature of May's map is that its adjoint $E_k \Rightarrow \Omega^k\Sigma^k$ is a map of monads.

\newglossaryentry{etam}{%
	name={\ensuremath{\eta_M}},
	description={May's natural transformation},
	type=symbols
}
The natural transformation $\gls{etam}$ is given by maps $S^k \wedge (\cC_k(i) \times_{\fS_i} X^{\times i})_+ \to S^k \wedge X_+$ defined as follows: a tuple $(v,e ; x_1, \ldots, x_i)$ of a $v \in ((0,1)^k)^+ \cong S^k$, a $e = (e_1,\ldots,e_i) \in \cC_k(i)$, and labels $x_1, \ldots, x_i \in X$, is mapped to the basepoint unless $v$ is contained in a cube $e_j$, in which case we normalize to identify that cube with $(0,1)^k$ and record the position of $v$ in this cube and the label $x_j$. (Note that this is different to the construction in Section \ref{sec:FreeAlgMaps}, where the analogous map was to the basepoint if $i>1$.) It is easy to verify that $\eta_M$, as well as giving $\Sigma^k$ the structure of a right $E_k$-functor, is a map of right $E_k$-functors. Hence its adjoint natural transformation $\eta_M^\vee \colon E_k \Rightarrow \Omega^k \Sigma^k$ of functors $\cat{Top} \to \cat{Top}$ is also a map of right $E_k$-functors.

Furthermore, $\Omega^k \Sigma^k$ also has the structure of a left $E_k$-functor and the natural transformation $\eta_M^\vee$ is one of left $E_k$-functors. The left $E_k$-functor structure on $\Omega^k \Sigma^k$ is described in Section 5 of \cite{GILS} and appeared in Example \ref{exam:loop-spaces}; it is induced by maps $\cC_k(i) \times (\Omega^k S^k)^i \to \Omega^k S^k$ defined as follows. Given an element $(e,f_1,\ldots,f_i) \in \cC_k(i) \times (\Omega^k S^k)^i$ we may form the pointed map $S^k \to S^k$
\[\theta(e,f_1,\ldots,f_k) \colon x \longmapsto \begin{cases}f_i(e^{-1}(x)) & \text{if $x \in \mr{im}(e_i)$,} \\
* & \text{otherwise.}\end{cases}\]

Thus we may produce a zig-zag of maps of non-unital $E_k$-algebras in $\cat{Top}$ as in Theorem 13.1 of \cite{GILS}
\begin{equation*}
\gR \longleftarrow \gr{B_\bullet(E_k,E_k,\gR)} \lra  \gr{B_\bullet(\Omega^k\Sigma^k,E_k,\gR)} \lra \Omega^k \gr{B_\bullet(\Sigma^k,E_k,\gR)},
\end{equation*}
where it may be helpful to observe that the leftmost map is a weak equivalence with a section induced by the inclusion $\mr{id} \times \gR \to B_0(E_k,E_k,\gR)$. We call the resulting zig-zag
\[\gR \overset{\simeq}\longleftarrow \gr{B_\bullet(E_k,E_k,\gR)} \lra \Omega^k \gr{B_\bullet(\Sigma^k,E_k,\gR)}\] \emph{group completion}.\index{group completion}

\begin{definition}We denote $\gr{B_\bullet(\Sigma^k,E_k,\gR)} \in \cat{Top}_\ast$ by $B^k_\mr{M} (\gR)$.\end{definition}

The following is consequence of Lemma \ref{lem:reedy-cofibrancy-simplicial} (i).

\begin{lemma}\label{lem:bkm-cofibrancy} If $\gR$ is cofibrant in $\cat{Top}$, then $B^k_\mr{M}(\gR)$ is cofibrant in $\cat{Top}_\ast$. The functor $B^k_\mr{M}(-) \colon \Alg_{E_k}(\cat{Top}) \to \cat{Top}_\ast$ preserves (trivial) cofibrations between $E_k$-algebras with cofibrant underlying objects.\end{lemma}

\begin{remark}\label{rem:group-completion-grouplike} May has proved that the group completion is a weak equivalence if and only if the monoid $\pi_0(\gR)$ is a group (\cite{GILS} proves this for $\pi_0(\gR)$ trivial and \cite{may-group-completion} addresses the case when $\pi_0(\gR)$ is a group).
\end{remark}

\subsubsection{Comparing deloopings}
We claim that $B^k_\mr{M}(-)$ and $\cB^{E_k}(-,\epsilon_\mr{gc})$ are naturally weakly equivalent when applied to $E_k$-algebras $\gR$ that are cofibrant in $\cat{Top}$, which is stated as Corollary \ref{cor:bek-vs-bkm} below. Note that
\[\cB^{E_k}_{0,\ldots,0}(\gE_k^+(X),\epsilon_\mr{gc}) = \cI_k(0, \ldots, 0) \times \gE_k^+(X).\]
This has a map to $S^k \wedge X_+$ given by forgetting the first factor, and then taking the composition
\[\gE_k^+(X) \subset \gE_k^+(X)_+ = S^0 \wedge \gE_k^+(X)_+ \lra S^k \wedge \gE_k^+(X)_+ \overset{\eta_M}\lra S^k \wedge X_+\]
where the middle map is induced by the inclusion $S^0 = (1/2,\ldots,1/2)_+ \subset S^k = ((0,1)^k)^+$. Unravelling this definition, the map is given by the basepoint unless the point $p_0 \coloneqq (1/2,\ldots,1/2) \in (0,1)^k$ lies in one of the cubes $e_j$, in which case we record $e_j^{-1}(p_0)$ and the label in $X$. It is easy to see that this coequalizes all face maps arriving at $\cB^{E_k}_{0,\ldots,0}(\gE_k^+(X),\epsilon_\mr{gc})$, so defines an augmentation
\[\cB^{E_k}_{\bullet,\ldots,\bullet}(\gE_k^+(X),\epsilon_\mr{gc}) \lra \Sigma^k X.\]
It is in order to define this map that we have passed from $B^{E_k}$ to the model $\cB^{E_k}$: in the former we have $B^{E_k}_{0,\ldots,0}(\gE_k^+(X),\epsilon_\mr{gc}) \simeq S^0$ so there is no useful augmentation.

By construction the augmentation maps the sub-object $\cB^{E_k}_{\bullet,\ldots,\bullet}(\bunit,\epsilon_\bunit)$ to the basepoint and thus upon geometric realization it factors over a natural transformation
\begin{equation}\label{eq:ComparingAug}
\tilde{\cB}^{E_k}(\gE_k^+(-),\epsilon_\mr{gc}) \Rightarrow \Sigma^k(-) \colon \cat{Top} \lra \cat{Top}_*,
\end{equation}
which is a natural transformation of right $E_k$-functors. The following is a well-known folklore result with a proof analogous to that of Theorem \ref{thm:CalcFree}, which we shall only sketch.

\begin{lemma}\label{lem:bek-vs-bkm} 
The map $\tilde{\cB}^{E_k}(\gE^+_k(X), \epsilon_\mr{gc}) \to \Sigma^k X$ is a weak equivalence if $X$ is cofibrant.
\end{lemma}

\begin{proof}[Sketch of proof]
As in the proof of Theorem \ref{thm:CalcFree}, we can obtain zig-zag of natural transformation analogous to (\ref{eqn:zigzag-bk-qek}) consisting of weak equivalences when $X$ is cofibrant, by comparing $\tilde{\cB}^{E_k}(\gE_k^+(X), \epsilon_\mr{gc})$ to a space of cubes in $\bR^k$ with labels in $X$ that can disappear outside $[0,1]^k$. By applying a scanning argument to this space, we may weakly deformation retract onto the subspace where at most a single cube appears, which is visibly weakly equivalent to $\Sigma^k X$.
\end{proof}

The conclusion of this section is as follows:

\begin{corollary}\label{cor:bek-vs-bkm}
There is a zig-zag of natural transformations 
	\[\tilde{B}^{E_k}(\gR^+,\epsilon_\mr{gc}) \Leftarrow \cdots \Rightarrow B^k_\mr{M}(\gR),\]
of functors $\Alg_{E_k}(\cat{Top}) \to \cat{Top}_{\ast}$, which consists of weak equivalences if $\gR$ is cofibrant in $\cat{Top}$.
\end{corollary}
\begin{proof}
Consider the following zig-zag 
\[\tilde{\cB}^{E_k}(\gR^+, \epsilon_\mr{gc}) \overset{\simeq}{\longleftarrow} \fgr{B_\bullet(\tilde{\cB}^{E_k}(\gE^+_k(-), \epsilon_\mr{gc}),E_k,\gR)} \lra  \fgr{B_\bullet(\Sigma^k,E_k,\gR)} \lra B^k_\mr{M}(\gR)\]
The left-hand map is a weak equivalence by commuting the two bar constructions, as $\tilde{\cB}^{E_k}((-)^+, \epsilon_\mr{gc}) \colon \Alg_{E_k}(\cat{Top}) \to \cat{Top}_*$ preserves weak equivalences between objects which are cofibrant in $\cat{Top}$. By Lemma \ref{lem:bek-vs-bkm} the natural transformation \eqref{eq:ComparingAug} is a weak equivalence when applied to a cofibrant object, so by Lemma \ref{lem:thick-geom-rel-cofibrations} the middle map is a weak equivalence. The right-hand map is the natural map from the thick to the thin geometric realization, which is a weak equivalence as long as the simplicial space $B_\bullet(\Sigma^k, E_k, \gR)$ is Reedy cofibrant by Lemma \ref{lem:thick-to-thin}, which it is by Lemma \ref{lem:reedy-cofibrancy-simplicial}.
\end{proof}

\subsubsection{Group completion by cell attachment} \index{group completion!by cell attachment} As an application of the preceding discussion, we explain an alternative construction of the group completion for $E_k$-algebras in $\cat{Top}$. As a consequence of Theorem \ref{thm:BarHomologyIndec}, up to weak equivalence $\tilde{B}^{E_k}$ takes cell attachments in $\Alg_{E_k}(\sfC)$ to cell attachments in $\sfC_\ast$ with a $k$-fold suspension. By considering May's delooping, we will obtain a similar result for $B^k_M $, which we need to prove our alternative group completion is indeed a group completion.

Any map $e \colon \partial D^d \to \gR$ induces a map $\Sigma^k e \colon \Sigma^k \partial D^d \to \Sigma^k \gR$, and since the latter are the $0$-simplices of the simplicial object defining $B_\mr{M}^k(\gR)$, we get an induced map
\[\Sigma^k e \colon \Sigma^k \partial D^d \lra B^k_\mr{M}(\gR).\]
Using Corollary \ref{cor:bek-vs-bkm}, we get an analogue for $B^k_\mr{M}$ of the result that $\tilde{B}^{E_k}(-)$ takes $E_k$-cell attachments to cell attachments. The difference is that $k$-fold bar construction used the canonical augmentation, whereas the lemma below uses the group completion augmentation.

\begin{lemma}\label{lem:bkm-cell-attachment}
Given a map $e \colon \partial D^d \to \gR$, the map $B^k_M(\gR) \cup_{\Sigma^k e} \Sigma^k D^d \to B^k_M(\gR \cup^{E_k}_e \gD^d)$ is a weak equivalence.
\end{lemma}

\begin{proof}
We shall first prove this when $\gR = \gE_k(X)$ with $X$ cofibrant. In that case we saw in Section \ref{sec:derived-cell-simplicial} that we may compute a derived cell attachment by taking the thick geometric realization of the simplicial object $[p] \mapsto \gE_k(E_k^{p}(X) \cup_{e} D^d)$ where $e \colon \partial  D^d \to E_k^{p}(X)$ is given by applying the unit transformation for the monad $E_k$ $p$ times. In this case, by an extra degeneracy argument as in Lemma \ref{lem:extra-degeneracy} to $B^k_\mr{M}(\fgr{\gE_k(E_k^{p}(X) \cup_{e} D^d)})$, we may compute $B^k_\mr{M}$ as the thick geometric realization of the simplicial object 
	\[[p] \mapsto \Sigma^k(E_k^{p}(X) \cup_e D^d) 
	\cong \Sigma^k(E_k^{p}(X)) \cup_{\Sigma^k e} \Sigma^k D^d\]
	and since geometric realization commutes with push-outs, we see this is in turn isomorphic to $B^k_M(\gE_k(X)) \cup_{\Sigma^k e} \Sigma^k D^d$.
	
	The general case follows by taking a free simplicial resolution $\gR_{\bullet} \to \gR$, in particular the thick monadic bar construction $B_\bullet(\gE_k,E_k,\gR) \to \gR$. Then we have a commutative diagram
	\[\begin{tikzcd} {\fgr{[p] \mapsto B^k_\mr{M}(\gE_k(E_k^p \gR) \cup^{E_k}_e \gD^d)}} \rar{\simeq} & {B^k_\mr{M}(\gR \cup^{E_k}_e \gD^d)} \\
	{\fgr{[p] \mapsto B^k_\mr{M}(\gE_k(E_k^p \gR)) \cup_{\Sigma^k e} \Sigma^k D^d}} \rar{\simeq} \uar{\simeq} & {B^k_\mr{M}(\gR) \cup_{\Sigma^k e} \Sigma^k D^d} \uar \end{tikzcd}\] 
	where the left vertical map is a weak equivalence by Lemma \ref{lem:thick-geom-rel-cofibrations} because it is a levelwise weak equivalence between levelwise cofibrant semi-simplicial objects. That it is a levelwise weak equivalence follows from the case of free $E_k$-algebras discussed above. For levelwise cofibrancy we remark that $- \cup^{E_k}_e \gD^d$ and $- \cup_{\Sigma^k e} \Sigma^k D^d$ preserve cofibrant objects in $\Alg_{E_k}(\cat{Top})$ and $\cat{Top}_\ast$ respectively, so that it suffices to prove that $B^k_\mr{M}$ preserves cofibrant objects, which follows from Lemma \ref{lem:bkm-cofibrancy}. We conclude that the right vertical map is a weak equivalence, proving the lemma.
\end{proof}

Suppose that $\gR$ is a $E_k$-algebra in $\cat{Top}$ with $\pi_0(\gR) \cong \bN_{>0}$ as a monoid, and let $\sigma_1 \in \gR$ be a choice of point in the path-component corresponding to $1 \in \bN_{>0} \cong \pi_0(\gR)$. Then we may construct a new $E_k^+$-algebra from the unitalization $\gR^+$ of $\gR$ by freely adding a $0$-cell $\sigma_{-1}$. The result $\gR^+ \cup^{E_k} \sigma_{-1}$ has $\pi_0$ isomorphic as a unital monoid to the free unital abelian monoid on two generators $\sigma_1$ and $\sigma_{-1}$. We then further attach a $1$-cell $\eta$ along a map $\partial D^1 \to \gR^+ \cup^{E_k} \sigma_{-1}$ sending $0$ to the unit in $\gR^+$ and $1$ to a point in the component corresponding to $\sigma_1 \sigma_{-1}$. We denote 
\[\gR^+[\pi_0^{-1}] \coloneqq  \gR^+  \cup^{E_k} \sigma_{-1}  \cup^{E_k} \eta.\]

\begin{lemma}
The map $B^{k}_\mr{M}(\gR^+) \to B^{k}_\mr{M}(\gR^+[\pi_0^{-1}])$ is a weak equivalence.
\end{lemma}

\begin{proof}Lemma \ref{lem:bkm-cell-attachment} tells us that up to homotopy $B^k_\mr{M}(\gR[\pi_0^{-1}])$ differs from $B^k_\mr{M}(\gR)$ by the attachment of a $k$- and a $(k+1)$-dimensional cell. An inspection of the attaching maps shows that these cells cancel, so that $B^k_\mr{M}(\gR) \to B^k_\mr{M}(\gR[\pi_0^{-1}])$ is a weak equivalence.
\end{proof}

We conclude that the map $\Omega^k B^k_\mr{M}(\gR^+) \to \Omega^k B^k_\mr{M}(\gR^+[\pi_0^{-1}])$ is also a weak equivalence. Since the monoid $\pi_0(\gR^+[\pi_0^{-1}]) \cong \bZ$ is in fact a group, we obtain the following from Remark \ref{rem:group-completion-grouplike}:

\begin{corollary}$\gR^+[\pi_0^{-1}]$ is weakly equivalent to $\Omega^k B^k_\mr{M} (\gR^+)$ as an $E_k$-algebra.
\end{corollary}
\section{Transferring vanishing lines}\label{sec:Transferring} 

In this section we explain under which circumstances vanishing lines for $E_k$-homology imply vanishing lines for $E_{k-1}$-homology (``transferring down'') or $E_{k+1}$-homology (``transferring up''). Transferring up using the bar constructions of the previous section, transferring down is the first application of our theory of $E_k$-cells. As before, $\sfC = \sfS^\sfG$ with $\sfS$ satisfying the axioms of Sections \ref{sec:axioms-of-cats} and \ref{sec:axioms-of-model-cats}.

\subsection{The bar spectral sequence}
The $k$-fold iterated bar construction $B^{E_k}(\gR,\epsilon)$ for an augmented $E_k$-algebra constructed in Section \ref{sec:iterated-bar-def} is a variation on the ordinary $k$-fold iterated bar construction. As such, there exists a bar spectral sequence\index{spectral sequence!bar} which computes $H_{*,*}(B^{E_k}(\gR,\epsilon))$ from $H_{*,*}(B^{E_{k-1}}(\gR,\epsilon))$. This will be used in Section \ref{sec:TrfUp} to show that if $\gR$ is an $E_k$-algebra whose $E_l$-homology for $l<k$ vanishes in a range of bidegrees, then the same is true for its $E_k$-homology. 

The setup is identical to that for the bar spectral sequence in Section \ref{sec:bar-ss}. Let $\cat{GrMod}_\bk$ denote the category of graded modules over a commutative ring $\bk$ with tensor product as usual involving a Koszul sign, and given a monoidal category $\sfG$ let $\cat{GrMod}_\bk^\sfG$ denote the category of functors with the Day convolution monoidal structure. As in Section \ref{sec:an-estimate}, we let $\bk[\bunit]$  denote the monoidal unit $(\bunit_\sfG)_*(\bk)$ in this category, given by the functor $g \mapsto \bk[\sfG(\bunit_\sfG,-)]$.

\begin{theorem}\label{thm:BarSS} Let $\gR$ be an augmented $E_k^+$-algebra which is cofibrant in $\sfC$. Then for each $\bk$-module $A$ there is a strongly convergent spectral sequence
	\[E^1_{g,p,q} = H_{g,q}(B^{E_{k-1}}(\gR,\epsilon)^{\otimes p};A) \Longrightarrow H_{g,p+q}(B^{E_k}(\gR,\epsilon);A)\]
	with differentials $d^r \colon E^r_{g,p,q} \to E^r_{g,p-r,q+r-1}$.

More generally, given a map $f \colon \gR \to \gS$ of augmented $E_k^+$-algebras which are cofibrant in $\sfC$ there is a strongly convergent spectral sequence
	\[\begin{tikzcd} E^1_{g,p,q} = H_{g,q}(B^{E_{k-1}}(\gS,\epsilon)^{\otimes p},B^{E_{k-1}}(\gR,\epsilon)^{\otimes p};A) \dar[Rightarrow] \\[-6pt] H_{g,p+q}(B^{E_k}(\gS,\epsilon),B^{E_k}(\gR,\epsilon);A).\end{tikzcd}\]
\end{theorem}

\begin{proof}There is a semi-simplicial object in $\sfC$
	\[X_\bullet \colon [p] \longmapsto \fgr{B^{E_k}_{p, \bullet \ldots, \bullet}(\gR,\epsilon)}\]
	given by forming the thick geometric realization in the last $(k-1)$ simplicial directions. This is levelwise cofibrant, by Lemma \ref{lem:thick-geom-rel-cofibrations}. By Lemma \ref{lem:bek-unit}, $X_0$ is weakly equivalent to $\bunit$. The object $X_1$ is isomorphic to the geometric realization of $\mathcal{P}_1(1) \times (B^{E_{k-1}}_{\bullet, \ldots, \bullet}(\gR,\epsilon))$, and using the fact that $\mathcal{P}_1(1)$ is contractible, we may conclude that this is weakly equivalent to $B^{E_{k-1}}(\gR,\epsilon)$. More generally there is a $(k-1)$-fold simplicial map
	\begin{equation}\label{eq:Segal}
	B^{E_k}_{p, \bullet \ldots, \bullet}(\gR,\epsilon) \lra \mathcal{P}_1(p) \times B^{E_{k-1}}_{\bullet, \ldots, \bullet}(\gR,\epsilon)^{\otimes p}
	\end{equation}
	induced by the inclusion
	\[\mathcal{P}_k(p_1, \ldots, p_k) \hookrightarrow \mathcal{P}_1(p_1) \times \mathcal{P}_{k-1}(p_2, \ldots, p_k)^{p_1}\]
	which remembers the grid inside each strip in the first simplicial direction. This map is a homotopy equivalence, so \eqref{eq:Segal} is a levelwise weak equivalence. As both objects are levelwise cofibrant, by Lemma \ref{lem:thick-geom-rel-cofibrations} we obtain an equivalence $X_p \simeq X_1^{\otimes p}$. The spectral sequence is then the geometric realization spectral sequence of Theorem \ref{thm:geom-rel-ss-thick} applied to the levelwise cofibrant simplicial object $X_\bullet$, using the equivalences $X_p \simeq X_1^{\otimes p}$ to identify the $E^1$-page.
	
	For the more general case, take the map of semi-simplicial objects
		\[\left([p] \mapsto \fgr{B^{E_k}_{p, \bullet \ldots, \bullet}(\gR,\epsilon)}\right) \lra[\left([p] \mapsto \fgr{B^{E_k}_{p, \bullet \ldots, \bullet}(\gS,\epsilon)}\right)\]
	induced by $f \colon \gR \to \gS$, and take the relative geometric realization spectral sequence.
\end{proof}
	
If $\bk$ is a field, or more generally if $H_{*,*}(B^{E_{k-1}}(\gR,\epsilon);\bk)$ consists of flat $\bk$-modules, there is an algebraic description of the $E_2$-page of the absolute version of the bar spectral sequence. To describe it, we use the Segal-like nature of the $k$-fold iterated bar construction in one of its $k$ directions to endow $H_{*,*}(B^{E_{k-1}}(\gR,\epsilon);\bk)$ with the structure of an augmented associative algebra in $\cat{GrMod}_\bk^\sfG$. 	
	
\begin{lemma} Let $\gR$ be as in Theorem \ref{thm:BarSS} and further suppose that $\sfG$ is a groupoid such that $G_x \times G_y \to G_{x \oplus y}$ is injective for all $x,y \in \sfG$. Then $H_{*,*}(B^{E_{k-1}}(\gR,\epsilon);\bk)$ has a natural structure of an augmented associative algebra in $\cat{GrMod}_\bk^\sfG$. 
	
If $H_{*,*}(B^{E_{k-1}}(\gR,\epsilon);\bk)$ consists of flat $\bk$-modules, then we may identify the $E_2$-page of the first spectral sequence in Theorem \ref{thm:BarSS} as
	\[E^2_{*, p,*} = \mathrm{Tor}_p^{H_{*,*}(B^{E_{k-1}}(\gR,\epsilon);\bk)}({\bk[\bunit]},{\bk[\bunit]}),\]
	with $\mr{Tor}$ formed in the category $\cat{GrMod}_\bk^\sfG$.
\end{lemma}	

\begin{proof}
	For the second part, when $\sfG$ is a groupoid such that $G_x \times G_y \to G_{x \oplus y}$ is injective for all $x,y \in \sfG$ we can apply edge homomorphisms of Lemma \ref{lem:KunnethFormula} (i) to get a map 
	\[H_{*,*}(X_1;\bk) \otimes H_{*,*}(X_1;\bk) \lra H_{*,*}(X_1^{\otimes 2};\bk) \cong H_{*,*}(X_2;\bk) \overset{(d_1)_*}\lra H_{*,*}(X_1;\bk),\]
	where $H_{*,*}(X_1;\bk) \cong H_{*,*}(B^{E_{k-1}}(\gR,\epsilon);\bk)$. This	induces a multiplication on $H_{*,*}(X_1;\bk)$, which may be seen to be associative by considering the face maps $X_3 \to X_1$. The two face maps $X_1 \to X_0$ are equal, and define an augmentation $H_{*,*}(X_1;\bk) \to H_{*,*}(X_0;\bk) =\bk[\bunit]$.
	
	When $H_{*,*}(B^{E_{k-1}}(\gR,\epsilon);\bk)$ consists of flat $\bk$-modules, the K\"unneth isomorphism of Lemma \ref{lem:KunnethFormula} (i) gives an isomorphism	\[H_{*,*}(X_p;\bk) \cong H_{*,*}(X_1^{\otimes p};\bk) \cong H_{*,*}(X_1;\bk)^{\otimes p}.\]
	In terms of this data, we make the identification
	\[H_{*,*}(X_p;\bk) \cong (\bunit_\sfG)_*(\bk) \otimes_{H_{*,*}(X_1;\bk)} (H_{*,*}(X_1;\bk)^{\otimes p+1} \otimes \bk[\bunit])\]
	and we recognize the chain complex $(E^1_{*,*,*}, d^1)$ as the result of applying the functor ${\bk[\bunit]} \otimes_{H_{*,*}(X_1;\bk)} - \colon \cat{GrMod}_\bk^\sfG \to \cat{GrMod}_\bk^\sfG$ levelwise to the canonical bar resolution of ${\bk[\bunit]}$ by free left $H_{*,*}(X_1;\bk)$-modules. This gives $E^2_{*, p, *} \cong \mathrm{Tor}_p^{H_{*,*}(X_1;\bk)}({\bk[\bunit]},{(\bk[\bunit]})$ as claimed.
\end{proof}

When $\gR$ is an augmented $E_{k+1}^+$-algebra, then it is in particular an augmented $E_k^+$-algebra so the second part of Theorem \ref{thm:BarSS} says that $H_{*,*}(B^{E_{k-1}}(\gR,\epsilon);\bk)$ is an augmented associative $\bk$-algebra. We shall shortly prove that this is in fact a \emph{commutative} one. (This should come as no surprise given Theorems \ref{thm:BarHomologyIndec} and \ref{thm:iterated-decomposables}, which endow the \emph{reduced} $E_{k-1}$-bar construction with an $E_2$-algebra structure.) Then we can combine the external tensor product with the multiplication map on $H_{*,*}(B^{E_{k-1}}(\gR,\epsilon);\bk)$ (which is a map of algebras if and only if it is commutative) to obtain a multiplication
	\[\begin{tikzcd}\mathrm{Tor}_*^{H_{*,*}(B^{E_{k-1}}(\gR,\epsilon);\bk)}({\bk[\bunit]},{\bk[\bunit]}) \otimes \mathrm{Tor}_*^{H_{*,*}(B^{E_{k-1}}(\gR,\epsilon);\bk)}({\bk[\bunit]},{\bk[\bunit]}) \dar \\ \mathrm{Tor}_*^{H_{*,*}(B^{E_{k-1}}(\gR,\epsilon);\bk) \otimes H_{*,*}(B^{E_{k-1}}(\gR,\epsilon);\bk)}({\bk[\bunit]} \otimes {\bk[\bunit]},{\bk[\bunit]}\otimes {\bk[\bunit]}) \dar \\   \mathrm{Tor}_*^{H_{*,*}(B^{E_{k-1}}(\gR,\epsilon);\bk)}({\bk[\bunit]},{\bk[\bunit]}),\end{tikzcd}\]
making $\mathrm{Tor}_*^{H_{*,*}(B^{E_{k-1}}(\gR,\epsilon);\bk)}({\bk[\bunit]},{\bk[\bunit]})$ a graded-commutative algebra with additional grading.
	
\begin{lemma}\label{lem:bar-ss-multiplicative} If $\gR$ is an augmented $E_{k+1}^+$-algebra, then $H_{*,*}(B^{E_{k-1}}(\gR,\epsilon);\bk)$ is an augmented commutative algebra. When $H_{*,*}(B^{E_{k-1}}(\gR,\epsilon);\bk)$ consists of flat $\bk$-modules, the first spectral sequence of Theorem \ref{thm:BarSS} is a spectral sequence of $\bk$-algebras.\end{lemma}

\begin{proof}
For the statement to be meaningful we must have $k-1 \geq 1$ and so $k+1 \geq 3$, which means that $\sfC$ must be symmetric monoidal. It shall be helpful to make the following general observation. For $r \leq k$, take the map
\begin{equation} \label{eqn:grid-c1} \mathcal{P}_k(p_1, \ldots, p_k) \lra \cC_r(p_1 \cdots p_r) \times \mathcal{P}_{k-r}(p_{r+1}, \ldots, p_k),\end{equation}
which considers the first $r$ grid directions as a collection of little $r$-cubes and remembers the remaining $(k-r)$ grid directions as a grid in an $(k-r)$-dimensional cube. We then define a $(k-r)$-fold semi-simplicial object $Y^{(k)}_{p_1,\ldots,p_r,\bullet,\ldots,\bullet}$ with $(p_{r+1},\ldots,p_k)$-simplices given by
\[\cC_r(p_1 \cdots p_r) \times \mathcal{P}_{k-r}(p_{r+1}, \ldots, p_k) \times G_{p_1,\ldots,p_k}(\epsilon),\]
with $G_{p_1,\ldots,p_k}(\epsilon)$ as in Definition \ref{def:kfold-bar-augmented}. As in that definition, the $i$th face map $d^j_i$ in the $j$th direction (here $r+1 \leq j \leq k$) is given by the face map of Definition \ref{def:pk} on the first factor. On the second factor, it is given by adjunction, by the map of simplicial sets
\begin{align*}
\cC_r(p_1 \cdots p_r) \times &\cP_{k-r}(p_{r+1},\ldots,p_k) \lra \cC_k(2) \\
&\qquad \overset{\alpha}\lra \mr{Map}_\sfC(G_{p_1, \ldots, p_k}(\epsilon), G_{p_1,\ldots, p_{j-1}, p_j-1, p_{j+1},\ldots,p_k}(\epsilon))
\end{align*}
with the first map given by $\{e\} \times \{t^j_i\} \mapsto \mr{id}_{I^r} \times \delta^j_i$, and the second map as in Definition \ref{def:kfold-bar-augmented}.

We next describe a $(k-r)$-fold simplicial map 
\begin{equation}\label{eqn:multiplication-ykr} Y^{(k)}_{p_1,\ldots,p_r,\bullet,\ldots,\bullet} \lra B^{E_{k-r}}_{\bullet,\ldots,\bullet}(\gR,\epsilon).\end{equation} 
On the first factor this is simply the projection $\cC_r(p_1 \cdots p_r) \times \mathcal{P}_{k-r}(p_{r+1}, \ldots, p_k) \to \mathcal{P}_{k-r}(p_{r+1}, \ldots, p_k)$. On the second factor, it is given by adjunction, by the map of simplicial sets
\begin{align*}
\cC_r(p_1 \cdots p_r) \times \cP_{k-r}(p_{r+1},\ldots,p_k) &\lra \cC_k(p_1 \cdots p_r) \overset{\beta}\lra \mr{Map}_\sfC(G_{p_1,\ldots, p_k}(\epsilon), G_{p_{r+1},\ldots,p_k}(\epsilon))\\
\{e\} \times \{t^j_i\} &\longmapsto \{e \times \mr{id}_{I^{k-r}}\},
\end{align*}
with $\beta$ given as follows: as long as $1 \leq q_j \leq p_j$ for all $j \geq r+1$ by the map
\begin{align*}\cC_k(p_1 \cdots p_r) \lra \cE_{\gR}(p) &= \mr{Map}_\sfC(\gR^{\otimes p_1 \cdots p_r}, \gR) \\
&\qquad = \mr{Map}_\sfC\left(\bigotimes_{j=1}^r \bigotimes_{i_j=1}^{p_j} B_{p_1,\ldots,p_r,p_{r+1},\ldots, p_k}^{i_1,\ldots,i_r,q_{r+1},\ldots,q_k}, B_{p_{r+1},\ldots,p_k}^{q_{r+1},\ldots,q_k}\right)\end{align*}
and the evident identity maps on the remaining factors. If for some $j \geq r+1$, $q_j$ is either $0$ or $p_j+1$, it is the same map but with $\gR$ replaced by $\bunit$.

We shall augment notation from the proof of Theorem \ref{thm:BarSS} to make the dependence on $k$ clear: $X^{(k)}_\bullet \coloneqq \fgr{B^{E_{k}}_{p,\bullet,\ldots,\bullet}(\gR,\epsilon)}$. Since $\gR$ is an $E_{k+1}^+$-algebra, we may consider the $E_{k+1}$-bar construction. We set $r=2$, $p_1=1$, and $p_2=2$, then take the geometric realization of \eqref{eqn:grid-c1} and \eqref{eqn:multiplication-ykr} to obtain maps
\[\fgr{B^{E_{k+1}}_{1,2,\bullet,\ldots,\bullet}(\gR,\epsilon)} \lra \fgr{Y^{(k+1)}_{1,2,\bullet,\ldots,\bullet}} \lra B^{E_{k-1}}(\gR,\epsilon).\]
The multiplication on $X_1^{(k)}$ can be recovered from this. To do so, we use the evident homotopy equivalences 
\begin{align*}\cP_k(2,p_3,\ldots,p_{k+1}) &\lra \cP_1(2) \times \cP_{k-1}(p_3,\ldots,p_{k+1})^2 \\ 
\cP_{k+1}(1,2,p_3,\ldots,p_{k+1}) &\lra \cP_{k}(2,p_3,\ldots,p_{k+1}), \\
\cC_2(1 \cdot 2) \times \cP_{k-1}(p_3,\ldots,p_{k+1}) &\lra \cC_2(2) \times \cP_{k-1}(p_3,\ldots,p_{k+1})^2, \\
\cP_{k}(1,p_3,\ldots,p_{k+1}) &\lra \cP_{k-1}(p_3,\ldots,p_{k+1}),\end{align*} 
to obtain the weak equivalences in the following homotopy-commutative diagram
\[\begin{tikzcd}\cP_1(2) \times (X^{(k)}_1)^{\otimes 2} \arrow[dd] & X^{(k)}_2 \rar{d_1} \lar[swap]{\simeq} & X^{(k)}_1 \arrow[equals]{d} \\
& \uar{\simeq}  \fgr{B^{E_{k+1}}_{1,2,\bullet,\ldots,\bullet}(\gR,\epsilon)} \dar & X^{(k)}_1 \dar{\simeq} \\
\cC_2(2) \times (X^{(k)}_1)^{\otimes 2} & \lar[swap]{\simeq}  Y^{k}_{1,2}  \rar &  B^{E_{k-1}}(\gR,\epsilon).\end{tikzcd}\]
Proving that the both squares commute up to homotopy is a simple matter of tracing through the various maps of grids and cubes.
 
Thus we have exhibited the multiplication on $X_1^{(k)}$ up to weak equivalence as arising from a choice of point in $\cC_2(2)$. It now remains to observe that the multiplication in reverse order similarly arises by picking another point in $\cC_2(2)$, and since $\cC_2(2)$ is path-connected these maps are homotopic.

\vspace{.5 em}

We showed in Theorem \ref{thm:BarSS} that there is a weak equivalence and multiplication
\[(X_1^{(k+1)})^{\otimes 2} \overset{\simeq}\longleftarrow X_2^{(k+1)} \lra X_1^{(k+1)}.\]
As $X_1^{(k+1)} \simeq B^{E_k}(\gR,\epsilon)$, it is this zigzag of maps that endows the bar spectral sequence with an algebra structure, which by construction converges to the $\bk$-algebra structure on $H_{*,*}(B^{E_{k}}(\gR,\epsilon);\bk)$. On the $E^1$-page it gives the map on canonical bar resolutions induced by the $E_1$-algebra structure in the remaining direction. This is homotopic to the $E_1$-algebra structure used to construct the product in the second part of Theorem \ref{thm:BarSS}, and hence gives the $\bk$-algebra structure on $\mr{Tor}$-groups discussed above.
\end{proof}

\subsection{Transferring vanishing lines up}\label{sec:TrfUp} Transferring vanishes lines up follows from our expression of derived $E_k$-indecomposables in terms of the iterated bar construction, the bar spectral sequence described in Theorem \ref{thm:BarSS}, and a K{\"u}nneth-type theorem.

\begin{theorem}\label{thm:TrfUp}
Let $\gR \in \Alg_{E_k}(\sfC)$, and $\rho \colon \sfG \to [-\infty,\infty]_\geq$ be an abstract connectivity such that $\rho*\rho \geq \rho$. If $l \leq k$ is such that $H^{E_l}_{g,d}(\gR)=0$ for $d < \rho(g)-l$, then $H^{E_k}_{g,d}(\gR)=0$ for $d < \rho(g)-l$ too.
\end{theorem}

\begin{proof}
	We claim that it is enough to consider the case $(l,k) = (k-1,k)$. To prove this claim we need to explain how the case $(l,l+1)$ provides the input for $(l+1,l+2)$, etc. We can use the case $(l,l+1)$ to prove that if $H^{E_l}_{g,d}(\gR)=0$ for $d < \rho(g)-l$, then $\smash{H^{E_{l+1}}_{g,d}}(\gR)=0$ for $d < \rho(g)-l$ too. This conclusion provides input for the case $(l+1,l+2)$ when we rewrite it as $H^{E_{l+1}}_{g,d}(\gR)=0$ for $d < \rho'(g)-(l+1)$ with $\rho' \coloneqq \rho+1$, which still satisfies $\rho' \ast \rho' \geq \rho'$.
	
	Let us from now assume that $l = k-1$. By Theorem \ref{thm:BarHomologyIndec}, we have equivalences
	\[\tilde{B}^{E_{k-1}}(\gR) \simeq S^{k-1} \wedge Q_\bL^{E_{k-1}}(\gR) \quad \text{and} \quad \tilde{B}^{E_k}(\gR) \simeq  S^k \wedge Q_\bL^{E_k}(\gR),\]
	so the assumption of the theorem is equivalent to saying that $\tilde{B}^{E_{k-1}}(\gR)$ is homologically $\rho$-connective, and our desired conclusion is equivalent to saying that $\tilde{B}^{E_{k}}(\gR)$ is homologically $(1+\rho)$-connective. We also have, by definition, homotopy cofibre sequences in $\sfC_*$
	\begin{align*}
		{B}^{E_{k-1}}(\bunit, \epsilon_\bunit)_+ \lra &{B}^{E_{k-1}}(\gR^+, \epsilon_{can})_+\lra \tilde{B}^{E_{k-1}}(\gR),\\
		{B}^{E_{k}}(\bunit, \epsilon_\bunit)_+ \lra &{B}^{E_{k}}(\gR^+,\epsilon_{can})_+\lra \tilde{B}^{E_{k}}(\gR).
	\end{align*}
	Let us write $\epsilon$ for either of the augmentations $\epsilon_\bunit$ or $\epsilon_{can}$.
	
	The relative version of the bar spectral sequence of Theorem \ref{thm:BarSS} starts form
	\[E^1_{g,p,q} = H_{g,q}(B^{E_{k-1}}(\gR^+,\epsilon)^{\otimes p}, {B}^{E_{k-1}}(\bunit,\epsilon)^{\otimes p})\]
	and converges strongly to $H_{g, p+q}(B^{E_{k}}(\gR^+,\epsilon), {B}^{E_{k}}(\bunit,\epsilon))=H_{g,p+q}(\tilde{B}^{E_{k}}(\gR))$. The assumption can be rephrased as saying that the map $B^{E_{k-1}}(\bunit,\epsilon) \to B^{E_{k-1}}(\gR^+,\epsilon)$ is homologically $\rho$-connective. As ${B}^{E_{k-1}}(\bunit,\epsilon) \simeq \bunit$ by Lemma \ref{lem:bek-unit}, which has homological connectivity given by the unit $\bunit_\text{conn} \in [-\infty,\infty]_\geq^\sfG$ as in (\ref{eqn:abstract-connectivity-unit}), the object $B^{E_{k-1}}(\gR^+,\epsilon)$ is $\inf(\bunit_\text{conn}, \rho)$-connective. By Corollary \ref{cor:connectivity-under-tensor2} the map ${B}^{E_{k-1}}(\bunit,\epsilon)^{\otimes p} \to {B}^{E_{k-1}}(\gR^+,\epsilon)^{\otimes p}$ is then homologically $(\inf(\bunit_\text{conn}, \rho)^{* p-1} * \rho)$-connective, and hence $\rho$-connective using the fact that $\rho*\rho \geq \rho$. Furthermore, it is $\infty$-connective if $p=0$, so $E^1_{g,p,q}$ vanishes if $p=0$ or if $q < \rho(g)$, so it vanishes for $p+q < 1+\rho(g)$. As this spectral sequence converges strongly to $H_{g, p+q}(\tilde{B}^{E_{k}}(\gR)) \cong H_{g,p+q-k}^{E_k}(\gR)$, the conclusion follows.
\end{proof}

The absolute case serves as input for the following relative version:

\begin{proposition}Let $f \colon \gR \to \gS$ be a map of $E_k$-algebras, and $\rho,\sigma \colon \sfG \to [-\infty,\infty]_\geq$ be abstract connectivities such that $\rho*\rho \geq \rho$ and $\rho*\sigma \geq \sigma$. If $l \leq k$ is such that $H^{E_l}_{g,d}(\gR)=0 = H^{E_l}_{g,d}(\gS)$ for $d < \rho(g)-l$ and $H^{E_l}_{g,d}(\gS,\gR) = 0$ for $d<\sigma(g)$, then $H^{E_k}_{g,d}(\gS, \gR)=0$ for $d < \sigma(g)$ too.
\end{proposition}

\begin{proof}
	By the same reasoning as in the proof of Theorem \ref{thm:TrfUp}, it suffices to consider the case $(l,k) = (k-1,k)$. As above, let us write $\epsilon$ for the canonical augmentation of $\gR^+$ or $\gS^+$. The relative version of the bar spectral sequence of Theorem \ref{thm:BarSS} starts from
	\[E^1_{g,p,q} = H_{g,q}(B^{E_{k-1}}(\gS^+,\epsilon)^{\otimes p}, {B}^{E_{k-1}}(\gR^+,\epsilon)^{\otimes p})\]
	and converges strongly to $H_{g, p+q}(B^{E_{k}}(\gS^+,\epsilon), {B}^{E_{k}}(\gR,^+\epsilon))$. Theorem \ref{thm:BarHomologyIndec} and Lemma \ref{lemma:bek-cofibrant} give that $B^{E_{k-1}}(\gR^+,\epsilon) \simeq \bk \oplus S^{k-1} \wedge Q^{E_{k-1}}_\bL(\gR)$ and similarly for $\gS$, so the map $B^{E_{k-1}}(\gR^+,\epsilon) \to B^{E_{k-1}}(\gS^+,\epsilon)$ is $(\sigma+k-1)$-connective. Its domain and target are $\inf(\bunit_\mr{conn},\rho)$-connective by Theorem \ref{thm:TrfUp}. 
	
	By Corollary \ref{cor:connectivity-under-tensor2} the map ${B}^{E_{k-1}}(\gR^+,\epsilon)^{\otimes p} \to {B}^{E_{k-1}}(\gS^+,\epsilon)^{\otimes p}$ is then homologically $(\inf(\bunit_\text{conn}, \rho)^{*p-1} * (\sigma+k-1))$-connective, and hence $(\sigma+k-1)$-connective. Furthermore, it is $\infty$-connective if $p=0$, so $E^1_{g,p,q}$ vanishes if $p=0$ or if $q < \sigma(g)+k-1$, so it vanishes for $p+q < \sigma(g) + k$. As this spectral sequence converges strongly to $H_{g, p+q}(B^{E_{k}}(\gS^+,\epsilon), {B}^{E_{k}}(\gR^+,\epsilon)) \cong H_{g,p+q-k}^{E_k}(\gS, \gR)$, the conclusion follows.
\end{proof}

\subsection{Transferring vanishing lines down} To transfer vanishing lines down, we use the theory of CW approximation that we have developed in Section \ref{sec:additive-case} and so we must assume Axiom \ref{axiom:Hurewicz}.

\begin{theorem}\label{thm:TrfDown}
Suppose that $\sfG$ is Artinian, let $\gR \in \Alg_{E_k}(\sfC)$ be reduced and 0-connective, $l \leq k$, and $\rho \colon \sfG \to [-\infty,\infty]_\geq$ be an abstract connectivity such that $\rho*\rho \geq \rho$ and $H^{E_k}_{g,d}(\gR)=0$ for $d < \rho(g) - l$. Then $H^{E_l}_{g,d}(\gR)=0$ for $d < \rho(g)-l$.
\end{theorem}

\begin{proof}
Firstly, the groupoid $\sfG$ and the operad $E_k$ satisfy the hypotheses of Lemma \ref{lem:TensorDetectNull}. The canonical morphism $\binit \to \gR$ is between 0-connective reduced $E_k$-algebras, so by Theorem \ref{thm:MinCellStr-additive} we may construct a CW approximation $\gZ \overset{\sim}\to \gR$, where $\gZ$ consists of $(g,d)$-cells with $d \geq \rho(g)-l$ and has skeletal filtration $\mr{sk}(\gZ) \in \Alg_{E_k}(\sfC^{\bZ_{\leq}})$. By Theorem \ref{thm:associated-graded-skeletal}, the associated graded $\grr(\mr{sk}(\gZ))$ of this filtration is given by $\gE_k(X)$, where $X$ is a wedge of $d_* S^{n, d}$'s with $d \geq \rho(g)-l$. The spectral sequence of Theorem \ref{thm:DerIndecSS} with $\cO = E_l$ takes the form
\[E^1_{g,p,q} = H_{g,p+q, p}^{E_l}(\gE_k(X)) \Longrightarrow H_{g,p+q}^{E_l}(\gR),\]
so, forgetting the internal grading, it is enough to show the vanishing of $H_{g,d}^{E_l}(\gE_k(X))$ for $d < \rho(g)-l$. To do this we use the weak equivalences
\[S^l \wedge Q_\bL^{E_l}(\gE_k(X)) \simeq \tilde{B}^{E_l}(\gE_k(X)) \simeq E_{k-l}(S^l \wedge X_+)\]
in $\sfC_*$ from Theorems \ref{thm:BarHomologyIndec} and \ref{thm:CalcFree}, so that it suffices to show that $E_{k-l}(S^l \wedge X_+)$ is homologically $\rho$-connective.

We have that $S^l \wedge X_+$ is homologically $\rho$-connective, so it follows from Lemma \ref{lem:connectivity-and-tensor-products} (i) that $(S^l \wedge X_+)^{\otimes p}$ is homologically $\rho^{* p}$-connective, so $\rho$-connective (as $\rho$ is lax monoidal). If $\sfC$ is $\infty$-monoidal, then we have
\[E_{k-l}(S^l \wedge X_+) = \bigvee_{i \geq 1} \cC_{k-l}(p) \times_{\fS_p} (E_{k-l}(S^l \wedge X_+))^{\otimes p},\]
and it follows from the homotopy orbit spectral sequence as in Section \ref{sec:homotopy-orbit-ss} that this is also homologically $\rho$-connective, as required. If $\sfC$ is $2$-monoidal, one needs to replace $\cC_{k-l}(p)$ by $\cC^{\cat{FB}_2}_{k-l}(p)$ and the symmetric group $\fS_p$ by the braid group $\beta_p$.
\end{proof}

By doing a more careful analysis we can occasionally relax the condition $d < \rho(g)-l$; we give the following theorem as an example of a general type of argument. 

\begin{theorem}\label{thm:TrfDownExplicit}
Let $\gR \in \Alg_{E_\infty}(\cat{sMod}_\bQ^\bN)$ be a reduced $E_\infty$-algebra in $\bN$-graded simplicial $\bQ$-modules such that $H_{*,0}(\gR^+) = \bQ[\sigma]$ with $\vert \sigma \vert = (1,0)$. If $H^{E_k}_{g,d}(\gR)=0$ for $d < 2(g-1)$ then $H^{E_1}_{g,d}(\gR)=0$ for $d < \frac{3}{2}(g-1)$.
\end{theorem}

This does not follow from Theorem \ref{thm:TrfDown}, as the assumed vanishing range for $E_k$-homology is $d < (2g-1)-1$, but $\rho(g) = 2g-1$ does not satisfy $\rho*\rho \geq \rho$.

\begin{proof}
Firstly, by transferring vanishing lines up we may suppose that $H^{E_\infty}_{g,d}(\gR)=0$ for $d < 2(g-1)$. As in the proof of Theorem \ref{thm:TrfDown}, by filtering a suitable CW approximation of $\gR$ we can reduce to the case $\gR = \gE_\infty(X)$ with $X$ a wedge of spheres such that $H_{g,d}(X)=0$ for $d < 2(g-1)$ and $H_{1,0}(X) = \bQ\{\sigma\}$.

We use the equivalences
\[S^1 \wedge Q^{E_1}_\bL(\gE_\infty(X))) \simeq \tilde{B}^{E_1}(\gE_\infty(X)) \simeq E_\infty(S^1 \wedge X)\]
from Theorems \ref{thm:BarHomologyIndec} and \ref{thm:CalcFree}. By F.\ Cohen's computations of the homology of free $E_k$-algebras, which shall be explained in Section \ref{sec:Cohen}, we have
\[H_{*,*}(E_\infty(S^1 \wedge X)) \cong \Lambda_\bQ(H_{*,*}(S^1 \wedge X)),\]
the free graded-commutative algebra on the rational homology of the suspension $S^1 \wedge X$, which we may write as $\Lambda_\bQ(\bQ\{s \sigma\}) \otimes A$ where $s\sigma$ has bidegree $(1,1)$ and $A$ is a free graded-commutative algebra on generators all of which have slope $\frac{d}{g} \geq \tfrac{3}{2}$, so $A$ is trivial in bidegrees $(g,d)$ with $d < \tfrac{3}{2}g$. It follows that $\smash{H^{E_1}_{g,d}}(\gR; \bQ)=0$ for $d < \frac{3}{2}(g-1)$ as required.
\end{proof}

\begin{remark}In fact, we may also do similar analyses with $\bF_p$-coefficients. Then we obtain the same vanishing range $H^{E_1}_{g,d}(\gR)=0$ for $d < \frac{3}{2}(g-1)$ as long as $p \geq 5$, and a lower range $d < \frac{4}{3}(g-1)$ for $p=3$. For $p=2$ one would need to know more information about the cell structures to improve upon Theorem \ref{thm:TrfDownExplicit}. To see this done in an example, see Theorem 10.2 of \cite{e2cellsIV}.\end{remark}
\section{Comparing algebra and module cells}\label{sec:comparing-algebra-and-module cells}

Given a map $f \colon \gR \to \gS$ of $E_k$-algebras, we may consider it as a map of $E_1$-algebras by neglect of structure and using the constructions of Section \ref{sec:modules} obtain a map $\overline{f} \colon \overline{\gR} \to \overline{\gS}$ of unital associative algebras. In this section we compare $H^{E_k}_{*,*}(\gS, \gR)$ (which measures the $E_k$-algebra cells necessary to construct $\gS$ from $\gR$) with $H^{\overline{\gR}}_{*,*}(\overline{\gS})$ (which measures the $\overline{\gR}$-module cells necessary to construct $\overline{\gS}$).

\subsection{Coproducts of $E_k$-algebras} In this section we discuss coproducts of unital $E_k$-algebras; we do so in general when $k=\infty$ and for free ones when $k<\infty$.

\subsubsection{Coproducts of $E_\infty$-algebras}
The coproduct of unital $E_\infty$-algebras is easy to describe. Recall from Proposition \ref{prop:algcc-symmetric-monoidal} that when $\sfC$ is symmetric monoidal and $\cC$ is an operad given by a $\infty$-symmetric sequence in $\cat{sSet}$, the category $\Alg_\cC(\sfC^{\sfG})$ inherits a symmetric monoidal structure $\otimes_\cC$ compatible with the forgetful functor, in the sense that there is a natural isomorphism $U^\cC(\gR \otimes_{\cC} \gS) \cong U^\cC(\gR) \otimes U^\cC(\gS)$.

Let $\gR$ and $\gS$ be $E_\infty^+$-algebras. The unit maps $\bunit \to \gR$ and $\bunit \to \gS$ are maps of $E_\infty^+$-algebras. Since there is a natural isomorphism $\gR \otimes_{E^+_\infty} \bunit \cong \gR$ of $E^+_\infty$-algebras, they induce maps
\[\gR \lra \gR \otimes_{E^+_\infty} \gS \longleftarrow \gS\]
and hence give a natural map $\gR \sqcup^{E_\infty^+} \gS \to \gR \otimes_{E^+_\infty} \gS$ of $E^+_\infty$-algebras. 

\begin{proposition}\label{prop:einty-coporoduct}
Let $\gR$ and $\gS$ be $E_\infty^+$-algebras in $\sfC= \sfS^\sfG$ for any $\sfG$, then
	\[\gR \sqcup^{\bL,E_\infty^+} \gS \lra \gR \otimes^\bL_{E^+_\infty} \gS\]
is a weak equivalence. Furthermore, if $\gR$ and $\gS$ are cofibrant in $\sfC$, then $\gR \otimes^\bL_{E^+_\infty} \gS \to \gR \otimes_{E^+_\infty} \gS$ is a weak equivalence.
\end{proposition}

\begin{proof}
For the first part, we must show that $\gR \sqcup^{E_\infty^+} \gS \to \gR \otimes_{E^+_\infty} \gS$ is a weak equivalence if $\gR$ and $\gS$ are cofibrant in $\Alg_{E^+_\infty}(\sfC)$.
	
We first prove this when $\gR$ and $\gS$ are free $E^+_\infty$-algebras on cofibrant objects in $\sfC$. The map in question is induced by a map of right $E^+_\infty$-module functors $\sfC^2 \to \sfC$
	\[\gE^+_\infty(-) \sqcup^{E^+_\infty} \gE^+_\infty(-) \cong \gE^+_\infty(- \sqcup -) \Longrightarrow \gE^+_\infty(-) \otimes_{E^+_\infty} \gE^+_\infty(-),\]
which we claim is a weak equivalence when evaluated on cofibrant objects of $\sfC$. We may verify this on underlying objects, where the natural transformation becomes a map
\[E^+_\infty(- \sqcup -) \Longrightarrow E^+_\infty(-) \otimes E^+_\infty(-)\]
which we shall now describe. First, note that $E^+_\infty(X \sqcup Y)$ may be described as the colimit as $k \to \infty$ of a coproduct over $n$ and $n'$ of terms
\[\mr{Emb}^\mr{rect}(\sqcup_{n+n'} I^k,I^k) \times_{\fS_{n} \times \fS_{n'}} X^{\otimes n} \otimes Y^{\otimes n'},\]
while $E^+_\infty(X) \otimes E^+_\infty(Y)$ may be described as the colimit as $k \to \infty$ of a coproduct over $n$ and $n'$ of terms
\[\left(\mr{Emb}^\mr{rect}(\sqcup_n I^k,I^k) \times \mr{Emb}^\mr{rect}(\sqcup_{n'} I^k,I^k)\right) \times_{\fS_{n} \times \fS_{n'}} X^{\otimes n} \otimes Y^{\otimes n'}.\]
The natural transformation is induced by the $\fS_n \times \fS_{n'}$-equivariant restriction map
\[\begin{tikzcd}\colim\limits_{k \to \infty} \mr{Emb}^\mr{rect}(\sqcup_{n+n'} I^k,I^k) \dar \\[-5pt] \colim\limits_{k \to \infty} \left(\mr{Emb}^\mr{rect}(\sqcup_n I^k,I^k) \times \mr{Emb}^\mr{rect}(\sqcup_{n'} I^k,I^k)\right),\end{tikzcd}\]
which is a weak equivalence (because both are contractible) between free $\fS_n \times \fS_{n'}$-spaces. If $X$ and $Y$ are cofibrant then $- \times_{\fS_{n} \times \fS_{n'}} X^{\otimes n} \otimes Y^{\otimes n'}$ is a left Quillen functor so preserves weak equivalences between free $\fS_n \times \fS_{n'}$-spaces, as required.

\vspace{1em}

Next we discuss the case that $\gR$ and $\gS$ are cofibrant in $\Alg_{E^\infty_+}(\sfC)$. In that case, we take the thick monadic bar resolutions to obtain functorial free simplicial resolutions $\gR_\bullet = \sigma_* B_\bullet(E_\infty^+,E_\infty^+,\gR)$ and $\gS_\bullet = \sigma_* B_\bullet(E_\infty^+,E_\infty^+,\gS)$, as $\gR$ and $\gS$ are in particular cofibrant in $\sfC$.

Then the levelwise coproduct $\gR_\bullet \sqcup^{E^+_\infty} \gS_\bullet$ is a bisimplicial object augmented over $\gR \sqcup^{E^+_\infty} \gS$, and the levelwise tensor product $\gR_\bullet \otimes_{E^+_\infty} \gS_\bullet$ is a bisimplicial object augmented over $\gR \otimes_{E^+_\infty} \gS$. We now remark that the geometric realization of a bisimplicial object may be obtained in two steps. Thus consider first for fixed $q$ the commutative diagram
\[\begin{tikzcd} {\gr{[p] \mapsto \gR_p \sqcup^{E^+_\infty} \gS_q}}_{E^+_\infty} \dar \rar & \gR \sqcup^{E^+_\infty} \gS_q \dar \\
{\gr{[p] \mapsto \gR_p \otimes_{E^+_\infty} \gS_q}}_{E^+_\infty} \rar & \gR \otimes_{E^+_\infty} \gS_q.\end{tikzcd}\]

To prove that the top horizontal map is a weak equivalence, we use that $- \sqcup^{E_\infty^+} \gS_q$, being a colimit, commutes with geometric realization. Thus it suffices to prove that $\gr{\gR_\bullet}_{E^+_\infty} \sqcup^{E_\infty^+} \gS_q \to \gR \sqcup^{E_\infty^+} \gS_q$ is a weak equivalence.  Since $\gS_q$ is cofibrant in $\Alg_{E^+_\infty}(\sfC)$, the functor $- \sqcup^{E_\infty^+} \gS_q$ preserves weak equivalences between cofibrant objects. That $\gr{\gR_\bullet}_{E^+_\infty}$ is cofibrant follows from Lemma \ref{lem:geom-rel-cofibrations}, and we assumed that $\gR$ is cofibrant.

There is a map $\gr{\gR_\bullet \otimes_{E_\infty^+} \gS_q}_{E^+_\infty} \to \gr{\gR_\bullet}_{E^+_\infty} \otimes_{E_\infty^+} \gS_q$ which is an isomorphism because geometric realization in $E_k^+$-algebras or in the underlying category coincide by part (iii) of Lemma \ref{lem:monad-functors-geomrel}, $U^{E_\infty^+}(-\otimes_{E_\infty^+} \gS_q) \cong U^{E_\infty^+}(-) \otimes U^{E_\infty^+}(\gS_q)$, and $-\otimes-$ commutes with geometric realization in each entry. Thus it suffices to prove that $\gr{\gR_\bullet}_{E^+_\infty} \otimes_{E_\infty^+} \gS_q \to \gR \otimes_{E_\infty^+} \gS_q$ is a weak equivalence. By the above observations it suffices to note that $- \otimes U^{E^+_\infty} \gS_q$ preserves weak equivalences between cofibrant objects.

The left vertical map is a levelwise weak equivalence by the discussion for free $E_\infty^+$-algebras, and we claim that both simplicial objects are Reedy cofibrant. This follows since $- \sqcup^{E_\infty^+} \gS_q$ and $-\otimes_{E_\infty^+} \gS_q$ commute with the formation of latching objects and preserve cofibrations between cofibrant objects since $\gS_q$ is cofibrant. Thus by Lemma \ref{lem:geom-rel-cofibrations}, the left vertical map is a weak equivalence between cofibrant objects.

A similar argument with $\gR$ taking the role of $\gS_q$ implies that $\gR \sqcup^{E^+_\infty} \gS \to \gR \otimes_{E^+_\infty} \gS$ is a weak equivalence.

\vspace{1em}

Finally, it remains to observe that the above argument tells us that the map
\[\gR \otimes_{E^+_\infty}^\bL \gS \to \gR \otimes_{E^+_\infty} \gS\]
is a weak equivalence if $\gR$ and $\gS$ are cofibrant in $\sfC$; whenever we discussed $\smash{\otimes_{E^+_\infty}}$ we applied $\smash{U^{E^+_\infty}}$ and argued only using cofibrancy in $\sfC$.
\end{proof}

Let us specialize this to the case that $\gR$ and $\gS$ are both free $E_\infty^+$-algebras.

\begin{corollary}\label{cor:CoprodEinfty}
For $A, B \in \sfC$ cofibrant, 0-connective and reduced there is a natural weak equivalence 
\[\overline{\gE_\infty(A \sqcup B)} \simeq \overline{\gE_\infty(A)} \otimes E_\infty^+( B)\]
of left $\overline{\gE_\infty(A)}$-modules.
\end{corollary}
\begin{proof}
The composition $E_\infty^+( B) \to E_\infty^+(A \sqcup B) \to \overline{\gE_\infty(A \sqcup B)}$ has target a left $\overline{\gE_\infty(A)}$-module, so extends to a morphism
\[\overline{\gE_\infty(A)} \otimes E_\infty^+( B) \lra \overline{\gE_\infty(A \sqcup B)}\]
of $\overline{\gE_\infty(A)}$-modules. To verify that it is an equivalence we may neglect the module structure, whereupon it is proven as in the first part of the proof of Proposition \ref{prop:einty-coporoduct}.\end{proof}

\subsubsection{Coproducts of $E_k$-algebras, $k < \infty$}

For unital $E_k$-algebras with $k < \infty$ there is not such a simple description of the coproduct as Proposition \ref{prop:einty-coporoduct}, but we do nonetheless have the following analogue of Corollary \ref{cor:CoprodEinfty}, which for many purposes suffices. The proof will use homology to detect weak equivalences, so we must assume that $\sfS$ satisfies Axiom \ref{axiom:Hurewicz}, and in particular that it is pointed.

\begin{theorem}\label{thm:AlgAsMod}
Suppose that $\sfS$ satisfies Axiom \ref{axiom:Hurewicz} and $\sfG$ is $(k+1)$-monoidal as well as Artinian. Then for $A, B \in \sfC$ cofibrant, 0-connective and reduced there is a natural weak equivalence 
\[\overline{\gE_k(A \vee B)} \simeq \overline{\gE_k(A)} \otimes E_k^+(E_1^+(S^{k-1} \wedge A) \otimes B)\]
of left $\overline{\gE_k(A)}$-modules.
\end{theorem}

\begin{remark}\label{rem:AlgAsModComputational}
When $\sfG$ is discrete, this theorem can proven at the level of homology by direct calculation, using the description (due to F.\ Cohen) of the homology of free $E_k$-algebras which we shall recall in Section \ref{sec:Cohen}. This is enough to prove Theorem \ref{thm:EkModuleHurewicz} and Corollary \ref{cor:BarEstimate} below in the case that $\sfG$ is discrete, which suffices for most applications. Rather than follow the admittedly difficult proof of Theorem \ref{thm:AlgAsMod} the reader may wish to go ahead to Section \ref{sec:Cohen}, after which they will be able to prove Theorem \ref{thm:AlgAsMod} at the level of homology themselves.
\end{remark}

\begin{remark}
The assumption that $\sfG$ be $(k+1)$-monoidal, and not merely $k$-monoidal as suffices to define $E_k$-algebras, may seem surprising but is in fact necessary for the statement of Theorem \ref{thm:AlgAsMod} to hold. This may be seen by taking $\sfG$ to be the free braided monoidal category on one generator $X$, $\sfS = \mathsf{sMod}_\bk$, and $A = B = X_*(\bk)$, then directly calculating both sides.
\end{remark}

We shall return to the proof of this theorem in Section \ref{sec:proofAlgAsMod}.

\subsection{Comparison in the $E_\infty$ case}

In the case of $E_\infty$-algebras, the comparison results we have in mind can be deduced from the following more precise theorem.

\begin{theorem}\label{thm:e2cellsIV-ss}
	Suppose that $\sfG$ is Artinian. Let $f \colon \gR \to \gS$ be a morphism in $\Alg_{E_\infty}(\sfC)$ between cofibrant, $0$-connective and reduced objects. Then there is a strongly convergent spectral sequence
	\[E^1_{g,p,q} = \widetilde{H}_{g,p+q,p}(E_\infty^+((-1)_* Q^{E_\infty}_\bL(\gS)/Q^{E_\infty}_\bL(\gR));A) \Longrightarrow H_{g,p+q}^{\overline{\gR}}(\overline{\gS};A)\]
with differentials $d^r \colon E^r_{g,p,q} \to E^r_{g,p-r,q+r-1}$.
\end{theorem}

To interpret this as a comparison result, suppose that $\sigma \colon \sfG \to [-\infty,\infty]_{\geq}$ is an abstract connectivity such that $\sigma * \sigma \geq \sigma$, and that $H_{g,d}^{E_\infty}(\gS, \gR)=0$ for $d < \sigma(g)$. That is, $Q^{E_\infty}_\bL(\gS)/Q^{E_\infty}_\bL(\gR)$ is homologically $\sigma$-connective. Then it follows from Lemma \ref{lem:connectivity-and-tensor-products} (i) that $(Q^{E_\infty}_\bL(\gS)/Q^{E_\infty}_\bL(\gR))^{\otimes p}$ is homologically $\sigma^{* p}$-connective, and hence from the homotopy orbits spectral sequence of Section \ref{sec:homotopy-orbit-ss} that 
\[\cC_\infty(p)_+ \wedge_{\fS_p} (Q^{E_\infty}_\bL(\gS)/Q^{E_\infty}_\bL(\gR))^{\otimes p}\]
is homologically $\sigma^{* p}$-connective too. 

In particular, as $\sigma * \sigma \geq \sigma$ it follows that the map
\[0_*\bunit  \lra E_\infty^+((-1)_* Q^{E_\infty}_\bL(\gS)/Q^{E_\infty}_\bL(\gR))\]
is $\sigma$-connective, which together with the spectral sequence of Theorem \ref{thm:e2cellsIV-ss} shows that $H_{g,d}^{\overline{\gR}}(\overline{\gS})=0$ for $d < \inf(\bunit_\text{conn}, \sigma)(g)$.

Going further, it also follows that the map
\[0_*\bunit \vee (-1)_*(Q^{E_\infty}_\bL(\gS)/Q^{E_\infty}_\bL(\gR)) \lra E_\infty^+((-1)_* Q^{E_\infty}_\bL(\gS)/Q^{E_\infty}_\bL(\gR))\]
is $(\sigma*\sigma)$-connective, which can be used to obtain a map $H_{g,d}^{\overline{\gR}}(\overline{\gS}, \overline{\gR}) \to H_{g,d}^{E_\infty}(\gS, \gR)$ which is an isomorphism for $d < (\sigma * \sigma)(g)$ and an epimorphism for $d < (\sigma * \sigma)(g)+1$. We do not state these as formal theorems here, as they follow from the more general Theorem \ref{thm:EkModuleHurewicz} and Corollary \ref{cor:BarEstimate} below.

\begin{lemma}\label{lem:modIndFreeAlg} 
Let $X \to Y$ be a cofibration between cofibrant, 0-connective and reduced objects in $\sfC$, inducing a morphism $\gE_\infty(X) \to \gE_\infty(Y)$ of non-unital $E_\infty$-algebras. Then there is a homology equivalence
	\[B\left(\bunit, \overline{\gE_\infty(X)}, \overline{\gE_\infty(Y)}\right) \simeq 
		 E_\infty^+(Y/X)\]
	where $Y/X$ denotes the cofibre of $X \to Y$.
\end{lemma}

\begin{proof}
	The map $\overline{\gE_\infty(Y)} \to \overline{\gE_\infty(Y/X)}$ induced by the quotient $Y \to Y/X$ gives an augmentation
	\begin{equation}\label{eq:AugMap}
		B_\bullet(\bunit, \overline{\gE_\infty(X)}, \overline{\gE_\infty(Y)}) \lra \overline{\gE_\infty(Y/X)},
	\end{equation}
	and we shall show that this is a homology equivalence after geometric realization.
	
	To do so, promote $X$ and $Y$ to filtered objects $fX$ and $fY$ by: $fX(i)=\ast$ for $i<0$, and $fX(i) = X$ for $i \geq 0$; $fY(i)=\ast$ for $i<0$, $fY(0) = X$, and $fY(i) = Y$ for $i > 0$, with the evident structure maps. Then $\mr{gr}(fX) = 0_*X$ and $\mr{gr}(fY) = 0_*X \vee 1_*(Y/X)$. The augmented simplicial object \eqref{eq:AugMap} is promoted to a filtered one, with associated graded given by
	\[B_\bullet(0_*\bunit, \overline{\gE_\infty(0_*X)}, \overline{\gE_\infty(0_*X \vee 1_*(Y/X))}) \lra \overline{\gE_\infty(1_*(Y/X))}.\]
	This augmentation is a split epimorphism. By Corollary \ref{cor:CoprodEinfty} with $A = 0_*X$ and $B = 1_*(Y/X)$ we have an equivalence 
	\[\overline{\gE_\infty(0_*X \vee 1_*(Y/X))} \simeq 
		\overline{\gE_\infty(0_*X)} \otimes {E^+_k(1_*(Y/X))}\]
	of left $\overline{\gE_\infty(0_*X)}$-modules, so that
\[B(0_*\bunit, \overline{\gE_\infty(0_*X)}, \overline{\gE_\infty(0_*X \vee 1_*(Y/X))}) \simeq 
		{E^+_\infty(1_*(Y/X))}.\]
	Thus this augmentation is a weak equivalence. It follows from the strongly convergent spectral sequence associated to this filtered object that \eqref{eq:AugMap} is a homology equivalence after geometric realization.
\end{proof}

\begin{proposition}\label{prop:CanMultFiltOnBar}
	Suppose that $\sfG$ is Artinian. Let $f \colon \gR \to \gS$ be a morphism in $\Alg_{E_\infty}(\sfC)$ between cofibrant, $0$-connective and reduced objects. Then there is a cofibrant descendingly filtered object with colimit $B(\bunit, \overline{\gR}, \overline{\gS})$, contractible homotopy limit, and whose associated graded is homology equivalent to
	\[E_\infty^+\left((-1)_* Q^{E_\infty}_\bL(\gS)/Q^{E_\infty}_\bL(\gR)\right),\]
	where $Q^{E_\infty}_\bL(\gS)/Q^{E_\infty}_\bL(\gR)$ denotes the homotopy cofibre of $Q^{E_\infty}_\bL(\gR) \to Q^{E_\infty}_\bL(\gS)$.
\end{proposition}

\begin{proof}
	We shall use the canonical multiplicative filtration described in Section \ref{sec:canon-mult-filtr}, given by the functor
	\begin{equation*}
		(-1)_*^\mr{alg}\colon \Alg_{E_\infty}(\sfC) \lra \Alg_{E_\infty}(\sfC^{\bZ_{\leq}})
	\end{equation*}
	left adjoint to evaluating at $-1 \in \bZ_\leq$. 
	
	There is an induced morphism 
	\[(-1)_*^\mr{alg} \gR \lra (-1)_*^\mr{alg} \gS,\]
	between cofibrant objects in $\Alg_{E_k}(\sfC^{\bZ_{\leq}})$,
	which can be rectified to a morphism of monoids $\overline{(-1)_*^\mr{alg} \gR} \to \overline{(-1)_*^\mr{alg} \gS}$ between objects which are cofibrant in $\sfC^{\bZ_{\leq}}$, by Lemma \ref{lem:tilde-properties} (i).  We may therefore form the cofibrant filtered object
	\begin{equation}
		B\left(0_* \bunit, \overline{(-1)_*^\mr{alg} \gR}, \overline{(-1)_*^\mr{alg} \gS}\right) \in \sfC^{\bZ_{\leq}}.\label{eq:19}
	\end{equation}
	This has colimit $B(\bunit, \overline{\gR}, \overline{\gS})$. In the proof of Theorem \ref{thm:CanMultFiltSS} we established the claim that for any coefficients $A$ and reduced $\gR$, $\bL C_*((-1)_*^\mr{alg}(\gR)(-a);A)$ is $c^{\ast a}$-connective, with the abstract connectivity $c$ as in \eqref{eqn:connectivity-c}. The same result will hold for $\gS$. Because $c \ast c \geq c$, this condition is preserved by tensor product. By an induction over the skeleta it is also preserved by geometric realization, so it holds for the bar construction too. As $\sfG$ is assumed Artinian it admits a rank functor $r \colon \sfG \to \bN_{\leq}$, and we have the estimate $c^{* a}(g) \geq a-r(g)$. It follows that~\eqref{eq:19} has contractible homotopy limit.
	
	The associated graded is
	\[B\left(0_* \bunit, \overline{\gE_\infty((-1)_*Q^{E_\infty}_\bL(\gR))}, \overline{\gE_\infty((-1)_*Q^{E_\infty}_\bL(\gS))}\right) \in \sfC^{\bZ_{=}},\]
	as taking associated graded commutes with $\overline{(-)}$ by Lemma \ref{lem:tilde-properties} (iii), and then using Proposition \ref{prop:associated-graded-of-canonical-filt}. Lemma \ref{lem:modIndFreeAlg} identifies this up to homology equivalence with the associated graded as described in the statement.
\end{proof}

\begin{proof}[Proof of Theorem \ref{thm:e2cellsIV-ss}]
This is the spectral sequence for the filtered object of Proposition \ref{prop:CanMultFiltOnBar}; as its homotopy limit is contractible the spectral sequence converges conditionally (see Theorem \ref{thm:SSAsc}).

The strong convergence is proved as in Theorem \ref{thm:CanMultFiltSS}, which yields the special case $\gR = \binit$ when we take $\cO = \cC_\infty$. We briefly recall the argument. Since $\gR$ and $\gS$ are reduced and $0$-connective, so is $Q^{E_\infty}_\bL(\gS)/Q^{E_\infty}_\bL(\gR)$, i.e.~it is $c$-connective. Then $E_\infty^+((-1)_* Q^{E_\infty}_\bL(\gS)/Q^{E_\infty}_\bL(\gR))(-a)$ is $c^{\ast a}$-connective, so that in terms of a rank functor $r \colon \sfG \to \bN_{\leq}$ we have $E^1_{g,p,q}=0$ for $p+q < -p-r(g)$. In this spectral sequence the differentials take the form
	\[d^r \colon E^r_{n,p,q} \lra E^r_{n,p-r,q+r-1}\]
so there are only finitely-many differentials exiting each position. Therefore the derived $E^\infty$-page vanishes and by \cite[Theorem 7.3]{Boardman} this spectral sequence actually converges strongly.
\end{proof}

\subsection{Comparison in general}

The comparison results explained after Theorem \ref{thm:e2cellsIV-ss} will, formulated correctly, hold for $E_k$-algebras, although nothing so simple as Lemma \ref{lem:modIndFreeAlg}  or Proposition \ref{prop:CanMultFiltOnBar} is available when $k < \infty$.

\begin{theorem}\label{thm:EkModuleHurewicz}
Suppose that $\sfS$ satisfies Axiom \ref{axiom:Hurewicz}, and that $\sfG$ is $(k+1)$-monoidal and Artinian. Let $\rho, \sigma \colon \sfG \to [-\infty,\infty]_{\geq}$ be abstract connectivities such that $\rho * \rho \geq \rho$, $\sigma * \sigma \geq \sigma$, and $\rho * \sigma \geq \sigma * \sigma$. If
\begin{enumerate}[\indent (i)]
\item $\gR \in \Alg_{E_k}(\sfC)$ is such that $H^{E_k}_{g,d}(\gR)=0$ for $d < \rho(g)-(k-1)$,
\item $f \colon \gR \to \gS$ is an $E_k$-algebra map such that $H^{E_k}_{g,d}(\gS, \gR)=0$ for $d <  \sigma(g)$, and
\item $\gR$ and $\gS$ are cofibrant in $\sfC$, $0$-connective, and reduced,
\end{enumerate}
then there is a map $H^{\overline{\gR}}_{g,d}(\overline{\gS}, \overline{\gR}) \to H^{E_k}_{g,d}(\gS, \gR)$ which is an isomorphism for $d < (\sigma * \sigma)(g)$, and an epimorphism for $d < (\sigma * \sigma)(g)+1$.
\end{theorem}

\begin{proof}
Without loss of generality, we may assume that $f \colon \gR \to \gS$ is a cofibration of cofibrant $E_k$-algebras. Similarly to the proof of Proposition \ref{prop:CanMultFiltOnBar}, we use the canonical multiplicative filtration, and so consider the induced cofibration $\overline{(-1)_*^\mr{alg} \gR} \to \overline{(-1)_*^\mr{alg} \gS}$ of associative monoid objects which are cofibrant in $\sfC^{\bZ_{\leq}}$. We then form the cofibrant filtered object
\[B \coloneqq B(0_* \bunit, \overline{(-1)_*^\mr{alg} \gR}, \overline{(-1)_*^\mr{alg} \gS}/\overline{(-1)_*^\mr{alg} \gR}) \in \mathsf{C}^{\bZ_{\leq}},\]
which has colimit $B(\bunit, \overline{\gR}, \overline{\gS}/\overline{\gR})$ whose homology groups are $H^{\overline{\gR}}_{*,*}(\overline{\gS}, \overline{\gR})$. This filtered object has contractible homotopy limit just as in the proof of Proposition \ref{prop:CanMultFiltOnBar}, using the assumption that $\sfG$ is Artinian. Its associated graded is
	\[\mr{gr}(B) = B\left(0_* \bunit, \overline{\gE_k((-1)_*Q^{E_k}_\bL(\gR))}, \overline{\gE_k((-1)_*Q^{E_k}_\bL(\gS))} /  \overline{\gE_k((-1)_*Q^{E_k}_\bL(\gR))} \right).\]
The graded object $\overline{\gE_k((-1)_*Q^{E_k}_\bL(\gS))} /  \overline{\gE_k((-1)_*Q^{E_k}_\bL(\gR))}$ is contractible in grading 0 and is $Q^{E_k}_\bL(\gS)/Q^{E_k}_\bL(\gR)$ in grading $-1$, so $\mr{gr}(B)$ has these properties too. Thus the spectral sequence
\[E^1_{g,p,q} = H_{g, p+q, p}(\mr{gr}(B)) \Longrightarrow H^{\overline{\gR}}_{g,p+q}(\overline{\gS}, \overline{\gR})\]
with differentials $d^r \colon E^r_{n,p,q} \to E^r_{n,p-r,q+r-1}$ provides an edge homomorphism $H^{\overline{\gR}}_{g,d}(\overline{\gS}, \overline{\gR}) \to H_{g, d, -1}(\mr{gr}(B)) = H_{g,d}^{E_k}(\gS, \gR)$, which is the required map. The connectivity property of this map then follows from this spectral sequence and the following vanishing result.
	
\vspace{1ex}
	
\noindent\textbf{Claim.} $H_{g,p+q,p}(\mr{gr}(B))=0$ for $p+q < (\sigma * \sigma)(g)$ and $p < -1$.

\vspace{1ex}
	 
As in the proof of Lemma \ref{lem:modIndFreeAlg} we may endow $Q^{E_k}_\bL(\gR)$ and $Q^{E_k}_\bL(\gS)$ with filtrations such that $\mr{gr}(Q^{E_k}_\bL(\gR)) = 0_*(Q^{E_k}_\bL(\gR))$ and $\mr{gr}(Q^{E_k}_\bL(\gS)) = 0_*(Q^{E_k}_\bL(\gR)) \vee 1_*(Q^{E_k}_\bL(\gS)/Q^{E_k}_\bL(\gR))$. This induces a further filtration of $\mr{gr}(B)$. By Theorem \ref{thm:AlgAsMod} $\overline{\gE_k((-1)_*\mr{gr}(Q^{E_k}_\bL(\gS)))} /  \overline{\gE_k((-1)_*\mr{gr}(Q^{E_k}_\bL(\gR)))}$ is equivalent to
\begin{align*} &\overline{\gE_k((-1,0)_*\mr{gr}(Q^{E_k}_\bL(\gR)))} \\
&\qquad\otimes E_k(E_1^+(S^{k-1} \wedge (-1,0)_*(Q^{E_k}_\bL(\gR))) \otimes (-1,1)_*(Q^{E_k}_\bL(\gS)/Q^{E_k}_\bL(\gR)) )\end{align*}
as a $\overline{\gE_k((-1,0)_*\mr{gr}(Q^{E_k}_\bL(\gR)))}$-module, and so we have
\[\mr{gr}(\mr{gr}(B)) \simeq E_k(E_1^+(S^{k-1} \wedge (-1,0)_*(Q^{E_k}_\bL(\gR))) \otimes (-1,1)_*(Q^{E_k}_\bL(\gS)/Q^{E_k}_\bL(\gR)))\]
as an object of $\sfC$ with two additional $\bZ$-gradings. Let us for now suppress the additional gradings from the notation.

By assumption $Q^{E_k}_\bL(\gS)/Q^{E_k}_\bL(\gR)$ is homologically $\sigma$-connective and $S^{k-1} \wedge Q^{E_k}_\bL(\gR)$ is homologically $\rho$-connective. In particular, using
$$E_1^+(S^{k-1} \wedge Q^{E_k}_\bL(\gR)) \simeq \bigvee_{p \geq 0} (S^{k-1} \wedge Q^{E_k}_\bL(\gR))^{\otimes p}$$
it follows from Lemma \ref{lem:connectivity-and-tensor-products} (i) (and the fact that $\rho*\rho \geq \rho$, which we shall use repeatedly) that the map
$$\bunit \lra E_1^+(S^{k-1} \wedge Q^{E_k}_\bL(\gR))$$
is homologically $\rho$-connective, its source is homologically $\bunit_\text{conn}$-connective, and its target is homologically $\inf(\bunit_\text{conn}, \rho)$-connective. From this it follows that the map
$$Q^{E_k}_\bL(\gS)/Q^{E_k}_\bL(\gR) \lra E_1^+(S^{k-1} \wedge Q^{E_k}_\bL(\gR)) \otimes (Q^{E_k}_\bL(\gS)/Q^{E_k}_\bL(\gR))$$
is homologically $(\rho * \sigma)$-connective, its source is homologically $\sigma$-connective, and its target is homologically $(\inf(\bunit_\text{conn}, \rho)*\sigma)$-connective. Using the definition
$$E_k^+(-) = \bigvee_{p \geq 0} \mathcal{C}_k(p)_+ \wedge_{\fS_p} (-)^{\otimes p},$$
Corollary \ref{cor:connectivity-of-tensor-power}, and the homotopy orbits spectral sequence of Section \ref{sec:homotopy-orbit-ss}, it follows that the map
$$E_k(Q^{E_k}_\bL(\gS)/Q^{E_k}_\bL(\gR)) \lra \mr{gr}(\mr{gr}(B))$$
is homologically $\inf_{a+b \geq 0}(\sigma^{* a} * (\rho * \sigma) * (\inf(\bunit_\text{conn}, \rho) * \sigma)^{* b})$-connective. Manipulating the various conditions on the abstract connectivities $\rho$ and $\sigma$, it follows that this map is homologically $\sigma * \sigma$-connective. On the other hand, the map
$$Q^{E_k}_\bL(\gS)/Q^{E_k}_\bL(\gR) \lra E_k(Q^{E_k}_\bL(\gS)/Q^{E_k}_\bL(\gR))$$
is also homologically $\sigma * \sigma$-connective, as $(Q^{E_k}_\bL(\gS)/Q^{E_k}_\bL(\gR))^{\otimes p}$ is for all $p \geq 2$. Recalling now the additional gradings, the composition 
$$(-1,1)_*(Q^{E_k}_\bL(\gS)/Q^{E_k}_\bL(\gR)) \lra \mr{gr}(\mr{gr}(B))$$
is the induced map on associated gradeds for the second filtration induced by the inclusion $(-1)_*( Q^{E_k}_\bL(\gS)/Q^{E_k}_\bL(\gR)) \to \mr{gr}(B)$ of the piece of grading $-1$, so the claim follows.
\end{proof}

We will typically apply this result via the following special case.

\begin{corollary}\label{cor:BarEstimate}
Suppose that $\sfS$ satisfies Axiom \ref{axiom:Hurewicz}, that $\sfG$ is $(k+1)$-monoidal and Artinian, and that $\rho \colon \sfG \to [-\infty,\infty]_{\geq}$ is an abstract connectivity so that $\rho * \rho \geq \rho$. If
\begin{enumerate}[\indent (i)]
\item $\gR \in \Alg_{E_k}(\sfC)$ is such that $H^{E_k}_{g,d}(\gR)=0$ for $d < \rho(g)-(k-1)$,
\item $f \colon \gR \to \gS$ is an $E_k$-algebra map such that $H^{E_k}_{g,d}(\gS, \gR)=0$ for $d <  \rho(g)$, and
\item $\gR$ and $\gS$ are cofibrant in $\sfC$, $0$-connective, and reduced,
\end{enumerate}
then we have $H^{\overline{\gR}}_{g,d}(\overline{\gS})=0$ for $d < \inf(\bunit_\mr{conn},\rho)(g)$. 

In addition, for an abstract connectivity $\mu$ such that $\mu*\rho \geq \mu$, if
\begin{enumerate}[\indent (i)]
\setcounter{enumi}{3}
\item $\gM$ is a left $\overline{\gR}$-module such that $H_{g,d}(\gM)=0$ for $d < \mu(g)$, and
\item $\gM$ is cofibrant in $\sfC$,
\end{enumerate}
then we have $H_{g,d}(B(\gM, \overline{\gR}, \overline{\gS}))=0$ for $d < \mu(g)$.
\end{corollary}

\begin{proof}
For the first part, we apply Theorem \ref{thm:EkModuleHurewicz} with $\sigma=\rho$. This yields a map $H^{\overline{\gR}}_{g,d}(\overline{\gS}, \overline{\gR}) \to H^{E_k}_{g,d}(\gS, \gR)$ which is an isomorphism for $d < (\rho*\rho)(g)$, and as the latter vanishes for $d < \rho(g) \leq (\rho*\rho)(g)$, it follows that $H^{\overline{\gR}}_{g,d}(\overline{\gS}, \overline{\gR}) =0$ for $d < \rho(g)$. 
From the long exact sequence
$$\cdots \lra H^{\overline{\gR}}_{g,d}(\overline{\gR}) \overset{\overline{f}_*}\lra H^{\overline{\gR}}_{g,d}(\overline{\gS}) \lra H^{\overline{\gR}}_{g,d}(\overline{\gS}, \overline{\gR}) \overset{\partial}\lra H^{\overline{\gR}}_{g,d-1}(\overline{\gR}) \lra \cdots$$
and the fact that $H^{\overline{\gR}}_{g,d}(\overline{\gR})$ is supported in bidegree $(g,d)=(0,0)$, it follows that $H^{\overline{\gR}}_{g,d}(\overline{\gS})=0$ for $d < \inf(\bunit_\mr{conn},\rho)(g)$. 

For the second part, we apply Theorem \ref{thm:MinCellStr-additive} to the operad $\cO$ having 1-ary operations given by $\overline{\gR}$ and no higher arity operations, whose algebras are right $\overline{\gR}$-modules, and the $\cO$-algebra morphism $\overline{f} : \overline{\gR} \to \overline{\gS}$. The augmentation $\epsilon: \overline{\gR} \to \bunit$ makes $\cO$ into an augmented non-unitary homologically connective $\Sigma$-cofibrant operad in $\sfC$. Furthermore, as $\gR$ is reduced, $\epsilon: \cO(1) = \overline{\gR} \to \bunit$ satisfies assumption (i) of Lemma \ref{lem:TensorDetectNull}. Thus Theorem \ref{thm:MinCellStr-additive} indeed applies, and gives a relative CW approximation $\overline{f}: \overline{\gR} \to \colim \mr{sk}(\overline{f}) \to \overline{\gS}$ where $\mr{sk}(\overline{f})$ has no $(g,d)$-cells with $d < \mu(g)$. This yields a filtered object
$$B(0_* \gM, 0_* \overline{\gR}, \mr{sk}(\overline{f}))$$
with colimit $B(\gM, \overline{\gR}, \mr{sk}(\overline{f})) \simeq B(\gM, \overline{\gR}, \overline{\gS})$, contractible limit, and associated graded (suppressing the additional grading) given by
$$B(\gM, \overline{\gR}, \overline{\gR} \vee \bigvee_{\alpha \in I} \overline{\gR} \otimes S^{g_\alpha, d_\alpha}) \simeq \gM \otimes (S^{0,0} \vee \bigvee_{\alpha \in I} S^{g_\alpha, d_\alpha})$$
with $d_\alpha \geq \mu(g_\alpha)$. By assumption $\gM$ is homologically $\mu$-connective, and $(S^{0,0} \vee \bigvee_{\alpha \in I} S^{g_\alpha, d_\alpha})$ is homologically $\inf(\bunit_\mr{conn},\rho)$-connective, so this has homological connectivity $\mu * \inf(\bunit_\mr{conn},\rho) = \inf(\mu, \mu * \rho) \geq \mu$. The spectral sequence for this filtered object then shows that $B(\gM, \overline{\gR}, \overline{\gS})$ has this homological connectivity too.
\end{proof}

\subsection{Proof of Theorem \ref{thm:AlgAsMod}}\label{sec:proofAlgAsMod}

In Section \ref{sec:construct-fk} below we will describe a map
\[f_k^+ \colon E_1^+(S^{k-1} \wedge A) \otimes B \overset{f_k}\lra E_k(A \vee B) \xrightarrow{\mr{inc}} E_k^+(A \vee B),\]
which only exists when the homotopy category of $\sfC$ is enriched in abelian groups (which follows from Axiom \ref{axiom:Hurewicz}). As the target is obtained by neglect of structure from an $E^+_k$-algebra this extends to a map $F_k \colon \gE_k^+(E_1^+(S^{k-1} \wedge A) \otimes B) \to \gE_k^+(A \vee B)$ of unital $E_k$-algebras, and hence to a map
\[\overline{F}_k \colon E_k^+(E_1^+(S^{k-1} \wedge A) \otimes B) \lra E_k^+(A \vee B) \lra \overline{E_k(A \vee B)}.\]
Furthermore, the map $\overline{\gE_k(A)} \to \overline{\gE_k(A \vee B)}$ makes $\overline{E_k(A \vee B)}$ into a left $\overline{\gE_k(A)}$-module, so $\overline{F}_k$ extends to a map
\[\alpha_k \colon \overline{E_k(A)} \otimes E_k^+(E_1^+(S^{k-1} \wedge A) \otimes B) \lra  \overline{E_k(A \vee B)}\]
of left $\overline{\gE_k(A)}$-modules. 

It is this map which we shall show is a weak equivalence. Though we have constructed it as a map of left $\overline{\gE_k(A)}$-modules, to show it is a weak equivalence we may forget this module structure.

If $\gR$ is a $E_1$-algebra then the inclusion $\gR^+ = \bunit \sqcup \gR \to \overline{\gR}$ is a weak equivalence by the formula (\ref{eq:TildeUnderlying}). In particular this applies to $\gR = \gE_k(X)$. Thus there is a homotopy commutative diagram
\begin{equation*}
\begin{tikzcd}
	E_k^+(A) \otimes E_k^+(E_1^+(S^{k-1} \wedge A) \otimes B) \rar{\simeq} \dar[swap]{E^+_k(\iota) \otimes F_k} \arrow[bend left=70,looseness=3]{dd}{\beta_k} & \overline{E_k(A)} \otimes E_k^+(E_1^+(S^{k-1} \wedge A) \otimes B) \arrow{dd}{\alpha_k} \\
	E_k^+(A \vee B) \otimes E_k^+(A \vee B) \dar[swap]{\mu} & \\
	E_k^+(A \vee B) \rar{\simeq} &  \overline{E_k(A \vee B)}
\end{tikzcd}
\end{equation*}
where $\mu$ is given by multiplication using a point in $\cC_1^+(2) \subset \cC_k^+(2)$. Denoting by $\beta_k$ the left composition in the diagram, 
we must show that this is an equivalence.

If $k=1$ then we will show this directly. If $k \geq 2$ then we have assumed that $\sfG$, and hence $\sfC$, is $\infty$-monoidal, so by Proposition \ref{prop:algcc-symmetric-monoidal} there is a $(k-1)$-monoidal structure on the category of $E_{k-1}^+$-algebras in $\sfC$. Using this to consider $E_k^+(A \vee B) \otimes E_k^+(A \vee B)$ as an $E_{k-1}^+$-algebra, our choice of multiplication as lying in $\cC_1^+(2)$ means that the map $\mu$ is an $E_{k-1}^+$-algebra map. Thus the map $\beta_k$ is a map of augmented $E^+_{k-1}$-algebras, so to show it is a weak equivalence it is enough, by Axiom \ref{axiom:Hurewicz} and Proposition \ref{prop:relwhitehead}, to take augmentation ideals $I(-)$ as in Section \ref{sec:non-unitary}, and show that $\smash{S^{k-1}\wedge Q^{E_{k-1}}_\bL(I(\beta_k))}$ is a weak equivalence. Applying Proposition \ref{prop:relwhitehead} uses our assumption that $A$ and $B$ are 0-connective and reduced.

We next apply the natural weak equivalence $S^{k-1} \wedge Q^{E_{k-1}}_\bL(-) \simeq \tilde{B}^{E_{k-1}}(-)$ of Theorem \ref{thm:BarHomologyIndec}. We shall apply these functors to the augmentation ideal of either a free $E^+_k$-algebra, or a tensor product of free $E_{k}^+$-algebras. In the first case we can use 
our calculation of the bar construction applied to a free algebra (from Section \ref{sec:bar-on-free-algebras}), and in the second case we will use the following lemma.

\begin{lemma}
There is a natural weak equivalence 
\[B^{E_k}(\gE^+_{n+k}(X) \otimes_{E_{k}^+} \gE^+_{n+k}(Y),\epsilon_\mr{can}) \simeq E_n^+(S^k \wedge X) \otimes E_n^+(S^k \wedge Y)\]
if $X$ and $Y$ are cofibrant.
\end{lemma}

\begin{proof}
By Proposition \ref{prop:algcc-symmetric-monoidal} the map $U^{E_k^+}(\gR \otimes_{E^+_k} \gS) \to U^{E_k^+}(\gR) \otimes U^{E_k^+}(\gS)$ is an isomorphism. Thus there is an isomorphism of $k$-fold semi-simplicial objects
	\begin{align*}B^{E_k}_{\bullet,\ldots,\bullet}&(\gE^+_{n+k}(X) \otimes_{E_{k}^+} \gE^+_{n+k}(Y),\epsilon_\mr{can}) \\
	&\cong B^{E_k}_{\bullet,\ldots,\bullet}(\gE^+_{n+k}(X),\epsilon_\mr{can}) \otimes B^{E_k}_{\bullet,\ldots,\bullet}(\gE^+_{n+k}(Y),\epsilon_\mr{can}).\end{align*}
	
	Both are restrictions of $k$-fold simplicial objects, which are Reedy cofibrant by Lemma \ref{lem:reedy-cofibrant-split-degeneracies} because their degeneracies are split. This implies that for the right hand side, the tensor product commutes with thick geometric realization up to weak equivalence:
	\begin{align*}B^{E_k}(\gE^+_{n+k}(X) \otimes_{E_{k}^+} &\gE^+_{n+k}(Y),\epsilon_\mr{can}) \\
&\cong \fgr{B^{E_k}_{\bullet,\ldots,\bullet}(\gE^+_{n+k}(X),\epsilon_\mr{can}) \otimes B^{E_k}_{\bullet,\ldots,\bullet}(\gE^+_{n+k}(Y),\epsilon_\mr{can})} \\
	&\simeq B(\gE^+_{n+k}(X),\epsilon_\mr{can}) \otimes B(\gE^+_{n+k}(Y),\epsilon_\mr{can}).\end{align*}
	
	Finally, we use Theorem \ref{thm:CalcFree} and the definition of $\tilde{B}^{E_k}$ in terms of $B^{E_k}$, to conclude that there is a natural weak equivalence $B^{E_k}(\gE^+_{n+k}(-),\epsilon_\mr{can}) \simeq E_n^+(S^k \wedge -)$ on cofibrant objects.
\end{proof}

In particular, we obtain weak equivalences
\begin{align*}S^0 \vee & \left( S^{k-1} \wedge Q^{E_{k-1}}_\bL(\gE_k^+(A) \otimes \gE_k^+(E_1^+(S^{k-1} \wedge A) \otimes B))) \right)\\
&\simeq E^+_1(S^{k-1} \wedge  A) \otimes E^+_1(S^{k-1} \wedge E_1^+(S^{k-1}\wedge A)\otimes B),\end{align*}
and
\[S^0 \vee \left( S^{k-1} \wedge Q^{E_{k-1}}_\bL(E_k^+(A \vee B)) \right) \simeq E^+_1(S^{k-1}\wedge (A \vee B)).\]

Thus for $k \geq 2$ checking that $\beta_k$ is a weak equivalence amounts to proving that the homotopy class
\[\begin{tikzcd}
E_1^+(S^{k-1} \wedge A) \otimes E_1^+(S^{k-1} \wedge E_1^+(S^{k-1}\wedge A) \otimes B) \dar{\mu(E_1^+(\iota) \otimes G_k)} \\ E_1^+(S^{k-1} \wedge (A \vee B))
\end{tikzcd}\]
is a weak equivalence, where we have used the computation of the bar construction on maps between free algebras (from Section \ref{sec:FreeAlgMaps}) to write $S^0\vee S^{k-1} \wedge Q^{E_{k-1}}_\bL(\beta_k)$ as $\mu(E_1^+(\iota) \wedge G_k)$, where the map $G_k$ is obtained by freely extending
\[\begin{tikzcd}
S^{k-1} \wedge E_1^+(S^{k-1} \wedge A) \otimes B \dar{S^{k-1} \wedge f_k} \arrow[bend right=70,looseness=2]{dd}[swap]{g_k} \\ S^{k-1} \wedge E_k(A \vee B) \dar{\eta^+_{A \vee B}} \\ E^+_1(S^{k-1} \wedge(A \vee B)),
\end{tikzcd}\]
to a map of $E^+_1$-algebras. Here we have used the natural transformation $\eta$ from Section \ref{sec:FreeAlgMaps} composed with the inclusion $E_1(-) \hookrightarrow E_1^+(-)$, which we denote $\eta^+$. We will show that this map $g_k$ may be identified up to homotopy with $f_1^+$ (with $A$ replaced by $S^{k-1}\wedge A$ and $B$ replaced by $S^{k-1}\wedge B$), so the map $\mu(E_1^+(\iota) \wedge G_k)$ is homotopic to the map $\mu(E^+_k(\iota) \wedge F_1) = \beta_1$, which reduces us to proving the case $k=1$. In this case we may work instead with the associative operad, where it will be a direct calculation.

\subsubsection{Constructing $f_k$}\label{sec:construct-fk}
Let us define a map
\[f'_k \colon E_1^+(S^{k-1} \wedge (A \vee B)) \otimes E_k(A \vee B) \lra E_k(A \vee B),\]
from which we obtain $f_k$ by precomposing with the product of inclusions
\[E_1^+(S^{k-1} \wedge A) \otimes B \lra E_1^+(S^{k-1} \wedge (A \vee B)) \otimes E_k(A \vee B)\]
and obtain $f_k^+$ by composing with the inclusion $E_k(A \vee B) \hookrightarrow E^+_k(A \vee B)$.
The domain and codomain of $f'_k$ are functors of $X = A \vee B$ and $f'_k$ will be a natural transformation of functors of $X$ whose definition need only refer to $X$, not $A$ and $B$ individually.

We define $f'_k$ by defining its adjoint
\[\hat{f}'_k \colon E_1^+(S^{k-1} \wedge X) \lra \Hom_\sfC(E_k(X), E_k(X))\]
using the internal hom object $\Hom_\sfC$ of the category $\sfC$. The target of $\hat{f}'_k$ is an associative unital monoid in $\sfC$, so in particular an $E_1^+$-algebra. Thus we can define $\hat{f}'_k$ as an $E_1^+$-map by describing a map
\[\mr{ad}_k \colon S^{k-1} \wedge X \lra \Hom_\sfC(E_k(X), E_k(X))\]
in $\sfC$, and freely extending it to a map of $E_1^+$-algebras.

To define $\mr{ad}_k$, we use that there is the morphism $s'_k \colon S^{k-1} \to S^{k-1}_+$ in $\sfC$ constructed as follows: there is a map $i_0 \colon S^0 \to S^{k-1}_+$ with retraction $r_0 \colon S^{k-1}_+ \to S^0$. As $\sfS$ is assumed to satisfy Axiom \ref{axiom:Hurewicz} it is in particular semistable, and it immediately follows that $\sfC = \sfS^\sfG$ is also semistable. Thus $\cat{Ho}(\sfC)$ is enriched in abelian groups (as described in Section \ref{sec:additive-case}). Using this enrichment we may form the difference $\mr{id}-i_0 \circ r_0 \colon S^{k-1}_+ \to S^{k-1}_+$. If we precompose this with $i_0$ we get $0$ in the homotopy category, so up to homotopy it factors over the homotopy cofiber $S^{k-1}$ as a map $s'_k \colon S^{k-1} \to S^{k-1}_+$. This induces a map 
\[s_k \coloneqq s'_k \wedge \mr{id}_X \colon  S^{k-1} \wedge X \longrightarrow S^{k-1}_+ \wedge X.\] 

If $k \geq 2$, so that we have assumed that $\sfG$ is symmetric monoidal, the map $\mr{ad}_k$ is defined as $\mu_k \circ s_k$ in terms of a map
\[\mu_k \colon S^{k-1}_+ \wedge X \lra \Hom_\sfC(E_k(X), E_k(X))\]
given by picking a homotopy equivalence $m_k \colon S^{k-1} \to \cC_k(2)$ (we give an explicit formula later) and then taking the adjoint of
\[\hat{\mu}_k \colon S^{k-1}_+ \wedge X \otimes E_k(X) \xrightarrow{m_k \otimes \mr{id}} \cC_{k}(2)_+ \wedge X \otimes E_k(X) \lra E_k(X)\]
with second map given by the canonical morphism $X \to E_k(X)$ and the $E_k$-algebra structure on $E_k(X)$.

If $k=1$, so that we have assumed that $\sfG$ is braided monoidal, then the map $\mr{ad}_k$ is defined as above but the map $\hat{\mu}_1$ is given on $S^{k-1} = S^0 = \{\pm 1\}$ by
\[\hat{\mu}_1\vert_{+1} \colon \{+1\}_+ \wedge X \otimes E_1(X) \xrightarrow{m_1 \otimes \mr{id}} \cC_{1}^{\cat{FB}_1}(2)_+ \wedge X \otimes E_1(X) \lra E_1(X)\]
and
\[\hat{\mu}_1\vert_{-1} \colon \{-1\}_+ \wedge X \otimes E_1(X) \xrightarrow{m_1 \otimes b} \cC_{1}^{\cat{FB}_1}(2)_+ \wedge E_1(X) \otimes X \lra E_1(X)\]
where $m_1 \colon \{*\} \to \cC_{1}^{\cat{FB}_1}(2)$ is a homotopy equivalence, $b$ denotes the braiding $\beta_{X, E_1(X)}$, and as above in both cases the second map is given by the canonical morphism $X \to E_1(X)$ and the $E_1$-algebra structure on $E_1(X)$. Note that if $\sfG$ is actually symmetric monoidal then the description given for the case $k \geq 2$ above still makes sense, and it agrees with this definition of $\hat{\mu}_k$.

\subsubsection{Comparing $f_1^+$ and $g_k$}

We shall now explain how the map $g_k$ may be identified with a special case of the map $f_1^+$. This only has content for $k \geq 2$, in which case $\sfG$ is assumed to be symmetric monoidal and we may use the simpler construction of $f_k$ in the previous section. More precisely we will show that the map $g_k$ for $A$ and $B$ is homotopic to the map $f_1^+$ for $S^{k-1} \wedge A$ and $S^{k-1} \wedge B$, after identifying the domains of these maps using the isomorphism
\[E_1^+(S^{k-1} \wedge A) \otimes S^{k-1} \wedge B \cong S^{k-1} \wedge E_1^+(S^{k-1} \wedge A) \otimes B\]
given by the symmetric monoidal structure preceded by the map induced by negation $w \mapsto \vec{1}-w$ on $S^{k-1} \wedge A = ((0,1)^{k-1})^+ \wedge A$.

To see this we let $X = A \vee B$ and consider the diagram
\begin{equation}\label{eq:ComparingF1Gk}
\begin{tikzcd}
E_1^+(S^{k-1} \wedge A) \otimes S^{k-1}\wedge B \rar{\cong} \ar{d}{\mr{inc} \otimes S^{k-1} \wedge \mr{inc}} & S^{k-1} \wedge E_1^+(S^{k-1} \wedge A) \otimes B \ar{d}{S^{k-1}\wedge\mr{inc} \otimes \mr{inc}}\\
E_1^+(S^{k-1} \wedge X) \otimes S^{k-1} \wedge E_k(X) \rar{\cong} \ar{dd}{\mr{id} \otimes \eta_{X}} \arrow[phantom, "\CircNum{!}"]{ddr} & S^{k-1}\wedge E_1^+(S^{k-1}\wedge X) \otimes E_k(X) \ar{d}{S^{k-1} \wedge f'_k}\\
&  S^{k-1} \wedge E_k(X) \ar{d}{\eta_{X}}\\
E_1^+(S^{k-1} \wedge X) \otimes E_1(S^{k-1}\wedge X) \ar{r}{f_1'} & E_1(S^{k-1} \wedge X) \ar{d}{\mr{inc}}\\
& E_1^+(S^{k-1}\wedge X),
\end{tikzcd}
\end{equation}
the horizontal isomorphisms being given by the symmetric monoidality preceded by the isomorphism induced by negation on $S^{k-1}\wedge A$ and $S^{k-1} \wedge X$. The composition along the right-hand edge is the definition of the map $g_k$, whereas the composition along the left-hand edge and bottom is, using that $\eta_{X} \circ S^{k-1}\wedge\mr{inc}$ agrees with the inclusion $\mr{inc} \colon S^{k-1} \wedge B \to E_1(S^{k-1}\wedge X)$, the definition of the map $f_1^+$. The top square of this diagram commutes, so it remains to show that $\CircNum{!}$ commutes up to homotopy.

\begin{lemma}
The square \emph{$\CircNum{!}$} commutes up to homotopy if it does so when precomposed with 
\[\mr{inc} \otimes \mr{id} \colon S^{k-1}\wedge X \otimes S^{k-1} \wedge E_k(X) \lra E_1^+(S^{k-1} \wedge X) \otimes S^{k-1} \wedge E_k(X).\]
\end{lemma}
\begin{proof}
There is a decomposition $E_1^+(S^{k-1} \wedge X) \simeq \bigvee_{n \geq 0} (S^{k-1} \wedge X)^{\otimes n}$, and it suffices to show that the squares $\CircNum{!}_n$ obtained by restricting to the $n$th summand commute up to (based) homotopy. On the $n=0$ summand the maps $S^{k-1} \wedge f_k'$ and $f_1'$ become homotopic to the identity, and so $\CircNum{!}_0$ indeed commutes as both directions give $\eta_{X}$. We have assumed that $\CircNum{!}_1$ commutes up to homotopy. 

For $n>1$ let us write $f'_k\vert_n$ for the restriction of $f'_k$ to the summand
\[(S^{k-1}\wedge X)^{\otimes n} \otimes E_k(X) \subset E_1^+(S^{k-1}\wedge X) \otimes E_k(X).\]
We then have homotopies
\[f_k'\vert_n \simeq f'_k\vert_{n-1} \circ (\mr{id}_{(S^{k-1}\wedge X)^{\otimes n-1}} \otimes f'_k\vert_1),\]
as the adjoint $\hat{f}'_k \colon E_1^+(S^{k-1}\wedge X) \to \Hom_\sfC(E_k(X), E_k(X))$ of $f_k'$ is by definition a morphism of $E_1^+$-algebras. The analogue holds for $f_1'$.

We may then form the diagram in Figure \ref{fig:Reduction} 
\begin{sidewaysfigure*}
	\centering
\[
\begin{tikzcd}[ampersand replacement=\&]
(S^{k-1}\wedge X)^{\otimes n}  \otimes S^{k-1}\wedge E_k(X) \ar{r}{\cong_1} \ar{d}{\mr{id}_{(S^{k-1}\wedge X)^{\otimes n}} \otimes \eta_{X}} \&[-8pt] (S^{k-1}\wedge X)^{\otimes n-1}  \otimes S^{k-1}\wedge (S^{k-1}\wedge X) \otimes E_k(X) \ar{d}[description]{\mr{id}_{(S^{k-1}\wedge X)^{\otimes n-1}}  \otimes S^{k-1}\wedge f_k'\vert_1} \ar{r}{\cong_2} \&[-4pt] S^{k-1}\wedge (S^{k-1}\wedge X)^{\otimes n}  \otimes E_k(X) \ar{d}[swap]{S^{k-1} \wedge \mr{id}_{(S^{k-1}\wedge X)^{\otimes n-1}} \otimes f'_k\vert_1}\\[30pt]
(S^{k-1}\wedge X)^{\otimes n}  \otimes E_1(S^{k-1}\wedge X) \ar{dd}{\mr{id}_{(S^{k-1}\wedge X)^{\otimes n-1}} \otimes f_1'\vert_1} \& (S^{k-1}\wedge X)^{\otimes n-1} \otimes S^{k-1}\wedge E_k(X) \ar{r}{\cong_3} \ar{dd}[description]{\mr{id}_{(S^{k-1}\wedge X)^{\otimes n-1}} \otimes \eta_{X}} \& S^{k-1}\wedge (S^{k-1}\wedge X)^{\otimes n-1} \otimes E_k(X) \ar{d}[swap]{S^{k-1}\wedge f'_k\vert_{n-1}}  \\[7pt]
 \& \& S^{k-1}\wedge E_k(X)  \ar{d}[swap]{\eta_{X}} \\[7pt]
(S^{k-1}\wedge X)^{\otimes n-1} \otimes E_1(S^{k-1}\wedge X) \ar[equals]{r} \& (S^{k-1}\wedge X)^{\otimes n-1} \otimes E_1(S^{k-1}\wedge X) \ar{r}{f_1'\vert_{n-1}} \& E_1(S^{k-1}\wedge X)
\end{tikzcd}
\]
	\caption{}
	\label{fig:Reduction}
\end{sidewaysfigure*}
where the horizontal isomorphisms are given by the symmetric monoidality preceded by: for $\cong_1$ the isomorphism induced by negation on the rightmost copy of $S^{k-1} \wedge X$; for $\cong_2$ and $\cong_3$ the isomorphism induced by negation on the $(n-1)$ leftmost copies of $S^{k-1} \wedge X$. The bottom right-hand square is $\CircNum{!}_{n-1}$, so by induction we may suppose this commutes up to homotopy. The top right-hand square commutes by symmetric monoidality and the fact that the isomorphisms $\cong_2$ and $\cong_3$ use the negation map in the same way. The left-hand square is obtained from $\CircNum{!}_1$ by applying $\mr{id}_{(S^{k-1}\wedge X)^{\otimes n-1}} \otimes -$, so commutes up to homotopy by assumption. The outer rectangle, by the observation above, is $\CircNum{!}_n$ and so commutes up to homotopy as required.
\end{proof}

To verify that the square $\CircNum{!}$ commutes up to homotopy when restricted to $S^{k-1}\wedge X \otimes S^{k-1}\wedge E_k(X)$, we observe that up to homotopy all the maps involved ($f'_k$, $f_1'$, and $\eta_X$) are induced by maps of symmetric sequences of pointed spaces. The only subtlety here involves the maps $s'_k \colon S^{k-1} \to S^{k-1}_+$ and $s'_1 \colon S^{0} \to S^{0}_+$, which are not defined at the level of spaces but use the enrichment of $\cat{Ho}(\sfC)$ in abelian groups. To address this we use that $k \geq 2$ to obtain maps of pointed spaces
\[s_k'' \colon S^{k-1} \wedge S^{k-1} \to S^{k-1} \wedge S^{k-1}_+ \quad \text{ and } \quad s_1'' \colon S^{k-1} \to S^{k-1} \wedge S^0_+\]
in the same way we formed $s'_k$ and $s'_1$, using that homotopy classes of maps out of a suspension obtain an abelian group structure. Using the tensoring of $\mathsf{C}$ over $\mathsf{Top}_*$, it follows that $s_k'' \wedge \mr{id}_X \simeq S^{k-1}\wedge s_k$ and $s_1'' \wedge \mr{id}_X \simeq S^{k-1}\wedge s_1$.

Translated to symmetric sequences, we must show that for each $i \geq 2$ the diagram
\[\begin{tikzcd}
(S^{k-1} \wedge \{\id\}_+) \wedge S^{k-1} \wedge \cC_k(i-1)_+ \ar{r}{\cong} \ar{d}{\mr{id} \otimes \eta_{i-1}} &[-17pt] S^{k-1} \wedge (S^{k-1} \wedge \{\id\}_+) \wedge \cC_k(i-1)_+ \ar{d}{s_k'' \wedge \mr{id}}\\
S^{k-1} \wedge \{\id\}_+ \wedge \cC_1(i-1)_+ \wedge (S^{k-1})^{\wedge i-1}  \ar{d}{s_1'' \wedge \mr{id}} & S^{k-1} \wedge S^{k-1}_+ \wedge \{\id\}_+ \wedge \cC_k(i-1)_+ \ar{d}{\mr{id} \wedge m_k \wedge \mr{id}}\\
S^{k-1} \wedge S^0_+ \wedge \{\id\}_+ \wedge \cC_1(i-1)_+ \wedge (S^{k-1})^{\wedge i-1}  \ar{d}{\mr{comm}} & S^{k-1} \wedge \cC_k(2)_+ \wedge \{\id\}_+ \wedge \cC_k(i-1)_+ \ar{d}{\mr{id} \wedge \mr{oper}_+}\\
\cC_1(i)_+ \wedge (S^{k-1})^{\wedge i}  & S^{k-1} \wedge  \cC_k(i)_+ \ar{l}[swap]{\eta_i}
\end{tikzcd}\]
commutes up to $\fS_{i-1}$-equivariant homotopy. Let us recall and explain these maps. The maps $\eta_i$ were defined in Section \ref{sec:FreeAlgMaps}, and in brief are given by
\[\eta_i(v, (e_1^1 \times e_1^{k-1}, \ldots, e_i^1 \times e_i^{k-1})) = ((e_1^1, \ldots, e_i^1),  ((e_1^{k-1})^{-1}(v), ..., (e_i^{k-1})^{-1}(v)))\]
 when the latter terms are defined, and $*$ otherwise; the map $\eta_{i-1}$ is analogous. The top horizontal isomorphism is induced by 
\begin{align*}
(S^{k-1} \wedge \{\id\}_+) \wedge S^{k-1} &\lra S^{k-1} \wedge (S^{k-1} \wedge \{\id\}_+)\\
(w, \id ; v) &\longmapsto (v; \vec{1}-w, \id).
\end{align*}
The map $\mr{comm}$ (which stands for ``commutator'') is given by
\begin{align*}
\mr{comm}(w, +1, \id, f_2, \ldots, f_{i}, t_2, \ldots, t_{i}) &= (\mu(\id, (f_2, \ldots, f_{i})), w, t_2, \ldots, t_{i})\\
\mr{comm}(w, -1, \id, f_2, \ldots, f_{i}, t_2, \ldots, t_{i}) &= (\bar{\mu}(\id, (f_2, \ldots, f_{i})), w, t_2, \ldots, t_{i}).
\end{align*}
for $\mu = (\iota_1, \iota_2), \bar{\mu}=(\bar{\iota}_1, \iota_2) \in \cC_1(2)$ defined using the little 1-cubes
\begin{align*}
\iota_1(x) = \tfrac{1}{5}x, \quad\quad \iota_2(x) = \tfrac{2}{5} + \tfrac{1}{5}x,  \quad\quad \bar{\iota}_1(x) = \tfrac{4}{5} + \tfrac{1}{5}x,
\end{align*}
as shown in Figure \ref{fig:MusAndCycle} (A). The map $\mr{oper} \colon \cC_k(2) \times \{\id\} \times \cC_k(i-1) \to \cC_k(i)$ given by the operadic composition (where we consider $\{\id\} \in \cC_k(1)$).

Up to homotopy we can take the map $m_k \colon S^{k-1} = \partial I^k \to \cC_k(2)$ to be given by $m_k(s_1, \ldots s_k) = (c_1(s_1, \ldots s_k), c_2)$ with
\begin{align*}
c_2(x_1, \ldots, x_k) &\coloneqq \tfrac{1}{5} (x_1, \ldots, x_k) + \tfrac{2}{5}\vec{1}\\
c_1(s_1, \ldots s_k)(x_1, \ldots, x_k) &\coloneqq \tfrac{1}{5} (x_1, \ldots, x_k) + \tfrac{4}{5}(s_1, s_2, \ldots, s_k)
\end{align*} as shown in Figure \ref{fig:MusAndCycle} (B).

\begin{figure}[h]
	\centering
	\begin{tikzpicture}
	
		\begin{scope}
		\node at (0,1) {(A)};
		
	\draw (0,0) -- (5,0);
	\draw[Mahogany!40!white,line width=1.6pt] (0,0) -- node[black,anchor=south] {$\iota_1$} ++ (1,0);
	\draw[Mahogany!40!white,line width=1.6pt] (2,0) -- node[black,anchor=south] {$\iota_2$} ++ (1,0);
	\node at (2.5,-0.5) {$\mu$};
	
		\draw (6,0) -- (11,0);
	\draw[Mahogany!40!white,line width=1.6pt] (10,0) -- node[black,anchor=south] {$\bar{\iota}_1$} ++ (1,0);
	\draw[Mahogany!40!white,line width=1.6pt] (8,0) -- node[black,anchor=south] {$\iota_2$} ++ (1,0);
		\node at (8.5,-0.5) {$\bar{\mu}$};
		\end{scope}
	
		\begin{scope}[yshift=-6.5cm, xshift=3cm]
		\node at (-1,4.5) {(B)};
		\draw[fill = Mahogany!10!white] (0.7,0) rectangle (1.7,1);
	\draw[fill = Mahogany!10!white] (2,2) rectangle (3,3);
	
	\draw (0,0) rectangle (5,5);
	
	\node at (1.2,0.5) {$c_1$};
	\node at (2.5,2.5) {$c_2$};
	
	\draw[dashed] (1.7,0.5) -- (4.5,0.5) -- (4.5,4.5) -- (0.5,4.5) -- (0.5,0.5) -- (0.7,0.5);

		\end{scope}

		\end{tikzpicture}
	\caption{(A) The maps $\mu$ and $\bar{\mu}$ used in the construction of $\mr{comm}$. (B) The map $m_k$.}
	\label{fig:MusAndCycle}
	\end{figure}

With this choice, and identifying $S^{k-1} = \partial I^k$, the map
\[\eta_i \circ (\mr{id} \wedge \mr{oper}_+) \circ (\mr{id} \wedge m_k \wedge \mr{id}) \colon S^{k-1} \wedge (\partial I^k)_+ \wedge \{\id\}_+ \wedge \cC_k(i-1)_+ \lra \cC_1(i)_+ \wedge (S^{k-1})^{\wedge i}\]
gives the basepoint whenever the second coordinate lies in $(\partial I^{k-1}) \times I \subset \partial I^k$, and so it factors uniquely over a map
\[\psi \colon S^{k-1} \wedge (S^{k-1}_0 \vee S^{k-1}_1) \wedge \{\id\}_+ \wedge \cC_k(i-1)_+ \lra \cC_1(i)_+ \wedge (S^{k-1})^{\wedge i}\]
with $S^{k-1}_j \coloneqq (I^{k-1} \times \{j\}) / (\partial I^{k-1} \times \{j\})$. The composition
\[S^{k-1} \wedge S^{k-1} \overset{s_k''}\lra S^{k-1} \wedge S^{k-1}_+ = S^{k-1} \wedge (\partial I^k)_+ \xrightarrow{\text{quot}} S^{k-1} \wedge (S_0^{k-1} \vee S_1^{k-1})\]
has degree $1$ on the first summand and degree $-1$ on the second. Thus the $\fS_{i-1}$-equivariant homotopy class of the clockwise composition is the difference, taken with respect to the co-$H$-space structure given by the first suspension coordinate, of the isomorphisms
\begin{align*}
(S^{k-1}_j \wedge \{\id\}_+) \wedge S^{k-1} \wedge \cC_k(i-1)_+ &\cong S^{k-1} \wedge (S^{k-1}_j \wedge \{\id\}_+) \wedge  \cC_k(i-1)_+\\
(w, \id ; v; e_2, \ldots, e_i) &\mapsto (v; \vec{1}-w, \id; e_2, \ldots, e_i).
\end{align*}
composed with the restrictions of $\psi$ to these two summands. On the other hand, the $\fS_{i-1}$-equivariant homotopy class of the anticlockwise composition is also given as the difference of two maps using the same co-$H$-space structure, namely the maps
\[\mr{comm}(-,+1, -)\circ\mr{id} \otimes \eta_{i-1} \quad \text{ and }\quad \mr{comm}(-,-1, -)\circ\mr{id} \otimes \eta_{i-1}.\]
We will show that these two pairs of maps are homotopic to each other.

The map $\psi$ on $S^{k-1} \wedge S^{k-1}_0  \wedge \{\id\}_+ \wedge \cC_k(i-1)_+$ is given at $(v; \vec{1}-w, \mr{id}; e_2, \ldots, e_{i})$ as follows. Suppose that $\eta_{i-1}(v; c_2 \circ e_2, \ldots, c_2 \circ e_{i}) = (f_2, \ldots, f_{i}; t_2, \ldots, t_{i})$ and $\eta_1(v; c_1(\vec{1}-w,0)) = (\iota_1, t_1)$. Then
\[\psi(v; w; \mr{id}; e_2, \ldots, e_{i}) = (\iota_1, f_2, \ldots, f_{i}; t_1, t_2, \ldots, t_{i}).\]
The point $t_1 \in S^{k-1} = ((0,1)^{k-1})^+ = (\bR/(-\infty,0] \cup [1,\infty))^{\wedge k-1}$ is represented by $t_1 = [5 v + 4(w - \vec{1})]$. As the remaining $t_i$ also tend to $*$ as $v$ does, 
this map is $\fS_{i-1}$-equivariantly homotopic to 
\[(v; w, \id; e_2, \ldots, e_{i}) \longmapsto (\iota_1, f_2, \ldots, f_{i}; w, t_2, \ldots, t_{i}),\]
which is $\fS_{i-1}$-equivariantly homotopic to
\[(v; w, \id; e_2, \ldots, e_{i}) \longmapsto (\mr{comm}(-,+1, -)\circ\mr{id} \otimes \eta_{i-1})(w, \mr{id}; v; e_2, \ldots, e_{i})\]
by a rescaling of the $t_2, \ldots, t_{i}$, as required.

Similarly, on $S^{k-1} \wedge S^{k-1}_1  \wedge \{\id\}_+ \wedge \cC_k(i-1)_+$ if $\eta_1(v, c_1(w,1)) = (\bar{\iota}_1, t_1)$ then we have
\[\psi(v; w, \mr{id}; e_2, \ldots, e_{i}) = (\bar{\iota}_1, f_2, \ldots, f_{i}; t_1, t_2, \ldots, t_{i}),\]
which analogously to the above is homotopic to
\[(v; w, \mr{id}; e_2, \ldots, e_{i}) \longmapsto (\mr{comm}(-,-1, -)\circ\mr{id} \otimes \eta_{i-1})(w, \mr{id}; v; e_2, \ldots, e_{i}).\]
This finishes the proof that \eqref{eq:ComparingF1Gk} commutes up to homotopy.

\subsubsection{The case $k=1$}\label{sec:module-cells-k1}
Recall that $\mr{Ass}^+$ denotes the unital associative operad, and that $\pi_0 \colon \cC_1^+ \to \mr{Ass}^+$ is a weak equivalent of operads. Thus we may replace the monad $E_1^+$ by the monad
\[\mr{Ass}^+(X) = \bigvee_{n=0}^\infty X^{\otimes n}.\]

We can define a bracket operation in terms of the multiplication $\cdot$ on $\mr{Ass}^+(X)$, the braiding on $\sfC$, and the enrichment of $\cat{Ho}(\sfC)$ in abelian groups, by
\begin{align*}
[\_, \_] \coloneqq (\_ \cdot \_) - (\_ \cdot \_) \circ \beta_{\mr{Ass}^+(X), \mr{Ass}^+(X)} \colon \mr{Ass}^+(X) \otimes \mr{Ass}^+(X) &\lra \mr{Ass}^+(X).
\end{align*}

The weak equivalence of operads $\pi_0 \colon \cC_1^+ \to \mr{Ass}^+$ induces a morphism of monads $E_1^+ \to \mr{Ass}^+$ which is a weak equivalence when evaluated on cofibrant objects. Tracing through the construction of $f_1$, we get a weakly equivalent map
\[\alpha_1 \colon \mr{Ass}^+(A) \otimes \mr{Ass}^+(\mr{Ass}^+(A) \otimes B) \lra \mr{Ass}^+(A \vee B)\]
of left $\mr{Ass}^+(A)$-modules. This is obtained from the map
\[h \colon \mr{Ass}^+(A) \otimes B = \bigvee_{n=0}^\infty A^{\otimes n} \otimes B \lra \mr{Ass}^+(A \vee B)\]
is given in terms of the bracket described above by $[\_,[\_,\cdots[\_,\_]]]$ on each summand by first extending to a map of unital associative algebras and then to a left $\mr{Ass}^+(A)$-module map. 

To understand this map, we use that the monoidal structure $\otimes$ commutes with colimits in each variable to see that
\[(A \vee B) \otimes X  \cong (A\otimes X) \vee (B\otimes X),\]
and hence can expand out $(A \vee B)^{\otimes n}$ as the coproduct of terms
\[A^{\otimes r} \otimes B \otimes A^{\otimes i_1} \otimes B \otimes A^{\otimes i_2} \otimes \cdots  \otimes B \otimes A^{\otimes i_k}\]
over all $(r; i_1, \ldots, i_k)$ with $r + k + \sum_{j=1}^k i_j=n$ and $r, i_j \geq 0$. This gives an isomorphism
\[\mr{Ass}^+(A \vee B)  \cong \mr{Ass}^+(A) \otimes \mr{Ass}^+(B \otimes \mr{Ass}^+(A))\]
of right left $\mr{Ass}^+(A)$-modules. 

With respect to this isomorphism, the map $\alpha_1$ is \emph{not} given by applying $\mr{Ass}^+(A) \otimes \mr{Ass}^+(-)$ to an isomorphism $\mr{Ass}^+(A) \otimes B \cong B \otimes \mr{Ass}^+(A)$. Rather the map of left $\mr{Ass}^+(A)$-module indecomposables,
\[Q^{\mr{Ass}^+(A)}(\alpha_1) \colon \mr{Ass}^+(\mr{Ass}^+(A) \otimes B) \lra \mr{Ass}^+(B \otimes \mr{Ass}^+(A))\]
is given by $\mr{Ass}^+(h')$ where $h' \colon \mr{Ass}^+(A) \otimes B \to B \otimes \mr{Ass}^+(A)$ is given on $A^{\otimes n} \otimes B$ by $(-1)^n$ times applying a certain braid (in the symmetric monoidal case it is given by $a_1 \otimes\cdots \otimes a_n \otimes b \mapsto (-1)^n b \otimes a_n \otimes \cdots \otimes a_1$)
so is an equivalence. 

As the source and target of $\alpha_1$ are cofibrant $\mr{Ass}^+(A)$-modules (as long as $A, B \in \sfC$ are cofibrant), their $\mr{Ass}^+(A)$-module indecomposables and derived indecomposables agree, we conclude that $\smash{Q^{\mr{Ass}^+(A)}_\bL(\alpha_1)}$ is an equivalence. By Corollary \ref{cor:AbsWhiteheadInfty}, this implies that $\alpha_1$ is an equivalence as required.

\section{$W_{k-1}$-algebras}\label{sec:Cohen}

In this paper and its sequels we will make use of computations of F. Cohen \cite[Part III]{CLM} to describe the homology of free unital $E_k$-algebras in the case that $\sfG$ is a \emph{discrete} symmetric monoidal groupoid. That is, $\sfG$ only has identity morphisms and hence simply encodes an additional grading. In a sense it is our explicit knowledge of the homology of free $E_k$-algebras that allows for some of the more intricate applications. In this section we give a careful description of these results, and from this we obtain a description of the homology of an $E_k$-cell attachment, in Section \ref{sec:homology-ekalg-coproduct}. We also study spectral sequences and operations on them in Sections \ref{sec:spectralsequences} and \ref{sec:SSStr}. As before, $\sfC = \sfS^\sfG$ with $\sfS$ satisfying the axioms of Sections \ref{sec:axioms-of-cats} and \ref{sec:axioms-of-model-cats}.

\subsection{Homology operations on $E_k$-algebras}\label{sec:homology-operations}

We wish to describe the collection of natural operations on the homology $H_{*,*}(\gR)$ of an $E_k$-algebra $\gR$ for $k \geq 2$. Recall that we have defined homology, in Section \ref{sec:homology} , in terms of a singular chain functor $C_* \colon \sfS \to \cat{A}$ where $\cat{A}$ is either the category $\cat{Ch}_\bk$ of chain complexes of $\bk$-modules, or the category $H\bk\text{-}\cat{Mod}$ of modules (in the category of symmetric spectra) over the Eilenberg--MacLane spectrum associated to the ring $\bk$. The axioms of such a singular chain functor give a lax monoidality which is a weak equivalence when evaluated on cofibrant objects, and say that when composed with $s \colon \cat{sSet} \to \sfS$ it agrees with either $C_*(-;\bk)$ or $H\bk \wedge \Sigma^\infty(-)_+$.

If $\gR$ is an $E_k$-algebra in $\sfC = \sfS^\sfG$ which is cofibrant in $\sfC$, with underlying object $R$ and structure map 
\[\alpha_R = \bigsqcup_{n \geq 1} \alpha_R(n) \colon E_k(R) = \bigsqcup_{n \geq 1} \cC_k(n) \times_{\fS_n} R^{\otimes n} \lra R,\]
then the lax monoidality of $C_*(-)$ gives $\fS_n$-equivariant maps
\[C_*(\cC_k(n)) \otimes C_*(R)^{\otimes n} \lra C_*(\cC_k(n) \times R^{\otimes n}) \xrightarrow{\,C_*(\alpha_R(n))\,} C_*(R).\]
On taking the quotient by the $\fS_n$-action these assemble to $\alpha_{C_*(R)} \colon E_k(C_*(R)) \to C_*(R)$ giving $C_*(R) \in \cat{A}^\sfG$ the structure of an $E_k$-algebra. The analogous discussion goes through for $E^+_k$-algebras. Therefore to define homology operations on $E_k$-algebras, it suffices to work in the category $\cat{A}^\sfG$.

Furthermore we have the following lemma which allows us to compute the homology of free $E_k$-algebras working just in $\cat{A}^\sfG$.
\begin{lemma}
Let $X \in \sfC$ be cofibrant. 
The natural map $C_*(X) \to C_*(E_k(X))$ has target an $E_k$-algebra, so extends to a map
\[E_k(C_*(X)) \lra C_*(E_k(X))\]
in $\cat{A}^\sfG$; this map is a weak equivalence. The analogous statement with $E^+_k$ holds too.
\end{lemma}
\begin{proof}
Both for $E_k$ and $E_k^+$ the map is given as a coproduct of the maps
\[C_*(\cC_k(n)) \otimes_{\fS_n} C_*(X)^{\otimes n} \lra C_*(\cC_k(n) \times_{\fS_n} X^{\otimes n})\]
so it is enough to see that these induce isomorphisms on homology. The $\fS_i$-actions on $\cC_k(n)$ and hence on $C_*(\cC_k(n))$ are free, so both quotients are in fact homotopy quotients. The homotopy orbit spectral sequence of Section \ref{sec:homotopy-orbit-ss} takes the form
\[E^1_{p,q} = H_p(\fS_n ; H_q(C_*(\cC_k(n)) \otimes C_*(X)^{\otimes n}))\]
in the source and 
\[E^1_{p,q} = H_p(\fS_n ; H_q(C_*(\cC_k(n) \otimes X^{\otimes n})))\]
in the target, but as the monoidality on $C_*$ is a weak equivalence on cofibrant objects the natural map between these is an isomorphism.
\end{proof}

\subsubsection{The product and Browder bracket}\label{sec:ProdBracket}
\newglossaryentry{product}{%
	name={\ensuremath{- \cdot -}},
	description={Product of a $W_{k-1}$-algebra},
	type=symbols
}
\newglossaryentry{bracket}{%
	name={\ensuremath{[- , -]}},
	description={Browder bracket of a $W_{k-1}$-algebra},
	type=symbols
}
Let $\gR$ be an $E_k$-algebra in $\cat{A}^\sfG$ with underlying object $R$. Recall that we assume that $\sfG$ is discrete. The simplest operations to define only make use of the map
\[\theta_2 \colon \cC_k(2) \times R \otimes R \lra \cC_k(2) \times_{\fS_2} (R \otimes R) \overset{\alpha_R(i)}\lra R.\]
The external product map provided by Lemma \ref{lem:KunnethFormula} (i) induces maps
\[(\theta_2)_*\colon H_{d}(\cC_k(2)) \otimes H_{g, q}(R) \otimes H_{g', q'}(R) \lra H_{g \oplus g', q+q'+d}(R).\]
Now the equivalence $S^{k-1} \overset{\sim} \to\cC_k(2)$ gives elements $u_0 \in H_{0}(\cC_k(2))$ and $u_{k-1} \in H_{k-1}(\cC_k(2))$. Using these we define the \emph{product}\index{product} $\gls{product}$ given by
\[(\theta_2)_*(u_0 \otimes - \otimes -) \colon H_{g, q}(R) \otimes H_{g', q'}(R) \lra H_{g \oplus g', q+q'}(R)\]
and the (Browder) \emph{bracket}\index{bracket} $\gls{bracket}$ given by 
\[(-1)^{(k-1)q+1}\cdot (\theta_2)_*(u_{k-1} \otimes - \otimes -) \colon H_{g, q}(R) \otimes H_{g', q'}(R) \lra H_{g \oplus g', q+q'+k-1}(R).\]
(See page 248 of \cite{CLM} to confirm this choice of sign.) We write $\mr{ad}(x)(y) \coloneqq [x,y]$.

\subsubsection{Araki--Kudo--Dyer--Lashof operations and ``top'' operations}

Now let $\bk = \bF_\ell$ be the finite field with $\ell$ elements, for a prime number $\ell$. Let $\gR$ be an $E_k$-algebra in $\cat{A}^\sfG$. The operations we describe below are constructed for $\cat{A}=\cat{Ch}_\bk$ in \cite[Chapter III]{CLM} (in fact they are written there for $\cat{Top}$, but the constructions are given on the chain level and also work when the chain complex does not arise as the singular chains on a space) and for $\cat{A} = H\bk\text{-}\cat{Mod}$ in \cite[III.\S 3]{Hinf} (in fact they are written there for the $\bk$-homology of an $\bS$-module, but go through for the homotopy of an $H\bk$-module). These definitions easily extend to $\sfG$-graded objects, since the functor $X \mapsto \bigsqcup_{g \in \sfG} X(g) \colon \cat{A}^{\sfG} \to \cat{A}$ is strong symmetric monoidal.

\newglossaryentry{qs}{%
	name={\ensuremath{Q^s}},
	description={Dyer--Lashof operation of a $W_{k-1}$-algebra},
	type=symbols
}
Suppose first that $\ell$ is odd. Then there are defined for $s \in \bZ$ and $2s-q < k-1$ \emph{Dyer--Lashof operations}\index{Dyer--Lashof operation}
	\[\gls{qs} \colon H_{g,q}(R) \lra H_{g^{\oplus \ell},q+2s(\ell-1)}(R),\]
and, for $2s-q < k-1$,
	\[\beta Q^s \colon H_{g,q}(R) \lra H_{g^{\oplus \ell},q+2s(\ell-1)-1}(R).\]
	
\begin{remark}
Note that $\beta Q^s$ does \emph{not} denote the composition of $Q^s$ with a Bockstein operation $\beta$ on the homology of $R$. Indeed, the chain complex of $\bF_\ell$-modules (or $H\bF_\ell$-module spectra) $R(g^{\oplus \ell})$ need not arise by reduction along $\bZ/\ell^2 \to \bZ/\ell = \bF_\ell$ so will not typically have Bockstein operations defined on it. However, if $R$ is obtained by reduction modulo $\ell$ from a chain complex of flat $\bZ$-modules then the Bockstein is defined and indeed $\beta Q^s = \beta \circ Q^s$, cf.\ \cite[Proposition 2.3 (v)]{MaySteenrod}. We imagine an analogous statement holds for spectra, but have not found a reference.
\end{remark}

\newglossaryentry{xi}{%
	name={\ensuremath{\xi}},
	description={Top operation of a $W_{k-1}$-algebra},
	type=symbols
}
For $q+(k-1)$ even there is defined a ``top'' operation \index{top operation}
	\[\gls{xi} \colon H_{g,q}(R) \lra H_{g^{\oplus \ell},\ell q+(k-1)(\ell-1)}(R)\]
and an associated operation
	\[\zeta  \colon H_{g,q}(R) \lra H_{g^{\oplus \ell},\ell q+(k-1)(\ell-1)-1}(R).\]

\begin{remark}
This operation is \emph{not} defined by $\zeta = \beta \circ \xi(-) - \mr{ad}^{\ell-1}(-)(\beta(-))$, as suggested on p.\ 217 of \cite{CLM}, but rather is defined on p.\ 248 of \cite{CLM}. Indeed, it is defined in situations where the Bockstein is not.
\end{remark}

Suppose now that $\ell=2$. Then there are defined for $s \in \bZ$ and $s-q < k-1$ Dyer--Lashof (rather, Araki--Kudo, but we keep to the former for uniformity) operations  
\[Q^s \colon H_{g,q}(R) \lra H_{g \oplus g,q+s}(R).\]
There is also defined a ``top'' operation 
\[\xi \colon H_{g,q}(R) \lra H_{g \oplus g,2q+(k-1)}(R).\]

\subsection{Relations among homology operations}\label{sec:RelationsAmongOps}

There are numerous relations among the operations described above, which appear in Theorems 1.1, 1.2 and 1.3 of \cite[Chapter III]{CLM} and in \cite[III.\S 3]{Hinf}. We will discuss below all those relations which do not involve Steenrod operations, which are not defined in our algebraic contexts. Consulting May \cite{MaySteenrod}, we have written out some relations which in \cite{CLM, Hinf} are left as implicit consequences of having Bockstein operations. For definiteness we consider the case of $E_k^+$-algebras.

In the following we work over $\bF_\ell$, which in the case $\ell=0$ denotes $\bQ$. We write $\cat{GrMod}_{\bF_\ell}$ for the category of graded $\bF_\ell$-vector spaces with graded tensor product given by $(V \otimes W)_k \coloneqq \bigoplus_{k'+k''=k} V_{k'} \otimes V_{k''}$ and symmetric braiding given by the Koszul sign rule $a \otimes b \mapsto (-1)^{|a||b|} b \otimes a$. We denote by $\cat{GrMod}_{\bF_\ell}^\sfG$ the category of functors $\sfG \to \cat{GrMod}_{\bF_\ell}$, equipped with the Day convolution symmetric monoidal structure.

\subsubsection{The case $\ell > 2$}\mbox{}

\noindent \textbf{Restricted $\lambda_{k-1}$-algebras.} The homology $H_*(\gR)$ equipped with the bracket $[-,-]$ forms a \emph{$(k-1)$-Lie algebra}, meaning that the bracket satisfies the following relations:
	\begin{enumerate}[(a)]	
		\item The bracket is linear in both entries.
		\item The bracket is symmetric up to a sign:
		\[[x,y] = (-1)^{|x||y|+1+(k-1)(|x|+|y|+1)} [y,x].\]
		\item The bracket satisfies the Jacobi identity up to sign:
		\begin{align*} 0=& (-1)^{(|x|+k-1)(|z|+k-1)}[x,[y,z]] \\
		&+ (-1)^{(|x|+k-1)(|y|+k-1)} [y,[z,x]] \\
		&+ (-1)^{(|y|+k-1)(|z|+k-1)} [z,[x,y]]. \end{align*}
		It also satisfies $[x,[x,x]] = 0$ (this follows from the above for $\ell > 3$, but for $\ell=3$ is new).
	\end{enumerate}

\noindent Considering also the operations $\xi$ and $\zeta$, we have the following relations:

	\begin{enumerate}[(a)] \setcounter{enumi}{3}
		\item The operation $\xi$ is not linear, but instead satisfies
		\[\xi(x+y) = \xi(x) + \xi(y) + \sum_{i=1}^{\ell-1} d^i_n(x)(y)\]
		with $d^i_n$ a certain iterated application of $\mr{ad}(x)$ and $\mr{ad}(y)$ to $x$ described on page 218 of \cite{CLM}. The operation $\zeta$ is linear \cite[Lemma 1.6]{WellingtonThesis}. For $\lambda \in \bF_\ell$ one has $\xi(\lambda x) = \lambda \xi(x)$ and $\zeta(\lambda x) = \lambda \zeta(x)$.
		\item The operations $\xi$ and $\zeta$ interact with the bracket as
		\[[x,\xi(y)] = \mr{ad}^\ell(y)(x) \qquad \text{and} \qquad [x,\zeta(y)] = 0.\]
	\end{enumerate}

\newglossaryentry{lk1}{%
	name={\ensuremath{L_{k-1}}},
	description={Free $L_{k-1}$-algebra functor},
	type=symbols
}
\newglossaryentry{dk}{%
	name={\ensuremath{D_{k}}},
	description={Free allowable Dyer--Lashof module functor},
	type=symbols
}
\newglossaryentry{vk1}{%
	name={\ensuremath{V_{k-1}}},
	description={Free allowable Dyer--Lashof algebra functor},
	type=symbols
}
\newglossaryentry{wk1}{%
	name={\ensuremath{W_{k-1}}},
	description={Free $W_{k-1}$ algebra functor},
	type=symbols
}

\noindent The resulting algebraic structure is called a \emph{restricted $\lambda_{k-1}$-algebra}. There is a \emph{free restricted $\lambda_{k-1}$-algebra} functor
	\[\gls{lk1} \colon \cat{GrMod}_{\bF_\ell}^\sfG \lra \Alg_{L_{k-1}}(\cat{GrMod}_{\bF_\ell}^\sfG)\]
	which is defined inductively as follows. Firstly, $L_0(V)$ is the free restricted Lie algebra generated by $V$, i.e.\ it is the smallest subobject of the tensor algebra $T(V)$, defined using the Day convolution formula, which contains $V$ and is closed under the bracket and under the forming of $\ell$th powers of even-degree elements. Secondly, $L_1(V) = s^{-1} L_0(sV) \oplus \zeta\cdot s^{-1}L_0(sV)_\mr{odd}$ where $s$ is the suspension defined by  $(sV)_{(g,i)} \coloneqq V_{(g,i+1)}$ (recall that $\zeta$ is only defined on odd-degree elements when $k-1=1$). For $i \geq 2$, $L_i(V) \coloneqq s^{-1} L_{i-1} (sV)$. On $L_0(V)$ the operations $\xi$ and $\zeta$, defined on even-degree elements, are respectively $x \mapsto x^\ell$ and $0$. On $L_1(V)$ they are defined as $\xi(x) = s^{-1} \xi(sx)$ (when $sx \in L_0(sV)$ has odd degree) and $\zeta(x) = \zeta \cdot x$ (when $sx \in L_0(sV)$ has odd degree). Note that as $\zeta$ always produces elements of even degree one can never form $\xi \zeta$ or $\zeta \zeta$, so this gives a complete description of these operations. For $i \geq 2$, we inductively define the operations by $\xi(x) = s^{-1} \xi(sx)$ and $\zeta(x) = -s^{-1}\zeta(sx)$ in terms of $L_i(V) = s^{-1} L_{i-1} (sV)$.\\
	
\noindent \textbf{Allowable Dyer--Lashof algebras.} The Dyer--Lashof operations satisfy the following relations:
\begin{enumerate}[(a')]
		\item The Dyer--Lashof operations are linear.
		\item \label{enum.dltoosmall} A Dyer--Lashof operation vanishes if the degree of $x$ is too large: $Q^s x = 0$ if $2s<|x|$ and $\beta Q^s x =0$ if $2s \leq |x|$.
		\item \label{enum.adem} The Dyer--Lashof operations satisfy the Adem relations. That is, if $r > \ell s$ we have
\[Q^rQ^s = \sum_{i} (-1)^{r+i} {\ell i -(\ell-1)s-i-1 \choose r-(\ell-1)s-i-1} Q^{r+s-i}Q^i\]
and
\[\beta Q^rQ^s = \sum_{i} (-1)^{r+i} {\ell i -(\ell-1)s-i-1 \choose r-(\ell-1)s-i-1} \beta Q^{r+s-i}Q^i,\]
and if $r \geq \ell s$ we have 
\begin{align*}
		Q^r \beta Q^s &= \sum_{i} (-1)^{r+i} {\ell i-(\ell-1)s-i \choose r-(\ell-1)s-i} \beta Q^{r+s-i}Q^i\\
 & \quad - \sum_i (-1)^{r+i} \binom{\ell i + (\ell-1)s-i-1}{r-(\ell-1)s-i} Q^{r+s-i} \beta Q^i
\end{align*}
and
\[\beta Q^r \beta Q^s = - \sum_i (-1)^{r+i} \binom{\ell i + (\ell-1)s-i-1}{r-(\ell-1)s-i} \beta Q^{r+s-i} \beta Q^i.\]
	\end{enumerate}	
\noindent The resulting algebraic structure is called a \emph{Dyer--Lashof module}. There is a free allowable Dyer--Lashof module functor
	\[\gls{dk} \colon \cat{GrMod}_{\bF_\ell}^\sfG \lra \cat{Alg}_{D_k}(\cat{GrMod}_{\bF_\ell}^\sfG)\]
	which is defined as follows: $D_k(V)$ is the quotient of the graded $\bF_\ell$-vector space generated by words in the Dyer--Lashof operations applied to elements of $V$, modulo the relations (a')-(c').

Considering the interaction of the Dyer--Lashof operations and the product, we have the following further relations:
	\begin{enumerate}[(a')] \setcounter{enumi}{3}
		\item \label{enum.dlequal} A Dyer--Lashof operation is an $\ell$-fold power in the critical degree: $Q^s x = x^\ell$ if $2 s = |x|$.
		\item \label{enum.dl1} Dyer--Lashof operation of non-zero degree vanish on the unit: if $1 \in H_{\bunit_\sfG,0}(\gR)$ is the identity element, then $Q^s 1 = 0$ if $s \neq 0$, and $ \beta Q^s 1=0$ for all $s$. 
		\item The Dyer--Lashof operations satisfy the Cartan formula:
		\[Q^s(xy) = \sum_{i+j=s} (Q^i x)(Q^j y)\]
and 
\[\beta Q^{s+1}(xy) = \sum_{i+j=s} (\beta Q^{i+1} x)(Q^j y) + (-1)^{|x|}(Q^i x)(\beta Q^{j+1} y).\]
	\end{enumerate}
\noindent The resulting algebraic structure is called an \emph{allowable Dyer--Lashof algebra}. There is a free allowable Dyer--Lashof algebra functor 
	\[\gls{vk1} \colon \cat{GrMod}_{\bF_\ell}^\sfG \lra \Alg_{V_{k-1}}(\cat{GrMod}_{\bF_\ell}^\sfG)\] 
	given by sending $V$ to the free graded-commutative algebra on $D_k(V)$ and taking the quotient by the ideal generated by $x^\ell - Q^s(x)$ for $|x| = 2s$. \\

	\noindent \textbf{$W_{k-1}$-algebras.} Finally we combine the allowable Dyer--Lashof algebra structure with the restricted $\lambda_{k-1}$-algebra structure. That means we need to describe the relations between $[-,-]$, $\xi$ and $\zeta$ on the one hand, and the product, $Q^s$ and $\beta Q^s$ on the other hand.
	
\begin{enumerate}[(a'')]
		\item The bracket is a derivation of the product in each variable, up to a sign: 
		\[[x,yz] = [x,y]z + (-1)^{|y|(k-1+|x|)} y[x,z].\]
		\item The bracket with the unit vanishes: $[1,x] = 0$ if $1 \in H_{\bunit_\sfG,0}(\gR)$ is the identity element.
		\item A bracket with a Dyer--Lashof operation vanishes: $[x,Q^s y] = 0 = [x,\beta Q^s y]$.
		\item To incorporate the operations $\xi$ and $\zeta$, one can observe that they behave something like Dyer--Lashof operations. Note by the restrictions on $s$ in the definitions of $Q^s$ and $\beta Q^s$, currently we have not defined $Q^{(|x|+k-1)/2}(x)$ and $\beta Q^{(|x|+k-1)/2}(x)$. However, the effect on bidegrees of $\xi$ and $\zeta$ coincides with such hypothetical operations. Hence by convention we shall define
\[Q^{(|x|+k-1)/2}(x) \coloneqq \xi \qquad \text{and} \qquad \beta Q^{(|x|+k-1)/2}(x) \coloneqq \zeta.\]
		Using this notation, $\zeta$ satisfies the relations (a')-(f'), and $\xi$ satisfies the relations (b')-(e') but (a') must be replaced as explained in (d) 
and (f') must be replaced by
\[\xi(xy) = \sum_{i+j=\frac{1}{2}(k-1+|x|+|y|)} (Q^ix)(Q^jy) + \sum_{0 \leq i,j \leq \ell} x^i y^j \Gamma_{ij}\]
with $\Gamma_{ij}$ a certain function of $x$ and $y$ defined on page 335 of \cite{CLM}. 
\end{enumerate}
	
The combined algebraic structure of an restricted $\lambda_{k-1}$-algebra and an allowable Dyer--Lashof algebra, satisfying these conditions, is called a \emph{$W_{k-1}$-algebra}\index{$W_{k-1}$-algebra}. There is a free $W_{k-1}$-algebra functor
\[\gls{wk1}  \colon \cat{GrMod}_{\bF_\ell}^\sfG \lra \Alg_{W_{k-1}}(\cat{GrMod}_{\bF_\ell}^\sfG)\]
given as follows: it is the quotient of $V_{k-1}(L_{k-1}(V))$ by the ideal generated by the relations (a'')-(d'').

Unwinding the definitions of $V_{k-1}$ and $L_{k-1}$ allows one to write down a basis of a free $W_{k-1}$-algebra on a $V \in \cat{GrMod}_{\bF_\ell}^\sfG$. Something like following description appears on page 227 of \cite{CLM} (see also \cite[Section I.1]{WellingtonBook}), but we have fixed a number of unfortunate typos. It is obtained by writing $Q^{(|x|+k-1)/2}(x) \coloneqq \xi(x)$ and $\beta Q^{(|x|+k-1)/2}(x) \coloneqq \zeta(x)$, and recalling that $L_{k-1}(V)$ is obtained by applying basic Lie words or applications of $\xi$ or $\zeta$ to basic Lie words, see e.g.\ \cite[Section II.2]{bourbakilie} on how to obtain these. The result is the free graded-commutative algebra generated by those $Q^I(y)$ (which includes the shorthand for $\xi$ and $\zeta$) such that
\begin{itemize}
	\item $y$ is a basic Lie word in a basis of $V$,
	\item $I = (\epsilon_1,s_1,\ldots,\epsilon_r,s_r)$ is admissible (i.e.\ $\ell s_j - \epsilon_j \geq s_{j-1}$ for $2 \leq j \leq r$), $e(I) + \epsilon_1 > |y|$ where $e(I) = 2s_1-\epsilon_1-\sum_{j=2}^r (2s_j(\ell-1)-\epsilon_j)$, and $2s_r \leq |y|+k-1$.
\end{itemize} 

\subsubsection{The case $\ell=2$}\mbox{}

\noindent \textbf{Restricted $\lambda_{k-1}$-algebras.} The bracket and $\xi$ satisfy the following relations, describing the structure of a \emph{restricted $\lambda_{k-1}$-algebra} (for $\ell = 2$).	
\begin{enumerate}[(a)]	
		\item The bracket is linear in both entries.
		\item The bracket satisfies $[x,x]=0$.
		\item The bracket satisfies the Jacobi identity:
		\[0= [x,[y,z]] + [y,[z,x]]+ [z,[x,y]].\]
		\item The operation $\xi$ is not linear, but satisfies
		\[\xi(x+y) = \xi(x) + \xi(y)+[y,x].\]
		\item The bracket and $\xi$ interact as
		\[[x,\xi(y)] = [y,[y,x]].\]
\end{enumerate}

\noindent  \textbf{Allowable Dyer--Lashof algebras.} The product and Dyer--Lashof operations satisfy the following relations, describing an \emph{allowable Dyer--Lashof algebra} (for $\ell = 2$).
	\begin{enumerate}[(a')]
		\item The Dyer--Lashof operations are linear.
		\item A Dyer--Lashof operation vanishes if the degree of $x$ is too large: $Q^s x = 0$ if $s<|x|$.
		\item The Dyer--Lashof operations satisfy the Adem relations. That is, if $r > 2s$, we have
		\[Q^rQ^s = \sum_{i} {2i-r \choose r-s-i-1} Q^{r+s-i}Q^i.\]
		\item A Dyer--Lashof operation is squaring in the critical degree: $Q^s x = x^2$ if $s = |x|$.
		\item A Dyer--Lashof operation of non-zero degree vanishes on the unit: $Q^s 1 = 0$ if $s \neq 0$ and $1 \in H_{\bunit_\sfG,0}(\gR)$ is the identity element.
		\item The Dyer--Lashof operations satisfy the Cartan formula:
		\[Q^s(xy) = \sum_{i+j=s} (Q^i x)(Q^j y).\]
	\end{enumerate}
	
\noindent \textbf{$W_{k-1}$-algebras.} Finally, the compatibility between the restricted $\lambda_{k-1}$-algebra and allowable Dyer--Lashof algebra structures is described by the following set of relations, leading to the structure of a \emph{$W_{k-1}$-algebra}\index{$W_{k-1}$-algebra} (for $\ell = 2$).
	
	\begin{enumerate}[(a'')]
		\item The bracket is a derivation of the product in each variable: $[x,yz] = [x,y]z + y[x,z]$.
		\item The bracket with the unit vanishes: $[1,x] = 0$ if $1 \in H_{\bunit_\sfG,0}(\gR)$ is the identity element.
		\item A bracket with a Dyer--Lashof operation vanishes: $[x,Q^s y] = 0$.
		\item To incorporate the operations $\xi$, one defines $Q^{|x|+k-1}(x) \coloneqq \xi$. Using this notation, the relations (b')-(e') hold, but (a') must be replaced by (d) 
and (f') must be replaced by
		\[\xi(xy) = \sum_{i+j=k-1+|x|+|y|} (Q^ix)(Q^jy) + x[x,y]y.\]
	\end{enumerate}

\noindent There is again a free $W_{k-1}$-algebra, using these definitions.

As above, for a $V \in \cat{GrMod}_{\bF_\ell}^\sfG$ we find (see \cite[p.\ 227]{CLM}) that $W_{k-1}(V)$ is the free commutative algebra generated by those $Q^I(y)$ (which includes the shorthand $Q^{|x|+k-1}(x) \coloneqq \xi(x)$) such that
\begin{itemize}
	\item $y$ is a basic Lie word in a basis of $V$,
	\item $I = (s_1,\ldots,s_r)$ is admissible (i.e.\ $2s_j \geq s_{j-1}$ for $2 \leq j \leq r$), $e(I) > |y|$ where $e(I) = s_1 - \sum_{j=2}^r s_j$, and $s_r \leq |y|+k-1$.
\end{itemize}

\subsubsection{The case $\ell=0$} 
In this case only the product and bracket are defined, and they satisfy the following relations:
\begin{enumerate}[(a)]	
	\item The bracket is linear in both entries.
	\item The bracket is symmetric up to a sign: $[x,y] = (-1)^{|x||y|+1+(k-1)(|x|+|y|+1)} [y,x]$.
	\item The bracket satisfies the Jacobi identity up to sign:
	\begin{align*} 0= (-1)^{(|x|+k-1)(|z|+k-1)}[x,[y,z]] & + (-1)^{(|x|+k-1)(|y|+(k-1))} [y,[z,x]] \\
	&+ (-1)^{(|y|+k-1)(|z|+k-1)} [z,[x,y]] \end{align*}
\end{enumerate}
\begin{enumerate}[(a'')]
	\item The bracket is a derivation up to a sign: $[x,yz] = [x,y]z + (-1)^{|y|(k-1+|x|)} y[x,z]$.
	\item The bracket with the unit vanishes: $[1,x] = 0$ if $1 \in H_{\bunit_\sfG,0}(\gR)$ is the identity element.
\end{enumerate}

This algebraic structure is also known as \emph{$(k-1)$-Gerstenhaber algebra}.\index{Gerstenhaber algebra} Let $W_{k-1}(V)$ denote the free $(k-1)$-Gerstenhaber algebra on a $V \in \cat{GrMod}_{\bF_\ell}^\sfG$. It may be constructed as a quotient of the free graded-commutative algebra on the free $(k-1)$-Lie algebra $L_{k-1}(V)$, by enforcing the above relations.

\subsubsection{Non-unital $E_k$-algebras}\label{sec:ReducedW}

We have so far discussed the relations which hold between homology operations on the homology of an $E_k^+$-algebra, but all of these operations are defined using only the $E_k$-algebra structure. The only relations which do not make sense for $E_k$-algebras are those involving the unit $1 \in H_{\bunit_\sfG, 0}(\gR)$. Let us define a variant of the construction $W_{k-1}$ adapted to $E_k$-algebras as follows. We have $W_{k-1}(0) = \bF_\ell[\bunit]$, where $\bF_\ell[\bunit]$ was shorthand for $(\bunit_\sfG)_*(\bF_\ell)$. As any object $V \in \cat{GrMod}_{\bF_\ell}^\sfG$ has a canonical morphism $z \colon V \to 0$, we may define
\[\tilde{W}_{k-1}(V) \coloneqq \mathrm{ker}(W_{k-1}(V) \to W_{k-1}(0))\]
as the kernel of this canonical augmentation. This inherits the structure of a monad.

\subsection{The homology of free $E_k$-algebras}\label{sec.homologyfreeekalgebras} 

The main theorems of F. Cohen's contribution to \cite{CLM} say that for $X \in \mathsf{Top}_\ast$, the homology groups $H_*(\gE_k^+(X);\bF_\ell)$ are the free $W_{k-1}$-algebra on $\tilde{H}_*(X;\bF_\ell)$. Here we shall show that the same is true for $X \in \cat{A}^\sfG$, and hence for $X \in \sfC = \sfS^\sfG$.

\begin{theorem}\label{thm.wkfreealg}
Let $\bk = \bF_\ell$ and $k \geq 2$. For an object $X \in \cat{A}^{\sfG}$ the natural map
	\[W_{k-1}({H}_{*, *}(X)) \lra H_{*, *}(\gE_k^+(X))\]
is an isomorphism of $W_{k-1}$-algebras, and the natural map
	\[\tilde{W}_{k-1}({H}_{*, *}(X)) \lra H_{*, *}(\gE_k(X))\]
is an isomorphism of $\tilde{W}_{k-1}$-algebras.
\end{theorem}
\begin{proof}
We may verify this after applying $X \mapsto \bigsqcup_{g \in \sfG} X(g) \colon \cat{A}^{\sfG} \to \cat{A}$, so we may forget about $\sfG$. Cohen has shown \cite[III Theorem 3.1]{CLM} that the first claim holds for any pointed space, and hence simplicial set, $X$ (in which case homology is to be interpreted as reduced homology). Thus it holds for all objects of $\cat{A}$ in the essential image of the functor
\begin{align*}
X \mapsto \Sigma^\infty X_+ \colon \cat{Ho}(\cat{sSet}) &\lra \cat{Ho}(H\bk\text{-}\cat{Mod}) \qquad \text{or}\\
X \mapsto C_*(X;\bk)  \colon \cat{Ho}(\cat{sSet}) &\lra \cat{Ho}(\cat{Ch}_\bk).
\end{align*}
As $\bk$ is a field, this essential image consists precisely of the 0-connective objects.

Now suppose that $X \in \cat{A}$ is bounded below, so that $S^{2N} \otimes X$ is 0-connective for some $N \gg 0$. The permutation action of $\fS_r$ on $(S^{2N})^{\otimes r}$ is homotopically trivial, i.e.\ there is a zig-zag of $\fS_r$-equivariant weak equivalences from $(S^{2N})^{\otimes r}$ to $S^{2Nr}$ with the trivial action, so there is a weak equivalence
\[\cC_k(r) \times_{\fS_r} (S^{2N} \otimes X)^{\otimes r} \simeq S^{2Nr} \otimes (\cC_k(r) \times_{\fS_r} X^{\otimes r}).\]
See \cite[VII.\S3]{Hinf} for a similar discussion. If we work in $\cat{A}^\bN$, where $X$ is placed in grading 1, the above says that
\[H_{r,d}(\gE^+_k(S^{2N} \otimes X);\bk) \cong H_{r, d-2Nr}(\gE_k^+(X);\bk),\]
and by the 0-connective case the left-hand side may be identified with the part of ${W}_{k-1}(H_{1,*}(S^{2N} \otimes X;\bk))$ of bidegree $(r,d)$. Considering the definition of ${W}_{k-1}$, one sees that this is isomorphic to the part of ${W}_{k-1}(H_{1,*}(X);\bk)$ of bidegree $(r,d-2Nr)$. In terms of the bases we have described, the isomorphism is given as follows: if $y$ is a basic Lie word of length $\rho$ in a basis of $H_{1,*}(S^{2N} \otimes X;\bk)$, and $y'$ is the corresponding basic Lie word in the corresponding basis of $H_{1,*}(X;\bk)$, then 
\[Q^{(\epsilon_1, s_1, \ldots, \epsilon_k, s_k)}(y) \mapsto Q^{(\epsilon_1, s_1-2N\ell^{k-1}\rho, \ldots\epsilon_{k-1}, s_{k-1}-2N\ell\rho, \epsilon_k, s_k-2N\rho)}(y').\]
Putting this together it follows that $W_{k-1}(H_*(X;\bk)) \to H_*(\gE_k^+(X);\bk)$ is an isomorphism when $X$ is bounded below. But it is then an isomorphism for general $X$, as any $X$ is a filtered colimit of bounded below objects, and both sides commute with filtered colimits.

The second claim follows from the first, using the decompositions $W_{k-1}(V) \cong \bF_\ell[\bunit] \oplus \tilde{W}_{k-1}(V)$ and $H_{*,*}(\gE_k^+(X)) \cong \bF_\ell[\bunit] \oplus H_{*,*}(\gE_k(X))$, which are natural.
\end{proof}

\subsection{Coproducts of $E_\infty$-algebras}\label{sec:homology-ekalg-coproduct}
Combining Proposition \ref{prop:einty-coporoduct} with Lemma \ref{lem:KunnethFormula} (i) gives the following.

\begin{corollary}\label{cor:homology-einfty-coproduct}
Let $\gR$ and $\gS$ be $E_\infty^+$-algebras in $\sfC= \sfS^\sfG$ with $\sfG$ discrete, cofibrant in $\sfC$. Then taking homology in a field we have an isomorphism 
	\[H_{*,*}(\gR \sqcup^{E_\infty^+} \gS) \cong H_{*,*}(\gR) \otimes H_{*,*}(\gS)\]
of objects of $\Alg_{W_{\infty}}(\mathsf{M}^{\sfG})$.
\end{corollary}

\begin{remark}
In the case of $E_k^+$-algebras with $k < \infty$ we do not know an algebraic formula for $H_{*,*}(\gR \sqcup^{E_k^+} \gS)$ in terms of the $W_{k-1}$-algebras $H_{*,*}(\gR)$ and $H_{*,*}(\gS)$. The na{\"i}ve guess, as the coproduct of $W_{k-1}$-algebras, is false even when $\gS = \gE_k^+(S^{g,d})$ and working over $\bQ$. We thank the referee for explaining this to us.
\end{remark}

\subsection{$E^1$-pages of spectral sequences}\label{sec:spectralsequences}

The cell attachment spectral sequence of Corollary \ref{cor:cell-attachment-spectral-sequence} and the skeletal spectral sequence of Corollary \ref{cor:skeletalss} apply to the $E_k^+$-operad and so can be used to compute the homology of cellular $E_k^+$-algebras. In this and the following section we develop some basic tools for such calculations in the case that $\sfG$ is a discrete groupoid: we describe the $E^1$-pages of these spectral sequences and discuss how differentials interact with the product, bracket, and Dyer--Lashof operations.

\subsubsection{The $E^1$-page of the skeletal spectral sequence}\label{sec:spectralsequences.skel}
 For the skeletal spectral sequence of a CW $E_k^+$-algebra $\gZ$ in $\sfC$ as described in Corollary \ref{cor:skeletalss}, the $E^1$-page is given by the homology of the free $E_k^+$-algebra on a wedge of spheres corresponding to the cells. When $\sfG$ is a discrete groupoid Theorem \ref{thm.wkfreealg} explains how F.\ Cohen's results give a description of the homology of free $E_k^+$-algebras with coefficients in $\bF_\ell$, where $\ell$ is either prime or $0$ (in which case $\bF_\ell$ means $\bQ$), in terms of the functor $W_{k-1}$. Thus in this case the skeletal spectral sequence takes the form
\[E^1_{*,*,*} \cong W_{k-1}\left(\bigoplus_{d \geq 0} \bigoplus_{a \in I_d}  \bF_{\ell}[g_\alpha,d,d]\right) \Longrightarrow H_{*,*}(\gZ;\bF_{\ell})\]
with $\bF_{\ell}[g_\alpha,d,d]$ replacing the more cumbersome notation $d_*\bF_{\ell}[g_\alpha,d]$.

\subsubsection{The $E^1$-page of the cell attachment spectral sequence for $E_\infty^+$-algebras} 

If $\gR$ is a $E_\infty^+$-algebra in $\sfC$ which is cofibrant in $\sfC$, then the cell attachment spectral sequence of Corollary \ref{cor:cell-attachment-spectral-sequence} for $\gR \cup^{E_\infty^+}_f \gD^{g,d}$ is given by
\[E^1_{g,p,q} \cong \tilde{H}_{g,p+q,p}(0_* \gR_+ \vee^{E_\infty^+} \gE_\infty^+(1_* S^{g,d})) \Longrightarrow H_{g,p+q}(\gR \cup^{E_\infty^+}_f \gD^{g,d}),\]
with the $E^1$-page the homology of an $E_\infty^+$-algebra in $\sfC_*^{\bN_{=}}$. Under the additional assumptions that $\gR$ is cofibrant in $\Alg_{E^+_\infty}(\sfC)$ and that $\sfG$ is a discrete groupoid, it follows from Corollary \ref{cor:homology-einfty-coproduct} that with coefficients in the field $\bF_\ell$ there is an isomorphism 
	\[\tilde{H}_{*,*,*}(0_* \gR_+ \vee^{E_\infty^+} \gE_\infty^+(1_* S^{g,d})) \cong 0_* H_{*,*}(\gR) \otimes_{\bF_\ell} W_\infty(\bF_\ell[g,d,1]).\]

\subsubsection{The $E^1$-page of the spectral sequence of Theorem \ref{thm:e2cellsIV-ss}}

Recall that for a map $f \colon \gR \to \gS$ between non-unital $E_\infty$-algebras, satisfying certain hypotheses, Theorem \ref{thm:e2cellsIV-ss} provides a strongly convergent spectral sequence
\[E^1_{g,p,q} = \widetilde{H}_{g,p+q,p}(E_\infty^+((-1)_* Q^{E_\infty}_\bL(\gS)/Q^{E_\infty}_\bL(\gR));A) \Longrightarrow H_{g,p+q}^{\overline{\gR}}(\overline{\gS};A).\]
When we take coefficients in the field $\bk = \bF_\ell$, the $E^1$-page can be described as
\[W_\infty((-1)_* H^{E_\infty}_{*,*}(\gS,\gR;\bF_\ell)).\]
In particular, for $\ell=0$ we get a free graded-commutative algebra on generators $(-1)_* H^{E_\infty}_{*,*}(\gB, \gA;\bQ)$.

\subsection{$W_{k-1}$-algebra structures on spectral sequences}\label{sec:SSStr}

In this section we consider homology with coefficients in a field $\bF$. Let $\gR \in \Alg_{E_k}(\sfC^{\bZ_{\leq}})$ be an $E_k$-algebra with an ascending filtration, and suppose that the underlying filtered object $U^{E_k}\gR \in \sfC^{\bZ_{\leq}}$ is cofibrant. Then Theorem \ref{thm:SSAsc} gives a spectral sequence
\[E^1_{g,p,q}(\gR) = \widetilde{H}_{g,p+q,p}(\grr(U^{E_k}\gR);\bF) \Longrightarrow H_{g,p+q}(\colim U^{E_k}\gR;\bF),\]
and in the following subsections we wish to discuss the additional structure present on this spectral sequence arising from $\gR$ being an $E_k$-algebra.

\subsubsection{Multiplicative structures}

In this subsection we consider homology with coefficients in any commutative ring $\bk$. We do not need to assume that $\sfG$ is a discrete groupoid, but we do assume that it satisfies the hypothesis of Lemma \ref{lem:KunnethFormula} (i) so that the external product is available.

\begin{theorem}
There are operations
\begin{align*}
	- \cdot_r - \colon E^r_{g_1,p_1,q_1}(\gR) \otimes_\bk E^r_{g_2,p_2,q_2}(\gR) &\lra E^r_{g_1\oplus g_2,p_1+p_2,q_1+q_2}(\gR)\\
	[-, -]_r \colon E^r_{g_1,p_1,q_1}(\gR) \otimes_\bk E^r_{g_2,p_2,q_2}(\gR) &\lra E^r_{g_1\oplus g_2,p_1+p_2,q_1+q_2+(k-1)}(\gR).
\end{align*}
For $x \in E^r_{g_1,p_1,q_1}(\gR)$ and $y \in E^r_{g_2,p_2,q_2}(\gR)$ these satisfy
\begin{align*}
	d^r(x \cdot_r y) &= d^r(x) \cdot_r y + (-1)^{p_1+q_1} x \cdot_r d^r(y)\\
	d^r([x , y]_r) &= [d^r(x) , y]_r + (-1)^{(k-1)+(p_1+q_1)} [x , d^r(y)]_r,
\end{align*}
and if $d^r(x) = 0 = d^r(y)$ so that they represent classes $\bar{x} \in E^{r+1}_{g_1,p_1,q_1}(\gR)$ and $\bar{y} \in E^{r+1}_{g_2,p_2,q_2}(\gR)$ then
\begin{align*}
	\bar{x} \cdot_{r+1} \bar{y} &= \overline{x \cdot_r y}\\
	[\bar{x}, \bar{y}]_{r+1} &= \overline{[x,y]_r}.
\end{align*}

On $E^1_{*,*,*}(\gR)$ and $E^\infty_{*,*,*}(\gR)$ these operations are induced by the product and bracket on $H_{*,*,*}(\grr(\gR);\bk)$ and $H_{*,*}(\gR;\bk)$ respectively.
\end{theorem}
\begin{proof}
As part of the $E_k$-algebra structure on $\gR$ we have a morphism
\[\theta_2 \colon \cC_k(2) \times \gR \otimes \gR \lra E_k(\gR) \lra \gR\]
between cofibrant objects in $\sfC^{\bZ_{\leq}}$ (we omit the notation $U^{E_k}$ for clarity from now on), which induces a map of spectral sequences. There is a morphism of spectral sequences
\[H_{*}(\cC_k(2);\bk)[0,0] \otimes_\bk E^r_{*,*,*}(\gR) \otimes_\bk E^r_{*,*,*}(\gR) \lra E^r_{*,*,*}(\cC_k(2) \times \gR \otimes \gR) ,\]
given (using the external product provided by Lemma \ref{lem:KunnethFormula} (i)) by a morphism of defining exact couples. Thus there is a map of spectral sequences
\[(\theta_2)_* \colon H_{*}(\cC_k(2);\bk)[0,0] \otimes_\bk E^r_{*,*,*}(\gR) \otimes_\bk E^r_{*,*,*}(\gR) \lra E^r_{*,*,*}(\gR).\]

Recall from Section \ref{sec:ProdBracket} that using the equivalence $S^{k-1} \overset{\sim}\to \cC_k(2)$ and the canonical classes $u_0 \in H_0(S^{k-1};\bk)$ and $u_{k-1} \in H_{k-1}(S^{k-1};\bk)$, the product and bracket are defined on $X = U^{E_k} \colim \gR$ by $(\theta_2)_*(u_0 \otimes - \otimes -)$ on $H_{g_1, p_1+q_1}(X) \otimes H_{g_2, p_2+q_2}(X)$ and by $(-1)^{(k-1)(p_1+q_1)+1}\cdot (\theta_2)_*(u_{k-1} \otimes - \otimes -)$ on $H_{g_1, p_1+q_1}(X) \otimes H_{g_2, p_2+q_2}(X)$. Therefore defining $- \cdot_r -$ on $E^r_{g_1,p_1,q_1}(\gR) \otimes_\bk E^r_{g_2,p_2,q_2}(\gR)$ by $(\theta_2)_*(u_0 \otimes - \otimes -)$, and defining $[-,-]_r$ on $E^r_{g_1,p_1,q_1}(\gR) \otimes_\bk E^r_{g_2,p_2,q_2}(\gR)$ by $(-1)^{(k-1)(p_1+q_1)+1}\cdot (\theta_2)_*(u_{k-1} \otimes - \otimes -)$, we have the desired properties.
\end{proof}

\subsubsection{Dyer--Lashof operations}

In this subsection we consider homology with coefficients in a prime field $\bF_\ell$ with $\ell>0$, and we assume that $\sfG$ is a discrete groupoid.

\begin{theorem}
Let $\gR \in \Alg_{E_\infty}(\sfC^{\bZ_{\leq}})$ be an $E_\infty$-algebra with an ascending filtration. If $x \in E^1_{g,p,q}(\gR)$ survives to $E^r_{g,p,q}(\gR)$ and $d^r([x]) = [y]$ for $y \in E^1_{g, p-r,q+r-1}(\gR)$, then
\begin{enumerate}[(i)]
\item $Q^s(x)$ survives to $E^{\ell r}_{*,*,*}$ and $d^{\ell r}([Q^s(x)])$ is represented by $Q^s(y)$,

\item $\beta Q^s(x)$ survives to $E^{\ell r}_{*,*,*}$ and $d^{\ell r}([\beta Q^s(x)])$ is represented by $-\beta Q^s(y)$.
\end{enumerate}

If $x \in E^1_{g,p,q}(\gR)$ survives to $E^\infty_{g,p,q}(\gR)$ and represents $z \in H_{g, p+q}(\colim\gR;\bF_\ell)$, then
\begin{enumerate}[(i)]
\item $Q^s(x)$ survives to $E^{\infty}_{*,*,*}$ and represents $Q^s(z)$,

\item $\beta Q^s(x)$ survives to $E^{\infty}_{*,*,*}$ and represents $\beta Q^s(z)$.
\end{enumerate}
\end{theorem}

\begin{proof}
In the case $\cat{A} = \mathsf{Ch}_{\bF_\ell}$ 
the first part of the theorem follows from the definition \cite[p.\ 7]{CLM} of the operations $Q^s$ and $\beta Q^s$ along with May's general approach to Steenrod operations, specifically \cite[Proposition 3.5]{MaySteenrod} with $f(a,b,c) = a$. The argument given there can be mimicked in $\cat{A}= H\bk\text{-}\cat{Mod}$, as in \cite[III\S 1]{Hinf}.

For the second part, that $Q^s(x)$ and $\beta Q^s(x)$ survive to $E^\infty_{*,*,*}$ follows from the first part. That they represent the claimed elements is immediate from the construction of the spectral sequence.
\end{proof}

The arguments of \cite{MaySteenrod} can be adapted to study $E_k$-algebras, where modified formulae can be obtained for how $Q^s$ and $\beta Q^s$ interact with differentials, as well as for how $\xi$ and $\zeta$ do. For example, for $\ell=2$ one finds that $d^r([\xi(x)])$ is represented by $[y,x]$. We will not develop those formulae, as we have---as yet---no need for them.

\subsection{Derived $\tilde{W}_{k-1}$-indecomposables}
If $\gR$ is an $E_k$-algebra which is cofibrant in $\sfC$, then one may use the equivalence $Q^{E_k}_\bL(\gR) \simeq B(\bunit, E_k, \gR)$ of Section \ref{sec:simplicial-formula-indecomposables}, the geometric realization spectral sequence of Theorem \ref{thm:geom-rel-ss}, and Theorem \ref{thm.wkfreealg} to obtain a spectral sequence
\[E^1_{*, p, *} = (\tilde{W}_{k-1})^p(H_{*,*}(\gR)) \Longrightarrow H^{E_k}_{*, *}(\gR).\]
One may consider the $E^1$-page to be obtained by taking the $\tilde{W}_{k-1}$-indecomposables $Q^{\tilde{W}_{k-1}}(-)$ of the canonical simplicial resolution $(\tilde{W}_{k-1})^{\bullet+1}(H_{*,*}(\gR)) \to H_{*,*}(\gR)$, so tautologically one has $E^2_{*, p, *} = (\bL_p Q^{\tilde{W}_{k-1}})(H_{*,*}(\gR))$, the (simplicially) derived $\tilde{W}_{k-1}$-indecomposables of the $\tilde{W}_{k-1}$-algebra $H_{*,*}(\gR)$.

We will not make use of this spectral sequence, but it has been studied in some detail by Richter--Ziegenhagen \cite{RichterZiegenhagen}, especially in even characteristic.


\part{A framework for examples} \label{part:framework}  In this final part we use the techniques developed in parts \ref{part:category}, \ref{part:homotopy}, and \ref{part:ek} to prove results related to homological stability, in the setting of $E_k$-algebras arising from monoidal groupoids. We show that in this setting the derived $E_1$-indecomposables can be computed in terms of a semi-simplicial set of splittings. We then describe the relationship to Koszul duality and give a generic homological stability results for both constant and local coefficients. As before, $\sfC = \sfS^\sfG$ with $\sfS$ satisfying the axioms of Sections \ref{sec:axioms-of-cats} and \ref{sec:axioms-of-model-cats}. Often we shall also assume Axiom \ref{axiom:Hurewicz} (which roughly says that in the category $\sfC$ homotopy groups detect weak equivalences, and there is a Hurewicz theorem which holds in all degrees), so that we have a Hurewicz theorem and CW approximation theorem in the category of $E_k$-algebras in $\sfC$.

\section{$E_k$-algebras from monoidal groupoids}\label{sec:algebras-from-groupoids}

In many applications of the theory developed in this paper, the $E_k$-algebras in question will arise as the classifying space of a monoidal groupoid; these are always $E_1$-algebras, but if the monoidal groupoid is braided, resp.\ symmetric,\index{monoidal groupoid} the result is an $E_2$-algebra, resp.\ $E_\infty$-algebra. In this section we shall explain the basic application of the theory developed in this paper to such examples. As before, $\sfC = \sfS^\sfG$ with $\sfS$ satisfying the axioms of Sections \ref{sec:axioms-of-cats} and \ref{sec:axioms-of-model-cats}.

\subsection{Constructing $E_k$-algebras}\label{sec:EkAlgFromGpd}
\newglossaryentry{gx}{%
	name={\ensuremath{G_x}},
	description={Automorphism group of object $x \in \sfG$},
	type=symbols
}
Let $(\mathsf{G}, \oplus, \bunit)$ be a $k$-monoidal groupoid, let $r \colon \sf{G} \to \bN$ be a monoidal functor, which we call the \emph{rank}\index{rank}, and suppose that $r^{-1}(0)$ consists precisely of those objects isomorphic to $\bunit_\sfG$. Recall that for an object $x \in \mathsf{G}$ we write $\gls{gx} \coloneqq \mathrm{Aut}_\mathsf{G}(x) = \sfG(x,x)$, and we make the following assumption:

\begin{assum}\label{assum:Unit}
$G_\bunit$ is trivial.
\end{assum}

We have defined the monoidal category $\mathsf{sSet}^{\mathsf{G}}$, and it has a canonical object $\underline{*}$ given by $\underline{*}(x)= *$ for all $x \in \cat{G}$. The object $\underline{*} \in \cat{sSet}^\sfG$ is terminal, so it has the structure of a $\cC$-algebra for any operad $\cC$ in simplicial sets (as its endomorphism operad is the terminal operad). In particular it has the structure of an $E_k^+$-algebra. Under Assumption \ref{assum:Unit} the unit $\bunit_{\cat{sSet}^\sfG}$ is the functor that takes the value $*$ on objects isomorphic to $\bunit_\sfG$ and $\varnothing$ otherwise, so we recognise that $\underline{*} = \bunit_{\cat{sSet}^\sfG} \sqcup \underline{*}_{>0}$ is the unitalisation of the $E_k$-algebra $\underline{*}_{>0}$ having
\[\underline{*}_{>0}(x) \coloneqq \begin{cases}
\varnothing & \text{ if $x\cong \bunit$,}\\
* & \text{ else},
\end{cases}\]
and we may find a cofibrant approximation
\[\gT \overset{\sim}\lra \underline{*}_{>0}\]
as an $E_k$-algebra in $\cat{sSet}^\sfG$. We may then form the left Kan extension
\[\gR \coloneqq r_*(\gT) \in \Alg_{E_k}(\mathsf{sSet}^{\N}).\]

As $\gT$ is a cofibrant $E_k$-algebra it is in particular cofibrant in $\mathsf{sSet}^{\mathsf{G}}$, and so $r_*(\gT)$ is a derived left Kan extension. Thus $U^{E_k}\gR \simeq \bL r_*(\underline{*}_{>0})$ and we have
\[(U^{E_k}\gR)(n) \simeq \begin{cases} \varnothing & \text{if $n=0$,} \\
\bigsqcup_{\substack{[x] \in \pi_0(\sfG) \\ r(x)=n}}BG_x & \text{if $n>0$,}\end{cases}\]
where the coproduct is over isomorphism classes of objects in $\sfG$ of rank $n$.

\subsection{$E_1$-splitting complexes}\label{sec:splittingcomplex}

In Definition \ref{def:derived-indecomposables} we defined the derived $E_1$-indecom\-posables $Q^{E_1}_\bL(\gR)$ of an $E_1$-algebra $\gR$, whose homology is the $E_1$-homology of $\gR$. This heuristically computes the generators, relations, etc., of $\gR$, and is used to bound the number of cells needed for a CW approximation of $\gR$ in the category of $E_1$-algebras. In Section \ref{sec:BarConstr} we proved it may be computed by a bar construction.

For the non-unital $E_k$-algebras $\gR$ arising as in the previous section, we wish to give a combinatorial model for this bar construction. To do so, we will make the following simplifying assumption (later we will explain how this can be omitted, at the expense of complicating the answer a little).

\begin{assum}\label{assum:MonoidalInj}
For all objects $x, y \in \cat{G}$, the homomorphism $- \oplus - \colon G_x \times G_y \to G_{x \oplus y}$ is injective.
\end{assum}

Let $\sfG_{r>0}$ denote the full subgroupoid of $\sfG$ on those objects $x$ with rank $r(x) > 0$, i.e.\ those objects not isomorphic to $\bunit_\sfG$.

\newglossaryentry{suspsplit}{%
	name={\ensuremath{\suspsplit}},
		description={Twice-suspended $E_1$-splitting complex},
		type=symbols
	}
\begin{definition}\label{def:YCx}
For $x \in \mathsf{G}$ let $\gls{suspsplit}_\bullet(x) \in \mathsf{ssSet}_*$ be the semi-simplicial pointed set with $p$-simplices given by
\[\suspsplit_p(x) \coloneqq \left(\colim_{x_1, \ldots, x_{p} \in \sfG_{r>0}^p} \sfG(x_1 \oplus \cdots \oplus x_{p}, x)\right)_+.\]
The face maps $d_0, d_p \colon  \suspsplit_p(x) \to \suspsplit_{p-1}(x)$ are the constant maps to the basepoint. For $0 < j < p$ the face map $d_j$ is induced by replacing $(x_j, x_{j+1})$ by $x_j \oplus x_{j+1}$. 
We write $\suspsplit(x) \coloneqq \fgr{\suspsplit_\bullet(x)} \in \mathsf{sSet}_*$ for its thick geometric realization into pointed simplicial sets.
\end{definition}

In Section \ref{sec:BarConstr} we defined an $E_1$-bar construction $\tilde{B}^{E_1}(\gR)$ for a non-unital $E_1$-algebra $\gR$ in a pointed category such as $\cat{sSet}_\ast$. Let us recall its definition. It is the thick geometric realization of a semi-simplicial space with $p$-simplices $\tilde{B}_p^{E_1}(\gR)$ given by the quotient of $\cP_1(p) \times (\gR^+)^{\wedge p}$ by the subobject consisting entirely of units. Here $\cP_1(p)$ is a contractible space of divisions of $[0,1]$ into $p+2$ intervals, and $\gR^+$ is the unitalization of $\gR$ obtained by formally adding a unit. Let $\gT$ be obtained from $\underline{\ast}_{>0}$ as in Section \ref{sec:EkAlgFromGpd}, with $k=1$.

\begin{proposition}\label{prop:YCxModel}
Under Assumption \ref{assum:MonoidalInj} there are $G_x$-equivariant homotopy equivalences\index{indecomposables!derived}
\[S^1 \wedge Q^{E_1}_\bL(\underline{*}_{>0})(x) \simeq S^1 \wedge Q^{E_1}(\gT)(x) \simeq  \suspsplit(x)\]
of pointed simplicial sets.
\end{proposition}

\begin{proof}
First note that all three terms are equivalent to $*$ if $x \cong \bunit$, so we may suppose that $x \not\cong \bunit$. As $\gT$ is cofibrant in $\sfC$ we may apply Theorem \ref{thm:BarHomologyIndec}, which gives an equivalence $S^1 \wedge Q^{E_1}_\bL(\gT)\simeq \tilde{B}^{E_1}(\gT_+)$ between the suspension of the derived $E_1$-indecomposables of $\gT$ and the $E_1$-bar construction recalled above. Here $\gT_+ \colon \sfG \to \cat{sSet}_*$ is the functor obtained by levelwise adding a disjoint basepoint. 

By definition of the Day convolution product we have
\begin{align*}\cP_1(p)& \times (\gT_+^+)^{\wedge p}(x) \\
&= \colim_{x_1, \ldots, x_p \in \sfG^p} \left(\cP_1(p) \times \sfG(x_1 \oplus \cdots \oplus x_p, x) \times \gT^+(x_1) \times \cdots \times \gT^+(x_p)\right)_+,\end{align*}
and as $\tilde{B}_p^{E_1}(\gT_+)$ is obtained by taking the quotient by the subobject consisting entirely of units, under our assumption that $x \not\cong \bunit_\sfG$ this is also $\tilde{B}_p^{E_1}(\gT_+)(x)$.

There is a natural transformation of functors $\sfG^p \to \cat{sSets}$
\begin{equation*}
\begin{tikzcd}
\cP_1(p)_+ \wedge \sfG(x_1 \oplus \cdots \oplus x_p, x)_+ \wedge \gT^+(x_1)_+ \wedge \cdots \wedge \gT^+(x_p)_+  \dar \\ \sfG(x_1 \oplus \cdots \oplus x_p, x)_+ 
\end{tikzcd}
\end{equation*}
given by $\cP_1(p)_+ \overset{\sim}\to S^0$ and $\gT^+(x_i)_+ \overset{\sim}\to S^0$. As $G_{x_i}$ acts freely on $\gT(x_i)$, the source is a cofibrant functor from $\sfG^p$ to $\cat{sSet}_*$, and under Assumption \ref{assum:MonoidalInj} the group $G_{x_1} \times \cdots \times G_{x_p}$ acts freely on $\sfG(x_1 \oplus \cdots \oplus x_p, x)$ so the target is also a cofibrant functor: this is therefore a weak equivalence between cofibrant functors, so a weak equivalence on colimits. Letting $Z^{E_1}_\bullet(x) \in \mathsf{ssSet}_*$ be the semi-simplicial pointed set with 
\[Z^{E_1}_p(x) \coloneqq \colim_{x_1, \ldots, x_{p} \in \sfG^p} \left( \sfG(x_1 \oplus \cdots \oplus x_{p}, x)\right)_+,\]
and face maps analogous to $\suspsplit_\bullet(x)$, this discussion determines a semi-simplicial map $\tilde{B}_\bullet^{E_1}(\gT_+)(x) \to Z^{E_1}_\bullet(x)$ which is a levelwise weak equivalence, and so a weak equivalence on geometric realization.

Now the semi-simplicial object $Z^{E_1}_\bullet(x)$ admits a system of degeneracies, by inserting copies of $\bunit_\sfG$ into a tuple $(x_1, \ldots, x_p)$, giving it the structure of a simplicial set. With these degeneracies, a face of a non-degenerate simplex is non-degenerate (as we have a monoidal rank functor $r \colon \sf{G} \to \bN$ such that $r(x)=0$ if and only if $x \cong \bunit_\sfG$), so the non-degenerate simplices form a sub-semi-simplicial set, and this is precisely $\suspsplit_\bullet(x)$. Thus the composition
\[\fgr{\suspsplit_\bullet(x)} \lra \fgr{Z^{E_1}_\bullet(x)} \lra \gr{Z^{E_1}_\bullet(x)}\]
is an isomorphism, and the second map is a weak equivalence, so the first map is also a weak equivalence. This gives a zig-zag
\[\suspsplit(x)= \fgr{\suspsplit_\bullet(x)} \overset{\sim}\lra \fgr{Z^{E_1}_\bullet(x)} \overset{\sim}\longleftarrow \fgr{\tilde{B}_\bullet^{E_1}(\gT_+)(x)} = \tilde{B}^{E_1}(\gT_+)(x)\]
of $G_x$-equivariant maps which are each weak equivalences.
\end{proof}

\begin{corollary}\label{cor:DerIndFromYCx}
	There is an equivalence
	\begin{equation*}
	S^1 \wedge Q^{E_1}_\bL(\gR)(n) 
	\simeq \bigvee_{\substack{[x] \in \pi_0(\sfG) \\ r(x)=n}} \suspsplit(x) \wedge_{G_x} (EG_x)_+.
	\end{equation*}
\end{corollary}
\begin{proof}
	We have $\gR = r_*(\gT)$, with $\gT \overset{\sim}\to \underline{*}_{>0}$ a cofibrant approximation as an $E_1$-algebra. Thus $\gR$ is a cofibrant $E_1$-algebra, so
	\[Q^{E_1}_\bL(\gR) \simeq Q^{E_1}(\gR) = Q^{E_1}(r_*(\gT)) \cong r_*(Q^{E_1}(\gT)).\]
	By Proposition \ref{prop:YCxModel} there is a weak equivalence $S^1 \wedge Q^{E_1}(\gT)(x) \simeq \suspsplit(x) \wedge (EG_x)_+$ of cofibrant $G_x$-spaces, so as required
	\[S^1 \wedge r_*(Q^{E_1}(\gT))(n) \simeq \bigvee_{\substack{[x] \in \pi_0(\sfG) \\ r(x)=n}} \suspsplit(x) \wedge_{G_x} (EG_x)_+. \qedhere\]
\end{proof}

By this corollary, the $E_1$-homology of $\gR$ can be interpreted in terms of the $G_x$-equivariant homology of $\suspsplit(x)$. The simplicial set $\suspsplit(x)$ has no homology in degrees above $r(x)$, because $\suspsplit_p(x)=*$ for $p>r(x)$, and the best possible situation is when it only has homology in this degree.

\newglossaryentry{stei}{%
	name={\ensuremath{{St}^{E_1}}},
		description={$E_1$-Steinberg module},
		type=symbols
	}
\begin{definition}\label{defn:StdConnEst}
	If the homology of $\suspsplit(x)$ is concentrated in degree $r(x)$ for every $x \in \mathsf{G}$, then we say that $(\mathsf{G}, \oplus, \bunit)$ \emph{satisfies the standard connectivity estimate}\index{standard connectivity estimate}, and call the $\bZ[G_x]$-module ${St}^{E_1}(x) \coloneqq \tilde{H}_{r(x)}(\suspsplit(x);\bZ)$ the associated \emph{$E_1$-Steinberg module}\index{$E_1$-Steinberg module}.
\end{definition}

In this case we have
\[H^{E_1}_{n,d}(\gR;\bZ) = H_d(Q^{E_1}_\bL(\gR)(n);\bZ) = \bigoplus_{\substack{[x] \in \pi_0(\sfG)\\r(x)=n}}H_{d-(n-1)}(G_x ; {St}^{E_1}(x))\]
so in particular $H^{E_1}_{n,d}(\gR;\bZ)=0$ for $d < n-1$. As $E_k$-homology may be computed as a $k$-fold bar construction, a bar spectral sequence may be used to transfer this vanishing line for $E_1$ to a vanishing line for $E_2$-homology. More precisely, by Theorem \ref{thm:TrfUp} it follows that if $\sfG$ is braided monoidal then $\smash{H^{E_2}_{n,d}}(\gR;\bZ)=0$ for $d < n-1$ too, and if $\mathsf{G}$ is symmetric monoidal then $\smash{H^{E_\infty}_{n,d}}(\gR;\bZ)=0$ for $d < n-1$ as well.

\begin{remark}
It is clear from the proof of Proposition \ref{prop:YCxModel} that Assumption \ref{assum:MonoidalInj} may be omitted if in the definition of $\suspsplit_p(x)$ one forms the homotopy colimit (of simplicial sets) rather than the colimit (of sets). This is analogous to the relaxation of a similar injectivity condition in \cite{RWW} obtained by Krannich \cite[\S 7.3]{Krannich}. In Section \ref{sec:splittingcomplexRevis} we will give a combinatorial model for $Q^{E_1}_\bL(\gT)(x)$ that does not use Assumption \ref{assum:MonoidalInj}.
\end{remark}

\begin{remark}If $\sfG$ is $k$-monoidal then there is an evident $k$-fold semi-simplicial pointed set generalising that of Definition \ref{def:YCx}, for which one can prove the analogue of Proposition \ref{prop:YCxModel} relating it to the derived $E_k$-indecomposables of $\underline{*}_{>0}$. If $\sfG$ is symmetric monoidal, these assemble into an infinite bar spectrum, which in Section \ref{sec:InfiniteBarSpectra} we have shown is equivalent to the suspension spectrum of $Q^{E_\infty}_\bL(\gR)$. In Section \ref{sec:infty-splitting-complex} we will give a combinatorial model for it.
\end{remark}

The semi-simplicial set $\suspsplit_\bullet(x)$ is visibly a double semi-simplicial suspension, where one suspension is formed on the left and one is formed on the right, of the following semi-simplicial set $\split_\bullet(x)$.   We record this in Lemma~\ref{lem:split-suspension} below.

\newglossaryentry{split}{%
	name={\ensuremath{\split}},
		description={$E_1$-splitting complex},
		type=symbols
	}
\begin{definition}\label{def:e1-splitting} 
	For $x \in \sfG$ let $\gls{split}_\bullet(x) \in \mathsf{ssSet}$ be the semi-simplicial set with $p$-simplices given by
	\[\split_p(x) \coloneqq \colim_{x_0, \ldots, x_{p+1} \in \sfG_{r>0}^{p+2}} \sfG(x_0\oplus \cdots \oplus x_{p+1}, x)\]
	and all face maps are given by the monoidal structure. This is called the \emph{$E_1$-splitting complex}\index{$E_1$-splitting complex}. We write $\split(x) \coloneqq \fgr{\split_\bullet(x)} \in \mathsf{sSet}$ for its thick geometric realization into simplicial sets.
\end{definition}

\begin{lemma}\label{lem:split-suspension} There is a $G_x$-equivariant homotopy equivalence $\suspsplit(x) \simeq \Sigma^2 \split(x)$ of pointed simplicial sets.\end{lemma}

We note that $\Sigma$ denotes the unreduced suspension, and in particular $\Sigma \varnothing = S^0$.  In the following remark we give a more concrete description of $\split(x)$ under mild conditions. For another perspective on it and its relationship to the derived decomposables, see Section \ref{sec:splittingcomplexRevis}.

\begin{remark}\label{rem:young-type}
We may give a more concrete description of $\split_\bullet(x)$ in terms of subgroups of $G_x = \sfG(x,x)$ which we call \emph{Young-type subgroups}.  An ordered tuple of objects $(x_0,\ldots,x_{p+1}) \in \sfG_{r>0}^{p+2}$ together with an isomorphism $\iota\colon x_0 \oplus \cdots \oplus x_{p+1} \to x$ defines an element of $S^{E_1}_p(x)$ which we also denote $\iota$.  Acting on $\iota$ defines a $G_x$-equivariant injection
  \begin{equation}\label{eq:1}
    G_x/G_{(x_0, \dots, x_{p+1})} \hookrightarrow S^{E_1}_p(x),
  \end{equation}
  where $G_{(x_0, \dots, x_{p+1})} < G_x$ denotes the image of the homomorphism
\[G_{x_0} \times \cdots \times G_{x_{p+1}} \lra G_x\]
induced by $\iota$, which is injective by Assumption \ref{assum:MonoidalInj}.  This image $G_{(x_0,\ldots,x_{p+1})}$ is the \emph{Young-type subgroup}\index{Young-type subgroup}.  Different choices of $\iota$ lead to conjugate subgroups.

The set $S^{E_1}_p(x)$ is the disjoint union of the images of the injections~\eqref{eq:1}, over  tuples $(x_0, \dots, x_{p+1}) \in \sfG^{p+2}_{r > 0}$, one in each isomorphism class, and isomorphisms $\iota \colon x_0 \oplus \dots \oplus x_{p+1} \to x$, one in each $G_x$-orbit.

We shall say that $(x'_0,\ldots,x'_{p'+1})$ is a \emph{refinement} of $(x_0,\ldots,x_{p+1})$ if there exists a surjective order-preserving map $\phi \colon [p'+1] \to [p+1]$ such that there exist isomorphisms $\bigoplus_{j \in \phi^{-1}(i)} x'_j \cong x_i$.  Choosing such isomorphisms for each $i \in [p+1]$ leads to a diagram of isomorphisms
\begin{equation}\label{eq:4}
  \begin{tikzcd}
    x'_0 \oplus \dots \oplus x'_{p'+1} \rar["\iota'"] \dar & x \arrow[d,equal]\\
    x_0 \oplus \dots \oplus x_{p+1} \rar["\iota"]  & x,
  \end{tikzcd}
\end{equation}
where $\iota'$ is defined by commutativity of the diagram.  In this situation, the Young-type subgroup $G_{(x'_0, \dots, x'_{p'+1})} < G_x$ is a subgroup of $G_{(x_0, \dots, x_{p+1})}$, and we obtain a $G_x$-equivariant diagram
\begin{equation}\label{eq:7}
  \begin{tikzcd}
    G_x/G_{(x'_0, \dots, x'_{p'+1})} \rar[hook] \dar[twoheadrightarrow] & S^{E_1}_{p'}(x)\dar["\theta"]\\
    G_x/G_{(x_0, \dots, x_{p+1})} \rar[hook]& S^{E_1}_p(x),
  \end{tikzcd}
\end{equation}
when the surjection $\phi\colon [p'+1] \to [p+1]$ is induced from $\theta\colon [p] \hookrightarrow[p']$.  If we instead choose $\iota'$ arbitrarily, the diagram~\eqref{eq:4} commutes only up to an element $g \in G_x$, in which case $G_{(x'_0, \dots, x'_{p'+1})}$ will only be conjugate to a subgroup of $G_{(x_0, \dots, x_{p+1})}$.  In this situation we still obtain a $G_x$-equivariant diagram of the form~\eqref{eq:7}, but the left horizontal map is induced by the element $g \in G_x$.

Let us also remark that when the monoidal structure on $\sfG$ is strict, and the underlying category is skeletal, the isomorphisms $\iota$, $\iota'$, and $\phi_i$ in the above discussion can all be chosen as the identity.
\end{remark}

\begin{example}\label{exam:n-e1}
Let us consider the case $\sfG = \bN$ with monoidal product given by addition, with monoidal rank functor $r \colon \bN \to \bN$ given by the identity. The resulting $E_1$-algebra $\gR$ is weakly equivalent to the non-unital associative algebra $\bN_{>0}$ with multiplication given by addition. Using Remark \ref{rem:young-type}, its $E_1$-splitting complex $S^{E_1}(n)$ is the thick geometric realization of the semi-simplicial set with $p$-simplices given by the \emph{ordered} sum decompositions $n_0+\cdots+n_{p+1} = n$ with each $n_i>0$.
	
If $n=1$ this is empty. If $n \geq 2$, it is contractible. To see this, note that $S^{E_1}(n)$ is homeomorphic to the geometric realization of the nerve of the poset with objects the \emph{ordered} sum decompositions $n_0+\cdots+n_{p+1} = n$ ordered by refinement. If $n \geq 2$, this has a terminal object $1+\cdots+1 = n$.

We conclude that $T^{E_1}(n) = \Sigma^2 S^{E_1}(n) \simeq S^1$ if $n=1$ and $\ast$ otherwise. Using Proposition \ref{prop:YCxModel} we see that $H^{E_1}_{n,d}(\bN_{>0};\bZ)$ is $\bZ$ if $(n,d) = (1,0)$ and $0$ otherwise. This also follows from the fact that $\bN_{>0}$ is weakly equivalent to the free $E_1$-algebra on $1_* \ast \in \cat{sSet}^\bN$.
\end{example}

\subsection{$E_k$-splitting complexes}\label{sec:ek-splitting-complex} \index{$E_k$-splitting complex}Here we give the generalization of Section \ref{sec:splittingcomplex} to $k \geq 2$. It is more convenient to generalize the simplicial set $Z^{E_1}_\bullet(x)$ rather than the semi-simplicial set $T^{E_1}_\bullet(x)$. Considering the finite set $\underline{p} = \{1,\ldots,p\}$ as a discrete category, we define $\cat{G}^{p_1 \cdots p_k}$ to be the category of functors $\cat{Fun}(\underline{p_1} \times \cdots \times \underline{p_k},\cat{G})$. We denote an object of this category as $\vec{x}$, which consists of a collection of objects $x_{i_1,\ldots,i_k}$ for $1 \leq i_j \leq p_j$.

\begin{definition}For $x \in \cat{G}$ define a $k$-fold simplicial pointed set  $Z^{E_k}_{\bullet,\ldots,\bullet}(x)$ with the pointed set of $(p_1,\ldots,p_k)$-simplices given by 
	\[Z^{E_k}_{p_1,\ldots,p_k}(x) \coloneqq \left(\underset{\vec{x} \in \cat{G}^{p_1\cdots p_k}}{\mr{colim}} \cat{G}(x_{1,\ldots,1} \oplus \cdots \oplus x_{p_1,\ldots,p_k},x) \right)_+.\]
	
	For $1 \leq i \leq k$, the face maps $d^i_0,d^i_{p_i}$ are the constant maps to the basepoint. The face maps $d^i_j$ for $0<j<p_i$ are induced by replacing each pair of objects $(x_{a_1,\ldots,a_{i-1},j,a_j,\ldots,a_p},x_{a_1,\ldots,a_{i-1},j+1,a_j,\ldots,a_p})$ by its sum. The degeneracy maps insert $\bunit$'s. It has a remaining $G_x$-action, and we write $Z^{E_k}(x) \coloneqq |Z_{\bullet,\ldots,\bullet}(x)| \in \cat{sSet}_\ast^{G_x}$ for its $k$-fold thin geometric realization.\end{definition}

Let $\gT$ be obtained from $\underline{\ast}_{>0}$ as in Section \ref{sec:EkAlgFromGpd}.  The analogue of Proposition \ref{prop:YCxModel} is a computation of the derived indecomposables $Q^{E_k}_\bL(\underline{\ast}_{>0}) \in \cat{sSet}_\ast^\cat{G}$ in terms of $Z^{E_k}(x)$ under Assumption \ref{assum:MonoidalInj}. Following the proof of Proposition \ref{prop:YCxModel} but using the $E_k$-bar construction in place of the $E_1$-bar construction, we get a $G_x$-equivariant homotopy equivalence $\tilde{B}^{E_k}(\gT)(x) \to Z^{E_k}(x)$ and conclude that:

\begin{proposition}Under Assumption \ref{assum:MonoidalInj}, there are $G_x$-equivariant homotopy equivalences
\[S^k \wedge Q^{E_k}_\bL(\underline{\ast}_{>0})(x) \simeq S^k \wedge Q^{E_k}(\gT)(x) \simeq Z^{E_k}(x).\]
pointed simplicial sets for $x \in \cat{G}_{>0}$.
\end{proposition}

Taking the derived pushforward along $r \colon \cat{G} \to \bN$, a straightforward adaptation of Corollary \ref{cor:DerIndFromYCx} gives us:

\begin{corollary}For $n \geq 1$, there is a weak equivalence
	\[S^k \wedge Q^{E_k}_\bL(\bL r_* \underline{\ast}_{>0})(n) \simeq \bigvee_{\substack{[x] \in \pi_0(\cat{G}) \\ r(x)=n}} Z^{E_k}(x) \hcoker G_x.\]
\end{corollary}

\begin{remark}
	As in Remark \ref{rem:young-type} we can identify the simplices of $Z^{E_k}_{\bullet,\ldots,\bullet}(x)$ with a disjoint union of cosets of Young-type subgroups. 
	
	There is such a Young-type subgroup for each isomorphism $\iota \colon x_1 \oplus \cdots \oplus x_n \overset{\sim}\lra x$ of objects $\cat{G}$; it is the subgroup $G_{(x_1,\cdots,x_n)} \leq G_x$ given by the image of the injective homomorphism $G_{x_1} \times \cdots \times G_{x_n} \hookrightarrow G_x$ induced by $\iota$. Under the above conditions, we can identify the sets which appear: 
	\[\underset{\vec{x} \in \cat{G}^{p_1\cdots p_k}}{\mr{colim}} \cat{G}(x_{1,\ldots,1} \oplus \cdots \oplus x_{p_1,\ldots,p_k},x)  \cong \bigsqcup_{\substack{\vec{x}\in \cat{G}^{p_1 \cdots p_k} \\
			\iota: \bigoplus x_{i_1,\ldots,i_k} \to x}} \frac{G_x}{G_{(x_{1,\ldots,1},\cdots,x_{p_1,\cdots,p_k})}},\]
	the indexing set running over tuples $\vec{x} \in \in \cat{G}^{p_1 \cdots p_k}$, one in each isomorphism class, and isomorphisms $\iota \colon\bigoplus x_{i_1,\ldots,i_k} \overset{\sim}\to x$, one in each $G_x$-orbit.
\end{remark}

\subsection{$E_\infty$-splitting complexes}\label{sec:infty-splitting-complex}

If $(\mathsf{G}, \oplus, \bunit)$ is a symmetric monoidal groupoid then the construction of Section \ref{sec:EkAlgFromGpd} provides an $E_\infty$-algebra $\gR$, which has derived $E_\infty$-indecomposables $Q^{E_\infty}_\bL(\gR)$ as in Definition \ref{def:derived-indecomposables}. Here we wish to give a combinatorial model for this, in the same way that we gave a combinatorial model for the derived $E_1$-indecomposables given in the previous section. However, the strategy in this case will be different from that of Section \ref{sec:splittingcomplex}; in Section \ref{sec:splittingcomplexRevis} we will explain how the ideas of this section can also be applied in the $E_1$ case.

\begin{definition}\label{def:e_infty-splitt-cat}
  For $x \in \cat{G}$, define the \emph{$E_\infty$-splitting category} $\sfS^{E_\infty}(x)$ as follows:\index{$E_\infty$-splitting category}
	\begin{itemize}\item Its objects are given by triples $([n], f, \phi)$ of a finite set $[n]\coloneqq \{1,2,\ldots, n\}$ with $n\geq 2$, a function $f \colon [n] \to \mathrm{ob}(\sfG_{r>0})$, and a morphism $\phi \colon \bigoplus_{\alpha \in [n]} f(\alpha) \to x$ in $\sfG$ (necessarily an isomorphism).
	\item A morphism $([n], f, \phi) \to ([n'], f', \phi')$ is the data of a surjection $e\colon [n] \to [n']$ and isomorphisms $\varphi_\alpha\colon f'(\alpha) \to \bigoplus_{\beta \in e^{-1}(\alpha)} f(\beta)$ for each $\alpha \in [n']$ such that the following diagram commutes
\begin{equation*}
\begin{tikzcd}
\bigoplus_{\alpha \in [n']} f'(\alpha) \dar[swap]{\bigoplus_{\alpha \in [n']} \varphi_\alpha}\rar{\phi'} & x\\
\bigoplus_{\alpha \in [n']} \bigoplus_{\beta \in e^{-1}(\alpha)} f(\beta) \ar{r}{\cong}& \bigoplus_{\beta \in [n]} f(\beta), \uar[swap]{\phi}
\end{tikzcd}
\end{equation*}
with bottom map the canonical identification of these two sums.
\item Composition is given by composing the surjections $e$, and composing the appropriate direct sums of the isomorphisms $\phi_\alpha$.
\end{itemize}
\end{definition}

The $E_\infty$-splitting category is similar in nature to the category of simplices of the semi-simplicial set $\split_\bullet (x)$ defined in Definition~\ref{def:e1-splitting} in the $E_1$ setting, but also incorporates morphisms that permute direct summands.  We discuss this connection in more depth in Section~\ref{sec:splittingcomplexRevis} below.

\newglossaryentry{splitinf}{%
	name={\ensuremath{\splitinf}},
		description={$E_\infty$-splitting complex},
		type=symbols
	}
\begin{definition}For $x \in \mathsf{G}$, we let the \emph{$E_\infty$-splitting complex}\index{$E_\infty$-splitting complex} $\gls{splitinf}(x) \in \mathsf{sSet}$ be the nerve of the \emph{$E_\infty$-splitting category} $\sfS^{E_\infty}(x)$.\end{definition}

As we shall see in Proposition~\ref{prop:DecEqualsSplittingCx}, the $E_\infty$-splitting complex gives a model for the derived $E_\infty$-decomposables of $\gT$ and hence $\underline{*}_{>0}$.

Let $\sfG^n \wr \fS_n$ denote the Grothendieck construction of the action of $\fS_n$ on the groupoid $\sfG^n$: it is a groupoid whose objects are given by tuples $(x_1, \ldots, x_n)$ of objects of $\sfG$, and a morphism from $(x_1, \ldots, x_n)$ to $(x'_1, \ldots, x'_n)$ is given by a permutation $\sigma \in \fS_n$ and a collection of morphisms $\varphi_i \colon x_{\sigma(i)} \to x'_{i}$ in $\sfG$ for $i=1,2,\ldots,n$. As $\oplus$ is symmetric monoidal, the functor $\oplus \colon \sfG^n \to \sfG$ extends to a functor $\pi_n \colon \sfG^n \wr \fS_n \to \sfG$.

If $X \in \cat{sSet}^\sfG$, then there is a functor $X^{\otimes n} \wr \fS_n \colon \sfG^n \wr \fS_n \to \cat{sSet}$ given on objects by
\[(X^{\otimes n} \wr \fS_n)(x_1, \ldots, x_n) = X(x_1) \otimes \cdots \otimes X(x_n)\]
and on a morphism $(\sigma, \{\varphi_i\})$ by
\[X(x_1) \otimes \cdots \otimes X(x_n) \overset{\hat{\sigma}}\lra  X(x_{\sigma(1)}) \otimes \cdots \otimes X(x_{\sigma(n)}) \xrightarrow{\prod_i \varphi_i}  X(x'_1) \otimes \cdots \otimes X(x'_n)\]
where $\hat{\sigma}$ denotes the permutation of the factors given by $\sigma$.

\begin{lemma}\label{lem:GpdModelForEinfty}
If $X \in \cat{sSet}^\sfG$ is cofibrant then there is a weak equivalence in $\cat{sSet}^\sfG$
\[\cC_\infty(n) \times_{\fS_n} X^{\otimes n} \simeq \bL(\pi_n)_*(X^{\otimes n} \wr \fS_n).\]
\end{lemma}
\begin{proof}
Consider the functor $C \colon \sfG^n \wr \fS_n \to \cat{sSet}$ given on objects by
\[C(x_1, \ldots, x_n) \coloneqq \cC_\infty(n) \times X(x_1) \otimes \cdots \otimes X(x_n)\]
and on a morphism $(\sigma, \{\varphi_i\})$ by
\begin{align*}
&\cC_\infty(n) \times X(x_1) \otimes \cdots \otimes X(x_n) \xrightarrow{\sigma \times \hat{\sigma}} \cC_\infty(n) \times X(x_{\sigma(1)}) \otimes \cdots \otimes X(x_{\sigma(n)})\\ &\qquad\qquad \xrightarrow{\mr{id} \times \prod_i \varphi_i} \cC_\infty(n) \times X(x'_1) \otimes \cdots \otimes X(x'_n).\end{align*}
The projection $\cC_\infty(n) \to *$ gives a natural transformation $C \to X^{\otimes n} \wr \fS_n$, which is an objectwise weak equivalence as $\cC_\infty(n) \simeq *$. As $X$ is cofibrant, each $X(x)$ has a free $G_{x}$-action, and certainly $\cC_\infty(n)$ has a free $\fS_n$-action. Thus $C(x_1, \ldots, x_n)$ has a free $\mr{Aut}_{X^n \wr \fS_n}(x_1, \ldots, x_n)$-action, and hence $C$ is cofibrant. Hence 
\[(\pi_n)_*(C) = \cC_\infty(n) \times_{\fS_n} X^{\otimes n}\]
computes the homotopy Kan extension of $X^{\otimes n} \wr \fS_n$, as required.
\end{proof}

This lemma shall be used to produce a model for the $E_\infty$-decomposables of $\gT$.

\begin{definition}
For $x \in \cat{G}$, let $\cat{GS}^{E_\infty}(x)$ be the category object internal to groupoids described as follows. 
\begin{itemize}
	\item The ``object'' groupoid $\cat{O}=\mr{ob}(\cat{GS}^{E_\infty}(x))$ has objects triples $([n], f, \phi)$ of a finite set $[n]\coloneqq \{1,2,\ldots, n\}$ with $n\geq 2$, a function $f \colon [n] \to \mr{ob}(\sfG_{r>0})$, and a morphism $\phi \colon \bigoplus_{\alpha \in [n]} f(\alpha) \to x$ in $\sfG$. A morphism from $([n], f, \phi)$ to $([n'], f', \phi')$ is the data of a bijection $b \colon [n] \to [n']$ and a collection of morphisms $\epsilon(\alpha) \colon f(\alpha) \to f'(b(\alpha))$ such that 
	\begin{equation*}
	\begin{tikzcd}
	\bigoplus_{\alpha \in [n]} f(\alpha) \dar[swap]{\bigoplus_{\alpha \in [n]} \epsilon(\alpha)}\rar{\phi} & x\\
	\bigoplus_{\alpha \in [n]} f'(b(\alpha)) \rar{b_*}& \bigoplus_{\beta \in [n']} f'(\beta) \uar[swap]{\phi'}
	\end{tikzcd}
	\end{equation*}
	commutes.
	\item The ``morphism'' groupoid $\cat{M}=\mr{mor}(\cat{GS}^{E_\infty}(x))$ has objects given by tuples $([n_0], f_0, \phi_0; [n_1], f_1, \phi_1; e, \{\varphi_\alpha\})$ of a pair of objects of the groupoid $\mathrm{ob}(\bar{\sfS}^{E_\infty}(x))$, a surjection $e \colon [n_0] \twoheadrightarrow [n_1]$, and isomorphisms $\varphi_\alpha \colon f_1(\alpha) \to \bigoplus_{\beta \in e^{-1}(\alpha)} f_0(\beta)$ for $\alpha \in [n_1]$. A morphism from such an object to another $([n_0'], f_0', \phi_0'; [n_1'], f_1', \phi_1'; e', \{\varphi_\alpha'\})$ is the data of morphisms in the groupoid $\cat{O}$
\[(b_i, \epsilon(\alpha)_i) \colon ([n_i], f_i, \phi_i) \lra ([n_i'], f_i', \phi_i')\]
	such that $b_1 \circ e = e' \circ b_0$ and
	\begin{equation*}
	\begin{tikzcd}
	f_1(\alpha) \rar{\varphi_\alpha} \dar[swap]{\epsilon(\alpha)_1} & \bigoplus_{\beta \in e^{-1}(\alpha)} f_0(\beta) \dar{\bigoplus \epsilon(\beta)_0}\\
	f_1'(b_1(\alpha)) \rar{\varphi'_{b_1(\alpha)}} & \bigoplus_{b_0(\beta) \in (e')^{-1}(b_1(\alpha))} f_0'(b_0(\beta))
	\end{tikzcd}
	\end{equation*}
	commutes for all $\alpha \in [n_1]$.
      \item The composition functor $\cat{M} \times_{\cat{O}} \cat{M}$ is defined similarly to Definition~\ref{def:e_infty-splitt-cat}: compose the surjections $e$ and compose the appropriate direct sums of the isomorphisms $\phi_\alpha$.
\end{itemize}
\end{definition}

Let us write $N\cat{GS}^{E_\infty}(x)$ for the simplicial set obtained as the nerve of this category, i.e.\ first form the simplicial category with objects $N_\bullet \cat{O}$ and morphisms $N_\bullet \cat{M}$, take the bisimplicial nerve of this simplicial category, then form the diagonal simplicial set.

\begin{lemma}\label{lem:DecEqualsSplittingGpd}
There are $G_x$-equivariant homotopy equivalences
\[\Dec_\bL^{E_\infty}(\gT)(x) \simeq N\cat{GS}^{E_\infty}(x)_+\]
of pointed simplicial sets.
\end{lemma}
\begin{proof}
We use the description $\Dec_\bL^{E_\infty}(\gT) \simeq B((E_\infty^{\geq 2})_+, E_\infty, \gT)$ given in Section \ref{sec:deriv-decomp}, where we have
\[B_p((E_\infty^{\geq 2})_+, E_\infty, \gT) = (E_\infty^{\geq 2})_+((E_\infty)^p (\gT)).\]
By Lemma \ref{lem:GpdModelForEinfty} we have,
\[E_\infty^{\geq r}(X) = \coprod_{n \geq r} \cC_\infty(n) \times_{\fS_n} X^{\otimes n} \overset{\simeq}\lra \bL\left(\coprod_{n \geq r} \pi_n\right)_*\left(\coprod_{n \geq r} X^{\otimes n} \wr \fS_n\right).\]
By iterating, this lets us describe $(E_\infty^{\geq 2})_+((E_\infty)^p (X))$ as the homotopy Kan extension along a certain functor $r_p \colon \cat{B}_p \to \sfG$ of a certain functor $X_p \colon \cat{B}_p \to \cat{sSet}$. 

Unwinding definitions, the category $\cat{B}_p$ is given as follows:
\begin{itemize}\item An object of  $\cat{B}_p$ consists of the following data: a collection of finite sets and surjections
	\[[n_0] \overset{e(1)}\longtwoheadleftarrow [n_{1}] \overset{e(2)}\longtwoheadleftarrow \cdots \overset{e(p)} \longtwoheadleftarrow [n_p]\]
	 with $n_0 \geq 2$, together with functions $f_j \colon [n_j] \to \mr{ob}(\sfG)$ for $0 \leq j \leq p$ and morphisms $\phi_{j-1}(i) \colon \bigoplus_{k \in e(j)^{-1}(i)} f_j(k) \to f_{j-1}(i)$ in $\sfG$ for $1 \leq j \leq p$. 
	\item A morphism from such an object to another, decorated with primes, is a collection of bijections $b_j \colon [n_j] \to [n'_j]$ and a collection of morphisms $\epsilon_j(i) \colon f_j(i) \to f'_j(b_j(i))$ in $\sfG$ for $0 \leq j \leq p$ and $i \in [n_j]$ which intertwine the $e$'s and $\phi$'s in the evident way.\end{itemize}

The functor $r_p \colon \cat{B}_p \to \sfG$ sends an object as described above to $\bigoplus_{i \in [n_0]} f_i(0)$ and a morphism as described above to
\[\bigoplus_{i \in [n_0]} f_0(i) \xrightarrow{\bigoplus \epsilon_0(i)} \bigoplus_{i \in [n_0]} f'_0(b_0(i)) \xrightarrow{(b_0)_*} \bigoplus_{k \in [n'_0]} f'_0(k).\]

The functor $X_p \colon \cat{B}_p \to \cat{sSet}$ can be expressed in terms of $X$ as follows. It sends an object as described above to $X(f_p(1)) \otimes \cdots \otimes X(f_p(n_p))$, and induces the evident operation on morphisms.

In the case $X=\gT \simeq \underline{*}_{> 0}$ we therefore obtain
\[B_p((E_\infty^{\geq 2})_+, E_\infty, \gT) \simeq \bL (r_p)_*(\gT_p) \simeq \bL (r_p)_*((\underline{*}_{> 0})_p).\]
The value of this object at $x$ may be described as the classifying space of the subgroupoid of $r_p/x$ in which all objects $f_i(j)$ which arise are required to lie in $\sfG_{r>0}$. This subgroupoid is recognisable as the groupoid $N_p(\cat{GS}^{E_\infty}(x))$ obtained as the $p$th stage of the nerve of the category object in groupoids $\cat{GS}^{E_\infty}(x)$. Let us write $\cat{sS}^{E_\infty}(x)$ for the category object in simplicial sets obtained by taking the nerve in the groupoid direction. The naturality with respect to face maps of the above discussion is easily seen, and we obtain a zig-zag of $G_x$-equivariant maps of simplicial objects
\[B_\bullet((E_\infty^{\geq 2})_+, E_\infty, \gT)(x) \lra \cdots \longleftarrow N_\bullet(\cat{sS}^{E_\infty}(x))\]
which are levelwise weak equivalences. Therefore, taking geometric realization gives the desired conclusion.
\end{proof}

We now wish to simplify our model for $\Dec_\bL^{E_\infty}(\gT)(x)$ for the classifying space of the category object in groupoids $\cat{GS}^{E_\infty}(x)$  and relate it to the classifying space of the category $\cat{S}^{E_\infty}(x)$, which is our definition of $\splitinf(x)$. We shall see in the proof that this amounts to the claim that we may discretize the objects and morphisms of $\cat{sS}^{E_\infty}(x)$ without affecting the homotopy type of its nerve.

\begin{proposition}\label{prop:DecEqualsSplittingCx}
There are $G_x$-equivariant homotopy equivalences
\[\Dec_\bL^{E_\infty}(\gT)(x) \simeq \splitinf(x)_+\]
of pointed simplicial sets.
\end{proposition}

From this proposition we may immediately deduce the following corollary, using the cofibre sequence
\[\Dec_\bL^{E_\infty}(\gT) \lra \fgr{B_\bullet(E_\infty,E_\infty,\gT)}_+ \lra Q^{E_\infty}_\bL(\gT)\]
given in \eqref{eq:DecCofSeq}, and the equivalence $\fgr{B_\bullet(E_\infty,E_\infty,\gT)} \simeq \gT \simeq \underline{*}_{>0}$.

\begin{corollary}\label{cor:IndecEqualsSuspSplittingCx}
There are $G_x$-equivariant homotopy equivalences
\[Q^{E_\infty}_\bL(\gT)(x) \simeq \Sigma \splitinf(x)\]
of pointed simplicial sets.
\end{corollary}

\begin{proof}[Proof of Proposition \ref{prop:DecEqualsSplittingCx}]
In view of Lemma \ref{lem:DecEqualsSplittingGpd}, we must show that there are $G_x$-equivariant homotopy equivalences $N\cat{GS}^{E_\infty}(x) \simeq N{\sfS}^{E_\infty}(x) = \splitinf(x)$ between the nerves of these two categories.

The category ${\sfS}^{E_\infty}(x)$ is obtained from $\cat{GS}^{E_\infty}(x)$ by taking the underlying sets of objects of $\cat{O}$ and $\cat{M}$. Let $\cat{TS}^{E_\infty}(x)$ be the category object internal to topological spaces obtained by forming the (thin) geometric realization of the groupoids $\cat{O}$ and $\cat{M}$ of objects and morphisms, denoted $O=|N_\bullet \cat{O}|$ and $M=|N_\bullet \cat{M}|$. Write $\delta O \subset O$ for the (discrete) subspace of 0-simplices, i.e.\ objects of $\cat{O}$, and similarly $\delta M \subset M$.

The combined source-target map
\[s \times t \colon \cat{M} \lra \cat{O} \times \cat{O}\]
is easily seen to be a fibration of groupoids, i.e.\ to induce a Kan fibration on nerves, so, as the geometric realization of a Kan fibration is a Serre fibration \cite{QuillenKanFib}, the combined source-target map $s \times t \colon M \to O \times O$ is a Serre fibration. This puts us in a position to apply \cite[Theorem 5.2]{EbertRWSx}, applied to the 0-connected map $i: \delta O \to O$. This produces a new topological category $\cat{TS}^{E_\infty}(x)^{\delta O}$ with space of objects $\delta O$ and space of morphisms given by the pullback
\begin{equation*}
\begin{tikzcd}
\mathrm{mor}(\cat{TS}^{E_\infty}(x)^{\delta O}) \rar \dar& M \dar{s \times t}\\
\delta O \times \delta O \rar{i \times i} & O \times O,
\end{tikzcd}
\end{equation*}
and shows that the inclusion $\cat{TS}^{E_\infty}(x)^{\delta O} \to \cat{TS}^{E_\infty}(x)$ induces an equivalence on classifying spaces. We may compute the above pullback as follows: it is the geometric realization of the nerve of the groupoid having the same objects as $\cat{M}$, but only those morphisms which map to identity morphisms under $s \times t$. By the definition of morphisms in $\cat{M}$ this is the discrete groupoid with the same objects as $\cat{M}$, i.e.\ $\mathrm{mor}(\cat{TS}^{E_\infty}(x)^{\delta O}) = \delta M$. Thus $\cat{TS}^{E_\infty}(x)^{\delta O}$ is just the category $\cat{S}^{E_\infty}(x)$, and the conclusion is that the inclusion, which is $G_x$-equivariant, induces a weak equivalence of topological spaces
\[|N_\bullet\cat{S}^{E_\infty}(x)| \overset{\sim}\lra |N_\bullet \cat{TS}^{E_\infty}(x)|.\]
Recognising this as the geometric realization of the map of simplicial sets $\splitinf(x) = N\cat{S}^{E_\infty}(x) \to N\cat{GS}^{E_\infty}(x)$, the desired conclusion follows.
\end{proof}

\begin{example}\label{exam:n-einfty} We continue Example \ref{exam:n-e1}. There we used that addition endows the discrete groupoid $\bN$ with a monoidal structure, which is of course a symmetric monoidal structure (with identity symmetries). Thus the $E_1$-algebra structure on $\gR$ extends to an $E_\infty$-algebra structure, weakly equivalent to the obvious non-unital commutative algebra structure on $\bN_{>0}$. Its $E_\infty$-homology is significantly more complicated than its $E_1$-homology, being related to the associated graded of the symmetric power filtration on $H\bZ$, as explained in the work of Arone--Lesh \cite{AroneLesh, AroneLesh2}.

In this case for $n \in \bN$, the $E_\infty$-splitting category $\cat{S}^{E_\infty}(n)$ may be described as follows. Objects are (ordered) tuples $(n_1, n_2, \ldots, n_k)$ of $k \geq 2$ natural numbers such that $n=\sum n_i$, and a morphism from such a tuple to $(n'_1, n'_2, \ldots, n'_{k'})$ is the data of a surjection $e \colon [k] \twoheadrightarrow [k']$ such that $n'_i = \sum_{j \in e^{-1}(i)} n_j$.

In this example we will explain the relationship between $\cat{S}^{E_\infty}(n)$ and the partition complex of $n$ in the sense of Arone--Dwyer \cite{AroneDwyer}.  Consider the analogous category $\cat{D}([n])$ whose objects are tuples $(S_1, S_2, \ldots, S_k)$ of $k \geq 2$ subsets of $[n]=\{1,2,\ldots,n\}$ such that $[n] = \sqcup_{i=1}^k S_i$, and a morphism from such a tuple to $(S'_1, S'_2, \ldots, S'_{k'})$ is the data of a surjection $e \colon [k] \to [k']$ such that $S'_i = \sqcup_{j \in e^{-1}(i)} S_j$. This category has a natural action of $\fS_n$, by permuting elements of $[n]$ and hence also its subsets. 

Sending a set to its cardinality defines a functor $\phi \colon \cat{D}([n]) \to \cat{S}^{E_\infty}(n)$, which strictly commutes with the $\fS_n$-action, and there is therefore a factorisation
\[B\phi \colon N\cat{D}([n]) \lra N\cat{D}([n])/\fS_n \overset{\varphi}\lra N\cat{S}^{E_\infty}(n)\]
of the induced map on simplicial nerves. It is easy to see that $\varphi$ is an isomorphism.

To describe the middle term we study the $\fS_n$-equivariant homotopy type of $N\cat{D}([n])$. First note that the object $(\{1\}, \{2\}, \ldots, \{n\})$ is initial in $\cat{D}([n])$, so $N\cat{D}([n]) \simeq *$ for $n \geq 2$ (and it is empty for $n=1$). A $p$-simplex is given by a sequence of surjections 
\[[k_0] \overset{e(1)}\longtwoheadleftarrow [k_1] \overset{e(2)}\longtwoheadleftarrow \cdots \overset{e(p)}\longtwoheadleftarrow [k_p]\] with $k_0 \geq 2$ and a tuple of sets $(S_1, S_2, \ldots, S_{k_p})$ decomposing $[n]$. A permutation of $[n]$ stabilises this simplex if and only if it preserves the sets $S_i$, so the stabiliser of this simplex is the Young subgroup $\fS_{S_1} \times \cdots \times \fS_{S_{k_p}} \leq \fS_n$. On the other hand for any decomposition $[n] = \sqcup_{i=1}^r X_i$ defining a Young subgroup $\fS_{X_1} \times \cdots \times \fS_{X_r} \leq \fS_n$, the fixed points for this subgroup consist of the simplices as above such that each $S_i$ is a union of $X_j$'s. We can identify this with $N_\bullet \cat{D}(\{X_1, X_2, \ldots, X_r\})$, the construction above applied to the set of parts $X_j$, so it is contractible as long as $r \geq 2$, and empty if $r=1$. We recognise the above properties as characterising the $\fS_n$-equivariant homotopy type of $E\mathcal{Y}$, the universal $\fS_n$-space whose isotropy is in the collection $\mathcal{Y}$ of all Young subgroups associated to \emph{proper} decompositions of $[n]$. Thus we have $N\cat{D}([n])/\fS_n = B\mathcal{Y}$.

Finally, by \cite[Proposition 7.3]{AroneDwyer} (put into our notation) there is a homotopy equivalence
\[S^1 \wedge \Sigma B\mathcal{Y} \simeq S^1 \wedge (S^n \wedge \Sigma P_n) \hcoker \fS_n,\]
where the homotopy orbits are formed in pointed simplicial sets. Here $P_n$ is the $n$th \emph{partition complex}, i.e.\ the nerve of the poset of partitions of $[n]$, with the discrete and indiscrete partitions removed. It is known that $\Sigma P_n$ is a wedge of $(n-2)$-spheres, and following \cite{AroneDwyer} we write $\mr{Lie}_n^* \coloneqq \mr{sign} \otimes \tilde{H}_{n-2}(\Sigma P_n;\bZ)$.

Putting the above together, and using Corollary \ref{cor:IndecEqualsSuspSplittingCx}, we have
\[S^1 \wedge Q^{E_\infty}_\bL(\bN_{>0})(n) \simeq S^1 \wedge\Sigma S^{E_\infty}(n) \simeq S^1 \wedge\Sigma N\cat{S}^{E_\infty}(n) \simeq S^1 \wedge (S^n \wedge \Sigma P_n) \hcoker \fS_n,\]
and so
\[H_{n,d}^{E_\infty}(\bN_{>0};\bk) \cong H_{d-2n+2}(\fS_n ;  \mr{Lie}_n^* \otimes \bk).\]
In particular we have $H_{n,d}^{E_\infty}(\bN_{>0};\bk)=0$ for $d< 2(n-1)$; a range twice as large as that given by the standard connectivity estimate. Furthermore, if $\bk$ is a finite field of characteristic $p$, it vanishes unless $n$ is a power of $p$ \cite[Theorem 1.1]{AroneDwyer}.
\end{example}

\subsection{$E_1$-splitting complexes revisited}\label{sec:splittingcomplexRevis} One may develop the analogue of the results of the previous section for a monoidal groupoid $(\mathsf{G}, \oplus, \bunit)$ too, which we outline here.

In this case there is an analogous \emph{$E_1$-splitting category} $\cat{S}^{E_1}(x)$ for $x \in \sfG$:\index{$E_1$-splitting category}
\begin{itemize}\item It has objects given by triples $([n],f,\phi)$ of a finite ordered set $[n] \coloneqq \{1<2<\ldots<n\}$ with $n \geq 2$, a function $f \colon [n] \to \mr{ob}(\sfG_{r>0})$ and a morphism $\phi \colon \bigoplus_{\alpha \in [n]} f(\alpha) \to x$ in $\sfG$ (necessarily an isomorphism).
		\item A morphism $([n],f,\phi) \to ([n'],f',\phi')$ is the data of an order-preserving surjection $e \colon [n] \twoheadrightarrow [n']$ and isomorphisms $\varphi_\alpha \colon f'(\alpha) \to \bigoplus_{\beta \in e^{-1}(\alpha)} f(\beta)$ for all $\alpha \in [n']$ such that the following diagram commutes
		\[\begin{tikzcd} \bigoplus_{\alpha \in [n']} f'(\alpha) \rar{\phi'} \dar[swap]{\bigoplus_{\alpha \in [n']} \varphi_\alpha} & x \\ 
		\bigoplus_{\alpha \in [n']} \bigoplus_{\beta \in e^{-1}(\alpha)} f(\beta) \rar{\cong} & \bigoplus_{\beta \in [n]} f(\beta) \uar[swap]{\phi},\end{tikzcd}\]
		with bottom map the canonical identification of these two sums.
\end{itemize}
The following analogues of Proposition \ref{prop:DecEqualsSplittingCx} and Corollary \ref{cor:IndecEqualsSuspSplittingCx} are established by the same arguments: there are $G_x$-equivariant equivalences
\[\Dec_\bL^{E_1}(\gT)(x) \simeq B\cat{S}^{E_1}(x)_+\]
and hence
\[Q^{E_1}_\bL(\gT)(x) \simeq \Sigma B\cat{S}^{E_1}(x).\]
These results hold without Assumption \ref{assum:MonoidalInj}.

This is related to the discussion of Section \ref{sec:splittingcomplex} as follows. There is a functor
\[\cat{S}^{E_1}(x) \lra \cat{Simp}(\split_\bullet(x))\]
to the poset of simplices of the semi-simplicial set $\split_\bullet(x)$ given by sending an object $([n],f,\phi)$ to the equivalence class of the element $\phi \in \sfG(f(1) \oplus \cdots \oplus f(n),x)$ in $\colim_{x_0, \ldots, x_{p+1} \in \sfG_{r>0}^{p+2}} \sfG(x_0\oplus \cdots \oplus x_{p+1}, x) = \split_p(x)$, and the evident map on morphisms. This functor is full and essentially surjective, and under Assumption \ref{assum:MonoidalInj} it is also faithful and so is an equivalence of categories. As equivalences of categories induce homotopy equivalences on classifying spaces, and the classifying space of the poset of simplices is the barycentric subdivision, we deduce that $\split(x) \simeq B\cat{S}^{E_1}(x)$. 
\section{Application to homological stability}\label{sec:homological-stability-applications}

One of the basic consequences of the theory developed so far is a homological stability theorem\index{homological stability} for the groups $G_x$ of automorphisms in a braided monoidal groupoid $\cat{G}$ satisfying the hypotheses described in Section \ref{sec:algebras-from-groupoids}: 
\begin{enumerate}[(i)]
	\item there is a monoidal functor $r \colon \sfG \to \bN$ with $r^{-1}(0)$ the objects isomorphic to $\bunit$, 
	\item Assumption \ref{assum:Unit} that $G_\bunit$ is trivial, and 
	\item Assumption \ref{assum:MonoidalInj} that for all objects $x, y \in \cat{G}$, the homomorphism $- \oplus - \colon G_x \times G_y \to G_{x \oplus y}$ is injective. 
\end{enumerate}	
As before, $\sfC = \sfS^\sfG$ with $\sfS$ satisfying the axioms of Sections \ref{sec:axioms-of-cats} and \ref{sec:axioms-of-model-cats}.

In this section we shall deduce a generic ``homological stability'' result based on the theory developed in this monograph.  We emphasize that our main motivation for this project is to give applications ``beyond homological stability'' (as in \cite{e2cellsII,e2cellsIII,e2cellsIV}).

\subsection{A generic homological stability result}

The basic generic homological stability result is as follows. It has surprisingly many applications, and in Section \ref{sec:outlook} we shall give one. A related result regarding the relative groups $H_d(G_{\sigma^{2d}}, G_{\sigma^{2d-1}};\bZ)$ has been given by Hepworth \cite{hepworth}.

\begin{theorem}\label{thm:StdStability}
Suppose that $(\mathsf{G}, \oplus, \bunit)$ is a braided monoidal groupoid as above which satisfies the standard connectivity estimate of Definition \ref{defn:StdConnEst}, and such that there is a unique $\sigma \in \mathsf{G}$ with $r(\sigma)=1$ up to isomorphism. Then up to isomorphism the objects of $\mathsf{G}$ are precisely the powers of $\sigma$, and $H_d(G_{\sigma^{n}}, G_{\sigma^{n-1}};\bZ)=0$ for $2d \leq n-1$.

In addition, if $\mathbbm{k}$ is a commutative ring such that the map $\sigma \cdot - \colon H_1(G_{\sigma};\mathbbm{k}) \to H_1(G_{\sigma^2};\mathbbm{k})$ is surjective, then in fact $H_d(G_{\sigma^{n}}, G_{\sigma^{n-1}};\mathbbm{k})=0$ for $3d \leq 2n-1$.
\end{theorem}

We will prove this by considering the $E_2$-algebra $\gR \in \Alg_{E_2}(\cat{sSet}^\bN)$ associated to $\cat{G}$ as in Section \ref{sec:EkAlgFromGpd}. As we are only interested in the homology of $\gR$, say with coefficients in a commutative ring $\mathbbm{k}$, there is no harm in applying the free $\mathbbm{k}$-module functor to consider $\gR_\mathbbm{k} \coloneqq \mathbbm{k}\gR \in \Alg_{E_2}(\mathsf{sMod}_\mathbbm{k}^\bN)$ as an $E_2$-algebra in the category of $\bN$-graded simplicial $\mathbbm{k}$-modules. The category $\mathsf{sMod}_\mathbbm{k}$ satisfies all the axioms of Sections \ref{sec:axioms-of-cats} and \ref{sec:axioms-of-model-cats}, as well as Axiom \ref{axiom:Hurewicz}, so all the tools we have developed so far may be applied. Moreover, this construction preserves derived indecomposables:

\begin{lemma}\label{lem:indec-commutes-with-k}There is a natural weak equivalence $Q^{E_k}_\bL(\gR_\bk) \simeq \bk[Q^{E_k}_\bL(\gR)]$ in $\cat{sMod}_\bk^\bN$.
\end{lemma}

\begin{proof}
This follows from the considerations of Section \ref{sec:deriv-indec-1} applied to the Quillen adjunction $\bk[-] \dashv U$ between $\cat{sSet}^\bN$ and $\cat{sMod}_\bk^\bN$, using the fact that the left adjoint $\bk[-]$ is strong monoidal so there is a natural isomorphism $(\bk[-],\phi_{\bk[-]})_*(\gR) \cong \gR_\bk$ by Lemma \ref{lem:l-left-right-iso}.
\end{proof}

Let us take for granted for the moment the part of Theorem \ref{thm:StdStability} that says that the objects of $\sfG$ up to isomorphism are $\sigma^n$ for $n \in \bN$. Recall that in Section \ref{sec:modules} we have defined an associative unital algebra $\overline{\gR}_\mathbbm{k}$, weakly equivalent to the unitalisation $\gR^+_\mathbbm{k}$ as an $E_1^+$-algebra, and using the adapter construction of Section \ref{sec:adapter} and the element $\sigma \in H_0(G_\sigma ; \mathbbm{k}) = H_0(G_1 ; \mathbbm{k}) = \pi_{1,0}(\gR) \cong \pi_{1,0}(\overline{\gR}_\mathbbm{k})$ we have described in Section \ref{sec:ApplyingAdapters} how to form a left $\overline{\gR}_\mathbbm{k}$-module $\overline{\gR}_\mathbbm{k}/\sigma$. By definition, the underlying homotopy type of $\overline{\gR}_\mathbbm{k}/\sigma$ is that of the homotopy cofibre of the composition
\begin{equation*}
  S^{1,0} \otimes \overline{\gR}_\mathbbm{k} \xrightarrow{\sigma \otimes \mathrm{id}} \overline{\gR}_\mathbbm{k} \otimes \overline{\gR}_\mathbbm{k}
  \xrightarrow{\mu} \overline{\gR}_\mathbbm{k},
\end{equation*}
with $\mu$ the multiplication, so that we have
\[H_{n,d}(\overline{\gR}_\mathbbm{k}/\sigma) \cong H_d(G_{\sigma^{n}}, G_{\sigma^{n-1}};\mathbbm{k}).\]
Thus Theorem \ref{thm:StdStability} will be a consequence of the following general theorem for $E_k$-algebras in the category of $\bN$-graded simplicial $\mathbbm{k}$-modules, using Lemma \ref{lem:indec-commutes-with-k} to translate the hypothesis on the $E_k$-homology of $\gR$ into one about $\gR_\bk$.

\begin{theorem}\label{thm:GeneralStab}
Let $k \geq 2$ and $\gR \in \Alg_{E_k}(\mathsf{sMod}_\mathbbm{k}^\bN)$ be a non-unital $E_k$-algebra such that $H_{\ast,0}(\overline{\gR}) = \bk[\sigma]$ with $\vert\sigma \vert = (1,0)$. If $H_{n,d}^{E_k}(\gR)=0$ for $d < n-1$, then $H_{n,d}(\overline{\gR}/\sigma)=0$ for $2d \leq n-1$. 

If in addition the map $\sigma \cdot - \colon H_{1,1}(\gR) \to H_{2,1}(\gR)$ is surjective then in fact $H_{n,d}(\overline{\gR}/\sigma)=0$ for $3d \leq 2n-1$, and $H_{2,1}^{E_k}(\gR)=0$.
\end{theorem}
\begin{proof}
We will make a cumulative sequence of reductions, until we can directly compute the homology of $\overline{\gR}/\sigma$ using Cohen's computation of the homology of free $E_k$-algebras, described in Section \ref{sec:Cohen}.

\vspace{1ex}

\noindent\textbf{Reduction 1.} {\it It is enough to consider the case $\gR = \gE_k(X)$ with $X$ a finite wedge of spheres such that $H_{n,d}(X)=0$ for $d < n-1$ and $H_{1,0}(X) = \mathbbm{k} \{\sigma\}$.}

\vspace{1ex}

Firstly, the groupoid $\sfG=\bN$ and the operad $\cC_k$ satisfy the hypotheses of Lemma \ref{lem:TensorDetectNull}. Every $E_2$-algebra in $\cat{sMod}_\bk^\bN$ is 0-connective, and $\gR$ is reduced (cf.~Definition \ref{defn:Reduced}) because $H_{0,0}(\gR^+;\bk) = \bk$ implies $H_{0,0}(\gR;\bk) = 0$. The canonical map $\binit \to \gR$ is thus between 0-connective reduced $E_k$-algebras, and thus our CW approximation theorem for $E_k$-algebras applies to it, Theorem \ref{thm:MinCellStr-additive}. Its conclusion is that there is a CW approximation $\gZ \overset{\sim}\to \gR$, where $\gZ$ is a CW $E_k$-algebra built out of cells in bidegrees $(n,d)$ with $d \geq n-1$ and a single $(1,0)$-cell $\sigma \colon S^{1,0} \to \gZ$. 

The CW $E_k$-algebra $\gZ$ has a skeletal filtration $\mr{sk}(\gZ) \in \Alg_{E_k}((\mathsf{sMod}_\mathbbm{k}^\bN)^{\bZ_{\leq}})$ with $(1,0)$-cell $\sigma$ has filtration 1, which is attached along a filtered map $\sigma \colon 1_*S^{1,0} \to \mr{sk}(\gZ)$. We may form the unital associative algebra $\overline{\mr{sk}(\gZ)}$ in $(\mathsf{sMod}_\mathbbm{k}^\bN)^{\bZ_{\leq}}$ and, using the adapter construction of Section \ref{sec:modules}, taking the quotient by $\sigma \cdot -$ gives rise to a left $\overline{\mr{sk}(\gZ)}$-module $\overline{\mr{sk}(\gZ)}/\sigma$. Upon taking colimit this recovers $\overline{\gZ}$ and $\overline{\gZ}/\sigma$, so there is a spectral sequence
\[E^1_{n,p,q} = {H}_{n,p+q,p}(\grr(\overline{\mr{sk}(\gZ)}/\sigma)) \Longrightarrow H_{n,p+q}(\overline{\gZ}/\sigma),\]
whose target is isomorphic to $H_{n,p+q}(\overline{\gR}/\sigma)$. 

Recall that $\overline{\mr{sk}(\gZ)}/\sigma$ is defined to be the homotopy pushout
  \begin{equation*}
    \begin{tikzcd}
      {0_*S^{1,0} \times B(A(\mr{sk}(\gZ)),\overline{\mr{sk}(\gZ)},\overline{\mr{sk}(\gZ)})} \arrow{d}{} \arrow{r}{\sigma''} &
      {B(A(\mr{sk}(\gZ)),\overline{\mr{sk}(\gZ)},\overline{\mr{sk}(\gZ)})} \arrow{d}{}\\
      {0_*D^{1,1} \times B(A(\mr{sk}(\gZ)),\overline{\mr{sk}(\gZ)},\overline{\mr{sk}(\gZ)})} \arrow{r}{} & {\overline{\mr{sk}(\gZ)}/\sigma}      
    \end{tikzcd}
  \end{equation*}
in the category of $\overline{\mr{sk}(\gZ)}$-modules, where the top map is multiplication by the map $\sigma \colon 0_*S^{1,0} \to \overline{\mr{sk}(\gZ)}$ using the special left $\overline{\mr{sk}(\gZ)}$-module structure on the adapter $A(\mr{sk}(\gZ))$. As $\grr(-)$ commutes with pushouts (as it is a left adjoint), bar constructions (as it is symmetric monoidal and preserves thick geometric realizations), $\overline{(-)}$ (by Lemma \ref{lem:tilde-properties} (iii)), and $A(-)$ (by the construction in Section \ref{sec:constr-adapt}), we find that $\grr(\overline{\mr{sk}(\gZ)}/\sigma) \cong \overline{\grr(\mr{sk}(\gZ))}/\sigma$.

By Theorem \ref{thm:associated-graded-skeletal} the associated graded of the skeletal filtration is given by
\[\grr(\mr{sk}(\gZ)) \simeq \gE_k\left(\bigoplus_{d \geq 0} \bigoplus_{\alpha \in I_d} S^{n_\alpha, d, d}_\bk\right),\]
with $S^{n_\alpha, d, d}_\bk$ shorthand for $d_* S^{n_\alpha,d}_\bk$. This is a free $E_k$-algebra on a wedge of spheres $X = \bigoplus_{d \geq 0} \bigoplus_{\alpha \in I_d} S^{n_\alpha, d, d}_\bk$ such that $d \geq n_\alpha-1$ for all $\alpha$, and there is a single sphere $\sigma$ of degree $(1,0,0)$. While this may not be a finite wedge of spheres, it is the colimit of its finite sub-wedges, and so $\overline{\grr(\mr{sk}(\gZ))}/\sigma \cong \overline{\gE_k(X)}/\sigma$ is a colimit of $\overline{\gE_k(X')}/\sigma$'s with $X' \subset X$ a finite sub-wedge.

\vspace{1ex}

\noindent\textbf{Reduction 2.} {\it It is enough to consider the case $\mathbbm{k}=\bZ$.}

\vspace{1ex}

Let us write $- \otimes_\bZ \mathbbm{k} \colon \mathsf{Mod}_\bZ \to \mathsf{Mod}_\mathbbm{k}$ for the base-change functor, which is symmetric monoidal and commutes with colimits, and use the same notation for the induced functor between categories of $\bN$-graded simplicial modules. Writing $S^{n,d}_\mathbbm{k} \in \mathsf{sMod}^\bN_\mathbbm{k}$ for sphere objects for now, we have $S^{n,d}_\mathbbm{k} = S^{n,d}_\bZ \otimes_\bZ \mathbbm{k}$. As base-change is symmetric monoidal and commutes with colimits, we recognize 
\[\gR = \gE_k\left(\bigoplus_{d \geq 0} \bigoplus_{\alpha \in I_d} S^{n_\alpha, d}_\mathbbm{k}\right) = \gE_k\left(\bigoplus_{d \geq 0} \bigoplus_{\alpha \in I_d} S^{n_\alpha, d}_\bZ\right) \otimes_\bZ \mathbbm{k}\]
as the base-change of $\gR_\bZ \coloneqq \gE_k(\bigoplus_{d \geq 0} \bigoplus_{\alpha \in I_d} S^{n_\alpha, d}_\bZ)$. Thus we have $\overline{\gR}/\sigma = (\overline{\gR}_\bZ/\sigma) \otimes_\bZ \mathbbm{k}$, so by the universal coefficient sequence
\[0 \lra H_{n,d}(\overline{\gR}_\bZ/\sigma) \otimes_\bZ \mathbbm{k} \lra  H_{n,d}(\overline{\gR}/\sigma) \lra \mathrm{Tor}^\bZ_1(H_{n,d-1}(\overline{\gR}_\bZ/\sigma), \mathbbm{k}) \lra 0\]
it is enough to establish the result for the case $\mathbbm{k}=\bZ$.

\vspace{1ex}

\noindent\textbf{Reduction 3.} {\it It is enough to consider the cases $\mathbbm{k}=\bF_\ell$ for all primes $\ell$.}

\vspace{1ex}

Suppose $\gR = \gE_k(X) \in \Alg_{E_k}(\mathsf{sMod}_\bZ^\bN)$ with $X$ a finite wedge of spheres such that $H_{n,d}(X)=0$ for $d < n-1$ and $H_{1,0}(X) = \bZ \{\sigma\}$. Then each $H_{n,d}(\gR)$ is a finitely-generated $\bZ$-module, so the same is true of each $H_{n,d}(\overline{\gR}/\sigma)$. Thus, if $H_{n,d}(\overline{\gR}/\sigma) \neq 0$ then there is a prime number $\ell$ such that $H_{n,d}(\overline{\gR}/\sigma) \otimes_\bZ \bF_\ell \neq 0$, and so by the universal coefficient sequence
\[0 \lra H_{n,d}(\overline{\gR}/\sigma) \otimes_\bZ \bF_\ell \lra  H_{n,d}(\overline{\gR}/\sigma \otimes_\bZ \bF_\ell) \lra \mathrm{Tor}^\bZ_1(H_{n,d-1}(\overline{\gR}/\sigma), \bF_\ell) \lra 0\]
we have $H_{n,d}(\overline{\gR}/\sigma \otimes_\bZ \bF_\ell) =H_{n,d}((\overline{\gR} \otimes_\bZ \bF_\ell)/\sigma) \neq 0$. Contrapositively, if $H_{n,d}((\overline{\gR} \otimes_\bZ \bF_\ell)/\sigma) = 0$ in a range of bidegrees for all primes $\ell$, then $H_{n,d}(\overline{\gR}/\sigma) = 0$ in that range of bidegrees too.

\vspace{2ex}

So let us consider the case $\mathbbm{k}=\bF_\ell$, and $\gR = \gE_k(X)$ with $X= S^{1,0} \sigma \oplus \bigoplus_{\alpha \in I} S^{n_\alpha,d_\alpha} x_\alpha$ a finite wedge of spheres such that $d_\alpha \geq n_\alpha -1$ and $d_\alpha > 0$ for each $\alpha$. In this case we can compute the homology of $\overline{\gR}/\sigma$, using the results of F.\ Cohen summarized in Theorem \ref{thm.wkfreealg}. Namely, $H_{*,*}(\gR^+) = W_{k-1}(H_{*,*}(X))$ is a free graded commutative algebra with generators $Q^I_\ell(y)$ where $y$ is a basic Lie word on $\sigma$ and the $x_\alpha$, satisfying certain conditions. Thus $H_{*,*}(\overline{\gR}/\sigma) \cong H_{*,*}(\gR^+)/(\sigma)$ is the free graded commutative algebra with the same generators except for $\sigma$ (though it should only be considered as an $H_{*,*}(\gR^+)$-module, not as a ring). Applying $\beta Q^s_\ell$ or $Q^s_\ell$ to an element increases its slope $\frac{d}{n}$ (in principle there can be elements of infinite slope when $n=0$, which is fine) and the bracket of two elements has larger slope than the smaller of the slopes of the two elements, so a generator of minimal slope is one of: (i) $[\sigma,\sigma]$, (ii) $\beta Q^1_\ell(\sigma)$ (or $Q^1_2(\sigma)$ if $\ell=2$), or (iii) $x_\alpha$ with $\vert x_\alpha \vert = (r, r-1)$ and $r \geq 2$. These have slope (i) $\frac{k-1}{2}$, (ii) $\tfrac{2(\ell-1)-1}{\ell}$ (or $\tfrac{1}{2}$ if $\ell=2$), and (iii) $\tfrac{r-1}{r}$, so all have slope $\geq \tfrac{1}{2}$. Thus $H_{n,d}(\overline{\gR}/\sigma)=0$ for $2d < n$.

\vspace{1ex}

For the addendum we will employ essentially the above argument but using a modified filtration. We have that
\[E_k(S^{1,0}_\bZ)(2) = \mathcal{C}_k(2) \otimes_{\fS_2} (S^{0}_\bZ)^{\otimes 2} \simeq \bZ \mathrm{Sing}_\bullet(\mathcal{C}_k(2)/\fS_2) \simeq \bZ \mathrm{Sing}_\bullet(\mathbb{RP}^{k-1}).\]
We will make use of the element
\[Q^1_\bZ(\sigma) \in H_{2,1}(\gE_k(S^{1,0}_\bZ \sigma)) = 
\begin{cases}
\bZ & \text{if $k=2$,}\\
\bZ/2 & \text{if $k>2$,}
\end{cases}\]
characterized by the properties that it reduces modulo 2 to $Q^1_2(\sigma)$ and, if $k=2$, it satisfies $2Q^1_\bZ(\sigma) = -[\sigma, \sigma]$ (the sign is due to Cohen's conventions for the bracket; note that when $k=2$ we have $Q^1_2(\sigma)=\xi(\sigma)$). We make a choice of a map $Q^1_\bZ(\sigma) \colon S^{2,1}_\bZ \to \gE_k(S^{1,0}_\bZ \sigma)$ representing this class. Under  $- \otimes_\bZ \mathbbm{k}$ it yields a map $Q^1_\bk(\sigma) \colon S^{2,1}_\bk \to \gE_k(S^{1,0}_\bk \sigma)$. As the map $\sigma \colon S^{1,0}_\mathbbm{k} \to \gR$ extends to a map $\gE_k(S^{1,0}_\mathbbm{k} \sigma) \to \gR$, we obtain a map $Q^1_\mathbbm{k}(\sigma) \colon S^{2,1}_\bk \to \gR$.

Let $\{x_\alpha\}_{\alpha \in I}$ be a set of generators of the $\mathbbm{k}$-module $H_{1,1}(\gR)$. Under the assumption that ${\sigma} \cdot - \colon H_{1,1}(\gR) \to H_{2,1}(\gR)$ is surjective, there is an ${x} \in H_{1,1}(\gR)$ such that $Q^1_\mathbbm{k}({\sigma}) = {\sigma} \cdot {x}$, and we may suppose that $x$ is one of the generators $x_\alpha$. There is thus an $E_k$-map
\[\gZ_{0} \coloneqq \gE_k\left(S^{1,0}_\mathbbm{k}{\sigma} \oplus \bigoplus_{\alpha \in I}S^{1,1}_\mathbbm{k}{x}_\alpha\right) \cup^{E_k}_{Q^1_\mathbbm{k}({\sigma})- {\sigma} \cdot x} \gD^{2,2}_\mathbbm{k}\rho \lra \gR,\]
given by a choice of nullhomotopy of the map $Q^1_\mathbbm{k}({\sigma})- {\sigma} \cdot x \colon S^{2,1}_\bk \to \gR$. 

\vspace{2ex}

\noindent\textbf{Claim:} We have $H_{2,1}^{E_k}(\gR, \gZ_0)=0 = H_{2,1}^{E_k}(\gR)$.

\begin{proof}[Proof of claim]
The long exact sequence on $E_k$-homology gives
\[0=H_{2,1}^{E_k}(\gZ_0) \lra H_{2,1}^{E_k}(\gR) \lra H_{2,1}^{E_k}(\gR,\gZ_0) \lra H_{2,0}^{E_k}(\gZ_0)=0,\]
so the two vanishing statements are equivalent. It also gives
\begin{equation*}
\begin{tikzcd}
     & 0 \rar & H_{1,2}^{E_k}(\gR) \rar \ar[draw=none]{d}[name=X, anchor=center]{} & H_{1,2}^{E_k}(\gR, \gZ_0)
                \ar[rounded corners,
                to path={ -- ([xshift=2ex]\tikztostart.east)
                	|- (X.center) \tikztonodes
                	-| ([xshift=-2ex]\tikztotarget.west)
                	-- (\tikztotarget)}]{dll} \\
             & {\bigoplus\limits_{\alpha \in I}\mathbbm{k}\{x_\alpha\}} \rar[two heads] & H_{1,1}^{E_k}(\gR) \rar \ar[draw=none]{d}[name=Y, anchor=center]{} & H_{1,1}^{E_k}(\gR, \gZ_0)
                \ar[rounded corners,
                to path={ -- ([xshift=2ex]\tikztostart.east)
                	|- (Y.center) \tikztonodes
                	-| ([xshift=-2ex]\tikztotarget.west)
                	-- (\tikztotarget)}]{dll} \\
             & \mathbbm{k}\{\sigma\} \rar{\sim} & H_{1,0}^{E_k}(\gR) \rar & H_{1,0}^{E_k}(\gR, \gZ_0) \rar& 0
\end{tikzcd}
\end{equation*}
so we have 
$H_{1,0}^{E_k}(\gR, \gZ_0)=H_{1,1}^{E_k}(\gR, \gZ_0)=0$ and hence by the Hurewicz theorem (Corollary \ref{cor:absolute-Hurewicz}) we have an isomorphism
\[H_{2,1}(\gR, \gZ_0) \overset{\sim}\lra H_{2,1}^{E_k}(\gR, \gZ_0).\]
(We have used that $\cC_k(1) \simeq *$ so that $\bk[\bunit] \otimes_{H_{*,0}(\cC_k(1);\bk)} -$ is the identity.) But the diagram
\begin{equation*}
\begin{tikzcd}
	H_{2,1}(\gZ_0) \rar & H_{2,1}(\gR) \rar & H_{2,1}(\gR, \gZ_0) \rar{0} & H_{2,0}(\gZ_0) \rar{\sim} & H_{2,0}(\gR)\\
	H_{1,1}(\gZ_0) \rar[two heads] \uar{\sigma} & H_{1,1}(\gR) \uar[two heads]{\sigma}
\end{tikzcd}
\end{equation*}
shows that $H_{2,1}(\gR, \gZ_0)=0$. Thus $H_{2,1}^{E_k}(\gR, \gZ_0)=0$ as claimed. 
\end{proof}

By again applying the CW approximation theorem for $E_k$-algebras, Theorem \ref{thm:MinCellStr-additive}, we can extend the map $\gZ_0 \to \gR$ to a relative CW approximation $\gZ_0 \to \gZ \overset{\sim}\to \gR$ having no $(2,1)$-cells (as $H^{E_k}_{2,1}(\gR, \gZ_0)=0$). The associated skeletal filtration $\mr{sk}(\gZ)$ has
\[\grr(\mr{sk}(\gZ)) \simeq 0_*(\gZ_0) \vee^{E_k} \gE_k(X')\]
where $X' \in (\cat{sMod}^\bN_\bk)^{\bN_{=}}$ is a wedge of $S^{n,d,d}$'s with $d \geq n-1$, $d>0$, and $(n,d) \neq (2,1)$. Taking underlying ungraded objects, as in Reduction 1 it is enough to prove the vanishing of the homology of $(\gZ_0 \vee^{E_k} \gE_k(\colim X'))^+/\sigma$ in the appropriate range of degrees, and we may suppose without loss of generality that $\colim X'$ is a finite wedge of $S^{n,d}$'s with $d \geq n-1$, $d>0$, and $(n,d) \neq (2,1)$.

Now we observe that $\gZ_0$ is obtained by base-change along $\bZ \to \mathbbm{k}$, as
\[\left(\gE_k\left(S^{1,0}_\bZ{\sigma} \oplus \bigoplus_{\alpha \in I}S^{1,1}_\bZ{x}_\alpha\right) \cup^{E_k}_{Q^1_\bZ({\sigma})- {\sigma} \cdot x} \gD^{2,2}_\bZ\rho\right) \otimes_\bZ \mathbbm{k},\]
and $\gE_k(X')$ for $X'$ a wedge of spheres is too. Thus, as in Reduction 2, it is enough to consider the case $\mathbbm{k}=\bZ$. Finally, as in Reduction 3 it is enough to consider the case $\mathbbm{k}=\bF_\ell$ for all primes $\ell$. 

\vspace{1ex}

We therefore consider the case $\mathbbm{k}=\bF_\ell$, and $\gR = \gZ_0 \vee^{E_k} \gE_k(X')$ with $X' = \bigoplus_{\alpha \in J} S^{n_\beta, d_\beta}x_\beta$ a finite wedge of spheres such that $d_\beta \geq n_\beta-1$ and $(n_\beta, d_\beta) \neq (2,1)$ for each $\beta \in J$. If we give $\gR$ the cell attachment filtration for the cell $\rho$, the associated filtration of $\overline{\gR}/\sigma$ has spectral sequence starting with
\[E^1_{n,p,q} = H_{n, p+q, p}\left(\overline{\gE_k\left(S^{1,0,0}_{\bF_\ell}\sigma \oplus \bigoplus_{\alpha \in I}S^{1,1,0}_{\bF_\ell}{x}_\alpha \oplus S^{2,2,1}_{\bF_\ell}\rho \oplus \bigoplus_{\beta \in J} S_{\bF_\ell}^{n_\beta, d_\beta,0}x_\beta\right)}/\sigma\right)\]
and converging to $H_{n, p+q}(\overline{\gR}/\sigma)$. As in the first case, this $E^1$-page may be written as the free graded commutative algebra on a certain set of generators, except $\sigma$. The $d^1$-differential satisfies
\[d^1(\rho) = Q^1_{\bF_\ell}(\sigma) - \sigma \cdot x \equiv Q^1_{\bF_\ell}(\sigma) \mod (\sigma).\]
By our characterisation of $Q^1_{\bF_\ell}(\sigma)$ we have $Q^1_{\bF_\ell}(\sigma) = Q^1_2(\sigma)$ if $\ell=2$, and $Q^1_{\bF_\ell}(\sigma) = -\tfrac{1}{2} [\sigma, \sigma]$ if $\ell$ is odd.

\vspace{1ex}

\noindent\textbf{Claim:} $E^2_{n,p+q,p} = 0$ vanishes for $\frac{p+q}{n}<\frac{2}{3}$.

\begin{proof}[Proof of claim] The groups $E^2_{*,*,*}$ are given by the homology of the chain complex $(E^1_{*,*,*},d^1) = (\Lambda_{\bF_\ell}(L/\langle \sigma \rangle), d^1)$ where $L$ is the trigraded vector space with homogeneous basis $Q^I_\ell(y)$ for $y$ a basic Lie word in $\{\sigma, x_\alpha, \rho, x_\beta\}$. To estimate these homology groups we introduce an additional ``computational'' filtration, which has the virtue of filtering away most of the $d^1$-differential. We let $F^\bullet E^1_{*,*,*}$ be the filtration in which $Q^1_{\bF_\ell}(\sigma)$ and $\rho$ are given filtration $0$, the remaining elements of a homogeneous basis are given filtration equal to their homological degree, and this filtration is extended multiplicatively. The $d^1$-differential preserves this filtration.
	
The associated graded $\grr(F^\bullet E^1_{\ast,\ast,\ast})$ is the tensor product of chain complexes
\[\left(\Lambda_{\bF_\ell}(\bF_\ell\{Q^1_{\bF_\ell}(\sigma),\rho\}),\delta\right) \otimes \left(\Lambda_{\bF_\ell}(L/\langle \sigma,Q^1_{\bF_\ell}(\sigma),\rho \rangle),0\right)\]
where $\delta(Q^1_{\bF_\ell}(\sigma)) = 0$ and $\delta(\rho) = Q^1_{\bF_\ell}(\sigma)$. This is the first page of a spectral sequence converging to $E^2_{\ast,\ast,\ast}$.

First note that all elements of $L/\langle \sigma,Q^1_{\bF_\ell}(\sigma),\rho \rangle$ have slope $\geq \frac{2}{3}$, so the right term consist of elements of slope $\geq \frac{2}{3}$. There are now two cases, the difference owing to the fact that $(Q^1_{\bF_\ell}(\sigma))^2 = 0$ if and only if $\ell$ is odd. In the case $\ell = 2$, the left term has homology $\Lambda_{\bF_2}(\bF_2\{\rho^2\})$ in which every element has slope 1, so the second page of this spectral sequence consists of elements of slope $\geq \frac{2}{3}$. In the case $\ell$ is odd, the non-zero class of lowest slope in the homology of the left term is $Q^1_{\bF_\ell}(\sigma) \cdot \rho^{\ell-1}$, of bidegree $(2\ell, 2\ell-1)$ so slope $\frac{2\ell-1}{2\ell} \geq \frac{5}{6}$. Thus in this case the second page of this spectral sequence consists of elements of slope $\geq \frac{2}{3}$ too. 
This implies the claim.
\end{proof}

The claim says that the $E^2$-page of a spectral sequence converging to $H_{n,d}(\overline{\gR}/\sigma)$ vanishes for $\frac{d}{n}<\frac{2}{3}$, hence so does $H_{n,d}(\overline{\gR}/\sigma)$.  This finishes the proof of Theorem~\ref{thm:GeneralStab}.  
\end{proof}

\begin{proof}[Proof of Theorem \ref{thm:StdStability}]
Let $\gR \in \Alg_{E_2}(\cat{sSet}^\bN)$ be the $E_2$-algebra associated to $\cat{G}$ as in Section \ref{sec:EkAlgFromGpd}. The powers of $\sigma$ are distinct because $r \colon \sfG \to \bN$ is a monoidal functor sending $\sigma$ to $1$. For contradiction, suppose that the objects of $\mathsf{G}$ are not precisely the powers of $\sigma$, and let $x \in \mathsf{G}$ be an object of minimal rank $r(x) \geq 2$ which is not a power of $\sigma$. Such an $x$ would define a 1-simplex $\frac{G_x}{G_x}$ of $\suspsplit(x)$, but there are no non-degenerate simplices of higher dimension. Thus $\suspsplit(x) \simeq S^1$ and hence $Q^{E_1}_\bL (\gR)(r(x))$ contains $\Sigma^\infty_+ BG_x$ as a summand, which is not 0-connected and so violates the standard connectivity estimate.

Now let $\gR_\bZ \coloneqq \bZ\gR \in \Alg_{E_2}(\cat{sMod}_\bZ^\bN)$. Then $H_{n,d}^{E_1}(\gR_\bZ) = H^{E_1}_{n,d}(\gR;\bZ)=0$ for $d < n-1$ by the standard connectivity estimate. We can transfer the vanishing line from $E_1$- to $E_2$-homology using Theorem \ref{thm:TrfUp} and conclude that $H^{E_2}_{n,d}(\gR_\bZ)=0$ for $d < n-1$ too. We explained above that there is an identification $H_{n,0}(\gR_\bZ^+) = \bZ\{\sigma^n\}$, and thus Theorem \ref{thm:GeneralStab} applies.
\end{proof}

In the subsequent papers we will employ the techniques developed so far to study the homology of families of groups arising from braided monoidal groupoids (such as mapping class groups, automorphism groups of free groups, general linear groups, unitary groups, and so on), as well as other moduli spaces which can be arranged to form $E_k$-algebras (such as classifying spaces of diffeomorphism groups, configuration spaces, and so on). Homological stability, as in Theorem \ref{thm:StdStability}, of such spaces is one application, but the methods of this paper allow us to prove other types of results about the homology of such spaces. We feel that the reader, having done so much work, is due an application of these techniques: luckily, there are two which are accessible with no further theory.

\subsection{Example: general linear groups of Dedekind domains} \label{sec:outlook} 

\index{homological stability!Dedekind domains} Let $\Lambda$ be a Dedekind domain, and $(\cat{P}_\Lambda, \oplus, 0)$ denote the symmetric monoidal category of finitely-generated projective $\Lambda$-modules under direct sum. Assigning to such a $\Lambda$-module its rank defines a functor $r \colon \cat{P}_\Lambda \to \bN$, and only the zero module has rank zero. The construction of Section \ref{sec:EkAlgFromGpd} therefore defines an $E_\infty$-algebra $\gR \in \Alg_{E_\infty}(\cat{sSet}^\bN)$ such that
\[H_{n,d}(\gR;\mathbbm{k}) = \bigoplus_{\substack{[P] \in \pi_0(\cat{P}_\Lambda) \\ \mr{rk}(P)=n}} H_d(\mr{GL}(P);\mathbbm{k})\]
is the sum of the homologies of the general linear groups of all (isomorphism classes of) projective $\Lambda$-modules of rank $n$. The realization of the semi-simplicial set $\split_\bullet(P)$ from Definition \ref{def:e1-splitting} is isomorphic to the realization of the  nerve of the ``split Tits poset" $S_\Lambda(P)$ introduced by Charney \cite{Charney}, which she has shown has the homotopy type of a wedge of $(r(P)-2)$-spheres (it is here the assumption that $\Lambda$ is a Dedekind domain is used). Thus $\suspsplit(P)$ is a wedge of $r(P)$-spheres, and so $\cat{P}_\Lambda$ satisfies the standard connectivity estimate.

To apply Theorem \ref{thm:StdStability} we require that $\cat{P}_\Lambda$ has a unique (isomorphism class of) object of rank 1, i.e.\ that the Dedekind domain $\Lambda$ has class number 1, so all finitely generated projective $\Lambda$-modules are free. In this case the first part of Theorem \ref{thm:StdStability} simply recovers van der Kallen's stability range for the groups $\mr{GL}_n(\Lambda)$ \cite[Theorem 4.11]{vdK}, but the power of our method becomes apparent with the second part of Theorem \ref{thm:StdStability}. Applied to Dedekind domains satisfying $H_1(\mr{GL}_2(\Lambda), \mr{GL}_1(\Lambda);\bZ)=0$, it tells us $H_d(\mr{GL}_n(\Lambda), \mr{GL}_{n-1}(\Lambda);\bZ)=0$ for $3d \leq 2n-1$. This property holds for example for rings of integers in any number field except $\bQ(\sqrt{-d})$ for $d \neq 1,2,3,7,11$ \cite{Cohn, Vaserstein}, so such rings of class number 1 enjoy this improved stability range. A less obvious application is that if $H_1(\mr{GL}_2(\Lambda), \mr{GL}_1(\Lambda);\bZ)$ is a finite group of order $N$, then $H_d(\mr{GL}_n(\Lambda), \mr{GL}_{n-1}(\Lambda);\bZ[\tfrac{1}{N}])=0$ for $3d \leq 2n-1$. This applies for example to $\Lambda=\bZ$ and $N=2$, where it seems to be new.\index{general linear groups!Dedekind domains}

\subsection{Example: general linear groups of $\bF_q$} \label{sec:FiniteFields} 

Let us specialize the example of the previous section to $\Lambda = \bF_q$ the finite field with $q=p^m$ elements. \index{homological stability!finite fields}\index{general linear groups!finite fields}

If $\ell$ is a prime number other than $p$ then Quillen has computed the homology groups $H_*(\mr{GL}_n(\bF_q);\bF_\ell)$, and in fact in the notation of the previous section he has computed $H_{*,*}(\gR;\bF_\ell)$ as a ring. To express the answer, let $r$ be the smallest positive integer such that $q^r \equiv 1 \mod \ell$. Then  Quillen \cite[Theorem 3]{quillenfinite} shows there is an isomorphism of bigraded rings
\[H_{*,*}(\gR;\bF_\ell) \cong \bF_\ell[\sigma, \xi_1, \xi_2, \ldots] \otimes \Lambda_{\bF_\ell}[\eta_1, \eta_2, \ldots]\]
where $\sigma$ has bidegree $(1,0)$, $\xi_i$ has bidegree $(r, 2ir)$, and $\eta_i$ has bidegree $(r, 2ir-1)$. Homological stability in this case is evident: one has $H_{*,*}(\gR/\sigma;\bF_\ell) \cong \bF_\ell[\xi_1, \xi_2, \ldots] \otimes \Lambda_{\bF_\ell}[\eta_1, \eta_2, \ldots]$ which vanishes in bidegrees $(n,d)$ with $\tfrac{d}{n} < 2-\tfrac{1}{r}$.

On the other hand when $\ell=p$ the homology groups $H_*(\mr{GL}_n(\bF_q);\bF_p)$ are not yet known. However, Quillen has shown \cite[Theorem 6]{quillenfinite} that $H_{n,d}(\gR;\bF_p)=0$ for $0 < d < m(p-1)$, and this vanishing range has been improved by Friedlander--Parshall \cite[Lemma A.1]{FriedlanderParshall} to $0 < d < m(2p-3)$. The free $E_\infty$-algebra $\gE_\infty(S^{1,0}\sigma)$ on a generator $\sigma$ of bidegree $(1,0)$ has $\bF_p$-homology class of smallest positive degree given by $\beta Q^1(\sigma)$ (or $Q^1(\sigma)$ if $p=2$) of degree $2p-3$. Thus the natural map $\gE_\infty(S^{1,0}\sigma) \to \gR$ is an isomorphism on $\bF_p$-homology in homological degrees $* < 2p-3$. Combining this with the vanishing line for $E_1$-, and hence $E_\infty$-, homology of $\gR$ in the previous section, we find that $H_{n,d}^{E_\infty}(\gR, \gE_\infty(S^{1,0}\sigma);\bF_p)=0$ for $\tfrac{d}{n} < \tfrac{2p-3}{2p-2}$. Consulting Cohen's calculations we have $H_{n,d}(\gE_\infty(S^{1,0}\sigma)/\sigma ; \bF_p)=0$ for $\tfrac{d}{n} < \tfrac{2p-3}{p}$, so using
\[\overline{\gR}/\sigma \simeq B(\overline{\gE_\infty(S^{1,0}\sigma)}/\sigma, \overline{\gE_\infty(S^{1,0}\sigma)}, \overline{\gR})\]
and Corollary \ref{cor:BarEstimate} we obtain $H_{n,d}(\gR/\sigma;\bF_p)=0$ for $\tfrac{d}{n} < \tfrac{2p-3}{2p-2}$. This establishes homological stability for $H_d(\mr{GL}_n(\bF_q);\bF_p)$ with slope $\tfrac{2p-3}{2p-2}$. In fact, Quillen has shown that for $q \neq 2$ these groups have homological stability with slope 1 \cite[1974-I, p.\ 10]{quillennotes}:  in \cite{e2cellsIII} we explain how our methods can be used to prove this result, which was also obtained by Sprehn and Wahl \cite[Theorem A]{SprehnWahl}.

\section{Local coefficients}\label{sec:local-coefficients} In this section we explain how to obtain results for local coefficients analogous to the generic homological stability results obtained in the previous section. In particular, we again assume that $\sfG$ satisfies the assumptions (i)--(iii) described in the beginning of Section \ref{sec:algebras-from-groupoids}.

\subsection{Coefficient systems}\label{sec:CoeffSys}

In the framework we have described in Sections \ref{sec:algebras-from-groupoids} and \ref{sec:homological-stability-applications}, it is easy to discuss homology of the collection of groups $G_x$ with coefficients in a collection of $\mathbbm{k}[G_x]$-modules. In fact it is no more difficult, and is technically convenient, to discuss hyperhomology of the groups $G_x$ with coefficients in simplicial $\mathbbm{k}[G_x]$-modules.

Let $\underline{\mathbbm{k}} \coloneqq \mathbbm{k}\underline{*} \in \cat{sMod}_\mathbbm{k}^\cat{G}$ denote the free $\mathbbm{k}$-module on $\underline{*} \in \cat{sSet}^\cat{G}$. As $\underline{*}$ has the structure of a unital commutative monoid, so does $\underline{\mathbbm{k}}$.

\begin{definition}
A \emph{coefficient system}\index{coefficient system} for $\sfG$ is a left $\underline{\mathbbm{k}}$-module, i.e.\ a functor $\gA\in \cat{sMod}_\mathbbm{k}^\cat{G}$ equipped with maps
\[\mu_{a,b} \colon \gA(b) \cong \underline{\mathbbm{k}}(a) \otimes_{\mathbbm{k}} \gA(b) \lra \gA(a \oplus b)\]
of simplicial $\mathbbm{k}$-modules which are appropriately associative.

We say a coefficient system is \emph{discrete} if each simplicial $\mathbbm{k}$-module $\gA(x)$ is actually just a $\mathbbm{k}$-module.
\end{definition}

As before, we take a cofibrant approximation $\gT \overset{\sim}\to \underline{*}_{>0}$ of non-unital $E_k$-algebras. This gives in particular a map $\overline{\gT} \to \overline{\underline{*}_{>0}}$ of unital associative algebras, and as $\underline{*}_{>0}$ is strictly associative the latter has a map of unital associative algebras to $\underline{*}$. Taking the associated simplicial $\bk$-modules this gives a map $\overline{\gT}_\mathbbm{k} \to \underline{\mathbbm{k}}$ of associative algebras, and so any coefficient system is also a left $\overline{\gT}_\mathbbm{k}$-module.

We may therefore find a cofibrant approximation $c\gA \overset{\sim}\to \gA$ as a left $\overline{\gT}_\mathbbm{k}$-module, which will be cofibrant in $\cat{sMod}_\mathbbm{k}^\cat{G}$, and hence the derived left Kan extension $\gR_\gA \coloneqq r_*(c\gA) \simeq \bL r_*(\gA)$ has the structure of a left $\overline{\gR}_\mathbbm{k} = r_*(\overline{\gT}_\mathbbm{k})$-module. Furthermore, by definition of the (derived) left Kan extension one has
\[H_{n,d}(\gR_\gA) = \bigoplus_{\substack{[x] \in \pi_0(\sfG)\\r(x)=n}} \bH_d(G_x ; \gA(x)),\]
the direct sum over the isomorphism classes of rank $n$ objects $x$ of the hyperhomologies of the groups $G_x$ with coefficients in the simplicial $\mathbbm{k}$-modules $\gA(x)$.

In particular for each $\sigma \in H_{1,0}(\gR_\mathbbm{k}) = H_{1,0}(\overline{\gR}_\mathbbm{k})$ there is a left multiplication map $\sigma \cdot - \colon S^{1,0} \otimes \gR_\gA \to \gR_\gA$, inducing a map
\[\sigma \cdot - \colon \bigoplus_{\substack{y \in \cat{G}\\r(x)=n-1}} \bH_d(G_y ; \gA(y)) \lra \bigoplus_{\substack{x \in \cat{G}\\r(x)=n}} \bH_d(G_x ; \gA(x)),\]
on homology, and one may ask for stability with respect to this map. Using the adapter as in Section \ref{sec:ApplyingAdapters}, this left multiplication map is equivalent to a map of left $\overline{\gR}_\mathbbm{k}$-modules, with homotopy cofibre $\gR_\gA/\sigma$, and so the homology stability of this map is equivalent to the vanishing of $H_{n,d}(\gR_\gA/\sigma)$ in a range of bidegrees. Our general result in this direction establishes stability in terms of the $\overline{\gR}_\mathbbm{k}$-module homology of $\gR_\gA$.

\begin{theorem}\label{thm:CoeffSysStab}
Let $\gA$ be a coefficient system for $\sfG$, and suppose that there are $\lambda \leq 1$ and $c$ such that $H_{n,d}(\overline{\gR}_\mathbbm{k}/\sigma)=0$ for $d < \lambda n$ and
$H_{n, d}^{\overline{\gR}_\mathbbm{k}}(\gR_\gA)=0$ for $d < \lambda(n-c)$. Then $H_{n,d}(\gR_\gA/\sigma)=0$ for $d < \lambda(n-c)$.
\end{theorem}
\begin{proof}
We will apply the theory of CW approximations developed in Section \ref{sec:additive-case} in the setting of left $\overline{\gR}_\mathbbm{k}$-modules in $\cat{sMod}_\bk^\bN$. As the groupoid $\bN$ is Artinian, and the operad $\cO$ which models left $\overline{\gR}_\mathbbm{k}$-modules has $H_{0, 0}(\cO(1)) = H_{0,0}(\overline{\gR}_\mathbbm{k}) = \bk$, because $r^{-1}(0)$ consists of objects isomorphic to $\bunit_\sfG$, condition (ii) of Lemma \ref{lem:TensorDetectNull} is satisfied. As before all simplicial $\bk$-modules are $0$-connective, so the canonical morphism $\binit \to \gR_\gA$ is between 0-connective $\cO$-algebras. Thus by Theorem \ref{thm:MinCellStr-additive} we may construct a CW approximation $\gZ \overset{\sim}\to \gR_\gA$, where $\gZ$ consists of $(n,d)$-cells with $d \geq \lambda(n-c)$. It has skeletal filtration $\mr{sk}(\gZ) \in \Alg_{\cO}((\mathsf{sMod}_\mathbbm{k}^\bN)^{\bZ_{\leq}})$ and by Theorem \ref{thm:associated-graded-skeletal} its associated graded is given by
\[\grr(\mr{sk}(\gZ)) \cong \bigvee_{\alpha \in I} S^{n_\alpha,d_\alpha, d_\alpha} \wedge \overline{\gR}_\mathbbm{k}\]
with $d_\alpha \geq \lambda(n_\alpha-c)$ for each $\alpha$. 

The filtered object $B(0_*(\overline{\gR}_\mathbbm{k}/\sigma), 0_*(\overline{\gR}_\mathbbm{k}), \mr{sk}(\gZ))$ has colimit $B(\overline{\gR}_\mathbbm{k}/\sigma, \overline{\gR}_\mathbbm{k}, \gZ) \simeq B(\overline{\gR}_\mathbbm{k}/\sigma, \overline{\gR}_\mathbbm{k}, \gR_\gA) \simeq \gR_\gA/\sigma$ and associated graded 
\[B(0_*(\overline{\gR}_\mathbbm{k}/\sigma), 0_*(\overline{\gR}_\mathbbm{k}), \grr(\mr{sk}(\gZ))) \simeq \bigvee_{\alpha \in I} S^{n_\alpha,d_\alpha, d_\alpha} \wedge \overline{\gR}_\mathbbm{k}/\sigma,\]
so gives a spectral sequence
\[\bigoplus_{\alpha \in I} H_{*,*, *}(S^{n_\alpha,d_\alpha,d_\alpha} \wedge \overline{\gR}_\mathbbm{k}/\sigma) \Longrightarrow H_{*,*}(\gR_\gA/\sigma).\]
By assumption $H_{n,d}(\overline{\gR}_\mathbbm{k}/\sigma)=0$ for $d < \lambda n$, and so $H_{n,d}(S^{n_\alpha,d_\alpha} \wedge \overline{\gR}_\mathbbm{k}/\sigma)=0$ for $d-d_\alpha < \lambda (n-n_\alpha)$. Thus the target of this spectral sequence vanishes in bidegrees $(n,d)$ satisfying
\[d < \lambda n + \min_{\alpha \in I}(d_\alpha - \lambda n_\alpha),\]
so satisfying $d < \lambda n -\lambda c$ as required.
\end{proof}

Once the principle behind this proof is understood, other qualitative (and quantitative) results suggest themselves. For example, under the weaker hypothesis
\[H_{n,d}^{\overline{\gR}_\mathbbm{k}}(\gR_\gA)=0 \text{ for } n \gg d\]
the same argument allows one to conclude that $H_{n,d}(\gR_\gA/\sigma)=0$ for $n \gg d$. 

\begin{remark}\label{rem:StabImpliesModuleCells}
The following argument shows that all coefficient systems $\gA$ which are known to enjoy homological stability do in fact have a vanishing line for $H_{n,d}^{\overline{\gR}_\mathbbm{k}}(\gR_\gA)$; in this sense Theorem \ref{thm:CoeffSysStab} is optimal.

Any $\overline{\gR}_\mathbbm{k}$-module $\gM$ may be descendingly filtered by its $\bN$-grading, giving an associated graded $\grr(\gM)$ which is isomorphic to $\gM$ in $\cat{sMod}_\mathbbm{k}^\bN$ but which has the trivial $\overline{\gR}_\mathbbm{k}$-module structure induced by the augmentation $\epsilon \colon \overline{\gR}_\mathbbm{k} \to \mathbbm{k}$. This induces a filtration of $B(\mathbbm{k} , \overline{\gR}_\mathbbm{k} , \gM)$ with associated graded $B(\mathbbm{k} , \overline{\gR}_\mathbbm{k} , \mathbbm{k}) \otimes \grr(\gM)$ and, suppressing the internal grading, a spectral sequence
\[E^1_{*,*} = (\mathbbm{k} \oplus \Sigma^{0,1}H^{E_1}_{*,*}(\gR)) \otimes H_{*,*}(\gM) \Longrightarrow H_{*,*}^{\overline{\gR}_\mathbbm{k}}(\gM).\]
Suppose that $\gA$ is a coefficient system which is known to enjoy homological stability, i.e.\ $H_{n, d}(\gR_\gA/\sigma)=0$ for $d < \lambda(n-c)$ for some $\lambda \leq 1$ and some $c$, and apply the above spectral sequence to the $\overline{\gR}_\mathbbm{k}$-module $\gR_\gA/\sigma$. Assuming that $\cat{G}$ satisfies the standard connectivity estimate, it follows that $E^1_{n, d}=0$ for $d < \lambda(n-c)$, and hence $\smash{H_{n,d}^{\overline{\gR}_\mathbbm{k}}}(\gR_\gA/\sigma)$ vanishes in this range too. But the left $\overline{\gR}_\mathbbm{k}$-module map $S^{1,0} \otimes \gR_\gA \to \gR_\gA$ constructed using the adapter is nullhomotopic on $\overline{\gR}_\mathbbm{k}$-module derived indecomposables, so
\[H_{n,d}^{\overline{\gR}_\mathbbm{k}}(\gR_\gA/\sigma) \cong H_{n,d}^{\overline{\gR}_\mathbbm{k}}(\gR_\gA) \oplus H_{n-1,d-1}^{\overline{\gR}_\mathbbm{k}}(\gR_\gA)\]
and hence $H_{n,d}^{\overline{\gR}_\mathbbm{k}}(\gR_\gA)=0$ for $d < \lambda(n-c)$.
\end{remark}

Applying Theorem \ref{thm:CoeffSysStab} requires an effective way of computing the homology groups $H_{n,d}^{\overline{\gR}_\mathbbm{k}}(\gR_\gA)$, or at least of proving their vanishing. We shall give one method to do so. 

If $\gA$ is discrete then 
\[H_{*,d}(c\gA) = \begin{cases}
\gA &\text{if }d=0,\\
0 & \text{else},
\end{cases} \qquad \text{and} \qquad H_{*,d}(\overline{\gT}_\mathbbm{k}) = \begin{cases}
\underline{\mathbbm{k}}&\text{if }d=0,\\
0 & \text{else}.
\end{cases}
\]
Now note that under Assumption \ref{assum:MonoidalInj} we may apply Lemma \ref{lem:KunnethFormula} (i). As $\underline{\mathbbm{k}}$ is objectwise flat, the K{\"u}nneth formula gives
\[H_{*,d}(B_p(\mathbbm{k} , \overline{\gT}_\mathbbm{k}, c\gA)) = H_{*,d}(\mathbbm{k} \otimes \overline{\gT}_\mathbbm{k}^{\otimes p} \otimes c\gA) \cong \begin{cases}
\mathbbm{k} \otimes\underline{\mathbbm{k}}^{\otimes p} \otimes \gA&\text{if }d=0,\\
0 & \text{else.}
\end{cases}\]
Hence the bar spectral sequence for $B(\mathbbm{k} , \overline{\gT}_\mathbbm{k}, c\gA)$ takes the form
\[E^2_{x,d, p} = \mathrm{Tor}^{\underline{\mathbbm{k}}}_{p, d}(\mathbbm{k} , \gA)(x) \Longrightarrow H_{x,d+p}(B(\mathbbm{k} , \overline{\gT}_\mathbbm{k}, c\gA))\]
and is supported along the line $d=0$ so collapses at $E^2$. 
By the equivalences
$Q_\bL^{\overline{\gR}_\mathbbm{k}}(\gR_\gA) \simeq B(\mathbbm{k} , \overline{\gR}_\mathbbm{k} , \gR_\gA) = r_* B(\mathbbm{k} , \overline{\gT}_\mathbbm{k}, c\gA)$ there is a spectral sequence
\begin{equation}\label{eq:ComputeModHomology} E^1_{n,p,q} = \bigoplus_{\substack{[x] \in \pi_0(\sfG) \\r(x)=n}} H_p(G_x ; \mathrm{Tor}^{\underline{\mathbbm{k}}}_{q}(\mathbbm{k} , \gA)(x)) \Longrightarrow H_{n, p+q}^{\overline{\gR}_\mathbbm{k}}(\gR_\gA),
\end{equation}
from which the following lemma is immediate.

\begin{lemma}\label{lem:ComputeModHomology}
If $\gA$ is a discrete coefficient system and $\lambda, \mu \in \bZ$ are such that for each $x \in \cat{G}$ we have $\mathrm{Tor}^{\underline{\mathbbm{k}}}_{d}(\mathbbm{k} , \gA)(x)=0$ for $d \leq \lambda r(x)+\mu$, then $H_{n, d}^{\overline{\gR}_\mathbbm{k}}(\gR_\gA)=0$ for $d \leq \lambda n+\mu$.
\end{lemma}

\begin{remark}\label{rem:CentralStabCx}
If $\cat{G}$ has objects $\bN$ and $\gA$ is a discrete coefficient system, considered as a $\overline{\gT}_\mathbbm{k}$-module, then $c\gA$ can be filtered as in Remark \ref{rem:StabImpliesModuleCells} by rank. Its associated graded $\mr{gr}(c\gA)$ is isomorphic to $c\gA$ in $\cat{sMod}_\bk^\sfG$ but has trivial module structure. We obtain from it a descending filtration of $B(\mathbbm{k} , \overline{\gT}_\mathbbm{k}, c\gA)$ with associated graded $B(\mathbbm{k} , \overline{\gT}_\mathbbm{k}, \mathbbm{k}) \owedge \mr{gr}(c\gA)$.

If $\cat{G}$ satisfies the standard connectivity estimate then
$B(\mathbbm{k} , \overline{\gT}_\mathbbm{k}, \mathbbm{k})(n)$ is a wedge of $n$-spheres, and its $n$th reduced homology is by definition the $E_1$-Steinberg module ${St}^{E_1}(n)$. Thus the spectral sequence associated to the above filtration has
\[E^1_{n,n, q} = \mathrm{Ind}_{G_{n+q} \times G_{-q}}^{G_n}({St}^{E_1}(n+q) \otimes \gA(-q)) \text{ for } q \leq 0\]
and all other terms zero: hence it collapses at $E^2$. This gives a chain complex $C_{n,*}(\gA) \coloneqq (E^1_{n,n,*-n}, d^1)$ computing $\mathrm{Tor}^{\underline{\mathbbm{k}}}_{*}(\mathbbm{k} , \gA)(n)$, quite different to that given by the bar resolution.
\end{remark}

\begin{example}
If $\cat{G}=\cat{FB}$ is the category of finite sets and bijections, then a discrete $\underline{\mathbbm{k}}$-module is precisely the datum of an FI-module in the sense of \cite{CEF}, and $\mathrm{Tor}^{\underline{\mathbbm{k}}}_{d}(\mathbbm{k} , -)(S)$ is precisely FI-homology.\index{FI-homology} By the Noetherian property of the category of FI-modules (when $\mathbbm{k}$ is Noetherian) \cite{CEFN}, if $\gA$ is a FI-module which is finitely-generated, then it is equivalent to a cellular $\underline{\mathbbm{k}}$-module with finitely-many cells of each dimension, so it satisfies $\mathrm{Tor}^{\underline{\mathbbm{k}}}_{d}(\mathbbm{k} , \gA)(S)=0$ for $\vert S \vert \gg d$. In fact, by the Castelnuovo--Mumford regularity property of the category of FI-modules \cite{ChurchEllenberg} a finitely-generated FI-module $\gA$ satisfies $\mathrm{Tor}^{\underline{\mathbbm{k}}}_{d}(\mathbbm{k} , \gA)(S)=0$ for $d \leq \vert S \vert +c$ for a constant $c$ (which may be determined in terms of the $d=0$ and $d=1$ pieces).

In this case the $E_1$-Steinberg modules are all given by the sign representation, and the chain complex $C_*(\gA)$ from Remark \ref{rem:CentralStabCx} is Putman's ``central stability chain complex" \cite{PutmanCong}, which is known to compute FI-homology.
\end{example}

\subsection{Polynomial coefficients}\label{sec:Polynomial}

There is a common source of examples of coefficient systems for which the results of the previous section can be applied, namely tensor powers of ``linear" coefficient systems.

In this section we will make use of the \emph{objectwise} tensor product $A \boxtimes B$ of functors $A, B \colon \cat{G} \to \cat{sMod}_\mathbbm{k}$, given by
\[(A \boxtimes B)(x) \coloneqq A(x) \otimes_\mathbbm{k} B(x).\]
If $A$ is a coefficient system and $S$ is a finite set, let us write $A^{\boxtimes S}$ for the objectwise $S$th tensor power of $A$, i.e.\
\[A^{\boxtimes S}(a) = A(a)^{\otimes_\mathbbm{k} S}.\]
The maps $\mu_{a,b}^{\otimes S} \colon A(b)^{\otimes S} \to A(a \oplus b)^{\otimes S}$ make $A^{\boxtimes S}$ into a coefficient system. Furthermore, the symmetric group $\fS_S$ acts on $A^{\boxtimes S}$ by maps of coefficient systems.

\begin{definition}\label{defn:LinCoeffSys}
Let $(\cat{G}, \oplus, \bunit)$ be a braided monoidal groupoid. A \emph{linear coefficient system}\index{coefficient system!linear} $\gL = (L, s)$ on $\cat{G}$ is a functor $L \colon \cat{G} \to \cat{Mod}_\mathbbm{k}$ equipped with a strong braided monoidality
\[s_{a,b} \colon L(a) \oplus L(b) \lra L(a \oplus b)\]
with respect to direct-sum of $\mathbbm{k}$-modules. This implies that $L(\bunit)=0$, and that $L(\beta_{b,a} \circ \beta_{a,b}) = \mathrm{Id}$. This yields a coefficient system in the usual sense by considering $\mathbbm{k}$-modules as discrete simplicial $\mathbbm{k}$-modules, and with the $\underline{\mathbbm{k}}$-module structure given by
\[\mu_{a,b} \colon \gL(b) = 0 \oplus \gL(b) \lra \gL(a) \oplus \gL(b) \overset{s_{a,b}}\lra \gL(a \oplus b).\]
\end{definition}

In Remark \ref{rem:lin-coeff-braided} we give an example showing that it is necessary for our purposes that $L$ is braided.

Our goal is to establish a vanishing line for $H^{\underline{\mathbbm{k}}}_{x,d}(\gL^{\boxtimes S})$ when $\gL$ is linear and $S$ is a finite set. By the ``multiplicative" philosophy we are expounding, it is advantageous to collect the coefficient systems $\gL^{\boxtimes S}$ into a single multiplicative object, as follows. Consider the functor
\begin{align*}
\gL^\boxtimes \colon \cat{G} \times \cat{FB} &\lra \cat{Mod}_\mathbbm{k} \subset \cat{sMod}_\mathbbm{k}\\
(x, S) & \longmapsto \begin{cases}
0 & \text{if } (x,S) \cong (\bunit, \varnothing)\\
 \gL(x)^{\otimes S}  & \text{else}.
\end{cases}
\end{align*}
Let us write $\underline{\mathbbm{k}} = \mathbbm{k} \oplus \underline{\mathbbm{k}}_{>0} \colon \cat{G} \to \cat{sMod}_\mathbbm{k}$, so $\underline{\mathbbm{k}}_{>0}$ is a nonunital commutative monoid, and consider the functor 
\[\pi_\cat{G}^*\underline{\mathbbm{k}}_{>0} \colon \cat{G} \times \cat{FB} \overset{\pi_\cat{G}}\lra \cat{G} \overset{\underline{\mathbbm{k}}_{>0}}\lra \cat{sMod}_\mathbbm{k}.\]
We may then write, tautologically, $\gL^\boxtimes = (\pi_\cat{G}^*\underline{\mathbbm{k}}_{>0}) \boxtimes \gL^\boxtimes$. 

We may calculate
\begingroup
\allowdisplaybreaks
\begin{align*}
(\gL^\boxtimes \otimes \gL^\boxtimes)(x,S) &= \colim_{\substack{a, b \in \cat{G}\\r(a), r(b)>0\\a \oplus b \to x}} \colim_{\substack{A, B \in \cat{FB}\\A \sqcup B \to S}}  \gL(a)^{\otimes A} \otimes  \gL(b)^{\otimes B},\\
\intertext{now note that holding $a$ and $b$ fixed, the innermost left Kan extension along $\sqcup \colon \cat{FB} \times \cat{FB} \to \cat{FB}$ is $S \mapsto (\gL(a) \oplus \gL(b))^{\otimes S}$ (a categorification of the Binomial Theorem) and hence isomorphic to $\gL(a \oplus b)^{\otimes S}$ via $s_{a,b}^{\otimes S}$, so}
 &= \colim_{\substack{a, b \in \cat{G}\\r(a), r(b)>0\\a \oplus b \to x}} \gL(a \oplus b)^{\otimes S}\\
&= \left(\colim_{\substack{a, b \in \cat{G}\\r(a), r(b)>0\\a \oplus b \to x}} \mathbbm{k}\right) \otimes_\mathbbm{k} \gL(x)^{\otimes S}\\
 &= (\underline{\mathbbm{k}}_{>0} \otimes \underline{\mathbbm{k}}_{>0})(x) \otimes_\mathbbm{k} \gL(x)^{\otimes S}.
\end{align*}
\endgroup
The conclusion is that there is an isomorphism $\gL^\boxtimes \otimes \gL^\boxtimes \cong ((\pi_\cat{G}^*\underline{\mathbbm{k}}_{>0}) \otimes (\pi_\cat{G}^*\underline{\mathbbm{k}}_{>0})) \boxtimes \gL^\boxtimes$; analogously we have isomorphisms $(\gL^\boxtimes)^{\otimes p} \cong (\pi_\cat{G}^*\underline{\mathbbm{k}}_{>0})^{\otimes p} \boxtimes \gL^\boxtimes$. In particular, the multiplication map $\mu_{\underline{\mathbbm{k}}_{>0}} \colon \underline{\mathbbm{k}}_{>0} \otimes \underline{\mathbbm{k}}_{>0} \to \underline{\mathbbm{k}}_{>0}$ defines a multiplication map \[\mu_{\gL^\boxtimes} \colon \gL^\boxtimes \otimes \gL^\boxtimes \lra \gL^\boxtimes\] making $\gL^\boxtimes$ into an associative algebra in $\cat{sMod}_\mathbbm{k}^{\cat{G} \times \cat{FB}}$. As $\underline{\mathbbm{k}}_{>0}$ is in fact commutative (in the sense that $\mu_{\underline{\mathbbm{k}}_{>0}} \circ \beta_{\underline{\mathbbm{k}}_{>0}, \underline{\mathbbm{k}}_{>0}} = \mu_{\underline{\mathbbm{k}}_{>0}}$, which makes sense even though $\cat{sMod}_\mathbbm{k}^{\cat{G} \times \cat{FB}}$ is only braided monoidal), and $\gL$ is a braided monoidal functor, one can check that this multiplication makes $\gL^\boxtimes$ into a non-unital commutative algebra, and therefore into an $E_2$-algebra.

\begin{theorem}\label{thm:VanLineLBox}
We have $S^1 \wedge Q^{E_1}_\bL(\gL^\boxtimes)(x,S) \simeq S^1 \wedge  Q^{E_1}_\bL(\underline{\mathbbm{k}}_{>0})(x) \otimes \gL(x)^{\otimes S}$. In particular, if $\cat{G}$ satisfies the standard connectivity estimate then
\[H^{E_1}_{x, S, d}(\gL^\boxtimes)=0 \text{ for } d < r(x)-1.\]
\end{theorem}
\begin{proof}
We may choose a cofibrant approximation $\gT_\gL \overset{\sim}\to \gL^\boxtimes$ as an $E_1$-algebra. As in the proof of Proposition \ref{prop:YCxModel}, there are equivalences
\[S^1 \wedge Q^{E_1}_\bL(\gL^\boxtimes) \simeq S^1 \wedge Q^{E_1}_\bL(\gT_\gL) \simeq \tilde{B}^{E_1}(\gT_\gL),\]
with $\tilde{B}^{E_1}_p(\gT_\gL)$ the quotient of $\mathcal{P}_1(p) \times (\gT_\gL^+)^{\otimes p}$ by the subobject consisting of units. By the same analysis as that which showed $(\gL^\boxtimes)^{\otimes p} \cong (\pi^*_\cat{G}\underline{\mathbbm{k}}_{>0})^{\otimes p} \boxtimes \gL^\boxtimes$, there are equivalences
\[(\gT_\gL^+)^{\otimes p} \overset{\sim}\lra (\pi^*_\cat{G}(\gT_\mathbbm{k})^+)^{\otimes p} \boxtimes \gT_\gL^+,\]
giving a map $\tilde{B}^{E_1}_\bullet(\gT_\gL) \to \pi_\sfG^*(\tilde{B}^{E_1}_\bullet(\gT_\bk)) \boxtimes \gT_\gL^+$ of semi-simplicial objects which is a levelwise weak equivalence. As both objects are levelwise cofibrant, and $-\boxtimes-$ commutes with geometric realization in each entry, this gives a weak equivalence between geometric realizations.
\end{proof}

\begin{corollary}\label{cor:PolyModuleCells}
Suppose $\cat{G}$ satisfies the standard connectivity estimate.
\begin{enumerate}[(i)]
\item If $\sfG$ is symmetric monoidal, then
\[\mathrm{Tor}^{\underline{\mathbbm{k}}}_{d}(\mathbbm{k}, \gL^{\boxtimes S})(x)=0 \text{ for } d < r(x)-\vert S \vert.\]

\item If $\sfG$ is  braided monoidal, then
\[H_{n,d}^{\overline{\gR}_\bk}(\gR_{\gL^{\boxtimes S}})=0 \text{ for } d < n-\vert S \vert.\]

\end{enumerate}
\end{corollary}

\begin{proof}
If $\sfG$ is symmetric monoidal then (ii) follows from (i) by Lemma \ref{lem:ComputeModHomology}. We will explain the proof of (i), then explain the modifications necessary to prove (ii).

Suppose that $\sfG$ is symmetric monoidal. The statement has no content if $r(x)=0$, so we may assume that $x \not\cong \bunit_\sfG$. Consider the $E_2$-algebra map $f \colon \pi_\sfG^*\underline{\mathbbm{k}}_{>0} \otimes \pi_\cat{FB}^*\mathbbm{k} \to \gL^{\boxtimes}$, given by the identity maps $\mathbbm{k} \to \gL(x)^{\otimes S}$ when $S = \varnothing$ and $r(x)>0$. The unit object $\mathbbm{k} \in \cat{sMod}_\mathbbm{k}^\cat{FB}$ is cofibrant, though not as an $E_2$-algebra. We choose cofibrant approximations $\gT_\mathbbm{k} \overset{\sim}\to \underline{\mathbbm{k}}_{>0}$ and $\gT_\gL \overset{\sim}\to \gL^\boxtimes$ as $E_2$-algebras, and can therefore lift $f$ to a map
\[F \colon \pi_\sfG^*\gT_\mathbbm{k} \otimes \pi_\cat{FB}^*\mathbbm{k} \lra \gT_\gL\] of $E_2$-algebras which are cofibrant in $\cat{sMod}_\mathbbm{k}^{\sfG \times \cat{FB}}$.

We will apply Corollary \ref{cor:BarEstimate} to the $E_2$-algebra map $F$, using the lax monoidal abstract connectivity $\rho(x, S) = \sup(0, r(x)-\vert S \vert)$. Note that $(\pi_\sfG^*\gT_\mathbbm{k} \otimes \pi_\cat{FB}^*\mathbbm{k})^+ \simeq \pi_\sfG^*(\gT_\mathbbm{k}^+) \otimes \pi_\cat{FB}^*\mathbbm{k}$ so we have
\[\bk \oplus \Sigma Q^{E_1}_\bL(\pi_\sfG^*\gT_\mathbbm{k} \otimes \pi_\cat{FB}^*\mathbbm{k}) \simeq B^{E_1}(\pi_\sfG^*\gT_\mathbbm{k}^+ \otimes \pi_\cat{FB}^*\mathbbm{k}) \simeq \pi_\sfG^*B^{E_1}(\gT_\mathbbm{k}^+) \otimes \pi_\cat{FB}^*B^{E_1}(\mathbbm{k}).\]
As $\bk \in \cat{sMod}^\cat{FB}_\bk$ is the unit, $B^{E_1}(\mathbbm{k}) \simeq \bk$ by Lemma \ref{lem:bek-unit}. Furthermore $B^{E_1}(\gT_\mathbbm{k}^+) \simeq \bk \oplus \Sigma Q^{E_1}_\bL(\gT_\mathbbm{k})$, so
\[\Sigma Q^{E_1}_\bL(\pi_\sfG^*\gT_\mathbbm{k} \otimes \pi_\cat{FB}^*\mathbbm{k}) \simeq \pi_\sfG^*(\Sigma Q^{E_1}_\bL(\gT_\mathbbm{k})) \otimes \pi_\cat{FB}^*(\bk).\]

We have assumed that $\cat{G}$ satisfies the standard connectivity estimate, $H^{E_1}_{d,x}(\gT_\mathbbm{k})=0$ for $d < r(x)-1$, and so $H^{E_1}_{x,S,d}(\pi_\sfG^*\gT_\mathbbm{k} \times \pi_\cat{FB}^*\mathbbm{k})=0$ unless $S=\varnothing$ and $d \geq r(x)-1$. In particular, it vanishes for $d < \rho(x, S)-1$. Thus by Theorem \ref{thm:TrfUp} the $E_2$-homology vanishes in the same range, which verifies hypothesis (i) of Corollary \ref{cor:BarEstimate}. 

To verify hypothesis (ii), note that the map $F$ is an equivalence evaluated at any $(x, \varnothing)$, so the map
\[H^{E_2}_{x, \varnothing, d}(\pi_\sfG^*\gT_\mathbbm{k} \times \pi_\cat{FB}^*\mathbbm{k}) \lra H^{E_2}_{x, \varnothing, d}(\gT_\gL)\]
is an isomorphism, and so $H^{E_2}_{x, \varnothing, d}(\gT_\gL, \pi_\sfG^*\gT_\mathbbm{k} \times \pi_\cat{FB}^*\mathbbm{k})=0$. On the other hand, for $S \neq \varnothing$ we have $H^{E_1}_{x,S,d}(\pi_\sfG^*\gT_\mathbbm{k} \times \pi_\cat{FB}^*\mathbbm{k}) = 0$ and so
\[H^{E_2}_{x, S, d}(\gL^\boxtimes) \overset{\sim}\lra H^{E_2}_{x, S, d}(\gT_\gL, \pi_\sfG^*\gT_\mathbbm{k} \times \pi_\cat{FB}^*\mathbbm{k}).\]
By Theorem \ref{thm:VanLineLBox} we have that $H^{E_1}_{x, S, d}(\gT_\gL)=H^{E_1}_{x, S, d}(\gL^\boxtimes)=0$ for $d < r(x)-1$, and so for $d < r(x) - \vert S \vert$ as $S \neq \varnothing$; these groups always vanish in negative degrees, so they vanish for $d < \rho(x, S) = \sup(0, r(x)-\vert S \vert)$. By Theorem \ref{thm:TrfUp} (applied using the abstract connectivity $\rho+1$, which is also lax monoidal), the same vanishing holds for $E_2$-homology of $\gT_\gL$, and so for $H^{E_2}_{x, S, d}(\gT_\gL, \pi_\sfG^*\gT_\mathbbm{k} \times \pi_\cat{FB}^*\mathbbm{k})$ as required.

By Corollary \ref{cor:BarEstimate} we conclude that 
\[H_{x,S,d}^{\overline{\pi_\sfG^*\gT_\mathbbm{k} \times \pi_\cat{FB}^*\mathbbm{k}}}(\overline{\gT_\gL})=0 \text{ for $d < \rho(x,S)$},\]
(recall that we have assumed that $x \not\cong \bunit_\sfG$). Finally, as
\[H_{*,d}(\overline{\pi_\sfG^*\gT_\mathbbm{k} \times \pi_\cat{FB}^*\mathbbm{k}}) = \begin{cases}
\pi_\sfG^*\underline{\mathbbm{k}} \times \pi_\cat{FB}^*\mathbbm{k}& \text{if $d=0$,}\\
0 & \text{otherwise,}
\end{cases}\]
is objectwise flat in $\cat{sMod}_\mathbbm{k}^{\sfG \times \cat{FB}}$, we may apply the K{\"u}nneth formula of Lemma \ref{lem:KunnethFormula} (i) to identify the $E^2$-page of the bar spectral sequence for
\[Q_\bL^{\overline{\pi_\sfG^*\gT_\mathbbm{k} \times \pi_\cat{FB}^*\mathbbm{k}}}(\overline{\gT_\gL}) \simeq B(\mathbbm{k}, \overline{\pi_\sfG^*\gT_\mathbbm{k} \times \pi_\cat{FB}^*\mathbbm{k}} , \overline{\gT_\gL})\]
as $E^2_{x, S,p,q} = \mathrm{Tor}_{p}^{\underline{\mathbbm{k}}}(\mathbbm{k}, \gL^{\boxtimes S})(x)$ for $q=0$ and zero otherwise. Thus the spectral sequence collapses, showing that 
\[\mathrm{Tor}_{d}^{\underline{\mathbbm{k}}}(\mathbbm{k}, \gL^{\boxtimes S})(x) \cong H_{x,S,d}^{\overline{\pi_\sfG^*\gT_\mathbbm{k} \times \pi_\cat{FB}^*\mathbbm{k}}}(\overline{\gT_\gL}),\]
from which the result follows.

Now, if $\sfG$ is only braided monoidal then we cannot apply Corollary \ref{cor:BarEstimate} to the $E_2$-algebra map $F$ as the category $\cat{sMod}_\bk^{\sfG \times \cat{FB}}$ is only braided monoidal, which is not sufficient for Corollary \ref{cor:BarEstimate}. However, forming the Kan extension along $r \times \mr{Id} \colon \sfG \times \cat{FB} \to \bN \times \cat{FB}$ gives an $E_2$-algebra map
\[F' \colon \pi_\bN^*(\gR_\bk) \otimes \pi^*_\cat{FB}(\bk) \lra \gR_\gL\]
in the symmetric monoidal category $\cat{sMod}_\bk^{\bN \times \cat{FB}}$. By applying the change-of-diagram-category spectral sequence from Section \ref{sec:change-of-diagram-ss} one checks that the above estimates descend to this map, so Corollary \ref{cor:BarEstimate} applies to it and shows that $H_{n,S,d}^{\overline{\pi_\sfG^*\gR_\mathbbm{k} \times \pi_\cat{FB}^*\mathbbm{k}}}(\overline{\gR_\gL})=0$ for $d < \sup(0, n-\vert S \vert)$. As above this translates into the required vanishing range.
\end{proof}

We do not know whether Corollary \ref{cor:PolyModuleCells} (i) holds when $\sfG$ is only braided monoidal, as Corollary \ref{cor:BarEstimate} does not apply. It would be interesting to know whether it is true nonetheless.

\begin{corollary}
If $\cat{G}$ is symmetric monoidal, has objects $\mathbb{N}$, and satisfies the standard connectivity estimate, and if $\gL$ is a linear coefficient system, then
\[\mathrm{Tor}^{\underline{\mathbbm{k}}}_{d}(\mathbbm{k}, \gL)(n)=0 \text{ for } d \neq n-1.\]
\end{corollary}
\begin{proof}
Remark \ref{rem:CentralStabCx} gives a chain complex $C_{n,d}(\gL)$ computing these groups. Now the fact that $C_{n,d}(\gL)=0$ for $d \geq n$ (as $\gL(0)=0$) and Corollary \ref{cor:PolyModuleCells} imply the result.
\end{proof}

\begin{remark}
There are more general notions of ``polynomial coefficient systems", cf.\ \cite{DwyerTwisted}, \cite{DVFaiblement}, \cite[\S 4.4]{RWW}. One would expect an analogue of Corollary \ref{cor:PolyModuleCells} for these, expressing a vanishing line for the $\overline{\gT}_\mathbbm{k}$-module homology of $\gA$ in terms of its degree, but we were not able to find a satisfactory such analogue. However, such coefficient systems are generally known to enjoy homological stability, and so by Remark \ref{rem:StabImpliesModuleCells} the $\overline{\gR}_\mathbbm{k}$-module homology of $\gR_\gA$ does generally have a vanishing line.
\end{remark}

\begin{remark}\label{rem:lin-coeff-braided} It is perhaps not obvious why we demanded that the functor $L$ is braided monoidal in Definition \ref{defn:LinCoeffSys}: let us give an example to show why this condition is necessary.
	
	Suppose that $\sfG=\bN$ with (symmetric) monoidal structure given by sum and the trivial braiding. Choose a $A \in \cat{Mod}_\bk$ and define $L \colon \bN \to \cat{Mod}_\bk$ by $L(n) = A^{\oplus n}$. This admits a strong monoidality with respect to direct sum of $\bk$-modules, but it is not braided: $L(\beta_{1,1}) \colon L(2) \to L(2)$ is trivial, as $\beta_{1,1}$ is, but $\beta_{A, A} \colon A \oplus A \to A \oplus A$ is not. One may compute that $\mr{Tor}^{\underline{\bk}}_0(\bk, \gL)(n) \cong A$ for all $n>0$, so the conclusion of Corollary \ref{cor:PolyModuleCells} (i) does not hold. The step that fails is that the multiplication on $\gL^\boxtimes$ is not commutative.
\end{remark}

\subsection{Example: general linear groups of $\bF_q$ with local coefficients} \label{sec:FiniteFieldsTwisted}   \index{homological stability!finite fields, with local coefficients}\index{general linear groups!finite fields}

Let us consider the example of Section \ref{sec:FiniteFields} with local coefficients. Taking $\Lambda = \bF_q$ there is a linear coefficient system
\[\gV \colon \cat{P}(\bF_q) \lra \cat{Mod}_{\bF_q}\]
given by sending an $\bF_q$-module to itself. The associated $\overline{\gR}_{\bF_q}$-module $\gR_{\gV^{\boxtimes S}}$ has
\[H_{n,d}(\gR_{\gV^{\boxtimes S}}) \cong H_d(\mr{GL}_n(\bF_q), (\bF_q^n)^{\otimes_{\bF_q} S}),\]
so a vanishing of the homology of $\gR_{\gV^{\boxtimes S}}/\sigma$ corresponds to homological stability with these local coefficients. Now we have $H_{n,d}^{\overline{\gR}_{\bF_q}}(\gR_{\gV^{\boxtimes S}})=0$ for $d < n - |S|$ by Corollary \ref{cor:PolyModuleCells} (ii), so applying Theorem \ref{thm:CoeffSysStab} with the estimate $H_{n,d}(\gR/\sigma;\bF_p)=0$ for $\tfrac{d}{n} < \tfrac{2p-3}{2p-2}$ of Section \ref{sec:FiniteFields} gives
\[H_{n,d}(\gR_{\gV^{\boxtimes S}}/\sigma)=0 \text{ for } d < \tfrac{2p-3}{2p-2}(n-|S|).\]

\section{Koszulity and connectivity}
\label{sec:koszul}

In this section we shall specialize to $\sfS = \cat{sMod}_\bk$ for a
commutative ring $\bk$ and $\sfC = \cat{sMod}_\bk^{\sfG}$ for an Artinian monoidal category $\sfG$, as in Definition \ref{defn:Artinian}. In particular, there is a monoidal rank functor $r \colon \sf{G} \to \bN$ such that $r(x)>0$ if $x$ is not $\oplus$-invertible. Our goal is explain how, for a quadratic algebra, the standard connectivity estimate of Definition~\ref{defn:StdConnEst} is equivalent to this algebra having the Koszul property. While this material will not be logically necessary for the applications we intend, it puts our work in perspective and explains how it relates to work of other authors, and may be clarifying for some readers.  The definitions and presentation in this section have been adapted to our setting from \cite{LodayVallette}. 

\subsection{Quadratic data}
\label{sec:definitions} In this section we associate to a quadratic datum as below both a quadratic algebra and coalgebra.

\begin{definition}
  A \emph{quadratic datum}\index{quadratic!datum} in $\sfC$ is a pair $(V,R)$, where $V \colon \sfG \to \Modk$ is a functor with $V(x) = 0$ if $r(x)=0$, regarded as a (simplicially constant) object of $\sfC$, and $R \subset V \otimes V$ is a subfunctor.
\end{definition}

\subsubsection{Quadratic algebras}
 To define the quadratic algebra, let us in this section write $A(V) \coloneqq \bigoplus_{n=1}^\infty V^{\otimes n}$ for the free associative non-unital algebra in $\smash{\Modk^\sfG}$ generated by $V$, and $\mr{Ass}$ for the non-unital associative operad (consistent with the use of $\mr{Ass}^+$ for the unital associative operad).

\begin{definition}
  The \emph{quadratic algebra}\index{quadratic!algebra} presented by a quadratic datum $(V,R)$ is the quotient $A(V) \to A(V,R)$, terminal among algebra homomorphisms which vanish when restricted to $R \subset V^{\otimes 2} \subset A(V)$. More explicitly, the vector space $A(V,R)(x)$ is defined by the exact sequence
  \begin{equation*}
    \bigoplus_{n=1}^{\infty} \bigoplus_{i=1}^{n-1} \left(V^{\otimes i-1} \otimes R \otimes V^{\otimes (n-1-i)}\right)(x) \lra 
    \bigoplus_{n=1}^{\infty} V^{\otimes n}(x) \lra A(V,R)(x) \lra 0.
  \end{equation*}
\end{definition}

It may be described as a pushout in the category of non-unital associative algebras
\begin{equation*}
  \begin{tikzcd}
    F^\mathrm{Ass}(R)\arrow{r}\arrow{d} & F^\mathrm{Ass}(V) \arrow{d}\\
    0 \arrow{r} & A(V,R),
  \end{tikzcd}
\end{equation*}
but we emphasize that this need not be a homotopy pushout, since the homotopy pushout would likely have non-trivial higher homotopy groups.

\subsubsection{Quadratic coalgebras}
To define the quadratic coalgebra, let us next recall the \emph{deconcatenation coproduct} on the object $C(V) = \oplus_{n=1}^\infty V^{\otimes n} \in \Modk^{\sfG}$.  For $a,b \geq 1$ and $a + b = n$, let
\begin{equation*}
  \Delta_{a,b} \colon V^{\otimes n} \overset{\cong}\lra (V^{\otimes a}) \otimes (V^{\otimes b}) \lra C(V) \otimes C(V)
\end{equation*}
be the canonical isomorphism composed with the maps induced from inclusions of the $a$th and $b$th direct summands.  Let $\Delta_n = \Delta_{1,n-1} + \dots + \Delta_{n-1,1} \colon V^{\otimes n} \to C(V) \otimes C(V)$ and assemble to a morphism $\Delta \colon C(V) \to C(V) \otimes C(V)$ from the infinite direct sum (i.e.\ coproduct). As usual, this makes $(C(V),\Delta)$ into an associative non-unital coalgebra. The assumption that $\sfG$ is Artinian implies that the canonical morphism $C(V) \to \prod_{n=1}^{\infty} V^{\otimes n}$ is an isomorphism, from which it is easily verified that the projection $C(V) \to V$ to the $n=1$ summand makes $C(V)$ into the \emph{cofree coalgebra}: it has the universal property that coalgebra maps into $C(V)$ are in natural bijection with linear maps into $V$.

\begin{definition}
  The \emph{quadratic coalgebra}\index{quadratic!coalgebra} of a quadratic datum $(V,R)$ is the subcoalgebra $C(V,R) \subset C(V)$ of the deconcatenation coalgebra, terminal among subalgebras whose projection to $V^{\otimes 2}/R$ vanishes.  Explicitly, it is given by the exact sequence
  \begin{equation*}
    \bigoplus_{n=1}^{\infty} \bigoplus_{i=1}^{n-1} (V^{\otimes i-1} \otimes (V^{\otimes 2}/R) \otimes V^{\otimes (n-1-i)})(x) \longleftarrow 
    \bigoplus_{n=1}^{\infty} V^{\otimes n}(x) \longleftarrow C(V,R)(x) \longleftarrow 0.
  \end{equation*}
\end{definition}

The fact that the (co)relations defining $A(V,R)$ and $C(V,R)$ are homogeneous implies that both objects may be lifted to $\bN$-graded objects, i.e.\ functors $\bN \times \sfG \to \Modk$, as follows.

\begin{definition}
  Let $sV \coloneqq 1_* V \in \Modk^{\bN \times \sfG} \subset \sfC^{\bN}$ and $s^2 R \coloneqq 2_* R \subset sV \otimes_{\sfC^\bN} sV$.  We obtain algebras $A(sV,s^2 R)$ and coalgebras $C(sV,s^2 R)$ in $\Modk^{\bN \times \sfG}$.  The underlying ungraded (co)algebras, i.e.\ those obtained by left Kan extension along the strong monoidal functor $\bN \to \ast$, are canonically isomorphic to $A(V,R)$ and $C(V,R)$.
\end{definition}

\subsection{Duality and the $E_1$-homology of quadratic algebras}
\label{sec:duality}

By Theorem \ref{thm:BarHomologyIndec}, up to a suspension the derived $E_1$-indecomposables of $\gA \in \Alg_{E_1}(\sfC)$ may be calculated using the $E_1$-bar construction. We now restrict to the special case and $\sfC = \cat{sMod}_\bk^\sfG$ and $U^{E_1}\gA = A \in \Modk^\sfG \subset \sfC$ is an (objectwise) discrete simplicial object. For $\gA = A(V,R)$  quadratic, we shall explain how to compute $H_{x,d}^{E_1}(\gA)$ for $d \geq r(x)-1$ in terms of $C(sV,s^2R)$.

If $A$ is discrete, then the endomorphism operad of $A$ is also discrete and so the $E_1$-algebra structure descends to an associative algebra structure. Thus the $E_1$-bar construction may be replaced by the standard bar construction model for $\bk \otimes^\mathbb{L}_{\gA^+} \bk$. Concretely, this is given by the $\bk$-linear chain complex in $\Modk^{\sfG}$
\begin{equation*}
   \cdots \xrightarrow{\partial} A ^{\otimes n} \xrightarrow{\partial} \cdots \xrightarrow{\partial} A ^{\otimes 2} \xrightarrow{\partial} A \to 0,
\end{equation*}
where we have removed the extra copy of $\bk$ in degree 0 to get a model for the suspended derived $E_1$-indecomposables. After forgetting the differential, the underlying $\bN$-graded object in $\Modk^{\sfG}$ is isomorphic to $C(sA)$.  Let us write $BA = (C(sA),\partial)$ for this chain complex; then $H^{E_1}_{g,d-1}(\gA) \cong \pi_d((BA)(g))$ for any such $\gA$.

Let us now assume that $A = A(V,R)$ comes from a quadratic datum $(V,R)$.  Let us also assume that $V(x) = 0$ unless $r(x) = 1$, and that $R(x) = 0$ unless $r(x) = 2$. Then we may deduce that the bar construction above is supported in homological degrees $\leq r(x)$, and hence that $H^{E_1}_{x,d}(A) = 0$ if $d \geq r(x)$. It also implies that $C(sV,s^2R)$ is concentrated in bidegrees $(x,d) \in \sfG \times \bN$ with $d = r(x)$.

The inclusion $V \to A$ of the generators as a direct summand induces split injections $V^{\otimes n} \to A^{\otimes n}$ for all $n \geq 1$ and in turn a split injection
\begin{equation}\label{eq:31}
  C(sV) \lra C(sA),
\end{equation}
as a map of functors $\bN \times \sfG \to \Modk$.  Upon identifying it with $BA$, the right hand side comes with a boundary homomorphism which makes it a model for the suspended derived indecomposables, but~\eqref{eq:31} is not a chain map (when the domain is given the trivial boundary homomorphism). However, as in \cite[Proposition 3.3.2]{LodayVallette}, one shows that it becomes one when restricted to $C(sV,s^2R)$.
\begin{proposition}\label{prop:Koszul}
  For $\sfG$, $r \colon \sfG \to \bN$, and $(V,R)$ a quadratic datum such that $V(x) = 0$ unless $r(x) = 1$, and that $R(x) = 0$ unless $r(x) = 2$, let $\gA = A(V,R) \in \Modk^{\sfG} \subset \sfC$.  Then the induced map
  \begin{equation}\label{eq:31a}
    C(sV,s^2 R) \lra BA
  \end{equation}
  is a chain map when the domain is given the trivial boundary homomorphism, and 
  induces an isomorphism
  \begin{equation}\label{eq:31b}
    C(sV,s^2R)(x,d) \lra H^{E_1}_{x,d-1}(\gA)
  \end{equation}
  for $d = r(x)$.
\end{proposition}

Hence, the $E_1$-homology of a quadratic algebra $A(V,R)$ in functors $\sfG \to \Modk \subset \cat{sMod}_\bk$ is completely known in all bidegrees $(x,d)$ with $d \geq r(x) - 1$: it vanishes for $d \geq r(x)$ and for $d = r(x) - 1$ is canonically isomorphic to the coalgebra presented by the same quadratic data. There may be non-zero $E_1$-homology in bidegrees $(x,d)$ with $d<r(x)-1$.

The (non-counital) coassociative coalgebra $\gC = C(sV,s^2 R)$ has
underlying functor $\sfG \to \cat{Mod}_\bk^{\bN} \subset \cat{sMod}_\bk$, where we
identify $\bN$-graded $\bk$-modules with chain complexes having trivial
boundary map, regarded as a full subcategory of simplicial
$\bk$-modules by the Dold--Kan functor.  By construction, $C(x)$ is
connected for all $x$, i.e.\ $\pi_0(C(x)) = 0$.  Hence the objectwise
desuspension as an $\bN$-graded vector space is again a functor
$s^{-1} C \colon \sfG \to \cat{Mod}_\bk^{\bN} \subset \cat{sMod}_\bk$.  (In
simplicial terms, we take based loops on the based Kan complexes
$C(g)$, resulting again in a simplicial $\bk$-module).  Now form the free associative algebra
\begin{equation*}
  \Omega(\gC) \coloneqq \bigoplus_{n=0}^\infty (s^{-1} C)^{\otimes n},
\end{equation*}
and define a degree-decreasing derivation $\partial \colon \Omega(\gC) \to \Omega(\gC)$ on the generating summand $(s^{-1}C)$ as the coproduct $(s^{-1}C) \to (s^{-1}C) \otimes (s^{-1}C)$ of the coalgebra (the desuspension makes this have degree $-1$).  The pair $(\Omega(\gC),\partial)$ may be regarded as a functor from $\sfG$ to $\bk$-linear chain complexes or, by applying the Dold--Kan functor, as a functor $\sfG \to \cat{sMod}_\bk$.  The fact that the Dold--Kan functor is lax monoidal ensures that the resulting functor
\begin{equation*}
  \Omega(\gC) \colon \sfG \lra \cat{sMod}_\bk
\end{equation*}
inherits the structure of an associative algebra.  This is the \emph{cobar construction} of $\gC$.

Dually to the homomorphism~\eqref{eq:31a} of coalgebras, we have a canonical morphism of associative (non-unital) algebras in $\cat{sMod}_\bk^\sfG$
\begin{equation}\label{eq:32}
  \Omega(\gC) \lra \gA.
\end{equation}
The following is a version of Koszul duality for algebras:\index{Koszul duality}

\begin{proposition}\label{prop:koszul-duality} 
  Let $(V,R)$ be a quadratic datum such that $V(x) = 0$ unless $r(x) = 1$, and that $R(x) = 0$ unless $r(x) = 2$, and $\gA = A(V,R)$ be the associated quadratic algebra and $\gC = C(sV,s^2 R)$ be the corresponding shifted coalgebra. Assume that $H^{E_1}_{x,d}(\gA) = 0$ for $d<r(x)-1$, then the morphism~(\ref{eq:32}) is a quasi-isomorphism.
\end{proposition}

\begin{proof}[Proof sketch]
This may be verified after taking derived $E_1$-indecomposables.  It is an easy exercise to verify that the derived $E_1$-indecomposables of $\Omega(\gC)$ is just $\gC$, concentrated in bidegrees $(x,d)$ with $d = r(x)-1$.  We have already seen that the assumptions imply that this is also the derived $E_1$-indecomposables of $\gA$, and one verifies that the map is isomorphic to the identity map of $\gC$.  See also \cite[Theorems 3.4.6 and 7.4.6]{LodayVallette}.
\end{proof}

\subsection{The fundamental example}\label{sec:fundamental-koszul-example}

	Let $\sfG$ be a monoidal groupoid with object set $\bN$ satisfying Assumptions \ref{assum:Unit} and \ref{assum:MonoidalInj}, and $r \colon \sfG \to \bN$ be the identity on objects. Recall that the automorphism groups of its object are denoted $G_x \coloneqq \mr{Aut}_\sfG(x) = \sfG(x,x)$.

	Let $\bk$ be a commutative ring and $\underline{\bk}_{>0}$ be the algebra object with $\underline{\bk}_{>0}(0) = 0$, $\underline{\bk}_{>0}(x) = \bk$ for $x > 0$, and algebra structure maps $\underline{\bk}_{>0}(x) \otimes_\bk \underline{\bk}_{>0}(x') \to \underline{\bk}_{>0}(x \oplus x')$ given by the multiplication in $\bk$. 
	
	Let $V(1) = \bk$ be the trivial representation of the group $G_1$. Then $(V \otimes V)(2)$ is the permutation representation associated to the action of $G_2$ on $G_2/(G_1 \times G_1)$, and hence comes with a map to the trivial representation by collapsing $G_2/(G_1 \times G_1)$ to a point. We define a $G_2$-representation $R(2)$ by the short exact sequence
	\begin{equation*}
	0 \lra R(2) \lra (V \otimes V)(2) \lra \bk \lra 0.
	\end{equation*}
	This gives quadratic data presenting a non-unital algebra $\gA = A(V,R) \colon \sfG \to \cat{Mod}_\bk$. Then the identity map $V(1) = \bk = \ul{\bk}_{>0}(1)$ extends to a homomorphism of associative algebras $\gA \to \ul{\bk}_{>0}$. 	
	
	Let us assume that $T^{E_1}(x)$ as in Definition \ref{def:e1-splitting}, is connected for $x \geq 3$. This is equivalent to the groups $G_x$ being generated by the image of the $(x-1)$-many maps $G_2 \to G_x$ obtained by applying the functors $(-) \oplus \id_1$ and $\id_1 \oplus (-)$ (note that the \emph{relations} of the groups $G_x$ play no direct role). Then the derived $E_1$-indecomposables of $\ul{\bk}_{>0}$ vanish in bidegree $(x,1)$ for $x > 2$ and no further relations are needed to present $\ul{\bk}_{>0}$. Hence $\gA \to \ul{\bk}_{>0}$ is an isomorphism of functors $\sfG \to \cat{Mod}_k$. 
	
	In this situation Proposition \ref{prop:Koszul} tells us how to compute $H^{E_1}_{x,d}(\ul{\bk}_{>0})$ for $d \geq r(x)-1$. If in addition $\sfG$ satisfies the standard connectivity estimate, the $E_1$-decomposables vanish for $d < r(x)-1$ and as in Proposition \ref{prop:koszul-duality} we have
	\[H^{E_1}_{x,d}(\ul{\bk}_{>0}) = \begin{cases} C(sV,s^2R)(x,d) & \text{if $d = r(x)-1$.} \\
	0 & \text{otherwise.}\end{cases}\]
	As we shall see in sequels to this paper, this applies to groupoids arising from mapping class groups, general linear groups, and more.\index{standard connectivity estimate!Koszul duality}
	
\begin{remark}
In this example the quasi-isomorphism~(\ref{eq:32}) spells out to just an acyclic chain complex
\begin{equation*}
	0 \lra \bk \lra C(x) \lra (C \otimes C)(x) \lra (C \otimes C \otimes C)(x)\lra \cdots
\end{equation*}
for all $x$ with $r(x) > 0$. We would like to compare this with the spectral sequence associated to the canonical multiplicative filtration, which was described in Theorem \ref{thm:CanMultFiltSS}, in the case that $\cO = E_1$. Applied to a cofibrant approximation $\gT \overset{\sim}\to \underline{\bk}_{>0}$ in $\Alg_{E_1}(\cat{sMod}^\sfG_\bk)$, and using that $\cC_1(1) \simeq *$ so that absolute and relative $E_1$-indecomposables are weakly equivalent, this takes the form
	\[E^1_{x,p,q} = H_{x,p+q,q}(\gE_1((-1)_* Q^{E_1}(\gT))) \Longrightarrow H_{x,p+q}(\gT).\]
Now $H_{x,d}(Q^{E_1}(\gT)) = H^{E_1}_{x,d}(\underline{\bk}_{>0})$ so if $\sfG$ satisfies the standard connectivity estimate then this vanishes for $d \neq r(x)-1$ and is, by definition, the $E_1$-Steinberg module $St^{E_1}(x)$ for $d=r(x)-1$. Including the additional grading of $(-1)_*Q^{E_1}(\gT)$ we have that $H_{x,p+q, q}((-1)_*Q^{E_1}(\gT))=0$ unless $q=-1$ and $p = r(x)$. But then as
	\[\gE_1((-1)_* Q^{E_1}(\gT)) \simeq \bigoplus_{n=1}^\infty ((-1)_* Q^{E_1}(\gT))^{\otimes n}\]
by the K{\"u}nneth theorem we have that $E^1_{x,p,q}=0$ for $p \neq r(x)$ or $q \geq 0$, and
	\[E^1_{x, r(x), q} \cong (St^{E_1})^{\otimes -q}(x)\]
for $q < 0$. The $d^1$-differential has the form $d^1 \colon (St^{E_1})^{\otimes -q}(x) \to (St^{E_1})^{\otimes -q+1}(x)$. There can be no higher differentials, for reasons of grading, so the spectral sequence collapses at $E^2$. We have that $H_{x,d}(\gT)=0$ for $d>0$ and is $\bk$ for $d=0$ (and $r(x) \neq 0$), which gives an acyclic chain complex
	\[0 \lra \bk \lra (St^{E_1})(x) \overset{d^1}\lra (St^{E_1})^{\otimes 2}(x) \overset{d^1}\lra (St^{E_1})^{\otimes 3}(x) \overset{d^1}\lra \cdots\]
for all $x$ with $r(x)>0$. Under the identification $C(x) \overset{\sim}\to H^{E_1}_{x, r(g)-1}(\underline{\bk}_{>0}) = St^{E_1}(x)$ of Proposition \ref{prop:Koszul} it seems inevitable that this acyclic complex is that of (\ref{eq:32}). We expect that this can be proved using the work of Ching--Harper \cite{ChingHarper}.
\end{remark}

\bibliographystyle{amsalpha}
\bibliography{biblio}

\newpage

\appendix
\printglossary[title=\listofsymbolsname,type=symbols,style=long4col,nogroupskip]

\printindex

\end{document}